\documentclass[12pt]{article}
\usepackage{fullpage}
\usepackage{epic}
\usepackage{eepic}
\usepackage[leqno]{amsmath}
\usepackage{amssymb,mathrsfs,bm}
\usepackage{hyperref}
\usepackage[all]{xy}
\usepackage{enumerate}
\usepackage{accents}

\usepackage{soul}
\usepackage{dsfont}
\newcommand{\vertiii}[1]{{\left\vert\kern-0.25ex\left\vert\kern-0.25ex\left\vert #1 
    \right\vert\kern-0.25ex\right\vert\kern-0.25ex\right\vert}}
\usepackage{tikz}
\usetikzlibrary{matrix,arrows}

\newcommand\souligner{}
\DeclareRobustCommand\souligner[1]{\underline{#1}}

\DeclareMathOperator{\Hom}{Hom}
\DeclareMathOperator{\Tr}{Trace}
\DeclareMathOperator{\Tra}{Tr}
\DeclareMathOperator{\No}{Norm}
\DeclareMathOperator{\Temp}{Temp}
\DeclareMathOperator{\Diff}{Diff}
\DeclareMathOperator{\Ima}{Im}
\DeclareMathOperator{\Ker}{Ker}
\DeclareMathOperator{\ad}{ad}
\DeclareMathOperator{\Ad}{Ad}
\DeclareMathOperator{\End}{End}
\DeclareMathOperator{\Lie}{Lie}
\DeclareMathOperator{\Nil}{Nil}
\DeclareMathOperator{\vol}{vol}
\DeclareMathOperator{\Supp}{Supp}
\DeclareMathOperator{\Gal}{Gal}
\DeclareMathOperator{\scusp}{scusp}
\DeclareMathOperator{\g}{\mathfrak{g}}
\DeclareMathOperator{\geom}{geom}
\DeclareMathOperator{\spec}{spec}
\DeclareMathOperator{\qc}{qc}
\DeclareMathOperator{\reg}{reg}
\DeclareMathOperator{\elli}{ell}
\DeclareMathOperator{\sreg}{sreg}
\DeclareMathOperator{\der}{der}
\DeclareMathOperator{\ssi}{ss}
\DeclareMathOperator{\tempe}{temp}
\DeclareMathOperator{\mini}{min}
\DeclareMathOperator{\quasid}{qd}
\DeclareMathOperator{\stab}{stab}
\DeclareMathOperator{\auxi}{aux}

\numberwithin{equation}{subsection}
\newcounter{keepeqno}
\newenvironment{num}
 {\setcounter{keepeqno}{\value{equation}}%
  \begin{list}{(\theequation)}{\usecounter{equation}}%
  \setcounter{equation}{\value{keepeqno}}}
 {\end{list}}

\usepackage[index]{glossaries}
\usepackage{glossary-mcols}
\renewcommand*{\indexname}{List of Notations}
\makenoidxglossaries

\newterm[name={\ensuremath{F}}]{F}
\newterm[name={\ensuremath{\delta(G)}}]{dg}
\newterm[name={\ensuremath{Z_G}}]{zg}
\newterm[name={\ensuremath{A_G}}]{Ag}
\newterm[name={\ensuremath{\mathcal{A}_G}}]{cAg}
\newterm[name={\ensuremath{\mathcal{A}^*_G}}]{cA*g}
\newterm[name={\ensuremath{H_G}}]{Hg}
\newterm[name={\ensuremath{\mathcal{A}_{G,F}}}]{cAgf}
\newterm[name={\ensuremath{\widetilde{\mathcal{A}}_{G,F}}}]{ctAgf}
\newterm[name={\ensuremath{\mathcal{A}^\vee_{G,F}}}]{cAvgf}
\newterm[name={\ensuremath{\widetilde{\mathcal{A}}^\vee_{G,F}}}]{ctAvgf}
\newterm[name={\ensuremath{\mathcal{A}^*_{G,F}}}]{cA*gf}
\newterm[name={\ensuremath{\mathcal{A}_{G,\mathbb{C}}}}]{cAgc}
\newterm[name={\ensuremath{\mathcal{A}^*_{G,\mathbb{C}}}}]{cA*gc}
\newterm[name={\ensuremath{R(A_M,P)}}]{RAMP}
\newterm[name={\ensuremath{H_P}}]{HP}
\newterm[name={\ensuremath{\mathcal{P}(M)}}]{PM}
\newterm[name={\ensuremath{\mathcal{L}(M)}}]{LM}
\newterm[name={\ensuremath{\mathcal{F}(M)}}]{FM}
\newterm[name={\ensuremath{\mathcal{A}_M^L}}]{AML}
\newterm[name={\ensuremath{Z_G(x)}}]{ZGx}
\newterm[name={\ensuremath{Z_G(X)}}]{ZGX}
\newterm[name={\ensuremath{G_x}}]{Gx}
\newterm[name={\ensuremath{G_X}}]{GX}
\newterm[name={\ensuremath{G_{\ssi}}}]{Gss}
\newterm[name={\ensuremath{G_{\reg}}}]{Greg}
\newterm[name={\ensuremath{\mathfrak{g}_{\ssi}}}]{gss}
\newterm[name={\ensuremath{\mathfrak{g}_{\reg}}}]{greg}
\newterm[name={\ensuremath{\mathcal{T}(G)}}]{TG}
\newterm[name={\ensuremath{G(F)_{\elli}}}]{Gell}
\newterm[name={\ensuremath{\mathfrak{g}(F)_{\elli}}}]{gell}
\newterm[name={\ensuremath{G_{\reg}(F)_{\elli}}}]{Gregell}
\newterm[name={\ensuremath{\mathfrak{g}_{\reg}(F)_{\elli}}}]{gregell}
\newterm[name={\ensuremath{D^G(x)}}]{DGx}
\newterm[name={\ensuremath{D^G(X)}}]{DGX}
\newterm[name={\ensuremath{W(G,M)}}]{WGM}
\newterm[name={\ensuremath{W(G,T)}}]{WGT}
\newterm[name={\ensuremath{{}^g f}}]{gf}
\newterm[name={\ensuremath{R}}]{R}
\newterm[name={\ensuremath{L}}]{L}
\newterm[name={\ensuremath{\mathcal{U}(\mathfrak{g})}}]{Ug}
\newterm[name={\ensuremath{\mathcal{Z}(\mathfrak{g})}}]{Zg}
\newterm[name={\ensuremath{S(\mathfrak{g})}}]{Sg}
\newterm[name={\ensuremath{S(\mathfrak{g}^*)}}]{Sg*}
\newterm[name={\ensuremath{\partial(u)}}]{partu}
\newterm[name={\ensuremath{I(\mathfrak{g})}}]{Ig}
\newterm[name={\ensuremath{I(\mathfrak{g}^*)}}]{Ig*}
\newterm[name={\ensuremath{\Diff(\mathfrak{g})}}]{Dg}
\newterm[name={\ensuremath{u^*}}]{u*}
\newterm[name={\ensuremath{z_T}}]{zT}
\newterm[name={\ensuremath{z_H}}]{zH}
\newterm[name={\ensuremath{p_T}}]{pT}
\newterm[name={\ensuremath{p_H}}]{pH}
\newterm[name={\ensuremath{u_T}}]{uT}
\newterm[name={\ensuremath{u_H}}]{uH}
\newterm[name={\ensuremath{\sigma_X}}]{sigmaX}
\newterm[name={\ensuremath{f^*\sigma_Y}}]{fsigmaY}
\newterm[name={\ensuremath{f_*\sigma_X}}]{fsigmaX}
\newterm[name={\ensuremath{\lVert g\rVert_G}}]{gnorm}
\newterm[name={\ensuremath{\lVert X\rVert_{\mathfrak{g}}}}]{Xnorm}
\newterm[name={\ensuremath{\sigma}}]{sigma}
\newterm[name={\ensuremath{X[<C]}}]{X<C}
\newterm[name={\ensuremath{X[\geqslant C]}}]{X>C}
\newterm[name={\ensuremath{A_{\mini}^+}}]{Amin+}
\newterm[name={\ensuremath{A^{\overline{Q},+}_{\mini}(\delta)}}]{AQ+}
\newterm[name={\ensuremath{C^\infty(X)}}]{CX}
\newterm[name={\ensuremath{C^\infty(M)}}]{CM}
\newterm[name={\ensuremath{C_c^\infty(X)}}]{CcX}
\newterm[name={\ensuremath{C_c^\infty(M)}}]{CcM}
\newterm[name={\ensuremath{\Diff^\infty(M)}}]{DiffM}
\newterm[name={\ensuremath{\Diff^\infty_{\leqslant k}(M)}}]{DiffkM}
\newterm[name={\ensuremath{\mathcal{D}'(M)}}]{D'M}
\newterm[name={\ensuremath{\mathcal{D}'(X)}}]{D'X}
\newterm[name={\ensuremath{T_F}}]{TF}
\newterm[name={\ensuremath{f_\lambda}}]{flambda}
\newterm[name={\ensuremath{T_\lambda}}]{Tlambda}
\newterm[name={\ensuremath{\mathcal{S}(\mathfrak{g}(F))}}]{S(g)}
\newterm[name={\ensuremath{\mathcal{S}(G(F))}}]{S(G)}
\newterm[name={\ensuremath{q_{N,u}}}]{qNu}
\newterm[name={\ensuremath{\mathcal{S}'(\mathfrak{g}(F))}}]{S'(g)}
\newterm[name={\ensuremath{q_{N,u,v}}}]{qNuv}
\newterm[name={\ensuremath{\widehat{f}}}]{fhat}
\newterm[name={\ensuremath{\widehat{T}}}]{That}
\newterm[name={\ensuremath{p_u}}]{pu}
\newterm[name={\ensuremath{u_p}}]{up}
\newterm[name={\ensuremath{\Xi^G}}]{XiG}
\newterm[name={\ensuremath{\mathcal{C}(G(F))}}]{C(G(F))}
\newterm[name={\ensuremath{p_d}}]{pd}
\newterm[name={\ensuremath{p_{u,v,d}}}]{puvd}
\newterm[name={\ensuremath{\mathcal{C}^w(G(F))}}]{Cw(G)}
\newterm[name={\ensuremath{\psi}}]{psi}
\newterm[name={\ensuremath{\nu(T)}}]{nuT}
\newterm[name={\ensuremath{\Nil(\mathfrak{g})}}]{Nilg}
\newterm[name={\ensuremath{\mu_V^*}}]{muV*}
\newterm[name={\ensuremath{\mu_V^\perp}}]{muVperp}
\newterm[name={\ensuremath{J_G(x,f)}}]{JGxf}
\newterm[name={\ensuremath{J_G(X,f)}}]{JGXf}
\newterm[name={\ensuremath{J_{\mathcal{O}}(f)}}]{JOf}
\newterm[name={\ensuremath{\Nil_{\reg}(\mathfrak{g})}}]{Nilregg}
\newterm[name={\ensuremath{\widehat{j}}}]{jhat}
\newterm[name={\ensuremath{\widehat{j}(\mathcal{O},.)}}]{jhatO}
\newterm[name={\ensuremath{\Gamma(H)}}]{GammaH}
\newterm[name={\ensuremath{\Gamma(\mathfrak{h})}}]{Gammah}
\newterm[name={\ensuremath{\Gamma_{\elli}(G)}}]{GammaellG}
\newterm[name={\ensuremath{\Gamma_{\reg}(G)}}]{GammaregG}
\newterm[name={\ensuremath{\Gamma_{\elli}(\mathfrak{g})}}]{Gammaellg}
\newterm[name={\ensuremath{\Gamma_{\reg}(\mathfrak{g})}}]{Gammaregg}
\newterm[name={\ensuremath{\lVert .\rVert_{\Gamma(\mathfrak{g})}}}]{normGammag}
\newterm[name={\ensuremath{G(F)_{ani}}}]{GFani}
\newterm[name={\ensuremath{\Gamma_{ani}(G)}}]{GammaaniG}
\newterm[name={\ensuremath{v_M(g)}}]{vMg}
\newterm[name={\ensuremath{J_M(x,f)}}]{JMxf}
\newterm[name={\ensuremath{v_L^Q(g)}}]{vLQg}
\newterm[name={\ensuremath{J_L^Q(x,f)}}]{JLQxf}
\newterm[name={\ensuremath{J_L^Q(X,f)}}]{JLQXf}
\newterm[name={\ensuremath{J_L(X,f)}}]{JLXf}
\newterm[name={\ensuremath{V_\pi^\infty}}]{Vpiinf}
\newterm[name={\ensuremath{\pi^\infty}}]{piinf}
\newterm[name={\ensuremath{\pi_\lambda}}]{pilambda}
\newterm[name={\ensuremath{i\mathcal{A}_{G,\pi}^\vee}}]{iAGpivee}
\newterm[name={\ensuremath{\overline{\pi}}}]{pibar}
\newterm[name={\ensuremath{\overline{\mathcal{H}_\pi}}}]{Hpibar}
\newterm[name={\ensuremath{\widehat{K}}}]{Khat}
\newterm[name={\ensuremath{d(\rho)}}]{drho}
\newterm[name={\ensuremath{\Delta_K}}]{DeltaK}
\newterm[name={\ensuremath{c(\rho)}}]{crho}
\newterm[name={\ensuremath{\lVert .\rVert_u}}]{normu}
\newterm[name={\ensuremath{\mathcal{H}_\pi^{-\infty}}}]{Hpi-inf}
\newterm[name={\ensuremath{\pi^{-\infty}}}]{pi-inf}
\newterm[name={\ensuremath{(.,.)}}]{scprod}
\newterm[name={\ensuremath{\End(\pi)}}]{Endpi}
\newterm[name={\ensuremath{\vertiii{T}}}]{normT}
\newterm[name={\ensuremath{\End(\pi)^\infty}}]{Endpiinf}
\newterm[name={\ensuremath{\vertiii{T}_{u,v}}}]{normTuv}
\newterm[name={\ensuremath{\theta_\pi}}]{thetapi}
\newterm[name={\ensuremath{T_{e,e'}}}]{Tee'}
\newterm[name={\ensuremath{\Temp(G)}}]{TempG}
\newterm[name={\ensuremath{\omega_\pi}}]{omegapi}
\newterm[name={\ensuremath{\Pi_2(G)}}]{Pi2G}
\newterm[name={\ensuremath{\Pi_2(G)/i\mathcal{A}_{G,F}^*}}]{Pi2GmodA}
\newterm[name={\ensuremath{d(\pi)}}]{dpi}
\newterm[name={\ensuremath{\chi_\pi}}]{chipi}
\newterm[name={\ensuremath{N^G(\pi)}}]{normpi}
\newterm[name={\ensuremath{i_P^G(\sigma)}}]{iPGsigma}
\newterm[name={\ensuremath{i_P^G(\sigma)^\infty}}]{iPGsigmainf}
\newterm[name={\ensuremath{i_M^G(\sigma)}}]{iMGsigma}
\newterm[name={\ensuremath{J_{P'\mid P}(\sigma_\lambda)}}]{JPPsigma}
\newterm[name={\ensuremath{j(\sigma_\lambda)}}]{jsigma}
\newterm[name={\ensuremath{R_{P'\mid P}(\sigma_\lambda)}}]{RPPsigma}
\newterm[name={\ensuremath{J_L^Q(\sigma,f)}}]{JLQsigmaf}
\newterm[name={\ensuremath{\mathcal{X}_{\tempe}(G)}}]{XtempG}
\newterm[name={\ensuremath{C^\infty(\mathcal{X}_{\tempe}(G),V)}}]{C(Xtemp,V)}
\newterm[name={\ensuremath{C^\infty(\mathcal{X}_{\tempe}(G))}}]{C(Xtemp)}
\newterm[name={\ensuremath{\mu(\pi)}}]{mupi}
\newterm[name={\ensuremath{C^\infty(\mathcal{X}_{\tempe}(G),\mathcal{E}(G))}}]{C(Xtemp,E)}
\newterm[name={\ensuremath{\mathcal{C}(\mathcal{X}_{\tempe}(G),\mathcal{E}(G))}}]{Cc(Xtemp,E)}
\newterm[name={\ensuremath{R_{\tempe}(G)}}]{RtempG}
\newterm[name={\ensuremath{\mathcal{X}_{\elli}(G)}}]{XellG}
\newterm[name={\ensuremath{D(\pi)}}]{Dpi}
\newterm[name={\ensuremath{R_{\elli}(G)}}]{RellG}
\newterm[name={\ensuremath{R_{ind}(G)}}]{RindG}
\newterm[name={\ensuremath{\mathcal{X}_{\elli}(G)/i\mathcal{A}_{G,F}^*}}]{XellmodA}
\newterm[name={\ensuremath{\mathcal{X}(G)}}]{XG}
\newterm[name={\ensuremath{C^\infty(\omega)^G}}]{ComegaG}
\newterm[name={\ensuremath{C^\infty(\Omega)^G}}]{COmegaG}
\newterm[name={\ensuremath{\mathcal{D}'(\omega)^G}}]{D'omegaG}
\newterm[name={\ensuremath{\mathcal{D}'(\Omega)^G}}]{D'OmegaG}
\newterm[name={\ensuremath{\mathcal{S}(\omega)}}]{Somega}
\newterm[name={\ensuremath{\mathcal{S}(\Omega)}}]{SOmega}
\newterm[name={\ensuremath{\Diff_{\leqslant n}^\infty(\omega)^G}}]{DiffnomegaG}
\newterm[name={\ensuremath{\Diff_{\leqslant n}^\infty(\Omega)^G}}]{DiffnOmegaG}
\newterm[name={\ensuremath{\Diff^\infty(\omega)^G}}]{DiffomegaG}
\newterm[name={\ensuremath{\Diff^\infty(\Omega)^G}}]{DiffOmegaG}
\newterm[name={\ensuremath{\mathcal{J}^\infty(\omega)}}]{Jinfomega}
\newterm[name={\ensuremath{\mathcal{J}^\infty(\Omega)}}]{JinfOmega}
\newterm[name={\ensuremath{\overline{\Diff^\infty(\omega)^G}}}]{DiffomegaGbar}
\newterm[name={\ensuremath{\overline{\Diff^\infty(\Omega)^G}}}]{DiffOmegaGbar}
\newterm[name={\ensuremath{\Diff_{\leqslant n}(\mathfrak{g})^G}}]{DiffngG}
\newterm[name={\ensuremath{\Diff(\mathfrak{g})^G}}]{DiffgG}
\newterm[name={\ensuremath{\mathcal{J}(\mathfrak{g})}}]{Jg}
\newterm[name={\ensuremath{\overline{\Diff(\mathfrak{g})^G}}}]{DiffgGbar}
\newterm[name={\ensuremath{\eta_X^G}}]{etaXG}
\newterm[name={\ensuremath{f_{X,\omega_X}}}]{fXomegaX}
\newterm[name={\ensuremath{T_{X,\omega_X}}}]{TXomegaX}
\newterm[name={\ensuremath{\eta_x^G}}]{etaxG}
\newterm[name={\ensuremath{f_{x,\Omega_x}}}]{fxOmegax}
\newterm[name={\ensuremath{p_X}}]{pX}
\newterm[name={\ensuremath{u_X}}]{uX}
\newterm[name={\ensuremath{z_x}}]{zx}
\newterm[name={\ensuremath{D_{X,\omega_X}}}]{DXomegaX}
\newterm[name={\ensuremath{D_{x,\Omega_x}}}]{DxOmegax}
\newterm[name={\ensuremath{D_T}}]{DT}
\newterm[name={\ensuremath{f_T}}]{fT}
\newterm[name={\ensuremath{j^G}}]{jG}
\newterm[name={\ensuremath{f_\omega}}]{fomega}
\newterm[name={\ensuremath{T_\omega}}]{Tomega}
\newterm[name={\ensuremath{D_\omega}}]{Domega}
\newterm[name={\ensuremath{f_\Omega}}]{fOmega}
\newterm[name={\ensuremath{T_\Omega}}]{TOmega}
\newterm[name={\ensuremath{D_\Omega}}]{DOmega}
\newterm[name={\ensuremath{f^{(P)}}}]{f(P)}
\newterm[name={\ensuremath{T_P}}]{TP}
\newterm[name={\ensuremath{i_M^G(T)}}]{iMGT}
\newterm[name={\ensuremath{G_{\sreg}}}]{Gsreg}
\newterm[name={\ensuremath{\souligner{\mathcal{X}}(G)}}]{underXG}
\newterm[name={\ensuremath{\souligner{\mathcal{X}}_{\elli}(G)}}]{underXellG}
\newterm[name={\ensuremath{\mathcal{X}^M(x)}}]{XMx}
\newterm[name={\ensuremath{QC(\omega)}}]{QComega}
\newterm[name={\ensuremath{QC_c(\omega)}}]{QCcomega}
\newterm[name={\ensuremath{QC(\Omega)}}]{QCOmega}
\newterm[name={\ensuremath{QC_c(\Omega)}}]{QCcOmega}
\newterm[name={\ensuremath{SQC(\mathfrak{g}(F))}}]{SQCgF}
\newterm[name={\ensuremath{SQC(\mathfrak{g}(\mathbb{R}))}}]{SQCgR}
\newterm[name={\ensuremath{q_{L,u}}}]{qLu}
\newterm[name={\ensuremath{q_{u,N}}}]{quN}
\newterm[name={\ensuremath{q_{L,z}}}]{qLz}
\newterm[name={\ensuremath{c_{\theta,\mathcal{O}}(x)}}]{cthetaOx}
\newterm[name={\ensuremath{c_\theta}}]{ctheta}
\newterm[name={\ensuremath{M_\lambda}}]{Mlambda}
\newterm[name={\ensuremath{c_{\pi,\mathcal{O}}(1)}}]{cpiO1}
\newterm[name={\ensuremath{f^P}}]{fP}
\newterm[name={\ensuremath{\varphi^P}}]{phiP}
\newterm[name={\ensuremath{\mathcal{C}_{\scusp}(G(F))}}]{CscuspG}
\newterm[name={\ensuremath{\mathcal{S}_{\scusp}(G(F))}}]{SscuspG}
\newterm[name={\ensuremath{\mathcal{S}_{\scusp}(\mathfrak{g}(F))}}]{Sscuspg}
\newterm[name={\ensuremath{\mathcal{S}_{\scusp}(\Omega)}}]{SscuspOmega}
\newterm[name={\ensuremath{\mathcal{S}_{\scusp}(\omega)}}]{Sscuspomega}
\newterm[name={\ensuremath{M(x)}}]{M(x)}
\newterm[name={\ensuremath{\theta_f}}]{thetaf}
\newterm[name={\ensuremath{M(X)}}]{M(X)}
\newterm[name={\ensuremath{\mathcal{C}_{\scusp}(\mathcal{X}_{\tempe}(G),\mathcal{E}(G))}}]{Cscusp(Xtemp,E)}
\newterm[name={\ensuremath{\widehat{\theta}_f}}]{thetafhat}
\newterm[name={\ensuremath{K^A_{f,f'}}}]{KAff'}
\newterm[name={\ensuremath{E}}]{E}
\newterm[name={\ensuremath{\overline{x}}}]{xbar}
\newterm[name={\ensuremath{N_{E/F}}}]{NEF}
\newterm[name={\ensuremath{\Tra_{E/F}}}]{TrEF}
\newterm[name={\ensuremath{\eta}}]{eta}
\newterm[name={\ensuremath{\overline{E}}}]{Ebar}
\newterm[name={\ensuremath{U(V)}}]{U(V)}
\newterm[name={\ensuremath{\mathfrak{u}(V)}}]{u(V)}
\newterm[name={\ensuremath{c(v,v')}}]{c(v,v')}
\newterm[name={\ensuremath{(W,V)}}]{(W,V)}
\newterm[name={\ensuremath{Z}}]{Z}
\newterm[name={\ensuremath{G}}]{G}
\newterm[name={\ensuremath{H}}]{H}
\newterm[name={\ensuremath{\xi}}]{xi}
\newterm[name={\ensuremath{(z_i)_{i=0,\pm 1,\ldots,\pm r}}}]{Zbasis}
\newterm[name={\ensuremath{P_V}}]{PV}
\newterm[name={\ensuremath{N}}]{N}
\newterm[name={\ensuremath{M_V}}]{MV}
\newterm[name={\ensuremath{P}}]{P}
\newterm[name={\ensuremath{M}}]{M}
\newterm[name={\ensuremath{\lambda}}]{lambda}
\newterm[name={\ensuremath{r}}]{r}
\newterm[name={\ensuremath{\overline{F}}}]{Fbar}
\newterm[name={\ensuremath{\lambda_F}}]{lambdaF}
\newterm[name={\ensuremath{\nu}}]{nu}
\newterm[name={\ensuremath{d}}]{d}
\newterm[name={\ensuremath{m}}]{m}
\newterm[name={\ensuremath{Z_+}}]{Z+}
\newterm[name={\ensuremath{Z_-}}]{Z-}
\newterm[name={\ensuremath{D}}]{D}
\newterm[name={\ensuremath{V_0}}]{V_0}
\newterm[name={\ensuremath{H_0}}]{H_0}
\newterm[name={\ensuremath{G_0}}]{G_0}
\newterm[name={\ensuremath{T}}]{T}
\newterm[name={\ensuremath{A}}]{A}
\newterm[name={\ensuremath{B(.,.)}}]{B}
\newterm[name={\ensuremath{\Hom_H(\pi,\xi)}}]{HomHpixi}
\newterm[name={\ensuremath{m(\pi)}}]{mpi}
\newterm[name={\ensuremath{\overline{P}_0}}]{P0bar}
\newterm[name={\ensuremath{A_0}}]{A_0}
\newterm[name={\ensuremath{A_0^+}}]{A0+}
\newterm[name={\ensuremath{\Delta_P}}]{DeltaP}
\newterm[name={\ensuremath{\Xi^{H\backslash G}}}]{XiHG}
\newterm[name={\ensuremath{\mathcal{P}_{H,\xi}}}]{PHxi}
\newterm[name={\ensuremath{\mathcal{L}_\pi}}]{Lpi}
\newterm[name={\ensuremath{L_\pi}}]{Lpid}
\newterm[name={\ensuremath{K(f,.)}}]{Kf}
\newterm[name={\ensuremath{J(f)}}]{Jf}
\newterm[name={\ensuremath{K^{\Lie}(f,.)}}]{KLief}
\newterm[name={\ensuremath{J^{\Lie}(f)}}]{JLief}
\newterm[name={\ensuremath{J_{\spec}(f)}}]{Jspecf}
\newterm[name={\ensuremath{J_{\auxi}(f,f')}}]{Jauxff'}
\newterm[name={\ensuremath{\mathcal{C}_d^w(G(F))}}]{CdwG}
\newterm[name={\ensuremath{d\pi}}]{dpimes}
\newterm[name={\ensuremath{\overline{P}}}]{Pbar}
\newterm[name={\ensuremath{\overline{N}}}]{Nbar}
\newterm[name={\ensuremath{\Xi}}]{Xi}
\newterm[name={\ensuremath{\mathfrak{h}^\perp}}]{hperp}
\newterm[name={\ensuremath{\Sigma}}]{Sigma}
\newterm[name={\ensuremath{\Lambda}}]{Lambda}
\newterm[name={\ensuremath{Q}}]{Q}
\newterm[name={\ensuremath{\Lambda'}}]{Lambda'}
\newterm[name={\ensuremath{\Sigma'}}]{Sigma'}
\newterm[name={\ensuremath{B(0,R)}}]{B0r}
\newterm[name={\ensuremath{\Gamma(\Sigma)}}]{GammaSigma}
\newterm[name={\ensuremath{\kappa_N}}]{kappaN}
\newterm[name={\ensuremath{J^{\Lie}_N}}]{JLieN}
\newterm[name={\ensuremath{\gamma_X}}]{gammaX}
\newterm[name={\ensuremath{X_\Sigma}}]{XSigma}
\newterm[name={\ensuremath{\kappa_{N,X}}}]{kappaNX}
\newterm[name={\ensuremath{J^{\Lie}_{N,T}}}]{JLieNT}
\newterm[name={\ensuremath{\mathfrak{t}(F)'[>\epsilon]}}]{t'>}
\newterm[name={\ensuremath{M_T}}]{MT}
\newterm[name={\ensuremath{\mathfrak{t}(F)'}}]{tF'}
\newterm[name={\ensuremath{\widetilde{v}(Y,.)}}]{tildevY}
\newterm[name={\ensuremath{J_{Y,T,N^{-b}}(f)}}]{JYTN-bf}
\newterm[name={\ensuremath{J^{\Lie}_{N,T,\epsilon}(f)}}]{JNTepsilonf}
\newterm[name={\ensuremath{(G_x,H_x,\xi_x)}}]{GxHxxix}
\newterm[name={\ensuremath{W'_x}}]{W'_x}
\newterm[name={\ensuremath{V'_x}}]{V'_x}
\newterm[name={\ensuremath{W''_x}}]{W''_x}
\newterm[name={\ensuremath{H''_x}}]{H''_x}
\newterm[name={\ensuremath{G''_x}}]{G''_x}
\newterm[name={\ensuremath{H'_x}}]{H'_x}
\newterm[name={\ensuremath{G'_x}}]{G'_x}
\newterm[name={\ensuremath{\Gamma(G,H)}}]{GammaGH}
\newterm[name={\ensuremath{\souligner{\mathcal{T}}}}]{undercalT}
\newterm[name={\ensuremath{T_{\natural}}}]{Tnat}
\newterm[name={\ensuremath{W''_T}}]{W''_T}
\newterm[name={\ensuremath{W(T)}}]{W(T)}
\newterm[name={\ensuremath{\mathcal{T}}}]{calT}
\newterm[name={\ensuremath{\Gamma(G_x,H_x)}}]{GammaGxHx}
\newterm[name={\ensuremath{\Gamma^{\Lie}(G,H)}}]{GammaLieGH}
\newterm[name={\ensuremath{G'_X}}]{G'_X}
\newterm[name={\ensuremath{G''_X}}]{G''_X}
\newterm[name={\ensuremath{H'_X}}]{H'_X}
\newterm[name={\ensuremath{H''_X}}]{H''_X}
\newterm[name={\ensuremath{W'_X}}]{W'_X}
\newterm[name={\ensuremath{V'_X}}]{V'_X}
\newterm[name={\ensuremath{W''_X}}]{W''_X}
\newterm[name={\ensuremath{\mathfrak{t}_\natural}}]{tnat}
\newterm[name={\ensuremath{\Delta(x)}}]{Deltax}
\newterm[name={\ensuremath{\Delta(X)}}]{DeltaX}
\newterm[name={\ensuremath{j_G^H(X)}}]{jGHX}
\newterm[name={\ensuremath{\Delta_x(y)}}]{Deltaxy}
\newterm[name={\ensuremath{\eta_{G,x}^H(y)}}]{etaGxHy}
\newterm[name={\ensuremath{m_{\geom}(\theta)}}]{mgeomtheta}
\newterm[name={\ensuremath{m_{x,\geom}(\theta_x)}}]{mxgeomthetax}
\newterm[name={\ensuremath{m_{\geom}^{\Lie}(\theta)}}]{mgeomLietheta}
\newterm[name={\ensuremath{m_{\geom}^L}}]{mgeomL}
\newterm[name={\ensuremath{J_{\geom}(f)}}]{Jgeomf}
\newterm[name={\ensuremath{m_{\geom}(\pi)}}]{mgeompi}
\newterm[name={\ensuremath{J_{\geom}^{\Lie}(f)}}]{JgeomLief}
\newterm[name={\ensuremath{J_{\qc}}}]{Jqc}
\newterm[name={\ensuremath{J_{\qc}^{\Lie}}}]{JqcLie}
\newterm[name={\ensuremath{\Gamma_F}}]{GammaF}
\newterm[name={\ensuremath{\sim_{\stab}}}]{stabconj}
\newterm[name={\ensuremath{G_\alpha}}]{Galpha}
\newterm[name={\ensuremath{e(G_\alpha)}}]{eGalpha}
\newterm[name={\ensuremath{Br_2(F)}}]{Br2F}
\newterm[name={\ensuremath{\widetilde{G}}}]{Gtilde}
\newterm[name={\ensuremath{\widetilde{U(V)}}}]{U(V)tilde}
\newterm[name={\ensuremath{(G_\alpha,H_\alpha,\xi_\alpha)}}]{GalphaHalphaxialpha}
\newterm[name={\ensuremath{W_\alpha}}]{Walpha}
\newterm[name={\ensuremath{V_\alpha}}]{Valpha}
\newterm[name={\ensuremath{W_F}}]{W_F}
\newterm[name={\ensuremath{L_F}}]{L_F}
\newterm[name={\ensuremath{{}^LG}}]{LG}
\newterm[name={\ensuremath{\Pi^G(\varphi)}}]{PiGphi}
\newterm[name={\ensuremath{\Pi^{G_\alpha}(\varphi)}}]{PiGalphaphi}
\newterm[name={\ensuremath{\theta_{\alpha,\varphi}}}]{thetaalphaphi}
\newterm[name={\ensuremath{\theta_{\varphi}}}]{thetaphi}
\newterm[name={\ensuremath{T_x}}]{T_x}
\newterm[name={\ensuremath{\widehat{\otimes}_p}}]{otimesp}
\newterm[name={\ensuremath{C_b^\infty(V,U)}}]{CbinfVU}

\newtheorem{conj}{Conjecture}[subsection]
\newtheorem{theo}[conj]{Theorem}
\newtheorem{them}{Theorem}
\newtheorem{lem}[conj]{Lemma}
\newtheorem{prop}[conj]{Proposition}
\newtheorem{cor}[conj]{Corollary}

\title{A local trace formula for the Gan-Gross-Prasad conjecture for unitary groups: the Archimedean case \protect}
\author{Rapha\"{e}l Beuzart-Plessis}
\begin{document}

\maketitle

\begin{abstract}
In this paper, we prove, following earlier work of Waldspurger \cite{Wa1}, \cite{Wa4} a sort of local relative trace formula which is related to the local Gan-Gross-Prasad conjecture for unitary groups over a local field $F$ of characteristic zero. As a consequence, we obtain a geometric formula for certain multiplicities $m(\pi)$ appearing in this conjecture and deduce from it a weak form of the local Gan-Gross-Prasad conjecture (multiplicity one in tempered L-packets). These results were already known over $p$-adic fields \cite{Beu1} and thus are only new when $F=\mathbb{R}$. 
\end{abstract}

\tableofcontents

\section*{Introduction}

Let $F$ be a local field of characteristic $0$ which is different from $\mathbb{C}$. So, $F$ is either a $p$-adic field (that is a finite extension of $\mathbb{Q}_p$) or $F=\mathbb{R}$. Let $E/F$ be a quadratic extension of $F$ (if $F=\mathbb{R}$, we have $E=\mathbb{C}$) and let $W\subset V$ be a pair of hermitian spaces having the following property: the orthogonal complement $W^\perp$ of $W$ in $V$ is odd-dimensional and its unitary group $U(W^\perp)$ is quasi-split. To such a pair (that we call an admissible pair, cf.\ Section \ref{section 6.2}), Gan, Gross and Prasad associate a triple $(G,H,\xi)$. Here, $G$ is equal to the product $U(W)\times U(V)$ of the unitary groups of $W$ and $V$, $H$ is a certain algebraic subgroup of $G$ and $\xi:H(F)\to \mathbb{S}^1$ is a continuous unitary character of the $F$-points of $H$. In the case where $\dim(W^\perp)=1$, we just have $H=U(W)$ embedded in $G$ diagonally and the character $\xi$ is trivial. For the definition in codimension greater than $1$, we refer the reader to Section \ref{section 6.2}. We call a triple like $(G,H,\xi)$ (constructed from an admissible pair $(W,V)$) a GGP triple.

\vspace{2mm}

\noindent Let $\pi$ be a tempered irreducible representation of $G(F)$. By this, we mean that $\pi$ is an irreducible unitary representation of $G(F)$ whose coefficients satisfy a certain growth condition (an equivalent condition is that $\pi$ belongs weakly to the regular representation of $G(F)$). We denote by $\pi^\infty$ the subspace of smooth vectors in $\pi$. This subspace is $G(F)$-invariant and carries a natural topology (if $F=\mathbb{R}$, this topology makes $\pi^\infty$ into a Fr\'echet space whereas if $F$ is $p$-adic the topology on $\pi^\infty$ doesn't play any role but in order to get a uniform treatment we endow $\pi^\infty$ with its finest locally convex topology). Following Gan, Gross and Prasad, we define a multiplicity $m(\pi)$ by

$$m(\pi)=\dim \Hom_H(\pi^\infty,\xi)$$

\noindent where $\Hom_H(\pi^\infty,\xi)$ denotes the space of continuous linear forms $\ell$ on $\pi^\infty$ satisfying the relation $\ell\circ\pi(h)=\xi(h)\ell$ for all $h\in H(F)$. By the main result of \cite{JSZ} (in the real case) and \cite{AGRS} (in the $p$-adic case) together with Theorem 15.1 of \cite{GGP}, we know that this multiplicity is always less or equal to $1$.

\vspace{2mm}

\noindent The main result of this paper extends this multiplicity one result to a whole $L$-packet of tempered representations of $G(F)$. This answers a conjecture of Gan, Gross and Prasad (Conjecture 17.1 of \cite{GGP}). Actually, the result is better stated if we consider more than one GGP triple at the same time. In any family of GGP triples that we are going to consider there is a distinguished one corresponding to the case where $G$ and $H$ are quasi-split over $F$. So, for convenience, we assume that the GGP triple $(G,H,\xi)$ we started with satisfies this condition. The other GGP triples that we need to consider may be called the pure inner forms of $(G,H,\xi)$. Those are naturally parametrized by the Galois cohomology set $H^1(F,H)$. A cohomology class $\alpha\in H^1(F,H)$ corresponds to a hermitian space $W_\alpha$ (up to isomorphism) of the same dimension as $W$. If we set $V_\alpha=W_\alpha\oplus^\perp W^\perp$, then $(W_\alpha,V_\alpha)$ is an admissible pair and thus gives rise to a new GGP triple $(G_\alpha,H_\alpha,\xi_\alpha)$. The pure inner forms of $(G,H,\xi)$ are exactly all the GGP triples obtained in this way.

\vspace{2mm}

\noindent Let $\varphi$ be a tempered Langlands parameter for $G$. According to the local Langlands correspondence (which is now known in all cases for unitary groups, cf.\ \cite{KMSW} and \cite{Mok}), this parameter determines an $L$-packet $\Pi^G(\varphi)$ consisting of a finite number of tempered representations of $G(F)$. Actually, this parameter also defines $L$-packets $\Pi^{G_\alpha}(\varphi)$ of tempered representations of $G_\alpha(F)$ for all $\alpha\in H^1(F,H)$. We can now state the main result of this paper as follows (cf.\ Theorem \ref{theorem 12.4.1}).

\vspace{2mm}

\begin{them}\label{theorem 1}
There exists exactly one representation $\pi$ in the disjoint union of $L$-packets

$$\displaystyle \bigsqcup_{\alpha\in H^1(F,H)}\Pi^{G_\alpha}(\varphi)$$

\noindent such that $m(\pi)=1$.
\end{them}

\vspace{2mm}

\noindent As we said, this answers in the affirmative a conjecture of Gan-Goss-Prasad (Conjecture 17.1 of \cite{GGP}). The analog of this theorem for special orthogonal groups has already been obtained by Waldspurger in the case where $F$ is $p$-adic \cite{Wa1}. In \cite{Beu1}, the author adapted the proof of Waldspurger to deal with unitary groups but again under the assumption that $F$ is $p$-adic. Hence, the only new result contained in Theorem \ref{theorem 1} is when $F=\mathbb{R}$. However, the proof we present here differs slightly from the original treatment of Waldspurger and we feel that this new approach is more amenable to generalizations. This is the main reason why we are including the $p$-adic case in this paper. Actually, it doesn't cost much: in many places, we have been able to treat the two cases uniformly and when we needed to make a distinction, it is often because the real case is more tricky.

\vspace{2mm}

\noindent As in \cite{Wa1} and subsequently \cite{Beu1}, Theorem \ref{theorem 1} follows from a formula for the multiplicity $m(\pi)$. This formula express $m(\pi)$ in terms of the Harish-Chandra character of $\pi$. Recall that, according to Harish-Chandra, there exists a smooth function $\theta_\pi$ on the regular locus $G_{\reg}(F)$ of $G(F)$ which is locally integrable on $G(F)$ and such that

$$\displaystyle \Tr\; \pi(f)=\int_{G(F)} \theta_\pi(x)f(x)dx$$

\noindent for all $f\in C_c^\infty(G(F))$ (here $C_c^\infty(G(F))$ denotes the space of smooth and compactly supported functions on $G(F)$). This function $\theta_\pi$ is obviously unique and is called the Harish-Chandra character of $\pi$. To state the formula for the multiplicity, we need to extend the character $\theta_\pi$ to a function

$$c_\pi\colon G_{\ssi}(F)\to \mathbb{C}$$

\noindent on the semi-simple locus $G_{\ssi}(F)$ of $G(F)$. If $x\in G_{\reg}(F)$, then $c_\pi(x)=\theta_\pi(x)$ but for a general element $x\in G_{\ssi}(F)$, $c_\pi(x)$ is in some sense the main coefficient of a certain local expansion of $\theta_\pi$ near $x$. For a precise definition of the function $c_\pi$, we refer the reader to Section \ref{section 4.5}, where we consider more general functions that we call quasi-characters and which are smooth functions on $G_{\reg}(F)$ sharing almost all of the good properties that characters of representations have. As we said, it is through the function $c_\pi$ that the character $\theta_\pi$ will appear in the multiplicity formula. The other main ingredient of this formula is a certain space $\Gamma(G,H)$ of semi-simple conjugacy classes in $G(F)$. For a precise definition of $\Gamma(G,H)$, we refer the reader to Section \ref{section 11.1}. Let us just say that $\Gamma(G,H)$ comes naturally equipped with a measure $dx$ on it and that this measure is not generally supported in the regular locus. For example, the trivial conjugacy class $\{1\}$ is an atom for this measure whose mass is equal to $1$. Apart from these two main ingredients (the function $c_\pi$ and the space $\Gamma(G,H)$), the formula for the multiplicity involves two normalizing functions $D^G$ and $\Delta$. Here, $D^G$ is the usual discriminant whereas $\Delta$ is some determinant function that is defined in Section \ref{section 11.1}. We can now state the formula for the multiplicity as follows (cf.\ Theorem \ref{theorem 11.2.2}).

\vspace{2mm}

\begin{them}\label{theorem 2}
For every irreducible tempered representation $\pi$ of $G(F)$, we have the equality

$$\displaystyle m(\pi)=\lim\limits_{s\to 0^+} \int_{\Gamma(G,H)} c_\pi(x)D^G(x)^{1/2}\Delta(x)^{s-1/2}dx$$
\end{them}

\vspace{2mm}

\noindent The integral in the right hand side of the equality above is absolutely convergent for all $s\in \mathbb{C}$ such that $Re(s)>0$ and moreover the limit as $s\to 0^+$ exists (cf.\ Proposition \ref{proposition 11.1.1}).

\vspace{2mm}

\noindent As we said, Theorem \ref{theorem 1} follows from Theorem \ref{theorem 2}. This is proved in the last chapter of this paper (Chapter \ref{section 12}). Let us fix a tempered Langlands parameter $\varphi$ for $G$. The main idea of the proof, the same as for Theorem 13.3 of \cite{Wa1}, is to show that the sum

\begin{align}\label{eq 0.1}
\displaystyle \sum_{\alpha\in H^1(F,H)}\sum_{\pi\in \Pi^{G_\alpha}(\varphi)} m(\pi)
\end{align}

\noindent when expressed geometrically through Theorem \ref{theorem 2} contains a lot of cancellations which roughly come from the character relations between the various stable characters associated to $\varphi$ on the pure inner forms of $G$. Once these cancellations are taken into account, the only remaining term is the term corresponding to the conjugacy class of the identity inside $\Gamma(G,H)$. By classical results of Rodier and Matumoto, this last term is related to the number of generic representations inside the quasi-split $L$-packet $\Pi^G(\varphi)$. By the generic packet conjecture, which is now known for unitary groups, we are able to show that this term is equal to $1$ and this immediately implies Theorem \ref{theorem 1}. Let us now explain in more detail how it works. Fix momentarily $\alpha\in H^1(F,H)$. Using Theorem \ref{theorem 2}, we can express the sum

$$\displaystyle \sum_{\pi\in \Pi^{G_\alpha}(\varphi)} m(\pi)$$

\noindent as

\begin{align}\label{eq 0.2}
\displaystyle \lim\limits_{s\to 0^+}\int_{\Gamma(G_\alpha,H_\alpha)} c_{\varphi,\alpha}(x) D^{G_\alpha}(x)^{1/2}\Delta(x)^{s-1/2}dx
\end{align}

\noindent where we have set $c_{\varphi,\alpha}=\sum_{\pi\in \Pi^{G_\alpha}(\varphi)} c_\pi$. One of the main properties of $L$-packets is that the sum of characters $\theta_{\varphi,\alpha}=\sum_{\pi\in \Pi^{G_\alpha}(\varphi)} \theta_\pi$ defines a function on $G_{\alpha,\reg}(F)$ which is stable, which here means that it is invariant by $G_{\alpha}(\overline{F})$-conjugation. In Section \ref{section 12.1}, we define a notion of strongly stable conjugacy for semi-simple elements of $G_\alpha(F)$. This definition of stable conjugacy differs from the usually accepted one (cf.\ \cite{Kott1}) and is actually stronger (hence the use of the word ``strongly"). The point of introducing such a notion is the following: it easily follows from the stability of $\theta_{\varphi,\alpha}$ that the function $c_{\varphi,\alpha}$ is constant on semi-simple strongly stable conjugacy classes. This allows us to further transform the expression \ref{eq 0.2} to write it as 

$$\displaystyle \lim\limits_{s\to 0^+}\int_{\Gamma_{\stab}(G_\alpha,H_\alpha)}\lvert p_{\alpha,\stab}^{-1}(x)\rvert c_{\varphi,\alpha}(x) D^{G_\alpha}(x)^{1/2}\Delta(x)^{s-1/2}dx$$

\noindent where $\Gamma_{\stab}(G_{\alpha},H_\alpha)$ denotes the space of strongly stable conjugacy classes in $\Gamma(G_\alpha,H_\alpha)$ and $p_{\alpha,\stab}$ stands for the natural projection $\Gamma(G_\alpha,H_\alpha)\twoheadrightarrow \Gamma_{\stab}(G_\alpha,H_\alpha)$ (thus $\lvert p_{\alpha,\stab}^{-1}(x)\rvert$ is just the number of conjugacy classes in $\Gamma(G_\alpha,H_\alpha)$ belonging to the strongly stable conjugacy class of $x$). Returning to the sum \ref{eq 0.1}, we can now write it as

\begin{align}\label{eq 0.3}
\displaystyle \sum_{\alpha\in H^1(F,H)} \lim\limits_{s\to 0^+} \int_{\Gamma_{\stab}(G_\alpha,H_\alpha)}\lvert p_{\alpha,\stab}^{-1}(x)\rvert c_{\varphi,\alpha}(x) D^{G_\alpha}(x)^{1/2}\Delta(x)^{s-1/2}dx
\end{align}

\noindent A second very important property of $L$-packets is that the stable character $\theta_{\varphi,\alpha}$ is related in a simple manner to the stable character $\theta_{\varphi,1}$ on the quasi-split form $G(F)$. More precisely, Kottwitz \cite{Kott2} has defined a sign $e(G_\alpha)$ such that we have $\theta_{\varphi,\alpha}(y)=e(G_\alpha)\theta_{\varphi,1}(x)$ as soon as $y\in G_{\alpha,\reg}(F)$ and $x\in G_{\reg}(F)$ are stably conjugate regular elements (i.e., are conjugate over the algebraic closure where $G_\alpha(\overline{F})=G(\overline{F})$). Once again, this relation extends to the functions $c_{\varphi,\alpha}$ and $c_{\varphi,1}$ and we have $c_{\varphi,\alpha}(y)=e(G_\alpha)c_{\varphi,1}(x)$ for all strongly stably conjugate elements $y\in G_{\alpha,ss}(F)$ and $x\in G_{\ssi}(F)$. As it happens, and contrary to the regular case, there might exist semi-simple elements in $G_\alpha(F)$ which are not strongly stably conjugate to any element of the quasi-split form $G(F)$. However, we can show that the function $c_{\varphi,\alpha}$ vanishes on such elements $x\in G_{\alpha,ss}(F)$. Therefore, these conjugacy classes don't contribute to the sum \ref{eq 0.3} and transferring the remaining terms to $G(F)$, we can express \ref{eq 0.3} as a single integral

$$\displaystyle \lim\limits_{s\to 0^+} \int_{\Gamma(G,H)} \left(\sum_{y\sim_{\stab} x} e(G_{\alpha(y)})\right) c_{\varphi,1}(x) D^G(x)^{1/2}\Delta(x)^{s-1/2}dx$$

\noindent where the sum

\begin{align}\label{eq 0.4}
\displaystyle \sum_{y\sim_{\stab} x} e(G_{\alpha(y)})
\end{align}

\noindent is over the conjugacy classes $y$ in the disjoint union $\bigsqcup_{\alpha\in H^1(F,H)}\Gamma(G_\alpha, H_\alpha)$ that are strongly stably conjugate to $x$ and $\alpha(y)\in H^1(F,H)$ denotes the only cohomology class such that $y$ lives in $\Gamma(G_{\alpha(y)},H_{\alpha(y)})$. There is a natural anisotropic torus $T_x\subset H$ associated to $x\in \Gamma_{\stab}(G,H)$ such that the set of conjugacy classes in $\bigsqcup_{\alpha\in H^1(F,H)}\Gamma(G_\alpha, H_\alpha)$ lying inside the strongly stable conjugacy class of $x$ is naturally in bijection with $H^1(F,T_x)$ (cf.\ Section \ref{section 12.5} for the definition of $T_x$). Moreover, for $y\in H^1(F,T_x)$, the cohomology class $\alpha(y)$ is just the image of $y$ via the natural map $H^1(F,T_x)\to H^1(F,H)$. Hence, the sum \ref{eq 0.4} equals

\begin{align}\label{eq 0.5}
\displaystyle \sum_{y\in H^1(F,T_x)} e(G_{\alpha(y)})
\end{align}

\noindent In order to further analyze this sum, we need to recall the definition of the sign $e(G_\alpha)$. In \cite{Kott2}, Kottwitz constructs a natural map $H^1(F,G)\to H^2(F,\{\pm 1\})=Br_2(F)$ from $H^1(F,G)$ to the $2$-torsion subgroup of the Brauer group of $F$. Since $F$ is either $p$-adic or real, we have an isomorphism $Br_2(F)\simeq \{\pm 1\}$. The sign $e(G_\alpha)$ for $\alpha\in H^1(F,H)$ is now just the image of $\alpha$ by the composition of this map with $H^1(F,H)\to H^1(F,G)$. Following Kottwitz's definition, it is not hard to see that the composition $H^1(F,T_x)\to H^1(F,G)\to Br_2(F)$ is a group homomorphism. Moreover, it turns out that for $x\neq 1$ this homomorphism is surjective and this immediately implies that for such an $x$ the sum \ref{eq 0.5} is zero. Going back to \ref{eq 0.3}, we are only left with the contribution of $1\in \Gamma(G,H)$ which is equal to

$$c_{\varphi,1}(1)$$

\noindent By a result of Rodier \cite{Ro} in the $p$-adic case and of Matumoto \cite{Mat} in the real case, the term $c_{\varphi,1}(1)$ has an easy interpretation in terms of Whittaker models. More precisely, this term equals the number of representations in the $L$-packet $\Pi^G(\varphi)$ having a Whittaker model, a representation being counted as many times as the number of types of Whittaker models it has, divided by the number of types of Whittaker models for $G(F)$. A third important property of $L$-packets is that $\Pi^G(\varphi)$ contains exactly one representation having a Whittaker model of a given type. It easily follows from this that $c_{\varphi,1}(1)=1$. Hence, the sum \ref{eq 0.1} equals $1$ and this ends our explanation of how Theorem \ref{theorem 2} implies Theorem \ref{theorem 1}.

\vspace{2mm}

\noindent The proof of Theorem \ref{theorem 2} is more involved and takes up most of this paper. It is at this point that our strategy differs from the one of Waldspurger. In what follows, we explain the motivations and the main steps of the proof of Theorem \ref{theorem 2}. Consider the unitary representation $L^2(H(F)\backslash G(F),\xi)$ of $G(F)$. It is the $L^2$-induction of the character $\xi$ from $H(F)$ to $G(F)$ and it consists in the measurable functions $\varphi\colon G(F)\to\mathbb{C}$ satisfying the relation $\varphi(hg)=\xi(h)\varphi(g)$ ($h\in H(F)$, $g\in G(F)$) almost everywhere and such that

$$\displaystyle \int_{H(F)\backslash G(F)} \lvert \varphi(x)\rvert^2dx<\infty$$

\noindent The action of $G(F)$ on $L^2(H(F)\backslash G(F),\xi)$ is given by right translation. Since the triple $(G,H,\xi)$ is of a very particular form, the direct integral decomposition of $L^2(H(F)\backslash G(F),\xi)$ only involves tempered representations and moreover an irreducible tempered representation $\pi$ of $G(F)$ appears in this decomposition if and only if $m(\pi)=1$. It is thus natural for our problem to study this big representation $L^2(H(F)\backslash G(F),\xi)$. A function $f\in C_c^\infty(G(F))$ naturally acts on this space by

$$\displaystyle \left(R(f)\varphi\right)(x)=\int_{G(F)} f(g)\varphi(xg)dg,\;\;\; \varphi\in L^2(H(F)\backslash G(F),\xi),\; x\in G(F)$$

\noindent Moreover, this operator $R(f)$ is actually a kernel operator. More precisely, we have

$$\displaystyle \left(R(f)\varphi\right)(x)=\int_{H(F)\backslash G(F)} K_f(x,y)\varphi(y)dy,\;\;\; \varphi\in L^2(H(F)\backslash G(F),\xi),\; x\in G(F)$$

\noindent where

$$\displaystyle K_f(x,y)=\int_{H(F)} f(x^{-1}hy)\xi(h)dh,\;\;\; x,y\in G(F)$$

\noindent is the kernel function associated to $f$. In order to study the representation $L^2(H(F)\backslash G(F),\xi)$, we would like to compute the trace of $R(f)$ (because for example it would give some informations about the characters of the representations appearing in $L^2(H(F)\backslash G(F),\xi)$). Formally, we may write

\begin{align}\label{eq 0.6}
\displaystyle ``\Tr\; R(f)=\int_{H(F)\backslash G(F)} K(f,x)dx"
\end{align}

\noindent where $K(f,x)=K_f(x,x)$, $x\in H(F)\backslash G(F)$, is the restriction of the kernel to the diagonal. Unfortunately, neither of the two sides of the equality \ref{eq 0.6} makes sense in general: the operator $R(f)$ is not generally of trace class and the integral of the right hand side is not usually convergent. The first main step towards the proof of Theorem \ref{theorem 2} is to prove that nevertheless the expression in the right hand side of \ref{eq 0.6} still makes sense for a wide range of functions $f$. A function $f\in C_c^\infty(G(F))$ is said to be strongly cuspidal if for every proper parabolic subgroup $P=MU$ of $G$, we have

$$\displaystyle \int_{U(F)} f(mu)du=0,\;\;\; \mbox{for all } m\in M(F)$$

\noindent In Chapter \ref{section 7}, we prove the following (see Theorem \ref{proposition 7.1.1}).

\vspace{2mm}

\begin{them}\label{theorem 3}
For every strongly cuspidal function $f\in C_c^\infty(G(F))$, the integral

$$\displaystyle \int_{H(F)\backslash G(F)} K(f,x)dx$$

\noindent is absolutely convergent.
\end{them}

\vspace{2mm}

\noindent We actually prove more: we show that the above integral is absolutely convergent for every strongly cuspidal function in the Harish-Chandra Schwartz space $\mathcal{C}(G(F))$ rather than just $C_c^\infty(G(F))$. This seemingly technical detail is in fact rather important since in the real case the author was only able to construct enough strongly cuspidal functions in the space $\mathcal{C}(G(F))$ and not in $C_c^\infty(G(F))$.

\vspace{2mm}

\noindent Once we have Theorem \ref{theorem 3}, we can consider the distribution

$$\displaystyle f \mapsto J(f)=\int_{H(F)\backslash G(F)} K(f,x)dx$$

\noindent which is defined on the subspace $\mathcal{C}_{\scusp}(G(F))$ of strongly cuspidal functions in $\mathcal{C}(G(F))$. The next two steps toward the proof of Theorem \ref{theorem 2} are to give two rather different expressions for the distribution $J(.)$. The first expansion that we prove is spectral. It involves a natural space $\mathcal{X}(G)$ of tempered representations. In fact, elements of $\mathcal{X}(G)$ are not really tempered representations but rather virtual tempered representations. The space $\mathcal{X}(G)$ is build up from Arthur's elliptic representations of the group $G$ and of all of its Levi subgroups. We refer the reader to Section \ref{section 2.7} for a precise definition of $\mathcal{X}(G)$. Here, we only need to know that $\mathcal{X}(G)$ comes equipped with a natural measure $d\pi$ on it. For all $\pi\in \mathcal{X}(G)$, Arthur has defined a weighted character

$$f\in \mathcal{C}(G(F))\mapsto J_{M(\pi)}(\pi,f)$$

\noindent Here, $M(\pi)$ denotes the Levi subgroup from which the representation $\pi$ originates (more precisely, $\pi$ is parabolically induced from an elliptic representation of $M(\pi)$). When $M(\pi)=G$, the distribution $J_G(\pi,.)$ simply reduces to the usual character of $\pi$, that is $J_G(\pi,f)=\Tr\;\pi(f)$. When $M(\pi)\neq G$, the definition of the distribution $J_{M(\pi)}(\pi,.)$ is more involved and actually depends on some auxiliary choices (a maximal compact subgroup $K$ of $G(F)$ and some normalization of intertwining operators). However, it can be shown that the restriction of $J_{M(\pi)}(\pi,.)$ to $\mathcal{C}_{\scusp}(G(F))$ doesn't depend on any of these choices. For $f\in \mathcal{C}_{\scusp}(G(F))$, we define a function $\widehat{\theta}_f$ on $\mathcal{X}(G)$ by

$$\widehat{\theta}_f(\pi)=(-1)^{a_{M(\pi)}}J_{M(\pi)}(\pi,f),\;\;\; \pi\in \mathcal{X}(G)$$

\noindent where $a_{M(\pi)}$ is the dimension of $A_{M(\pi)}$ the maximal central split subtorus of $M(\pi)$. The spectral expansion of the distribution $J(.)$ now reads as follows (cf.\ Theorem \ref{theorem 9.1.1}):

\vspace{2mm}

\begin{them}\label{theorem 4}
For every strongly cuspidal function $f\in \mathcal{C}_{\scusp}(G(F))$, we have

$$\displaystyle J(f)=\int_{\mathcal{X}(G)} D(\pi)\widehat{\theta}_f(\pi)m(\pi)d\pi$$
\end{them}

\vspace{2mm}

\noindent The factor $D(\pi)$ appearing in the formula above is a certain determinant function which comes from Arthur's definition of elliptic representations. What is really important in the above spectral expansion of $J(.)$ is the appearance of the abstractly defined multiplicity $m(\pi)$. Its presence is due to the existence of an explicit description of the space $\Hom_H(\pi^\infty,\xi)$. More precisely, for $\pi$ an irreducible tempered representation of $G(F)$, we may define a certain hermitian form $\mathcal{L}_\pi$ on $\pi^\infty$ by

$$\displaystyle \mathcal{L}_\pi(e,e')=\int_{H(F)}^* (e,\pi(h)e')\xi(h)dh,\;\;\; e,e'\in \pi^\infty$$

\noindent The above integral is not necessarily absolutely convergent and needs to be regularized (cf.\ Section \ref{section 8.1}), it is why we put a star at the top of the integral sign. In any case, $\mathcal{L}_\pi$ is continuous and satisfies the intertwining relation

$$\displaystyle \mathcal{L}_\pi(\pi(h)e,\pi(h')e')=\xi(h)\overline{\xi(h')}\mathcal{L}_\pi(e,e'),\;\;\; e,e'\in \pi^\infty,\; h,h'\in H(F)$$

\noindent In particular, we see that for all $e'\in \pi^\infty$ the linear form $e\in \pi^\infty\mapsto \mathcal{L}_\pi(e,e')$ belongs to $\Hom_H(\pi^\infty,\xi)$. Hence, if $\mathcal{L}_\pi$ is not zero so is $m(\pi)$. In Chapter \ref{section 8}, we prove that the converse is also true. Namely, we have (cf.\ Theorem \ref{theorem 8.2.1})

\vspace{2mm}

\begin{them}\label{theorem 5}
For every irreducible tempered representation $\pi$ of $G(F)$, we have

$$\mathcal{L}_\pi\neq 0\Leftrightarrow m(\pi)\neq 0$$ 
\end{them}

\vspace{2mm}

\noindent This theorem has already been established in \cite{Beu1} when $F$ is $p$-adic (Th\'eor\`eme 14.3.1 of \cite{Beu1}). An analogous result for special orthogonal groups was proved previously by Waldspurger in \cite{Wa4} (Proposition 5.7) and then reproved in a different manner by Y. Sakellaridis and A. Venkatesh in \cite{SV} (Theorem 6.4.1) in a more general setting but under the additional assumption that the group is split. The proof given in \cite{Beu1} followed closely the treatment of Sakellaridis and Venkatesh whereas here we have been able to give an uniform proof in both the $p$-adic and the real case which is closer to the original work of Waldspurger.

\vspace{2mm}

\noindent As already explained, Theorem \ref{theorem 5} is a crucial step in the proof of the spectral expansion (Theorem \ref{theorem 4}). Actually, once Theorem \ref{theorem 5} is established, Theorem \ref{theorem 4} essentially reduces to the spectral expansion of Arthur's local trace formula \cite{A1} together with an argument allowing us to switch two integrals. This step is carried out in Chapter \ref{section 9}.

\vspace{2mm}

\noindent We now come to the geometric expansion of $J(.)$. It involves again the space of conjugacy classes $\Gamma(G,H)$ that appears in the formula for the multiplicity (Theorem \ref{theorem 2}). The other main ingredient is a function $c_f\colon G_{\ssi}(F)\to\mathbb{C}$ that is going to take the role played by the function $c_\pi$ in the multiplicity formula. The definition of $c_f$ involves the weighted orbital integrals of Arthur. Recall that for every Levi subgroup $M$ of $G$ and all $x\in M(F)\cap G_{\reg}(F)$, Arthur has defined a certain distribution

$$f\in \mathcal{C}(G(F))\mapsto J_M(x,f)$$

\noindent called a weighted orbital integral. If $M=G$, it simply reduces to the usual orbital integral

$$\displaystyle J_G(x,f)=\int_{G_x(F)\backslash G(F)} f(g^{-1}xg)dg,\;\;\; f\in \mathcal{C}(G(F))$$

\noindent When $M\neq G$, the distribution $J_M(x,.)$ depends on the choice of a maximal compact subgroup $K$ of $G(F)$. However, as for weighted characters, the restriction of $J_M(x,.)$ to the subspace $\mathcal{C}_{\scusp}(G(F))$ of strongly cuspidal functions doesn't depend on such a choice. For $f\in \mathcal{C}_{\scusp}(G(F))$, we define a function $\theta_f$ on $G_{\reg}(F)$ by

$$\theta_f(x)=(-1)^{a_{M(x)}}J_{M(x)}(x,f),\;\;\; x\in G_{\reg}(F)$$

\noindent where $M(x)$ denotes the minimal Levi subgroup of $G$ containing $x$ and $a_{M(x)}$ denotes, as before, the dimension of $A_{M(x)}$ the maximal central split subtorus of $M(x)$. The function $\theta_f$ is invariant and we can show that it shares a lot of the good properties that characters of representations have. It is what we call a quasi-character (cf.\ Chapter \ref{section 4}). In particular, as for characters, there is a natural extension of $\theta_f$ to a function

$$c_f\colon G_{\ssi}(F)\to \mathbb{C}$$

\noindent and we can now state the geometric expansion of $J(.)$ as follows (cf.\ Theorem \ref{theorem 11.2.1}).

\vspace{2mm}

\begin{them}\label{theorem 6}
For all strongly cuspidal functions $f\in \mathcal{C}_{\scusp}(G(F))$, we have

$$\displaystyle J(f)=\lim\limits_{s\to 0^+} \int_{\Gamma(G,H)} c_f(x)D^G(x)^{1/2}\Delta(x)^{s-1/2}dx$$
\end{them}

\vspace{2mm}

\noindent Once again, the expression of the right hand side of the equality above is absolutely convergent for all $s\in \mathbb{C}$ such that $Re(s)>0$ and the limit as $s\to 0^+$ exists (cf.\ Proposition \ref{proposition 11.1.1}).

\vspace{2mm}

\noindent It is from the equality between the two expansions of Theorem \ref{theorem 4} and Theorem \ref{theorem 6} that we deduce the formula for the multiplicity (Theorem \ref{theorem 2}). To be more precise, we first prove the spectral expansion (Theorem \ref{theorem 4}) and then we proceed to show the geometric expansion (Theorem \ref{theorem 6}) and the formula for the multiplicity (Theorem \ref{theorem 2}) together in a common inductive proof. The main reason for proceeding this way and not in a more linear order is that we use the spectral expansion together with the multiplicity formula for some ``smaller" GGP triples in order to show that the distribution $J(.)$ is supported in the elliptic locus $G(F)_{\elli}$ of $G(F)$. This fact is used crucially in the proof of Theorem \ref{theorem 6} and the author was not able to give an independent proof of it.

\vspace{2mm}

\noindent We now give a quick description of the content of each chapter. The fist two chapters are mainly intended to set up the notations, fix some normalizations and remind the reader of some well-known results. In particular, the second chapter contains the basic material we will be using on tempered representations. It includes a strong statement of the Harish-Chandra Plancherel theorem sometimes called matricial Paley-Wiener theorem which in the $p$-adic case is due to Harish-Chandra \cite{Wa2} and in the real case is due to Arthur \cite{A2}. In Chapter \ref{section 3}, we recall the Harish-Chandra technique of descent. There are two: descent from the group to the centralizer of one of its semi-simple elements (semi-simple descent) and descent from the group to its Lie algebra. In both cases, the descent takes the form of a map between some function spaces. We will be particularly concerned by the behavior of invariant differential operator (in the real case) under these two types of descent and we collect the relevant results there. The last Section (\ref{section 3.4}) is devoted to a third type of descent that we may call parabolic descent. However, we will be mainly interested in the dual of this map which allows us to ``induce" invariant distributions of Levi subgroups. In Chapter \ref{section 4}, we define the notion of quasi-characters and develop the main features of those functions that in many ways looks like characters. The main results of this chapter are, in the $p$-adic case, already contained in \cite{Wa1} and so we focus mainly on the real case. Chapter \ref{section 5} is devoted to the study of strongly cuspidal functions. In particular, it is in this chapter that we define the functions $\widehat{\theta}_f$ and $\theta_f$ to which we alluded before. This chapter also contains a version of Arthur's local trace formula for strongly cuspidal function (Theorems \ref{theorem 5.5.1} and \ref{theorem 5.5.2}). The proof of these two theorems, which are really just slight variation around Arthur's local trace formula, will appear elsewhere \cite{Beu2}. In Chapter \ref{section 6}, we define the GGP triples, the multiplicity $m(\pi)$ and we show some estimates that will be needed in the proof of the main theorems. Chapter \ref{section 8} is devoted to the proof of Theorem \ref{theorem 5}. In Chapter \ref{section 7}, we establish Theorem \ref{theorem 3} as well as an analog for the Lie algebra $\mathfrak{g}(F)$ of $G(F)$. This allows us to define two distributions $J(.)$ and $J^{\Lie}(.)$ on the group $G(F)$ and its Lie algebra respectively. Chapter \ref{section 9} concentrates on the spectral expansion of the distribution $J(.)$ (Theorem \ref{theorem 4}). As already explained, the main ingredient in the proof is Theorem \ref{theorem 5}. In Chapter \ref{section 10}, we establish some ``spectral" expansion for the Lie algebra analog $J^{\Lie}$. More precisely, we express $J^{\Lie}(f)$ in terms of weighted orbital integrals of $\widehat{f}$, the Fourier transform of the function $f$. Chapter \ref{section 11} contains the proofs of Theorem \ref{theorem 6} and Theorem \ref{theorem 2}. As we said, these two theorems are proved together. The last chapter, Chapter \ref{section 12}, is devoted to the proof of the main result of this paper (Theorem \ref{theorem 1}) following the outline given above. Finally, I collected in two appendices some definitions and results that are used throughout the text. Appendix \ref{section A} is concerned with locally convex topological vector spaces and particularly smooth and holomorphic maps taking values in such spaces whereas Appendix \ref{section B} contains some general estimates.

\bigskip

\noindent\ul{\textbf{Acknowledgements}}: I would like to thank Patrick Delorme for providing me with the unpublished work of Arthur \cite{A8}. I am particularly indebted to Jean-Loup Waldspurger for his very careful readings and comments on earlier versions of this paper. It will, I think, be obvious to the reader how much this paper owes to the pioneering work of Waldspurger and I would also like to thank him for sharing so generously his insights on this problem over the past few years. I thank Hang Xue for helpful discussions about the material to be found in Chapter \ref{section 8}. I am grateful to the referees for a thorough reading of a first version of this paper allowing to correct many inaccuracies and for the numerous suggestions to make the text more readable. This work was supported by the Gould Fund and by the National Science Foundation under agreement No. DMS-1128155. Any opinions, findings and conclusions or recommendations expressed in this material are those of the author and do not necessarily reflect the views of the National Science Foundation.

\section{Preliminaries}\label{section 1}

This is a preparatory chapter. We mainly set up notations, conventions and recall some standard results from the literature. In more details, Section \ref{section 1.1} fixes general notations, Section \ref{section 1.2} is devoted to a certain notion of norm on algebraic varieties over local fields due to Kottwitz \cite{Kott3} that we will use extensively, in Section \ref{section 1.3} we prove some useful estimates that will be needed later, in Section \ref{section 1.4} we introduce the most common spaces of functions that will appear in this paper, Section \ref{section 1.5} discusses the very important Harish-Chandra Schwartz space of functions and its basic properties, Section \ref{section 1.6} explains our normalizations of measures, in Section \ref{section 1.8} we introduce some spaces of conjugacy classes that we equip with topologies and measures, in Section \ref{section 1.7} we set up notations for orbital integrals and recall some of their properties, finally in Sections \ref{section 1.9} and \ref{section 1.10} we recall Arthur's notions of {\em $(G,M)$-family} and {\em weighted orbital integral} respectively.

\subsection{General notation and conventions}\label{section 1.1}

\noindent Throughout this paper, we fix a field \gls{F} which is either $p$-adic (i.e., a finite extension of $\mathbb{Q}_p$) or $\mathbb{R}$ the field of real numbers. We denote by $\lvert .\rvert$ the normalized absolute value on $F$ i.e., for every Haar measure $dx$ on $F$ we have $d(ax)=\lvert a\rvert dx$, for all $a\in F$. We fix once and for all an algebraic closure $\gls{Fbar}$ for $F$ and let $\gls{GammaF}=\Gal(\overline{F}/F)$ be the corresponding absolute Galois group. We will also denote by $\lvert .\rvert$ the unique extension of the absolute value to $\overline{F}$. All varieties, schemes, algebraic groups will be assumed, unless otherwise specified, to be defined over $F$. Moreover we will identify any algebraic variety $X$ defined over $F$ or $\overline{F}$ with its set of $\overline{F}$-points. For $G$ a locally compact separable group (for example the $F$-points of an algebraic group defined over $F$), we will usually denote by $d_L g$ (resp.\ $d_R g$) a left (resp.\ a right) Haar measure on $G$. If the group is unimodular then we will usually denote both by $dg$. Finally $\delta_G$ will stand for the modular character of $G$ that is defined by $d_L(gg'^{-1})=\delta_G(g')d_Lg$ for all $g\in G$.

\vspace{2mm}

\noindent We fix, until the end of Chapter \ref{section 5}, a connected reductive group $G$ over $F$. We denote by $\mathfrak{g}$ its Lie algebra and by

$$G\times \mathfrak{g}\to \mathfrak{g}$$
$$(g,X)\mapsto gXg^{-1}$$

\noindent the adjoint action. A sentence like ``Let $P=MU$ be a parabolic subgroup of $G$" will mean as usual that $P$ is defined over $F$, $U$ is its unipotent radical and that $M$ is a Levi component of $P$ defined over $F$. We define an integer \gls{dg} by

$$\delta(G)=\dim(G)-\dim(T)$$

\noindent where $T$ is any maximal torus of $G$. It is also the dimension of any regular conjugacy class in $G$.

\vspace{2mm}

\noindent Let us recall some of the usual objects attached to $G$. We shall denote by \gls{zg} the center of $G$ and by \gls{Ag} its split component. We define the real vector space

$$\gls{cAg}=\Hom(X^*(G),\mathbb{R})$$

\noindent and its dual

$$\mathcal{A}_G^*=X^*(G)\otimes \mathbb{R}$$

\noindent where $X^*(G)$ stands for the module of $F$-rational characters of $G$. We have a natural homomorphism

$$\gls{Hg}\colon G(F)\to \mathcal{A}_G$$

\noindent given by

$$\displaystyle \langle \chi,H_G(g)\rangle=\log\left(\lvert \chi(g)\rvert\right),\;\; g\in G(F),\; \chi\in X^*(G)$$

\noindent Set $\gls{cAgf}=H_G(G(F))$ and $\gls{ctAgf}=H_G(A_G(F))$. In the real case, we have $\widetilde{\mathcal{A}}_{G,F}=\mathcal{A}_{G,F}=\mathcal{A}_G$. In the $p$-adic case, $\widetilde{\mathcal{A}}_{G,F}$ and $\mathcal{A}_{G,F}$ are both lattices inside $\mathcal{A}_G$. We also set $\gls{cAvgf}=\Hom(\mathcal{A}_{G,F},2\pi\mathbb{Z})$ and $\gls{ctAvgf}=\Hom(\widetilde{\mathcal{A}}_{G,F},2\pi\mathbb{Z})$. In the $p$-adic case, $\widetilde{\mathcal{A}}_{G,F}^\vee$ and $\mathcal{A}_{G,F}^\vee$ are this time lattices inside $\mathcal{A}_G^*$ whereas in the real case we have $\widetilde{\mathcal{A}}_{G,F}^\vee=\mathcal{A}_{G,F}^\vee=0$. We set $\gls{cA*gf}=\mathcal{A}^*_G/\mathcal{A}_{G,F}^\vee$ and we identify $i\mathcal{A}_{G,F}^*$ with the group of unitary unramified characters of $G(F)$ by mean of the pairing $(\lambda,g)\in i\mathcal{A}_{G,F}^*\times G(F)\mapsto e^{\langle \lambda, H_G(g)\rangle}$. We shall also denote by $\gls{cAgc}$ and $\gls{cA*gc}$ the complexifications of $\mathcal{A}_G$ and $\mathcal{A}^*_G$.

\vspace{2mm}

\noindent Let $P=MU$ be a parabolic subgroup of $G$. Of course, the previous constructions apply to $M$. We will denote by $\gls{RAMP}$ the set of roots of $A_M$ in the unipotent radical of $P$. If $K$ is a maximal compact subgroup of $G(F)$ which is special in the $p$-adic case, we have the Iwasawa decomposition $G(F)=M(F)U(F)K$. We may then choose maps $m_P\colon G(F)\to M(F)$, $u_P\colon G(F)\to U(F)$ and $k_P\colon G(F)\to K$ such that $g=m_P(g)u_P(g)k_P(g)$ for all $g\in G(F)$. We then extend the homomorphism $H_M$ to a map $\gls{HP}\colon G(F)\to \mathcal{A}_M$ by setting $H_P(g)=H_M(m_P(g))$. This extension depends of course on the maximal compact $K$ but its restriction to $P(F)$ doesn't and is given by $H_P(mu)=H_M(m)$ for all $m\in M(F)$ and all $u\in U(F)$. By a Levi subgroup of $G$ we mean a subgroup of $G$ which is the Levi component of some parabolic subgroup of $G$. We will also use Arthur's notation: if $M$ is a Levi subgroup of $G$, then we denote by $\gls{PM}$, $\gls{LM}$ and $\gls{FM}$ the finite sets of parabolic subgroups admitting $M$ as a Levi component, of Levi subgroups containing $M$ and of parabolic subgroups containing $M$ respectively. If $M\subset L$ are two Levi subgroups, we set $\gls{AML}=\mathcal{A}_M/\mathcal{A}_L$. We have a canonical decomposition

$$\mathcal{A}_M=\mathcal{A}_L\oplus \mathcal{A}_M^L$$

\noindent and its dual

$$\mathcal{A}_M^*=\mathcal{A}_L^*\oplus \left(\mathcal{A}_M^L\right)^*$$

\noindent If $H$ is an algebraic group, we shall denote by $H^0$ its neutral connected component. For $x\in G$ (resp.\ $X\in \mathfrak{g}$), we denote by $\gls{ZGx}$ (resp.\ $\gls{ZGX}$) the centralizer of $x$ (resp.\ $X$) in $G$ and by $\gls{Gx}=Z_G(x)^0$ (resp.\ $\gls{GX}=Z_G(X)^0$) the neutral component of the centralizer. Recall that if $X\in \mathfrak{g}$ is semi-simple then $Z_G(X)=G_X$. We will denote by $\gls{Gss}$ and $\gls{Greg}$ (resp.\ $\gls{gss}$ and $\gls{greg}$) the subsets of semi-simple and regular semi-simple elements in $G$ (resp.\ in $\mathfrak{g}$). For any subset $A$ of $G(F)$ (resp.\ of $\mathfrak{g}(F)$), we will denote by $A_{\reg}$ the intersection $A\cap G_{\reg}(F)$ (resp.\ $A\cap \mathfrak{g}_{\reg}(F)$) and by $A_{\ssi}$ the intersection $A\cap G_{\ssi}(F)$ (resp.\ $A\cap \mathfrak{g}_{\ssi}(F)$). We will usually denote by $\gls{TG}$ a set of representatives for the conjugacy classes of maximal tori in $G$. Recall that a maximal torus $T$ of $G$ is said to be elliptic if $A_T=A_G$. Elliptic maximal tori always exist in the $p$-adic case but not necessarily in the real case. An element $x\in G(F)$ will be said to be {\em elliptic} if it belongs to some elliptic maximal torus (in particular it is semi-simple). Similarly, an element $X\in \mathfrak{g}(F)$ will be said to be {\em elliptic} if it belongs to the Lie algebra of some elliptic maximal torus. We will denote by $\gls{Gell}$ and $\gls{gell}$ the subsets of elliptic elements in $G(F)$ and $\mathfrak{g}(F)$ respectively. We will also set $\gls{Gregell}=G(F)_{\elli}\cap G_{\reg}(F)$ and $\gls{gregell}=\mathfrak{g}(F)_{\elli}\cap \mathfrak{g}_{\reg}(F)$. For all $x\in G_{\ssi}(F)$ (resp.\ all $X\in \mathfrak{g}_{\ssi}(F)$), we set

$$\displaystyle \gls{DGx}=\left\lvert \det(1-\Ad(x))_{\mid \mathfrak{g}/\mathfrak{g}_x}\right\rvert \;\;\; \left(\mbox{resp. } \gls{DGX}=\left\lvert \det \;\ad(X)_{\mid \mathfrak{g}/\mathfrak{g}_X}\right\rvert\right)$$

\noindent If a group $H$ acts on a set $X$ and $A$ is a subset of $X$, we shall denote by $\No_H(A)$ the normalizer of $A$ in $H$. For every Levi subgroup $M$ and every maximal torus $T$ of $G$, we will denote by $\gls{WGM}$ and $\gls{WGT}$ the Weyl groups of $M(F)$ and $T(F)$ respectively, that is

$$\displaystyle W(G,M)=\No_{G(F)}(M)/M(F)\;\;\mbox{ and }\;\; W(G,T)=\No_{G(F)}(T)/T(F)$$

\noindent If $f$ is a function on either $G(F)$ or $\mathfrak{g}(F)$, for all $g\in G(F)$ we will denote by $\gls{gf}$ the function $f\circ \Ad(g)$. We shall denote by $\gls{R}$ and $\gls{L}$ the natural actions of $G(F)$ on functions on $G(F)$ given by right and left translation respectively. That is $\left(R(g)f\right)(\gamma)=f(\gamma g)$ and $\left(L(g)f\right)(\gamma)=f(g^{-1}\gamma)$ for every function $f$ and all $g,\gamma\in G(F)$.

\vspace{2mm}

\noindent Assume that $F=\mathbb{R}$. Then, we will denote by $\gls{Ug}$ the enveloping algebra of $\mathfrak{g}=\mathfrak{g}(\mathbb{C})$ and by $\gls{Zg}$ its center. The right and left actions of $G(F)$ on smooth functions on $G(F)$ of course extend to $\mathcal{U}(\mathfrak{g})$. We still denote by $R$ and $L$ these actions. For $z\in \mathcal{Z}(\mathfrak{g})$, we will simply set $zf=R(z)f$ ($=L(z^*)f$ with a notation introduced below) for all $f\in C^\infty(G(F))$. We will also denote by $\gls{Sg}$ and $\gls{Sg*}$ the symmetric algebras of $\mathfrak{g}$ and $\mathfrak{g}^*$ respectively. We will identify $S(\mathfrak{g}^*)$ with the algebra of complex-valued polynomial functions on $\mathfrak{g}$ and we will identify $S(\mathfrak{g})$ with the algebra of differential operators on $\mathfrak{g}$ with constant coefficients. More precisely for $u\in S(\mathfrak{g})$, we shall denote by $\gls{partu}$ the corresponding differential operator. We denote by $\gls{Ig}$ and $\gls{Ig*}$ the subalgebras of $G$-invariant elements in $S(\mathfrak{g})$ and $S(\mathfrak{g}^*)$ respectively. We will also need the algebra $\gls{Dg}$ of differential operators with polynomial coefficients on $\mathfrak{g}(F)$. We will denote by $u\mapsto \gls{u*}$ the unique $\mathbb{C}$-algebra automorphism of both $\mathcal{U}(\mathfrak{g})$ and $S(\mathfrak{g})$ that sends every $X\in \mathfrak{g}$ to $-X$. We then have

$$\displaystyle \int_{G(F)} \left(R(u)f_1\right)(g) f_2(g)dg=\int_{G(F)} f_1(g) \left(R(u^*)f_2\right)(g) dg,\;\mbox{ for all } f_1,f_2\in C_c^\infty(G(F)), u\in \mathcal{U}(\mathfrak{g})$$

$$\displaystyle \int_{\mathfrak{g}(F)} \left(\partial(u)f_1\right)(X) f_2(X)dX=\int_{\mathfrak{g}(F)} f_1(X) \left(\partial(u^*)f_2\right)(X) dX,\;\mbox{ for all } f_1,f_2\in C_c^\infty(\mathfrak{g}(F)), u\in S(\mathfrak{g})$$

\noindent For each maximal torus $T$ of $G$, the Harish-Chandra homomorphism provides us with an isomorphism

$$\mathcal{Z}(\mathfrak{g})\simeq S(\mathfrak{t})^{W(G_{\mathbb{C}},T_{\mathbb{C}})}$$
$$z\mapsto \gls{zT}$$

\noindent where $W(G_{\mathbb{C}},T_{\mathbb{C}})$ denotes the Weyl group of $T_{\mathbb{C}}$ in $G_{\mathbb{C}}$. Assume that $H$ is a connected reductive subgroup of $G$ of the same rank as $G$, for example $H$ can be a Levi subgroup or the connected centralizer of a semi-simple element. Let $T\subset H$ be a maximal torus. Since $W(H_{\mathbb{C}},T_{\mathbb{C}})\subset W(G_{\mathbb{C}},T_{\mathbb{C}})$, the Harish-Chandra isomorphism for $T$ induces an injective homomorphism

$$\mathcal{Z}(\mathfrak{g})\hookrightarrow \mathcal{Z}(\mathfrak{h})$$
$$z\mapsto \gls{zH}$$

\noindent which is in fact independent of $T$ and such that the extension $\mathcal{Z}(\mathfrak{h})/\mathcal{Z}(\mathfrak{g})$ is finite. Similarly, over the Lie algebra we have isomorphisms

$$I(\mathfrak{g}^*)\simeq S(\mathfrak{t}^*)^{W(G_{\mathbb{C}},T_{\mathbb{C}})}\;\;\; I(\mathfrak{g})\simeq S(\mathfrak{t})^{W(G_{\mathbb{C}},T_{\mathbb{C}})}$$
$$p\mapsto \gls{pT}\;\;\;\;\;\;\;\;\;\;\;\;\;\;\;\;\;\;\;\; u\mapsto \gls{uT}\;\;\;\;\;\;$$

\noindent The first one is just the ``restriction to $\mathfrak{t}$" homomorphism, the second one may be deduced from the first one once we choose a $G$-equivariant isomorphism $\mathfrak{g}\simeq \mathfrak{g}^*$ (but it doesn't depend on such a choice). As for the group, if $H$ is a connected reductive subgroup of $G$ of the same rank as $G$, we deduce from these isomorphisms two injective homomorphisms

$$I(\mathfrak{g}^*)\hookrightarrow I(\mathfrak{h}^*)\;\;\; I(\mathfrak{g})\hookrightarrow I(\mathfrak{h})$$
$$\;\;\;\;\;\;\;\;\;\;\;p\mapsto \gls{pH}\;\;\;\;\;\;\;\;\;\;\;\; u\mapsto \gls{uH}\;\;\;\;\;\;$$

\noindent which are such that the extensions $I(\mathfrak{h}^*)/I(\mathfrak{g}^*)$ and $I(\mathfrak{h})/I(\mathfrak{g})$ are finite. Also the two isomorphisms $\mathcal{Z}(\mathfrak{g})\simeq S(\mathfrak{t})^{W(G_{\mathbb{C}},T_{\mathbb{C}})}$ and $I(\mathfrak{g})\simeq S(\mathfrak{t})^{W(G_{\mathbb{C}},T_{\mathbb{C}})}$ induce an isomorphism $\mathcal{Z}(\mathfrak{g})\simeq I(\mathfrak{g})$ that we shall denote by $z\mapsto u_z$.

\vspace{2mm}

\noindent We will also adopt the following slightly imprecise but convenient notation. If $f$ and $g$ are positive functions on a set $X$, we will write

\begin{center}
$f(x)\ll g(x)$ for all $x\in X$
\end{center}

\noindent and we will say that $f$ is essentially bounded by $g$, if there exists a $c>0$ such that

\begin{center}
$f(x)\leqslant cg(x)$, for all $x\in X$
\end{center}

\noindent We will also say that $f$ and $g$ are equivalent and we will write

\begin{center}
$f(x)\sim g(x)$ for all $x\in X$
\end{center}

\noindent if both $f$ is essentially bounded by $g$ and $g$ is essentially bounded by $f$.

\vspace{2mm}

\subsection{Reminder of norms on algebraic varieties}\label{section 1.2}

All along this paper, we will assume that $\mathfrak{g}=\mathfrak{g}(\overline{F})$ has been equipped with a (classical) norm $\lvert .\rvert_{\mathfrak{g}}$, that is a map $\lvert .\rvert_{\mathfrak{g}}:\mathfrak{g}\to \mathbb{R}_+$ satisfying $\lvert \lambda X\rvert_{\mathfrak{g}}=\lvert \lambda\rvert. \lvert X\rvert_{\mathfrak{g}}$, $\lvert X+Y\rvert_{\mathfrak{g}}\leqslant \lvert X\rvert_{\mathfrak{g}} +\lvert Y\rvert_{\mathfrak{g}}$ and $\lvert X\rvert_{\mathfrak{g}}=0$ if and only if $X=0$ for all $\lambda\in F$ and $X,Y\in \mathfrak{g}$. For any $R>0$, we will denote by $\gls{B0r}$ the closed ball of radius $R$ centered at the origin in $\mathfrak{g}(F)$.

\vspace{2mm}

\noindent We will make an heavy use of the notion of norm on varieties over local field introduced by Kottwitz in \cite{Kott3}. Actually, we will use a slight variation of Kottwitz's norms that is more convenient for us and that we will call log-norms because these are essentially logarithms of Kottwitz's norms. For the convenience of the reader, we will recall here the definitions and main features of these log-norms.

\vspace{2mm}

\noindent First, an {\em abstract log-norm} on a set $X$ is just a real-valued function $x\mapsto \sigma(x)$ on $X$ such that $\sigma(x)\geqslant 1$, for all $x\in X$. For two abstract log-norms $\sigma_1$ and $\sigma_2$ on $X$, we will say that $\sigma_2$ {\em dominates} $\sigma_1$ if
 
$$\sigma_1(x)\ll \sigma_2(x)$$

\noindent for all $x\in X$ in which case we shall write $\sigma_1\ll \sigma_2$. We will say that $\sigma_1$ and $\sigma_2$ are {\em equivalent} if each of them dominates the other and in this case we will write $\sigma_1\sim \sigma_2$.

\vspace{2mm}

\noindent Let $X$ be an affine algebraic variety over $\overline{F}$ and denote by $\mathcal{O}(X)$ its ring of regular functions. Choosing a set of generators $f_1,\ldots,f_m$ of the $\overline{F}$-algebra $\mathcal{O}(X)$, we can define an abstract log-norm $\gls{sigmaX}$ on $X$ by setting

$$\sigma_X(x)=1+\log\left(max\{1,\lvert f_1(x)\rvert,\ldots,\lvert f_m(x)\rvert\}\right)$$

\noindent for all $x\in X$. The equivalence class of $\sigma_X$ doesn't depend on the choice of $f_1,\ldots,f_m$ and by a {\em log-norm on $X$} we will mean any abstract log-norm in this equivalence class. Note that if $U$ is the principal Zariski open subset of $X$ defined by the non-vanishing of $Q\in \mathcal{O}(X)$, then we have

$$\sigma_U(x)\sim \sigma_X(x)+\log\left(2+\lvert Q(x)\rvert^{-1}\right)$$

\noindent for all $x\in U$.

\vspace{2mm}

\noindent More generally, let $X$ be any algebraic variety over $\overline{F}$. Choose a finite covering $\left(U_i\right)_{i\in I}$ of $X$ by open affine subsets and fix log-norms $\sigma_{U_i}$ on $U_i$, $i\in I$. Then

$$\sigma_X(x)=\inf\{\sigma_{U_i}(x);i\in I \; \mbox{ such that } x\in U_i\}$$

\noindent defines an abstract log-norm on $X$ the equivalence class of which doesn't depend on the various choices. An abstract log-norm in this equivalence class will be just called a {\em log-norm on $X$}.

\vspace{2mm}

\noindent Let $f\colon X\to Y$ be a morphism of algebraic varieties over $\overline{F}$ and $\sigma_Y$ be a log-norm on $Y$. We define the abstract log-norm $\gls{fsigmaY}$ on $X$ by

$$f^*\sigma_Y(x)=\sigma_Y(f(x))$$

\noindent for all $x\in X$. The following lemma will be used without further notice throughout the text (cf.\ Proposition 18.1(1) of \cite{Kott3})

\begin{lem}\label{lemma 1.2.1}
Let $\sigma_X$ be a log-norm on $X$. Then $f^*\sigma_Y\ll \sigma_X$. If $f$ is moreover a finite morphism (in particular if it is a closed immersion), then $f^*\sigma_Y\sim \sigma_X$.
\end{lem}

\noindent Let this time $f\colon X\to Y$ be a morphism of algebraic varieties over $F$ and let $\sigma_X$ be a log-norm on $X$ (but we will only consider its restriction to $X(F)$). Define an abstract log-norm $\gls{fsigmaX}$ on $\Ima(X(F)\to Y(F))$ by

$$f_*\sigma_X(y)=\inf_{x\in X(F);\; f(x)=y} \sigma_X(x)$$

\noindent Let $\sigma_Y$ be a log-norm on $Y$. By the previous lemma, $f_*\sigma_X$ dominates $\sigma_Y$ (as abstract log-norms on $\Ima(X(F)\to Y(F))$). We say that $f$ has {\em the norm descent property} if $\sigma_Y$ and $f_*\sigma_X$ are equivalent as abstract log-norms on $\Ima(X(F)\to Y(F))$. Of course, if $f^*\sigma_Y$ and $\sigma_X$ are equivalent, then $f$ has the norm descent property, and so this is the case in particular if $f$ is finite. We state here the basic facts we will be using regarding to the norm descent property (cf.\ \cite{Kott3} Proposition 18.2).

\begin{lem}\label{lemma 1.2.2}
\begin{enumerate}[(i)]

\item The norm descent property is local on the basis. In other words if $f\colon X\to Y$ is a morphism of algebraic varieties over $F$ and $\left(U_i\right)_{i\in I}$ is a finite covering by Zariski-open subsets of $Y$ defined over $F$, then $f$ has the norm descent property if and only if each of the $f_i\colon  f^{-1}(U_i)\to U_i$, $i\in I$, has the norm descent property.

\item If $f$ admits a section, then it has the norm descent property.
\end{enumerate}
\end{lem}

\noindent We will also need the following nontrivial result (cf.\ \cite{Kott3} Proposition 18.3).

\begin{prop}\label{proposition 1.2.1}
Let $G$ be a connected reductive group over $F$ and $T$ an $F$-subtorus of $G$. Then the morphism $G\to T\backslash G$ has the norm descent property.
\end{prop}

\noindent In order not to confuse the reader, we remark that this last proposition is a statement of equivalence of log-norms on $T(F)\backslash G(F)$ (which is the image of $G(F)$ in $(T\backslash G)(F)$) and not on the full of $(T\backslash G)(F)$ (which in general can be slightly bigger).

\vspace{2mm}

\noindent We will assume that all algebraic varieties $X$ in this article (be they defined over $F$ or $\overline{F}$) are equipped with a log-norm $\sigma_X$. Note that if $X=V$ is a vector space over $\overline{F}$, then we may take

$$\sigma_V(v)= \log\left(2+\lvert v\rvert\right),\;\; v\in V$$

\noindent where $\lvert .\vert$ is a classical norm on $V$. We will usually assume that it is the case. Also, we will denote $\sigma_G$ simply by $\gls{sigma}$ and all closed subvarieties of $G$ will be equipped with the log-norm obtain by restriction of $\sigma$. Note that we have

$$\sigma(xy)\ll \sigma(x)+\sigma(y)\ll \sigma(x)\sigma(y)$$

\noindent for all $(x,y)\in G\times G$. It follows from the last proposition that we may assume, and we will throughout the paper, that we have

\begin{align}\label{eq 1.2.1}
\displaystyle \sigma_{T\backslash G}(g)=\inf_{t\in T(F)}\sigma(tg)
\end{align}

\noindent for all $g\in G(F)$. We will also need the following

\vspace{3mm}

\begin{num}
\item\label{eq 1.2.2} For every maximal torus $T\subset G$, we have
$$\sigma(g^{-1}Xg)+\log\left(2+D^G(X)^{-1}\right)\sim \sigma_{\mathfrak{g}}(X)+\sigma_{T\backslash G}(g)+\log\left(2+D^G(X)^{-1}\right)$$
for all $X\in \mathfrak{t}_{\reg}(F)$ and all $g\in G(F)$.
\end{num}

\vspace{3mm}

\noindent Indeed, this follows from Lemma \ref{lemma 1.2.1} and the fact that the regular map $T\backslash G\times \mathfrak{t}_{\reg}\to \mathfrak{g}_{\reg}$, $(g,X)\mapsto g^{-1}Xg$ is finite.

\vspace{2mm}

\noindent As we said, our log-norms are essentially the logarithm of Kottwitz's norms. For $G$ and its Lie algebra $\mathfrak{g}$, it will be convenient at some points to work with norms instead of log-norms. We therefore set

$$\displaystyle \gls{gnorm}=e^{\sigma(g)},\;\;\; g\in G$$

$$\displaystyle \gls{Xnorm}=e^{\sigma_{\mathfrak{g}}(X)},\;\;\; X\in \mathfrak{g}$$

\noindent Let $X$ be an algebraic variety over $F$ on which a log-norm $\sigma_X$ has been chosen. We will also use the following notation

$$\displaystyle \gls{X<C}:=\{x\in X(F);\; \sigma_X(x)<C \}$$

$$\displaystyle \gls{X>C}:=\{x\in X(F);\; \sigma_X(x)\geqslant C \}$$

\noindent for all $C>0$. We have the two following estimates. The first one is easy to prove and the second one is due to Harish-Chandra (cf.\ Theorem 9 p.37 of \cite{Va}).

\vspace{3mm}

\begin{num}
\item\label{eq 1.2.3} Let $V$ be a finite dimensional $F$-vector space and $Q$ be a polynomial function on $V$. Set $V'=\{Q\neq 0\}$ and let $\omega_V\subset V$ be a relatively compact subset. Then, for all $k\geqslant 0$, there exists $\delta>0$ such that
$$\displaystyle \int_{V'[\leqslant \epsilon]\cap \omega_V} \log\left(2+\lvert Q(X)\rvert^{-1}\right)^k dX\ll \epsilon^\delta$$
for all $\epsilon>0$.

\item\label{eq 1.2.4} Let $T\subset G$ be a maximal torus. Then, for all $N_0>0$ there exists $N>0$ such that
$$\displaystyle D^G(X)^{1/2}\int_{T(F)\backslash G(F)} \lVert g^{-1}Xg\rVert_{\mathfrak{g}}^{-N}dg\ll \lVert X\rVert_{\mathfrak{g}}^{-N_0}$$
for all $X\in \mathfrak{t}_{\reg}(F)$.
\end{num}

\subsection{A useful lemma}\label{section 1.3}

\noindent Fix a minimal parabolic subgroup $\overline{P}_{\mini}=M_{\mini}\overline{U}_{\mini}$ of $G$ and let $A_{\mini}=A_{M_{\mini}}$ be the maximal split central subtorus of $M_{\mini}$. Set

$$\displaystyle \gls{Amin+}=\{a\in A_{\mini}(F); \lvert \alpha(a)\rvert\geqslant 1 \; \forall \alpha\in R(A_{\mini},\overline{P}_{\mini})\}$$

\noindent Let $\overline{Q}=LU_{\overline{Q}}$ be a parabolic subgroup containing $\overline{P}_{\mini}$, where $L$ is the unique Levi component such that $M_{\mini}\subseteq L$. For all $\delta>0$, we define

$$\displaystyle \gls{AQ+}:=\{a\in A_{\mini}^+; \lvert\alpha(a)\rvert\geqslant e^{\delta\sigma(a)}\; \forall \alpha\in R(A_{\mini},U_{\overline{Q}})\}$$

\noindent Let $Q=LU_Q$ be the parabolic subgroup opposite to $\overline{Q}$ with respect to $L$.

\begin{lem}\label{proposition 1.3.1}
\begin{enumerate}[(i)]
\item Let $\epsilon>0$ and $\delta>0$. Then, we have an inequality
$$\displaystyle \sigma(a)\ll \sup\left(\sigma(g),\sigma(a^{-1}ga)\right)$$
for all $a\in A_{\mini}^{\overline{Q},+}(\delta)$ and all $g\in G(F)\smallsetminus \left(\overline{Q}\bigl[ <\epsilon\sigma(a) \bigr] aU_Q\bigl[ <\epsilon\sigma(a) \bigr] a^{-1}\right)$.

\item Let $0<\delta'<\delta$ and $c_0>0$. Then for $\epsilon>0$ sufficiently small, we have
$$aU_Q\left[<\epsilon \sigma(a)\right]a^{-1}\subseteq \exp\left(B(0,c_0e^{-\delta'\sigma(a)})\cap \mathfrak{u}_Q(F)\right)$$
for all $a\in A_{\mini}^{\overline{Q},+}(\delta)$.
\end{enumerate}

\end{lem}

\vspace{2mm}

\noindent\ul{Proof}: 
\begin{enumerate}[(i)]
\item Let us set

$$W^G=W(G,M_{\mini}),\; W^{L}=W(L,M_{\mini})$$

\noindent for the Weyl groups of $M_{\mini}$ in $G$ and $L$ respectively. We have the following decomposition:

\begin{align}\label{eq 1.3.1}
\displaystyle G=\bigcup_{w\in W^{L}\backslash W^G} \overline{Q}U_Qw
\end{align}

\noindent Indeed, if we let $P_{\mini}=M_{\mini}U_{\mini}$ be the parabolic subgroup opposite to $\overline{P}_{\mini}$ with respect to $M_{\mini}$, we have $\overline{P}_{\mini}U_{\mini}\subseteq \overline{Q}U_Q$. Thus \ref{eq 1.3.1} is clearly a consequence of the decomposition

$$\displaystyle G=\bigcup_{w\in W^G}\overline{P}_{\mini}U_{\mini}w$$

\noindent which itself follows from the Bruhat decomposition

$$\displaystyle G=\bigsqcup_{w\in W^G}\overline{P}_{\mini}wU_{\mini}$$

\noindent and the fact that $\overline{P}_{\mini}wU_{\mini}\subseteq \overline{P}_{\mini}U_{\mini}w$, for all $w\in W^G$.

\vspace{2mm}

\noindent Fix a set $\mathcal{W}\subset W^G$ of representatives of the left $W^L$-cosets in $W^G$ and assume (as we may) that $1\in \mathcal{W}$. Set $\mathcal{U}_w=\overline{Q}U_Qw$ for all $w\in \mathcal{W}$. These are affine open subsets of $G$ which are all naturally isomorphic to $\overline{Q}\times U_Q$. For all $w\in\mathcal{W}$, we may define a log-norm $\sigma_w$ on $\mathcal{U}_w$ by

\begin{align}\label{eq 1.3.2}
\displaystyle \sigma_w(\overline{q}uw)=\sup\{\sigma(\overline{q}),\sigma(u)\}
\end{align}

\noindent for all $\overline{q}\in \overline{Q}$ and all $u\in U_Q$. By \ref{eq 1.3.1} the family $\left(\mathcal{U}_w\right)_{w\in \mathcal{W}}$ is a Zariski open cover of $G$ and so we have

\begin{align}\label{eq 1.3.3}
\displaystyle \sigma(g)\sim\inf\{\sigma_w(g);\; w\in \mathcal{W} \mbox{ such that } g\in \mathcal{U}_w\}
\end{align}

\noindent for all $g\in G$. 

\vspace{2mm}

\noindent Obviously we may assume that $\epsilon$ as small as we want: if we replace $\epsilon$ by another positive constant $\epsilon'<\epsilon$ then the assertion of the proposition becomes stronger. Moreover, the following is easy to see

\vspace{3mm}

\begin{num}
\item\label{eq 1.3.4} If $\epsilon$ is sufficiently small (depending on $\delta$), there exists a bounded subset $C_{\overline{Q}}\subseteq \overline{Q}$ such that
$$a^{-1}\overline{Q}\bigl[ <\epsilon\sigma(a) \bigr] a\subseteq C_{\overline{Q}}L(F)$$
for all $a\in A_{\mini}^{\overline{Q},+}(\delta)$.
\end{num}

\vspace{3mm}

\noindent Henceforth, we will assume $\epsilon$ sufficiently small so that it satisfies \ref{eq 1.3.4}. For all $w\in\mathcal{W}$ and $c>0$, we set

$$\displaystyle \mathcal{U}_w[<c]=\{g\in\mathcal{U}_w;\; \sigma_w(g)<c\}$$

\noindent We now show the following

\vspace{3mm}

\begin{num}
\item\label{eq 1.3.5} Let $w\in\mathcal{W}$. Then, we have
$$\sigma(a^{-1}ga)\sim \sigma_w(a^{-1}ga)$$
for all $a\in A_{\mini}^{\overline{Q},+}(\delta)$ and all $g\in\mathcal{U}_w[<\epsilon \sigma(a)]$.
\end{num}

\vspace{3mm}

\noindent Indeed, by \ref{eq 1.3.4} we have

$$\displaystyle a^{-1}\mathcal{U}_w[<\epsilon \sigma(a)]a\subseteq C_{\overline{Q}}L(F)U_Q(F)w$$

\noindent for all $a\in A_{\mini}^{\overline{Q},+}(\delta)$ and moreover $\sigma$ and $\sigma_w$ are equivalent on $C_{\overline{Q}}L(F)U_Q(F)w$. The point \ref{eq 1.3.5} follows.

\vspace{2mm}

\noindent By \ref{eq 1.3.3}, there is an inequality

$$\sigma(a)\ll \sigma(g)$$

\noindent for all $a\in A_{\mini}$ and all $g\in G\backslash \bigcup_{w\in\mathcal{W}} \mathcal{U}_w[<\epsilon \sigma(a)]$. Combining this with \ref{eq 1.3.5}, we see that the estimate of the proposition is a consequence of the following claims:

\vspace{2mm}

\begin{num}
\item\label{eq 1.3.6} We have
$$\sigma(a)\ll\sigma_1(a^{-1}ga)$$
for all $a\in A_{\mini}^{\overline{Q},+}(\delta)$ and all $g\in \mathcal{U}_1[<\epsilon \sigma(a)]\big\backslash \left(\overline{Q}\left[<\epsilon \sigma(a)\right] aU_Q\left[<\epsilon \sigma(a)\right]a^{-1}\right)$.
\end{num}

\vspace{2mm}

\begin{num}
\item\label{eq 1.3.7} If $\epsilon$ is sufficiently small, then for all $w\in\mathcal{W}$ such that $w\neq 1$, we have an inequality
$$\sigma(a)\ll\sigma_w(a^{-1}ga)$$
for all $a\in A_{\mini}^{\overline{Q},+}(\delta)$ and all $g\in\mathcal{U}_w[<\epsilon \sigma(a)]$.
\end{num}

\vspace{2mm}

\noindent Claim \ref{eq 1.3.6} is a simple consequence of the definition of the log-norm $\sigma_1$ and the inclusion
$$\displaystyle \mathcal{U}_1[<\epsilon \sigma(a)]\big\backslash \left(\overline{Q}\left[<\epsilon \sigma(a)\right] aU_Q\left[<\epsilon \sigma(a)\right]a^{-1}\right)\subseteq \overline{Q}(F)\left(U_Q(F)\backslash aU_Q\left[<\epsilon \sigma(a)\right]a^{-1}\right)$$

\noindent We now prove \ref{eq 1.3.7}. Fix $w\in \mathcal{W}$ such that $w\neq 1$. For all $g=\overline{q}uw\in \mathcal{U}_w$ with $\overline{q}\in \overline{Q}$, $u\in U_Q$ and for all $a\in A_{\mini}$, we have

$$a^{-1}ga=a^{-1}\overline{q}w(a)\left(w(a)^{-1}uw(a)\right)w$$

\noindent Here $a^{-1}\overline{q}w(a)\in \overline{Q}$ and $w(a)^{-1}uw(a)\in U_Q$. Thus by \ref{eq 1.3.2} , we have

$$\sigma_w\left(a^{-1}ga\right)\geqslant \sigma\left(a^{-1}\overline{q}w(a)\right)$$

\noindent For all $c>0$, we have $\mathcal{U}_w[<c]=\overline{Q}[<c]U_Q[<c]w$. Consequently, to prove \ref{eq 1.3.7} it is sufficient to establish the following

\vspace{3mm}

\begin{num}
\item\label{eq 1.3.8} If $\epsilon$ is sufficiently small,we have
$$\sigma(a)\ll\sigma\left(a^{-1}\overline{q}w(a)\right)$$
for all $a\in A_{\mini}^{\overline{Q},+}(\delta)$ and all $\overline{q}\in\overline{Q}[<\epsilon \sigma(a)]$.
\end{num}

\vspace{3mm}

\noindent We have an inequality

$$\lvert H_{\overline{Q}}(\overline{q})\rvert \ll \sigma(\overline{q})$$

\noindent for all $\overline{q}\in \overline{Q}$. Hence,

\[\begin{aligned}
\displaystyle \left\lvert H_{\overline{Q}}\left(w(a)\right)-H_{\overline{Q}}(a)\right\rvert-\left\lvert H_{\overline{Q}}(\overline{q})\right\rvert & \leqslant \left\lvert H_{\overline{Q}}\left(w(a)\right)-H_{\overline{Q}}(a)+ H_{\overline{Q}}(\overline{q})\right\rvert \\
 & =\left\lvert H_{\overline{Q}}\left(a^{-1}\overline{q}w(a)\right)\right\rvert \\
 & \ll \sigma\left(a^{-1}\overline{q}w(a)\right)
\end{aligned}\]

\noindent for all $a\in A_{\mini}$ and all $\overline{q}\in \overline{Q}$. Therefore, to prove \ref{eq 1.3.8} it suffices to establish that

\begin{align}\label{eq 1.3.9}
\displaystyle \sigma(a)\ll \left\lvert H_{\overline{Q}}\left(w(a)\right)-H_{\overline{Q}}(a)\right\rvert,\; \mbox{ for all }a\in A_{\mini}^{\overline{Q},+}(\delta)
\end{align}

\noindent Set $\mathcal{A}_{\mini}=\mathcal{A}_{M_{\mini}}$ and $H_{\mini}=H_{M_{\mini}}$. Denote by $\Delta\subseteq R(A_{\mini},\overline{P}_{\mini})$ the subset of simple roots and by $\Delta^\vee\subseteq \mathcal{A}_{\mini}$, $\widehat{\Delta}\subseteq \mathcal{A}^*_{\mini}$ the corresponding sets of simple coroots and fundamental weights respectively. The coroot and fundamental weight associated to $\alpha\in \Delta$ will be denoted by $\alpha^\vee$ and $\varpi_\alpha$ respectively. Set $\Delta_{\overline{Q}}=\Delta\cap R(A_{\mini},U_{\overline{Q}})$. We have

$$\displaystyle \langle\varpi_\alpha,H_{\overline{Q}}(a)\rangle=\langle\varpi_\alpha,H_{\mini}(a)\rangle$$

\noindent for all $\alpha\in \Delta_{\overline{Q}}$ and all $a\in A_{\mini}$. Therefore, to get \ref{eq 1.3.9} it is sufficient to establish the following claim:

\begin{num}
\item\label{eq 1.3.10} There exists $\alpha\in \Delta_{\overline{Q}}$ such that
$$\displaystyle \sigma(a)\ll \left\langle\varpi_\alpha,H_{\mini}(a)-H_{\mini}(w(a))\right\rangle$$
for all $a\in A_{\mini}^{\overline{Q},+}(\delta)$.
\end{num}

\noindent Let $w=s_{\alpha_1}\ldots s_{\alpha_k}$ be a minimal decomposition of $w$ as a product of distinct simple reflections (thus with $\alpha_i\in \Delta$ for all $i$). Let $X\in \mathcal{A}_{\mini}$. We have the following identity which is easy to establish by induction

$$\displaystyle X-w(X)=\sum_{i=1}^k \langle\alpha_i,X\rangle s_{\alpha_1}\ldots s_{\alpha_{i-1}}(\alpha_i^\vee)$$

\noindent Hence for $\alpha\in \Delta$, we have

$$\displaystyle \left\langle\varpi_\alpha,X-w(X)\right\rangle=\sum_{i=1}^k \langle\alpha_i,X\rangle\langle\varpi_\alpha,s_{\alpha_1}\ldots s_{\alpha_{i-1}}(\alpha_i^\vee)\rangle$$

\noindent Notice that $\langle\varpi_\alpha,s_{\alpha_1}\ldots s_{\alpha_{i-1}}(\alpha_i^\vee)\rangle\geqslant 0$ for all $1\leqslant i\leqslant k$. Consequently, if $X\in\mathcal{A}_{\mini}^+$, meaning that $\langle\alpha,X\rangle\geqslant 0$ for all $\alpha\in \Delta$, we have

$$\left\langle\varpi_\alpha,X-w(X)\right\rangle\geqslant \langle\alpha_i,X\rangle\langle\varpi_\alpha,s_{\alpha_1}\ldots s_{\alpha_{i-1}}(\alpha_i^\vee)\rangle$$

\noindent for all $1\leqslant i\leqslant k$. Since $w\notin W^{L}$, there exists $i$ such that $\alpha_i\in \Delta_{\overline{Q}}$. Let $i$ be the minimal such index and set $\alpha=\alpha_i$. By the above, we have

$$\left\langle\varpi_\alpha,X-w(X)\right\rangle\geqslant \langle\alpha,X\rangle\langle\varpi_\alpha,s_{\alpha_1}\ldots s_{\alpha_{i-1}}(\alpha^\vee)\rangle$$

\noindent for all $X\in \mathcal{A}_{\mini}^+$. Since $\alpha_j\neq \alpha$ for $1\leqslant j\leqslant i-1$, we have

$$\langle\varpi_\alpha,s_{\alpha_1}\ldots s_{\alpha_{i-1}}(\alpha^\vee)\rangle=\langle\varpi_\alpha,\alpha^\vee\rangle=1$$

\noindent So finally

$$\left\langle\varpi_\alpha,X-w(X)\right\rangle\geqslant \langle\alpha,X\rangle$$

\noindent for all $X\in \mathcal{A}_{\mini}^+$. Now \ref{eq 1.3.10} follows immediately since by definition of $A_{\mini}^{\overline{Q},+}(\delta)$, we have

$$\langle\alpha, H_{\mini}(a)\rangle\geqslant \delta \sigma(a)$$

\noindent for all $a\in A_{\mini}^{\overline{Q},+}(\delta)$.

\item If $\sigma(a)\leqslant\epsilon^{-1}$ then the left hand side is empty and there is nothing to prove (as $\sigma(g)\geqslant 1$ for all $g\in G(F)$). So we shall only prove the inclusion for $\sigma(a)>\epsilon^{-1}$. There exists $\alpha>0$ such that
 
$$\lvert \log(u)\rvert_{\g} \leqslant e^{\alpha\sigma(u)}$$
 
\noindent for all $u\in U_Q(F)$. Also, there exists $\beta>0$ such that

$$\left\lvert aXa^{-1}\right\rvert_{\g}\leqslant \beta e^{-\delta \sigma(a)}\left\lvert X\right\rvert_{\g}$$

\noindent for all $X\in \mathfrak{u}_Q(F)$ and all $a\in A_{\mini}^{\overline{Q},+}(\delta)$. It follows that for a given $\epsilon>0$, we have

\[\begin{aligned}
\lvert \log\left(aua^{-1}\right)\rvert_{\g}=\lvert a\log(u)a^{-1}\rvert_{\g} & \leqslant \beta e^{-\delta\sigma(a)} \lvert log(u)\rvert_{\g} \\
 & \leqslant \beta e^{\left(\alpha \epsilon-\delta\right)\sigma(a)} \\
 & =\beta e^{\left(\alpha \epsilon-\delta+\delta'\right)\sigma(a)} e^{-\delta'\sigma(a)} \\
\end{aligned}\]

\noindent for all $a\in A_{\mini}^{\overline{Q},+}(\delta)$ and all $u\in U_Q\left[<\epsilon \sigma(a)\right]$. Now, it suffices to choose $\epsilon$ sufficiently small such that we have $\beta e^{\left(\alpha \epsilon-\delta+\delta'\right)\sigma(a)}\leqslant c_0$ for all $a\in A_{\mini}(F)$ such that $\sigma(a)>\epsilon^{-1}$ $\blacksquare$
\end{enumerate}

\subsection{Common spaces of functions}\label{section 1.4}

\noindent Let $X$ be a locally compact Hausdorff totally disconnected topological space and let $M$ be a real smooth manifold. In this paper the adjective smooth will have two meanings: a function from $X$ to a topological vector space $E$ is {\em smooth} if it is locally constant whereas a function from $M$ to a topological space $E$ is {\em smooth} if it is weakly $C^\infty$ in the sense of Appendix \ref{section A.3}. We shall denote by $\gls{CX}$ and $\gls{CM}$ the spaces of all smooth complex-valued functions on $X$ and $M$ respectively and by $\gls{CcX}$, $\gls{CcM}$ the subspaces of compactly-supported functions. We equip $C_c^\infty(M)$ and $C^\infty(M)$ with their usual locally convex topology. Then $C^\infty(M)$ is a Fr\'echet space whereas $C_c^\infty(M)$ is an LF space. We endow $C_c^\infty(X)$ with its finest locally convex topology. If $X$ admits a countable basis of open subsets, $C_c^\infty(X)$ is also an LF space (since it admits a countable basis). Restrictions to compact-open subsets induces an isomorphism

$$C^\infty(X)\simeq \varprojlim_{\mathcal{K}} C_c^\infty(\mathcal{K})$$

\noindent where $\mathcal{K}$ runs through the compact-open subsets of $X$. We shall endow $C^\infty(X)$ with the projective limit topology relative to this isomorphism. We also denote by $\gls{DiffM}$ the space of all smooth differential operators on $M$ which are globally of finite order and by $\gls{DiffkM}$, $k\in \mathbb{N}$, the subspace of smooth differential operators of order less than $k$. $\Diff^\infty_{\leqslant k}(M)$ carries a natural locally convex topology and if $M$ is countable at infinity, it is a Fr\'echet space. We endow $\Diff^\infty(M)$ with the direct limit topology relative to the natural isomorphism

$$\displaystyle \Diff^\infty(M)=\varinjlim_{k} \Diff^\infty_{\leqslant k}(M)$$

\noindent (so that if $M$ is countable at infinity, $\Diff^\infty(M)$ is an LF space). We denote by $\gls{D'M}$ and $\gls{D'X}$ the topological duals of $C_c^\infty(M)$ and $C_c^\infty(X)$ respectively and we call them the spaces of distributions on $M$ and $X$ respectively. If we have fixed a regular Borel measure $dm$ on $M$ (resp.\ $dx$ on $X$), then for every locally integrable function $F$ on $M$ (resp.\ on $X$), we will denote by $\gls{TF}$ the associated distribution on $M$ (resp.\ on $X$) i.e., we have

$$\displaystyle \langle T_F,f\rangle=\int_M F(m)f(m)dm,\;\;\; f\in C_c^\infty(M)$$
$$\displaystyle \left(\mbox{resp. }\langle T_F,f\rangle=\int_X F(x)f(x)dx,\;\;\; f\in C_c^\infty(X)\right)$$

\vspace{2mm}

\noindent Let $V$ be a finite dimensional $F$-vector space. Then, for all $f\in C^\infty(V)$ and all $\lambda\in F^\times$, we will denote by $\gls{flambda}$ the function defined by

$$f_\lambda(v)=f(\lambda^{-1}v),\;\;\; v\in V$$

\noindent We extend this action of $F^\times$ to the space of distributions $\mathcal{D}'(V)$ by setting $\langle \gls{Tlambda},f\rangle= \lvert \lambda\rvert^{\dim(V)}\langle T,f_{\lambda^{-1}}\rangle$ for all $T\in \mathcal{D}'(V)$, all $f\in C_c^\infty(V)$ and all $\lambda\in F^\times$. Moreover, we will say that a distribution $T\in \mathcal{D}'(V)$ is {\em homogeneous of degree $d$} if

$$T_\lambda=\lvert \lambda\rvert^{-d} T$$

\noindent for every $\lambda\in (F^{\times})^2$ (where $(F^\times)^2$ denotes the set of squares in $F^\times$).

\vspace{2mm}

\noindent We will also need the Schwartz spaces $\gls{S(g)}$ and $\gls{S(G)}$. If $F$ is $p$-adic, we have $\mathcal{S}(\mathfrak{g}(F))=C_c^\infty(\mathfrak{g}(F))$ and $\mathcal{S}(G(F))=C_c^\infty(G(F))$. Assume that $F=\mathbb{R}$. Then, $\mathcal{S}(\mathfrak{g}(F))$ is the space of all functions $f\in C^\infty(\mathfrak{g}(F))$ such that

$$\displaystyle \gls{qNu}(f)=\sup_{X\in\mathfrak{g}(F)} \lVert X\rVert_{\mathfrak{g}}^{N}\left\lvert \left(\partial(u)f\right)(X)\right\rvert<\infty$$

\noindent for all $N\geqslant 1$ and all $u\in S(\mathfrak{g})$. We endow $\mathcal{S}(\mathfrak{g}(F))$ with the topology defined by the semi-norms $q_{N,u}$, for all $N\geqslant 1$ and all $u\in S(\mathfrak{g})$. It is a Fr\'echet space. The natural inclusion $C_c^\infty(\mathfrak{g}(F))\subset \mathcal{S}(\mathfrak{g}(F))$ is continuous with dense image. We will say that a distribution $T$ on $\mathfrak{g}(F)$ is {\em tempered} if it extends to a continuous linear form on $\mathcal{S}(\mathfrak{g}(F))$. We denote by $\gls{S'(g)}$ the space of tempered distributions on $\mathfrak{g}(F)$.

\vspace{2mm}

\noindent Similarly, we define $\mathcal{S}(G(F))$ to be the space of all functions $f\in C^\infty(G(F))$ such that

$$\displaystyle \gls{qNuv}(f)=\sup_{x\in G(F)} \lVert x\rVert_G^{N}\left\lvert \left(L(u)R(v)f\right)(x)\right\rvert<\infty$$

\noindent for all $N\geqslant 1$ and all $u,v\in \mathcal{U}(\mathfrak{g})$. We endow $\mathcal{S}(G(F))$ with the topology defined by the semi-norms $q_{N,u,v}$, for all $N\geqslant 1$ and all $u,v\in \mathcal{U}(\mathfrak{g})$. It is also a Fr\'echet space.

\vspace{2mm}

\noindent Assume that a non-degenerate symmetric bilinear form $B$ and a measure have been fixed on $\mathfrak{g}(F)$ (this will be done in Section \ref{section 1.6}). Then we define the Fourier transform on $\mathcal{S}(\mathfrak{g}(F))$ by

$$\displaystyle \gls{fhat}(X)=\int_{\mathfrak{g}(F)} f(Y) \psi(B(X,Y)) dY, \;\;\; f\in \mathcal{S}(\mathfrak{g}(F)), X\in \mathfrak{g}(F)$$

\noindent and we extend this definition to tempered distributions by setting

$$\langle \gls{That},f\rangle=\langle T,\widehat{f}\rangle$$

\noindent for all $T\in \mathcal{S}'(\mathfrak{g}(F))$ and all $f\in \mathcal{S}(\mathfrak{g}(F))$. In the real case there exist two isomorphisms $S(\mathfrak{g})\simeq S(\mathfrak{g}^*)$, $u\mapsto \gls{pu}$, and $S(\mathfrak{g}^*)\simeq S(\mathfrak{g})$, $p\mapsto \gls{up}$, such that

$$\widehat{\partial(u)T}=p_u\widehat{T} \mbox{ and } \widehat{pT}=\partial(u_p)\widehat{T}$$

\noindent for all $T\in \mathcal{S}'(\mathfrak{g}(F))$, all $u\in S(\mathfrak{g})$ and all $p\in S(\mathfrak{g}^*)$. Note that if $B$ is chosen to be $G$-invariant, then these isomorphisms restrict to give isomorphisms $I(\mathfrak{g})\simeq I(\mathfrak{g}^*)$ and $I(\mathfrak{g}^*)\simeq I(\mathfrak{g})$.

\subsection{Harish-Chandra Schwartz space}\label{section 1.5}

\noindent We will denote by $\gls{XiG}$ the Harish-Chandra function. Let us recall its definition. Let $P_{\mini}$ be a minimal parabolic subgroup of $G$ and let $K$ be a maximal compact subgroup of $G(F)$ which is special in the $p$-adic case. Then we have $G(F)=P_{\mini}(F)K$ (Iwasawa decomposition). Consider the (smooth and normalized) induced representation (cf.\ Section \ref{section 2.3})

$$\displaystyle i^G_{P_{\mini}}\left(1\right)^\infty:=\{e\in C^\infty(G(F));\; e(pg)=\delta_{\mini}(p)^{1/2}e(g)\; \forall p\in P_{\mini}(F), g\in G(F)\}$$

\noindent that we equip with the scalar product

$$\displaystyle (e,e')=\int_K e(k)\overline{e'(k)} dk,\;\;\; e,e'\in i^G_{P_{\mini}}(1)^\infty$$

\noindent Let $e_K\in i^G_{P_{\mini}}\left(1\right)^\infty$ be the unique function such that $e_K(k)=1$ for all $k\in K$. Then the Harish-Chandra function $\Xi^G$ is defined by

$$\displaystyle \Xi^G(g)=\left(i_{P_{\mini}}^G(1,g)e_K,e_K\right),\; g\in G(F)$$

\noindent Of course, the function $\Xi^G$ depends on the various choices we made, but this doesn't matter because different choices would yield equivalent functions and the function $\Xi^G$ will only be used to give estimates. The next proposition summarizes the main properties of the function $\Xi^G$ that we will need. We indicate references for these after the statement.

\begin{prop}\label{proposition 1.5.1}
\begin{enumerate}[(i)]

\item Set

$$M_{\mini}^+=\{m\in M_{\mini}(F);\; \lvert \alpha(m)\rvert\leqslant 1 \;\forall \alpha\in R(A_{M_{\mini}},P_{\mini})\}$$

\noindent Then, there exists $d>0$ such that

$$\delta_{P_{\mini}}(m)^{1/2}\ll \Xi^G(m)\ll \delta_{P_{\mini}}(m)^{1/2}\sigma(m)^d$$

\noindent for all $m\in M_{\mini}^+$.

\item Let $m_{P_{\mini}}\colon G(F)\to M_{\mini}(F)$ be any map such that $g\in m_{P_{\mini}}(g)U_{\mini}(F)K$ for all $g\in G(F)$. Then, there exists $d>0$ such that

$$\Xi^G(g)\ll \delta_{P_{\mini}}(m_{P_{\mini}}(g))^{1/2}\sigma(g)^d$$

\noindent for all $g\in G(F)$.

\item Let $P=MU$ be a parabolic subgroup that contains $P_{\mini}$. Let $m_P\colon G(F)\to M(F)$ be any map such that $g\in m_{P}(g)U(F)K$ for all $g\in G(F)$. Then, we have

$$\displaystyle \Xi^G(g)=\int_K \delta_P(m_P(kg))^{1/2} \Xi^M(m_P(kg))dk$$

\noindent for all $g\in G(F)$.

\item Let $P=MU$ be a parabolic subgroup of $G$. Then, for all $d>0$, there exists $d'>0$ such that

$$\displaystyle \delta_P(m)^{1/2}\int_{U(F)} \Xi^G(mu)\sigma(mu)^{-d'}du\ll \Xi^M(m)\sigma(m)^{-d}$$

\noindent for all $m\in M(F)$.

\item There exists $d>0$ such that the integral

$$\displaystyle\int_{G(F)} \Xi^G(g)^2 \sigma(g)^{-d} dg$$

\noindent is convergent.

\item (Doubling principle) We have the equality

$$\int_K \Xi^G(g_1kg_2)dk=\Xi^G(g_1)\Xi^G(g_2)$$

\noindent for all $g_1,g_2\in G(F)$.
\end{enumerate}
\end{prop}

\vspace{2mm}

\noindent \ul{Proof}: Most of these are due to Harish-Chandra. A convenient reference is \cite{Wa2} in the $p$-adic case (see Lemme II.1.1 for (i), Lemme II.4.4 for (ii), Lemme II.1.6 for (iii), Proposition II.4.5 for (iv), Lemme II.1.5 for (v) and Lemme II.1.3 for (vi)) and \cite{Va} in the real case (see Theorem 30 p.339 for (i), Proposition 16(iv) p.329 for (iii), Theorem 23 p.360 for (iv), Proposition 31 p.340 for (v) and Proposition 16(iii) p.329 for (vi)) except concerning point (ii) of the Proposition for which we refer the reader to \cite{HC1} Lemma 85 and Corollary 1 p.108. $\blacksquare$

\vspace{2mm}

\noindent Using the function $\Xi^G$ we can define the {\em Harish-Chandra Schwartz space} $\gls{C(G(F))}$ as follows. For every function $f\in C(G(F))$ and all $d\in\mathbb{R}$, we set

$$\displaystyle \gls{pd}(f):=\sup_{g\in G(F)} \lvert f(g)\rvert \Xi^G(g)^{-1}\sigma(g)^d$$

\noindent If $F$ is $p$-adic then 

$$\displaystyle\mathcal{C}(G(F))=\bigcup_{K'}\mathcal{C}_{K'}(G(F))$$

\noindent where $K'$ runs through the open-compact subgroups of $G(F)$ and $\mathcal{C}_{K'}(G(F))$ is the space of functions $f\in C\left(K'\backslash G(F)/K'\right)$ such that $p_d(f)<\infty$ for all $d>0$. We endow the spaces $\mathcal{C}_{K'}(G(F))$ with the topology defined by the semi-norms $(p_d)_{d>0}$. These are Fr\'echet spaces and we equip $\mathcal{C}(G(F))$ with the direct limit topology. Thus $\mathcal{C}(G(F))$ is an LF space in this case.

\vspace{2mm}

\noindent If $F=\mathbb{R}$ then $\mathcal{C}(G(F))$ is by definition the space of all $f\in C^\infty(G(F))$ such that

$$\displaystyle \gls{puvd}(f):=p_d(R(u)L(v)f)<\infty$$

\noindent for all $d>0$ and all $u,v\in \mathcal{U}(\mathfrak{g})$. We equip $\mathcal{C}(G(F))$ with the topology defined by the semi-norms $p_{u,v,d}$, for all $u,v\in\mathcal{U}(\mathfrak{g})$ and all $d>0$. In this case, $\mathcal{C}(G(F))$ is a Fr\'echet space.

\vspace{2mm}

\begin{lem}\label{lemma 1.5.1bis}
Assume that $F=\mathbb{R}$. Let $f\in \mathcal{C}(G(F))$, $d>0$ and $B\subset \mathfrak{g}(F)$ be a compact. Then, we have
$$\displaystyle \left\lvert f(e^Xge^Y)-f(g)\right\rvert \ll \Xi^G(g) \sigma(g)^{-d} \left(\lvert X\rvert_{\mathfrak{g}}+\lvert Y\rvert_{\mathfrak{g}}\right)$$
for all $g\in G(F)$ and all $X,Y\in B$.
\end{lem}

\vspace{2mm}

\noindent\ul{Proof}: We have
\[\begin{aligned}
\displaystyle f(e^Xge^Y)-f(g) & =f(e^Xge^Y)-f(ge^Y)+f(ge^Y)-f(g) \\
 & =\int_0^1 (L(-X)f)(e^{tX}ge^Y)dt+\int_0^1 (R(Y)f)(ge^{tY}) dt
\end{aligned}\]
for all $g\in G(F)$ and all $X,Y\in B$. Hence,
\[\begin{aligned}
\displaystyle \lvert f(e^Yge^X)-f(g)\rvert & \leqslant \int_0^1 \lvert (L(-X)f)(e^{tX}ge^Y)\rvert dt+\int_0^1 \lvert (R(Y)f)(ge^{tY})\rvert dt \\
 & \leqslant p_d(L(-Y)f)\int_0^1 \Xi^G(e^{tX}ge^Y)\sigma(e^{tX}ge^Y)^{-d} dt+p_d(R(Y)f)\int_0^1 \Xi^G(ge^{tY})\sigma(ge^{tY})^{-d} dt \\
 & \ll \left[p_d(L(-Y)f)+ p_d(R(Y)f)\right] \Xi^G(g) \sigma(g)^{-d} \\
 & \leqslant \sup_{Z\in \mathfrak{g}(F); \lvert Z\rvert_{\mathfrak{g}}=1} \left[ p_d(L(Z)f)+ p_d(R(Z)f)\right] \left(\lvert X\rvert_{\mathfrak{g}}+\lvert Y\rvert_{\mathfrak{g}}\right) \Xi^G(g) \sigma(g)^{-d}
\end{aligned}\]
for all $g\in G(F)$ and all $X,Y\in B$. The lemma follows. $\blacksquare$

\vspace{2mm}

\noindent We now define what we will call the {\em weak Harish-Chandra Schwartz space} $\gls{Cw(G)}$. This topological space is important since it is the natural home for coefficients of tempered representations. Moreover, fixing a Haar measure on $G(F)$, this is precisely the smooth part of the space of tempered distributions. Again, the definition of $\mathcal{C}^w(G(F))$ differs in the $p$-adic and the real case.

\vspace{2mm}

\noindent If $F$ is $p$-adic, we have

$$\displaystyle \mathcal{C}^w(G(F))=\bigcup_{K'}\mathcal{C}_{K'}^w(G(F))$$

\noindent here again $K'$ runs through the open-compact subgroups of $G(F)$ and

$$\displaystyle \mathcal{C}^w_{K'}(G(F))=\bigcup_{d>0} \mathcal{C}^w_{K',d}(G(F))$$

\noindent where $\mathcal{C}^w_{K',d}(G(F))$ denotes the space of functions $f\in C\left(K'\backslash G(F)/K'\right)$ such that $p_{-d}(f)<\infty$. Equipped with the norm $p_{-d}$, $\mathcal{C}^w_{K',d}(G(F))$ is a Banach space. We endow $\mathcal{C}^w_{K'}(G(F))$ and $\mathcal{C}^w(G(F))$ with the direct limit topologies. These are LF spaces. We will also set

$$\displaystyle \gls{CdwG}=\bigcup_{K'} \mathcal{C}^w_{K',d}(G(F))$$

\noindent and we will equip this space with the direct limit topology. It is also an LF space.

\vspace{2mm}

\noindent If $F=\mathbb{R}$, we have

$$\displaystyle \mathcal{C}^w(G(F))=\bigcup_{d>0} \mathcal{C}^w_d(G(F))$$

\noindent where $\gls{CdwG}$ denotes the space of functions $f\in C^\infty(G(F))$ such that

$$\displaystyle p_{u,v,-d}(f):=p_{-d}\left(R(u)L(v)f\right)<\infty$$

\noindent for all $u,v\in\mathcal{U}(\mathfrak{g})$. We equip $\mathcal{C}^w_d(G(F))$ with the topology defined by the semi-norms $p_{u,v,-d}$, for all $u,v\in \mathcal{U}(\mathfrak{g})$. It is a Fr\'echet space. Finally, we endow $\mathcal{C}^w(G(F))$ with the direct limit topology so that it becomes an LF space.

\vspace{2mm}

\noindent In any case, the natural inclusion $\mathcal{C}(G(F))\subseteq \mathcal{C}^w(G(F))$ is continuous and we have the following

\begin{num}
\item\label{eq 1.5.1} $\mathcal{C}(G(F))$ is dense in $\mathcal{C}^w(G(F))$.
\end{num}

\noindent Indeed, we may even prove that $C_c^\infty(G(F))$ is dense in $\mathcal{C}^w(G(F))$. For all $t>0$, denote by $\kappa_t$ the characteristic function of $\{g\in G(F);\; \sigma(g)<t\}$. Let $\varphi\in C_c^\infty(G(F))$ be any positive function such that $\displaystyle \int_{G(F)}\varphi(g)dg=1$ and set $\varphi_t=\varphi\ast\kappa_t\ast\varphi$ for all $t>0$. Then we leave to the reader the task to prove that for all $f\in \mathcal{C}^w(G(F))$ we have

$$\displaystyle \lim\limits_{t\to \infty} \varphi_tf=f$$

\noindent in $\mathcal{C}^w(G(F))$. This proves the claim.

\vspace{2mm}

\noindent We end this section with a lemma that will be useful for us. The second part of this lemma gives a criterion for a function taking values in $\mathcal{C}^w(G(F))$ to be smooth.

\vspace{2mm}

\begin{lem}\label{lemma 1.5.1}
\begin{enumerate}[(i)]
\item Let $d>0$ and let $\nu$ be a continuous semi-norm on $\mathcal{C}_d^w(G(F))$. Then

\begin{enumerate}[(a)]

\item In the $p$-adic case, for all $\varphi_1,\varphi_2\in C_c^\infty(G(F))$, there exists a continuous semi-norm $\nu_{\varphi_1,\varphi_2}$ on $\mathcal{C}_d^w(G(F))$ such that

$$\displaystyle \nu\left(R(\varphi_1)L(\varphi_2)R(g_1)L(g_2)f\right)\leqslant \nu_{\varphi_1,\varphi_2}(f) \Xi^G(g_1)\Xi^G(g_2)\sigma(g_1)^d\sigma(g_2)^d$$

\noindent for all $f\in \mathcal{C}^w_d(G(F))$ and all $g_1,g_2\in G(F)$.

\item In the real case, there exists $k\geqslant 0$ (which depends on $\nu$) such that for all $\varphi_1,\varphi_2\in C_c^k(G(F))$, there exists a continuous semi-norm $\nu_{\varphi_1,\varphi_2}$ on $\mathcal{C}_d^w(G(F))$ such that

$$\displaystyle \nu\left(R(\varphi_1)L(\varphi_2)R(g_1)L(g_2)f\right)\leqslant \nu_{\varphi_1,\varphi_2}(f) \Xi^G(g_1)\Xi^G(g_2)\sigma(g_1)^d\sigma(g_2)^d$$

\noindent for all $f\in \mathcal{C}^w_d(G(F))$ and all $g_1,g_2\in G(F)$.

\end{enumerate}
\item Let $V$ be a real vector space and let $\varphi\colon V\times G(F)\to \mathbb{C}$ be a function such that

\begin{enumerate}[(a)]
\item In the $p$-adic case: for all $g\in G(F)$ the function $\lambda\in V\mapsto \varphi(\lambda,g)$ is smooth and there exists a compact-open subgroup $K'$ of $G(F)$ such that for all $\lambda\in V$ the function $\varphi(\lambda,.)$ is $K'$-biinvariant.

\item In the real case: for all $\lambda\in V$, the function $g\in G(F)\mapsto \varphi(\lambda,g)$ is smooth and for all $u,v\in \mathcal{U}(\mathfrak{g})$ and all $g\in G(F)$ the function $\lambda\in V \mapsto \left(R(u)L(v)\varphi\right)(\lambda,g)$ is smooth.

\item In the $p$-adic case: for every differential operator with constant coefficients $D\in S(V)$, there exist two constants $C,d>0$ such that

$$\left\lvert \left(D\varphi\right)(\lambda,g)\right\rvert\leqslant C\Xi^G(g)\sigma(g)^d$$

\noindent for all $g\in G(F)$ and all $\lambda\in V$.

\item In the real case: for every differential operator with constant coefficients $D\in S(V)$, there exists $d>0$ such that for all $u,v\in \mathcal{U}(\mathfrak{g})$ there is a positive continuous function $C_{u,v}(.)$ on $V$ such that

$$\left\lvert \left(DR(u)L(v)\varphi\right)(\lambda,g)\right\rvert\leqslant C_{u,v}(\lambda)\Xi^G(g)\sigma(g)^d$$

\noindent for all $g\in G(F)$ and all $\lambda\in V$.
\end{enumerate}

\noindent Then, the map $\lambda\mapsto \varphi(\lambda,.)$ takes value in $\mathcal{C}^w(G(F))$ and defines a smooth function from $V$ to $\mathcal{C}^w(G(F))$.
\end{enumerate}
\end{lem}

\vspace{2mm}

\noindent\ul{Proof}: 

\begin{enumerate}[(i)]

\item We will only prove (b), the proof of (a) being similar and easier. We may assume without loss of generality that $\nu=p_{u,v,-d}$ for some $u,v\in \mathcal{U}(\mathfrak{g})$. Set $k=\max\left(\deg(u),\deg(v)\right)$. Then, we have

$$\displaystyle \nu\left(R(\varphi_1)L(\varphi_2)R(g_1)L(g_2)f\right)=p_{-d}\left(R\left(u\varphi_1\right)L\left(v\varphi_2\right)R(g_1)L(g_2)f\right)$$

\noindent for all $\varphi_1,\varphi_2\in C_c^k(G(F))$, all $g_1,g_2\in G(F)$ and all $f\in \mathcal{C}^w_d(G(F))$, where $u\varphi_1$ and $v\varphi_2$ stand for $L(u)\varphi_1$ and $L(v)\varphi_2$ respectively. Hence, we may assume that $\nu=p_{-d}$. Let $\varphi_1,\varphi_2\in C_c(G(F))$. Then, we have

$$\displaystyle \left(R(\varphi_1)L(\varphi_2)R(g_1)L(g_2)f\right)(g)=\int_{G(F)\times G(F)} \varphi_1(\gamma_1)\varphi_2(\gamma_2)f\left(g_2^{-1}\gamma_2^{-1}g\gamma_1g_1\right)d\gamma_1 d\gamma_2$$

\noindent for all $f\in \mathcal{C}^w_d(G(F))$ and all $g,g_1,g_2\in G(F)$. Since $\sigma(xy)\ll \sigma(x)\sigma(y)$ for all $x,y\in G(F)$, it follows that

\[\begin{aligned}
\displaystyle & \left\lvert \left(R(\varphi_1)L(\varphi_2)R(g_1)L(g_2)f\right)(g)\right\rvert\ll \\
 & p_{-d}(f)\sigma(g)^d\sigma(g_1)^d\sigma(g_2)^d\int_{G(F)\times G(F)} \left\lvert\varphi_1(\gamma_1)\right\rvert \left\lvert\varphi_2(\gamma_2)\right\rvert \Xi^G\left(g_2^{-1}\gamma_2^{-1}g\gamma_1g_1\right) \sigma(\gamma_2)^d\sigma(\gamma_1)^d d\gamma_1 d\gamma_2
\end{aligned}\]

\noindent for all $f\in \mathcal{C}^w_d(G(F))$ and all $g,g_1,g_2\in G(F)$. Moreover, by the doubling principle (Proposition \ref{proposition 1.5.1}(vi)), we have

$$\displaystyle \int_{G(F)\times G(F)} \left\lvert\varphi_1(\gamma_1)\right\rvert \left\lvert\varphi_2(\gamma_2)\right\rvert \Xi^G\left(g_2^{-1}\gamma_2^{-1}g\gamma_1g_1\right) \sigma(\gamma_2)^d\sigma(\gamma_1)^d d\gamma_1 d\gamma_2\ll \Xi^G(g)\Xi^G(g_1)\Xi^G(g_2)$$

\noindent for all $g,g_1,g_2\in G(F)$. So finally, we get

$$\displaystyle p_{-d}\left(R(\varphi_1)L(\varphi_2)R(g_1)L(g_2)f\right)\ll p_{-d}(f)\Xi^G(g_1)\Xi^G(g_2)\sigma(g_1)^d\sigma(g_2)^d$$

\noindent for all $f\in \mathcal{C}^w_d(G(F))$ and all $g_1,g_2\in G(F)$ and this ends the proof of (i).

\item Assume first that $F$ is $p$-adic. Let $K'$ be as in (a). Then the condition (c) implies that for all $k\geqslant 0$ there exists $d>0$ such that $\lambda\mapsto \varphi(\lambda,.)$ defines a strongly $C^k$ map from $V$ to $\mathcal{C}^w_{d,K'}(G(F))$ and the result follows.

\vspace{2mm}

\noindent Assume now that $F=\mathbb{R}$. Then, by the condition (d), for all $u\in \mathcal{U}(\mathfrak{g})$ and for all $D\in S(V)$ the function

$$(\lambda,g)\in V\times G(F)\mapsto \left(DR(u)\right)\varphi(\lambda,g)$$

\noindent is locally bounded. It follows that $\varphi$ is smooth (as a function on $V\times G(F)$). In particular, for all $u,v\in \mathcal{U}(\mathfrak{g})$ and all $D\in S(V)$, we have

$$DR(u)L(v)\varphi=R(u)L(v)D\varphi$$

\noindent Let $k\geqslant 0$ be an integer. It now follows from (d) that  there exists $d>0$ such that $\left(D\varphi\right)(\lambda,.)\in \mathcal{C}^w_d(G(F))$ for all $\lambda\in V$ and all $D\in S(V)$ of degree less than $k$. From this we easily deduce, using (d) again, that the map $\lambda\mapsto \varphi(\lambda,.)$ defines a strongly $C^k$ map from $V$ to $\mathcal{C}^w_d(G(F))$. The result follows. $\blacksquare$
\end{enumerate}

\subsection{Measures}\label{section 1.6}

\noindent We fix once and for all a (unitary) continuous non-trivial additive character $\gls{psi}\colon F\to \mathbb{S}^1$ and we equip $F$ with the autodual Haar measure with respect to $\psi$. We also fix a Haar measure $d^\times t$ on $F^\times$ to be $\lvert t\rvert^{-1}dt$ where $dt$ is the Haar measure on $F$ that we just fixed.

\vspace{2mm}

\noindent Fix a $G(F)$-invariant nondegenerate bilinear form $B$ on $\mathfrak{g}(F)$. If $F=\mathbb{R}$, we choose $B$ so that for every maximal compact subgroup $K$ of $G(F)$ the restriction of $B$ to $\mathfrak{k}(F)$ is negative definite and the restriction to $\mathfrak{k}(F)^\perp$ (the orthogonal of $\mathfrak{k}(F)$ with respect to $B$) is positive definite. We endow $\mathfrak{g}(F)$ with the autodual measure with respect to $B$, it is the only Haar measure $dX$ on $\mathfrak{g}(F)$ such that the Fourier transform

$$\displaystyle \widehat{f}(Y)=\int_{\mathfrak{g}(F)}f(X)\psi\left(B(X,Y)\right)dX,\;\;\; f\in \mathcal{S}(\mathfrak{g}(F))$$

\noindent satisfies $\widehat{\widehat{f\, }}(X)=f(-X)$. We equip $G(F)$ with the unique Haar measure such that the exponential map has a Jacobian equal to $1$ at the origin. Similarly, for every $F$-algebraic subgroup $H$ of $G$ such that the restriction of $B(.,.)$ to $\mathfrak{h}(F)$ is non-degenerate, we equip $\mathfrak{h}(F)$ with the autodual measure with respect to $B$ and we lift this measure to $H(F)$ by means of the exponential map. This fixes for example the Haar measures on the Levi subgroups of $G$ as well as on the maximal subtori of $G$. For other subgoups of $G(F)$, for example unipotent radicals of parabolic subgroups of $G$, we fix an arbitrary Haar measure on the Lie algebra and we lift it to the group, again using the exponential map.

\vspace{2mm}

\noindent For every Levi subgroup $M$ of $G$, we equip $\mathcal{A}_M$ and $i\mathcal{A}_M^*$ with Haar measures as follows. In the real case we choose any measures whereas the $p$-adic case, we choose the unique Haar measures such that $meas\left(\mathcal{A}_M/\widetilde{\mathcal{A}}_{M,F}\right)=1$ and $meas\left(i\mathcal{A}^*_M/i\widetilde{\mathcal{A}}^\vee_{M,F}\right)=1$. 

\vspace{2mm}

\noindent Let $T$ be a maximal subtorus of $G$. Besides the Haar measure $dt$ that has been fixed above on $T(F)$, we will need another Haar measure that we shall denote by $d_ct$. First, we define a Haar measure $d_ca$ on $A_T(F)$ as follows. If $F$ is $p$-adic, it is the unique Haar measure such that the maximal compact subgroup of $A_T(F)$ is of measure $1$. In the real case, $d_ca$ is the unique Haar measure such that the surjective homomorphism $H_T\colon A_T(F)\to \mathcal{A}_T$ is locally measure preserving (note that $\mathcal{A}_T$ coincide with $\mathcal{A}_M$ for a certain Levi subgroup $M$ so that a Haar measure has already been fixed on $\mathcal{A}_T$). Finally, in both cases $d_ct$ is the unique Haar measure on $T(F)$ such that the quotient measure $d_ct/d_ca$ gives $T(F)/A_T(F)$ the measure $1$. To avoid confusions, we shall only use the Haar measure $dt$ but we need to introduce the only factor $\gls{nuT}>0$ such that $d_ct=\nu(T)dt$.

\vspace{2mm}

\noindent Denote by $\gls{Nilg}$ the set of nilpotent orbits in $\mathfrak{g}(F)$. Let $\mathcal{O}\in \Nil(\mathfrak{g})$. Then, for all $X\in \mathcal{O}$ the bilinear map $(Y,Z)\mapsto B(Y,[X,Z])$ descends to a non-degenerate symplectic form on $\mathfrak{g}(F)/\mathfrak{g}_X(F)$ that is the tangent space of $\mathcal{O}$ at $X$. This defines on $\mathcal{O}$ a structure of symplectic $F$-analytic manifold. Using the Haar measure on $F$, this equips $\mathcal{O}$ with a natural ``autodual" measure. This measure is obviously $G(F)$-invariant.

\vspace{2mm}

\noindent The following considerations will be useful for Chapter \ref{section 10} only. Let $V$ be an $\overline{F}$-subspace of $\mathfrak{g}=\mathfrak{g}(\overline{F})$. Even if $V$ is not defined over $F$ we can talk of {\em Haar measures} on $V$: these are elements of $(\bigwedge^{max}V)^*\smallsetminus\{0\}$ modulo multiplication by an element of norm $1$ in $\overline{F}$ (Indeed, in the case $V$ is defined over $F$ the autodual additive measure on $F$ allows to interpret such a class as a Haar measure on $V(F)$). Assume that a Haar measure $\mu_V$ has been fixed on $V$ (for example one of the measures that we fixed above). There is a natural notion of {\em dual Haar measure} $\gls{muV*}$ on $V^*$: it is the unique Haar measure on $V^*$ such that the image of $\mu_V\otimes \mu_V^*$ by the natural pairing $(\bigwedge^{max}V)^*\otimes (\bigwedge^{max} V^*)^*\to \overline{F}$ is of norm $1$. Let $V^\perp$ be the orthogonal of $V$ with respect to $B$. Then, we may associate to $\mu_V$ a Haar measure $\mu_V^\perp$ on $V^\perp$ as follows. Using the form $B$ we have a natural isomorphism $(\bigwedge^{max}\mathfrak{g})^*\simeq (\bigwedge^{max}V^*)^*\otimes (\bigwedge^{max}V^\perp)^*$. Then $\gls{muVperp}$ is the unique Haar measure on $V^\perp$ such that via this isomorphism we have $\mu_{\mathfrak{g}}=\mu_V^*\otimes \mu_V^\perp$ (modulo a scalar of norm $1$) where $\mu_{\mathfrak{g}}$ denotes the autodual Haar measure on $\mathfrak{g}$ that we fixed above. If $V$ is defined over $F$, we have the formula

\begin{align}\label{eq 1.6.1}
\displaystyle \int_{V(F)} f(v)d\mu_V(v)=\int_{V^\perp(F)} \widehat{f}(v^\perp) d\mu^\perp_V(v^\perp)
\end{align}

\noindent for all $f\in\mathcal{S}(\mathfrak{g}(F))$. We easily check that

\begin{align}\label{eq 1.6.2}
(\mu_V^\perp)^\perp=\mu_V
\end{align}

\noindent Also, if we have a decomposition $\mathfrak{g}=V_1\oplus V_2$ and two Haar measures $\mu_{V_1}$, $\mu_{V_2}$ on $V_1$ and $V_2$ such that

$$\displaystyle \mu_{\mathfrak{g}}=\mu_{V_1}\otimes \mu_{V_2}$$

\noindent then we also have the equality

\begin{align}\label{eq 1.6.3}
\displaystyle \mu_{\mathfrak{g}}=\mu_{V_1}^\perp\otimes \mu_{V_2}^\perp
\end{align}

\noindent relative to the decomposition $\mathfrak{g}=V_1^\perp\oplus V_2^\perp$.

\vspace{2mm}

\noindent Finally, suppose that $V$ and $W$ are $\overline{F}$-subspaces of $\mathfrak{g}$ and that $T\colon V\simeq W$ is a linear isomorphism. Then $T$ induces an isomorphism $(\bigwedge^{max}T)^*\colon (\bigwedge^{max}V)^*\simeq (\bigwedge^{max}W)^*$ and if $\mu_V$ is a Haar measure on $V$ then we will denote by $T_*\mu_V$ the image of $\mu_V$ by this isomorphism (a measure on $W$). Notice that if $V=W$ then $T_*\mu_V=\lvert \det(T)\rvert \mu_V$.

\subsection{Spaces of conjugacy classes and invariant topology}\label{section 1.8}

\noindent If $H$ is a connected linear algebraic group defined over $F$, we will denote by $\gls{GammaH}$ the set of semi-simple conjugacy classes in $H(F)$. Thus, we have a natural projection

$$H_{\ssi}(F)\twoheadrightarrow \Gamma(H)$$

\noindent and we endow $\Gamma(H)$ with the quotient topology. Then, $\Gamma(H)$ is Hausdorff and locally compact. Moreover for every connected linear algebraic group $H'$ over $F$ and every embedding $H'\hookrightarrow H$ the induced map $\Gamma(H')\to \Gamma(H)$ is continuous and proper. We define similarly the space $\gls{Gammah}$ of semi-simple conjugacy classes in $\mathfrak{h}(F)$. This space satisfies similar properties.

\vspace{2mm}

\noindent We will say of a subset $A\subseteq \mathfrak{h}(F)$ (resp.\ $A\subseteq H(F)$) that it is {\em completely $H(F)$-invariant} if it is $H(F)$-invariant and if moreover for all $X\in A$ (resp.\ $g\in A$) its semi-simple part $X_s$ (resp.\ $g_s$) also belongs to $A$. Closed invariant subsets are automatically completely $H(F)$-invariant. We easily check that the completely $H(F)$-invariant open subsets of $\mathfrak{h}(F)$ (resp.\ of $H(F)$) define a topology. We will call it the {\em invariant topology}. This topology coincides with the pull-back of the topology on $\Gamma(\mathfrak{h})$ (resp.\ on $\Gamma(H)$) just defined by the natural map

$$\mathfrak{h}(F)\to \Gamma(\mathfrak{h})\;\; \left(\mbox{resp.\ } H(F)\to \Gamma(H)\right)$$

\noindent which associates to $X\in \mathfrak{h}(F)$ (resp.\ $g\in H(F)$) the conjugacy class of the semi-simple part of $X$ (resp.\ of $g$). In particular, we have the following property which will be used many times implicitly in that paper: If $\omega\subseteq \mathfrak{h}(F)$ (resp.\ $\Omega\subseteq H(F)$) is a completely $H(F)$-invariant open subset and $\omega'\subseteq \omega$ (resp.\ $\Omega'\subseteq \Omega$) is invariant open and contains $\omega_{\ssi}$ (resp.\ $\Omega_{\ssi}$) then $\omega'=\omega$ (resp.\ $\Omega'=\Omega$). We will say of an invariant subset $L\subseteq \mathfrak{h}(F)$ (resp.\ $L\subseteq H(F)$) that it is {\em compact modulo conjugation} if it is closed and if there exists a compact subset $\mathcal{K}\subseteq \mathfrak{h}(F)$ (resp.\ $\mathcal{K}\subseteq H(F)$) such that $L=\mathcal{K}^G$, it is equivalent to ask that $L$ is completely $H(F)$-invariant and that for every maximal torus $T\subset H$ the intersection $L\cap\mathfrak{t}(F)$ (resp.\ $L\cap T(F)$) is compact, it is also equivalent to the fact that $L$ is compact for the invariant topology.

\vspace{2mm}

\noindent We will denote by $\gls{GammaellG}$ and $\gls{GammaregG}$ (resp.\ $\gls{Gammaellg}$ and $\gls{Gammaregg}$) the subsets of elliptic and regular conjugacy classes in $\Gamma(G)$ (resp.\ in $\Gamma(\mathfrak{g})$) respectively. The subset $\Gamma_{\elli}(G)$ (resp.\ $\Gamma_{\elli}(\mathfrak{g})$) is closed in $\Gamma(G)$ (resp.\ in $\Gamma(\mathfrak{g})$) whereas $\Gamma_{\reg}(G)$ (resp.\ $\Gamma_{\reg}(\mathfrak{g})$) is an open subset of $\Gamma(G)$ (resp.\ of $\Gamma(\mathfrak{g})$). Let $\mathcal{T}(G)$ be a set of representatives for the conjugacy classes of maximal tori in $G$. We equip $\Gamma(G)$ and $\Gamma(\mathfrak{g})$ with the unique regular Borel measures $dx$ and $dX$ such that

$$\displaystyle \int_{\Gamma(G)} \varphi_1(x)dx=\sum_{T\in \mathcal{T}(G)} \lvert W(G,T)\rvert^{-1}\int_{T(F)} \varphi_1(t) dt$$

$$\displaystyle \int_{\Gamma(\mathfrak{g})} \varphi_2(X)dX=\sum_{T\in \mathcal{T}(G)} \lvert W(G,T)\rvert^{-1}\int_{\mathfrak{t}(F)} \varphi_2(X) dX$$

\noindent for all $\varphi_1\in C_c(\Gamma(G))$ and all $\varphi_2\in C_c(\Gamma(\mathfrak{g}))$. We have the Weyl integration formula

$$\displaystyle \int_{G(F)} f(g) dg=\int_{\Gamma(G)} D^G(x)^{1/2} J_G(x,f)dx \;\; \left(\mbox{resp. } \int_{\mathfrak{g}(F)} f(X)dX=\int_{\Gamma(\mathfrak{g})} D^G(X)^{1/2} J_G(X,f)dX\right)$$

\noindent for all $f\in \mathcal{S}(G(F))$ (resp.\ for all $f\in \mathcal{S}(\mathfrak{g}(F))$). We deduce from this and the local boundedness of normalized orbital integrals the following fact

\vspace{3mm}

\begin{num}
\item\label{eq 1.8.1} The function $x\mapsto D^G(x)^{-1/2}$ (resp.\ $X\mapsto D^G(X)^{-1/2}$) is locally integrable on $G(F)$ (resp.\ on $\mathfrak{g}(F)$).
\end{num}

\vspace{3mm}

\noindent We define an abstract norm $\gls{normGammag}$ on $\Gamma(\mathfrak{g})$ as follows. Fix a set of tori $\mathcal{T}(G)$ as above. Then, we define $\lVert.\rVert_{\Gamma(\mathfrak{g})}$ by

$$\displaystyle \lVert X\rVert_{\Gamma(\mathfrak{g})}=\inf_{X'} \; 1+\lvert X'\rvert ,\;\;\; X\in \Gamma(\mathfrak{g})$$

\noindent where the infimum is taken over the set of $X'\in \bigsqcup_{T\in \mathcal{T}(G)} \mathfrak{t}(F)$ that belong to the conjugacy class of $X$. We will need the two following estimates

\vspace{3mm}

\begin{num}
\item\label{eq 1.8.2} For all $k\geqslant 0$ and for all $N>0$ sufficiently large, the integral
$$\displaystyle\int_{\Gamma(\mathfrak{g})} \log\left(2+D^G(X)^{-1}\right)^k \lVert X\rVert_{\Gamma(\mathfrak{g})}^{-N}dX$$
is absolutely convergent.

\item\label{eq 1.8.3} Assume that $F=\mathbb{R}$. Then, for all $N>0$, there exists a continuous semi-norm $\nu_N$ on $\mathcal{S}(\mathfrak{g}(F))$ such that
$$\displaystyle \left\lvert J_G(X,f)\right\rvert\leqslant \nu_N(f) \lVert X\rVert_{\Gamma(\mathfrak{g})}^{-N}$$
for all $f\in \mathcal{S}(\mathfrak{g}(F))$.
\end{num}

\vspace{3mm}

\noindent We say that an element $x\in G(F)$ is {\em anisotropic} if it is regular semi-simple and $G_x(F)$ is compact. We will denote by $\gls{GFani}$ the subset of anisotropic elements and by $\gls{GammaaniG}$ the set of anisotropic conjugacy classes in $G(F)$. We equip $\Gamma_{ani}(G)$ with the quotient topology relative to the natural projection $G(F)_{ani}\twoheadrightarrow \Gamma_{ani}(G)$. Let $\mathcal{T}_{ani}(G)$ be a set of representatives for the $G(F)$-conjugacy classes of maximal anisotropic tori of $G$ (a torus $T$ is {\em anisotropic} if $T(F)$ is compact). We equip $\Gamma_{ani}(G)$ with the quotient topology and we endow it with the unique regular Borel measure such that

$$\displaystyle \int_{\Gamma_{ani}(G)}\varphi(x)dx=\sum_{T\in\mathcal{T}_{ani}(G)}\lvert W(G,T)\rvert^{-1}\nu(T)\int_{T(F)} \varphi(t)dt$$

\noindent for all $\varphi\in C_c(\Gamma_{ani}(G))$, where the factor $\nu(T)$ has been defined in Section \ref{section 1.6}. Note that if $A_G\neq 1$ then $\Gamma_{ani}(G)=\emptyset$.

\subsection{Orbital integrals and their Fourier transforms}\label{section 1.7}

\noindent For $x\in G_{\reg}(F)$ (resp.\ $X\in \mathfrak{g}_{\reg}(F)$), we define the {\em normalized orbital integral} at $x$ (resp.\ at $X$) by

$$\displaystyle \gls{JGxf}=D^G(x)^{1/2}\int_{G_x(F)\backslash G(F)} f(g^{-1}xg)dg,\;\;\; f\in \mathcal{C}(G(F))$$

\begin{center}
(resp.\ $\displaystyle \gls{JGXf}=D^G(X)^{1/2}\int_{G_X(F)\backslash G(F)} f(g^{-1}Xg)dg,\;\;\; f\in \mathcal{S}(\mathfrak{g}(F))$)
\end{center}

\noindent the integral being absolutely convergent for all $f\in \mathcal{C}(G(F))$ (resp.\ for all $f\in \mathcal{S}(\mathfrak{g}(F))$). This defines a tempered distribution $J_G(x,.)$ (resp.\ $J_G(X,.)$) on $G(F)$ (resp.\ on $\mathfrak{g}(F)$). For all $f\in \mathcal{C}(G(F))$ (resp.\ $f\in \mathcal{S}(\mathfrak{g}(F))$), the function $x\in G_{\reg}(F)\mapsto J_G(x,f)$ (resp.\ $X\in \mathfrak{g}_{\reg}(F)\mapsto J_G(X,f)$) is locally bounded on $G(F)$ (resp.\ on $\mathfrak{g}(F)$).

\vspace{2mm}

\noindent Similarly, for $\mathcal{O}\in \Nil(\mathfrak{g})$, we define the orbital integral on $\mathcal{O}$ by

$$\displaystyle \gls{JOf}=\int_{\mathcal{O}}f(X)dX,\;\;\; f\in \mathcal{S}(\mathfrak{g}(F))$$

\noindent We have

\begin{align}\label{eq 1.7.1}
\displaystyle J_{\mathcal{O}}(f_{\lambda})=\lvert \lambda\rvert^{\dim(\mathcal{O})/2} J_{\mathcal{O}}(f)
\end{align}

\noindent for all $\mathcal{O}\in \Nil(\mathfrak{g})$ and all $\lambda\in F^{\times 2}$ (recall that $f_\lambda(X)=f(\lambda^{-1}X)$). Denote by $\gls{Nilregg}$ the subset of regular nilpotent orbits in $\mathfrak{g}(F)$. This set is empty unless $G$ is quasi-split in which case we have $\dim(\mathcal{O})=\delta(G)$ for all $\mathcal{O}\in \Nil_{\reg}(\mathfrak{g})$. By the above equality, the distributions $J_{\mathcal{O}}$ for $\mathcal{O}\in \Nil_{\reg}(\mathfrak{g})$ are all homogeneous of degree $\delta(G)/2-\dim(\mathfrak{g})$. This characterizes the distributions $J_{\mathcal{O}}$, $\mathcal{O}\in \Nil_{\reg}(\mathfrak{g})$, among the invariant distributions supported in the nilpotent cone. More precisely, we have

\vspace{3mm}

\begin{num}
\item\label{eq 1.7.2} The invariant distributions on $\mathfrak{g}(F)$ supported in the nilpotent cone and homogeneous of degree $\delta(G)/2-\dim(\mathfrak{g})$ are precisely linear combinations of the distributions $J_{\mathcal{O}}$ for $\mathcal{O}\in \Nil_{\reg}(\mathfrak{g})$.
\end{num}

\vspace{3mm}

\noindent This follows from Lemma 3.3 of \cite{HCDS} in the $p$-adic case and from Corollary 3.9 of \cite{BV} in the real case (there is a sign error in this last reference, $n-\alpha$ should be replaced by $\alpha-n$ and the inequality $\alpha\geqslant \frac{n-r}{2}$ should be replaced by $\alpha\leqslant \frac{n-r}{2}$).

\vspace{2mm}

\noindent According to Harish-Chandra, there exists a unique smooth function $\gls{jhat}$ on $\mathfrak{g}_{\reg}(F)\times \mathfrak{g}_{\reg}(F)$ which is locally integrable on $\mathfrak{g}(F)\times \mathfrak{g}(F)$ such that

$$\displaystyle J_G(X,\widehat{f})=\int_{\mathfrak{g}(F)} \widehat{j}(X,Y) f(Y)dY$$

\noindent for all $X\in \mathfrak{g}_{\reg}(F)$ and all $f\in \mathcal{S}(\mathfrak{g}(F))$. We have the following control on the size of $\widehat{j}$:

\vspace{3mm}

\begin{num}
\item\label{eq 1.7.3} The function $(X,Y)\in \mathfrak{g}_{\reg}(F)\times \mathfrak{g}_{\reg}(F)\mapsto D^G(Y)^{1/2}\widehat{j}(X,Y)$ is globally bounded.
\end{num}

\vspace{3mm}

\noindent cf.\ Theorem 7.7 and Lemma 7.9 of \cite{HCDS} in the $p$-adic case and Proposition 9 p.112 of \cite{Va} in the real case. We will need the following property regarding to the non-vanishing of the function $\widehat{j}$

\vspace{3mm}

\begin{num}
\item\label{eq 1.7.4} Assume that $G$ admits elliptic maximal tori. Then, for all $Y\in \mathfrak{g}_{\reg}(F)$ there exists $X\in \mathfrak{g}_{\reg}(F)_{\elli}$ such that $\widehat{j}(X,Y)\neq 0$.
\end{num}

\vspace{3mm}

\noindent In the $p$-adic case, this follows from Theorem 9.1 and Lemma 9.6 of \cite{HCDS} whereas in the real case, it is a consequence of Theorem 4 p.104 and Theorem 11 p.126 of \cite{Va}.

\vspace{2mm}

\noindent Similarly, for any nilpotent orbit $\mathcal{O}\in \Nil(\mathfrak{g})$, there exists a smooth function $\gls{jhatO}$ on $\mathfrak{g}_{\reg}(F)$ which is locally integrable on $\mathfrak{g}(F)$ such that

$$\displaystyle J_{\mathcal{O}}(\widehat{f})=\int_{\mathfrak{g}(F)} \widehat{j}(\mathcal{O},X)f(X)dX$$

\noindent for all $f\in \mathcal{S}(\mathfrak{g}(F))$.  We know that the function $(D^G)^{1/2}\widehat{j}(\mathcal{O},.)$ is locally bounded on $\mathfrak{g}(F)$ (Theorem 6.1 of \cite{HCDS} in the $p$-adic case and Theorem 17 p.63 of \cite{Va} in the real case). By \ref{eq 1.7.1}, the functions $\widehat{j}(\mathcal{O},.)$ satisfy the following homogeneity property

\begin{align}\label{eq 1.7.5}
\displaystyle \widehat{j}(\mathcal{O},\lambda X)=\lvert \lambda\rvert^{-\dim(\mathcal{O})/2} \widehat{j}(\lambda\mathcal{O},X)
\end{align}

\noindent for all $\mathcal{O}\in \Nil(\mathfrak{g})$, all $X\in \mathfrak{g}_{\reg}(F)$ and all $\lambda\in F^\times$. Recall also that for every nilpotent orbit $\mathcal{O}\in \Nil(\mathfrak{g})$, we have $\lambda \mathcal{O}=\mathcal{O}$ for all $\lambda\in F^{\times 2}$.

\subsection{$(G,M)$-families}\label{section 1.9}

\noindent We collect here some useful facts from Arthur's theory of $(G,M)$-families as developed for example in \cite{A3} Section 17.

\noindent Let $M$ be a Levi subgroup of $G$ and $V$ a locally convex topological vector space. A {\em $(G,M)$-family} with values in $V$ is a family $(c_P)_{P\in \mathcal{P}(M)}$ of smooth functions on $i\mathcal{A}_M^*$ taking values in $V$ such that for all adjacent parabolic subgroups $P,P'\in \mathcal{P}(M)$, the functions $c_P$ and $c_{P'}$ coincide on the hyperplane supporting the wall that separates the positive chambers for $P$ and $P'$. Arthur associates to any $(G,M)$-family $(c_P)_{P\in \mathcal{P}(M)}$ (taking values in $V$) an element $c_M$ of $V$ as follows. The function

$$\displaystyle c_M(\lambda)=\sum_{P\in \mathcal{P}(M)} c_P(\lambda) \theta_P(\lambda)^{-1}$$

\noindent where

$$\displaystyle \theta_P(\lambda)=meas\left(\mathcal{A}^G_M/\mathbb{Z}\Delta_P^\vee\right)^{-1} \prod_{\alpha\in \Delta_P} \lambda(\alpha^\vee),\;\;\; P\in \mathcal{P}(M)$$

\noindent extends to a smooth function on $i\mathcal{A}_M^*$ and we have $c_M=c_M(0)$. Here, $\Delta_P$ denotes the set of simple roots of $A_M$ in $P$, $\Delta_P^\vee$ denotes the corresponding set of simple coroots and for every $\alpha\in \Delta_P$ we have denoted by $\alpha^\vee$ the corresponding simple coroot. For all $P\in \mathcal{P}(M)$, Arthur also constructs an element $c_P'\in V$ from the $(G,M)$-family $(c_P)_{P\in\mathcal{P}(M)}$. This element is the value at $\lambda=0$ of the function

$$\displaystyle c'_P(\lambda)=\sum_{P\subset Q} (-1)^{a_P-a_Q} c_P(\lambda_Q) \widehat{\theta}_P^Q(\lambda)^{-1}\theta_Q(\lambda_Q)^{-1}$$

\noindent where the sum is over the parabolic subgroups $Q=L_QU_Q$ containing $P$, $\theta_Q$ is defined as above, $\lambda_Q$ denotes the projection of $\lambda$ onto $i\mathcal{A}_{L_Q}^*$, and

$$\displaystyle \widehat{\theta}_P^Q(\lambda)=meas\left(\mathcal{A}_M^{L_Q}/\mathbb{Z}\left(\widehat{\Delta}^Q_P\right)^\vee\right)^{-1}\prod_{\alpha\in \Delta_P^Q} \lambda(\varpi_\alpha^\vee)$$

\noindent where this time $\Delta_P^Q$ denotes the set of simple roots of $A_M$ in $P\cap L_Q$, $\widehat{\Delta}_P^Q$ denotes the corresponding set of simple coweights and for every $\alpha\in \Delta_P^Q$ we have denoted by $\varpi^\vee_\alpha$ the corresponding simple coweight.

\vspace{2mm}

\noindent Let $L\in \mathcal{L}(M)$ and $Q=L_QU_Q\in \mathcal{F}(L)$. Starting from a $(G,M)$-family $(c_P)_{P\in \mathcal{P}(M)}$, we can construct a $(L_Q,L)$-family $(c_R^Q)_{R\in \mathcal{P}^{L_Q}(L)}$ as follows: for all $R\in \mathcal{P}^{L_Q}(L)$ and all $\lambda\in i\mathcal{A}_L^*$, we set $c_R^Q(\lambda)=c_P(\lambda)$ where $P$ is any parabolic subgroup in $\mathcal{P}(M)$ such that $P\subset Q(R)=RU_Q$. Applying the previous formal procedure to this new $(L_Q,L)$-family, we obtain an element $c^Q_L\in V$. We will usually simply set $c_L=c^G_L$. Notice that using the $(G,L)$-families $(c_Q)_{Q\in \mathcal{P}(L)}$, $L\in \mathcal{L}(M)$, we may define as above elements $c'_Q\in V$ for all $Q\in \mathcal{F}(M)$.

\vspace{2mm}

\noindent Assume now that $V$ is equipped with a continuous multiplication $V\times V\to V$ making it into a $\mathbb{C}$-algebra. Starting from two $(G,M)$-families $(c_P)_{P\in \mathcal{P}(M)}$ and $(d_P)_{P\in \mathcal{P}(M)}$ we may form they product $((cd)_P)_{P\in\mathcal{P}(M)}$, given by $(cd)_P=c_Pd_P$, which is again a $(G,M)$-family. We have the following splitting formulas (cf.\ Lemma 17.4 and Lemma 17.6 of \cite{A3})

\begin{align}\label{eq 1.9.1}
\displaystyle (cd)_M=\sum_{Q\in\mathcal{F}(M)} c^Q_Md'_Q
\end{align} 

\noindent and

\begin{align}\label{eq 1.9.2}
\displaystyle (cd)_M=\sum_{L_1,L_2\in\mathcal{L}(M)} d_M^G(L_1,L_2) c^{Q_1}_Md^{Q_2}_M 
\end{align}

\noindent where in the second formula $Q_1\in \mathcal{P}(L_1)$, $Q_2\in \mathcal{P}(L_2)$ are parabolic subgroups that depend on the choice of a point $X\in\mathcal{A}_M$ in general position and $d_M^G(L_1,L_2)$ is a non-negative real number that is nonzero if and only if $\mathcal{A}_{L_1}^G\oplus\mathcal{A}_{L_2}^G=\mathcal{A}_M^G$. Moreover, we have $d^G_M(G,M)=d^G_M(M,G)=1$. Starting from only one $(G,M)$-family, we also have the following descent formula (cf.\ Lemma 17.5 of \cite{A3})

\begin{align}\label{eq 1.9.3}
\displaystyle c_L=\sum_{L'\in\mathcal{L}(M)} d_M^G(L,L') c^{Q'}_M
\end{align}

\vspace{2mm}

\noindent A {\em $(G,M)$-orthogonal set} is a family $(Y_P)_{P\in \mathcal{P}(M)}$ of points in $\mathcal{A}_M$ such that for all adjacent parabolic subgroups $P,P'\in \mathcal{P}(M)$ there exists a real number $r_{P,P'}$ such that $Y_P-Y_{P'}=r_{P,P'}\alpha^\vee$, where $\alpha$ is the unique root of $A_M$ that is positive for $P$ and negative for $P'$. If moreover we have $r_{P,P'}\geqslant 0$ for all adjacent $P,P'\in\mathcal{P}(M)$, then we say that the family is {\em positive}. Obviously, if $(Y_P)_{P\in\mathcal{P}(M)}$ is a $(G,M)$-orthogonal set, then the family $(c_P)_{P\in\mathcal{P}(M)}$ defined by $c_P(\lambda)=e^{\lambda(Y_P)}$ is a $(G,M)$-family. If the family $(Y_P)_{P\in\mathcal{P}(M)}$ is positive then there is an easy interpretation for the number $c_M$: it is the volume in $\mathcal{A}_M^G$ of the convex hull of the set $\{Y_P,\; P\in \mathcal{P}(M)\}$.

\subsection{Weighted orbital integrals}\label{section 1.10}

\noindent Let $M$ be a Levi subgroup of $G$. Choose a maximal compact subgroup $K$ of $G(F)$ that is special in the $p$-adic case. Recall that using $K$, we may construct for every $P\in \mathcal{P}(M)$ a map

$$H_P\colon G(F)\to \mathcal{A}_M$$

\noindent (cf.\ Section \ref{section 1.1}). For every $g\in G(F)$, the family $(H_P(g))_{P\in \mathcal{P}(M)}$ is a positive $(G,M)$-orthogonal set. Hence, it defines a $(G,M)$-family $(v_P(g,.))_{P\in \mathcal{P}(M)}$ and the number $\gls{vMg}$ associated to this $(G,M)$-family is just the volume in $\mathcal{A}_M^G$ of the convex hull of the $H_P(g)$, $P\in \mathcal{P}(M)$. The function $g\mapsto v_M(g)$ is obviously invariant on the left by $M(F)$ and on the right by $K$.

\vspace{2mm}

\noindent Let $x\in M(F)\cap G_{\reg}(F)$. Then, for $f\in \mathcal{C}(G(F))$, we define the {\em weighted orbital integral} of $f$ at $x$ to be

$$\displaystyle \gls{JMxf}=D^G(x)^{1/2}\int_{G_x(F)\backslash G(F)} f(g^{-1}xg)v_M(g)dg$$

\noindent (note that the above expression is well-defined since $G_x\subset M$). The integral above is absolutely convergent and this defines a tempered distribution $J_M(x,.)$ on $G(F)$. More generally, we have seen in the last section how to associate to the $(G,M)$-family $(v_P(g,.))_{P\in \mathcal{P}(M)}$ complex numbers $\gls{vLQg}$ for all $L\in \mathcal{L}(M)$ and all $Q\in \mathcal{F}(L)$. This allows us to define tempered distributions $J_L^G(x,.)$ on $G(F)$ for all $L\in \mathcal{L}(M)$ and all $Q\in \mathcal{F}(L)$ by setting

$$\displaystyle \gls{JLQxf}=D^G(x)^{1/2}\int_{G_x(F)\backslash G(F)} f(g^{-1}xg)v_L^Q(g)dg,\;\;\; f\in \mathcal{C}(G(F))$$

\noindent The functions $x\in M(F)\cap G_{\reg}(F)\mapsto J_L^Q(x,f)$ are easily seen to be $M(F)$-invariant.

\vspace{2mm}

\noindent Let $X\in \mathfrak{m}(F)\cap \mathfrak{g}_{\reg}(F)$. We define similarly {\em weighted orbital integrals} $J_L^Q(X,.)$, $L\in \mathcal{L}(M)$, $Q\in \mathcal{F}(L)$. These are tempered distributions on $\mathfrak{g}(F)$ given by

$$\displaystyle \gls{JLQXf}=D^G(X)^{1/2}\int_{G_X(F)\backslash G(F)} f(g^{-1}Xg)v_L^Q(g)dg,\;\;\; f\in \mathcal{S}(\mathfrak{g}(F))$$

\noindent When $Q=G$, we will simply set $J_L^G(X,f)=\gls{JLXf}$. For all $L\in \mathcal{L}(M)$ and all $Q\in \mathcal{F}(L)$, we have an inequality $v_L^Q(g)\ll \sigma_{M\backslash G}(g)$ for all $g\in G(F)$. Using \ref{eq 1.2.2} and \ref{eq 1.2.4}, we easily deduce the following

\vspace{3mm}

\begin{num}
\item\label{eq 1.10.1} Assume that $F=\mathbb{R}$. Then, there exists $k\geqslant 0$ such that for all $N\geqslant 0$ there exists a continuous semi-norm $\nu_N$ on $\mathcal{S}(\mathfrak{g}(F))$ such that
$$\displaystyle \left\lvert J_L^Q(X,f)\right\rvert\leqslant \nu_N(f) \log\left(2+D^G(X)^{-1}\right)^k \lVert X\rVert_{\Gamma(\mathfrak{g})}^{-N}$$
for all $f\in \mathcal{S}(\mathfrak{g}(F))$.
\end{num}

\vspace{3mm}

\noindent We will need the following lemma regarding to the behavior of weighted orbital integrals under the action of invariant differential operators  (cf.\ Proposition 11.1 and Lemma 12.4 of \cite{A6}).

\begin{lem}\label{lemma 1.10.1}
Assume that $F=\mathbb{R}$ and let $T\subset M$ be a maximal torus. Then, we have
\begin{enumerate}[(i)]
\item For all $f\in \mathcal{C}(G(F))$ the function $x\in T_{\reg}(F)\mapsto J_M(x,f)$ is smooth and for all $z\in \mathcal{Z}(\mathfrak{g})$, there exist smooth differential operators $\partial^L_M(.,z_L)$ on $T_{\reg}(F)$ for all $L\in \mathcal{L}(M)\backslash \{M\}$ such that

$$\displaystyle J_M(x,zf)-z_TJ_M(x,f)=\sum_{\substack{L\in \mathcal{L}(M) \\ L\neq M}} \partial^L_M(x,z_L)J_L(x,f)$$

\noindent for all $f\in \mathcal{C}(G(F))$ and all $z\in T_{\reg}(F)$.

\item For all $f\in \mathcal{S}(\mathfrak{g}(F))$ the function $X\in \mathfrak{t}_{\reg}(F)\mapsto J_M(X,f)$ is smooth and for all $u\in I(\mathfrak{g})$, there exist smooth differential operators $\partial^L_M(.,u_L)$ on $\mathfrak{t}_{\reg}(F)$ for all $L\in \mathcal{L}(M)\backslash \{M\}$ such that

$$\displaystyle J_M(X,\partial(u)f)-\partial(u_T)J_M(X,f)=\sum_{\substack{L\in \mathcal{L}(M) \\ L\neq M}} \partial^L_M(X,u_L)J_L(X,f)$$

\noindent for all $f\in \mathcal{S}(\mathfrak{g}(F))$ and all $X\in \mathfrak{t}_{\reg}(F)$.
\end{enumerate}
\end{lem}

\section{Representations}\label{section 2}

This chapter contains some background on representations of $G(F)$ that will be used extensively in the rest of the paper. Here is a more precise description of the content of each section. In Section \ref{section 2.1}, we collect some basic facts on smooth representations. In Section \ref{section 2.2}, we recall the fundamental notion of {\em tempered representations} as well as some important properties of those. Sections \ref{section 2.3} and \ref{section 2.4} concern parabolic induction of smooth representations and (normalized) intertwining operators on them. These are used in Section \ref{section 2.5} to define, following Arthur, {\em weighted characters} which are distributions on the group $G(F)$ generalizing the usual characters and are spectral counterparts to the weighted orbital integrals defined in Section \ref{section 1.10}. In Section \ref{section 2.6}, we recall two fundamental results of harmonic analysis on $G(F)$ which are the {\em matricial Paley-Wiener} theorem and the {\em Harish-Chandra Plancherel formula}. They together give a full spectral decomposition of the Harish-Chandra Schwartz space $\mathcal{C}(G(F))$ and an inversion formula allowing to recover a function from its {\em Fourier transform}. Both are due to Harish-Chandra (a convenient reference being \cite{Wa2} in the $p$-adic case) except for the matricial Paley-Wiener in the Archimedean case which was proved by Arthur \cite{A8}. Finally, in Section \ref{section 2.7} we collect some facts on the so-called {\em elliptic representations} in the sense of Arthur \cite{A4}.
 
\subsection{Smooth representations, Elliptic regularity}\label{section 2.1}

\noindent Recall that a {\em continuous representation} of $G(F)$ is a pair $(\pi,V_\pi)$ where $V_\pi$ is a locally convex topological vector space and $\pi\colon G(F)\to Gl(V_\pi)$ is a morphism such that the resulting action

$$G(F)\times V_\pi\to V_\pi$$
$$(g,v)\mapsto \pi(g)v$$

\noindent is continuous. If $V_\pi$ is complete or even quasi-complete, we get an action of $\left(C_c(G(F)),\ast\right)$ on $V_\pi$ given by

$$\displaystyle \pi(f)v=\int_{G(F)} f(g)\pi(g)vdg,\;\; f\in C_c(G(F)), v\in V_\pi$$

\noindent A vector $v\in V_\pi$ is said to be {\em smooth} if the orbit map

$$\gamma_v\colon  g\in G(F)\mapsto \pi(g)v\in V_\pi$$

\noindent is smooth (i.e., it is locally constant in the $p$-adic case and weakly infinitely differentiable in the real case). We will denote by $\gls{Vpiinf}$ the subspace of smooth vectors. This subspace is $G(F)$-invariant and, if $F=\mathbb{R}$, it carries a natural action of $\mathcal{U}(\mathfrak{g})$. These two actions will be denoted by $\gls{piinf}$ or even by $\pi$ is there is no risk of confusion. For all $f\in C_c^\infty(G(F))$, the image of $\pi(f)$ is included in $V_\pi^\infty$. A continuous representation $(\pi,V_\pi)$ is said to be {\em smooth} if $V_\pi=V_\pi^\infty$.

\vspace{2mm}

\noindent Let $(\pi,V_\pi)$ be a smooth representation of $G(F)$. In the $p$-adic case we always have

$$\displaystyle \pi(C_c^\infty(G(F)))V_\pi=V_\pi$$

\noindent In the real case it is not always true. By a theorem of Dixmier-Malliavin \cite{DM}, it is at least true when $V_\pi$ is a Fr\'echet space. For example $\mathcal{C}(G(F))$ is a smooth Fr\'echet representation of $G(F)$ for the action given by left translation. Hence, we have a factorization

\begin{align}\label{eq 2.1.1}
\displaystyle \mathcal{C}(G(F))=C_c^\infty(G(F))\ast \mathcal{C}(G(F))
\end{align}

\vspace{2mm}

\noindent where $\ast$ denotes the convolution operator. Assume that $F=\mathbb{R}$ and let $H$ be an algebraic subgroup of $G$. Fix a basis $X_1,\ldots,X_h$ of $\mathfrak{h}(F)$ and set

$$\displaystyle \Delta_H=1-X_1^2-\ldots-X_h^2\in \mathcal{U}(\mathfrak{h})$$

\noindent The differential operator $R(\Delta_H)$ on $H(F)$ is elliptic. Hence, by elliptic regularity (cf.\ \cite{BK} Lemma 3.7), for every integer $m$ such that $2m>\dim(H)$, there exists a function $\varphi_1\in C_c^{2m-\dim(H)-1}(H(F))$ and a function $\varphi_2\in C_c^\infty(H(F))$ such that

\begin{align}\label{eq 2.1.2}
\displaystyle \varphi_1\ast \Delta_H^m+\varphi_2=\delta_1
\end{align}

\noindent where $\delta_1$ denotes the Dirac distribution at the identity and $\Delta_H^m$ is viewed as a distribution supported at the origin. It follows in particular that for every smooth representation $(\pi,V_\pi)$ of $G(F)$, we have

$$\displaystyle \pi(\varphi_1)\pi(\Delta_H^m)+\pi(\varphi_2)=Id_{V_\pi}$$

\subsection{Unitary and tempered representations}\label{section 2.2}

\noindent Recall that a {\em unitary representation} of $G(F)$ is a continuous representation $(\pi,\mathcal{H}_\pi)$ of $G(F)$ on a Hilbert space $\mathcal{H}_\pi$ such that for all $g\in G(F)$ the operator $\pi(g)$ is unitary. A unitary representation $(\pi,\mathcal{H}_\pi)$ is {\em irreducible} if $\mathcal{H}_\pi$ is nonzero and has no nontrivial closed $G(F)$-invariant subspace. We will only consider unitary representations that are of finite length. Such representations are finite direct sums of irreducible unitary representations. To avoid multiple repetitions of the words ``finite length" we will henceforth say ``unitary representation" to mean ``unitary representation of finite length". There is an action of $i\mathcal{A}_G^*$ on unitary representations given by $(\lambda,\pi)\mapsto \gls{pilambda}$ where $\pi_\lambda$ acts on the same space as $\pi$ and $\pi_\lambda(g)=e^{\lambda(H_G(g))}\pi(g)$, for all $g\in G(F)$. We will denote by $\gls{iAGpivee}$ the stabilizer of $\pi$ for this action. Notice that we always have $i\mathcal{A}^\vee_{G,F}\subset i\mathcal{A}_{G,\pi}^\vee\subset i\widetilde{\mathcal{A}}^\vee_{G,F}$. For $(\pi,\mathcal{H}_\pi)$ an unitary representation, we will denote by $(\gls{pibar},\gls{Hpibar})$ the {\em complex-conjugate representation} which identifies naturally (using the scalar product on $\mathcal{H}_\pi$) to the dual representation. 

\vspace{2mm}

\noindent Let us fix a compact maximal subgroup $K$ of $G(F)$. We will denote by $\gls{Khat}$ the set of equivalence classes of irreducible representations of $K$. For $\rho\in \widehat{K}$, we will denote by $\gls{drho}$ its dimension. For $(\pi,\mathcal{H}_\pi)$ a unitary representation of $G(F)$ and $\rho\in \widehat{K}$, we will denote by $\mathcal{H}_\pi(\rho)$ the $\rho$-isotypic component of $\mathcal{H}_\pi$. Every irreducible unitary representation $(\pi,\mathcal{H}_\pi)$ of $G(F)$ is {\em admissible} in the sense that

$$\dim \mathcal{H}_\pi(\rho)<\infty$$

\noindent for all $\rho\in \widehat{K}$. In the real case, we even have

\begin{align}\label{eq 2.2.1}
\dim \mathcal{H}_\pi(\rho)\leqslant d(\rho)^2
\end{align}

\noindent for all $\rho\in \widehat{K}$. Still in the real case, let us choose a basis $X_1,\ldots,X_n$ of $\mathfrak{k}(\mathbb{R})$ such that $B(X_i,X_j)=-\delta_{i,j}$ for $i,j=1,\ldots,n$ (recall that we choose the bilinear form $B$ such that $B_{\mid \mathfrak{k}}$ is negative definite) and set $\gls{DeltaK}=1-X_1^2-\ldots-X_n^2\in \mathcal{U}(\mathfrak{k})$. Then, $\Delta_K$ is in the center $\mathcal{Z}(\mathfrak{k})$ of $\mathcal{U}(\mathfrak{k})$ and doesn't depend on the basis chosen. It follows that for all $\rho\in \widehat{K}$, $\Delta_K$ acts by a scalar $\gls{crho}$ on the space of $\rho$. We always have $c(\rho)\geqslant 1$ and there exists $k\geqslant 1$ such that the sum

$$\displaystyle \sum_{\rho\in \widehat{K}}c(\rho)^{-k}$$

\noindent converges absolutely. Moreover, there exists $\ell\geqslant 1$ such that $d(\rho)\leqslant c(\rho)^\ell$ for all $\rho\in \widehat{K}$. Hence the sum

$$\displaystyle \sum_{\rho\in \widehat{K}}d(\rho)^2c(\rho)^{-k}$$

\noindent is convergent for $k$ sufficiently large. By \ref{eq 2.2.1}, it follows that for $k\geqslant 1$ sufficiently large there exists $C_k>0$ such that

\begin{align}\label{eq 2.2.2}
\displaystyle \sum_{\rho\in \widehat{K}}c(\rho)^{-k}\dim \mathcal{H}_\pi(\rho)<C_k
\end{align}

\noindent for every unitary irreducible representation $(\pi,\mathcal{H}_\pi)$ of $G(F)$.

\vspace{2mm}

\noindent Let $(\pi,\mathcal{H}_\pi)$ be a unitary representation. We endow the subspace of smooth vectors $\mathcal{H}^\infty_\pi$ with its own locally convex topology which is defined as follows. In the $p$-adic case, $\mathcal{H}_\pi^\infty$ is equipped with its finest locally convex topology. If $F=\mathbb{R}$, we endow $\mathcal{H}_\pi^\infty$ with the topology defined by the semi-norms

$$\displaystyle \lVert e\rVert_u=\lVert \pi^\infty(u)e\rVert,\;\;\; e\in \mathcal{H}_\pi^\infty$$

\noindent for all $u\in\mathcal{U}(\mathfrak{g})$, where $\lVert.\rVert$ is the norm derived from the scalar product on $\mathcal{H}_\pi$. In this case, $\mathcal{H}_\pi^\infty$ is a Fr\'echet space. A vector $e\in \mathcal{H}_\pi$ is smooth if and only if it is smooth for the $K$-action. Moreover, if $F=\mathbb{R}$, the semi-norms $\gls{normu}$, $u\in \mathcal{U}(\mathfrak{k})$, already generate the topology on $\mathcal{H}_\pi^\infty$. More precisely, the topology on $\mathcal{H}_\pi^\infty$ is generated by the family of semi-norms $(\lVert.\rVert_{\Delta_K^n})_{n\geqslant 0}$ and we have $\lVert .\rVert_{\Delta_K^m}\leqslant\lVert.\rVert_{\Delta_K^n}$ for $m\leqslant n$.

\vspace{2mm}

\noindent We will denote by $\gls{Hpi-inf}$ the topological dual of $\mathcal{H}_\pi^\infty$ which following our convention of Appendix \ref{section A} is equipped with the strong topology (in the $p$-adic case this is just the algebraic dual of $\mathcal{H}_\pi$ with the weak topology on it) and by $\gls{pi-inf}$ the natural representation of $G(F)$ on that space. It is a continuous representation. The scalar product on $\mathcal{H}_\pi$ gives a natural embedding $\overline{\mathcal{H}_\pi}\subset \mathcal{H}_\pi^{-\infty}$ and we have

$$\displaystyle \pi^{-\infty}\left(C_c^\infty(G(F))\right)\mathcal{H}_\pi^{-\infty}\subseteq \overline{\mathcal{H}_\pi^\infty}$$

\noindent (we even have an equality by Dixmier-Malliavin). We will always use the slight abuse of notation of denoting by $\pi$, $\overline{\pi}$, $\pi^\infty$ and $\pi^{-\infty}$ both the representations and the spaces on which these representations act. Also, we will always denote by $\gls{scprod}$ the scalar product on a given unitary representation (linear in the first variable) and by $\lVert . \rVert$ the induced norm.

\vspace{2mm}

\noindent Let again $\pi$ be a unitary representation of $G(F)$. Then we will denote by $\gls{Endpi}$ the space of continuous endomorphisms of the space of $\pi$. It is naturally a Banach space for the operator-norm

$$\displaystyle \gls{normT}=\sup_{\lVert e\rVert=1} \lVert Te\rVert,\;\;\; T\in \End(\pi)$$

\noindent Moreover $\End(\pi)$ is a continuous representation of $G(F)\times G(F)$ for the action given by left and right translations. We will denote by $\gls{Endpiinf}$ the subspace of smooth vectors and we will equip it with its own locally convex topology as follows. In the $p$-adic case we endow this space with its finest locally convex topology whereas in the real case we equip it with the topology defined by the semi-norms

$$\displaystyle \gls{normTuv}=\vertiii{\pi(u)T\pi(v)}, \;\;\; u,v\in \mathcal{U}(\mathfrak{g}), T\in \End(\pi)^\infty$$

\noindent in which case it is a Fr\'echet space. Once again, the topology on $\End(\pi)^\infty$ is generated by the semi-norms $(\vertiii{.}_{\Delta_K^n,\Delta_K^n})_{n\geqslant 1}$ and we have $\vertiii{.}_{\Delta_K^m,\Delta_K^m}\leqslant \vertiii{.}_{\Delta_K^n,\Delta_K^n}$ for $m\leqslant n$. Every $T\in \End(\pi)^\infty$ is traceable and for all $f\in C_c^\infty(G(F))$, the operator $\pi(f)$ belongs to $\End(\pi)^\infty$. Moreover, by Harish-Chandra, there exists a smooth function $\gls{thetapi}$ on $G_{\reg}(F)$ which is locally integrable on $G(F)$ such that

$$\displaystyle \Tr(\pi(f))=\int_{G(F)} \theta_\pi(g) f(g)dg$$

\noindent for all $f\in C_c^\infty(G(F))$. We call $\theta_\pi$ the character of $\pi$.

\vspace{2mm}

\noindent We have a natural embedding $\pi^\infty\otimes \overline{\pi^\infty}\subset \End(\pi)^\infty$ which sends $e\otimes e'$ to the operator $\gls{Tee'}$ given by

$$e_0\in \pi\mapsto (e_0,e')e$$

\noindent In the $p$-adic case, we even have an equality $\End(\pi)^\infty=\pi^\infty\otimes \overline{\pi^\infty}$  whereas in the real case $\pi^\infty\otimes\overline{\pi^\infty}$ is only a dense subspace of $\End(\pi)^\infty$ and we have $\End(\pi)^\infty=\pi^\infty\widehat{\otimes}_p\overline{\pi^\infty}$ where $\widehat{\otimes}_p$ denotes the projective topological tensor product (cf.\ Appendix \ref{section A.5}). It is easy to infer from this description that any $T\in \End(\pi)^\infty$ extends to a continuous linear map $T\colon \overline{\pi^{-\infty}}\to \pi^\infty$, the extension being necessarily unique since $\pi^\infty$ is dense in $\overline{\pi^{-\infty}}$. This induces a natural linear map

$$\displaystyle \End(\pi)^\infty\to \Hom(\overline{\pi^{-\infty}},\pi^\infty)$$

\noindent which is continuous, where we equip $\Hom(\overline{\pi^{-\infty}},\pi^\infty)$ with the strong topology.

\vspace{2mm}

\noindent We will say that a unitary representation $\pi$ is {\em tempered} if for all $e,e'\in\pi^\infty$ we have an inequality

\begin{align}\label{eq 2.2.3}
\displaystyle \left\lvert (\pi(g)e,e')\right\rvert\ll \Xi^G(g)
\end{align}

\noindent for all $g\in G(F)$. This inequality extends to $\End(\pi)^\infty$ in the sense that for all $T\in \End(\pi)^\infty$, we have an inequality 

\begin{align}\label{eq 2.2.4}
\displaystyle \left\lvert \Tr(\pi(g)T) \right\rvert\ll \Xi^G(g)
\end{align}

\noindent for all $g\in G(F)$. In the $p$-adic case, this is well-known and follows from \cite{CHH} Theorem 2 as the function $g\mapsto \Tr(\pi(g)T)$ is a finite sum of coefficients of $\pi$. If $F=\mathbb{R}$ and $\pi$ is moreover irreducible (which we can assume without loss of generality), we actually have a more precise inequality which follows from \cite{Sun}. Indeed, from {\it loc. cit.} and \ref{eq 2.2.2} we easily infer that there exists $n\geqslant 0$ and $C>0$ such that for every irreducible tempered representation $\pi$ and for all $T\in \End(\pi)^\infty$, we have

\begin{align}\label{eq 2.2.5}
\displaystyle \left\lvert \Tr(\pi(g)T)\right\rvert \leqslant C \Xi^G(g)\vertiii{T}_{\Delta_K^n,\Delta_K^n}
\end{align}

\noindent for all $g\in G(F)$. In particular, we have

\begin{align}\label{eq 2.2.6}
\displaystyle \lvert (\pi(g)e,e')\rvert\leqslant C \Xi^G(g)\lVert e\rVert_{\Delta_K^n} \lVert e'\rVert_{\Delta_K^n}
\end{align}

\noindent for all $e,e'\in \pi^\infty$ and all $g\in G(F)$ (still assuming that $\pi$ is an irreducible tempered representation).

\noindent Twists by unitary characters preserve tempered representations. We will denote by $\gls{TempG}$ the set of isomorphism classes of irreducible tempered representations. If $\pi$ is a tempered representation, then we may extend the action of $C_c^\infty(G(F))$ on $\pi^\infty$ to an action of $\mathcal{C}(G(F))$ by setting

$$\displaystyle (\pi(f)e,e')=\int_{G(F)} f(g) (\pi(g)e,e') dg$$

\noindent for all $e,e'\in \mathcal{C}(G(F))$ and all $f\in \mathcal{C}(G(F))$. Note that the vector $\pi(f)v$ a priori belongs to $\overline{\pi^{-\infty}}$ (in the real case this follows from \ref{eq 2.2.6}), the fact that it actually belongs to $\pi^\infty$ follows from the factorization \ref{eq 2.1.1}. This factorization also implies that we have $\pi(f)\in \End(\pi)^\infty$ for all $f\in \mathcal{C}(G(F))$. In the real case, it is easy to infer from \ref{eq 2.2.5} that there exists a continuous semi-norm $\nu$ on $\mathcal{C}(G(F))$ such that

\begin{align}\label{eq 2.2.7}
\displaystyle \lVert \pi(f)e\rVert\leqslant \nu(f) \lVert e\rVert
\end{align}

\noindent for all tempered representations $\pi$ of $G(F)$, all $e\in \pi$ and all $f\in \mathcal{C}(G(F))$.

\vspace{2mm}

\noindent Let $\pi$ be an irreducible unitary representation of $G(F)$. By Schur's lemma $Z(G)(F)$ acts by a unitary character on $\pi$. We call it the {\em central character} of $\pi$ and we denote it by $\gls{omegapi}$. We say that an irreducible unitary representation $\pi$ is {\em square-integrable} if for all $e,e'\in \pi$ the function

$$g\in G(F)/A_G(F)\mapsto \lvert(\pi(g)e,e')\rvert$$

\noindent is square-integrable. We will denote by $\gls{Pi2G}$ the set of isomorphism classes of square-integrable representations of $G(F)$. Square-integrable representations are obviously preserved by unramified twists. We will denote by $\gls{Pi2GmodA}$ the set of orbits for this action. For $\pi\in \Pi_2(G)$, we define the {\em formal degree} $\gls{dpi}$ of $\pi$ to be the only positive real number such that

$$\displaystyle \int_{G(F)/A_G(F)} (\pi(g)e_0,e_0')(e_1,\pi(g)e_1')dg=d(\pi)^{-1} (e_0,e_1')(e_1,e_0')$$

\noindent for all $e_0,e_0',e_1,e_1'\in \pi$. Square-integrable representations are tempered, hence we have an inclusion $\Pi_2(G)\subseteq \Temp(G)$.

\vspace{2mm}

\noindent Assume now that $F=\mathbb{R}$. Let $\pi\in \Temp(G)$. Recall that $\mathcal{Z}(\mathfrak{g})$ denotes the center of the enveloping algebra $\mathcal{U}(\mathfrak{g})$. By Schur's lemma, $\mathcal{Z}(\mathfrak{g})$ acts by a character on $\pi^\infty$. This is the {\em infinitesimal character} of $\pi$. We will denote it by $\gls{chipi}$. It is convenient to introduce a norm $\pi\mapsto \gls{normpi}$ on $\Temp(G)$ as follows. Fix a maximal torus $T\subset G$. We have the Harish-Chandra isomorphism

$$\mathcal{Z}(\mathfrak{g})\simeq S(\mathfrak{t})^{W(G_{\mathbb{C}},T_{\mathbb{C}})}$$

\noindent Hence, the set of characters of $\mathcal{Z}(\mathfrak{g})$ gets identified with $\mathfrak{t}^*/W(G_{\mathbb{C}},T_{\mathbb{C}})$. Fix an hermitian norm $\lvert.\rvert$ on $\mathfrak{t}^*$ which is $W(G_{\mathbb{C}},T_{\mathbb{C}})$-invariant. Then, we set

$$N^G(\pi)=1+\lvert \chi_\pi\rvert$$

\noindent for all $\pi\in \Temp(G)$. Note that although the definition of $N^G(\pi)$ depends on some choices, two different choices would give two norms that are equivalent. Since the norm $N^G(.)$ will only be used for the purpose of estimates, the precise choices involved in its definition won't really matter and we will always assume implicitly that such choices have been made. We extend the norm $N^G(.)$ to all tempered representations by

$$\pi=\pi_1\oplus\ldots\oplus\pi_k\mapsto N^G(\pi)=\max\left( N^G(\pi_1),\ldots,N^G(\pi_k)\right)$$

\noindent where $\pi_1,\ldots,\pi_k$ are irreducible tempered representations of $G(F)$ (recall that all our unitary representations have finite length). Later, we will need the following inequality

\begin{num}
\item\label{eq 2.2.8} There exists an integer $k\geqslant 1$ such that
$$d(\pi)\ll N^G(\pi)^k$$
for all $\pi\in \Pi_2(G)$.
\end{num}

\subsection{Parabolic induction}\label{section 2.3}

\noindent Let $P=MU$ be a parabolic subgroup of $G$ and $\sigma$ a tempered of $M(F)$. We extend $\sigma$ to a representation of $P(F)$ trivial on $U(F)$. We will denote by $\gls{iPGsigma}$ the {\em unitary parabolic induction} of $\sigma$. It is a tempered representation of $G(F)$. The space on which $i_P^G(\sigma)$ acts may be described as the completion of the space of continuous functions $e\colon G(F)\to \sigma$ satisfying $e(mug)=\delta_P(m)^{1/2}\sigma(m)e(g)$ for all $m\in M(F)$, $u\in U(F)$ and $g\in G(F)$ for the topology defined by the scalar product

$$\displaystyle (e,e')=\int_{P(F)\backslash G(F)} (e(g),e'(g)) dg$$

\noindent The action of $G(F)$ is given by right translation. The smooth subspace $\gls{iPGsigmainf}$ of $i_P^G(\sigma)$ is exactly the space of smooth functions $e\colon G(F)\to \sigma^\infty$ satisfying the equality $e(mug)=\delta_P(m)^{1/2}\sigma(m)\varphi(g)$ for all $m\in M(F)$, $u\in U(F)$ and $g\in G(F)$. The isomorphism class of $i_P^G(\sigma)$ only depends on $M$ and $\sigma$ and not on $P$. When we only consider this representation modulo isomorphism, we will denote it by $\gls{iMGsigma}$.

\vspace{2mm}

\noindent If $F=\mathbb{R}$, recall that in the previous section we have introduced a norm $N^G$ on the set of (isomorphism classes of) tempered representations of $G(F)$. Of course, this construction also applies to $M$ and yields a norm $N^M$ on $\Temp(M)$. It is not hard to see that we may choose $N^M$ in such a way that

\begin{align}\label{eq 2.3.1}
\displaystyle N^G(i_M^G(\sigma))=N^M(\sigma)
\end{align}

\noindent for all $\sigma\in \Temp(M)$.

\vspace{2mm}

\noindent Let $K$ be a maximal compact subgroup of $G(F)$ that is special in the $p$-adic case. By the Iwasawa decomposition, restriction to $K$ defines a $K$-equivariant isomorphism between $i_P^G(\sigma)$ and $i_{K_P}^K(\sigma_{K_P})$, where $K_P=K\cap P(F)$ and $\sigma_{K_P}$ denotes the restriction of $\sigma$ to $K_P$. This isomorphism restricts to an isomorphism of topological vector spaces $i_P^G(\sigma)^\infty\simeq i_{K_P}^K(\sigma_{K_P})^{\infty}$. Note that if $\lambda\in i\mathcal{A}_M^*$ then $(\sigma_{\lambda})_{K_P}=\sigma_{K_P}$. Hence we get topological isomorphisms $i_P^G(\sigma_\lambda)^\infty\simeq i_{K_P}^K(\sigma_{K_P})^{\infty}$ for all $\lambda\in i\mathcal{A}_M^*$.

\begin{lem}\label{lemma 2.3.1}
\begin{enumerate}[(i)]
\item For every tempered representation $\pi$ of $G(F)$, the linear map

$$\Tra_\pi\colon  \End(\pi)^\infty\to \mathcal{C}^w(G(F))$$

$$T\mapsto \left(g\mapsto \Tr(\pi(g)T)\right)$$

\noindent is continuous.

\item Let $P=MU$ be a parabolic subgroup of $G$ and $\sigma\in \Pi_2(M)$. Let $K$ be a maximal compact subgroup of $G(F)$ which is special in the $p$-adic case. Set $\pi_K=i_{P\cap K}^K(\sigma_{\mid P\cap K})$ and $\pi_\lambda=i_P^G(\sigma_\lambda)$ for all $\lambda\in i\mathcal{A}_M^*$. Consider the isomorphism $\pi_\lambda^\infty\simeq \pi_K^\infty$, $\lambda\in i\mathcal{A}_M^*$, given by restriction to $K$ as an identification. Then the map

$$i\mathcal{A}_M^*\to \Hom\left(\End(\pi_K)^\infty, \mathcal{C}^w(G(F))\right)$$

$$\lambda\mapsto \Tra_{\pi_\lambda}$$

\noindent is smooth.
\end{enumerate}
\end{lem}

\vspace{2mm}

\noindent\ul{Proof}: (i) follows from \ref{eq 2.2.5}. We prove (ii). Let us denote by $(.,.)$ the scalar product on $\pi_K$ given by

$$\displaystyle (e,e')=\int_K (e(k),e'(k))dk$$

\noindent (where the scalar product inside the integral is the scalar product on $\sigma$). Using this scalar product we have an inclusion $\pi_K^\infty\otimes \overline{\pi_K^\infty}\subseteq \End(\pi_K)^\infty$ which induces an identification $\End(\pi_K)^\infty=\pi_K^\infty\widehat{\otimes}_p \overline{\pi_K^\infty}$. Hence, by \ref{eq A.5.1}, it suffices to show that for all $e,e'\in \pi_K^\infty$ the map

$$\displaystyle \lambda\in i\mathcal{A}_M^*\mapsto \left(g\mapsto (\pi_\lambda(g)e,e')\right)\in \mathcal{C}^w(G(F))$$

\noindent is smooth. Fix two vectors $e,e'\in \pi_K^\infty$ and set

$$\displaystyle \varphi(\lambda,g)=(\pi_\lambda(g)e,e')$$

\noindent for all $g\in G(F)$ and all $\lambda\in i\mathcal{A}_M^*$. We are going to apply Lemma \ref{lemma 1.5.1}(ii) to this function $\varphi$. We need to check the various hypothesis of this lemma. The condition (a) in the $p$-adic case, is obvious. Let us show that, in the real case, the condition (b) holds. Let $u,v\in \mathcal{U}(\mathfrak{g})$. Then, we have

\begin{align}\label{eq 2.3.2}
\displaystyle \left(R(u)L(v)\varphi\right)(\lambda,g)=\left(\pi_\lambda(g)\pi_\lambda(u)e,\pi_\lambda(v)e'\right)
\end{align}

\noindent for all $\lambda\in i\mathcal{A}_M^*$ and all $g\in G(F)$ and it suffices to show that the two maps

\begin{align}\label{eq 2.3.3}
\displaystyle \lambda\in i\mathcal{A}_M^*\mapsto \pi_\lambda(g)\pi_\lambda(u)e\in \pi_K^\infty
\end{align}

\begin{align}\label{eq 2.3.4}
\displaystyle \lambda\in i\mathcal{A}_M^*\mapsto \pi_\lambda(v)e'\in \pi_K^\infty
\end{align}

\noindent are smooth. Denote by $\mathfrak{k}^\perp$ the orthogonal of $\mathfrak{k}$ in $\mathfrak{g}$ for the form $B(.,.)$. Fix a basis $X_1,\ldots,X_k$ of $\mathfrak{k}(\mathbb{R})$ such that $B(X_i,X_j)=-\delta_{i,j}$ for $i,j=1,\ldots,k$ and fix a basis $Y_1,\ldots,Y_p$ of $\mathfrak{k}^\perp(\mathbb{R})$ such that $B(Y_i,Y_j)=\delta_{i,j}$ for $i,j=1,\ldots,p$. Set

$$\displaystyle \Delta_K=1-X_1^2-\ldots-X_k^2\in \mathcal{U}(\mathfrak{k})$$

$$\displaystyle \Delta_G=\Delta_K-Y_1^2-\ldots-Y_p^2\in \mathcal{U}(\mathfrak{g})$$

\noindent Using elliptic regularity \ref{eq 2.1.2}, we easily see that the smoothness of \ref{eq 2.3.3} and \ref{eq 2.3.4} follows from the following claim

\vspace{3mm}

\begin{num}
\item\label{eq 2.3.5} For all $f\in C_c(G(F))$, the maps
$$\lambda\in i\mathcal{A}_M^*\mapsto \pi_\lambda(\Delta_G)\in \End(\pi_K^\infty)$$
$$\lambda\in i\mathcal{A}_M^*\mapsto \pi_\lambda(f)\in \End(\pi_K^\infty)$$
are smooth.
\end{num}

\vspace{3mm}

\noindent The difference $C_G=2\Delta_K-\Delta_G$ is in the center $\mathcal{Z}(\mathfrak{g})$ of $\mathcal{U}(\mathfrak{g})$ and the map

$$\displaystyle \lambda\in i\mathcal{A}_M^*\mapsto \chi_{\pi_\lambda}(C_G)$$

\noindent is easily seen to be smooth. This shows the first point of \ref{eq 2.3.5} since the map $\lambda\in i\mathcal{A}_M^*\mapsto \pi_\lambda(\Delta_K)=\pi_K(\Delta_K)$ is constant. For the second point, we first notice that

\begin{align}\label{eq 2.3.6}
\displaystyle \pi_\lambda(f)=\int_{G(F)}f(g)m_\lambda(g)\pi_0(g)dg
\end{align}

\noindent for all $\lambda\in i\mathcal{A}_M^*$, where $m_\lambda(g)$ is the operator that multiply a function $e_0\in \pi_K^\infty$ by the function

$$\displaystyle k\in K\mapsto e^{\langle \lambda, H_P(kg)\rangle}$$

\noindent where $H_P\colon G(F)\to \mathcal{A}_M$ is the extension of $H_M$ to $G(F)$ associated to $K$. The function $g\in G(F)\mapsto H_P(g)$ is easily seen to be smooth. It follows that the map $\lambda\in i\mathcal{A}_M^*\mapsto m_\lambda(g)\in \End(\pi_K^\infty)$ is smooth and its derivatives are easy to compute. By \ref{eq 2.3.6} and the general theorem of differentiation under the integral sign, we deduce that the map $\lambda\in i\mathcal{A}_M^*\mapsto \pi_\lambda(f)$ is smooth. This ends the proof of \ref{eq 2.3.5}.

\vspace{2mm}

\noindent Next we need to check the conditions (c) (in the $p$-adic case) and (d) (in the real case) of Lemma \ref{lemma 1.5.1}. Let $D=\partial(\lambda_1)\ldots\partial(\lambda_n)$ where $\lambda_1,\ldots,\lambda_n\in i\mathcal{A}_M^*$. In the real case, by what we just saw, the functions $\lambda\mapsto \pi_\lambda(u)e$ and $\lambda\mapsto \pi_\lambda(v)e'$ are smooth for all $u,v\in \mathcal{U}(\mathfrak{g})$. Hence, in both the $p$-adic and the real case it suffices to show the existence of a continuous semi-norm $\nu$ on $\pi_K^\infty$ such that

\begin{align}\label{eq 2.3.7}
\displaystyle \left\lvert D_\lambda (\pi_\lambda(g)e_0,e_1)\right\rvert\leqslant \nu(e_0)\nu(e_1) \Xi^G(g)\sigma(g)^n
\end{align}

\noindent for all $\lambda\in i\mathcal{A}_M^*$, all $g\in G(F)$ and all $e_0,e_1\in \pi_K^\infty$. We have

$$\displaystyle (\pi_\lambda(g)e_0,e_1)=\int_K e^{\langle \lambda, H_M(m_P(kg))\rangle} \delta_P(m_P(kg))^{1/2}\left(\sigma\left(m_P(kg)\right)e_0\left(k_P(kg)\right),e_1(k)\right)dk$$

\noindent for all $\lambda\in i\mathcal{A}_M^*$, all $g\in G(F)$, all $e_0,e_1\in \pi_K^\infty$ and where $m_P\colon G(F)\to M(F)$ is as before and $k_P\colon G(F)\to K$ is any map such that $g\in m_P(g)U(F)k_P(g)$ for all $g\in G(F)$. It follows that

\[\begin{aligned}
\displaystyle D_\lambda(\pi_\lambda(g)e_0,e_1)=\int_K \prod_{i=1}^n \left\langle \lambda_i,H_M(m_P(kg))\right\rangle & e^{\langle \lambda, H_M(m_P(kg))\rangle} \delta_P(m_P(kg))^{1/2} \\
 & \left(\sigma\left(m_P(kg)\right)e_0\left(k_P(kg)\right),e_1(k)\right)dk
\end{aligned}\]

\noindent for all $\lambda\in i\mathcal{A}_M^*$, all $g\in G(F)$ and all $e_0,e_1\in \pi_K^\infty$. Obviously we have an inequality

$$\displaystyle \prod_{i=1}^n \left\lvert\left\langle \lambda_i,H_M(m_P(kg))\right\rangle\right\rvert\ll \sigma(g)^n$$

\noindent for all $g\in G(F)$ and all $k\in K$. Hence,

$$\displaystyle \left\lvert D_\lambda(\pi_\lambda(g)e_0,e_1)\right\rvert \ll \sigma(g)^n\int_K \delta_P(m_P(kg))^{1/2}\left\lvert\left(\sigma\left(m_P(kg)\right)e_0\left(k_P(kg)\right),e_1(k)\right)\right\rvert dk$$

\noindent for all $\lambda\in i\mathcal{A}_M^*$, all $g\in G(F)$ and all $e_0,e_1\in \pi_K^\infty$. By \ref{eq 2.2.3} and \ref{eq 2.2.6} there exists a continuous semi-norm $\nu_\sigma$ on $\sigma^\infty$ such that $\left\lvert (\sigma(m)v_0,v_1)\right\rvert\leqslant \nu_\sigma(v_0)\nu_\sigma(v_1)\Xi^M(m)$ for all $v_0,v_1\in \sigma^\infty$ and all $m\in M(F)$. It follows that

\[\begin{aligned}
\displaystyle \left\lvert D_\lambda(\pi_\lambda(g)e_0,e_1)\right\rvert & \ll \sup_{k\in K}\left[\nu_\sigma\left(e_0(k)\right)\right] \sup_{k\in K}\left[\nu_\sigma\left(e_1(k)\right)\right]\sigma(g)^n\int_K \delta_P(m_P(kg))^{1/2} \Xi^M(m_P(kg))dk \\
 & =\sup_{k\in K}\left[\nu_\sigma(e_0(k))\right] \sup_{k\in K}\left[\nu_\sigma(e_1(k))\right]\Xi^G(g)\sigma(g)^n
\end{aligned}\]

\noindent \noindent for all $\lambda\in i\mathcal{A}_M^*$, all $g\in G(F)$ and all $e_0,e_1\in \pi_K^\infty$, where in the last equality we used Proposition \ref{proposition 1.5.1}(iii). By the standard Sobolev inequality the semi-norm

$$e\mapsto \sup_{k\in K}\left[\nu_\sigma\left(e(k)\right)\right]$$

\noindent is continuous on $\pi_K^\infty$. This proves \ref{eq 2.3.7} and ends the proof of the lemma. $\blacksquare$

\subsection{Normalized intertwining operators}\label{section 2.4}

\noindent Let $M$ be a Levi subgroup of $G$, $\sigma$ be a tempered representation of $M(F)$ and fix $K$ a maximal compact subgroup of $G(F)$ which is special in the $p$-adic case. The definition of the representations $i_P^G(\sigma_\lambda)^\infty$, $P\in \mathcal{P}(M)$, actually still makes sense for any $\lambda\in \mathcal{A}_{M,\mathbb{C}}^*$ (although these are not anymore unitary representations in general). Moreover, restriction to $K$ still induces $K$-equivariant isomorphisms $i_P^G(\sigma_\lambda)^\infty\simeq i_{K_P}^K(\sigma_{K_P})^\infty$ for all $\lambda\in \mathcal{A}_{M,\mathbb{C}}^*$ and all $P\in \mathcal{P}(M)$ (recall that $K_P=P(F)\cap K$).

\vspace{2mm}

\noindent Let $P=MU,P'=MU'\in \mathcal{P}(M)$. For $Re(\lambda)$ in a certain open cone, the expression

$$\displaystyle \left(\gls{JPPsigma}e\right)(g):=\int_{\left(U(F)\cap U'(F)\right)\backslash U'(F)} e(u'g) du'$$

\noindent is absolutely convergent for all $e\in i_P^G(\sigma_\lambda)^\infty$ and defines a $G(F)$-equivariant continuous linear map

$$\displaystyle J_{P'\mid P}(\sigma_\lambda)\colon i_P^G(\sigma_\lambda)^\infty\to i_{P'}^G(\sigma_\lambda)^\infty$$

\noindent Via the isomorphisms $i_P^G(\sigma_\lambda)^\infty\simeq i_{K_P}^K(\sigma_{K_P})^\infty$ and $i_{P'}^G(\sigma_\lambda)\simeq i_{K_{P'}}^K(\sigma_{K_{P'}})^\infty$, we can view the map $\lambda\mapsto J_{P'\mid P}(\sigma_\lambda)$ as taking values in $\Hom\left(i_{K_P}^K(\sigma_{K_P})^\infty,i_{K_{P'}}^K(\sigma_{K_{P'}})^\infty\right)$ the space of continuous linear maps between $i_{K_P}^K(\sigma_{K_P})^\infty$ and $i_{K_{P'}}^K(\sigma_{K_{P'}})^\infty$. This function admits a meromorphic continuation to $\mathcal{A}_{M,\mathbb{C}}^*$ (see \cite{Wall2} Theorem 10.1.6 in the Archimedean case and \cite{Wa2} Th\'eor\`eme IV.1.1 in the $p$-adic case).

\vspace{2mm}

\noindent Let $\overline{P}=M\overline{U}\in \mathcal{P}(M)$ be the parabolic subgroup opposite to $P=MU$. Assume that the Haar measures $du$ and $d\overline{u}$ on $U(F)$ and $\overline{U}(F)$ have been normalized so that $dg=\delta_P(m)^{-1}dudmd\overline{u}$ where $dg$ and $dm$ denotes the Haar measures on $dg$ and $dm$ respectively. Then the meromorphic function $\lambda\mapsto \gls{jsigma}=J_{P\mid \overline{P}}(\sigma_\lambda)J_{\overline{P}\mid P}(\sigma_\lambda)$ is scalar-valued and doesn't depend on the choice of $P$. Moreover this function takes on $i\mathcal{A}_M^*$ positive real values (including $\infty$). We will need the following:

\vspace{3mm}

\begin{num}
\item\label{eq 2.4.1} Assume that $F=\mathbb{R}$. Then, there exists an integer $k\geqslant 1$ such that
$$j(\sigma)^{-1}\ll N^M(\sigma)^k$$
for all $\sigma\in \Temp(M)$.
\end{num}

\vspace{3mm}

\noindent We will also need {\em normalized intertwining operators}. In the Archimedean case, such normalizations have been defined and extensively studied by Knapp and Stein in \cite{KSI} and \cite{KSII}. However, in this paper we shall prefer Arthur's normalization \cite{A5} which better fits our purposes and takes care of both the Archimedean and non-Archimedean cases. The general construction is as follows. There exist complex-valued meromorphic functions $\lambda\mapsto r_{P'\mid P}(\sigma_\lambda)$, for all $P,P'\in \mathcal{P}(M)$, such that if we set
$$\displaystyle \gls{RPPsigma}=r_{P'\mid P}(\sigma_\lambda)^{-1} J_{P'\mid P}(\sigma_\lambda),\;\;\;P,P'\in\mathcal{P}(M), \lambda\in \mathcal{A}_{M,\mathbb{C}}^*$$
these operators satisfy the conditions (R1)-(R8) of the Section 2 of \cite{A5}. The most important conditions for us will be the following:

\vspace{2mm}

\begin{num}
\item\label{eq 2.4.2} $R_{P''\mid P'}(\sigma_\lambda)R_{P'\mid P}(\sigma_\lambda)=R_{P''\mid P}(\sigma_\lambda)$ for all $P,P',P''\in \mathcal{P}(M)$;

\item\label{eq 2.4.3} $\lambda\mapsto R_{P'\mid P}(\sigma_\lambda)$ is holomorphic and unitary on $i\mathcal{A}_M^*$ for all $P',P\in \mathcal{P}(M)$;

\item\label{eq 2.4.5} If $\lambda\in i\mathcal{A}_G^*$, then $R_{P'\mid P}(\sigma_\lambda)=R_{P'\mid P}(\sigma)_{\lambda}$ via the natural isomorphisms $i_P^G(\sigma_\lambda)\simeq i_P^G(\sigma)_\lambda$ and $i_{P'}^G(\sigma_\lambda)\simeq i_{P'}^G(\sigma)_\lambda$;

\item\label{eq 2.4.6} For $P,P'\in \mathcal{P}(M)$, if $Q=LU_Q$ denotes the parabolic subgroup generated by $P$ and $P'$ then $R_{P'\mid P}(\sigma_\lambda)=i_Q^G\left(R_{P'\cap L\mid P\cap L}(\sigma_\lambda)\right)$ via the isomorphisms of induction by stages $i_P^G(\sigma_\lambda)\simeq i_Q^G\left(i_{P\cap L}^{L}(\sigma_\lambda)\right)$ and $i_{P'}^G(\sigma_\lambda)\simeq i_Q^G\left(i_{P'\cap L}^{L}(\sigma_\lambda)\right)$;

\item\label{eq 2.4.7} For all $g\in G(F)$ and all $P,P'\in \mathcal{P}(M)$, we have $R_{gP'g^{-1}\mid gPg^{-1}}(g\sigma g^{-1})=A_{P'}(g)R_{P'\mid P}(\sigma)A_P(g)^{-1}$ where $A_P(g)$ is the isomorphism $i_P^G(\sigma)\simeq i_{gPg^{-1}}^G(g\sigma g^{-1})$ given by $(A_P(g)e)(\gamma)=e(g^{-1}\gamma)$ (and $A_{P'}(g)$ is defined similarly).

\item\label{eq 2.4.8} Assume that $F=\mathbb{R}$. Then, for every differential operator with constant coefficients $D$ on $i\mathcal{A}_M^*$ and for all $P,P'\in \mathcal{P}(M)$, there exist $k,r\geqslant 1$ such that
$$\displaystyle \left\lVert D_\lambda R_{P'\mid P}(\sigma_\lambda)e\right\rVert \leqslant \lVert i_P^G(\sigma_\lambda,\Delta_K^r) e\rVert N^M(\sigma_\lambda)^k$$
for all $\sigma\in \Temp(M)$, all $\lambda\in i\mathcal{A}_M^*$ and all $e\in i_P^G(\sigma_\lambda)^\infty$. 
\end{num}

\noindent Properties \ref{eq 2.4.2}, \ref{eq 2.4.3} and \ref{eq 2.4.7} correspond to Arthur's conditions (R3), (R4) and (R6) respectively. Property \ref{eq 2.4.5} is a direct consequence of the requirement (r.1) p.171 of \cite{A5} on the normalizing factors $r_{P'\mid P}(\sigma_\lambda)$ (in that they only depend on the projection of $\lambda$ to $(\mathcal{A}^G_{M,\mathbb{C}})^*$). The identity \ref{eq 2.4.6} also follows from the same requirement (r.1) of \cite{A5} on normalizing factors together with the analogous property of unnormalized intertwining operators. Finally, the last condition \ref{eq 2.4.8} is a consequence of Lemma 2.1 of \cite{A5}.

\subsection{Weighted characters}\label{section 2.5}

\noindent We keep the notation and assumptions of the previous section : $M$ is a Levi subgroup of $G$, $\sigma$ a tempered representation of $M(F)$ and $K$ a maximal compact subgroup of $G(F)$ that is special in the $p$-adic case. Fix $P\in \mathcal{P}(M)$. For all $P'\in \mathcal{P}(M)$, we may consider the function $\mathcal{R}_{P'}(\sigma,P)$ on $i\mathcal{A}_M^*$ defined by

$$\displaystyle \mathcal{R}_{P'}(\lambda,\sigma,P)=R_{P'\mid P}(\sigma)^{-1}R_{P'\mid P}(\sigma_{\lambda})$$

\noindent The family $\left(\mathcal{R}_{P'}(\sigma,P)\right)_{P'\in \mathcal{P}(M)}$ is a $(G,M)$-family taking values in $\End(i_{K_P}^K(\sigma_{K_P})^\infty)$ (\cite{A9} p.43). Following Arthur, we may associate to this family operators $\mathcal{R}_L^Q(\sigma,P)$ in $\End(i_{K_P}^K(\sigma_{K_P})^\infty)$  for all $L\in \mathcal{L}(M)$ and all $Q\in \mathcal{F}(L)$ (cf.\ Section \ref{section 1.9}). Then, for all $f\in \mathcal{C}(G(F))$ all $L\in \mathcal{L}(M)$ and all $Q\in\mathcal{F}(L)$, we set

$$\displaystyle \gls{JLQsigmaf}=\Tr(\mathcal{R}_L^Q(\sigma,P)i_P^G(\sigma,f))$$

\noindent The trace is well-defined in the $p$-adic case since then $i_P^G(\sigma,f)$ is a finite rank operator. To see that it is also well-defined in the real case, we may proceed as follows: by the factorization \ref{eq 2.1.1} and by linearity, we may assume that $f=f_1\ast f_2$ where $f_1\in C_c^\infty(G(F))$ and $f_2\in \mathcal{C}(G(F))$. Then $\mathcal{R}_L^Q(\sigma,P)i_P^G(\sigma,f_1)$ extends continuously to an endomorphism of $i_P^G(\sigma)$ and since $i_P^G(\sigma,f_2)$ is traceable so is $\mathcal{R}_L^Q(\sigma,P)i_P^G(\sigma,f)=\mathcal{R}_L^Q(\sigma,P)i_P^G(\sigma,f_1)i_P^G(\sigma,f_2)$.

\vspace{2mm}

\noindent  This defines a family of tempered distributions $\left(J_L^Q(\sigma,.)\right)_{L,Q}$ on $G(F)$ which doesn't depend on $P$ but depends on $K$ and the way we normalized the intertwining operators. Note that if $L=Q=G$, this reduces to the usual character, that is

$$J^G_G(\sigma,f)=\Tr(i_M^G(\sigma,f))$$

\noindent for all $f\in \mathcal{C}(G(F))$.

\begin{lem}\label{lemma 2.5.1}
Assume $F=\mathbb{R}$. Let $L\in \mathcal{L}(M)$ and $Q\in\mathcal{F}(L)$. Then, for all $k\geqslant 0$, there exists a continuous semi-norm $\nu_k$ on $\mathcal{C}(G(F))$ such that

$$\displaystyle \left\lvert J_L^Q(\sigma,f)\right\rvert\leqslant \nu_k(f) N^M(\sigma)^{-k}$$

\noindent for all $\sigma\in \Temp(M)$.
\end{lem}

\vspace{2mm}

\noindent\ul{Proof}: For all $z\in \mathcal{Z}(\mathfrak{g})$, we have

$$J_L^Q(\sigma,zf)=\chi_\sigma(z_M)J_L^Q(\sigma,f)$$

\noindent for all $\sigma\in \Temp(M)$ and all $f\in \mathcal{C}(G(F))$. Moreover, there exist $z_1,\ldots,z_n\in \mathcal{Z}(\mathfrak{g})$ such that

$$\lvert \chi_\sigma(z_1)\rvert+\ldots+\lvert \chi_\sigma(z_n)\rvert\geqslant N^M(\sigma)$$

\noindent for all $\sigma\in \Temp(M)$. Consequently, we only need to show the following

\vspace{3mm}

\begin{num}
\item\label{eq 2.5.1} There exists $k\geqslant 0$ and a continuous semi-norm $\nu$ on $\mathcal{C}(G(F))$ such that
$$\left\lvert J_L^Q(\sigma,f)\right\rvert\leqslant \nu(f)N^M(\sigma)^k$$
for all $\sigma\in \Temp(M)$ and all $f\in \mathcal{C}(G(F))$.
\end{num}

\vspace{3mm}

\noindent It follows from \ref{eq 2.4.8} that we may find two integers $k\geqslant 0$ and $r\geqslant 0$ such that

$$\displaystyle \left\lVert \mathcal{R}_L^Q(\sigma,P)e\right\rVert\leqslant \lVert i_P^G(\sigma,\Delta_K^r)e\rVert N^M(\sigma)^k$$

\noindent for all $\sigma\in \Temp(M)$ and all $e\in i_P^G(\sigma)^\infty$.

\vspace{2mm}

\noindent Let $\sigma\in \Temp(M)$ and $f\in \mathcal{C}(G(F))$. For all $\rho\in \widehat{K}$, let us fix an orthonormal basis $\mathcal{B}_\rho(\sigma)$ of $i_P^G(\sigma)(\rho)$ (the $\rho$ isotypic component of $i_P^G(\sigma)$). Then, for every integer $\ell\geqslant 1$, we have

\[\begin{aligned}
\displaystyle \left\lvert J_L^Q(\sigma,f)\right\rvert & = \left\lvert \Tr\left(\mathcal{R}_M(\sigma,P)i_P^G(\sigma,f)\right) \right\rvert \\
 & \leqslant \sum_{\rho\in \widehat{K}} \sum_{e\in \mathcal{B}_\rho(\sigma)} \left\lvert \left( \mathcal{R}_M(\sigma,P)i_P^G(\sigma,f)e,e\right)\right\rvert \\
 & \leqslant N^M(\sigma)^k \sum_{\rho\in \widehat{K}} \sum_{e\in \mathcal{B}_\rho(\sigma)} \left\lVert i_P^G(\sigma,L(\Delta_K^r)f)e\right\rVert \\
 & =N^M(\sigma)^k \sum_{\rho\in \widehat{K}} \sum_{e\in \mathcal{B}_\rho(\sigma)} c(\rho)^{-\ell}\left\lVert i_P^G(\sigma,L(\Delta_K^r)R(\Delta_K^\ell)f)e\right\rVert
\end{aligned}\]

\noindent By \ref{eq 2.2.7}, there exists a continuous semi-norm $\nu_\ell$ on $\mathcal{C}(G(F))$ which doesn't depend on $\sigma$ and such that the sum above is bounded by

$$\displaystyle \nu_\ell(f)N^M(\sigma)^k\sum_{\rho\in \widehat{K}} \sum_{e\in \mathcal{B}_\rho(\sigma)} c(\rho)^{-\ell}$$

\noindent and by \ref{eq 2.2.2}, if $\ell$ is sufficiently large the sum above is absolutely convergent and bounded by a constant independent of $\sigma$. This shows \ref{eq 2.5.1} and ends the proof of the lemma. $\blacksquare$

\subsection{Matricial Paley-Wiener theorem and Plancherel-Harish-Chandra theorem}\label{section 2.6}

\noindent Let us define $\gls{XtempG}$ to be the set of isomorphism classes of tempered representations of $G(F)$ which are of the form $i_M^G(\sigma)$ where $M$ is a Levi subgroup of $G$ and $\sigma\in \Pi_2(M)$ is a square-integrable representation. According to Harish-Chandra two such representations $i_M^G(\sigma)$ and $i_{M'}^G(\sigma')$ are isomorphic if and only if the pairs $(M,\sigma)$ and $(M',\sigma')$ are conjugate under $G(F)$. Let $\mathcal{M}$ be a set of representatives for the conjugacy classes of Levi subgroups of $G$. Then, $\mathcal{X}_{\tempe}(G)$ is naturally a quotient of

$$\displaystyle \widetilde{\mathcal{X}}_{\tempe}(G)=\bigsqcup_{M\in \mathcal{M}} \bigsqcup_{\mathcal{O}\in \Pi_2(M)/i\mathcal{A}_{M,F}^*} \mathcal{O}$$

\noindent which has a natural structure of real smooth manifold since each orbit $\mathcal{O}\in \{\Pi_2(M)\}$ is a quotient of $i\mathcal{A}_{M,F}^*$ by a finite subgroup hence is naturally a real smooth manifold. We equip $\mathcal{X}_{\tempe}(G)$ with the quotient topology. Note that the connected components of $\mathcal{X}_{\tempe}(G)$ are the image of unramified classes $\mathcal{O}\in \{\Pi_2(M)\}$, $M\in \mathcal{M}$. The following is due to Harish-Chandra (cf.\ Theorem VIII.1.2 of \cite{Wa2})

\vspace{3mm}

\begin{num}
\item\label{eq 2.6.1} If $F$ is $p$-adic, then for every compact-open subgroup $K\subset G(F)$ the set 
$$\displaystyle \{\pi\in \mathcal{X}_{\tempe}(G); \; \pi^K\neq 0\}$$
is relatively compact in $\mathcal{X}_{\tempe}(G)$ (i.e., is contained in the union of a finite number of components).
\end{num}

\vspace{3mm}

\noindent Let $V$ be a locally convex topological vector space. We will say of a function $f\colon  \mathcal{X}_{\tempe}(G)\to V$ that it is {\em smooth} if the pullback of $f$ to $\widetilde{\mathcal{X}}_{\tempe}(G)$ is a smooth function. We will denote by $\gls{C(Xtemp,V)}$ the space of smooth functions on $\mathcal{X}_{\tempe}(G)$ taking values in $V$. We will also simply set $\gls{C(Xtemp)}=C^\infty(\mathcal{X}_{\tempe}(G),\mathbb{C})$.

\vspace{2mm}

\noindent We define a regular Borel measure $\gls{dpimes}$ on $\mathcal{X}_{\tempe}(G)$ by requesting that

$$\displaystyle \int_{\mathcal{X}_{\tempe}(G)} \varphi(\pi)d\pi=\sum_{M\in \mathcal{M}} \left\lvert W(G,M)\right\rvert^{-1}\sum_{\mathcal{O}\in \Pi_2(M)/i\mathcal{A}_{M,F}^*} [i\mathcal{A}^\vee_{M,\sigma}:i\mathcal{A}^\vee_{M,F}]^{-1} \int_{i\mathcal{A}^*_{M,F}} \varphi(i_P^G(\sigma_\lambda)) d\lambda$$

\noindent for all $\varphi\in C_c(\mathcal{X}_{\tempe}(G))$, where for all $M\in \mathcal{M}$ we have fixed $P\in \mathcal{P}(M)$ and for all $\mathcal{O}\in \{\Pi_2(M)\}$ we have fixed a base-point $\sigma\in \mathcal{O}$.

\vspace{2mm}

\noindent For all $\pi=i_M^G(\sigma)\in \mathcal{X}_{\tempe}(G)$, we set $\gls{mupi}=d(\sigma)j(\sigma)^{-1}$. This quantity really only depends on $\pi$ since another pair $(M',\sigma')$ yielding $\pi$, where $M'$ is a Levi subgroup and $\sigma'\in \Pi_2(M')$, is $G(F)$-conjugate to $(M,\sigma)$.

\vspace{2mm}

\noindent Assume that $F=\mathbb{R}$. Recall that we defined in Section \ref{section 2.2} a norm $N^G$ on the set of (isomorphism classes of) tempered representations of $G(F)$. By \ref{eq 2.2.8}, \ref{eq 2.4.1} and \ref{eq 2.3.1}, there exists $k\geqslant 0$ such that $\mu(\pi)\ll N^G(\pi)^k$ for all $\pi\in \mathcal{X}_{\tempe}(G)$. The following basic estimate will be used several times:

\vspace{3mm}

\begin{num}
\item\label{eq 2.6.2} There exists an integer $k\geqslant 1$ such that the integral
$$\displaystyle \int_{\mathcal{X}_{\tempe}(G)}N^G(\pi)^{-k}d\pi$$
is absolutely convergent.
\end{num}

\vspace{3mm}

\noindent We are now going to define a space of functions $\gls{C(Xtemp,E)}$. The elements of that space are certain assignments $T\colon \pi\in \Temp(G)\mapsto T_\pi\in \End(\pi)^\infty$ (notice that for all $\pi\in \Temp(G)$, the space $\End(\pi)^\infty$ is well-defined up to a unique isomorphism). First, we extend such an assignment to all (isomorphism classes of) tempered representations by $\pi=\pi_1\oplus\ldots\oplus \pi_k\mapsto T_\pi=T_{\pi_1}\oplus\ldots\oplus T_{\pi_k}\in \End(\pi)^\infty$ where the $\pi_i$'s are irreducible. Let us fix a maximal compact subgroup $K$ of $G(F)$ which is special in the $p$-adic case. We may now define $C^\infty(\mathcal{X}_{\tempe}(G),\mathcal{E}(G))$ as the space of functions $\pi\in \Temp(G)\mapsto T_\pi\in \End(\pi)^\infty$ such that for any parabolic subgroup $P=MU$ and for all $\sigma\in \Pi_2(M)$, setting $\pi_K=i_{P\cap K}^K(\sigma_{\mid P\cap K})$ and $\pi_\lambda=i_P^G(\sigma_\lambda)$ for all $\lambda\in i\mathcal{A}_M^*$, the function

$$\displaystyle \lambda\in i\mathcal{A}_M^*\mapsto T_{\pi_\lambda}\in \End(\pi_\lambda)^\infty\simeq \End(\pi_K)^{\infty}$$

\noindent is smooth.

\vspace{2mm}

\noindent We define a subspace $\gls{Cc(Xtemp,E)}$ of $C^\infty(\mathcal{X}_{\tempe}(G),\mathcal{E}(G))$ as follows. This is the subspace of sections $T\in C^\infty(\mathcal{X}_{\tempe}(G),\mathcal{E}(G))$ such that

\vspace{2mm}

\begin{itemize}
\renewcommand{\labelitemi}{$\bullet$}
\item in the $p$-adic case: $\Supp(T)=\overline{\{\pi\in \mathcal{X}_{\tempe}(G); \; T_\pi\neq 0\}}$ is compact (i.e., is contained in a finite union of connected components);

\item in the real case: for every parabolic subgroup $P=MU$ and for every differential operator with constant coefficient $D$ on $i\mathcal{A}_M^*$ the function

$$\displaystyle DT\colon \sigma\in \Pi_2(M)\mapsto D_\lambda\left(\lambda\mapsto T_{i_P^G(\sigma_\lambda)}\in \End(i_{P\cap K}^K(\sigma))^{\infty}\right)$$

\noindent has the property that

$$\displaystyle p_{D,u,v,k}(T)=\sup_{\sigma\in \Pi_2(M)} \vertiii{ (DT)_\sigma}_{u,v}N(\sigma)^k<\infty$$

\noindent for all $u,v\in \mathcal{U}(\mathfrak{k})$ and all $k\in \mathbb{N}$.
\end{itemize}

\vspace{2mm}

\noindent We equip the space $\mathcal{C}(\mathcal{X}_{\tempe}(G),\mathcal{E}(G))$ with a locally convex topology as follows. If $F=\mathbb{R}$, it is the topology defined by the semi-norms $p_{D,u,v,k}$ for all $D,u,v$ and $k$ as above. If $F$ is $p$-adic, we remark that $\mathcal{C}(\mathcal{X}_{\tempe}(G),\mathcal{E}(G))$ is naturally a subspace of

\begin{align}\label{eq 2.6.3}
\displaystyle \bigoplus_{M\in \mathcal{M}}\bigoplus_{\mathcal{O}\in\{\Pi_2(M)\}} C^\infty\left(i\mathcal{A}_{M,F}^*, \End(i_{K_P}^K(\sigma_{K_P}))^\infty\right)
\end{align}

\noindent where for all $M\in \mathcal{M}$ we have fixed a parabolic subgroup $P\in \mathcal{P}(M)$ and for all $\mathcal{O}\in \{\Pi_2(M)\}$ we have fixed a base-point $\sigma\in \mathcal{O}$. The spaces $C^\infty\left(i\mathcal{A}_{M,F}^*, \End(i_{K_P}^K(\sigma_{K_P}))^\infty\right)$ have natural locally convex topologies. We endow the space \ref{eq 2.6.3} with the direct sum topology and $\mathcal{C}(\mathcal{X}_{\tempe}(G),\mathcal{E}(G))$ with the subspace topology.

\vspace{2mm}

\noindent We will need the following strong version of the Harish-Chandra Plancherel formula also called matricial Paley-Wiener theorem (cf.\ Theorem VII.2.5 and Theorem VIII.1.1 of \cite{Wa2} in the $p$-adic case and \cite{A2}, \cite{A8} in the real case).

\vspace{3mm}

\begin{theo}\label{theorem 2.6.1}
\begin{enumerate}[(i)]
\item The map $f\in \mathcal{C}(G)\mapsto \left(\pi\in \Temp(G)\mapsto \pi(f)\in \End(\pi)^\infty\right)$ induces a topological isomorphism $\mathcal{C}(G)\simeq \mathcal{C}(\mathcal{X}_{\tempe}(G),\mathcal{E}(G))$.

\item The inverse of that isomorphism is given by sending $T\in \mathcal{C}(\mathcal{X}_{\tempe}(G),\mathcal{E}(G))$ to the function $f_T$ defined by

$$\displaystyle f_T(g)=\int_{\mathcal{X}_{\tempe}(G)} \Tr\left(\pi(g^{-1})T_\pi\right) \mu(\pi) d\pi$$
\end{enumerate}
\end{theo}

\vspace{2mm}

\noindent\ul{Remark}: The last integral above is absolutely convergent by \ref{eq 2.2.5} and \ref{eq 2.7.2}.

\subsection{Elliptic representations and the space $\mathcal{X}(G)$}\label{section 2.7}

\noindent Denote by $\gls{RtempG}$ the space of complex virtual tempered representations of $G(F)$, that is $R_{\tempe}(G)$ is the complex vector space with basis $\Temp(G)$. We may extend almost all our constructions to virtual representations. In particular:

\vspace{2mm}

\begin{itemize}
\renewcommand{\labelitemi}{$\bullet$}
\item We extend the action of $i\mathcal{A}_{G,F}^*$ by linearity to $R_{\tempe}(G)$;

\item Let $M$ be a Levi subgroup of $G$. Then the functor $i_M^G$ extends by linearity to give a linear map $i_M^G(.):R_{\tempe}(M)\to R_{\tempe}(G)$. Also, we extend the weighted character $\sigma\mapsto J_L^Q(\sigma,.)$ ($L\in \mathcal{L}(M)$, $Q\in\mathcal{F}(L)$) of Section \ref{section 2.5} by linearity to $R_{\tempe}(M)$;

\item If $F=\mathbb{R}$, we extend the norm $N^G$ to $R_{\tempe}(G)$ by

$$N^G(\lambda_1 \pi_1+\ldots+\lambda_k\pi_k)=\max\left(\lvert \lambda_1\rvert N^G(\pi_1),\ldots,\lvert \lambda_k\rvert N^G(\pi_k)\right)$$

\item We will denote by $\pi\mapsto \overline{\pi}$ the unique conjugate-linear extension of $\pi\in \Temp(G)\mapsto \overline{\pi}\in \Temp(G)$ to $R_{\tempe}(G)$.
\end{itemize}

\vspace{2mm}

\noindent In \cite{A4}, Arthur defines a set $T_{\elli}(G)$ of virtual tempered representations of $G(F)$, that we will denote by $\gls{XellG}$ in this paper. The elements of $\mathcal{X}_{\elli}(G)$ are actually well-defined only up to a scalar of module $1$. That is, we have $\mathcal{X}_{\elli}(G)\subset R_{\tempe}(G)/\mathbb{S}^1$. These are the so-called {\em elliptic representations}. Let us recall their definition. Let $P=MU$ be a parabolic subgroup of $G$ and $\sigma\in \Pi_2(M)$. For all $g\in G(F)$, we define $g\sigma$ to be the representation of $gM(F)g^{-1}$ given by $(g\sigma)(m')=\sigma(g^{-1}m'g)$ for all $m'\in gM(F)g^{-1}$. Denote by $\No_{G(F)}(\sigma)$ the subgroup of elements $g\in \No_{G(F)}(M)$ such that $g\sigma\simeq \sigma$ and set $W(\sigma)=\No_{G(F)}(\sigma)/M(F)$. Fix $P\in \mathcal{P}(M)$. Then, we may associate to every $w\in W(\sigma)$ an unitary endomorphism $R_P(w)$ of the representation $i_P^G(\sigma)$ that is well-defined up to a scalar of module $1$ as follows. Choose a lift $\widetilde{w}\in \No_{G(F)}(\sigma)$ of $w$ and an unitary endomorphism $A(\widetilde{w})$ of $\sigma$ such that $\sigma(\widetilde{w}^{-1}m\widetilde{w})=A(\widetilde{w})^{-1}\sigma(m)A(\widetilde{w})$ for all $m\in M(F)$. We define the operator $R_P(w):i_P^G(\sigma)^\infty\to i_P^G(\sigma)^\infty$ as the composition $R_{P\mid wPw^{-1}}(\sigma)\circ \mathcal{A}(\widetilde{w})\circ \lambda(\widetilde{w})$, where

\begin{itemize}
\item $\lambda(\widetilde{w})$ is the isomorphism $i_P^G(\sigma)\simeq i^G_{wPw^{-1}}(\widetilde{w}\sigma)$ given by $(\lambda(\widetilde{w})e)(g)=e(\widetilde{w}^{-1}g)$;

\item $\mathcal{A}(\widetilde{w})$ is the isomorphism $i^G_{wPw^{-1}}(\widetilde{w}\sigma)\simeq i^G_{wPw^{-1}}(\sigma)$ given by $(\mathcal{A}(\widetilde{w})e)(g)=A(\widetilde{w})e(g)$;

\item $R_{P\mid wPw^{-1}}(\sigma): i^G_{wPw^{-1}}(\sigma)^\infty\to i^G_{P}(\sigma)^\infty$ is the normalized intertwining operator defined in Section \ref{section 2.4}.
\end{itemize}

\noindent We immediately check that $R_P(w)$ is $G(F)$-equivariant and that it depends on all the choices ($\widetilde{w}$, $A(\widetilde{w})$ and the normalization of the intertwining operator $R_{P\mid wPw^{-1}}(\sigma)$) only up to a scalar of module $1$. We associate to any $w\in W(\sigma)$ a virtual tempered representation $i_M^G(\sigma,w)$, well-defined up to a scalar of module $1$, by setting

$$\displaystyle i_M^G(\sigma,w)=\sum_{\lambda\in \mathbb{C}} \lambda\; i_P^G(\sigma,w,\lambda)$$

\noindent where for all $\lambda\in \mathbb{C}$, $i_P^G(\sigma,w,\lambda)$ denotes the subrepresentation of $i_P^G(\sigma)$ where $R_P(w)$ acts by multiplication by $\lambda$ (as is indicated in the notation this definition doesn't depend on the choice of $P$). Let $W_0(\sigma)$ be the subgroup of elements $w\in W(\sigma)$ such that $R_P(w)$ is a scalar multiple of the identity and let $W(\sigma)_{\reg}$ be the subgroup of elements $w\in W(\sigma)$ such that $\mathcal{A}_M^w=\mathcal{A}_G$. We will say that the virtual representation $i_M^G(\sigma,w)$, $w\in W(\sigma)$, is {\em elliptic} if $W_0(\sigma)=\{1\}$ and $w\in W(\sigma)_{\reg}$. The set $\mathcal{X}_{\elli}(G)$ is the set of all virtual elliptic representations (well-defined up to multiplication by a scalar of module $1$) that are obtained in this way. Let $\pi\in \mathcal{X}_{\elli}(G)$ and write $\pi=i_M^G(\sigma,w)$ with $M$, $\sigma$ and $w\in W(\sigma)_{\reg}$ as before. Then we set

$$\displaystyle \gls{Dpi}=\lvert \det(1-w)_{\mathcal{A}_M^G}\rvert^{-1} \lvert W(\sigma)_w\rvert^{-1}$$

\noindent where $W(\sigma)_w$ denotes the centralizer of $w$ in $W(\sigma)$. This number doesn't depend on the particular choice of $M$, $\sigma$ and $w$ representing $\pi$ because any other choice yielding $\pi$ will be $G(F)$-conjugate to $(M,\sigma,w)$. The set $\mathcal{X}_{\elli}(G)$ satisfies the following important property. Denote by $\gls{RellG}$ the subspace of $R_{\tempe}(G)$ generated by $\mathcal{X}_{\elli}(G)$ and denote by $\gls{RindG}$ the subspace of $R_{\tempe}(G)$ generated by the image of all the linear maps $i_M^G:R_{\tempe}(M)\to R_{\tempe}(G)$ for $M$ a proper Levi subgroup of $G$. Then we have the decomposition

\begin{align}\label{eq 2.7.1}
\displaystyle R_{\tempe}(G)=R_{ind}(G)\oplus R_{\elli}(G)
\end{align}

\noindent The set $\mathcal{X}_{\elli}(G)$ is invariant under unramified twists. We will denote by $\gls{XellmodA}$ the set of unramified orbits in $\mathcal{X}_{\elli}(G)$. Also, we will denote by $\gls{underXellG}$ the inverse image of $\mathcal{X}_{\elli}(G)$ in $R_{\tempe}(G)$. This set is invariant under multiplication by $\mathbb{S}^1$.

\vspace{2mm}

\noindent We define $\gls{XG}$ to be the subset of $R_{\tempe}(G)/\mathbb{S}^1$ consisting of virtual representations of the form $i_M^G(\sigma)$ where $M$ is a Levi subgroup of $G$ and $\sigma\in \mathcal{X}_{\elli}(M)$. Also, we will denote by $\gls{underXG}$ the inverse image of $\mathcal{X}(G)$ in $R_{\tempe}(G)$. Hence, the fibers of the natural projection $\underline{\mathcal{X}}(G)\to \mathcal{X}(G)$ are all isomorphic to $\mathbb{S}^1$. Let $\mathcal{M}$ be a set of representatives for the conjugacy classes of Levi subgroups of $G$. Then, $\mathcal{X}(G)$ is naturally a quotient of

$$\bigsqcup_{M\in \mathcal{M}} \bigsqcup_{\mathcal{O}\in \mathcal{X}_{\elli}(M)/i\mathcal{A}_{M,F}^*} \mathcal{O}$$

\noindent This defines, as for $\mathcal{X}_{\tempe}(G)$, a structure of topological space on $\mathcal{X}(G)$. We also define a regular Borel measure $d\pi$ on $\mathcal{X}(G)$ by requesting that

$$\displaystyle \int_{\mathcal{X}(G)} \varphi(\pi) d\pi=\sum_{M\in \mathcal{M}} \left\lvert W(G,M)\right\rvert^{-1}\sum_{\mathcal{O}\in \mathcal{X}_{\elli}(M)/i\mathcal{A}_{M,F}^*} [i\mathcal{A}^\vee_{M,\sigma}:i\mathcal{A}^\vee_{M,F}]^{-1} \int_{i\mathcal{A}^*_{M,F}} \varphi(i_M^G(\sigma_\lambda))d\lambda$$

\noindent for all $\varphi\in C_c(\mathcal{X}_{\elli}(G))$, where we have fixed a base point $\sigma\in \mathcal{O}$ for every orbit $\mathcal{O}\in \mathcal{X}_{\elli}(M)/i\mathcal{A}_{M,F}^*$. We will need the following:

\vspace{3mm}

\begin{num}
\item\label{eq 2.7.2} If $F=\mathbb{R}$, there exists an integer $k\geqslant 0$ such that the integral
$$\int_{\mathcal{X}(G)} N^G(\pi)^{-k} d\pi$$
is absolutely convergent.
\end{num}

\vspace{3mm}

\noindent Finally, we extend the function $\pi\mapsto D(\pi)$ to $\mathcal{X}(G)$ by setting $D(\pi)=D(\sigma)$ for $\pi=i_M^G(\sigma)$, where $M$ is a Levi subgroup and $\sigma\in \mathcal{X}_{\elli}(M)$.

\section{Harish-Chandra descent}\label{section 3}

In this chapter, we collect some well-known facts concerning Harish-Chandra's technique of descent which are scattered over the literature. In this paper, we will use three types of descent. First, there is the semi-simple descent for the group or its Lie algebra which allows to localize functions or distributions near a semi-simple conjugacy class. This is the object of Section \ref{section 3.2}. Then, there is the descent from the group to its Lie algebra which is covered in Section \ref{section 3.3}. Finally, in Section \ref{section 3.4} we will discuss Harish-Chandra's notion of parabolic descent or more precisely its dual form which allows to parabolically induce invariant distributions. In Section \ref{section 3.1}, we set up notations and record basic properties of spaces of invariant functions/distributions and algebras of differential operators acting on those.

\subsection{Invariant analysis}\label{section 3.1}

\noindent Let $\omega\subseteq \mathfrak{g}(F)$ (resp.\ $\Omega\subseteq G(F)$) be a completely $G(F)$-invariant open subset (see Section \ref{section 1.8} for this notion). We will denote by $\gls{ComegaG}$ (resp.\ $\gls{COmegaG}$) the space of smooth and $G(F)$-invariant functions on $\omega$ (resp.\ $\Omega$). It is a closed subspace of $C^\infty(\omega)$ (resp.\ $C^\infty(\Omega)$) and we endow it with the induced locally convex topology. The notation $\gls{D'omegaG}$ (resp.\ $\gls{D'OmegaG}$) will stand for the space of $G(F)$-invariant distributions on $\omega$ (resp.\ $\Omega$). We will also denote by $\gls{Somega}$ (resp.\ $\gls{SOmega}$) the space of all functions $f\in \mathcal{S}(\mathfrak{g}(F))$ (resp.\ $f\in \mathcal{S}(G(F))$) such that $\overline{\Supp(f)^G}\subseteq \omega$ (resp.\ $\overline{\Supp(f)^G}\subseteq \Omega$). 

\vspace{2mm}

\noindent Assume now that $F=\mathbb{R}$. For each integer $n\geqslant 0$, we will denote by $\gls{DiffnomegaG}$ (resp.\ $\gls{DiffnOmegaG}$) the space of smooth invariant differential operators on $\omega$ (resp.\ on $\Omega$) that are of order less than $n$. It is a closed subspace of $\Diff^\infty_{\leqslant n}(\omega)$ (resp.\ of $\Diff^\infty_{\leqslant n}(\Omega)$) and we endow it with the induced locally convex topology. We will set

$$\displaystyle \gls{DiffomegaG}=\bigcup_{n\geqslant 0} \Diff_{\leqslant n}^\infty(\omega)^G \; (\mbox{resp.\ } \gls{DiffOmegaG}=\bigcup_{n\geqslant 0} \Diff_{\leqslant n}^\infty(\Omega)^G)$$

\noindent and we equip this space with the direct limit topology. We define

$$\displaystyle \gls{Jinfomega}=\{D\in \Diff^\infty(\omega)^G; \; DT=0 \; \forall T\in \mathcal{D}'(\omega)^G\}$$

$$\displaystyle \bigg(\mbox{resp.\ }\gls{JinfOmega}=\{D\in \Diff^\infty(\Omega)^G; \; DT=0 \; \forall T\in \mathcal{D}'(\Omega)^G\}\bigg)$$

\noindent and set

$$\displaystyle \gls{DiffomegaGbar}=\Diff^\infty(\omega)^G/\mathcal{J}^\infty(\omega) \; \left(\mbox{resp.\ }\gls{DiffOmegaGbar}=\Diff^\infty(\Omega)^G/\mathcal{J}^\infty(\Omega)\right)$$

\noindent Note that we have a natural action of $\overline{\Diff^\infty(\omega)^G}$ (resp.\ $\overline{\Diff^\infty(\Omega)^G}$) on $\mathcal{D}'(\omega)^G$ (resp.\ $\mathcal{D}'(\Omega)^G$).

\vspace{2mm}

\noindent For all $n\geqslant 0$, we will also denote by $\gls{DiffngG}$ the space of invariant differential operators with polynomial coefficients on $\mathfrak{g}(F)$ of order less than $n$ and we will set

$$\displaystyle \gls{DiffgG}=\bigcup_{n\geqslant 0} \Diff_{\leqslant n}(\mathfrak{g})^G$$

\noindent Note that $I(\mathfrak{g})$ and $I(\mathfrak{g}^*)$ are both naturally subalgebras of $\Diff(\mathfrak{g})^G$. We define

$$\displaystyle \gls{Jg}=\{D\in \Diff(\mathfrak{g})^G; \; DT=0 \; \forall T\in \mathcal{D}'(\mathfrak{g}(F))^G\}$$

\noindent and we will denote by $\gls{DiffgGbar}$ the quotient $\Diff(\mathfrak{g})^G/\mathcal{J}(\mathfrak{g})$.

\vspace{2mm}

\begin{prop}\label{proposition 3.1.1}
\begin{enumerate}[(i)]
\item Let $\omega\subseteq \mathfrak{g}(F)$ (resp.\ $\Omega\subseteq G(F)$) be a completely $G(F)$-invariant open subset. Then, there exists a sequence $(\omega_n)_{n\geqslant 1}$ (resp.\ $(\Omega_n)_{n\geqslant 1}$) of completely $G(F)$-invariant open subsets of $\omega$ (resp.\ $\Omega$) such that

\begin{itemize}
\renewcommand{\labelitemi}{$\bullet$}
\item $\displaystyle \omega=\bigcup_{n\geqslant 1} \omega_n$ (resp.\ $\displaystyle \Omega=\bigcup_{n\geqslant 1} \Omega_n$);

\item For all $n\geqslant 1$, $\overline{\omega}_n$ (resp.\ $\overline{\Omega}_n$) is compact modulo conjugation and included in $\omega_{n+1}$ (resp.\ $\Omega_{n+1}$).
\end{itemize}

\item Let $(\omega_i)_{i\in I}$ (resp.\ $(\Omega_i)_{i\in I}$) be a family of completely $G(F)$-invariant open subsets  and $L\subseteq \mathfrak{g}(F)$ (resp.\ $L\subseteq G(F)$) an invariant compact modulo conjugation subset such that

$$\displaystyle L\subseteq \bigcup_{i\in I}\omega_i\;\; \left(\mbox{resp.\ } L\subseteq \bigcup_{i\in I}\Omega_i\right)$$

\noindent Then, there exist a finite subset $J\subseteq I$ and functions $\varphi_j\in C^\infty(\omega_j)^G$ (resp.\ $\varphi_j\in C^\infty(\Omega_j)^G$) such that

\begin{itemize}
\renewcommand{\labelitemi}{$\bullet$}
\item For all $j\in J$, $0\leqslant \varphi_j\leqslant 1$ and $\Supp_{\omega_j}(\varphi_j)$ (resp.\ $\Supp_{\Omega_j}(\varphi_j)$) is compact modulo conjugation.

\item $\displaystyle \sum_{j\in J}\varphi_j=1$ on some invariant neighborhood of $L$.
\end{itemize}

\item Let $L\subseteq \mathfrak{g}(F)$ (resp.\ $L\subseteq G(F)$) be invariant and compact modulo conjugation. Then, there exist two constants $c>0$ and $m\geqslant 1$ and a compact subset $\mathcal{K}\subseteq \mathfrak{g}(F)$ (resp.\ $\mathcal{K}\subseteq G(F)$) such that for all $X\in L$ (resp.\ $g\in L$) there exists an element $\gamma\in G(F)$ satisfying the two following conditions

\begin{itemize}
\renewcommand{\labelitemi}{$\bullet$}
\item $\lVert \gamma\rVert\leqslant c\lVert X\rVert^m$ (resp.\ $\lVert \gamma\rVert\leqslant c\lVert g\rVert^m$);

\item $\gamma^{-1}X\gamma\in \mathcal{K}$ (resp.\ $\gamma^{-1}g\gamma\in \mathcal{K}$).
\end{itemize}

\item Let $\varphi\in C^\infty(\mathfrak{g}(F))^G$ (resp.\ $\varphi\in C^\infty(G(F))^G$) be compactly supported modulo conjugation. Then, multiplication by $\varphi$ preserves $\mathcal{S}(\mathfrak{g}(F))$ (resp.\ preserves $\mathcal{S}(G(F))$), that is: for all $f\in \mathcal{S}(\mathfrak{g}(F))$ (resp.\ for all $f\in \mathcal{S}(G(F))$), we have $\varphi f\in \mathcal{S}(\mathfrak{g}(F))$ (resp.\ $\varphi f\in \mathcal{S}(G(F))$).

\item Assume that $F=\mathbb{R}$ and let $\omega\subseteq \mathfrak{g}(F)$ be a completely $G(F)$-invariant open subset. Then, for every integer $n\geqslant 0$, there exists a finite family $\{D_1,\ldots,D_k\}\subseteq \Diff_{\leqslant n}(\omega)^G$ such that the linear map

$$\displaystyle \left(C^\infty(\omega)^G\right)^k\to \Diff^\infty_{\leqslant n}(\omega)^G$$
$$(\varphi_1,\ldots,\varphi_k)\mapsto \varphi_1D_1+\ldots+\varphi_kD_k$$

\noindent is a topological isomorphism (in particular the $C^\infty(\omega)^G$-module $\Diff^\infty_{\leqslant n}(\omega)^G$ is free of finite rank).

\item Still assuming that $F=\mathbb{R}$. the $\mathbb{C}$-algebra $\overline{\Diff(\mathfrak{g})^G}$ is generated by the image of $I(\mathfrak{g})$ and $I(\mathfrak{g}^*)$.
\end{enumerate}
\end{prop}

\vspace{2mm}

\noindent\ul{Proof}: (i) is contained in Lemma 2.2.1 of \cite{Bou2} and Lemma 2.2.2 of \cite{Bou1} whereas (ii) follows from Lemma 2.3.1 of \cite{Bou2} and Lemma 2.3.1 of \cite{Bou1}.

\begin{enumerate}[(i)]
\setcounter{enumi}{2}
\item We prove it for the group the proof for the Lie algebra being similar and easier. Let $P_{\mini}=M_{\mini}U_{\mini}$ be a minimal parabolic subgroup of $G$ and $A_{\mini}=A_{M_{\mini}}$ be the split part of the center of $M_{\mini}$. Set

$$\displaystyle A_{\mini}^+=\{a\in A_{\mini}^+;\; \lvert\alpha(a)\rvert\leqslant 1\; \forall \alpha\in R(A_{\mini},P_{\mini})\}$$

\noindent Then by the Cartan decomposition, there exists a compact subset $C_G\subseteq G(F)$ such that

$$G(F)=C_G A_{\mini}^+C_G$$

\noindent Let $L\subseteq G(F)$ be invariant and compact modulo conjugation and fix a compact subset $\mathcal{K}_G\subseteq G(F)$ such that $L=(\mathcal{K}_G)^G$. Replacing $\mathcal{K}_G$ by $(\mathcal{K}_G)^{C_G}$ if necessary, we may assume that $L=(\mathcal{K}_G)^{A_{\mini}^+C_G}$. Since $C_G$ is compact, it suffices to show the following

\vspace{3mm}

\begin{num}
\item\label{eq 3.1.1} There exists a constant $c>0$ such that for all $g\in \mathcal{K}_G$ and all $a\in A_{\mini}^+$, there exists $a'\in A_{\mini}^+$ satisfying the three following conditions
$${a'}^{-1}a\in A_{\mini}^+,\;\; \sigma(a')\leqslant c\sigma(a^{-1}ga) \mbox{ and } \sigma (a'a^{-1}ga{a'}^{-1})\leqslant c$$
\end{num}

\vspace{3mm}

\noindent We prove this by induction on $\dim(G)$, the case where $G$ is a torus being trivial. Denote by $\overline{P}_{\mini}=M_{\mini}\overline{U}_{\mini}$ the parabolic subgroup opposite to $P_{\mini}$ with respect to $M_{\mini}$. Choose $\delta>0$ and set

$$\displaystyle A_{\mini}^{\overline{Q},+}(\delta)=\{a\in A_{\mini}^+;\; \lvert \alpha(a)\rvert\geqslant e^{\delta\sigma(a)} \; \forall \alpha\in R(A_{\mini},U_{\overline{Q}})\}$$

\noindent for every parabolic subgroup $\overline{Q}=M_{\overline{Q}}U_{\overline{Q}}\supseteq \overline{P}_{\mini}$. If $\delta$ is chosen sufficiently small, and we will assume that it is so in what follows, the complement of

$$\displaystyle \bigcup_{\overline{P}_{\mini}\subseteq\overline{Q}\neq G} A_{\mini}^{\overline{Q},+}(\delta)$$

\noindent in $A_{\mini}^+$ is compact. Hence, to get \ref{eq 3.1.1} it suffices to prove the following for every parabolic subgroup $\overline{P}_{\mini}\subseteq \overline{Q}\neq G$

\vspace{3mm}

\begin{num}
\item\label{eq 3.1.2} There exists $c=c_{\overline{Q}}>0$ such that for all $g\in \mathcal{K}_G$ and all $a\in A_{\mini}^{\overline{Q},+}(\delta)$, there exists $a'\in A_{\mini}^+$ satisfying the three following conditions
$${a'}^{-1}a\in A_{\mini}^+,\;\; \sigma(a')\leqslant c\sigma(a^{-1}ga) \mbox{ and } \sigma (a'a^{-1}ga{a'}^{-1})\leqslant c$$
\end{num}

\vspace{3mm}

\noindent Fix such a parabolic subgroup $\overline{Q}=MU_{\overline{Q}}$, where $M$ is the only Levi component of $\overline{Q}$ containing $M_{\mini}$, and let $Q=MU_Q$ be the parabolic subgroup opposite to $\overline{Q}$ with respect to $M$. Fix also $\epsilon>0$, that we will assume sufficiently small in what follows, and set

$$\displaystyle \mathcal{K}^{\overline{Q}}_{\epsilon,a}=\mathcal{K}_G\cap \left(U_{\overline{Q}}[<\epsilon \sigma(a)] M[<\epsilon \sigma(a)] aU_{Q}[<\epsilon \sigma(a)]a^{-1}\right)$$

\noindent for all $a\in A_{\mini}(F)$. By Lemma \ref{proposition 1.3.1}(i) and since $\mathcal{K}_G$ is compact, there exists $c_0>0$ such that $\sigma(a)\leqslant c_0\sigma(a^{-1}ga)$ for all $a\in A_{\mini}^{\overline{Q},+}(\delta)$ and all $g\in \mathcal{K}_G\backslash \mathcal{K}^{\overline{Q}}_{\epsilon,a}$. Hence, we only need to prove the existence of $c\geqslant c_0$ such that \ref{eq 3.1.2} holds for all $a\in A_{\mini}^{\overline{Q},+}(\delta)$ and all $g\in \mathcal{K}^{\overline{Q}}_{\epsilon,a}$ (otherwise, we just take $a'=a$). Choosing $\epsilon>0$ sufficiently small, we may assume that the subsets $aU_{Q}[<\epsilon \sigma(a)]a^{-1}$ remain uniformly bounded as $a$ varies in $A_{\mini}^{\overline{Q},+}(\delta)$. Then, there exists compact subsets $\mathcal{K}_{\overline{U}}\subseteq U_{\overline{Q}}(F)$, $\mathcal{K}_M\subseteq M(F)$ and $\mathcal{K}_U\subseteq U_{Q}(F)$ such that

$$\displaystyle \mathcal{K}^{\overline{Q}}_{\epsilon,a} \subseteq \mathcal{K}_{\overline{U}}\mathcal{K}_M\mathcal{K}_U$$

\noindent for all $a\in A_{\mini}^{\overline{Q},+}(\delta)$. Since the subsets $a^{-1}\mathcal{K}_{\overline{U}}a$ remain uniformly bounded for $a\in A_{\mini}^+$, to get \ref{eq 3.1.2} we only need to show the following

\vspace{3mm}

\begin{num}
\item\label{eq 3.1.3} There exists $c>0$ such that for all $m\in \mathcal{K}_M$, all $u\in \mathcal{K}_U$ and all $a\in A_{\mini}^+$, there exists $a'\in A_{\mini}^+$ satisfying the three following conditions
$${a'}^{-1}a\in A_{\mini}^+,\;\; \sigma(a')\leqslant c\sigma(a^{-1}mua)\;\mbox{ and } \sigma(a'a^{-1}mua{a'}^{-1})\leqslant c$$
\end{num}

\vspace{3mm}

\noindent Since $\sigma(mu)\sim \sigma(m)+\sigma(u)$ for all $m\in M(F)$ and all $u\in U_Q(F)$ and $\sigma(aua^{-1})\ll \sigma(u)$ for all $u\in U_Q(F)$ and all $a\in A_{\mini}^+$, the last claim will follow from the combination of the two next facts

\vspace{3mm}

\begin{num}
\item\label{eq 3.1.4} There exists $c_U>0$ such that for all $u\in \mathcal{K}_U$ and all $a\in A_{\mini}^+$, there exists $a'\in A_{\mini}^+$ satisfying the three following conditions
$$\displaystyle {a'}^{-1}a\in A_{\mini}^+,\;\; \sigma(a')\leqslant c_U\sigma(a^{-1}ua)\;\mbox{ and } \sigma(a'a^{-1}ua{a'}^{-1})\leqslant c_U$$

\item\label{eq 3.1.5} There exists $c_M>0$ such that for all $m\in \mathcal{K}_M$ and all $a\in A_{\mini}^+$, there exists $a'\in A_{\mini}^+$ satisfying the three following conditions
$$\displaystyle {a'}^{-1}a\in A_{\mini}^+,\;\; \sigma(a')\leqslant c_M\sigma(a^{-1}ma)\;\mbox{ and } \sigma(a'a^{-1}ma{a'}^{-1})\leqslant c_M$$
\end{num}

\vspace{3mm}

\noindent The proof of \ref{eq 3.1.4} is easy and left to the reader. The point \ref{eq 3.1.5} follows from the induction hypothesis applied to $M$. Indeed, by the induction hypothesis there exists $c_M>0$ such that for all $m\in \mathcal{K}_M$ and all $a\in A_{\mini}^+$, there exists $a'\in A_{\mini}^{M,+}$ such that

\begin{align}\label{eq 3.1.6}
\displaystyle {a'}^{-1}a\in A_{\mini}^{M,+},\;\; \sigma(a')\leqslant c_M\sigma(a^{-1}ma) \mbox{ and } \sigma(a'a^{-1}ma{a'}^{-1})\leqslant c_M
\end{align}

\noindent where we have set

$$\displaystyle A_{\mini}^{M,+}=\{a\in A_{\mini}(F);\; \lvert \alpha(a)\rvert \leqslant 1\; \forall \alpha\in R(A_{\mini},M\cap P_{\mini})\}$$

\noindent Denote by $\Delta\subseteq R(A_{\mini},P_{\mini})$ the subset of simple roots and set $\Delta_Q=\Delta\cap R(A_{\mini},U_Q)$. Let $A_M^G=\{a\in A_M;\; \chi(a)=1\; \forall \chi\in X^*(G)\}$ and $A_{\mini}^{\Delta_Q}=\{a\in A_{\mini};\; \alpha(a)=1\; \forall \alpha\in \Delta_Q\}$. The multiplication map $A_M^G\times A_{\mini}^{\Delta_Q}\to A_{\mini}$ is an isogeny. Hence, since $A_M^G(F)$ is in the center of $M(F)$, up to increasing the constant $c_M$ we see that for all $m\in \mathcal{K}_M$ and all $a\in A_{\mini}^+$, we can find $a'\in A_{\mini}^{\Delta_Q,+}=A_{\mini}^{\Delta_Q}(F)\cap A_{\mini}^{M,+}$ satisfying \ref{eq 3.1.6}. But obviously $A_{\mini}^{\Delta_Q,+}\subseteq A_{\mini}^+$ and for all $a\in A_{\mini}^+$ and $a'\in A_{\mini}^{\Delta_Q,+}$ the first condition of \ref{eq 3.1.6} is equivalent to ${a'}^{-1}a\in A_{\mini}^+$. This proves \ref{eq 3.1.5} and ends the proof of (iii).

\item This is clear in the $p$-adic case. Assume that $F=\mathbb{R}$. Then the result will follow at once from the following fact which is an easy consequence of (iii)

\vspace{3mm}

\begin{num}
\item\label{eq 3.1.7} For all $u\in I(\mathfrak{g})$ (resp.\ $u\in \mathcal{U}(\mathfrak{g})$), there exists $N\geqslant 1$ such that we have
$$\displaystyle \left\lvert \left(\partial(u)\varphi\right)(X)\right\rvert\ll \lVert X\rVert^N\; (\mbox{resp.\ } \left\lvert \left(R(u)\varphi\right)(g)\right\rvert\ll \lVert g\rVert^N)$$
for all $X\in \mathfrak{g}(F)$ (resp.\ for all $g\in G(F)$).
\end{num}

\vspace{3mm}

\item By Corollaire 3.7 of \cite{Bou3}, the $C^\infty(\omega)^G$-module $\Diff^\infty_{\leqslant n}(\omega)^G$ is free of finite rank and we can find a basis consisting of elements of $\Diff_{\leqslant n}(\mathfrak{g})^G$. Fix such a basis $(D_1,\ldots,D_k)$ then the linear map

$$\displaystyle \left( C^\infty(\omega)^G\right)^k\to \Diff^\infty_{\leqslant n}(\omega)^G$$
$$\displaystyle (\varphi_1,\ldots,\varphi_k)\mapsto \varphi_1D_1+\ldots+\varphi_kD_k$$

\noindent is continuous and bijective. Since both $\left( C^\infty(\omega)^G\right)^k$ and $\Diff^\infty_{\leqslant n}(\omega)^G$ are Fr\'echet spaces, by the open mapping theorem, this is a topological isomorphism.

\item This follows from Theorem 1, Theorem 3 and Theorem 5 of \cite{LS} $\blacksquare$

\end{enumerate}

\subsection{Semi-simple descent}\label{section 3.2}

\noindent Let $X\in \mathfrak{g}_{\ssi}(F)$. An open subset $\omega_X\subseteq \mathfrak{g}_X(F)$ will be called {\em $G$-good} if it is completely $G_X(F)$-invariant and if moreover the map

\begin{align}\label{eq 3.2.1}
\displaystyle \omega_X\times^{G_X(F)} G(F)\to \mathfrak{g}(F)
\end{align}
$$\displaystyle (Y,g)\mapsto g^{-1}Yg$$

\noindent induces an $F$-analytic isomorphism between $\omega_X\times^{G_X(F)} G(F)$ and $\omega_X^G$, where $\omega_X\times^{G_X(F)} G(F)$ (the {\em contracted product}) denotes the quotient of $\omega_X\times G(F)$ by the free $G_X(F)$-action given by

$$\displaystyle g_X\cdot(Y,g)=(g_XYg_X^{-1},g_Xg),\;\;\; g_X\in G_X(F), (Y,g)\in \omega_X\times G(F)$$

\noindent The Jacobian of the map \ref{eq 3.2.1} at $(Y,g)\in \omega_X\times^{G_X(F)} G(F)$ is equal to

$$\displaystyle \gls{etaXG}(Y)=\left\lvert \det \ad(Y)_{\mid \mathfrak{g}/\mathfrak{g}_X}\right\rvert$$

\noindent It follows that an open subset $\omega_X\subseteq \mathfrak{g}_X(F)$ is $G$-good if and only if the following conditions are satisfied

\vspace{2mm}

\begin{itemize}
\renewcommand{\labelitemi}{$\bullet$}
\item $\omega_X$ is completely $G_X(F)$-invariant;
\item For all $Y\in \omega_X$, we have $\eta^G_X(Y)\neq 0$
\item For all $g\in G(F)$, the intersection $g^{-1}\omega_Xg\cap \omega_X$ is nonempty if and only if $g\in G_X(F)$.
\end{itemize}

\vspace{2mm}

\noindent Let $\omega_X\subseteq \mathfrak{g}_X(F)$ be a $G$-good open neighborhood of $X$ and set $\omega=\omega_X^G$. Then $\omega$ is completely $G(F)$-invariant (since $\omega_X$ is completely $G_X(F)$-invariant). Moreover, the completely $G(F)$-invariant open subsets obtained in this way form a basis of neighborhood for $X$ in the invariant topology. We have the integration formula

\begin{align}\label{eq 3.2.2}
\displaystyle\int_{\omega}f(Y)dY=\int_{G_X(F)\backslash G(F)} \int_{\omega_X}f(g^{-1}Yg)\eta^G_X(Y) dYdg
\end{align}

\noindent for all $f\in L^1(\omega)$. For every function $f$ defined on $\omega$, we will denote by $\gls{fXomegaX}$ the function on $\omega_X$ given by $f_{X,\omega_X}(Y)=\eta_X^G(Y)^{1/2}f(Y)$. The map $f\mapsto f_{X,\omega_X}$ induces topological isomorphisms

$$\displaystyle C^\infty(\omega)^G\simeq C^\infty(\omega_X)^{G_X}\;\;\; C^\infty(\omega_{\reg})^G\simeq C^\infty(\omega_{X,\reg})^{G_X}$$

\noindent (Note that $\omega_X\cap \mathfrak{g}_{X,\reg}(F)=\omega_X\cap\mathfrak{g}_{\reg}(F)$ so that the notation $\omega_{X,\reg}$ is unambiguous). We also have an isomorphism 

$$\displaystyle \mathcal{D}'(\omega)^G\simeq \mathcal{D}'(\omega_X)^{G_X}$$
$$\displaystyle T\mapsto \gls{TXomegaX}$$

\noindent where for $T\in\mathcal{D}'(\omega)^G$, $T_{X,\omega_X}$ is the unique $G_X(F)$-invariant distribution on $\omega_X$ such that

$$\displaystyle \langle T,f\rangle=\int_{G_X(F)\backslash G(F)} \langle T_{X,\omega_X},({}^g f)_{X,\omega_X}\rangle dg$$

\noindent for all $f\in C_c^\infty(\omega)$. If $f$ is a locally integrable and invariant function on $\omega$, then by the integration formula \ref{eq 3.2.2}, $f_{X,\omega_X}$ is also locally integrable and we have

$$\displaystyle \left(T_f\right)_{X,\omega_X}=T_{f_{X,\omega_X}}$$

\vspace{2mm}

\noindent Let $x\in G_{\ssi}(F)$. We define similarly the notion of {\em $G$-good} open subset of $G_x(F)$. More precisely, an open subset $\Omega_x\subseteq G_x(F)$ is $G$-good if it is completely $Z_G(x)(F)$-invariant and if moreover the map

\begin{align}\label{eq 3.2.3}
\displaystyle \Omega_x\times^{Z_G(x)(F)} G(F)\to G(F)
\end{align}
$$\displaystyle (y,g)\mapsto g^{-1}yg$$

\noindent induces an $F$-analytic isomorphism $\Omega_x\times^{Z_G(x)(F)} G(F)\simeq \Omega_x^G$, where this time $\Omega_x\times^{Z_G(x)(F)} G(F)$ denotes the quotient of $\Omega_x\times G(F)$ by the free $Z_G(x)(F)$ action given by

$$\displaystyle g_x\cdot(y,g)=(g_xyg_x^{-1},g_xg),\;\;\; g_x\in Z_G(x)(F), (y,g)\in \Omega_x\times G(F)$$

\noindent The Jacobian of the map \ref{eq 3.2.3} at $(y,g)$ is given by

\begin{align}\label{eq 3.2.4}
\displaystyle \gls{etaxG}(y)=\left\lvert \det\left(1-\Ad(y)\right)_{\mid \mathfrak{g}/\mathfrak{g}_x}\right\rvert
\end{align}

\noindent Thus, an open subset $\Omega_x\subseteq G_x(F)$ is $G$-good if and only if the following conditions are satisfied

\vspace{2mm}

\begin{itemize}
\renewcommand{\labelitemi}{$\bullet$}
\item $\Omega_x$ is completely $Z_G(x)(F)$-invariant;
\item For all $y\in \Omega_x$, we have $\eta^G_x(y)\neq 0$
\item For all $g\in G(F)$, the intersection $g^{-1}\Omega_xg\cap \Omega_x$ is nonempty if and only if $g\in Z_G(x)(F)$.
\end{itemize}

\vspace{2mm}

\noindent Let $\Omega_x\subseteq G_x(F)$ be a $G$-good open subset and set $\Omega=\Omega_x^G$. We have the integration formula

\begin{align}\label{eq 3.2.5}
\displaystyle \displaystyle \int_{\Omega} f(y) dy & =\int_{Z_G(x)(F)\backslash G(F)}\int_{\Omega_x} f(g^{-1}yg)\eta^G_x(y)dydg  \\
\nonumber &  =[Z_G(x)(F):G_x(F)]^{-1}\int_{G_x(F)\backslash G(F)}\int_{\Omega_x} f(g^{-1}yg)\eta^G_x(y)dydg
\end{align}

\noindent for all $f\in L^1(\Omega)$. For every function $f$ on $\Omega$, we will denote by $\gls{fxOmegax}$ the function on $\Omega_x$ defined by $f_{x,\Omega_x}(y)=\eta^G_x(y)^{1/2} f(y)$. Again, the map $f\mapsto f_{x,\Omega_x}$ induces topological isomorphisms

$$\displaystyle C^\infty(\Omega)^G\simeq C^\infty(\Omega_x)^{Z_G(x)}\;\;\; C^\infty(\Omega_{\reg})^G\simeq C^\infty(\Omega_{x,\reg})^{Z_G(x)}$$

\noindent which extend to an isomorphism

$$\displaystyle \mathcal{D}'(\Omega)\simeq \mathcal{D}'(\Omega_x)^{Z_G(x)}$$

\noindent Finally if $F=\mathbb{R}$, we may also descend invariant differential operators. The result is summarized in the next lemma. We just need to introduce first some notation. Recall that for all $X\in\mathfrak{g}_{\ssi}(F)$ and all $x\in G_{\ssi}(F)$, the Harish-Chandra isomorphisms induces injective $\mathbb{C}$-algebra homomorphisms

$$\displaystyle I(\mathfrak{g}^*)\to I(\mathfrak{g}_X^*)\;\;\; I(\mathfrak{g})\to I(\mathfrak{g}_X)\;\;\; \mathcal{Z}(\mathfrak{g})\to \mathcal{Z}(\mathfrak{g}_x)$$
$$\displaystyle \;\; p\mapsto p_{G_X}\;\;\;\;\;\;\;\;\;\; u\mapsto u_{G_X}\;\;\;\;\;\;\;\;\;\;\;\; z\mapsto z_{G_x}$$

\noindent (cf.\ Section \ref{section 1.1}). We shall denote these homomorphisms simply by $p\mapsto \gls{pX}$, $u\mapsto \gls{uX}$ and $z\mapsto \gls{zx}$ respectively. Note that the image of $z\mapsto z_x$ is included in $\mathcal{Z}(\mathfrak{g}_x)^{Z_G(x)}$.

\vspace{3mm}

\begin{lem}\label{lemma 3.2.1}
Assume that $F=\mathbb{R}$. 

\begin{enumerate}[(i)]
\item Let $X\in \mathfrak{g}_{\ssi}(F)$ and let $\omega_X\subseteq \mathfrak{g}_X(F)$ be a $G$-good open subset. Set $\omega=\omega_X^G$. Then, there exists a unique topological isomorphism

$$\displaystyle \overline{\Diff^\infty(\omega)^G}\simeq \overline{\Diff^\infty(\omega_X)^{G_X}}$$
$$\displaystyle D\mapsto \gls{DXomegaX}$$

\noindent such that for all $T\in \mathcal{D}'(\omega)^G$ and all $D\in \Diff^\infty(\omega)^G$, we have

\begin{align}\label{eq 3.2.6}
\displaystyle \left(DT\right)_{X,\omega_X}=D_{X,\omega_X}T_{X,\omega_X}
\end{align}

\noindent Moreover, we have

\begin{align}\label{eq 3.2.7}
\displaystyle \left(\partial(u)\right)_{X,\omega_X}=\partial(u_X)\;\;\; \mbox{ and }\;\;\; (p)_{X,\omega_X}=p_X
\end{align}

\noindent for all $u\in I(\mathfrak{g})$ and all $p\in I(\mathfrak{g}^*)$. In particular, by Proposition \ref{proposition 3.1.1}(vi), the image of $\overline{\Diff(\mathfrak{g})^G}$ by the map $D\mapsto D_{X,\omega_X}$ lies in $\overline{\Diff(\mathfrak{g}_X)^{G_X}}$.

\item Let $x\in G_{\ssi}(F)$ and let $\Omega_x\subseteq G_x(F)$ be a $G$-good open subset. Set $\Omega=\Omega_x^G$. Then, there exists a unique topological isomorphism

$$\displaystyle \overline{\Diff^\infty(\Omega)^G}\simeq \overline{\Diff^\infty(\Omega_x)^{Z_G(x)}}$$
$$\displaystyle D\mapsto \gls{DxOmegax}$$

\noindent such that for all $T\in \mathcal{D}'(\Omega)^G$ and all $D\in \Diff^\infty(\Omega)^G$, we have

\begin{align}\label{eq 3.2.8}
\displaystyle \left(DT\right)_{x,\Omega_x}=D_{x,\Omega_x}T_{x,\Omega_x}
\end{align}

\noindent Moreover, we have

\begin{align}\label{eq 3.2.9}
\displaystyle (z)_{x,\Omega_x}=z_x
\end{align}

\noindent for all $z\in \mathcal{Z}(\mathfrak{g})$.

\end{enumerate}

\end{lem}

\vspace{2mm}

\noindent\ul{Proof}:
\begin{enumerate}[(i)]
\item Since $T\mapsto T_{X,\omega_X}$ is an isomorphism $\mathcal{D}'(\omega)^G\simeq \mathcal{D}'(\omega_X)^{G_X}$, for all $D\in \Diff^\infty(\omega)^G$ there exists at most one operator $D_{X,\omega_x}\in \overline{\Diff^\infty(\omega_X)^{G_X}}$ such that the relation \ref{eq 3.2.6} is satisfied for all $T\in \mathcal{D}'(\omega)^G$. Moreover, such an operator is constructed in Theorem 11 p.30 of \cite{Va} (although in this reference only analytic differential operators are considered, the construction applies equally well to smooth differential operators). This yields an injective linear map

$$\displaystyle \overline{\Diff^\infty(\omega)^G}\to \overline{\Diff^\infty(\omega_X)^{G_X}}$$
$$\displaystyle D\mapsto D_{X,\omega_X}$$

\noindent and we need to prove that it is a topological isomorphism. Since both $\overline{\Diff^\infty(\omega)^G}$ and $\overline{\Diff^\infty(\omega_X)^{G_X}}$ are LF spaces, by the open mapping theorem, we only need to construct a continuous right inverse to the previous linear map. Actually, we are going to construct a continuous linear map

$$\Diff^\infty(\omega_X)^{G_X}\to \Diff^\infty(\omega)^G$$
$$D\mapsto D_\omega$$

\noindent such that

$$\left(\overline{D_\omega}\right)_{X,\omega_X}=\overline{D}$$

\noindent for all $D\in \Diff^\infty(\omega_X)^{G_X}$, where we have denoted by $\overline{D}$ and $\overline{D_\omega}$ the image of $D$ and $D_\omega$ in $\overline{\Diff^\infty(\omega_X)^{G_X}}$ and $\overline{\Diff^\infty(\omega)^{G}}$ respectively. The construction is as follows. A smooth differential operator $D$ on $\omega_X$ (resp.\ on $\omega$) may be seen as map $Y\in\omega_X\mapsto D_Y\in S(\mathfrak{g}_X)$ (resp.\ $Y\in\omega\mapsto D_Y\in S(\mathfrak{g})$) which has its image in a finite dimensional subspace and is smooth, the action of $D$ on smooth functions being given by

$$(Df)(Y)=(\partial(D_Y)f)(Y)$$

\noindent for all $f\in C^\infty(\omega_X)$ (resp.\ $f\in C^\infty(\omega)$) and all $Y\in \omega_X$ (resp.\ $Y\in\omega$). Let $D\in \Diff^\infty(\omega_X)^{G_X}$ be an invariant and smooth differential operator on $\omega_X$. We first associate to $D$ a differential operator $D^\natural_\omega$ on $\omega$ by setting

$$D^\natural_{\omega,Y}=g^{-1}D_{Y'}g\in S(\mathfrak{g})$$

\noindent for all $Y\in \omega$, where $g\in G(F)$ and $Y'\in \omega_X$ are any elements such that $Y=g^{-1}Y'g$. Since $D$ is invariant and $\omega_X$ is $G$-good, this definition doesn't depend on the choice of $g$ and $Y'$ and $D^\natural_\omega$ is a smooth invariant differential operator on $\omega$. The function $\eta_X^G\in C^\infty(\omega_X)^{G_X}$ uniquely extends to a smooth and invariant function on $\omega$. Still denoting by $\eta_X^G$ this extension, we now set

$$D_\omega=\left(\eta^G_X\right)^{-1/2}\circ D^\natural_\omega\circ \left(\eta_X^G\right)^{1/2}$$

\noindent for all $D\in \Diff^\infty(\omega_X)^G$. Then it is easy to see that the linear map $D\mapsto D_\omega$ has all the desired properties.

\vspace{2mm}

\noindent The second equality of \ref{eq 3.2.7} is obvious whereas the first one follows from Theorem 15 p.30 of \cite{Va}.

\item Once again, there exists at most one linear map

\begin{align}\label{eq 3.2.10}
\displaystyle \overline{\Diff^\infty(\Omega)^G}\to \overline{\Diff^\infty(\Omega_x)^{Z_G(x)}}
\end{align}
$$D\mapsto D_{x,\Omega_x}$$

\noindent such that the relation \ref{eq 3.2.8} is satisfied for all $D\in \Diff^\infty(\Omega)^G$ and all $T\in \mathcal{D}'(\Omega)^G$ and such a linear map, if it exists, is necessarily injective. The construction of Proposition 4 p.224 of \cite{Va} proves the existence of such a map, where once again the extension of the construction of that reference from the analytic case to the smooth case is straightforward. Moreover, we may construct explicitly, analogously to what has been done in the proof of (i), a right continuous inverse to \ref{eq 3.2.10}. Since both $\overline{\Diff^\infty(\Omega)^G}$ and $\overline{\Diff^\infty(\Omega_x)^{Z_G(x)}}$ are LF spaces, this proves by the open mapping theorem that \ref{eq 3.2.10} is a topological isomorphism. The equality \ref{eq 3.2.9} follows from Theorem 12 p.229 of \cite{Va}. $\blacksquare$
\end{enumerate}

\vspace{2mm}

\noindent Consider the particular case where $G_X=T$, $X\in\mathfrak{g}_{\ssi}(F)$, is a maximal torus. Then, if $F=\mathbb{R}$, the lemma provides us with a morphism $\Diff(\mathfrak{g})^G\to \Diff(\mathfrak{t})$ that we shall denote by $D\mapsto \gls{DT}$ in this particular case. For $\omega\subseteq \mathfrak{g}(F)$ a completely $G(F)$-invariant open subset and $f$ an invariant function on $\omega_{\reg}$, we will denote by $\gls{fT}$ the function on $\mathfrak{t}(F)\cap \omega_{\reg}$ given by

$$f_T(Y)=D^G(Y)^{1/2}f(Y),\;\;\; Y\in \mathfrak{t}(F)\cap \omega_{\reg}$$

\noindent Hence, we have

$$(Df)_T=D_Tf_T$$

\noindent for all $f\in C^\infty(\omega_{\reg})^G$ and all $D\in \Diff(\mathfrak{g})^G$. Note that an invariant function $f$ on $\omega_{\reg}$ is smooth if and only if $f_T$ is a smooth function for every maximal torus $T\subset G$.

\vspace{2mm}

\noindent\ul{Remark}: We can extend the definition of $G$-good open subsets to the case where $G$ is not necessarily reductive (but is still a connected linear algebraic group over $F$). The definition is as follows. Let $x\in G_{\ssi}(F)$, an open subset $\Omega_x\subseteq G_x(F)$ is {\em $G$-good} if it is invariant by translation by $G_{x,u}(F)$, where $G_{x,u}$ denotes the unipotent radical of $G_x$, and if moreover its image in $G_x(F)/G_{x,u}(F)=(G/G_u)_x(F)$ is a $G/G_u$-good open subset, where this time $G_u$ denotes the unipotent radical of $G$. If $\Omega_x\subseteq G_x(F)$ is a $G$-good open subset, then the map \ref{eq 3.2.3} still induces an $F$-analytic isomorphism onto $\Omega=\Omega_x^G$ whose Jacobian at $(y,g)\in \Omega_x\times^{Z_G(x)(F)}G(F)$ is again given by the formula \ref{eq 3.2.4}. It follows in particular that the integration formula \ref{eq 3.2.5} is still valid in this more general setting.

\subsection{Descent from the group to its Lie algebra}\label{section 3.3}

\noindent We will say of an open subset $\omega\subseteq \mathfrak{g}(F)$ that it is {\em $G$-excellent}, if it satisfies the following conditions

\vspace{2mm}

\begin{itemize}
\renewcommand{\labelitemi}{$\bullet$}

\item $\omega$ is completely $G(F)$-invariant and relatively compact modulo conjugation;

\item The exponential map is defined on $\omega$ and induces an $F$-analytic isomorphism between $\omega$ and $\Omega=\exp(\omega)$.
\end{itemize}

\vspace{2mm}

\noindent For all $X\in \mathfrak{z}_G(F)$ (in particular $X=0$), the $G$-excellent open subsets containing $X$ form a basis of neighborhoods of $X$ for the invariant topology.

\vspace{2mm}

\noindent Let $\omega\subseteq \mathfrak{g}(F)$ be a $G$-excellent open subset and set $\Omega=\exp(\omega)$. The Jacobian of the exponential map

$$\exp: \omega\to \Omega$$
$$X\mapsto e^X$$

\noindent at $X\in \omega_{\ssi}$ is given by

\begin{align}\label{eq 3.3.1}
\displaystyle \gls{jG}(X)=D^G(e^X)D^G(X)^{-1}
\end{align}

\noindent Hence, we have the integration formula

\begin{align}\label{eq 3.3.2}
\displaystyle \int_{\Omega} f(g) dg=\int_{\omega}f(e^X)j^G(X)dX
\end{align}

\noindent for all $f\in L^1(\Omega)$. For every function $f$ on $\Omega$, we will denote by $\gls{fomega}$ the function on $\omega$ defined by $f_\omega(X)=j^G(X)^{1/2}f(e^X)$. The map $f\mapsto f_\omega$ induces topological isomorphisms

$$\displaystyle C^\infty(\Omega)\simeq C^\infty(\omega)\;\;\;\;\; C^\infty(\Omega_{\reg})\simeq C^\infty(\omega_{\reg})$$

\noindent We will also denote by

$$\displaystyle \mathcal{D}'(\Omega)\simeq \mathcal{D}'(\omega)$$
$$\displaystyle T\mapsto \gls{Tomega}$$

\noindent the isomorphism defined by the relations

$$\displaystyle \langle T_\omega,f_\omega\rangle=\langle T,f\rangle$$

\noindent for all $T\in \mathcal{D}'(\Omega)$ and all $f\in C_c^\infty(\Omega)$. By \ref{eq 3.3.2}, if $f$ is a locally integrable function on $\Omega$, we have

$$(T_f)_\omega=T_{f_\omega}$$

\noindent There also exists a unique topological isomorphism

$$\displaystyle \Diff^\infty(\Omega)\simeq \Diff^\infty(\omega)$$
$$D\mapsto \gls{Domega}$$

\noindent such that $(DT)_\omega=D_\omega T_\omega$ for all $D\in \Diff^\infty(\Omega)$ and all $T\in \mathcal{D}'(\Omega)$. By Theorem 14 p.231 of \cite{Va}, we have

\begin{align}\label{eq 3.3.3}
\displaystyle (z)_\omega=\partial(u_z)
\end{align}

\noindent for all $z\in \mathcal{Z}(\mathfrak{g})$ (recall that we are denoting by $z\mapsto u_z$ the Harish-Chandra isomorphism $\mathcal{Z}(\mathfrak{g})\simeq I(\mathfrak{g})$).

\vspace{2mm}

\noindent We will denote by $f\mapsto \gls{fOmega}$, $T\mapsto \gls{TOmega}$ and $D\mapsto \gls{DOmega}$ the inverse of the previous isomorphisms. So for example $f_\Omega(g)=j^G(\log(g))^{-1/2}f(\log(g))$ for all $f\in C^\infty(\omega)$ and all $g\in \Omega$, where $\log:\Omega\to \omega$ denotes the inverse of $\exp$.

\vspace{2mm}

\noindent The exponential map actually also induces an isomorphism between the corresponding Schwartz spaces. This is the object of the next lemma.

\vspace{2mm}

\begin{lem}\label{lemma 3.3.1}
Let $\omega\subseteq \mathfrak{g}(F)$ be a $G$-excellent open subset and set $\Omega=\exp(\omega)$. Then, the map $f\mapsto f_\omega$ induces a linear isomorphism

$$\mathcal{S}(\Omega)\simeq \mathcal{S}(\omega)$$
\end{lem}

\vspace{2mm}

\noindent\ul{Proof}: This is clear in the $p$-adic case. We assume from now on that $F=\mathbb{R}$. We need to prove the two following facts

\vspace{3mm}

\begin{num}
\item\label{eq 3.3.4} For all $f\in \mathcal{S}(\Omega)$, the function $f_\omega$ belongs to $\mathcal{S}(\omega)$.

\item\label{eq 3.3.5} For all $f\in \mathcal{S}(\omega)$, the function $f_\Omega$ belongs to $\mathcal{S}(\Omega)$.
\end{num}

\vspace{3mm}

\noindent We will prove \ref{eq 3.3.4}, the proof of \ref{eq 3.3.5} being analog (it suffices to replace $\exp$ by $\log$  and $\log$ by $\exp$ in what follows). Let $f\in \mathcal{S}(\Omega)$. Since $\omega$ is relatively compact modulo conjugation, $\overline{\Supp(f)^G}$ is compact modulo conjugation. Hence, there exists a compact subset $\mathcal{K}_0\subseteq \Omega$ such that $\overline{\Supp(f)^G}=\mathcal{K}_0^G$. Then, we have

$$\displaystyle \Supp(f_\omega)^G\subseteq \exp^{-1}\left(\mathcal{K}_0\right)^G$$

\noindent where of course $\exp^{-1}\left(\mathcal{K}_0\right)\subseteq \omega$ is compact. Since $\omega$ is completely $G(F)$-invariant, the closure of $\exp^{-1}\left(\mathcal{K}_0\right)^G$ in $\mathfrak{g}(F)$ is still contained in $\omega$. Hence, we have

$$\displaystyle \overline{\Supp(f_\omega)^G}\subseteq \omega$$

\noindent and it only remains to show that $f_\omega$ belongs to $\mathcal{S}(\mathfrak{g}(F))$ i.e., that $f_\omega$ and all its derivatives are rapidly decreasing. Set $L_G=\overline{\Supp(f)^G}$ and $L_{\mathfrak{g}}=\overline{\Supp(f_\omega)^G}$ so that $L_G$ and $L_{\mathfrak{g}}$ are invariant and compact modulo conjugation subsets of $G(F)$ and $\mathfrak{g}(F)$ respectively and $\exp$ realizes an homeomorphism $L_{\mathfrak{g}}\simeq L_G$. We start by proving that

\begin{align}\label{eq 3.3.6}
\displaystyle \sigma_{\mathfrak{g}}(X)\sim \sigma(e^X),\mbox{ for all }X\in L_{\mathfrak{g}}
\end{align}

\noindent Since $L_G$ is compact modulo conjugation, by Proposition \ref{proposition 3.1.1}(iii), we may find a compact subset $\mathcal{K}\subseteq L_G$ and two maps

$$g\in L_G\mapsto \gamma_g\in G(F)$$
$$g\in L_G\mapsto g_c\in \mathcal{K}$$

\noindent such that

$$g=\gamma_g^{-1}g_c\gamma_g \mbox{ and } \sigma(\gamma_g)\ll \sigma(g)$$

\noindent for all $g\in L_G$. Since $\log(\mathcal{K})$ is compact, we have

$$\sigma_{\mathfrak{g}}(X)=\sigma_{\mathfrak{g}}\left(\gamma_{e^X}^{-1}\log\left((e^X)_c\right)\gamma_{e^X}\right)\ll \sigma(\gamma_{e^X})\ll \sigma(e^X)$$

\noindent for all $X\in L_{\mathfrak{g}}$. This proves one half of \ref{eq 3.3.6}. The other half can be proved similarly, using Proposition \ref{proposition 3.1.1}(iii) for the Lie algebra rather than for the group.

\vspace{2mm}

\noindent The function $j^G$ is bounded on $\omega$ (since it is an invariant function on $\mathfrak{g}(F)$ and $\omega$ is relatively compact modulo conjugation) and hence, it already follows from \ref{eq 3.3.6} that the function $f_\omega$ is rapidly decreasing. Let $u\in S(\mathfrak{g})$. We want to show that the function $\partial(u)f_\omega$ is rapidly decreasing. Since we have

$$\partial(u)f_\omega=\left(\partial(u)_\Omega f\right)_\omega$$

\noindent by what we just saw it suffices to prove that the function $\partial(u)_\Omega f$ is rapidly decreasing (as a function on $G(F)$). For all $D\in \Diff^\infty(\Omega)$ and all $g\in \Omega$, let us denote by $D_g\in\mathcal{U}(\mathfrak{g})$ the unique element such that

$$(Df')(g)=\left(L(D_g)f'\right)(g)$$

\noindent for all $f'\in C^\infty(\Omega)$. Obviously, if $u$ is of degree $k$, then we have $\partial(u)_{\Omega,g}\in \mathcal{U}_{\leqslant k}(\mathfrak{g})$ for all $g\in \Omega$. Let us fix a classical norm $\lvert .\rvert$ on $\mathcal{U}_{\leqslant k}(\mathfrak{g})$. Then, since $f$ is a Schwartz function, to show that $\partial(u)_\Omega f$ is rapidly decreasing we only need to prove that there exist $c>0$ and $m\geqslant 1$ such that

\begin{align}\label{eq 3.3.7}
\displaystyle \left\lvert \partial(u)_{\Omega,g}\right\rvert\leqslant c\lVert g\rVert^m
\end{align}

\noindent for all $g\in L_G$. We have the following easy to check equality

$$\displaystyle \partial(u)_{\Omega,\gamma^{-1}g\gamma}=\gamma^{-1}\partial(\gamma u\gamma^{-1})_{\Omega,g}\gamma$$

\noindent for all $g\in \Omega$ and all $\gamma\in G(F)$. Let us introduce a compact subset $\mathcal{K}\subseteq L_G$ and functions $g \mapsto \gamma_g$ and $g\mapsto g_c$ as before. Then, by the previous identity, for all $g\in \Omega$ we have

$$\displaystyle \partial(u)_{\Omega,g}=\gamma_g^{-1}\partial(\gamma_g u\gamma_g^{-1})_{\Omega,g_c}\gamma_g$$

\noindent for all $g\in L_G$. The inequality \ref{eq 3.3.7} is now easy to deduce from this (note that for all $v\in S(\mathfrak{g})$, the function $g\mapsto \partial(v)_{\Omega,g}$ is bounded on $\mathcal{K}$ and for all $g\in \Omega$ the function $v\in S(\mathfrak{g})\mapsto \partial(v)_{\Omega,g}$ is linear). $\blacksquare$

\vspace{4mm}

\noindent\ul{Remark}: We can extend the definition of $G$-excellent open subsets to the case where $G$ is not necessarily reductive (but is still a connected algebraic group over $F$). The definition is as follows: an open subset $\omega\subseteq \mathfrak{g}(F)$ is {\em $G$-excellent} if it is invariant by $\mathfrak{g}_u(F)$, where $\mathfrak{g}_u$ denotes the Lie algebra of the unipotent radical of $G$, and if moreover its image in $\mathfrak{g}(F)/\mathfrak{g}_u(F)=\left( \mathfrak{g}/\mathfrak{g}_u\right)(F)$  is a $G/G_u$-excellent open subset. If $\omega\subseteq \mathfrak{g}(F)$ is a $G$-excellent subset, then the exponential map still induces an $F$-analytic isomorphism between $\omega$ and $\Omega=\exp(\omega)$ whose Jacobian at $X\in \omega_{\ssi}$ is again given by the same formula \ref{eq 3.3.1}. It follows in particular that the integration formula \ref{eq 3.3.2} is still valid in this more general setting.

\subsection{Parabolic induction of invariant distributions}\label{section 3.4}

\noindent Let $M$ be a Levi subgroup of $G$. Choose a parabolic subgroup $P=MU\in \mathcal{P}(M)$ and a maximal compact subgroup $K$ of $G(F)$ which is special in the $p$-adic case. Fix a Haar measure on $K$ such that

$$\displaystyle \int_{G(F)} f(g)dg=\int_{M(F)}\int_{U(F)}\int_K f(muk)dkdudm$$

\noindent for all $f\in C(G(F))$. We define a continuous linear map

$$C_c^\infty(G(F))\to C_c^\infty(M(F))$$
$$f\mapsto \gls{f(P)}$$

\noindent by setting

$$\displaystyle f^{(P)}(m)=\delta_P(m)^{1/2}\int_K\int_{U(F)} f(k^{-1}muk)dudk$$

\noindent Dually, this defines a linear map

$$\mathcal{D}'(M(F))\to \mathcal{D}'(G(F))$$
$$T\mapsto \gls{TP}$$

\noindent which is uniquely determined by the relations

$$\langle T_P,f\rangle=\langle T,f^{(P)}\rangle$$

\noindent for all $T\in \mathcal{D}'(M(F))$ and all $f\in C_c^\infty(G(F))$. For $T\in \mathcal{D}'(M(F))^M$ an invariant distribution, the distribution $T_P$ is also invariant and doesn't depend on the choices of $P$ or $K$. In this case, we shall denote this distribution by $\gls{iMGT}$ and call it the {\em parabolic induction from $M$ to $G$ of $T$}. If $F=\mathbb{R}$, we have

\begin{align}\label{eq 3.4.1}
zi_M^G(T)=i_M^G(z_MT)
\end{align}

\noindent for all $T\in \mathcal{D}'(M(F))^M$ and all $z\in \mathcal{Z}(\mathfrak{g})$ (Recall that $z\mapsto z_M$ denotes the homomorphism $\mathcal{Z}(\mathfrak{g})\to \mathcal{Z}(\mathfrak{m})$ deduced from the Harish-Chandra isomorphism, cf.\ Section \ref{section 1.1}).

\vspace{2mm}

\noindent If $T=T_{F_M}$ where $F_M$ is an invariant and locally integrable function on $M(F)$, then the distribution $i_M^G(T)$ is also representable by an invariant locally integrable function $F_G$ on $G(F)$. We shall also write $F_G=i_M^G(F_M)$. The function $F_G$ admits the following description in terms of $F_M$. Let us denote by $\gls{Gsreg}$ the subset of strongly regular elements in $G$ (i.e., $x\in G_{\sreg}$ if $Z_G(x)$ is a torus) and let us fix for each $x\in G_{\sreg}(F)$ a set $\gls{XMx}$ of representatives for the $M(F)$-conjugacy classes of elements in $M(F)$ that are $G(F)$-conjugate to $x$. Then, we have the equality

\begin{align}\label{eq 3.4.2}
\displaystyle D^G(x)^{1/2}F_G(x)=\sum_{y\in\mathcal{X}^M(x)} D^M(y)^{1/2} F_M(y)
\end{align}

\noindent for almost all $x\in G_{\sreg}(F)$.

\vspace{2mm}

\noindent If $\sigma\in \Temp(M)$ and $\pi=i_M^G(\sigma)$, then we have

\begin{align}\label{eq 3.4.3}
\theta_\pi=i_M^G(\theta_\sigma)
\end{align}

\noindent We define similarly an induction map $i_M^G$ from the space of invariant distributions on $\mathfrak{m}(F)$ to the space of invariant distributions on $\mathfrak{g}(F)$. Once again, if $F_M$ is an invariant and locally integrable function on $\mathfrak{m}(F)$, the distribution $i_M^G(F_M)$ is also representable by an invariant function $F_G$ and we have

\begin{align}\label{eq 3.4.4}
\displaystyle D^G(X)^{1/2}F_G(X)=\sum_{Y\in\mathcal{X}^M(X)} D^M(Y)^{1/2} F_M(Y)
\end{align}

\noindent for almost all $X\in \mathfrak{g}_{\sreg}(F)$, where this time $\mathcal{X}^M(X)$ is a set of representatives for the $M(F)$-conjugacy classes of elements in $\mathfrak{m}(F)$ that are $G(F)$-conjugate to $X$. We have for example

$$i_M^G\left(\widehat{j}^M(X,.)\right)=\widehat{j}^G(X,.)$$

\noindent for all $X\in \mathfrak{m}(F)\cap \mathfrak{g}_{\reg}(F)$. In particular, if $G$ is quasi-split, $B$ is a Borel subgroup of $G$ and $T_{\quasid}\subset B$ is a maximal torus, it follows from \ref{eq 3.4.4} that we have

\begin{align}\label{eq 3.4.5}
\displaystyle D^G(Y)^{1/2}\widehat{j}^G(X_{\quasid},Y)=\left\{
    \begin{array}{ll}
        \sum_{w\in W(G,T_{\quasid})} \psi\left(B(X_{\quasid},wY)\right) & \mbox{ if } Y\in \mathfrak{t}_{\quasid,\reg}(F) \\
        0 & \mbox{ if } Y\notin \mathfrak{t}_{\quasid,\reg}(F)^G
    \end{array}
\right.
\end{align}

\noindent for all $X_{\quasid}\in \mathfrak{t}_{\quasid,\reg}(F)$ and all $Y\in \mathfrak{g}_{\reg}(F)$. Still assuming that $G$ is quasi-split,we have

$$\displaystyle i_{T_{\quasid}}^G(1)=\sum_{\mathcal{O}\in \Nil_{\reg}(\mathfrak{g})}\widehat{j}(\mathcal{O},.)$$

\noindent from which it follows that

\begin{align}\label{eq 3.4.6}
\displaystyle D^G(X)^{1/2}\sum_{\mathcal{O}\in \Nil_{\reg}(\mathfrak{g})}\widehat{j}(\mathcal{O},X)=\left\{
    \begin{array}{ll}
        \left\lvert W(G,T_{\quasid})\right\rvert & \mbox{ if } X\in \mathfrak{t}_{\quasid,\reg}(F) \\
        0 & \mbox{ if } X\notin \mathfrak{t}_{\quasid,\reg}(F)^G
    \end{array}
\right.
\end{align}

\noindent for all $X\in \mathfrak{g}_{\reg}(F)$. On the maximal torus $T_{\quasid}$, we even have the following more precise equality

\begin{align}\label{eq 3.4.7}
\displaystyle D^G(X)^{1/2}\widehat{j}(\mathcal{O},X)=\lvert W(G,T_{\quasid})\rvert\lvert \Nil_{\reg}(\mathfrak{g})\rvert^{-1}
\end{align}

\noindent for all $X\in \mathfrak{t}_{\quasid,\reg}(F)$ and all $\mathcal{O}\in \Nil_{\reg}(\mathfrak{g})$. Indeed, by \ref{eq 3.4.6} it suffices to show that for all $\mathcal{O}_1,\mathcal{O}_2\in \Nil_{\reg}(\mathfrak{g})$ we have $\widehat{j}(\mathcal{O}_1,X)=\widehat{j}(\mathcal{O}_2,X)$ for all $X\in \mathfrak{t}_{\quasid,\reg}(F)$. Fix $\mathcal{O}_1,\mathcal{O}_2\in \Nil_{\reg}(\mathfrak{g})$. Then, there exists $g_{ad}\in G_{ad}(F)$ such that $g_{ad}^{-1}\mathcal{O}_1g_{ad}=\mathcal{O}_2$, where $G_{ad}$ denotes the adjoint group of $G$. Up to multiplying $g_{ad}$ by an element in $\Ima\left(G(F)\to G_{ad}(F)\right)$, we may assume that $g_{ad}^{-1}Bg_{ad}=B$ and $g_{ad}^{-1}T_{\quasid}g_{ad}=T_{\quasid}$. But then $g_{ad}$ belongs to $T_{\quasid,ad}$ the image of $T_{\quasid}$ in $G_{ad}$. Hence,we have

$$\displaystyle \widehat{j}(\mathcal{O}_1,X)=\widehat{j}(g_{ad}^{-1}\mathcal{O}_1g_{ad},g_{ad}^{-1}Xg_{ad})=\widehat{j}(\mathcal{O}_2,X)$$

\noindent for all $X\in \mathfrak{t}_{\quasid,\reg}(F)$ and this proves the claim.

\section{Quasi-characters}\label{section 4}

The goal of this chapter is to define and establish some crucial properties of what we call {\em quasi-characters} on the group $G(F)$ and its Lie algebra. These are invariant functions which, in some sense, ``locally look like a character". In the $p$-adic case, the notion is due to Waldspurger \cite{Wa1} and in Section \ref{section 4.1} we recall, following {\it loc. cit.}, the definition of quasi-characters and their main properties in this case. The definition in the real case is more technical and is the object of Sections \ref{section 4.2} (for the Lie algebra) and \ref{section 4.4} (for the group). In Section \ref{section 4.3}, we show (still in the real case) that quasi-characters are locally asymptotic to linear combinations of Fourier transforms of regular nilpotent orbital integrals (as do usual characters). In Section \ref{section 4.5}, we associate to any quasi-character $\theta$ a function $c_\theta$ on the set of semi-simple conjugacy classes and study some of its properties. This should be regarded as a regularization of the quasi-character at non-regular elements (where it is not usually defined) and is simply given by averaging the coefficients in the local expansions of $\theta$. In Section \ref{section 4.6}, we study homogeneous distributions on spaces of quasi-characters of the Lie algebra and prove some automatic continuity result for them. In Section \ref{section 4.7}, we study the effect of parabolic induction in the sense of Section \ref{section 3.4} on quasi-characters. Finally, in Section \ref{section 4.8} we collect some properties of characters of tempered representations (as special cases of quasi-characters) and more precisely we recall a well-known link between the coefficients of the local expansion of such a character at $1$ and the existence of Whittaker model for the corresponding representation (a result due to Rodier \cite{Ro} in the $p$-adic case and Matumoto \cite{Mat} in the real case).

\subsection{Quasi-characters when $F$ is $p$-adic}\label{section 4.1}

\noindent In this section we assume that $F$ is $p$-adic. The definition and basic properties of quasi-characters in this case have been established in \cite{Wa1}. We recall them now. Let $\omega\subseteq \mathfrak{g}(F)$ be a completely $G(F)$-invariant open subset. A {\em quasi-character} on $\omega$ is a $G(F)$-invariant smooth function $\theta:\omega_{\reg}\to \mathbb{C}$ satisfying the following condition: for all $X\in \omega_{\ssi}$, there exists $\omega_X\subseteq \mathfrak{g}_X(F)$ a $G$-good open neighborhood of $X$ such that $\omega_X^G\subseteq \omega$ and coefficients $c_{\theta,\mathcal{O}}(X)$ for all $\mathcal{O}\in \Nil(\mathfrak{g}_X)$ such that we have

$$\displaystyle \theta(Y)=\sum_{\mathcal{O}\in \Nil(\mathfrak{g}_X)}c_{\theta,\mathcal{O}}(X)\widehat{j}(\mathcal{O},Y)$$

\noindent for all $Y\in \omega_{X,\reg}$. Note that if $\theta$ is a quasi-character on $\omega$ and $f\in C^\infty(\omega)^G$ then $f\theta$ is also quasi-character on $\omega$. We will denote by $\gls{QComega}$ the space of all quasi-characters on $\omega$ and by $\gls{QCcomega}$ the subspace of quasi-characters on $\omega$ whose support (in $\omega$) is compact modulo conjugation. We will endow $QC_c(\omega)$ with its finest locally convex topology. Note that we have a natural isomorphism

$$QC(\omega)\simeq\varprojlim_{\omega'}QC_c(\omega')$$

\noindent where $\omega'$ runs through the completely $G(F)$-invariant open subsets of $\omega$ that are compact modulo conjugation, the maps $QC(\omega)\to QC_c(\omega')$ being given by $\theta\mapsto \mathbf{1}_{\omega'}\theta$. We equip $QC(\omega)$ with the projective limit topology relative to this isomorphism. To unify notation with the real case, we will also set $\gls{SQCgF}=QC_c(\mathfrak{g}(F))$ and we will call elements of that space Schwartz quasi-characters on $\mathfrak{g}(F)$.

\vspace{2mm}

\noindent Let $\Omega\subseteq G(F)$ be a completely $G(F)$-invariant open subset. A {\em quasi-character} on $\Omega$ is a $G(F)$-invariant smooth function $\theta:\Omega_{\reg}\to \mathbb{C}$ satisfying the following condition: for all $x\in \Omega_{\ssi}$, there exists $\omega_x\subseteq \mathfrak{g}_x(F)$ a $G_x$-excellent open neighborhood of $0$ such that $\left(x\exp(\omega_x)\right)^G\subseteq \Omega$ and coefficients $\gls{cthetaOx}$ for all $\mathcal{O}\in \Nil(\mathfrak{g}_x)$ such that we have the equality

$$\displaystyle \theta(xe^Y)=\sum_{\mathcal{O}\in \Nil(\mathfrak{g}_x)}c_{\theta,\mathcal{O}}(x)\widehat{j}(\mathcal{O},Y)$$

\noindent for all $Y\in \omega_{x,\reg}$. As before, we will denote by $\gls{QCOmega}$ the space of quasi-characters on $\Omega$ and by $\gls{QCcOmega}$ the subspace of quasi-characters that are compactly supported modulo conjugation. We again endow $QC_c(\Omega)$ with its finest locally convex topology and $QC(\Omega)$ with the projective limit topology relative to the natural isomorphism

$$QC(\Omega)\simeq \varprojlim_{\Omega'}QC_c(\Omega')$$

\noindent where $\Omega'$ runs through the completely $G(F)$-invariant open subsets of $\Omega$ that are compact modulo conjugation.

\vspace{2mm}

\begin{prop}\label{proposition 4.1.1}
\begin{enumerate}[(i)]
\item For all $X\in\mathfrak{g}_{\reg}(F)$, $\widehat{j}(X,.)$ is a quasi-character on $\mathfrak{g}(F)$. For all $\mathcal{O}\in \Nil(\mathfrak{g})$, $\widehat{j}(\mathcal{O},.)$ is a quasi-character on $\mathfrak{g}(F)$. For every irreducible admissible representation $\pi$ of $G(F)$, the character $\theta_\pi$ is a quasi-character on $G(F)$.

\item For all $\theta\in QC(G(F))$ (resp.\ $\theta\in QC(\mathfrak{g}(F))$) the function $(D^G)^{1/2}\theta$ is locally bounded.

\item The Fourier transform preserves $SQC(\mathfrak{g}(F))$ in the following sense: for all $\theta\in SQC(\mathfrak{g}(F))$, there exists $\widehat{\theta}\in SQC(\mathfrak{g}(F))$ such that $\widehat{T_\theta}=T_{\widehat{\theta}}$. Moreover, for all $\theta\in SQC(\mathfrak{g}(F))$, we have the equality

$$\displaystyle \widehat{\theta}=\int_{\Gamma(\mathfrak{g})}D^G(X)^{1/2}\theta(X)\widehat{j}(X,.)dX$$

\noindent the integral being absolutely convergent in $QC(\mathfrak{g}(F))$.

\item Let $\omega\subset \mathfrak{g}(F)$ be a $G$-excellent open subset. Set $\Omega=\exp(\omega)$. Then, the linear map

$$\theta\mapsto \theta_\omega$$

\noindent induces topological isomorphisms $QC(\Omega)\simeq QC(\omega)$ and $QC_c(\Omega)\simeq QC_c(\omega)$.

\item Let $X\in \mathfrak{g}_{\ssi}(F)$ and let $\omega_X\subseteq \mathfrak{g}_X(F)$ be a $G$-good open neighborhood of $X$. Set $\omega=\omega_X^G$. Then, the linear map

$$\theta\mapsto \theta_{X,\omega_X}$$

\noindent induces topological isomorphisms $QC(\omega)\simeq QC(\omega_X)$ and $QC_c(\omega)\simeq QC_c(\omega_X)$.

\item Let $x\in G_{\ssi}(F)$ and let $\Omega_x\subseteq G_x(F)$ be a $G$-good open neighborhood of $x$. Set $\Omega=\Omega_x^G$. Then, the linear map

$$\theta\mapsto \theta_{x,\Omega_x}$$

\noindent induces topological isomorphisms $QC(\Omega)\simeq QC(\Omega_x)^{Z_G(x)(F)}$ and $QC_c(\Omega)\simeq QC_c(\Omega_x)^{Z_G(x)(F)}$.

\end{enumerate}
\end{prop}

\vspace{3mm}

\noindent\ul{Proof}:

\begin{enumerate}[(i)]
\item The first part follows from Theorem 4.2 of \cite{Wa1} and the second from Theorem 16.2 of \cite{HCDS}.

\item This follows from the fact that the functions $X\mapsto D^G(X)^{1/2} \widehat{j}(\mathcal{O},X)$ are locally bounded for all $\mathcal{O}\in \Nil(\mathfrak{g})$ (cf.\ Section \ref{section 1.7}).

\item By Theorem 4.2 of \cite{Wa1}, the Fourier transform $\widehat{T}_{\theta}$ of a compactly supported modulo conjugation quasi-character $\theta$ is representable by a quasi-character $\widehat{\theta}$. To see that $\widehat{\theta}$ is again compactly supported modulo conjugation we may appeal to Lemma 6.1 and Proposition 6.4 of \cite{Wa1}. Indeed, by Proposition 6.4 of {\it loc. cit} there exists a strongly cuspidal function $f\in C_c^\infty(\mathfrak{g}(F))$ such that $\theta=\theta_f$ (cf.\ Chapter \ref{section 5} for the definition of strongly cuspidal functions and of the associated quasi-character $\theta_f$). Now, by Lemma 6.1 of {\it loc. cit} we have $\widehat{\theta}_f=\theta_{\widehat{f}}$ where $\widehat{f}$ denotes the usual Fourier transform of the function $f$ (again a strongly cuspidal function). But, by its very definition, the quasi-character $\theta_{\widehat{f}}$ is clearly compactly supported modulo conjugation. Hence, so is $\widehat{\theta}$. Finally, we sketch quickly the proof of the integral formula for $\widehat{\theta}$ since we are going to prove an analogous result over $\mathbb{R}$ (cf.\ Lemma \ref{lemma 4.2.2}(iii)). By Weyl's integration formula and the definition of the functions $\widehat{j}(X,.)$, $X\in \mathfrak{g}_{\reg}(F)$, for all $\varphi\in C_c^\infty(\mathfrak{g}(F))$ we have

\begin{align}\label{eq 4.1.1}
\displaystyle \int_{\mathfrak{g}(F)}\widehat{\theta}(Y)\varphi(Y)dY & =\int_{\mathfrak{g}(F)} \theta(X)\widehat{\varphi}(X)dX \\
\nonumber & =\int_{\Gamma(\mathfrak{g})}D^G(X)^{1/2} \theta(X) J_G(X,\widehat{\varphi})dX \\
\nonumber & =\int_{\Gamma(\mathfrak{g})}D^G(X)^{1/2} \theta(X)\int_{\mathfrak{g}(F)}\widehat{j}(X,Y)\varphi(Y)dYdX \\
\nonumber & =\int_{\mathfrak{g}(F)}\left(\int_{\Gamma(\mathfrak{g})}D^G(X)^{1/2} \theta(X)\widehat{j}(X,Y)dX\right) \varphi(Y)dY 
\end{align}

\noindent By (ii), the function $(D^G)^{1/2}\theta$ is locally bounded. Hence, by \ref{eq 1.7.3} and the fact that the function $(D^G)^{1/2}\theta$ has compact support modulo conjugation, the function

$$\displaystyle Y\in \mathfrak{g}_{\reg}(F)\mapsto \int_{\Gamma(\mathfrak{g})}D^G(X)^{1/2} \left\lvert \theta(X)\right\rvert \left\lvert \widehat{j}(X,Y)\right\rvert dX$$

\noindent is well-defined (i.e., absolutely convergent) and locally essentially bounded by $(D^G)^{-1/2}$. Consequently, by \ref{eq 1.8.1}, the expression \ref{eq 4.1.1} above is absolutely convergent as a double integral. This justifies the above computation and moreover shows that we have an equality

$$\displaystyle \widehat{\theta}(Y)=\int_{\Gamma(\mathfrak{g})}D^G(X)^{1/2} \theta(X)\widehat{j}(X,Y)dX$$

\noindent almost everywhere. To conclude, it suffices to prove that the integral

$$\displaystyle \int_{\Gamma(\mathfrak{g})}D^G(X)^{1/2} \theta(X)\widehat{j}(X,.)dX$$

\noindent is absolutely convergent in $QC(\mathfrak{g}(F))$. By definition of the topology on this space, it suffices to show that for every compact modulo conjugation open subset $\omega\subseteq \mathfrak{g}(F)$ the integral

$$\displaystyle \int_{\Gamma(\mathfrak{g})}D^G(X)^{1/2} \theta(X)\widehat{j}(X,.)_{\mid \omega}dX$$

\noindent is absolutely convergent in $QC_c(\omega)$. But by Howe conjecture, the space spanned by the quasi-characters $\widehat{j}(X,.)_{\mid \omega}$, $X\in \Supp(\theta)\cap \mathfrak{g}_{\reg}(F)$, is finite dimensional and so the absolute convergence of the above integral reduces to the pointwise absolute convergence already established.

\item and (v) are obvious from the definitions. $\blacksquare$
\end{enumerate}

\subsection{Quasi-characters on the Lie algebra for $F=\mathbb{R}$}\label{section 4.2}

\noindent In this section and until the end of Section \ref{section 4.4}, we assume that $F=\mathbb{R}$. Let $\omega\subseteq \mathfrak{g}(\mathbb{R})$ be a completely $G(\mathbb{R})$-invariant open subset. A {\em quasi-character} on $\omega$ is a function $\theta\in C^\infty(\omega_{\reg})^G$ which satisfies the two following conditions

\begin{itemize}
\renewcommand{\labelitemi}{$\bullet$}
\item For all $u\in I(\mathfrak{g})$, the function $(D^G)^{1/2}\partial(u)\theta$ is locally bounded on $\omega$ (so that by \ref{eq 1.8.1}, the function $\partial(u)\theta$ is locally integrable on $\omega$);

\item For all $u\in I(\mathfrak{g})$, we have the following equality of distributions on $\omega$

$$\partial(u)T_\theta=T_{\partial(u)\theta}$$

\end{itemize}

\noindent Notice that the notion of quasi-character is local for the invariant topology: if $\theta\in C^\infty(\omega_{\reg})^G$ then $\theta$ is a quasi-character on $\omega$ if and only if for all $X\in \omega_{\ssi}$ there exists $\omega'\subseteq \omega$ a completely $G$-invariant open neighborhood of $X$ such that $\theta_{\mid \omega'}$ is a quasi-character on $\omega'$. We will say that a quasi-character $\theta$ on $\omega$ is {\em compactly supported} if its support (in $\omega$) is compact modulo conjugation. Finally, a {\em Schwartz quasi-character} is a quasi-character $\theta$ on $\mathfrak{g}(\mathbb{R})$ such that for all $u\in I(\mathfrak{g})$ and for any integer $N\geqslant 1$, we have an inequality

$$D^G(X)^{1/2}\lvert \partial(u)\theta(X)\rvert \ll \lVert X\rVert_{\Gamma(\mathfrak{g})}^{-N}$$

\noindent for all $X\in \mathfrak{g}_{\reg}(\mathbb{R})$. Note that a compactly supported quasi-character is automatically a Schwartz quasi-character. 

\vspace{2mm}

\noindent Any invariant distribution $T$ on some completely $G(\mathbb{R})$-invariant open subset $\omega\subseteq \mathfrak{g}(\mathbb{R})$ such that $\dim(I(\mathfrak{g})T)<\infty$ is the distribution associated to a quasi-character on $\omega$. This follows from the representation theorem of Harish-Chandra on the Lie algebra, cf.\ Theorem 28 p.95 of \cite{Va}. In particular, the functions $\widehat{j}(X,.)$, $X\in\mathfrak{g}_{\reg}(\mathbb{R})$, and the functions $\widehat{j}(\mathcal{O},.)$, $\mathcal{O}\in \Nil(\mathfrak{g})$, are quasi-characters on $\mathfrak{g}(\mathbb{R})$. In our study of quasi-characters, we will need the following lemma which reduces essentially to Proposition 11 p.159 of \cite{Va} using semi-simple descent to maximal tori (cf. the remark after Lemma \ref{lemma 3.2.1}).

\vspace{3mm}

\begin{lem}\label{lemma 4.2.1}
Let $\omega\subseteq \mathfrak{g}(\mathbb{R})$ be a completely $G(\mathbb{R})$-invariant open subset. Let $J\subset I(\mathfrak{g})$ be a subalgebra such that the extension $I(\mathfrak{g})/J$ is finite. Let us define the following topological vector spaces

\begin{itemize}
\renewcommand{\labelitemi}{$\bullet$}
\item $L^\infty_{loc}(\omega,(D^G)^{1/2},\Diff)^G$ is the space of all invariant functions $\theta\in C^\infty(\omega_{\reg})^G$ such that

$$q_{L,D}(\theta)=\sup_{X\in L_{\reg}} D^G(X)^{1/2}\lvert D\theta(X)\rvert<\infty$$

\noindent for all $D\in \Diff(\mathfrak{g})^G$ and each invariant compact modulo conjugation subset $L\subseteq \omega$. We equip $L^\infty_{loc}(\omega,(D^G)^{1/2},\Diff)^G$ with the topology defined by the semi-norms $q_{L,D}$ for all $D$ and $L$ as before.

\item $L^\infty_{loc}(\omega,(D^G)^{1/2},I)^G$ is the space of all invariant functions $\theta\in C^\infty(\omega_{\reg})^G$ such that

$$q_{L,u}(\theta)<\infty$$

\noindent for all $u\in I(\mathfrak{g})$ and each invariant compact modulo conjugation subset $L\subseteq \omega$. We equip $L^\infty_{loc}(\omega,(D^G)^{1/2},I)^G$ with the topology defined by the semi-norms $q_{L,u}$ for all $u$ and $L$ as before.

\item $L^\infty_{loc}(\omega,(D^G)^{1/2},J)^G$ is the space of all invariant functions $\theta\in C^\infty(\omega_{\reg})^G$ such that 

$$q_{L,u}(\theta)<\infty$$

\noindent for all $u\in J$ and each invariant compact modulo conjugation subset $L\subseteq \omega$. We equip $L^\infty_{loc}(\omega,(D^G)^{1/2},J)^G$ with the topology defined by the semi-norms $q_{L,u}$ with $L$ and $u$ as before.

\item $L^\infty_{loc}(\omega,(D^G)^{1/2},T)^G$ is the space of all invariant functions $\theta\in C^\infty(\omega_{\reg})^G$ such that

$$q_{T,u,L_T}=\sup_{X\in L_{T,\reg}}\lvert \partial(u)\theta_T(X)\rvert<\infty$$

\noindent for every maximal torus $T\subset G$, all $u\in S(\mathfrak{t})$ and every compact subset $L_T\subset \omega_T$ (recall that $\omega_T=\omega\cap \mathfrak{t}(\mathbb{R})$ and $\theta_T(X)=D^G(X)^{1/2}\theta(X)$ for all $X\in \omega_{T,\reg}$). We equip $L^\infty_{loc}(\omega,(D^G)^{1/2},T)^G$ with the topology defined by the semi-norms $q_{T,u,L_T}$ with $T$, $u$ and $L_T$ as before.

\end{itemize}

\noindent Then, we have the following equalities of topological vector spaces

$$L^\infty_{loc}(\omega,(D^G)^{1/2},\Diff)^G=L^\infty_{loc}(\omega,(D^G)^{1/2},I)^G=L^\infty_{loc}(\omega,(D^G)^{1/2},J)^G=L^\infty_{loc}(\omega,(D^G)^{1/2},T)^G$$
\end{lem}

\vspace{5mm}

\begin{prop}\label{proposition 4.2.1}
\begin{enumerate}
\item Let $\theta\in C^\infty(\mathfrak{g}_{\reg}(\mathbb{R}))^G$ and assume that there exists $k\geqslant 0$ such that for all $N\geqslant 1$ we have an inequality

$$D^G(X)^{1/2}\lvert \theta(X)\rvert \ll \log\left(2+D^G(X)^{-1}\right)^k \lVert X\rVert_{\Gamma(\mathfrak{g})}^{-N}$$

\noindent for all $X\in\mathfrak{g}_{\reg}(F)$. Then

\begin{enumerate}[(i)]
\item  The function $\theta$ is locally integrable, the distribution $T_\theta$ is tempered and there exists a quasi-character $\widehat{\theta}$ on $\mathfrak{g}(\mathbb{R})$ such that

$$\widehat{T_\theta}=T_{\widehat{\theta}}$$

\noindent Moreover, we have

$$\displaystyle \widehat{\theta}(Y)=\int_{\Gamma(\mathfrak{g})}D^G(X)^{1/2}\theta(X)\widehat{j}(X,Y)dX$$

\noindent for all $Y\in \mathfrak{g}_{\reg}(\mathbb{R})$, the integral being absolutely convergent, and the function $X\in\mathfrak{g}_{\reg}(\mathbb{R})\mapsto D^G(X)^{1/2} \widehat{\theta}(X)$ is (globally) bounded.

\item Assume moreover that $\theta$ is a quasi-character and that for all $u\in I(\mathfrak{g})$ there exists $k\geqslant 0$ such that for all $N\geqslant 1$ we have an inequality

$$D^G(X)^{1/2}\lvert \partial(u)\theta(X)\rvert \ll \log\left(2+D^G(X)^{-1}\right)^k \lVert X\rVert_{\Gamma(\mathfrak{g})}^{-N}$$

\noindent for all $X\in\mathfrak{g}_{\reg}(F)$. Then the function $\widehat{\theta}$ is a Schwartz quasi-character and so is $\theta$.

\item The Fourier transform preserves the space of Schwartz quasi-characters on $\mathfrak{g}(\mathbb{R})$, that is: for every Schwartz quasi-character $\theta$ on $\mathfrak{g}(\mathbb{R})$, the distribution $T_{\theta}$ is tempered and there exists a Schwartz quasi-character $\widehat{\theta}$ such that $\widehat{T_\theta}=T_{\widehat{\theta}}$.

\item For every Schwartz quasi-character $\theta$ on $\mathfrak{g}(\mathbb{R})$ and for all $D\in \Diff(\mathfrak{g})^G$, the function $D\theta$ is a Schwartz quasi-character and we have $DT_{\theta}=T_{D\theta}$. 
\end{enumerate}

\item Let $\omega\subseteq \mathfrak{g}(\mathbb{R})$ be a completely $G(\mathbb{R})$-invariant open subset and let $\theta\in C^\infty(\omega_{\reg})^G$. Then

\begin{enumerate}[(i)]
\item Let $X\in \mathfrak{g}_{\ssi}(\mathbb{R})$ and $\omega_X\subseteq \mathfrak{g}_X(\mathbb{R})$ be a $G$-good open neighborhood of $X$. Assume that $\omega=\omega_X^G$. Then $\theta$ is a quasi-character on $\omega$ if and only if $\theta_{X,\omega_X}$ is a quasi-character on $\omega_X$.

\item Let $J\subseteq I(\mathfrak{g})$ be a subalgebra such that the extension $I(\mathfrak{g})/J$ is finite. Assume that $\theta$ satisfies the two following conditions

\begin{itemize}
\renewcommand{\labelitemi}{$\bullet$}
\item For all $u\in J$, the function $(D^G)^{1/2}\partial(u)\theta$ is locally bounded on $\omega$;

\item For all $u\in J$, we have the equality of distributions on $\omega$

$$\partial(u)T_{\theta}=T_{\partial(u)\theta}$$
\end{itemize}

\noindent Then, $\theta$ is a quasi-character on $\omega$.

\item Assume that $\theta$ is a quasi-character on $\omega$. Then for all $D\in \Diff^\infty(\omega)^G$ the function $D\theta$ is also a quasi-character on $\omega$ and we have the following equality of distributions on $\omega$

$$DT_{\theta}=T_{D\theta}$$

\end{enumerate}
\end{enumerate}
\end{prop}

\vspace{3mm}

\noindent\ul{Proof}:
\begin{enumerate}

\item 

\begin{enumerate}[(i)]

\item First, note that the function $\theta$ is locally integrable and the distribution $T_{\theta}$ is tempered by \ref{eq 1.8.2} and \ref{eq 1.8.3}. Let $\varphi\in C_c^\infty(\mathfrak{g}(\mathbb{R}))$. Then, by the Weyl integration formula, we have

$$\displaystyle \int_{\mathfrak{g}(\mathbb{R})}\theta(X)\widehat{\varphi}(X) dX=\int_{\Gamma(\mathfrak{g})} D^G(X)^{1/2}\theta(X) J_G(X,\widehat{\varphi}) dX$$

\noindent Moreover, by definition of the function $\widehat{j}(.,.)$, we have

$$\displaystyle J_G(X,\widehat{\varphi})=\int_{\mathfrak{g}(\mathbb{R})} \widehat{j}(X,Y)\varphi(Y) dY$$

\noindent for all $X\in \mathfrak{g}_{\reg}(\mathbb{R})$. Hence, we get

\begin{align}\label{eq 4.2.1}
\displaystyle \int_{\mathfrak{g}(\mathbb{R})}\theta(X)\widehat{\varphi}(X) dX=\int_{\Gamma(\mathfrak{g})}\int_{\mathfrak{g}(\mathbb{R})} D^G(X)^{1/2}\theta(X) \widehat{j}(X,Y)\varphi(Y) dY dX
\end{align}

\noindent If we introduce an absolute value inside the double integral above we get an expression which by \ref{eq 1.7.3} is essentially bounded by

$$\displaystyle \int_{\Gamma(\mathfrak{g})} D^G(X)^{1/2} \lvert\theta(X)\rvert dX \int_{\mathfrak{g}(\mathbb{R})} D^G(Y)^{-1/2}\lvert \varphi(Y)\rvert dY$$

\noindent This product is finite by \ref{eq 1.8.1} and \ref{eq 1.8.2}. Hence the double integral \ref{eq 4.2.1} is absolutely convergent. Switching the two integrals, we get

$$\displaystyle \int_{\mathfrak{g}(\mathbb{R})} \theta(X)\widehat{\varphi}(X) dX=\int_{\mathfrak{g}(\mathbb{R})}\widehat{\theta}(Y)\varphi(Y) dY$$

\noindent where

$$\displaystyle \widehat{\theta}(Y)=\int_{\Gamma(\mathfrak{g})} D^G(X)^{1/2}\theta(X)\widehat{j}(X,Y)dX$$

\noindent for all $Y\in \mathfrak{g}_{\reg}(\mathbb{R})$, the integral being absolutely convergent. This shows that $\widehat{T_{\theta}}$ is represented by the locally integrable function $\widehat{\theta}$. It follows from \ref{eq 1.7.3} and \ref{eq 1.8.2}, that the function $(D^G)^{1/2}\widehat{\theta}$ is globally bounded. Let us now show that $\widehat{\theta}$ is a quasi-character. For this, it suffices to show that for all $u\in I(\mathfrak{g})$ the distribution $\partial(u)T_{\widehat{\theta}}$ is representable by a function which is locally essentially bounded by $(D^G)^{-1/2}$. Let $u\in I(\mathfrak{g})$. Since $T_{\widehat{\theta}}=\widehat{T_\theta}$, the distribution $\partial(u)T_{\widehat{\theta}}$ is the Fourier transform of $p_uT_\theta=T_{p_u\theta}$. But it is not hard to see that $p_u\theta$ satisfies the same hypothesis as $\theta$ and so its Fourier transform is also representable by a function which is (globally) essentially bounded by $(D^G)^{-1/2}$. This shows that $\widehat{\theta}$ is indeed a quasi-character.

\item Let $u\in I(\mathfrak{g})$ and $p\in I(\mathfrak{g}^*)$. Then the function $p\partial(u)\theta$ satisfies the same hypothesis as $\theta$ in 1.(i) and we have $T_{p\partial(u)\theta}=p\partial(u)T_\theta$ (since $\theta$ is a quasi-character). Hence, by 1.(i), the Fourier transform of $p\partial(u)T_{\theta}$ is representable by a function, which is necessarily $\partial(u_p)(p_u\widehat{\theta})$, and moreover the function $(D^G)^{1/2}\partial(u_p)(p_u\widehat{\theta})$ is globally bounded. Let us show

\begin{num}
\item\label{eq 4.2.2} For all $u\in I(\mathfrak{g})$ and all $p\in I(\mathfrak{g}^*)$, the function
$$X\in \mathfrak{g}_{\reg}(\mathbb{R})\mapsto D^G(X)^{1/2}\lvert p(X)\rvert\lvert\partial(u)\widehat{\theta}(X)\rvert$$
is bounded.
\end{num}

Let $u$ and $p$ be as above. Since $p$ is bounded (in absolute value) by an element in $I(\mathfrak{g}^*)$ which is positive and bounded by below on $\mathfrak{g}(\mathbb{R})$ (just take $1+p\overline{p}$), we may assume that $p$ has this property. By what we just saw, for every integer $k\geqslant 1$ the function $(D^G)^{1/2}\partial(u)(p^k\widehat{\theta})$ is bounded. Consider the endomorphisms $R(p), L(p)$ and $\ad(p)$ of $\Diff(\mathfrak{g})$ given by $R(p)D=Dp$, $L(p)D=pD$ and $\ad(p)=L(p)-R(p)$. They all commute with each other and $\ad(p)$ is locally nilpotent. Hence, there exists an integer $M\geqslant 1$ such that $\ad(p)^M(\partial(u))=0$. It follows that for every integer $n\geqslant M$, we have

\[\begin{aligned}
\displaystyle p^n\partial(u)=L(p)^n(\partial(u)) & =\left(R(p)+\ad(p)\right)^n\partial(u) \\
 & =\sum_{k=0}^{M-1}\binom{n}{k} (R(p)^{n-k}\ad(p)^k)(\partial(u)) \\
 & =\left(\sum_{k=0}^{M-1} \binom{n}{k} \ad(p)^k(\partial(u))p^{M-1-k}\right) p^{n+1-M}
\end{aligned}\]

\noindent The last sum above stays in a subspace of dimension less than $M$ as $n$ varies. It easily follows that we may find two integers $n,m\geqslant 1$ and scalars $\lambda_1,\ldots,\lambda_m\in \mathbb{C}$ such that

$$p^{n+m}\partial(u)=\lambda_1p^{n+m-1}\partial(u)p+\ldots+p^n\partial(u)p^m$$

\noindent Because the functions $(D^G)^{1/2}\partial(u)(p^k\widehat{\theta})$, $1\leqslant k\leqslant m$, are all globally bounded, we get an inequality

$$D^G(X)^{1/2}\lvert p(X)\rvert^{n+m} \lvert\partial(u)\widehat{\theta}(X)\rvert\ll \lvert p(X)\rvert^{n+m-1}+\ldots+\lvert p(X)\rvert^{n}$$

\noindent for all $X\in\mathfrak{g}_{\reg}(\mathbb{R})$. Since $\lvert p\rvert$ is bounded by below, the last sum above is essentially bounded by $\lvert p(X)\rvert^{n+m-1}$. Then, after dividing by $\lvert p(X)\rvert^{n+m-1}$, we obtain \ref{eq 4.2.2}.

\vspace{2mm}

It is easy to see that we may find $p\in I(\mathfrak{g}^*)$ such that $\lVert X\rVert_{\Gamma(\mathfrak{g})}\ll \lvert p(X)\rvert$, for all $X\in \mathfrak{g}(\mathbb{R})$. Hence, it follows from \ref{eq 4.2.2} and 1.(i) that $\widehat{\theta}$ is a Schwartz quasi-character. It implies in particular that $\widehat{\theta}$ satisfies the same condition as $\theta$, hence its Fourier transform, which is just the function $X\mapsto \theta(-X)$, is also a Schwartz quasi-character. This shows that $\theta$ is itself a Schwartz quasi-character.

\item This follows directly from 1.(ii) (Note that a Schwartz quasi-character $\theta$ satisfies the assumptions of 1.(ii)).

\item Denote by $\mathcal{A}$ the subalgebra of $\Diff(\mathfrak{g})^G$ consisting of the operators $D$ such that for every Schwartz quasi-character $\theta$, the function $D\theta$ is also a Schwartz quasi-character and $DT_{\theta}=T_{D\theta}$. We want to show that $\mathcal{A}=\Diff(\mathfrak{g})^G$. Obviously we have $I(\mathfrak{g})\subset \mathcal{A}$. Since the Fourier transform preserves the space of Schwartz quasi-characters, it easily follows that $I(\mathfrak{g}^*)$ is also included in $\mathcal{A}$. By Proposition \ref{proposition 3.1.1}(vi), it follows that $\mathcal{A}=\Diff(\mathfrak{g})^G$.

\end{enumerate}

\item We are going to prove 2.(i), 2.(ii) and 2.(iii) by induction on $\dim(G)$. If $\dim(G)=1$, then $G$ is a torus and everything is obvious. We henceforth assume that 2.(i), 2.(ii) and 2.(iii) hold for every connected reductive groups $G'$ with $\dim(G')<\dim(G)$.

\vspace{2mm}

We first establish 2.(i). The direction $\theta_{X,\omega_X}$ quasi-character $\Rightarrow$ $\theta$ quasi-character is easy using Lemma \ref{lemma 3.2.1}(i). So assume that $\theta$ is a quasi-character. If $G_X=G$, there is nothing to prove. If $G_X\neq G$, by Lemma \ref{lemma 3.2.1}(i), we see that $\theta_{X,\omega_X}$ satisfies 2.(ii) for $J\subseteq I(\mathfrak{g}_X)$ the image of $I(\mathfrak{g})$ by the morphism $u\mapsto u_X$. By the induction hypothesis, it follows that $\theta_{X,\omega_X}$ is indeed a quasi-character and this ends the proof of 2.(i).

\vspace{2mm}

We now prove 2.(ii) and 2.(iii) together. That is, we take $\theta$ and $J$ as in 2.(ii) and we are going to prove that for all $D\in \Diff^\infty(\omega)^G$ the function $D\theta$ is a quasi-character on $\omega$ and that we have the equality $DT_{\theta}=T_{D\theta}$ of distributions on $\omega$. By definition, this amounts to showing the two following facts

\begin{num}
\item\label{eq 4.2.3} For all $D\in \Diff^\infty(\omega)^G$ the function $(D^G)^{1/2}D\theta$ is locally bounded on $\omega$.

\item\label{eq 4.2.4} For all $D\in \Diff^\infty(\omega)^G$ we have the equality $DT_{\theta}=T_{D\theta}$ of distributions on $\omega$.
\end{num}

\vspace{3mm}

\noindent By Proposition \ref{proposition 3.1.1}(v), we are immediately reduced to proving \ref{eq 4.2.3} and \ref{eq 4.2.4} only for $D\in \Diff(\mathfrak{g})^G$. By Lemma \ref{lemma 4.2.1}, we already know that for all $D\in \Diff(\mathfrak{g})^G$ the function $(D^G)^{1/2}D\theta$ is locally bounded on $\omega$. Hence, we only need to show the following

\vspace{3mm}

\begin{num}
\item\label{eq 4.2.5} For all $D\in \Diff(\mathfrak{g})^G$, we have the equality $DT_{\theta}=T_{D\theta}$ of distributions on $\omega$.
\end{num}

\vspace{3mm}

\noindent Let $D\in \Diff(\mathfrak{g})^G$. The question is local for the invariant topology i.e., we only need to prove that the equality holds near every $X\in \omega_{\ssi}$. So let $X\in \omega_{\ssi}$. Assume first that $G_X\neq G$. In this case, we can use semi-simple descent and the induction hypothesis on $G_X$. More precisely, let $\omega_X\subseteq \mathfrak{g}_X(\mathbb{R})$ be a $G$-good open neighborhood of $X$ such that $\omega_X^G\subseteq \omega$. Then, by Lemma \ref{lemma 3.2.1}(i), the function $\theta_{X,\omega_X}$ satisfies the assumptions of 2.(ii) on $\omega_X$ with $J$ replace by its image in $I(\mathfrak{g}_X)$ via the morphism $u\mapsto u_X$. Since the two extensions $I(\mathfrak{g}_X)/I(\mathfrak{g})$ and $I(\mathfrak{g})/J$ are finite, so is $I(\mathfrak{g}_X)/J$. Hence, applying 2.(ii) to $G_X$, we see that $\theta_{X,\omega_X}$ is a quasi-character on $\omega_X$. Then applying 2.(iii) to $G_X$, we see that we have the equality of distributions $D_XT_{\theta_{X,\omega_X}}=T_{D_X\theta_{X,\omega_X}}$ on $\omega_X$. By Lemma \ref{lemma 3.2.1}(i), this implies that the equality of distributions $DT_\theta=T_{D\theta}$ holds on $\omega_X^G$. This proves the claim \ref{eq 4.2.5} near $X$ in this case. Now assume that $G_X=G$. This condition is equivalent to $X\in \mathfrak{z}_G(\mathbb{R})\cap \omega$, where $\mathfrak{z}_G$ denotes the Lie algebra of $Z_G$. Let $G_{\der}$ be the derived subgroup of $G$ and let $\mathfrak{g}_{\der}$ denote its Lie algebra. We have the decompositions $\mathfrak{g}=\mathfrak{z}_G\oplus \mathfrak{g}_{\der}$ and $I(\mathfrak{g})=S(\mathfrak{z}_G)I(\mathfrak{g}_{\der})$. The question being local at $X$, we may replace $\omega$ by any completely $G(\mathbb{R})$-invariant open neighborhood of $X$ that is contained in $\omega$. In particular, we may assume without loss of generality that $\omega=\omega_{\mathfrak{z}}\times\omega_{\der}$ where $\omega_{\mathfrak{z}}\subseteq \mathfrak{z}_G(\mathbb{R})$ and $\omega_{\der}\subseteq \mathfrak{g}_{\der}(\mathbb{R})$ are open and completely $G(\mathbb{R})$-invariant. Note that we have $\omega_{\reg}=\omega_{\mathfrak{z}}\times \omega_{\der,\reg}$. Let $f_{\der}\in C^\infty(\omega_{\der})^G$ be such that

\vspace{2mm}

\begin{itemize}
\renewcommand{\labelitemi}{$\bullet$}
\item $\Supp(\omega_{\der})$ is compact modulo conjugation;
\item $f_{\der}=1$ near $0$.
\end{itemize}

\vspace{2mm}

\noindent Then, we claim that the function $f_{\der}\theta$ satisfies the same hypothesis as $\theta$, that is:

\vspace{3mm}

\begin{num}
\item\label{eq 4.2.6} for all $u\in J$, the function $(D^G)^{1/2}\partial(u)(f_{\der}\theta)$ is locally bounded on $\omega$ and we have the equality $\partial(u)T_{f_{\der}\theta}=T_{\partial(u)(f_{\der}\theta)}$ of distributions on $\omega$.
\end{num}

\vspace{3mm}

\noindent This is true near $\omega_{\mathfrak{z}}\times \{0\}$ since on some neighborhood of it $f_{\der}\theta$ coincide with $\theta$. For $X\in \omega_{\ssi}\backslash \omega_{\mathfrak{z}}\times\{0\}$, we can use semi-simple descent and the induction hypothesis (note that $G_X\neq G$) to show that \ref{eq 4.2.6} holds near $X$. Indeed, we already saw that there exists $\omega_X\subseteq \omega$ a $G$-good open neighborhood of $X$ such that $\theta_{X,\omega_X}$ is a quasi-character on $\omega_X$. By the property 2.(iii) applied to $G_X$, it follows that $(f_{\der})_{X,\omega_X}\theta_{X,\omega_X}$ is also a quasi-character on $\omega_X$. Hence, by 2.(i), the function $f_{\der}\theta$ is a quasi-character on $\omega_X^G$ and so a fortiori \ref{eq 4.2.6} is satisfied on $\omega_X^G$.

\vspace{2mm}

\noindent Since $f_{\der}\theta$ and $\theta$ coincide near $X$, we may replace $\theta$ by $f_{\der}\theta$ (recall that we want to show that the equality $DT_\theta=T_{D\theta}$ holds near $X$). Doing this, the function $\theta$ will now satisfy the following additional assumption

\vspace{3mm}

\begin{num}
\item\label{eq 4.2.7} There exists a $G(\mathbb{R})$-invariant compact modulo conjugation subset $L_{\der}\subset \omega$ such that $\Supp(\theta)\subseteq \omega_{\mathfrak{z}}\times L_{\der}$.
\end{num}

\vspace{3mm}

\noindent Let us now show the following

\vspace{3mm}

\begin{num}
\item\label{eq 4.2.8} For all $u\in S(\mathfrak{z}_G)$ and all $v\in J$, we have $\partial(uv)T_{\theta}=T_{\partial(uv)\theta}$ on $\omega$.
\end{num}

\vspace{3mm}

\noindent Indeed, for all $u\in S(\mathfrak{z}_G)$, all $v\in J$ and all $\varphi\in C_c^\infty(\omega)$, we have

\[\begin{aligned}
\displaystyle \int_{\omega} (\partial(v^*u^*)\varphi)(Y)\theta(Y)dY & =\int_{\omega} (\partial(u^*)\varphi)(Y)(\partial(v)\theta)(Y)dY \\
 & =\int_{\omega_{\der,\reg}}\int_{\omega_{\mathfrak{z}}}(\partial(u^*)\varphi)(Y_{\mathfrak{z}}+Y_{\der})(\partial(v)\theta)(Y_{\mathfrak{z}}+Y_{\der})dY_{\mathfrak{z}}dY_{\der} \\
 & =\int_{\omega_{\der,\reg}}\int_{\omega_{\mathfrak{z}}}\varphi(Y_{\mathfrak{z}}+Y_{\der})(\partial(uv)\theta)(Y_{\mathfrak{z}}+Y_{\der})dY_{\mathfrak{z}}dY_{\der} \\
 & =\int_{\omega} \varphi(Y)(\partial(uv)\theta)(Y)dY
\end{aligned}\]

\noindent where in the first equality we have used the equality $\partial(v)T_{\theta}=T_{\partial(v)\theta}$, in the third equality we have used the fact that the function $Y_{\mathfrak{z}}\in \omega_{\mathfrak{z}}\mapsto (\partial(v)\theta)(Y_{\mathfrak{z}}+Y_{\der})$ is smooth for all $Y_{\der}\in \omega_{\der,\reg}$ and in the fourth equality we have used the fact that the function $\partial(uv)\theta$ is locally integrable (since it is locally bounded by $(D^G)^{-1/2}$). This proves \ref{eq 4.2.8}. Up to replacing $J$ by $S(\mathfrak{z}_G)J$, we may now assume that $S(\mathfrak{z}_G)\subseteq J$. 

\vspace{2mm}

\noindent Choose $f_{\mathfrak{z}}\in C_c^\infty(\omega_{\mathfrak{z}})$ such that $f_{\mathfrak{z}}=1$ near $X$. Then, the function $f_{\mathfrak{z}}\theta$ coincides near $\theta$ with $X$. By 1.(iv), we thus only need to show that $f_{\mathfrak{z}}\theta$ is a Schwartz quasi-character on $\mathfrak{g}(\mathbb{R})$. Since, by \ref{eq 4.2.7}, the support of this function in $\omega$ is compact modulo conjugation, we even only need to prove that it is a quasi-character. Actually, we are going to prove this for all $f_{\mathfrak{z}}\in C^\infty_c(\omega_{\mathfrak{z}})$. For all $N\geqslant 1$, the function $\lVert .\rVert_{\Gamma(\mathfrak{g})}^N (D^G)^{1/2}f_{\mathfrak{z}}\theta$ is globally bounded (since it is locally bounded, invariant and compactly supported modulo conjugation). Hence, the function $f_{\mathfrak{z}}\theta$ satisfies the assumption of 1.(i) (with $k=0$) and so we know that the distribution $T_{f_{\mathfrak{z}}\theta}$ is tempered and that its Fourier transform is representable by a quasi-character $\widehat{f_{\mathfrak{z}}\theta}$ on $\mathfrak{g}(\mathbb{R})$ which is globally bounded by $(D^G)^{-1/2}$. We claim that we have the following

\vspace{3mm}

\begin{num}
\item\label{eq 4.2.9} For all $p\in I(\mathfrak{g}^*)$, the function 
$$Y\in \mathfrak{g}_{\reg}(\mathbb{R})\mapsto D^G(Y)^{1/2}p(Y)\widehat{f_{\mathfrak{z}}\theta}(Y)$$
is bounded.
\end{num}

\vspace{3mm}

\noindent Since $u\mapsto p_u$ is an isomorphism $I(\mathfrak{g})\simeq I(\mathfrak{g}^*)$ and $I(\mathfrak{g})=S(\mathfrak{z}_G)I(\mathfrak{g}_{\der})$, we only need to prove \ref{eq 4.2.9} when $p$ is a product $p_{u_{\mathfrak{z}}}p_{u_{\der}}$ with $u_{\mathfrak{z}}\in S(\mathfrak{z}_G)$ and $u_{\der}\in I(\mathfrak{g}_{\der})$. Note that the function $p_{u_{\mathfrak{z}}}(\widehat{f_{\mathfrak{z}}\theta})$ is the Fourier transform of $\partial(u_{\mathfrak{z}})f_{\mathfrak{z}}T_{\theta}$. Since $S(\mathfrak{z}_G)\subset J$ and the function $f_{\mathfrak{z}}\theta$ is supported in $\omega$, this distribution is represented by a function of the form 

$$\displaystyle \sum_{i=1}^N (\partial(u_i)f_{\mathfrak{z}})(\partial(v_i)\theta)$$

\noindent with $u_i,v_i\in S(\mathfrak{z}_G)$, $1\leqslant i\leqslant N$. Note that the function $\partial(v_i)\theta$ satisfies the same assumptions as $\theta$ (including \ref{eq 4.2.7}) and that the functions $\partial(u_i)f_{\mathfrak{z}}$ belong to $C_c^\infty(\omega_{\mathfrak{z}})$. Hence $\partial(u_{\mathfrak{z}})f_{\mathfrak{z}}\theta$ is a sum of functions of the same type than $f_{\mathfrak{z}}\theta$. It follows that we only need to prove \ref{eq 4.2.9} for $p=p_u$ with $u\in I(\mathfrak{g}_{\der})$. Since $p$ is bounded (in absolute value) by an element of $I(\mathfrak{g}_{\der}^*)$ that is positive and bounded by below on $\mathfrak{g}(\mathbb{R})$, we may assume that $p$ has this property. Because the extension $I(\mathfrak{g})/J$ is finite, we may find an integer $n\geqslant 1$ and elements $v_0,\ldots,v_{n-1}\in J$ such that

$$u^n=v_{n-1}u^{n-1}+\ldots+v_1u+v_0$$

\noindent The function $p_u^n(\widehat{f_{\mathfrak{z}}\theta})$ is the Fourier transform of

$$\partial(u^n)f_{\mathfrak{z}}T_\theta=f_{\mathfrak{z}}\partial(u^n)T_\theta=\sum_{i=0}^{n-1} f_{\mathfrak{z}}\partial(u^i)\partial(v_i)T_\theta=\sum_{i=0}^{n-1} \partial(u^i)f_{\mathfrak{z}}T_{\partial(v_i)\theta}$$

\noindent Note that, although the equality $\partial(v_i)T_\theta=T_{\partial(v_i)\theta}$ only holds on $\omega$, we can use it here as if it holds everywhere since the support of $f_{\mathfrak{z}}\theta$ in $\mathfrak{g}(\mathbb{R})$ is contained in $\omega$. Applying the Fourier transform to this equality, we obtain

\begin{align}\label{eq 4.2.10}
\displaystyle p_u^n\widehat{T}_{f_{\mathfrak{z}}\theta}=\sum_{i=0}^{n-1} p_u^i\widehat{T}_{f_{\mathfrak{z}}\partial(v_i)\theta}
\end{align}

\noindent Note that for all $1\leqslant i\leqslant n-1$, the function $f_{\mathfrak{z}}\partial(v_i)\theta$ is globally bounded by $(D^G)^{-1/2}$ and compactly supported modulo conjugation. Consequently, its Fourier transform is also representable by a function that is globally bounded by $(D^G)^{-1/2}$ (by 1.(i)). Hence, by \ref{eq 4.2.10}, we get an inequality

$$D^G(X)^{1/2}\lvert p_u(X)\rvert^n \lvert\widehat{f_{\mathfrak{z}}\theta}(X)\rvert \ll 1+\lvert p_u(X)\rvert+\ldots+\lvert p_u(X)\rvert^{n-1}$$

\noindent for all $X\in \mathfrak{g}_{\reg}(\mathbb{R})$. Since $\lvert p_u\rvert$ is bounded by below, we easily deduce that the function $(D^G)^{1/2}p_u(\widehat{f_{\mathfrak{z}}\theta})$ is bounded. This proves \ref{eq 4.2.9}.

\vspace{2mm}

\noindent Since there exists $p\in I(\mathfrak{g}^*)$ such that $\lVert X\rVert_{\Gamma(\mathfrak{g})}\ll \lvert p(X)\rvert$ for all $X\in \mathfrak{g}(\mathbb{R})$, it follows from \ref{eq 4.2.9} that the function $\widehat{f_{\mathfrak{z}}\theta}$ satisfies the assumptions of 1.(i). Hence, its Fourier transform, which is the function $Y\mapsto (f_{\mathfrak{z}}\theta)(-Y)$, is a quasi-character. This proves the claim that $f_{\mathfrak{z}}\theta$ is a quasi-character and ends the proof of 2.(ii) and 2.(iii). $\blacksquare$
\end{enumerate}

\vspace{3mm}

\noindent Let $\omega\subseteq \mathfrak{g}(\mathbb{R})$ be a completely $G(\mathbb{R})$-invariant open subset. We will denote by $\gls{QComega}$ the space of quasi-characters on $\omega$ and by $\gls{QCcomega}$ the subspace of compactly supported quasi-characters on $\omega$. If $L\subset \omega$ is invariant and compact modulo conjugation, we will also denote by $QC_L(\omega)\subset QC_c(\omega)$ the subspace of quasi-characters with support in $L$. Finally, we will denote by $\gls{SQCgR}$ the space of Schwartz quasi-characters on $\mathfrak{g}(\mathbb{R})$. \\

\noindent We will endow these spaces with locally convex topologies as follows. For $L\subset \omega$ as before and $u\in I(\mathfrak{g})$, we define a semi-norm $\gls{qLu}$ on $QC(\omega)$ by

$$\displaystyle q_{L,u}(\theta)=\sup_{X\in L_{\reg}} D^G(X)^{1/2}\lvert \partial(u)\theta(X)\rvert,\;\;\; \theta\in QC(\omega)$$

\noindent Then, we equip $QC_L(\omega)$ with the topology defined by the semi-norms $(q_{L,u})_{u\in I(\mathfrak{g})}$ and $QC(\omega)$ with the topology defined by the semi-norms $(q_{L,u})_{L,u}$ where $L$ runs through the invariant compact modulo conjugation subsets of $\omega$ and $u$ runs through $I(\mathfrak{g})$. We have a natural isomorphism

$$\displaystyle QC_c(\omega)\simeq\varinjlim\limits_{L} QC_L(\omega)$$

\noindent and we endow $QC_c(\omega)$ with the direct limit topology. Finally, we put on $SQC(\mathfrak{g}(\mathbb{R}))$ the topology defined by the semi-norms

$$\displaystyle \gls{quN}(\theta)=\sup_{X\in \mathfrak{g}_{\reg}(\mathbb{R})} \lVert X\rVert^N_{\Gamma(\mathfrak{g})}D^G(X)^{1/2}\lvert \partial(u)\theta(X)\rvert,\;\;\; \theta\in SQC(\mathfrak{g}(\mathbb{R}))$$

\noindent where $u$ runs through $I(\mathfrak{g})$ and $N$ runs through all positive integers.

\vspace{2mm}

\begin{lem}\label{lemma 4.2.2}
Let $\omega\subseteq \mathfrak{g}(\mathbb{R})$ be a completely $G(\mathbb{R})$-invariant open subset and $L\subseteq \omega$ be invariant and compact modulo conjugation. Then

\begin{enumerate}[(i)]
\item $QC(\omega)$, $QC_L(\omega)$ and $SQC(\mathfrak{g}(\mathbb{R}))$ are Fr\'echet spaces whereas $QC_c(\omega)$ is an LF space. The inclusions $QC_c(\omega)\subset QC(\omega)$ and $SQC(\mathfrak{g}(\mathbb{R}))\subset QC(\mathfrak{g}(\mathbb{R}))$ are continuous and moreover $QC_c(\omega)$, $QC_L(\omega)$ and $QC(\omega)$ are nuclear spaces.

\item Let $X\in \mathfrak{g}_{\ssi}(\mathbb{R})$ and $\omega_X\subseteq \mathfrak{g}_X(\mathbb{R})$ be a $G$-good open neighborhood of $X$. Assume that $\omega=\omega_X^G$. Then the linear map

$$\theta\mapsto \theta_{X,\omega_X}$$

\noindent induces topological isomorphisms $QC(\omega)\simeq QC(\omega_X)$ and $QC_c(\omega)\simeq QC_c(\omega_X)$.

\item The Fourier transform $\theta\mapsto \widehat{\theta}$ is a continuous linear automorphism of $SQC(\mathfrak{g}(\mathbb{R}))$ and for all $\theta\in SQC(\mathfrak{g}(F))$, we have

$$\displaystyle \widehat{\theta}=\int_{\Gamma(\mathfrak{g})} D^G(X)^{1/2} \theta(X)\widehat{j}(X,.)dX$$

\noindent the integral above being absolutely convergent in $QC(\mathfrak{g}(F))$. Moreover, for all $D\in \Diff(\mathfrak{g})^G$ the linear map

$$\theta\in SQC(\mathfrak{g}(\mathbb{R}))\mapsto D\theta\in SQC(\mathfrak{g}(\mathbb{R}))$$

\noindent is continuous.

\item The two bilinear maps

$$\Diff^\infty(\omega)^G\times QC_c(\omega)\to QC_c(\omega)\;\;\; \Diff^\infty(\omega)^G\times QC(\omega)\to QC(\omega)$$
$$(D,\theta)\mapsto D\theta$$

\noindent are separately continuous.

\item $QC_c(\omega)$ is dense in $QC(\omega)$ and $QC_c(\mathfrak{g}(\mathbb{R}))$ is dense in $SQC(\mathfrak{g}(\mathbb{R}))$.
\end{enumerate}
\end{lem}

\vspace{2mm}

\noindent\ul{Proof}:

\begin{enumerate}[(i)]

\item The claim about inclusions is obvious. Using Proposition \ref{proposition 3.1.1}(i), it is clear that the topologies on $QC(\omega)$, $QC_L(\omega)$ and $SQC(\omega)$ may be defined by a countable number of semi-norms and that $QC_c(\omega)$ is the direct limit of a countable family $(QC_{L_n}(\omega))_{n\geqslant 1}$ where $(L_n)_{n\geqslant 1}$ is an increasing sequence of invariant compact modulo conjugation subsets of $\omega$. Moreover, $QC_L(\omega)$ is a closed subspace of $QC(\omega)$ and the topology on $QC_L(\omega)$ is induced from the one on $QC(\omega)$. Hence, it suffices to show that $QC(\omega)$ and $SQC(\mathfrak{g}(\mathbb{R}))$ are complete and that $QC(\omega)$ is nuclear. Let $(\theta_n)_{n\geqslant 1}$ be a Cauchy sequence in $SQC(\mathfrak{g}(\mathbb{R}))$. Then, it is also a Cauchy sequence in $QC(\mathfrak{g}(\mathbb{R}))$. If this sequence admits a limit $\theta$ in $QC(\mathfrak{g}(\mathbb{R}))$, then it is clear that $\theta$ belongs to $SQC(\mathfrak{g}(\mathbb{R}))$ and that it is also a limit of the sequence $(\theta_n)_{n\geqslant 0}$ in $SQC(\mathfrak{g}(\mathbb{R}))$. Hence, we are only left with proving that $QC(\omega)$ is nuclear and complete. Let $\mathcal{T}(G)$ be a set of representatives for the $G(\mathbb{R})$-conjugacy classes of maximal tori in $G$. For all $T\in \mathcal{T}(G)$ set $\omega_T=\omega\cap \mathfrak{t}(\mathbb{R})$, $\omega_{T,\reg}=\omega\cap \mathfrak{t}_{\reg}(\mathbb{R})$ and define $C^\infty_b(\omega_{T,\reg},\omega_T)$ to be the space of all smooth functions $f:\omega_{T,\reg}\to \mathbb{C}$ such that $\partial(u)f$ is locally bounded in $\omega_T$ for all $u\in S(\mathfrak{t})$. We endow this space with the topology defined by the semi-norms

$$q_{T,u,L_T}(f)=\sup_{X\in L_{T,\reg}} \left\lvert (\partial(u)f)(X)\right\rvert,\;\;\; f\in C_b^\infty(\omega_{T,\reg},\omega_T)$$

\noindent for all $u\in S(\mathfrak{t})$ and every compact subset $L_T\subset \omega_T$. Then by Lemma \ref{lemma A.5.1}, the spaces $C^\infty_b(\omega_{T,\reg},\omega_T)$, $T\in \mathcal{T}(G)$, are all nuclear Fr\'echet spaces (note that since $\omega_{T,\reg}$ is the complement in $\omega_T$ of a finite union of subspaces of $\mathfrak{t}(\mathbb{R})$, the pair $(\omega_{T,\reg},\omega_T)$ trivially satisfies the assumption of Lemma \ref{lemma A.5.1}). Moreover, by Lemma \ref{lemma 4.2.1}, the linear map

$$\theta\mapsto (\theta_T)_{T\in\mathcal{T}(G)}$$

\noindent where $\theta_T(X)=D^G(X)^{1/2}\theta(X)$ for all $X\in \omega_{T,\reg}$, induces a closed embedding

$$QC(\omega)\hookrightarrow \bigoplus_{T\in\mathcal{T}(G)} C^\infty_b(\omega_{T,\reg},\omega_T)$$

\noindent The result follows.

\item By Proposition \ref{proposition 4.2.1}(i), the map $\theta\mapsto \theta_{X,\omega_X}$ induces linear isomorphisms $QC(\omega)\simeq QC(\omega_X)$ and $QC_c(\omega)\simeq QC_c(\omega_X)$. Moreover, the inverses of these isomorphisms are easily seen to be continuous. By the open mapping theorem, it follows that these are indeed topological isomorphisms.

\item First we prove the claim about the Fourier transform. By the closed graph theorem and Proposition \ref{proposition 4.2.1}(i), it is sufficient to prove that for all $\theta\in SQC(\mathfrak{g}(\mathbb{R}))$ the integral

\begin{align}\label{eq 4.2.11}
\displaystyle\int_{\Gamma(\mathfrak{g})} D^G(X)^{1/2} \theta(X) \widehat{j}(X,.)dX
\end{align}

\noindent is absolutely convergent in $QC(\mathfrak{g}(\mathbb{R}))$ and that the linear map

$$SQC(\mathfrak{g}(\mathbb{R}))\to QC(\mathfrak{g}(\mathbb{R}))$$
$$\displaystyle \theta\mapsto \int_{\Gamma(\mathfrak{g})} D^G(X)^{1/2} \theta(X) \widehat{j}(X,.)dX$$

\noindent is continuous. Let $L\subset \mathfrak{g}(\mathbb{R})$ be invariant and compact modulo conjugation and let $u\in I(\mathfrak{g})$. Since $\partial(u)\widehat{j}(X,.)=p_u(-X)\widehat{j}(X,.)$, by \ref{eq 1.7.3} we have

\[\begin{aligned}
\displaystyle \int_{\Gamma(\mathfrak{g})} D^G(X)^{1/2} \lvert \theta(X)\rvert q_{L,u}(\widehat{j}(X,.))dX & \ll \int_{\Gamma(\mathfrak{g})} D^G(X)^{1/2} \lvert p_u(-X)\theta(X)\rvert dX \\
 & \leqslant q_{1,N}(\theta) \int_{\Gamma(\mathfrak{g})} \lvert p_u(-X)\rvert \lVert X\rVert_{\Gamma(\mathfrak{g})}^{-N}dX
\end{aligned}\]

\noindent for all $\theta\in SQC(\mathfrak{g}(\mathbb{R}))$ and all $N\geqslant 1$. There exists $N_0\geqslant 0$ such that $\lvert p_u(-X)\rvert\ll \lVert X\rVert_{\Gamma(\mathfrak{g})}^{N_0}$ for all $X\in \mathfrak{g}(F)$ and so by \ref{eq 1.8.2}, the last integral above is absolutely convergent for $N$ sufficiently large (depending on $u$). This proves the convergence and the continuity in $\theta\in SQC(\mathfrak{g}(\mathbb{R}))$ of the integral \ref{eq 4.2.11}.

\vspace{2mm}

\noindent We now prove the claim about differential operators. Let us denote by $\mathcal{A}$ the subalgebra of differential operators $D\in \Diff(\mathfrak{g})^G$ that induce a continuous endomorphism of $SQC(\mathfrak{g}(\mathbb{R}))$. It is obvious that $I(\mathfrak{g})\subset \mathcal{A}$. Since the Fourier transform exchanges the actions of $I(\mathfrak{g}^*)$ with the action of $I(\mathfrak{g})$, it follows that we also have $I(\mathfrak{g}^*)\subset \mathcal{A}$. Hence, by Proposition \ref{proposition 3.1.1}(vi), we have $\mathcal{A}=\Diff(\mathfrak{g})^G$.

\item By the closed graph theorem, we only need to prove that the bilinear map

$$\Diff^\infty(\omega)^G\times QC(\omega)\to QC(\omega)$$
$$(D,\theta)\mapsto D\theta$$

\noindent is separately continuous. Since for all $u\in I(\mathfrak{g})$ multiplication by $\partial(u)$ is a continuous endomorphism of $\Diff^\infty(\omega)^G$, it suffices to prove that the bilinear map above is separately continuous when the target space is equipped with the topology defined by the semi-norms $(q_{L,1})_L$. By definition of the topology on $\Diff^\infty(\omega)^G$, it is even sufficient to prove the following

\vspace{3mm}

\begin{num}
\item\label{eq 4.2.12} For each invariant and compact modulo conjugation subset $L\subseteq \omega$ and all $n\geqslant 1$, there exist a continuous semi-norm $\nu_{n,L}$ on $\Diff_{\leqslant n}^\infty(\omega)^G$ and a continuous semi-norm $\mu_{n,L}$ on $QC(\omega)$ such that
$$q_{L,1}(D\theta)\leqslant \nu_{n,L}(D)\mu_{n,L}(\theta)$$
for all $D\in \Diff_{\leqslant n}^\infty(\omega)^G$ and all $\theta\in QC(\omega)$.
\end{num}

\vspace{3mm}

\noindent By Proposition \ref{proposition 3.1.1}(v), to prove \ref{eq 4.2.12}, we only need to show that for all $D\in \Diff(\mathfrak{g})^G$ the linear map

$$QC(\omega)\to QC(\omega)$$
$$\theta\mapsto D\theta$$

\noindent is continuous. This follows from Lemma \ref{lemma 4.2.1}.

\item The density of $QC_c(\omega)$ in $QC(\omega)$ follows from the existence, for every invariant compact modulo conjugation subset $L\subseteq \omega$, of a function $\varphi_L\in C^\infty(\omega)^G$ that is compactly supported modulo conjugation and such that $\varphi_L=1$ on some neighborhood of $L$. Let now $\varphi\in C^\infty(\mathfrak{g}(\mathbb{R}))^G$ be compactly supported modulo conjugation and such that $\varphi=1$ in some neighborhood of $0$. Let us set $\varphi_t(X)=\varphi(t^{-1}X)$ for all $t>0$ and all $X\in \mathfrak{g}(\mathbb{R})$. The density of $QC_c(\mathfrak{g}(\mathbb{R}))$ in $SQC(\mathfrak{g}(\mathbb{R}))$ will follow from the following claim

\vspace{3mm}

\begin{num}
\item\label{eq 4.2.13} For all $\theta\in SQC(\mathfrak{g}(\mathbb{R}))$, we have
$$\displaystyle \lim\limits_{t\to\infty}\varphi_t\theta=\theta$$
in $SQC(\mathfrak{g}(\mathbb{R}))$.
\end{num}

\vspace{3mm}

\noindent We need to see that for all $u\in I(\mathfrak{g})$ and all integers $N\geqslant 0$, we have

$$\lim\limits_{t\to \infty}q_{u,N}(\theta-\varphi_t\theta)=0$$

\noindent Of course, it is sufficient to deal with the case where $u$ is homogeneous. We henceforth fix an element $u\in I(\mathfrak{g})$ which is homogeneous. Let $\omega\subset\mathfrak{g}(\mathbb{R})$ be a completely $G(\mathbb{R})$-invariant open neighborhood of $0$ on which $\varphi$ equals $1$. It is not hard to see that

$$\lVert X\rVert_{\Gamma(\mathfrak{g})}^{-1}\ll t^{-1}$$

\noindent for all $t>0$ and all $X\in \mathfrak{g}(\mathbb{R})-t\omega$. Since for all $t>0$ the function $\theta-\varphi_t\theta$ is supported in $\mathfrak{g}(\mathbb{R})-t\omega$, it follows that for all integers $N,k\geqslant 0$, we have an inequality

$$q_{u,N}(\theta-\varphi_t\theta)\ll t^{-k}q_{u,N+k}(\theta-\varphi_t\theta)$$

\noindent for all $t>0$. Hence, it suffices to show the existence of an integer $d_u$ such that for every integer $N\geqslant 1$, we have an inequality

\begin{align}\label{eq 4.2.14}
q_{u,N}(\varphi_t\theta)\ll t^{d_u}
\end{align}

\noindent for all $t>1$. We define an action of $\mathbb{R}_+^*$ on $\Diff^\infty(\mathfrak{g}(\mathbb{R}))$ as follows: for all $(t,D)\in \mathbb{R}_+^*\times \Diff^\infty(\mathfrak{g}(\mathbb{R}))$, we define $D_t\in \Diff^\infty(\mathfrak{g}(\mathbb{R}))$ by

$$D_tf=(Df_{t^{-1}})_t$$

\noindent for all $f\in C^\infty(\mathfrak{g}(\mathbb{R}))$, where as before $f_t(X)=f(t^{-1}X)$ for all $t\in \mathbb{R}_+^*$ and all $f\in C^\infty(\mathfrak{g}(\mathbb{R}))$. Let $N\geqslant 1$ be an integer. Then, we have

\begin{align}\label{eq 4.2.15}
q_{u,N}(\varphi_t\theta)=q_N(\partial(u)(\varphi_t\theta))=t^{-\deg(u)}q_N((\partial(u)\circ \varphi)_t\theta)
\end{align}

\noindent for all $t>0$, where we have set $q_N=q_{1,N}$ and where $\partial(u)\circ \varphi$ denotes the differential operator obtained by composing $\partial(u)$ with the multiplication by $\varphi$. By Proposition \ref{proposition 3.1.1}(v), there exists an integer $n\geqslant 1$, functions $\varphi_1,\ldots,\varphi_n\in C^\infty(\mathfrak{g}(\mathbb{R}))^G$ and operators $D_1,\ldots,D_n\in \Diff(\mathfrak{g})^G$ such that

$$\partial(u)\circ \varphi=\varphi_1D_1+\ldots+\varphi_nD_n$$

\noindent Of course, we may assume that the functions $\varphi_1,\ldots,\varphi_n$ are compactly supported modulo conjugation and that for all $1\leqslant i\leqslant n$ there exists an integer $d_i$ such that $(D_i)_t=t^{d_i}D_i$, for all $t\in\mathbb{R}_+^*$. It then easily follows from \ref{eq 4.2.15} that

$$\displaystyle q_{u,N}(\varphi_t\theta)\ll\sum_{i=1}^n t^{d_i-\deg(u)} q_N(D_i\theta)$$

\noindent for all $t>0$. Since the functions $D_i\theta$, $1\leqslant i\leqslant n$ are Schwartz quasi-characters, the semi-norms $q_N(D_i\theta)$ are finite. We deduce that \ref{eq 4.2.14} holds for $d_u=\max_{1\leqslant i\leqslant n}(d_i)-\deg(u)$. This proves \ref{eq 4.2.13} and ends the proof of (v). $\blacksquare$
\end{enumerate}

\subsection{Local expansions of quasi-characters on the Lie algebra when $F=\mathbb{R}$}\label{section 4.3}

\noindent Let $\omega\subseteq \mathfrak{g}(\mathbb{R})$ be a completely $G(\mathbb{R})$-invariant open subset and $\theta$ a quasi-character on $\omega$. Recall that, for $T\subseteq G$ a maximal torus, we denote by $\theta_T$ the function on $\omega\cap \mathfrak{t}_{\reg}(\mathbb{R})$ defined by

$$\theta_T(X)=D^G(X)^{1/2}\theta(X),\;\; X\in \omega\cap \mathfrak{t}_{\reg}(\mathbb{R})$$

\vspace{2mm}

\begin{lem}\label{lemma 4.3.1}
\begin{enumerate}[(i)]
\item Let $T\subseteq G$ be a maximal torus. Then, for any connected component $\Gamma\subseteq \omega\cap\mathfrak{t}_{\reg}(\mathbb{R})$ and all $u\in S(\mathfrak{t})$ the function $X\in \Gamma\mapsto (\partial(u)\theta_T)(X)$ extends continuously to the closure of $\Gamma$ in $\omega\cap \mathfrak{t}(\mathbb{R})$.

\item Let $X\in \omega_{\ssi}$. Then, there exist constants $c_{\theta,\mathcal{O}}(X)$ for $\mathcal{O}\in \Nil_{\reg}(\mathfrak{g}_X)$ such that

$$\displaystyle D^{G}(X+Y)^{1/2}\theta(X+Y)=D^{G}(X+Y)^{1/2}\sum_{\mathcal{O}\in \Nil_{\reg}(\mathfrak{g}_X)} c_{\theta,\mathcal{O}}(X)\widehat{j}(\mathcal{O},Y)+O(\lvert Y\rvert)$$

\noindent for all $Y\in \mathfrak{g}_{X,\reg}(\mathbb{R})$ sufficiently near $0$.
\end{enumerate}
\end{lem}

\vspace{2mm}

\noindent\ul{Proof}:

\begin{enumerate}[(i)]
\item This follows directly from Lemma \ref{lemma 4.2.1} and the mean value theorem (Note that $\mathfrak{t}_{\reg}(\mathbb{R})$ is the complement in $\mathfrak{t}(\mathbb{R})$ of a finite union of subspaces).

\item By Lemma \ref{lemma 4.2.2}(ii), we are immediately reduced to the case where $X\in \mathfrak{z}_G(\mathbb{R})$. By translation, we may even assume that $X=0$. Define a function $\theta_0$ on $\mathfrak{g}_{\reg}(\mathbb{R})$ by

$$D^G(X)^{1/2}\theta_0(X)=\lim\limits_{t\to 0^+} D^G(tX)^{1/2}\theta(tX)$$

\noindent for all $X\in \mathfrak{g}_{\reg}(\mathbb{R})$. Note that for $t>0$ sufficiently small, we have $tX\in \omega_{\reg}$ and that the limit exists by (i). Moreover, the function $\theta_0$ is invariant and homogeneous of degree $-\delta(G)/2$, that is

\begin{align}\label{eq 4.3.1}
\theta_0(tX)=t^{-\delta(G)/2}\theta_0(X)
\end{align}

\noindent for all $X\in \mathfrak{g}_{\reg}(\mathbb{R})$ and all $t>0$. Also, the function $(D^G)^{1/2}\theta_0$ is (globally) bounded since $(D^G)^{1/2}\theta$ is bounded near $0$. By (i) again, for any maximal torus $T\subset G$ and any connected component $\Gamma\subseteq \mathfrak{g}_{\reg}(\mathbb{R})$, we have

$$D^G(X)^{1/2}\theta(X)=D^G(X)^{1/2}\theta_0(X)+O(\lvert X\rvert)$$

\noindent as $X\in \Gamma\cap \omega$ goes to $0$. Since there are only finitely many conjugacy classes of maximal tori, that for all of them $\mathfrak{t}_{\reg}(\mathbb{R})$ has only finitely many connected components and that $\lvert X\rvert\ll \lvert g^{-1}Xg\rvert$ for all $g\in G(\mathbb{R})$ and all $X\in\mathfrak{t}(\mathbb{R})$, it follows that for some completely $G(\mathbb{R})$-invariant neighborhood $\omega'\subseteq \omega$ of $0$, we have

\begin{align}\label{eq 4.3.2}
D^G(X)^{1/2}\theta(X)=D^G(X)^{1/2}\theta_0(X)+O(\lvert X\rvert)
\end{align}

\noindent for all $X\in \omega'_{\reg}$. It remains to show that the function $\theta_0$ is a linear combination of functions $\widehat{j}(\mathcal{O},.)$ for $\mathcal{O}\in \Nil_{\reg}(\mathfrak{g})$. Since the function $(D^G)^{1/2}\theta_0$ is bounded, the function $\theta_0$ is locally integrable and so it defines a distribution $T_{\theta_0}$. Let $I^+(\mathfrak{g})\subset I(\mathfrak{g})$ denote the subalgebra of elements without constant term. We first establish the following

\vspace{3mm}

\begin{num}
\item\label{eq 4.3.3} For all $u\in I^+(\mathfrak{g})$, we have $\partial(u)T_{\theta_0}=0$.
\end{num}

\vspace{3mm}

\noindent Let $u\in I^+(\mathfrak{g})$ be homogeneous of degree $d>0$. It follows easily from \ref{eq 4.3.1} that the distribution $\partial(u)T_{\theta_0}$ is homogeneous of degree $-d-\delta(G)/2$ in the following sense: for all $f\in C_c^\infty(\mathfrak{g}(\mathbb{R}))$, we have

\begin{align}\label{eq 4.3.4}
(\partial(u)T_{\theta_0})(f_t)=t^{\dim(\mathfrak{g})-d-\delta(G)/2}(\partial(u)T_{\theta_0})(f)
\end{align}

\noindent for all $t>0$ and where $f_t(X)=f(t^{-1}X)$. Since the function $(D^G)^{1/2}\partial(u)\theta$ is locally bounded and $(D^G)^{-1/2}$ locally integrable, it is easy to see that for all $f\in C_c^\infty(\mathfrak{g}(\mathbb{R}))$, we have

\begin{align}\label{eq 4.3.5}
\left\lvert (\partial(u)T_{\theta})(f_t)\right\rvert=\left\lvert T_{\partial(u)\theta}(f_t)\right\rvert\ll t^{\dim(\mathfrak{g})-\delta(G)/2}
\end{align}

\noindent for all $t>0$ sufficiently small (in particular so that $\Supp(f_t)\subseteq \omega$). Set $R=T_{\theta-\theta_0}$. It is a distribution on $\omega$. By \ref{eq 4.3.2}, we know that $X\mapsto \lvert X\rvert^{-1} D^G(X)^{1/2}\lvert \theta(X)-\theta_0(X)\rvert$ is bounded on some neighborhood of $0$. It follows easily that for all $f\in C_c^\infty(\mathfrak{g}(\mathbb{R}))$, we have

\begin{align}\label{eq 4.3.6}
\lvert (\partial(u)R)(f_t)\rvert\ll t^{1+\dim(\mathfrak{g})-d-\delta(G)/2}
\end{align}

\noindent for all $t>0$ sufficiently small. Since we have $\partial(u)R=\partial(u)T_\theta-\partial(u)T_{\theta_0}$, the equality (\ref{eq 4.3.4} and the inequalities \ref{eq 4.3.5} and \ref{eq 4.3.6} cannot be compatible unless $\partial(u)T_{\theta_0}=0$. This proves \ref{eq 4.3.3}.

\vspace{2mm}

By Lemma 2.2 of \cite{BV}, every homogeneous distribution is tempered. Hence, $T_{\theta_0}$ is tempered. Consider its Fourier transform $\widehat{T_{\theta_0}}$. It is an invariant and homogeneous distribution of degree $-\dim(\mathfrak{g})-\delta(G)/2$. Also, by \ref{eq 4.3.3}, we have $p_u\widehat{T_{\theta_0}}=0$ for all $u\in I^+(\mathfrak{g})$. Since the common zero locus of the polynomials $p_u$ for all $u\in I^+(\mathfrak{g})$ is the nilpotent cone $\mathcal{N}$ of $\mathfrak{g}(\mathbb{R})$, it follows that $\Supp(\widehat{T_{\theta_0}})\subseteq \mathcal{N}$. Hence, by \ref{eq 1.7.2}, $\widehat{T}_{\theta_0}$ is a linear combination of the distributions $J_{\mathcal{O}}$ for $\mathcal{O}\in \Nil_{\reg}(\mathfrak{g})$ and we are done. $\blacksquare$

\end{enumerate}

\subsection{Quasi-characters on the group when $F=\mathbb{R}$}\label{section 4.4}

\noindent In this section we still assume that $F=\mathbb{R}$. Let $\Omega\subseteq G(\mathbb{R})$ be a completely $G(\mathbb{R})$-invariant open subset. A {\em quasi-character} on $\Omega$ is a function $\theta\in C^\infty(\Omega_{\reg})^G$ that satisfies the two following conditions

\vspace{2mm}

\begin{itemize}
\renewcommand{\labelitemi}{$\bullet$}
\item For all $z\in \mathcal{Z}(\mathfrak{g})$, the function $(D^G)^{1/2}z\theta$ is locally bounded on $\Omega$ (so that by \ref{eq 1.8.1} the function $z\theta$ is locally integrable on $\Omega$);

\item For all $z\in \mathcal{Z}(\mathfrak{g})$, we have the following equality of distributions on $\Omega$

$$zT_\theta=T_{z\theta}$$
\end{itemize}

\noindent We say that a quasi-character $\theta$ is {\em compactly supported} if it is compactly supported modulo conjugation. We will denote by $\gls{QCOmega}$ the space of quasi-characters on $\Omega$ and by $\gls{QCcOmega}$ the subspace of compactly supported quasi-characters. If $L\subset \Omega$ is invariant and compact modulo conjugation we also introduce the subspace $QC_L(\Omega)\subset QC_c(\Omega)$ of quasi-characters supported in $L$. We endow $QC_L(\Omega)$ with the topology defined by the semi-norms

$$\displaystyle \gls{qLz}(\theta)=\sup_{x\in L_{\reg}} D^G(x)^{1/2}\lvert z\theta(x)\rvert$$

\noindent for all $z\in \mathcal{Z}(\mathfrak{g})$ and we equip $QC(\Omega)$ with the topology defined by the semi-norms $(q_{L,z})_{L,z}$ where $L$ runs through the invariant compact modulo conjugation subsets of $\Omega$ and $z$ runs through $\mathcal{Z}(\mathfrak{g})$. Finally, we put on $QC_c(\Omega)$ the inductive limit topology relative to the natural isomorphism

$$QC_c(\Omega)\simeq \varinjlim_{L}QC_L(\Omega)$$

\vspace{2mm}

\noindent As in the case of the Lie algebra, by the representation theorem of Harish-Chandra, any invariant distribution $T$ defined on some completely $G(\mathbb{R})$-invariant open subset $\Omega\subseteq G(\mathbb{R})$ such that $\dim(\mathcal{Z}(\mathfrak{g})T)<\infty$ is the distribution associated to a quasi-character on $\Omega$. In particular, for every admissible irreducible representation $\pi$ of $G(\mathbb{R})$ the character $\theta_\pi$ of $\pi$ is a quasi-character on $G(\mathbb{R})$.

\vspace{2mm}

\begin{prop}\label{proposition 4.4.1}
Let $\Omega\subseteq G(\mathbb{R})$ be a completely $G(\mathbb{R})$-invariant open subset and let $L\subseteq \Omega$ be invariant and compact modulo conjugation. Then

\begin{enumerate}[(i)]
\item Let $\omega\subseteq \mathfrak{g}(\mathbb{R})$ be a $G$-excellent open subset and assume that $\Omega=\exp(\omega)$. Then the linear map

$$\theta\mapsto \theta_\omega$$

\noindent induces topological isomorphisms $QC(\Omega)\simeq QC(\omega)$ and $QC_c(\Omega)\simeq QC_c(\omega)$.

\item $QC(\Omega)$ and $QC_L(\Omega)$ are nuclear Fr\'echet spaces and $QC_c(\Omega)$ is a nuclear LF space.

\item Let $x\in G_{\ssi}(\mathbb{R})$ and $\Omega_x\subseteq G_x(\mathbb{R})$ be a $G$-good open neighborhood of $x$. Assume that $\Omega=\Omega_x^G$. Then the linear map

$$\theta\mapsto \theta_{x,\Omega_x}$$

\noindent induces topological isomorphisms $QC(\Omega)\simeq QC(\Omega_x)^{Z_G(x)(F)}$ and $QC_c(\Omega)\simeq QC_c(\Omega_x)^{Z_G(x)(F)}$.

\item For all $D\in \Diff^\infty(\Omega)^G$ and all $\theta\in QC(\Omega)$, the function $D\theta$ is a quasi-character on $\Omega$ and we have $DT_\theta=T_{D\theta}$. Moreover, the two bilinear maps

$$\Diff^\infty(\Omega)^G\times QC(\Omega)\to QC(\Omega),\;\;\; \Diff^\infty(\Omega)^G\times QC_c(\Omega)\to QC_c(\Omega)$$
$$(D,\theta)\mapsto D\theta$$

\noindent are separately continuous.

\item If $G=G_1\times G_2$ with $G_1$ and $G_2$ two reductive connected groups over $F$, and $\Omega=\Omega_1\times \Omega_2$ where $\Omega_1\subseteq G_1(\mathbb{R})$ (resp.\ $\Omega_2\subseteq G_2(\mathbb{R})$) is a completely $G_1(\mathbb{R})$-invariant (resp.\ completely $G_2(\mathbb{R})$-invariant) open subset then there is a canonical isomorphism of topological vector spaces

$$QC(\Omega)\simeq QC(\Omega_1)\widehat{\otimes}_p QC(\Omega_2)$$

\item For all $\theta\in QC(\Omega)$ and all $x\in G_{\ssi}(\mathbb{R})$, there exist constants $\gls{cthetaOx}\in \mathbb{C}$, for $\mathcal{O}\in \Nil_{\reg}(\mathfrak{g}_x)$, such that

$$\displaystyle D^G(xe^Y)^{1/2}\theta(xe^Y)=D^G(xe^Y)^{1/2}\sum_{\mathcal{O}\in \Nil_{\reg}(\mathfrak{g}_x)}c_{\theta,\mathcal{O}}(x)\widehat{j}(\mathcal{O},Y)+O(\lvert Y\rvert)$$

\noindent for all $Y\in \mathfrak{g}_{x,\reg}(\mathbb{R})$ sufficiently near $0$.

\end{enumerate}
\end{prop}

\vspace{2mm}

\noindent\ul{Proof}: (i) is a straightforward consequence of \ref{eq 3.3.3}. (vi) follows from (i), (iii) and Lemma \ref{lemma 4.3.1}(ii). Before proving (ii)-(v), we need to show the following

\vspace{5mm}

\begin{num}
\item\label{eq 4.4.1} Let $J\subseteq \mathcal{Z}(\mathfrak{g})$ be a subalgebra such that the extension $\mathcal{Z}(\mathfrak{g})/J$ is finite. Assume that $\theta\in C^\infty(\Omega_{\reg})^G$ satisfies the two following conditions

\begin{itemize}
\renewcommand{\labelitemi}{$\bullet$}
\item For all $z\in J$, the function $(D^G)^{1/2}z\theta$ is locally bounded on $\Omega$;

\item For all $z\in J$, we have the equality of distributions on $\Omega$

$$zT_{\theta}=T_{z\theta}$$
\end{itemize}

Then, $\theta$ is a quasi-character on $\Omega$.
\end{num}

\vspace{5mm}

\noindent The proof is by induction on $\dim(G)$, the case of a torus being obvious. Let $\theta$ and $J$ be as in \ref{eq 4.4.1}. We need to show that $\theta$ is a quasi-character near every semi-simple point $x\in \Omega_{\ssi}$. If $G_x\neq G$, then we can use semi-simple descent (cf.\ Lemma \ref{lemma 3.2.1}(ii)) and the induction hypothesis for $G_x$ to conclude. Assume that $G_x=G$ i.e., $x\in Z_G(F)$, then translating $\theta$ and $\Omega$ by $x$, we may assume that $x=1$. But then the result follows from (i), \ref{eq 3.3.3} and the analogous result for the Lie algebra (Proposition \ref{proposition 4.2.1}.2.(ii)).

\vspace{2mm}

\noindent We may now proceed to the proof of (ii)-(v).

\vspace{2mm}

\begin{enumerate}[(i)]
\setcounter{enumi}{1}
\item Using Proposition \ref{proposition 3.1.1}(i), it is easy to see that the topology on $QC(\Omega)$ is defined by a countable number of semi-norms and that $QC_c(\Omega)$ is the direct limit of a countable family $(QC_{L_n}(\Omega))_{n\geqslant 1}$ where $(L_n)_{n\geqslant 1}$ is an increasing sequence of invariant and compact modulo conjugation subsets of $\Omega$. Moreover, $QC_L(\Omega)$ is a closed subspace of $QC(\Omega)$ and its topology is induced from the one on $QC(\Omega)$. Hence, we only need to show that $QC(\Omega)$ is nuclear and complete. For every $x\in \Omega_{\ssi}$, choose $\omega_x\subseteq \mathfrak{g}_x(F)$ a $G_x$-excellent open subset such that $\Omega_x=x\exp(\omega_x)\subseteq G_x(F)$ is a $G$-good open neighborhood of $x$ such that $\Omega_x^G\subseteq \Omega$. The linear map

$$\displaystyle QC(\Omega)\to \prod_{x\in \Omega_{\ssi}} QC(\Omega_x^G)$$
$$\displaystyle \theta\mapsto \left(\theta_{\mid \Omega_x^G}\right)_{x\in \Omega_{\ssi}}$$

\noindent is a closed embedding. Hence, it suffices to show that the spaces $QC(\Omega_x^G)$, $x\in \Omega_{\ssi}$, are nuclear and complete. Let $x\in \Omega_{\ssi}$ and consider the map

\begin{align}\label{eq 4.4.2}
\theta\in QC(\Omega_x^G)\mapsto \left(\theta_{x,\Omega_x}\right)_{\omega_x}\in C^\infty(\omega_{x,\reg})^{Z_G(x)}
\end{align}

\noindent Set $J=\{u_{z_x};\; z\in \mathcal{Z}(\mathfrak{g})\}$. It is a subalgebra of $I(\mathfrak{g}_x)$ with the property that the extension $I(\mathfrak{g}_x)/J$ is finite. By \ref{eq 4.4.1}, (i) and Proposition \ref{proposition 4.2.1}.2.(ii), we see that the linear map \ref{eq 4.4.2} induces a linear isomorphism $QC(\Omega_x^G)\simeq QC(\omega_x)^{Z_G(x)}$. By Lemma \ref{lemma 4.2.1}, it is even a topological isomorphism. Hence, by Lemma \ref{lemma 4.2.2}(i), $QC(\Omega_x^G)$ is a nuclear Fr\'echet space and this ends the proof of (ii).

\item Once again, by \ref{eq 4.4.1} (where we take $J=\{z_x,\; z\in \mathcal{Z}(\mathfrak{g})\}$) and Lemma \ref{lemma 3.2.1}(i), the linear map $\theta\mapsto \theta_{x,\Omega_x}$ induces linear isomorphisms $QC(\Omega)\simeq QC(\Omega_x)^{Z_G(x)}$ and $QC_c(\Omega)\simeq QC_c(\Omega_x)^{Z_G(x)}$. The inverse of these isomorphisms are obviously continuous. Hence, by the open-mapping theorem these are topological isomorphisms.

\item The first part of (iv) follows easily from (i), (iii) and Proposition \ref{proposition 4.2.1}.2.(iii). Choose, as in the proof of (i), for every $x\in \Omega_{\ssi}$ a $G_x$-excellent open subset $\omega_x\subseteq \mathfrak{g}_x(F)$ such that $\Omega_x=x\exp(\omega_x)\subseteq G_x(F)$ is $G$-good and $\Omega_x^G\subseteq \Omega$. Then, by the closed graph theorem, it suffices to show that for every $x\in \Omega_{\ssi}$ the bilinear map

$$\Diff^\infty(\Omega_x^G)^G\times QC(\Omega_x^G)\to QC(\Omega_x^G)$$
$$(D,\theta)\mapsto D\theta$$

\noindent is separately continuous. But by (i), (iii) and Lemma \ref{lemma 3.2.1}(ii), we have topological isomorphisms

$$QC(\Omega_x^G)\simeq QC(\omega_x)^{Z_G(x)}$$
$$\theta\mapsto \left(R(x)\theta_{x,\Omega_x}\right)_{\omega_x}$$

$$\overline{\Diff^\infty(\Omega_x^G)^G}\simeq \overline{\Diff^\infty(\omega_x)^{Z_G(x)}}$$
$$D\mapsto \left(R(x)D_{x,\Omega_x}\right)_{\omega_x}$$

\noindent Hence, we are reduced to show that the bilinear map

$$\Diff^\infty(\omega_x)^{Z_G(x)} \times QC(\omega_x)^{Z_G(x)}\to QC(\omega_x)^{Z_G(x)}$$
$$(D,\theta)\mapsto D\theta$$

\noindent is separately continuous. This follows from Lemma \ref{lemma 4.2.2}(iv).

\item The natural bilinear map

$$QC(\Omega_1)\times QC(\Omega_2)\to QC(\Omega)$$
$$(\theta_1,\theta_2)\mapsto \left[ (g_1,g_2)\mapsto \theta_1(g_1)\theta_2(g_2)\right]$$

\noindent is continuous. Hence it extends to a continuous linear map

$$QC(\Omega_1)\widehat{\otimes}_pQC(\Omega_2)\to QC(\Omega)$$

\noindent and we would like to show that this is a topological isomorphism. By the open mapping theorem and Proposition \ref{proposition A.5.1} of the appendix, this amounts to proving that for all $\theta\in QC(\Omega)$ the two following conditions are satisfied

\vspace{3mm}

\begin{num}
\item\label{eq 4.4.3} For all $g_1\in \Omega_{1,\reg}$, the function $g_2\in \Omega_{2,\reg}\mapsto \theta(g_1,g_2)$ belongs to $QC(\Omega_2)$;

\item\label{eq 4.4.4} For all $\lambda\in QC(\Omega_2)'$, the function $g_1\mapsto \lambda(\theta(g_1,.))$ belongs to $QC(\Omega_1)$.
\end{num}

\vspace{3mm}

\noindent The first condition is easy to check and left to the reader. Let $\theta\in QC(\Omega)$. In order to prove that the condition \ref{eq 4.4.4} is satisfied, it is obviously sufficient to establish the following

\vspace{3mm}

\begin{num}
\item\label{eq 4.4.5} The function $g_1\in \Omega_{1,\reg}\mapsto \theta(g_1,.)\in QC(\Omega_2)$ is smooth, for all $z_1\in \mathcal{Z}(\mathfrak{g}_1)$ and every invariant and compact modulo conjugation subset $L_1\subseteq \Omega_1$, the set 
$$\{D^{G_1}(g_1)^{1/2}\left(R_1(z_1)\theta\right)(g_1,.),\; g_1\in L_{1,\reg}\}$$
is bounded in $QC(\Omega_2)$ and for all $\varphi_1\in C_c^\infty(\Omega_1)$, we have the equality
$$\displaystyle \int_{\Omega_1}\left(R_1(z_1)\theta\right)(g_1,.)\varphi_1(g_1) dg_1=\int_{\Omega_1} \theta(g_1,.) \left(z_1^*\varphi_1\right)(g_1)dg_1$$
in $QC(\Omega_2)$.
\end{num}

\vspace{3mm}

\noindent (the index $1$ in $R_1(z_1)$ is here to emphasize that we are deriving in the first variable). Note that if the function $g_1\in \Omega_{1,\reg}\mapsto \theta(g_1,.)\in QC(\Omega_2)$ is smooth then for all $u\in \mathcal{U}(\mathfrak{g})$, the derivative $R(u)\left(g_1\mapsto \theta(g_1,.)\right)(g_1)\in QC(\Omega_2)$ is necessarily equal to $\left(R_1(u)\theta\right)(g_1,.)$, which is why we are using the function $\left(R_1(z_1)\theta\right)(g_1,.)$ above. The second claim in \ref{eq 4.4.5} is obvious and the last equality of \ref{eq 4.4.5} need only to be checked after applying to it the continuous linear forms

$$\displaystyle \theta_2\in QC(\Omega_2)\mapsto \int_{\Omega_2}\theta_2(g_2)\varphi_2(g_2)dg_2,\;\;\; \varphi_2\in C_c^\infty(\Omega_{2,\reg})$$

\noindent (because these linear forms separate elements of $QC(\Omega_2)$) where then it is an obvious consequence of $\theta$ being a quasi-character. Hence, we are only left with proving that the map $g_1\in \Omega_{1,\reg}\mapsto \theta(g_1,.)\in QC(\Omega_2)$ is smooth. This fact follows easily from the next claim

\vspace{3mm}

\begin{num}
\item\label{eq 4.4.6} For all $u\in \mathcal{U}(\mathfrak{g}_1)$ and all $g_1\in \Omega_{1,\reg}$, the function $(R_1(u)\theta)(g_1,.)$ is a quasi-character on $\Omega_2$ and for every compact subset $\mathcal{K}_1\subseteq \Omega_{1,\reg}$, the family
$$\{(R_1(u)\theta)(g_1,.),\; g_1\in \mathcal{K}_1\}$$
is bounded in $QC(\Omega_2)$.
\end{num}

\vspace{3mm}

\noindent Let $u\in \mathcal{U}(\mathfrak{g}_1)$. In order to get \ref{eq 4.4.6}, we only need to check the following

\vspace{3mm}

\begin{num}
\item\label{eq 4.4.7} For all $z_2\in \mathcal{Z}(\mathfrak{g}_2)$, all compact subsets $\mathcal{K}_1\subseteq \Omega_{1,\reg}$ and each invariant and compact modulo conjugation subset $L_2\subseteq \Omega_2$, there exists $C>0$ such that
$$\displaystyle D^{G_2}(g_2)^{1/2}\left\lvert \left(R_2(z_2)R_1(u)\theta\right)(g_1,g_2)\right\rvert\leqslant C$$
for all $(g_1,g_2)\in \mathcal{K}_1\times L_{2,\reg}$.
\end{num}

\vspace{3mm}

\noindent and

\vspace{3mm}

\begin{num}
\item\label{eq 4.4.8} For all $\varphi_2\in C_c^\infty(\Omega_2)$ and for all $z_2\in \mathcal{Z}(\mathfrak{g}_2)$, we have the equality
$$\displaystyle \int_{\Omega_2}\left(R_2(z_2)R_1(u)\theta\right)(g_1,g_2) \varphi_2(g_2) dg_2=\int_{\Omega_2}\left(R_1(u)\theta\right)(g_1,g_2)\left(R(z_2^*)\varphi_2\right)(g_2)dg_2$$
for all $g_1\in \Omega_{1,\reg}$.
\end{num}

\vspace{3mm}

\noindent If \ref{eq 4.4.7} is satisfied, then both sides of the equality \ref{eq 4.4.8} are continuous in $g_1\in \Omega_{1,\reg}$. Hence, the equality need only to be checked after integrating it against a function $\varphi_1\in C_c^\infty(\Omega_{1,\reg})$ where it is again an easy consequence of $\theta$ being a quasi-character. Let us prove \ref{eq 4.4.7}. Fix $z_2\in \mathcal{Z}(\mathfrak{g}_2)$ and an invariant and compact modulo conjugation subset $L_2\subseteq \Omega_2$. Then the functions $\left(R_2(z_2)\theta\right)(.,g_2)$, $g_2\in L_{2,\reg}$, are all quasi-characters on $\Omega_1$ and the family

$$\{D^{G_2}(g_2)^{1/2}\left(R_2(z_2)\theta\right)(.,g_2), \; g_2\in L_{2,\reg}\}$$

\noindent is bounded in $QC(\Omega_1)$. Hence, to get \ref{eq 4.4.7} it suffices to see that for every compact subset $\mathcal{K}_1\subseteq \Omega_{1,\reg}$, the linear forms

$$\theta_1\in QC(\Omega_1)\mapsto \left(R(u)\theta_1\right)(g_1),\;\; g_1\in \mathcal{K}_1$$

\noindent form a bounded subset of $QC(\Omega_1)'$. This follows from example from (iii) (where we take points $x\in \Omega_1$ which are regular). This ends the proof of (v). $\blacksquare$
\end{enumerate}

\subsection{Functions $c_{\theta}$}\label{section 4.5}

\noindent We henceforth drop the condition that $F=\mathbb{R}$ so that $F$ can be $p$-adic as well. Let $\theta$ be a quasi-character on $G(F)$. Then, for all $x\in G_{\ssi}(F)$ we have a local expansion

$$\displaystyle D^G(xe^X)^{1/2}\theta(xe^X)=D^{G}(xe^X)^{1/2}\sum_{\mathcal{O}\in \Nil_{\reg}(\mathfrak{g}_x)} c_{\theta,\mathcal{O}}(x)\widehat{j}(\mathcal{O},X)+O(\lvert X\rvert)$$

\noindent for all $X\in \mathfrak{g}_{x,\reg}(F)$ sufficiently near $0$ (in the $p$-adic case, this follows from the fact that $D^G(X)^{1/2}\widehat{j}(\mathcal{O},X)=O(\lvert X\rvert)$ near $0$ for all $\mathcal{O}\in \Nil(\mathfrak{g})\backslash \Nil_{\reg}(\mathfrak{g})$). It follows from the homogeneity property of the functions $\widehat{j}(\mathcal{O},.)$ and their linear independence that the coefficients $c_{\theta,\mathcal{O}}(x)$, $\mathcal{O}\in \Nil_{\reg}(\mathfrak{g}_x)$, are uniquely defined. We set

$$\displaystyle \gls{ctheta}(x)=\frac{1}{\lvert \Nil_{\reg}(\mathfrak{g}_x)\rvert}\sum_{\mathcal{O}\in \Nil_{\reg}(\mathfrak{g}_x)}c_{\theta,\mathcal{O}}(x)$$

\noindent for all $x\in G_{\ssi}(F)$. This defines a function

$$\displaystyle c_\theta: G_{\ssi}(F)\to \mathbb{C}$$

\noindent Similarly, to any quasi-character $\theta$ on $\mathfrak{g}(F)$ we associate a function

$$c_\theta: \mathfrak{g}_{\ssi}(F)\to \mathbb{C}$$

\vspace{2mm}

\begin{prop}\label{proposition 4.5.1}
\begin{enumerate}
\item Let $\theta$ be a quasi-character on $G(F)$ and let $x\in G_{\ssi}(F)$. Then

\begin{enumerate}[(i)]

\item If $G_x$ is not quasi-split then $c_{\theta}(x)=0$.

\item Assume that $G_x$ is quasi-split. Let $B_x\subset G_x$ be a Borel subgroup and $T_{\quasid,x}\subset B_x$ be a maximal torus (both defined over $F$). Then, we have

$$D^G(x)^{1/2}c_\theta(x)=\lvert W(G_x,T_{\quasid,x})\rvert^{-1}\lim\limits_{x'\in T_{\quasid,x}(F)\to x} D^G(x')^{1/2}\theta(x')$$

\noindent (in particular, the limit exists).

\item The function $(D^G)^{1/2}c_{\theta}$ is locally bounded on $G(F)$. More precisely, for any invariant and compact modulo conjugation subset $L\subset G(F)$, there exists a continuous semi-norm $\nu_L$ on $QC(G(F))$ such that

$$\sup_{x\in L_{\ssi}} D^G(x)^{1/2}\lvert c_\theta(x)\rvert \leqslant \nu_L(\theta)$$

\noindent for all $\theta\in QC(G(F))$.

\item Let $\Omega_x\subseteq G_x(F)$ be a $G$-good open neighborhood of $x$. Then, we have

$$D^G(y)^{1/2}c_\theta(y)=D^{G_x}(y)^{1/2}c_{\theta_{x,\Omega_x}}(y)$$

\noindent for all $y\in \Omega_{x,ss}$.
\end{enumerate}

\item Let $\theta$ be a quasi-character on $\mathfrak{g}(F)$ and let $X\in \mathfrak{g}_{\ssi}(F)$. Then

\begin{enumerate}[(i)]

\item If $G_X$ is not quasi-split then $c_{\theta}(X)=0$.

\item Assume that $G_X$ is quasi-split. Let $B_X\subset G_X$ be a Borel subgroup and $T_{\quasid,X}\subset B_X$ be a maximal torus (both defined over $F$). Then, we have

$$D^G(X)^{1/2}c_\theta(X)=\lvert W(G_X,T_{\quasid,X})\rvert^{-1}\lim\limits_{X'\in \mathfrak{t}_{\quasid,X}(F)\to X} D^G(X')^{1/2}\theta(X')$$

\noindent (in particular, the limit exists).

\item The function $(D^G)^{1/2}c_{\theta}$ is locally bounded on $\mathfrak{g}(F)$. More precisely, for any invariant and compact modulo conjugation subset $L\subset \mathfrak{g}(F)$, there exists a continuous semi-norm $\nu_L$ on $QC(\mathfrak{g}(F))$ such that

$$\sup_{X\in L_{\ssi}} D^G(X)^{1/2}\lvert c_\theta(X)\rvert \leqslant \nu_L(\theta)$$

\noindent for all $\theta\in QC(\mathfrak{g}(F))$.

\item For all $\lambda\in F^\times$ let $\gls{Mlambda} \theta$ be the quasi-character defined by $(M_\lambda\theta)(X)=\lvert \lambda\rvert^{-\delta(G)/2}\theta(\lambda^{-1}X)$ for all $X\in \mathfrak{g}_{\reg}(F)$. Then we have

$$D^G(X)^{1/2}c_{M_\lambda \theta}(X)=D^G(\lambda^{-1}X)^{1/2}c_\theta(\lambda^{-1}X)$$

\noindent for all $X\in \mathfrak{g}_{\ssi}(F)$ and all $\lambda\in F^\times$.

\item Assume that $G$ is quasi-split. Let $B\subset G$ be a Borel subgroup and $T_{\quasid}\subset B$ be a maximal torus (both defined over $F$). Then, for all $X\in \mathfrak{t}_{\quasid,\reg}(F)$, we have

$$c_{\widehat{j}(X,.)}(0)=1$$

\end{enumerate}
\end{enumerate}
\end{prop}

\vspace{2mm}

\noindent\ul{Proof}: 1.(i) and 2.(i) are obvious since for $G_x$ not quasi-split, $\Nil_{\reg}(\mathfrak{g}_x)$ is empty. 1.(ii) and 2.(ii) follow easily from \ref{eq 3.4.7} whereas 1.(iii), 1.(iv), 2.(iii) and 2.(iv) are direct consequences of the preceding points and the fact that the function $(D^G)^{1/2}\theta$ is locally bounded by a continuous semi-norm on $QC(G(F))$ (resp.\ on $QC(\mathfrak{g}(F))$) for all $\theta\in QC(G(F))$ (resp.\ for all $\theta\in QC(\mathfrak{g}(F))$). Finally, 2.(v) follows from 2.(ii) and \ref{eq 3.4.5}. $\blacksquare$

\subsection{Homogeneous distributions on spaces of quasi-characters}\label{section 4.6}

\noindent For all $\lambda\in F^\times$, let us denote by $M_\lambda$ the continuous operator on $QC_c(\mathfrak{g}(F))$ (resp.\ on $SQC(\mathfrak{g}(F))$, resp.\ on $QC(\mathfrak{g}(F))$) given by

$$M_\lambda \theta=\lvert \lambda\rvert^{-\delta(G)/2}\theta_\lambda,\;\;\; \theta\in QC_c(\mathfrak{g}(F))\; (\mbox{resp.\ } \theta\in SQC(\mathfrak{g}(F)), \mbox{ resp.\ } \theta\in QC(\mathfrak{g}(F)))$$

\noindent (recall that $\theta_\lambda(X)=\theta(\lambda^{-1}X)$ for all $X\in\mathfrak{g}_{\reg}(F)$).

\begin{prop}\label{proposition 4.6.1}
Let $\lambda\in F^\times$ be such that $\lvert \lambda\rvert\neq 1$. Then, we have the following:
\begin{enumerate}[(i)]
\item For all $\theta\in QC_c(\mathfrak{g}(F))$ (resp.\ $\theta\in SQC(\mathfrak{g}(F))$) such that $c_{\theta,\mathcal{O}}(0)=0$ for all $\mathcal{O}\in \Nil_{\reg}(\mathfrak{g})$ and for any integer $d\geqslant 1$, there exist $\theta_1,\theta_2\in QC_c(\mathfrak{g}(F))$ (resp.\ $\theta_1,\theta_2\in SQC(\mathfrak{g}(F))$) such that

\begin{itemize}
\renewcommand{\labelitemi}{$\bullet$}

\item $\theta=(M_\lambda-1)^d\theta_1+\theta_2$;

\item $0\notin \Supp(\theta_2)$.
\end{itemize}

\item Let $\ell$ be a continuous linear form on $QC_c(\mathfrak{g}(F))$ such that

$$\ell(M_\lambda \theta)=\ell(\theta)$$

\noindent for all $\theta\in QC_c(\mathfrak{g}(F))$. Then $\ell$ extends by continuity to $SQC(\mathfrak{g}(F))$.
\end{enumerate}
\end{prop}

\vspace{2mm}

\noindent\ul{Proof}: 

\begin{enumerate}[(i)]
\item Let $\lambda\in F^\times$ such that $\lvert \lambda\rvert\neq 1$. Since $M_\lambda M_{\lambda^{-1}}=Id$, we may assume that $\lvert \lambda\rvert >1$. Denote by $QC_0(\mathfrak{g}(F))$ the space of quasi-characters $\theta\in QC(\mathfrak{g}(F))$ such that $c_{\theta,\mathcal{O}}(0)=0$ for all $\mathcal{O}\in \Nil_{\reg}(\mathfrak{g}(F))$. It is a closed subspace of $QC(\mathfrak{g}(F))$. Obviously, we only need to prove that $(M_\lambda-1)$ is a linear bijection of $QC_0(\mathfrak{g}(F))$ onto itself. This will follow from the next claim

\vspace{3mm}

\begin{num}
\item\label{eq 4.6.1} For all $\theta\in QC_0(\mathfrak{g}(F))$, the series
$$\displaystyle \sum_{n=0}^\infty (M_\lambda)^n \theta$$
converges in $QC_0(\mathfrak{g}(F))$.
\end{num}

\vspace{3mm}

\noindent Let $\theta\in QC_0(\mathfrak{g}(F))$. Assume first that $F$ is $p$-adic. Then, we need to show that for any open invariant and compact modulo conjugation subset $\omega\subseteq \mathfrak{g}(F)$, the series

\begin{align}\label{eq 4.6.2}
\displaystyle \sum_{n=0}^\infty \mathbf{1}_\omega(M_\lambda)^n \theta
\end{align}

\noindent converges in $QC_c(\omega)$. In some invariant neighborhood of $0$, we have

$$\displaystyle \theta(X)=\sum_{\mathcal{O}\in \Nil(\mathfrak{g})}c_{\theta,\mathcal{O}}(0)\widehat{j}(\mathcal{O},X)$$

\noindent Hence, for $n$ sufficiently large, by \ref{eq 1.7.5}, we have

\[\begin{aligned}
\displaystyle (M_\lambda)^n\theta(X) & =\lvert \lambda\rvert^{-n\delta(G)/2}\sum_{\mathcal{O}\in \Nil(\mathfrak{g})}c_{\theta,\mathcal{O}}(0)\widehat{j}(\mathcal{O},\lambda^{-n}X) \\
 & =\sum_{\mathcal{O}\in \Nil(\mathfrak{g})}\lvert \lambda\rvert^{n(\dim(\mathcal{O})-\delta(G)/2)}c_{\theta,\mathcal{O}}(0)\widehat{j}(\lambda^n\mathcal{O},X)
\end{aligned}\]

\noindent for all $X\in \omega_{\reg}$. This shows that the family of quasi-characters $\{\mathbf{1}_\omega(M_\lambda)^n\theta;\; n\geqslant 1\}$ generates a finite dimensional space. Moreover, by the hypothesis made on $\theta$ and since $\dim(\mathcal{O})<\delta(G)/2$ for any nilpotent orbit $\mathcal{O}$ that is not regular, we see that the series \ref{eq 4.6.2} converges in that finite dimensional space.

\vspace{2mm}

\noindent Assume now that $F=\mathbb{R}$. Then, we need to show that for any invariant and compact modulo conjugation subset $L\subseteq\mathfrak{g}(F)$ and all $u\in I(\mathfrak{g})$, the series

$$\displaystyle \sum_{n=0}^\infty q_{L,u}((M_\lambda)^n \theta)$$

\noindent converges. Actually, we are going to show that

\begin{align}\label{eq 4.6.3}
q_{L,u}(M_{\lambda^n}\theta)\ll \lvert \lambda\rvert^{-n}
\end{align}

\noindent for all $n\geqslant 0$. Obviously, we may assume that $u$ is homogeneous. We distinguish two cases. First assume that $\deg(u)>0$. Enlarging $L$ if necessary, we may assume that $\lambda^{-1}L\subseteq L$. Then, for all $n\geqslant 1$, we have

\[\begin{aligned}
q_{L,u}(M_{\lambda^n}\theta) & =\lvert \lambda\rvert^{-n\deg(u)} \sup_{X\in L_{\reg}} D^G(\lambda^{-n}X)^{1/2} \lvert\left(\partial(u)\theta\right)(\lambda^{-n}X)\rvert \\
 & =\lvert \lambda\rvert^{-n\deg(u)} \sup_{X\in \lambda^{-n}L_{\reg}} D^G(X)^{1/2} \lvert\left(\partial(u)\theta\right)(X)\rvert \\
 & \leqslant \lvert \lambda\rvert^{-n\deg(u)} q_{L,u}(\theta)
\end{aligned}\]

\noindent and this proves \ref{eq 4.6.3} in this case. Now assume that $u=1$. Since $c_{\theta,\mathcal{O}}(0)=0$ for all $\mathcal{O}\in \Nil_{\reg}(\mathfrak{g})$, we have $D^G(X)^{1/2}\theta(X)=O(\lvert X\rvert)$ for $X$ in some invariant neighborhood of $0$. Hence, $L$ being compact modulo conjugation, we have

\[\begin{aligned}
q_{L,1}(M_{\lambda^n}\theta)=\sup_{X\in L_{\reg}}D^G(\lambda^{-n}X)^{1/2} \lvert \theta(\lambda^{-n}X)\rvert\ll \lvert \lambda\rvert^{-n}
\end{aligned}\]

\noindent for all $n\geqslant 0$. This proves \ref{eq 4.6.3} in this case too and this ends the proof of \ref{eq 4.6.1}.

\item There is nothing to prove in the $p$-adic case (since $SQC(\mathfrak{g}(F))=QC_c(\mathfrak{g}(F))$) so we assume $F=\mathbb{R}$. Fix $\lambda\in F^\times$ such that $\lvert \lambda\rvert<1$. Consider the series

$$\displaystyle S(\theta):=\sum_{n=0}^\infty (M_\lambda)^n\theta$$

\noindent for all $\theta\in SQC(\mathfrak{g}(F))$. Then, it is not hard to prove as above that this series converges in $QC(\mathfrak{g}(F)-0)$ and that this defines a continuous linear map $S:SQC(\mathfrak{g}(F))\to QC(\mathfrak{g}(F)-0)$. Let $\varphi\in C^\infty(\mathfrak{g}(F))^G$ be compactly supported modulo conjugation and such that $\varphi=1$ near $0$. Then, we claim that

\vspace{3mm}

\begin{num}
\item\label{eq 4.6.4} For all $\theta\in SQC(\mathfrak{g}(F))$, the quasi-character $(1-\varphi)S(\theta)$ is a Schwartz quasi-character.
\end{num} 

\vspace{3mm}

\noindent Let $\theta\in SQC(\mathfrak{g}(F))$. By Proposition \ref{proposition 4.2.1} 2.(iii), $(1-\varphi)S(\theta)$ is a quasi-character on $\mathfrak{g}(F)$. Denote by $L$ the support of $\varphi$. Let $u\in I(\mathfrak{g})$ and $N\geqslant 1$. Since $(1-\varphi)S(\theta)=S(\theta)$ outside $L$, we only need to show that

\begin{align}\label{eq 4.6.5}
D^G(X)^{1/2}\lvert (\partial(u)S(\theta))(X)\rvert\ll \lVert X\rVert_{\Gamma(\mathfrak{g})}^{-N}
\end{align}

\noindent for all $X\in \mathfrak{g}_{\reg}(F)-L$. Of course, we may assume that $u$ is homogeneous and $N>\deg(u)$. We have

\[\begin{aligned}
\displaystyle D^G(X)^{1/2}\lvert (\partial(u)S(\theta))(X)\rvert & \leqslant \sum_{n=0}^\infty D^G(X)^{1/2} \lvert \left(\partial(u)(M_{\lambda^n}\theta)\right)(X)\rvert \\
 & =\sum_{n=0}^\infty \lvert \lambda\rvert^{-n\deg(u)}D^G(\lambda^{-n}X)^{1/2} \lvert \left(\partial(u)\theta\right)(\lambda^{-n}X)\rvert \\
 & \leqslant q_{N,u}(\theta)\sum_{n=0}^\infty \lvert \lambda\rvert^{-n\deg(u)}\lVert \lambda^{-n}X\rVert_{\Gamma(\mathfrak{g})}^{-N}
\end{aligned}\]

\noindent for all $X\in \mathfrak{g}_{\reg}(F)-L$. Since $L$ is an invariant neighborhood of $0$, we have an inequality $\lVert X\rVert_{\Gamma(\mathfrak{g})}\ll \lvert \mu\rvert \lVert \mu^{-1}X\rVert_{\Gamma(\mathfrak{g})}$ for all $X\in \mathfrak{g}(F)-L$ and all $\mu\in F^\times$ such that $\lvert \mu\rvert\leqslant 1$. Hence, the last sum above is essentially bounded by

$$\displaystyle \left(\sum_{n=1}^\infty \lvert \lambda\rvert^{n(N-\deg(u))}\right)\lVert X\rVert_{\Gamma(\mathfrak{g})}^{-N}$$

\noindent for all $X\in \mathfrak{g}_{\reg}(F)-L$. Since we are assuming that $N>\deg(u)$ and $\lvert \lambda\rvert<1$, this last term is finite and this shows \ref{eq 4.6.5}. This ends the proof of \ref{eq 4.6.4}.

\vspace{2mm}

\noindent Consider the linear map

$$L_\varphi:SQC(\mathfrak{g}(F))\to QC_c(\mathfrak{g}(F))$$
$$L_\varphi(\theta)=\theta+(M_\lambda-1)\left[(1-\varphi)S(\theta)\right]$$

\noindent It indeed takes value in $QC_c(\mathfrak{g}(F))$ since $(M_\lambda-1)\left[(1-\varphi)S(\theta)\right]=-\theta$ outside $\Supp(\varphi)$. Moreover, it is continuous by Lemma \ref{lemma 4.2.2}(iv) and the closed graph theorem. Let $\ell$ be a continuous linear form on $QC_c(\mathfrak{g}(F))$ that satisfies $\ell(M_\lambda \theta)=\ell(\theta)$ for all $\theta\in QC_c(\mathfrak{g}(F))$. Then, we can extend $\ell$ to a continuous linear form on $SQC(\mathfrak{g}(F))$ by setting

$$\ell(\theta)=\ell(L_\varphi(\theta))$$

\noindent for all $\theta\in SQC(\mathfrak{g}(F))$. It is indeed an extension since $(1-\varphi)S(\theta)\in QC_c(\mathfrak{g}(F))$ for all $\theta\in QC_c(\mathfrak{g}(F))$. $\blacksquare$
\end{enumerate}

\subsection{Quasi-characters and parabolic induction}\label{section 4.7}

\noindent Let $M$ be a Levi subgroup of $G$. Let us fix for all $x\in G_{\ssi}(F)$ a system of representatives $\mathcal{X}^M(x)$ for the $M(F)$-conjugacy classes of elements in $M(F)$ that are $G(F)$-conjugate to $x$. Recall that in Section \ref{section 3.4}, we have defined a parabolic induction morphism $i_M^G:\mathcal{D}'(M(F))^M\to \mathcal{D}'(G(F))^G$ which sends distributions representable by locally integrable invariant functions on $M(F)$ to distributions representable by locally integrable invariant functions on $G(F)$. In particular, for any quasi-character $\theta^M$ on $M(F)$, $i_M^G(\theta^M)$ is a well-defined locally integrable invariant function on $G(F)$.

\begin{prop}\label{proposition 4.7.1}
Let $\theta^M$ be a quasi-character on $M(F)$. Then,

\begin{enumerate}[(i)]
\item $i_M^G(\theta^M)$ is a quasi-character.

\item Let $\theta=i_M^G(\theta^M)$. Then, we have

$$\displaystyle D^G(x)^{1/2}c_{\theta}(x)=[Z_G(x)(F):G_x(F)]\sum_{y\in \mathcal{X}^M(x)} [Z_M(y)(F):M_y(F)]^{-1}D^M(y)^{1/2}c_{\theta^M}(y)$$

\noindent for all $x\in G_{\ssi}(F)$.
\end{enumerate}
\end{prop}

\vspace{2mm}

\noindent\ul{Proof}:

\begin{enumerate}[(i)]
\item If $F$ is $p$-adic, this is proved in \cite{Wa1}. Assume that $F=\mathbb{R}$. We already know (cf.\ Section \ref{section 3.4}) that $i_M^G(\theta^M)$ is representable by a smooth function $\theta$ on $G_{\sreg}(F)$. It is not hard to see using \ref{eq 3.4.2}, that $\theta$ extends to a smooth function on $G_{\reg}(F)$ and that the function $(D^G)^{1/2}\theta$ is locally bounded. Moreover, by \ref{eq 3.4.1} and since $\theta^M$ is a quasi-character, for all $z\in \mathcal{Z}(\mathfrak{g})$ the distribution $zT_\theta$ is represented by the function $i_M^G(z_M\theta^M)$. It easily follows that $\theta$ is a quasi-character.

\item Assume first that $G_x$ is not quasi-split. Then, by Proposition \ref{proposition 4.5.1}.1.(i), both sides are easily seen to be zero (notice that for $y\in M_{\ssi}(F)$, $M_y$ is a Levi subgroup of $G_y$). Assume now that $G_x$ is quasi-split. Let us fix $B_x\subset G_x$ a Borel subgroup and $T_{\quasid,x}\subset B_x$ a maximal torus (both defined over $F$). Then, by Proposition \ref{proposition 4.5.1}.1.(ii), we have the formula

\begin{align}\label{eq 4.7.1}
\displaystyle D^G(x)c_{\theta}(x)=\lvert W(G_x,T_{\quasid,x})\rvert^{-1}\lim\limits_{x'\in T_{\quasid,x}(F)\to x}D^G(x')^{1/2}\theta(x')
\end{align}

\noindent Moreover, by \ref{eq 3.4.2}, we know that for all $x'\in T_{\quasid,x}(F)\cap G_{\sreg}(F)$ we have

\begin{align}\label{eq 4.7.2}
D^G(x')^{1/2}\theta(x')=\sum_{y'\in \mathcal{X}^M(x')} D^M(y')^{1/2}\theta^M(y')
\end{align}

\noindent For all $y\in \mathcal{X}^M(x)$, the group $G_y$ is quasi-split and $M_y$ is one of its Levi subgroup. Hence, we may fix for all $y\in \mathcal{X}^M(x)$ a Borel subgroup $B_y\subset G_y$ and a maximal torus $T_{\quasid,y}\subset B_y$ (both again defined over $F$) such that $T_{\quasid,y}\subset M_y$. Let us also fix for all $y\in \mathcal{X}^M(x)$ an element $g_y\in G(F)$ such that $g_y^{-1}xg_y=y$ and $g_y^{-1}T_{\quasid,x}g_y=T_{\quasid,y}$. We claim the following

\vspace{3mm}

\begin{num}
\item\label{eq 4.7.3} Let $x'\in T_{\quasid,x}(F)\cap G_{\sreg}(F)$. Then, for all $y'\in \mathcal{X}^M(x')$, there exist $y\in \mathcal{X}^M(x)$ and $g\in \No_{Z_G(y)(F)}(T_{\quasid,y})$ such that $y'$ and $g^{-1}g_y^{-1}x'g_yg$ are $M(F)$-conjugate. Moreover, for all $y_1,y_2\in \mathcal{X}^M(x)$ and all $g_i\in \No_{Z_G(y_i)(F)}(T_{\quasid,y_i})$, $i=1,2$, the elements $g_1^{-1}g_{y_1}^{-1}x'g_{y_1}g_1$ and $g_2^{-1}g_{y_2}^{-1}x'g_{y_2}g_2$ are $M(F)$-conjugate if and only if $y_1=y_2$ and $g_2\in g_1\No_{Z_M(y_1)(F)}(T_{\quasid,y_1})$.
\end{num}

\vspace{3mm}

\noindent Let $x'\in T_{\quasid,x}(F)\cap G_{\sreg}(F)$ and $y'\in \mathcal{X}^M(x')$. Choose $\gamma\in G(F)$ such that $y'=\gamma^{-1}x'\gamma$. Then the centralizer of $y'$ in $G$ is the maximal torus $\gamma^{-1}T_{\quasid,x}\gamma$ which is contained in $M$. It follows that $\gamma^{-1}x\gamma\in M(F)$. By definition of $\mathcal{X}^M(x)$, up to translating $\gamma$ by an element of $M(F)$ (and conjugating $y'$ by the same element), we may assume that there exists $y\in \mathcal{X}^M(x)$ such that $y=\gamma^{-1}x\gamma$. Then $\gamma^{-1}B_x\gamma$ is a Borel subgroup of $G_y$ and $\gamma^{-1}T_{\quasid,x}\gamma$ is a maximal torus of $\gamma^{-1}B_x\gamma$ which is the centralizer of $\gamma^{-1}x'\gamma$ and so is contained in $M_y$. Hence, up to translating $\gamma$ by an element of $M_y(F)$ we may further assume that $\gamma^{-1}T_{\quasid,x}\gamma=T_{\quasid,y}$. Consider the element $g=g_y^{-1}\gamma$. It centralizes $y$ and normalizes $T_{\quasid,y}$. It follows that $g\in \No_{Z_G(y)(F)}(T_{\quasid,y})$ and this proves the first part of the claim. Let $y_1$, $y_2$, $g_1$ and $g_2$ be as in the second part of the claim and assume that $g_1^{-1}g_{y_1}^{-1}x'g_{y_1}g_1$ and $g_2^{-1}g_{y_2}^{-1}x'g_{y_2}g_2$ are $M(F)$-conjugate. These two elements are conjugate by $m=g_1^{-1}g_{y_1}^{-1}g_{y_2}g_2$. Since the centralizer of $g_1^{-1}g_{y_1}^{-1}x'g_{y_1}g_1$ is the torus $T_{\quasid,y_1}$ which is contained in $M$, we have $m\in M(F)$. But, we easily check that $my_2m^{-1}=y_1$. By definition of $\mathcal{X}^M(x)$, it follows that $y_1=y_2$ and hence $m=g_1^{-1}g_2\in \No_{Z_G(y_1)(F)}(T_{\quasid,y_1})\cap M(F)=\No_{Z_M(y_1)(F)}(T_{\quasid,y_1})$ and this ends the proof of \ref{eq 4.7.3}.

\vspace{2mm}

Let us fix for all $y\in \mathcal{X}^M(x)$, a set $N^G_M(y)$ of representatives for the left cosets of $\No_{Z_M(y)(F)}(T_{\quasid,y})$ in $\No_{Z_G(y)(F)}(T_{\quasid,y})$. By \ref{eq 4.7.3}, we may assume that for all $x'\in T_{\quasid,x}(F)\cap G_{\sreg}(F)$ we have

$$\mathcal{X}^M(x')=\{g^{-1}g_y^{-1}x'g_yg;\;\; y\in \mathcal{X}^M(x)\; g\in N^G_M(y)\}$$

\noindent Hence, by \ref{eq 4.7.2}, we get

\begin{align}\label{eq 4.7.4}
\displaystyle D^G(x')^{1/2}\theta(x')=\sum_{y\in \mathcal{X}^M(x)}\sum_{N^G_M(y)} D^M(g^{-1}g_y^{-1}x'g_yg)^{1/2}\theta^M(g^{-1}g_y^{-1}x'g_yg)
\end{align}

\noindent for all $x'\in T_{\quasid,x}(F)\cap G_{\sreg}(F)$. Notice that for all $y\in \mathcal{X}^M(x)$ and all $g\in N^G_M(y)$, we have $g^{-1}g_y^{-1}T_{\quasid,x}g_yg=T_{\quasid,y}$ and $g^{-1}g_y^{-1}x'g_yg\to y$ as $x'\to x$. Hence, taking the limit in \ref{eq 4.7.4} as $x'\to x$, we get, by \ref{eq 4.7.1} and Proposition \ref{proposition 4.5.1}.1.(ii),

$$\displaystyle D^G(x)^{1/2}c_{\theta}(x)=\lvert W(G_x,T_{\quasid,x})\rvert^{-1}\sum_{y\in \mathcal{X}^M(x)} \lvert N_M^G(y)\rvert \lvert W(M_y,T_{\quasid,y})\rvert D^M(y)^{1/2}c_{\theta^M}(y)$$

\noindent To conclude, it suffices now to show that for all $y\in \mathcal{X}^M(x)$, we have

\begin{align}\label{eq 4.7.5}
\lvert W(G_x,T_{\quasid,x})\rvert^{-1}\lvert N_M^G(y)\rvert \lvert W(M_y,T_{\quasid,y})\rvert=[Z_G(x)(F):G_x(F)][Z_M(y)(F):M_y(F)]^{-1}
\end{align}

\noindent Let $y\in\mathcal{X}^M(x)$. By definition, we have

\[\begin{aligned}
\lvert N^G_M(y)\rvert & =\left\lvert \No_{Z_G(y)(F)}(T_{\quasid,y})/\No_{Z_M(y)(F)}(T_{\quasid,y})\right\rvert \\
 & =\left\lvert \No_{Z_G(y)(F)}(T_{\quasid,y})/T_{\quasid,y}(F)\right\rvert \times\left\lvert \No_{Z_M(y)(F)}(T_{\quasid,y})/T_{\quasid,y}(F)\right\rvert^{-1}
\end{aligned}\]

\noindent Since the pairs $(x,T_{\quasid,x})$ and $(y,T_{\quasid,y})$ are $G(F)$-conjugate, we have

$$\left\lvert \No_{Z_G(y)(F)}(T_{\quasid,y})/T_{\quasid,y}(F)\right\rvert=\left\lvert \No_{Z_G(x)(F)}(T_{\quasid,x})/T_{\quasid,x}(F)\right\rvert$$

\noindent Hence, the left hand side of \ref{eq 4.7.5} is the product of the two following terms

$$\left\lvert W(G_x,T_{\quasid,x})\right\rvert^{-1} \left\lvert \No_{Z_G(x)(F)}(T_{\quasid,x})/T_{\quasid,x}(F)\right\rvert$$

\noindent and

$$\left\lvert W(M_y,T_{\quasid,y})\right\rvert \left\lvert \No_{Z_M(y)(F)}(T_{\quasid,y})/T_{\quasid,y}(F)\right\rvert^{-1}$$

\noindent Let us look at the first term above. Every coset of $G_x(F)$ in $Z_G(x)(F)$ contains an element that normalizes $T_{\quasid,x}$. Indeed, for all $z\in Z_G(x)(F)$ the torus $zT_{\quasid,x}z^{-1}$ is also a maximal torus of $G_x$ which is contained in a Borel subgroup, hence there exists $g\in G_x(F)$ such that $gzT_{\quasid,x}z^{-1}g^{-1}=T_{\quasid,x}$. It now follows easily from this that we have

\begin{align}\label{eq 4.7.6}
\left\lvert W(G_x,T_{\quasid,x})\right\rvert^{-1} \left\lvert \No_{Z_G(x)(F)}(T_{\quasid,x})/T_{\quasid,x}(F)\right\rvert=[Z_G(x)(F):G_x(F)]
\end{align}

\noindent We similarly prove that 

\begin{align}\label{eq 4.7.7}
\left\lvert W(M_y,T_{\quasid,y})\right\rvert \left\lvert \No_{Z_M(y)(F)}(T_{\quasid,y})/T_{\quasid,y}(F)\right\rvert^{-1}=[Z_M(y)(F):M_y(F)]^{-1}
\end{align}

\noindent Now \ref{eq 4.7.5} follows from \ref{eq 4.7.6} and \ref{eq 4.7.7} and this ends the proof of the proposition. $\blacksquare$
\end{enumerate}

\subsection{Quasi-characters associated to tempered representations and Whittaker datas}\label{section 4.8}

\noindent Recall that a {\em Whittaker datum} for $G$ is a $G(F)$-conjugacy class of pairs $(U_B,\xi_B)$ where $U_B$ is the unipotent radical of a Borel subgroup $B$ of $G$, defined over $F$, and $\xi_B:U_B(F)\to \mathbb{S}^1$ is a generic character on $U_B(F)$ (generic means that the stabilizer of $\xi_B$ in $B(F)$ coincides with $Z_G(F)U_B(F)$). Of course, Whittaker data for $G$ exist if and only if $G$ is quasi-split. Using the bilinear form $B(.,.)$ and the additive character $\psi$, we can define a bijection $\mathcal{O}\mapsto (U_{\mathcal{O}},\xi_{\mathcal{O}})$ between $\Nil_{\reg}(\mathfrak{g})$ and the set of Whittaker data for $G$ as follows. Let $\mathcal{O}\in \Nil_{\reg}(\mathfrak{g})$. Pick $Y\in \mathcal{O}$ and extend it to an $\mathfrak{sl}_2$-triple $(Y,H,X)$. Then $X$ is a regular nilpotent element and hence belongs to exactly one Borel subalgebra $\mathfrak{b}_{\mathcal{O}}$ of $\mathfrak{g}$ that is defined over $F$. Let $B_{\mathcal{O}}$ be the corresponding Borel subgroup and $U_{\mathcal{O}}$ be its unipotent radical. The assignment $u\in U_{\mathcal{O}}(F)\mapsto \xi_{\mathcal{O}}(u)=\psi\left(B(Y,\log(u))\right)$ defines a generic character on $U_{\mathcal{O}}(F)$. Moreover, the $G(F)$-conjugacy class of $(U_{\mathcal{O}},\xi_{\mathcal{O}})$ only depends on $\mathcal{O}$ and this defines the desired bijection.

\vspace{2mm}

\noindent Let $\pi$ be a tempered irreducible representation of $G(F)$. For $\mathcal{O}\in \Nil_{\reg}(F)$, we will say that $\pi$ has a {\em Whittaker model of type $\mathcal{O}$} if there exists a nonzero continuous linear form $\ell:\pi^\infty\to \mathbb{C}$ such that $\ell\circ \pi(u)=\xi_{\mathcal{O}}(u)\ell$ for all $u\in U_{\mathcal{O}}(F)$. Recall that the character $\theta_\pi$ of $\pi$ is a quasi-character on $G(F)$. For all $\mathcal{O}\in \Nil_{\reg}(\mathfrak{g})$, we set $\gls{cpiO1}=c_{\theta_\pi,\mathcal{O}}(1)$. In Section \ref{section 2.7}, we have defined a space $\underline{\mathcal{X}}(G)$ of virtual tempered representations of $G(F)$. The character of a virtual representation is defined by linearity.

\vspace{2mm}

\begin{prop}\label{proposition 4.8.1}
\begin{enumerate}[(i)]
\item Let $\pi\in \Temp(G)$. Then, for all $\mathcal{O}\in \Nil_{\reg}(\mathfrak{g})$, we have

$$\displaystyle c_{\pi,\mathcal{O}}(1)=\left\{
    \begin{array}{ll}
        1 & \mbox{ if } \pi \mbox{ has a Whittaker model of type } \mathcal{O} \\
        0 & \mbox{ otherwise}
    \end{array}
\right.
$$

\item The map

$$\pi\in \underline{\mathcal{X}}(G)\mapsto \theta_\pi\in QC(G(F))$$

\noindent is smooth. Moreover, if $F=\mathbb{R}$, for every continuous semi-norm $\nu$ on $QC(G(F))$, there exists an integer $k\geqslant 0$ such that

$$\nu(\theta_{\pi})\ll N(\pi)^k$$

\noindent for all $\pi\in \mathcal{X}(G)$.
\end{enumerate}
\end{prop}

\vspace{2mm}

\noindent\ul{Proof}: 
\begin{enumerate}[(i)]
\item In the case where $F$ is $p$-adic and $G$ is split, this is due to Rodier \cite{Ro}. The same proof works equally well for general quasi-split groups and is contained in the more general results of \cite{MW}. Finally, when $F=\mathbb{R}$ it is a theorem of Matumoto (\cite{Mat} Theorem C).

\item The first part is easy to prove. Indeed, this amounts to showing that for every Levi subgroup $M$ of $G$ and for all $\sigma\in \underline{\mathcal{X}}_{\elli}(M)$ the map

$$\lambda\in i\mathcal{A}_M^*\mapsto \theta_{\pi_\lambda}\in QC(G(F))$$

\noindent is smooth, where we have set $\pi_\lambda=i_M^G(\sigma_\lambda)$ for all $\lambda\in i\mathcal{A}_M^*$. By \ref{eq 3.4.3}, we have $\theta_{\pi_\lambda}=i_M^G(\theta_{\sigma_\lambda})$ for all $\lambda\in i\mathcal{A}_M^*$. The linear map

$$i_M^G:QC(M(F))\to QC(G(F))$$

\noindent is easily seen to be continuous using \ref{eq 3.4.1} and \ref{eq 3.4.2}, whereas the map

$$\lambda\in i\mathcal{A}_M^*\mapsto \theta_{\sigma_\lambda}\in QC(M(F))$$

\noindent is obviously smooth. This handles the first part of the proposition.

\vspace{2mm}

\noindent Let us now prove the second part of the proposition. So assume that $F=\mathbb{R}$. Let $\mathcal{M}$ be a set of representatives for the conjugacy classes of Levi subgoups of $G$. Recall that by definition of $\mathcal{X}(G)$, we have

$$\displaystyle \mathcal{X}(G)=\bigcup_{M\in \mathcal{M}} i_M^G(\mathcal{X}_{\elli}(M))$$

\noindent Let $M\in \mathcal{M}$ be such that $M\neq G$. Then, by induction, we may assume that the result is true for $M$. As we just saw, the linear map $i_M^G:QC(M(F))\to QC(G(F))$ is continuous. It immediately follows that for every continuous semi-norm $\nu$ on $QC(G(F))$, there exists $k\geqslant 0$ such that

$$\nu(\theta_\pi)\ll N(\pi)^k$$

\noindent for all $\pi\in i_M^G(\mathcal{X}_{\elli}(M))$. Combining these inequalities for all $M\in \mathcal{M}$, $M\neq G$, we are left with proving the inequality of the proposition only for $\pi\in \mathcal{X}_{\elli}(G)$, that is

\vspace{3mm}

\begin{num}
\item\label{eq 4.8.1} For each continuous semi-norm $\nu$ on $QC(G(F))$, there exists $k\geqslant 0$ such that
$$\displaystyle \nu(\theta_\pi)\ll N(\pi)^k$$
for all $\pi\in \mathcal{X}_{\elli}(G)$.
\end{num} 

\vspace{3mm}

\noindent For all $z\in \mathcal{Z}(\mathfrak{g})$ we have $z\theta_\pi=\chi_\pi(z)\theta_\pi$ and there exists $k_0\geqslant 0$ such that $\left\lvert \chi_\pi(z)\right\rvert\ll N(\pi)^{k_0}$ for all $\pi\in \mathcal{X}_{\elli}(G)$. Hence, it clearly suffices to prove the existence of $C>0$ such that

\begin{align}\label{eq 4.8.2}
\displaystyle \sup_{x\in G_{\reg}(F)}D^G(x)^{1/2} \left\lvert \theta_\pi(x)\right\rvert\leqslant C
\end{align}

\noindent for all $\pi\in \mathcal{X}_{\elli}(G)$. Harish-Chandra has completely described the characters $\theta_\pi$ of elliptic representations $\pi\in \mathcal{X}_{\elli}(G)$. More precisely, let $T$ be a maximal torus of $G$ which is elliptic (if such a torus doesn't exist then $G$ has no elliptic representations). Denote by $T(\mathbb{R})^*$ the group of continuous unitary characters of $T(\mathbb{R})$. Then, to each element $b^*\in T(\mathbb{R})^*$, Harish-Chandra associates a certain function $\theta_{b^*}$ on $G_{\reg}(\mathbb{R})$ (cf.\ Theorem 24 p.261 of \cite{Va}). These are invariant eigendistributions on $G(\mathbb{R})$ and some of them might be equal to zero. Moreover, by Theorem 24 (c) p.261 of \cite{Va}, there exists $C>0$ such that

\begin{align}\label{eq 4.8.3}
\displaystyle \sup_{x\in G_{\reg}(F)}D^G(x)^{1/2} \left\lvert \theta_{b^*}(x)\right\rvert\leqslant C
\end{align}

\noindent for all $b^*\in T(\mathbb{R})^*$. Now let $\pi\in \mathcal{X}_{\elli}(G)$. Recall that $\pi$ is a linear combination of constituents of a certain induced representation $i_M^G(\sigma)$, $M$ a Levi subgroup of $G$ and $\sigma\in \Pi_2(M)$. Let us denote, as in Section \ref{section 2.7}, by $W(\sigma)$ the stabilizer of $\sigma$ in $W(G,M)$. Then, by the first equality after Theorem 13 of \cite{HC2}, there exists $b^*\in T(\mathbb{R})^*$ such that the equality

$$\theta_\pi=\lvert W(\sigma)\rvert \theta_{b^*}$$

\noindent holds up to a scalar of module one (recall that $\pi$ itself is defined up to such a scalar). Of course the term $\lvert W(\sigma)\rvert$ is bounded independently of $\pi$. Hence, \ref{eq 4.8.2} follows from \ref{eq 4.8.3} and this ends the proof of the proposition. $\blacksquare$

\end{enumerate}

\section{Strongly cuspidal functions}\label{section 5}

This chapter is devoted to the study of the so-called {\em strongly cuspidal} functions; a notion that we borrow from the work of Waldspurger \cite{Wa1}. These are functions $f$ on the group $G(F)$ or its Lie algebra satisfying a certain geometric condition; namely that for every proper parabolic subgroup $P=MU$ the function on $M(F)$ defined by integration over $U(F)$ vanishes identically. Their importance stems from the fact that the simple local trace formulas to be developed in Chapter \ref{section 7} to \ref{section 11} are functionals on the space of strongly cuspidal functions.

In Section \ref{section 5.1}, we define strongly cuspidal functions and derive their basic properties. Following Waldspurger, we study the weighted orbital integrals of such functions in Section \ref{section 5.2}. In Sections \ref{section 5.3} and \ref{section 5.4}, we give a spectral characterization of the strongly cuspidal functions and study their weighted characters. Section \ref{section 5.5} recalls the local trace formula of Arthur in the particular case where one of the test functions is strongly cuspidal (it then takes a particularly nice form). Section \ref{section 5.6} contains a very important construction that allows to associate to any strongly cuspidal function a quasi-character in the sense of Chapter \ref{section 4} (this construction is also due to Waldspurger \cite{Wa1}). Finally in Section \ref{section 5.7}, we prove a technical proposition which allows to replace a strongly cuspidal function $f$ by another one with the same associated quasi-character but whose semi-simple descents to elliptic elements is again strongly cuspidal. This proposition will play a crucial role in the proof of the geometric expansions of our local trace formulas (Theorem \ref{theorem 11.2.1} and Theorem \ref{theorem 11.2.3}).
 
\subsection{Definition, first properties}\label{section 5.1}

\noindent For every parabolic subgroup $P=MU$ of $G$, we define continuous linear maps

$$\displaystyle f\in\mathcal{C}(G(F))\mapsto \gls{fP}\in\mathcal{C}(M(F))$$

$$\displaystyle \varphi\in \mathcal{S}(\mathfrak{g}(F))\mapsto \gls{phiP}\in \mathcal{S}(\mathfrak{m}(F))$$

\noindent by setting

$$\displaystyle f^P(m)=\delta_P(m)^{1/2}\int_{U(F)}f(mu)du \mbox{ and } \varphi^P(X)=\int_{\mathfrak{u}(F)}\varphi(X+N)dN$$

\noindent We will say that a function $f\in \mathcal{C}(G(F))$ (resp.\ $\varphi\in \mathcal{S}(\mathfrak{g}(F))$) is {\em strongly cuspidal} if $f^P=0$ (resp.\ $\varphi^P=0$) for every proper parabolic subgroup $P$ of $G$. We will denote by $\gls{CscuspG}$, $\gls{SscuspG}$ and $\gls{Sscuspg}$ the subspaces of strongly cuspidal functions in $\mathcal{C}(G(F))$, $\mathcal{S}(G(F))$ and $\mathcal{S}(\mathfrak{g}(F))$ respectively. More generally, if $\Omega\subseteq G(F)$ (resp.\ $\omega\subseteq \mathfrak{g}(F)$) is a completely $G(F)$-invariant open subset, we will set $\gls{SscuspOmega}=\mathcal{S}(\Omega)\cap \mathcal{S}_{\scusp}(G(F))$ (resp.\ $\gls{Sscuspomega}=\mathcal{S}(\omega)\cap \mathcal{S}_{\scusp}(\mathfrak{g}(F))$) (the subspaces $\mathcal{S}(\Omega)$ and $\mathcal{S}(\omega)$ have been defined in Section \ref{section 3.1}). In the real case, we have $(zf)^P=z_Mf^P$ and $(\partial(u)\varphi)^P=\partial(u_M)\varphi^P$ for all $f\in \mathcal{C}(G(F))$, all $z\in \mathcal{Z}(\mathfrak{g})$, all $\varphi\in \mathcal{S}(\mathfrak{g}(F))$ and all $u\in I(\mathfrak{g})$. Hence, the action of $\mathcal{Z}(\mathfrak{g})$ preserves the spaces $\mathcal{C}_{\scusp}(G(F))$, $\mathcal{S}_{\scusp}(G(F))$ and $\mathcal{S}_{\scusp}(\Omega)$ and the action of $I(\mathfrak{g})$ preserves the spaces $\mathcal{S}_{\scusp}(\mathfrak{g}(F))$ and $\mathcal{S}_{\scusp}(\omega)$.

\vspace{2mm}

\noindent Let $f\in \mathcal{C}(G(F))$. By the usual variable change we see that

$$\displaystyle f^P(m)=D^G(m)^{1/2}D^M(m)^{-1/2}\int_{U(F)} f(u^{-1}mu)du$$

\noindent for every parabolic subgroup $P=MU$ and all $m\in M(F)\cap G_{\reg}(F)$. Hence, an equivalent condition for $f$ to be strongly cuspidal is that the integral

$$\displaystyle \int_{U(F)} f(u^{-1}mu)du$$

\noindent is zero for every proper parabolic subgroup $P=MU$ and all $m\in M(F)\cap G_{\reg}(F)$. Similarly, a function $\varphi\in \mathcal{S}(\mathfrak{g}(F))$ is strongly cuspidal if and only if the integral

$$\displaystyle \int_{U(F)} \varphi(u^{-1}Xu)du$$

\noindent is zero for every proper parabolic subgroup $P=MU$ and all $X\in \mathfrak{m}(F)\cap \mathfrak{g}_{\reg}(F)$. It follows from these descriptions and Proposition \ref{proposition 3.1.1}(iv) that for every function $\varphi\in C^\infty(G(F))^G$ (resp.\ $\varphi\in C^\infty(\mathfrak{g}(F))^G$) which is compactly supported modulo conjugation, multiplication by $\varphi$ preserves $\mathcal{S}_{\scusp}(G(F))$ (resp.\ $\mathcal{S}_{\scusp}(\mathfrak{g}(F))$). Moreover, choosing $\varphi\in C^\infty(\mathfrak{g}(F))^G$ compactly supported modulo conjugation and such that $\varphi=1$ in a neighborhood of $0$, it is not hard to see that $\lim\limits_{\lambda\to\infty}\varphi_\lambda f=f$ in $\mathcal{S}(\mathfrak{g}(F))$ for all $f\in \mathcal{S}_{\scusp}(\mathfrak{g}(F))$ (recall that $\varphi_\lambda(X)=\varphi(\lambda^{-1}X)$). Hence, we have

\vspace{3mm}

\begin{num}
\item\label{eq 5.1.1} The subspace of functions $f\in \mathcal{S}_{\scusp}(\mathfrak{g}(F))$ that are compactly supported modulo conjugation is dense in $\mathcal{S}_{\scusp}(\mathfrak{g}(F))$.
\end{num}

\vspace{3mm}

\noindent Let $\omega\subseteq \mathfrak{g}(F)$ be a $G$-excellent open subset and set $\Omega=\exp(\omega)$. Then the map $f\mapsto f_\omega$ induces an isomorphism

$$\mathcal{S}_{\scusp}(\Omega)\simeq \mathcal{S}_{\scusp}(\omega)$$

\noindent (this follows from Lemma \ref{lemma 3.3.1}). Finally, we leave to the reader the simple task of checking that the Fourier transform preserves $\mathcal{S}_{\scusp}(\mathfrak{g}(F))$.

\subsection{Weighted orbital integrals of strongly cuspidal functions}\label{section 5.2}

\noindent Let $M$ be a Levi subgroup of $G$. Recall that in Section \ref{section 1.10}, we have defined a family of tempered distributions $J_L^Q(x,.)$ on $G(F)$ for all $x\in M(F)\cap G_{\reg}(F)$, all $L\in\mathcal{L}(M)$ and all $Q\in\mathcal{F}(L)$. We have also defined tempered distributions $J_L^Q(X,.)$ on $\mathfrak{g}(F)$ for all $X\in \mathfrak{m}(F)\cap \mathfrak{g}_{\reg}(F)$, all $L\in\mathcal{L}(M)$ and all $Q\in\mathcal{F}(L)$. These distributions depended on the choice of a maximal compact subgroup $K$ which is special in the $p$-adic case.

\begin{lem}\label{lemma 5.2.1}
Let $f\in \mathcal{C}_{\scusp}(G(F))$ (resp.\ $f\in \mathcal{S}_{\scusp}(\mathfrak{g}(F))$) be a strongly cuspidal function and fix $x\in M(F)\cap G_{\reg}(F)$ (resp.\ $X\in \mathfrak{m}(F)\cap \mathfrak{g}_{\reg}(F)$). Then

\begin{enumerate}[(i)]

\item For all $L\in \mathcal{L}(M)$ and all $Q\in \mathcal{F}(L)$, if $L\neq M$ or $Q\neq G$, we have

$$J^Q_L(x,f)=0 \;\;\; \mbox{(resp.\ } J^Q_L(X,f)=0 \mbox{)}$$

\item The weighted orbital integral $J_M^G(x,f)$ (resp.\ $J_M^G(X,f)$) doesn't depend on the choice of $K$;

\item If $x\notin M(F)_{\elli}$ (resp.\ $X\notin \mathfrak{m}(F)_{\elli}$), we have

$$J_M^G(x,f)=0\;\;\; \mbox{(resp.\ } J_M^G(X,f)=0\mbox{)}$$

\item For all $y\in G(F)$, we have

$$J^G_{yMy^{-1}}(yx y^{-1},f)=J_M^G(x,f)\;\;\; \mbox{(resp.\ } J^G_{yMy^{-1}}(yXy^{-1},f)=J_M^G(X,f)\mbox{)}$$ 
\end{enumerate}
\end{lem}

\vspace{2mm}

\noindent\ul{Proof}: This is proved in \cite{Wa1} when $F$ is $p$-adic. The proof works equally well for $F=\mathbb{R}$. $\blacksquare$

\vspace{2mm}

\noindent For all $x\in G_{\reg}(F)$, let us denote by $\gls{M(x)}$ the centralizer of $A_{G_x}$ in $G$. It is the minimal Levi subgroup of $G$ containing $x$. Let $f\in\mathcal{C}_{\scusp}(G(F))$. Then, we set

$$\gls{thetaf}(x)=(-1)^{a_G-a_{M(x)}}\nu(G_x)^{-1}D^G(x)^{-1/2} J_{M(x)}^G(x,f)$$

\noindent for all $x\in G_{\reg}(F)$ (where we recall that $\nu(G_x)$ is a normalizing constant making the Haar measure on the torus $G_x(F)$ of total mass $1$, see \S \ref{section 1.6}). By the point (iv) of the above lemma, the function $\theta_f$ is invariant. We define similarly an invariant function $\theta_f$ on $\mathfrak{g}_{\reg}(F)$ for all $f\in \mathcal{S}_{\scusp}(\mathfrak{g}(F))$, by setting

$$\gls{thetaf}(X)=(-1)^{a_G-a_{M(X)}}\nu(G_X)^{-1}D^G(X)^{-1/2} J_{M(X)}^G(X,f)$$

\noindent for all $X\in \mathfrak{g}_{\reg}(F)$, where $\gls{M(X)}$ denotes the centralizer of $A_{G_X}$ in $G$.

\vspace{2mm}

\begin{lem}\label{lemma 5.2.2}
Assume that $F=\mathbb{R}$. Then,
\begin{enumerate}[(i)]
\item For all $f\in \mathcal{S}_{\scusp}(\mathfrak{g}(F))$, $\theta_f$ is a smooth function on  $\mathfrak{g}_{\reg}(F)$ and we have $\partial(u)\theta_f=\theta_{\partial(u)f}$ for all $u\in I(\mathfrak{g})$. Moreover, there exists $k\geqslant 0$ such that for all $N\geqslant 1$ there exists a continuous semi-norm $\nu_N$ on $\mathcal{S}_{\scusp}(\mathfrak{g}(F))$ such that

$$D^G(X)^{1/2}\lvert \theta_f(X)\rvert\leqslant \nu_N(f) \log\left(2+D^G(X)^{-1}\right)^k\lVert X\rVert_{\Gamma(\mathfrak{g})}^{-N}$$

\noindent for all $X\in \mathfrak{g}_{\reg}(F)$ and all $f\in \mathcal{S}_{\scusp}(\mathfrak{g}(F))$.

\item For all $f\in \mathcal{C}_{\scusp}(G(F))$, $\theta_f$ is a smooth function on $G_{\reg}(F)$ and we have $z\theta_f=\theta_{zf}$ for all $z\in \mathcal{Z}(\mathfrak{g})$.
\end{enumerate}
\end{lem}

\vspace{2mm}

\noindent\ul{Proof}: By semi-simple descent (Lemma \ref{lemma 3.2.1}), the first point of (i) and (ii) follow directly from Lemma \ref{lemma 1.10.1} and the point (i) of the last lemma. The estimates in (i) is a direct consequence of \ref{eq 1.10.1}. $\blacksquare$

\subsection{Spectral characterization of strongly cuspidal functions}\label{section 5.3}

\noindent Let $P=MU$ be a parabolic subgroup of $G$ and $\sigma$ a tempered representation of $M(F)$. We have a natural isomorphism $\End(i_P^G(\sigma))^\infty\simeq i_{P\times P}^{G\times G}(\End(\sigma))^\infty$ sending a function $\varphi\in i_{P\times P}^{G\times G}(\End(\sigma))^\infty$ to the operator

$$\displaystyle e\in i_P^G(\sigma)\mapsto \left(g\mapsto \int_{P(F)\backslash G(F)} \varphi(g,g')e(g')dg'\right)$$

\noindent Let $f\in \mathcal{C}(G(F))$. A direct computation shows that the operator $i_P^G(\sigma,f)\in \End(i_P^G(\sigma))^\infty$ corresponds to the function $i_P^G(\sigma,f)(.,.)\in i_{P\times P}^{G\times G}(\End(\sigma))^\infty$ given by

$$i_P^G(\sigma,f)(g,g')=\sigma\left[ \left(L(g)R(g')f\right)^P\right],\;\;\; g,g'\in G(F)$$

\noindent In particular, we have $\sigma(f^P)=i_P^G(\sigma,f)(1,1)$. Since the function $f^P$ is zero if and only if it acts trivially on every representation in $\mathcal{X}_{\tempe}(M)$ (by Theorem \ref{theorem 2.6.1}), we deduce that

\vspace{3mm}

\begin{num}
\item\label{eq 5.3.1} A function $f\in \mathcal{C}(G(F))$ is strongly cuspidal if and only if for every proper parabolic subgroup $P=MU$ and all $\sigma\in \mathcal{X}_{\tempe}(M)$, we have $i_P^G(\sigma,f)(1,1)=0$.
\end{num}

\vspace{3mm}

\noindent Recall that in Section \ref{section 2.6}, we have defined a topological space $\mathcal{C}(\mathcal{X}_{\tempe}(G),\mathcal{E}(G))$ of smooth sections $\pi\in \mathcal{X}_{\tempe}(G)\mapsto T_{\pi}\in \End(\pi)^\infty$ and that the map that associates to $f\in \mathcal{C}(G(F))$ its Fourier transform $\pi\in \mathcal{X}_{\tempe}(G)\mapsto \pi(f)$ induces a topological isomorphism $\mathcal{C}(G(F))\simeq \mathcal{C}(\mathcal{X}_{\tempe}(G),\mathcal{E}(G))$ (Theorem \ref{theorem 2.6.1}). Let us denote by $\gls{Cscusp(Xtemp,E)}$ the image by this isomorphism of $\mathcal{C}_{\scusp}(G(F))$. Then, we have the following

\vspace{2mm}

\begin{lem}\label{lemma 5.3.1}
\begin{enumerate}[(i)]
\item The subspace $\mathcal{C}_{\scusp}(\mathcal{X}_{\tempe}(G),\mathcal{E}(G))$ is stable by multiplication by functions $\varphi\in C_c^\infty(\mathcal{X}_{\tempe}(G))$;
\item The subspace of functions $f\in \mathcal{C}_{\scusp}(G(F))$ having a Fourier transform $\pi\in \mathcal{X}_{\tempe}(G)\mapsto \pi(f)$ that is compactly supported is dense in $\mathcal{C}_{\scusp}(G(F))$.
\end{enumerate} 
\end{lem}

\vspace{2mm}

\noindent\ul{Proof}: 
\begin{enumerate}[(i)]
\item This follows directly from the above characterization of strongly cuspidal functions;

\item There is nothing to say on the $p$-adic case since every function $f\in \mathcal{C}(G(F))$ has a compactly supported Fourier transform. In the real case this follows from (i) once we observe that there exists a sequence $(\varphi_N)_{N\geqslant 1}$ of functions in $C_c^\infty(\mathcal{X}_{\tempe}(G))$ such that

$$\lim\limits_{N\to \infty} \varphi_NT=T$$

\noindent for all $T\in \mathcal{C}(\mathcal{X}_{\tempe}(G),\mathcal{E}(G))$. $\blacksquare$
\end{enumerate}

\subsection{Weighted characters of strongly cuspidal functions}\label{section 5.4}

\noindent Let $M$ be a Levi subgroup of $G$ and $\sigma$ a tempered representation of $M(F)$. Recall that in Section \ref{section 2.5}, we have defined tempered distributions $J_L^Q(\sigma,.)$ on $G(F)$ for all $L\in\mathcal{L}(M)$ and all $Q\in\mathcal{F}(L)$. These distributions depended on the choice of a maximal compact subgroup $K$ which is special in the $p$-adic case and also on the way we normalize intertwining operators (cf.\ Section \ref{section 2.4}).

\begin{lem}\label{lemma 5.4.1}
Let $f\in \mathcal{C}(G(F))$ be a strongly cuspidal function.

\begin{enumerate}[(i)]

\item For all $L\in \mathcal{L}(M)$ and all $Q\in \mathcal{F}(L)$, if $L\neq M$ or $Q\neq G$, then we have

$$J^Q_L(\sigma,f)=0$$

\item The weighted character $J_M^G(\sigma,f)$ doesn't depend on the choice of $K$ or on the way we normalized the intertwining operators;

\item If $\sigma$ is induced from a proper parabolic subgroup of $M$ then

$$J_M^G(\sigma,f)=0$$

\item For all $x\in G(F)$, we have

$$J^G_{xMx^{-1}}(x\sigma x^{-1},f)=J_M^G(\sigma,f)$$ 
\end{enumerate}
\end{lem}

\vspace{3mm}

\noindent\ul{Proof}:

\begin{enumerate}[(i)]

\item First we do the case where $Q=SU_Q$ is different from $G$. Following the definition, we see that

$$J_L^Q(\sigma,f)=J_L^Q(i_M^L(\sigma),f)$$

\noindent Hence, we may assume without loss of generality that $L=M$. We will treat the natural isomorphisms $i_P^G(\sigma_\lambda)\simeq i_{K_P}^K(\sigma_{K_P})$ for $P\in \mathcal{P}(M)$ and $\lambda\in i\mathcal{A}_M^*$, where $K_P=K\cap P(F)$, as identifications. Choose $P\in \mathcal{P}(M)$ such that $P\subset Q$. We have

$$J_M^Q(\sigma,f)=\Tr(\mathcal{R}_M^Q(\sigma,P)i_P^G(\sigma,f))$$

\noindent where the operator $\mathcal{R}_M^Q(\sigma,P)\in \End(i_{K_P}^K(\sigma_{K_P})^\infty)$ is associated with the $(S,M)$-family

$$\mathcal{R}_R^Q(\lambda,\sigma,P)=R_{Q(R)\mid P}(\sigma)^{-1} R_{Q(R)\mid P}(\sigma_\lambda),\;\; R\in \mathcal{P}^{S}(M), \lambda\in i\mathcal{A}_M^*$$

\noindent where $Q(R)=RU_Q\in \mathcal{P}(M)$. Let $K_S$ be the projection of $K_Q=K\cap Q(F)$ onto $S(F)$ and for all $R\in \mathcal{P}^S(M)$, set $K_R=K_S\cap R(F)$. Then $K_S$ is a maximal compact subgroup of $S(F)$ that is special in the $p$-adic case. Hence, we have isomorphisms $i_R^S(\sigma_\lambda)\simeq i_{K_R}^{K_S}(\sigma_{K_R})$ for all $\lambda\in i\mathcal{A}_M^*$ and all $R\in \mathcal{P}^S(M)$. Also, for all $R\in\mathcal{P}^S(M)$, we have the isomorphism of induction by stages $i_{K_{Q(R)}}^K(\sigma_{K_{Q(R)}})\simeq i_{K_Q}^K(i_{K_R}^{K_S}(\sigma_{K_R}))$. In all what follows, we will treat these isomorphisms as identifications. Setting $P_S=P\cap S$, by \ref{eq 2.4.6} we have the equality $R_{Q(R)\mid P}(\sigma_\lambda)=i_{K_Q}^K(R_{R\mid P_S}(\sigma_\lambda))$ for all $\lambda\in i\mathcal{A}_M^*$ and all $R\in \mathcal{P}^S(M)$ (meaning that the $K$-homomorphism $R_{Q(R)\mid P}(\sigma_\lambda)$ is deduced from the $K_S$-homomorphism $R_{R\mid P_S}(\sigma_\lambda):i_{K_{P_S}}^{K_S}(\sigma_{K_{P_S}})^\infty\to i_{K_R}^{K_S}(\sigma_{K_R})^\infty$ by functoriality). We deduce immediately that

\begin{align}\label{eq 5.4.1}
\mathcal{R}_M^Q(\sigma,P)=i_{K_Q}^K\left( \mathcal{R}^S_M(\sigma,P_S)\right)
\end{align}

\noindent where $\mathcal{R}^S_M(\sigma,P_S)\in \End(i_{K_{P_S}}^{K_S}(\sigma_{K_{P_S}})^\infty)=\End(i_{P_S}^S(\sigma)^\infty)$ is associated with the $(S,M)$-family

$$\mathcal{R}_R^S(\lambda,\sigma,P_S)=R_{R\mid P_S}(\sigma)^{-1}R_{R\mid P_S}(\sigma_\lambda),\;\; R\in \mathcal{P}^S(M), \lambda \in i\mathcal{A}_M^*$$

\noindent Recall that we have a natural isomorphism $\End(i_P^G(\sigma))^\infty\simeq i_{Q\times Q}^{G\times G}\left( \End(i_{P_S}^S(\sigma))\right)^\infty$ (cf.\ Section \ref{section 5.3}) which sends the operator $i_P^G(\sigma,f)$ to the function

$$i_P^G(\sigma,f)(g_1,g_2)=i_{P_S}^S\left(\sigma,(L(g_1)R(g_2)f)^Q\right),\;\;\; g_1,g_2\in G(F)$$

\noindent It follows from this and \ref{eq 5.4.1} that the action of the operator $\mathcal{R}_M^Q(\sigma,P)i_P^G(\sigma,f)$ on $i_{K_P}^K(\sigma_{K_P})^\infty\simeq i_{K_Q}^K(i_{K_{P_S}}^{K_S}(\sigma_{K_{P_S}}))^\infty$ is given by

$$\displaystyle \left(\mathcal{R}_M^Q(\sigma,P)i_P^G(\sigma,f)e\right)(k)=\int_{K_Q\backslash K} \mathcal{R}_M^S(\sigma,P_S)i_{P_S}^S(\sigma, (L(k)R(k')f)^Q)e(k')dk'$$

\noindent for all $e\in i_{K_Q}^K(i_{K_{P_S}}^{K_S}(\sigma_{K_{P_S}}))^\infty$ and all $k\in K$. We may now write

\[\begin{aligned}
\displaystyle J_M^Q(\sigma,f) & =\Tr\left(\mathcal{R}_M^Q(\sigma,P)i_P^G(\sigma,f)\right) \\
 & =\int_K \Tr\left(\mathcal{R}_M^S(\sigma,P_S)i_{P_S}^S(\sigma,({}^kf)^Q)\right) dk=0
\end{aligned}\]

\noindent This proves the vanishing (i) in the case $Q\neq G$. Assume now $L\neq M$ but $Q=G$, applying the descent formula \ref{eq 1.9.3}, we see that

$$J_L^G(\sigma,f)=\sum_{L'\in \mathcal{L}(M)} d_M^G(L,L') J_M^{Q'}(\sigma,f)$$

\noindent By what we just saw, the terms in that sum corresponding to $L'\neq G$ vanish. Since $L\neq M$, we also have $d_M^G(L,G)=0$. Hence all terms in the sum above are zero.

\item First we prove the independence in $K$. Let $\widetilde{K}$ be another maximal compact subgroup that is special in the $p$-adic case. Let $P\in \mathcal{P}(M)$. Using $\widetilde{K}$ instead of $K$, we may define another $(G,M)$-family $(\widetilde{\mathcal{R}}_{P'}(\sigma,P))_{P'\in \mathcal{P}(M)}$ taking values in $\End(i_P^G(\sigma)^\infty)$. We deduce from this $(G,M)$-family another weighted character $\widetilde{J}_M^G(\sigma,.)$. For all $P'\in \mathcal{P}(M)$ and all $\lambda\in i\mathcal{A}_M^*$, we have a chain of natural isomorphisms

$$i_{P'}^G(\sigma)\simeq i_{\widetilde{K}_{P'}}^{\widetilde{K}}(\sigma_{\widetilde{K}_{P'}})\simeq i_{P'}^G(\sigma_\lambda)\simeq i_{K_{P'}}^K(\sigma)\simeq i_{P'}^G(\sigma)$$

\noindent where $\widetilde{K}_{P'}=K\cap P'(F)$. We will denote by $I_{P'}(\lambda,\sigma): i_{P'}^G(\sigma)\simeq i_{P'}^G(\sigma)$ their composition and we set $D_{P'}(\lambda,\sigma,P)=R_{P'\mid P}(\sigma)^{-1} I_{P'}(\lambda,\sigma)R_{P'\mid P}(\sigma)$. Then the family $(D_{P'}(\sigma,P))_{P'\in\mathcal{P}(M)}$ is a $(G,M)$-family (taking values in $\End(i_P^G(\sigma))$) and we have

$$\widetilde{\mathcal{R}}_{P'}(\lambda,\sigma,P)=D_{P'}(\lambda,\sigma,P) \mathcal{R}_{P'}(\lambda,\sigma,P)I_P(\lambda,\sigma)$$

\noindent for all $P'\in \mathcal{P}(M)$ and all $\lambda\in i\mathcal{A}_M^*$. We remark that the term $\lambda\mapsto I_P(\lambda,\sigma)$ doesn't depend on $P'$ and satisfies $I_P(0,\sigma)=Id$. Applying the splitting formula \ref{eq 1.9.1}, we get

$$\displaystyle\widetilde{\mathcal{R}}_{M}(\sigma,P)=\sum_{Q\in\mathcal{F}(M)} D'_Q(\sigma,P)\mathcal{R}_M^Q(\sigma,P)$$

\noindent Hence, we have

$$\displaystyle \widetilde{J}_M^G(\sigma,f)=\sum_{Q\in\mathcal{F}(M)} \Tr\left(D'_Q(\sigma,P)\mathcal{R}_M^Q(\sigma,P)i_P^G(\sigma,f)\right)$$

\noindent and the term indexed by $Q=G$ is precisely $J_M^G(\sigma,f)$ (the weighted character defined using $K$). Consequently, it suffices to show the following

\begin{align}\label{eq 5.4.2}
\Tr\left(D'_Q(\sigma,P)\mathcal{R}_M^Q(\sigma,P)i_P^G(\sigma,f)\right)=0
\end{align}

\noindent for all $Q\in \mathcal{F}(M)$ such that $Q\neq G$. Fix such a parabolic subgroup $Q=SU_Q$. Notice that for all $P'\in \mathcal{P}(M)$, we have

$$D'_Q(\sigma,P)=R_{P'\mid P}(\sigma)^{-1}D'_Q(\sigma,P')R_{P'\mid P}(\sigma)$$

$$\mathcal{R}_M^Q(\sigma,P)=R_{P'\mid P}(\sigma)^{-1}\mathcal{R}_M^Q(\sigma,P')R_{P'\mid P}(\sigma)$$

$$i_P^G(\sigma,f)=R_{P'\mid P}(\sigma)^{-1}i_{P'}^G(\sigma,f)R_{P'\mid P}(\sigma)$$

\noindent so that the trace \ref{eq 5.4.2} doesn't change if we replace $P$ by $P'$. Hence, we may assume without loss of generality that $P\subset Q$. The operator $D'_Q(\sigma,P)$ then only depends on the function $\lambda\mapsto D_P(\lambda,\sigma,P)=I_P(\lambda,\sigma)$. We now use again the isomorphism $i_P^G(\sigma)\simeq i_{K_P}^K(\sigma_{K_P})$ as an identification. Direct computation shows that

$$\left(I_P(\lambda,\sigma)e\right)(k)=e^{\langle \lambda,H_M(\widetilde{m}_P(k))\rangle} e(k)$$

\noindent for all $e\in i_{K_P}^K(\sigma_{K_P})$, all $\lambda\in i\mathcal{A}_M^*$ and all $k\in K$, where $\widetilde{m}_P:G(F)\to M(F)$ is any map such that $\widetilde{m}_P(g)^{-1}g\in U(F)\widetilde{K}$ for all $g\in G(F)$. It follows easily that there exists a smooth function $d'_Q(\sigma,P):K\to \mathbb{C}$ such that

$$\left(D'_Q(\sigma,P)e\right)(k)=d'_Q(k,\sigma,P)e(k)$$

\noindent for all $e\in i_{K_P}^K(\sigma_{K_P})$ and all $k\in K$. We have the isomorphism of induction by stage $i_{K_P}^K(\sigma_{K_P})^\infty\simeq i_{K_Q}^K(i_{P_S}^{S}(\sigma)_{K_Q})^\infty$ and we saw during the proof of (i) that $\mathcal{R}_M^Q(\sigma,P)$ is obtained by functoriality from the $K_Q$-endomorphism $\mathcal{R}_M^S(\sigma,P_S)$ of $i_{P_S}^{S}(\sigma)^\infty$. The image of $i_P^G(\sigma,f)$ via the natural isomorphism $\End(i_P^G(\sigma))^\infty\simeq i_{K_Q\times K_Q}^{K\times K}\left(\End(i_{P_S}^{S}\sigma)_{K_Q}\right)^\infty$ is the function given by

$$i_P^G(\sigma,f)(k_1,k_2)=i_{P_S}^S\left(\sigma,(L(k_1)R(k_2)f)^Q\right)$$

\noindent for all $k_1,k_2\in K$. Hence, we see that the operator $D_Q'(\sigma,P)\mathcal{R}_M^Q(\sigma,P)i_P^G(\sigma,f)$ acting on $i_{K_P}^K(\sigma)\simeq i_{K_Q}^K(i_{P_S}^{S}(\sigma)_{K_Q})$ is given by

\[\begin{aligned}
\displaystyle \big(D_Q'(\sigma,P) & \mathcal{R}_M^Q(\sigma,P)i_P^G(\sigma,f)e\big)(k) \\
 & =\int_{K_Q\backslash K} d'_Q(k,\sigma,P) \mathcal{R}_M^S(\sigma,P_S)i_{P_S}^S\left(\sigma,(L(k)R(k')f)^Q\right)e(k') dk'
\end{aligned}\]

\noindent for all $e\in i_{K_Q}^K(i_{P_S}^{S}(\sigma)_{K_Q})$ and all $k\in K$. Consequently, we have

\[\begin{aligned}
\displaystyle \Tr & \left(D_Q'(\sigma,P)\mathcal{R}_M^Q(\sigma,P)i_P^G(\sigma,f)\right) \\
 & =\int_{K_Q\backslash K} d'_Q(k,\sigma,P) \Tr\left(\mathcal{R}_M^S(\sigma,P_S)i_{P_S}^S\left(\sigma,({}^kf)^Q\right)\right) dk=0
\end{aligned}\]

\noindent This proves the vanishing \ref{eq 5.4.2} and ends the proof of the independence in $K$. \\

We now prove that $J_M^G(\sigma,f)$ does not depend on the way we normalized the intertwining operators. Assume we choose different normalization factors $\lambda\mapsto \widetilde{r}_{P'\mid P}(\sigma_\lambda)$, $P,P'\in \mathcal{P}(M)$, yielding new normalized intertwining operators $\widetilde{R}_{P'\mid P}(\sigma_\lambda)$ ($\lambda\in i\mathcal{A}_M^*$, $P,P'\in \mathcal{P}(M)$). Using these new normalized intertwining operators, we construct new $(G,M)$-families $(\widetilde{\mathcal{R}}_{P'}(\sigma,P))_{P'\in\mathcal{P}(M)}$ ($P\in\mathcal{P}(M)$) from which we derive a new weighted character $\widetilde{J}_M^G(\sigma,f)$. Fix $P\in \mathcal{P}(M)$. For all $P'\in \mathcal{P}(M)$, the quotient $r_{P'\mid P}(\sigma_\lambda)\widetilde{r}_{P'\mid P}(\sigma_\lambda)^{-1}$ is well-defined and nonzero for all $\lambda\in i\mathcal{A}_M^*$. Let us set

$$d_{P'}(\lambda,\sigma,P)=r_{P'\mid P}(\sigma)^{-1}\widetilde{r}_{P'\mid P}(\sigma)r_{P'\mid P}(\sigma_\lambda)\widetilde{r}_{P'\mid P}(\sigma_\lambda)^{-1}$$

\noindent for all $P'\in \mathcal{P}(M)$ and all $\lambda\in i\mathcal{A}_M^*$. Then the family $(d_{P'}(\sigma,P))_{P'\in \mathcal{P}(M)}$ is a scalar-valued $(G,M)$-family and we have

$$\widetilde{\mathcal{R}}_{P'}(\lambda,\sigma,P)=d_{P'}(\lambda,\sigma,P)\mathcal{R}_{P'}(\lambda,\sigma,P)$$

\noindent for all $P'\in \mathcal{P}(M)$ and all $\lambda\in i\mathcal{A}_M^*$. Hence by the splitting formula \ref{eq 1.9.1} and the definition of the weighted characters, we have

$$\displaystyle \widetilde{J}_M^G(\sigma,f)=\sum_{Q\in \mathcal{F}(M)} d'_Q J_M^Q(\sigma,f)$$

\noindent The term indexed by $Q=G$ in the above sum is equal to $J_M^G(\sigma,f)$ whereas by (i) all the other terms vanish. Hence we have the equality $\widetilde{J}_M^G(\sigma,f)=J_M^G(\sigma,f)$. This proves indeed that $J_M^G(\sigma,f)$ doesn't depend on the way we normalized the intertwining operators.

\item Assume that there exists a proper Levi subgroup $M_1\subset M$ and a tempered representation $\sigma_1$ of $M_1(F)$ such that $\sigma=i_{M_1}^M(\sigma_1)$. Following the definition, we have

$$J_M^G(\sigma,f)=J_M^G(\sigma_1,f)$$

\noindent and by (i) the right hand side above is zero.

\item By (ii), we may assume that to define $J_{xMx^{-1}}^G(x\sigma x^{-1},f)$ we have used the maximal compact subgroup $xKx^{-1}$ and normalization factors given by $r_{xPx^{-1}\mid xP'x^{-1}}((x\sigma x^{-1})_\lambda)=r_{P\mid P'}(\sigma_\lambda)$ (for $P,P'\in \mathcal{P}(M)$ and $\lambda\in i\mathcal{A}_M^*$). Then by ``transport de structure", we have the equality

$$J_{xMx^{-1}}^G(x\sigma x^{-1},f)=J_M^G(\sigma,f)$$

\noindent $\blacksquare$
\end{enumerate}

\vspace{2mm}

\noindent In Section \ref{section 2.7}, we have defined a set $\underline{\mathcal{X}}(G)$ of virtual tempered representations of $G(F)$. Let $\pi\in\underline{\mathcal{X}}(G)$. Then, there exists a pair $(M,\sigma)$ where $M$ is a Levi subgroup of $G$ and $\sigma\in \underline{\mathcal{X}}_{\elli}(M)$ such that $\pi=i_M^G(\sigma)$. We set

$$\displaystyle \gls{thetafhat}(\pi)=(-1)^{a_G-a_M} J_M^G(\sigma,f)$$

\noindent for all $f\in\mathcal{C}_{\scusp}(G(F))$ (recall that the weighted character $J_M^G(\sigma,.)$ is extended by linearity to all virtual tempered representations). This definition makes sense by the point (iv) of the lemma above since the pair $(M,\sigma)$ is well-defined up to conjugacy.

\vspace{2mm}

\begin{lem}\label{lemma 5.4.2}
\begin{enumerate}[(i)]
\item If $F$ is $p$-adic, then for every compact-open subgroup $K\subseteq G(F)$, there exists a compact subset $\Omega_K\subset \mathcal{X}(G)$ and a continuous semi-norm $\nu_K$ on $\mathcal{C}_K(G(F))$ such that

$$\lvert \widehat{\theta}_f(\pi)\rvert\leqslant \nu_K(f)\mathbf{1}_{\Omega_K}(\pi)$$

\noindent for all $f\in \mathcal{C}_{\scusp,K}(G(F))$ and all $\pi\in \mathcal{X}(G)$.

\item If $F=\mathbb{R}$, then for every integer $k\geqslant 1$ there exists a continuous semi-norm $\nu_k$ on $\mathcal{C}_{\scusp}(G(F))$ such that

$$\lvert \widehat{\theta}_f(\pi)\rvert\leqslant \nu_k(f)N(\pi)^{-k}$$

\noindent for all $f\in \mathcal{C}_{\scusp}(G(F))$ and all $\pi\in\mathcal{X}(G)$.
\end{enumerate}
\end{lem}

\vspace{2mm}

\noindent\ul{Proof}: The point (i) follows from \ref{eq 2.2.3} and \ref{eq 2.6.1} whereas the point (ii) is a consequence of Lemma \ref{lemma 2.5.1} together with \ref{eq 2.3.1}. $\blacksquare$

\subsection{The local trace formulas for strongly cuspidal functions}\label{section 5.5}

\noindent Let us set

$$\displaystyle \gls{KAff'}(g_1,g_2)=\int_{G(F)} f(g_1^{-1}gg_2)f'(g)dg$$

\noindent for all $f,f'\in \mathcal{C}(G(F))$ and all $g_1,g_2\in G(F)$. The integral above is absolutely convergent by Proposition \ref{proposition 1.5.1}(v). We also define

$$\displaystyle \gls{KAff'}(x,x)=\int_{\mathfrak{g}(F)} f(x^{-1}Xx)f'(X)dX$$

\noindent for all $f,f'\in \mathcal{S}(\mathfrak{g}(F))$ and all $x\in A_G(F)\backslash G(F)$.

\vspace{2mm}

\noindent The two theorems below are slight variations around the local trace formula of Arthur (cf.\ \cite{A1}) and its version for Lie algebras due to Waldspurger (cf.\ \cite{Wa3}). The proof of these two theorems will appear elsewhere \cite{Beu2}.

\vspace{3mm}

\begin{theo}\label{theorem 5.5.1}
\begin{enumerate}[(i)]

\item For all $d\geqslant 0$, there exists $d'\geqslant 0$ and a continuous semi-norm $\nu_{d,d'}$ on $\mathcal{C}(G(F))$ such that

$$\lvert K_{f,f'}^A(g_1,g_2)\rvert\leqslant \nu_{d,d'}(f)\nu_{d,d'}(f')\Xi^G(g_1)\sigma(g_1)^{-d}\Xi^G(g_2)\sigma(g_2)^{d'}$$

\noindent and 

$$\lvert K_{f,f'}^A(g_1,g_2)\rvert\leqslant \nu_{d,d'}(f)\nu_{d,d'}(f')\Xi^G(g_1)\sigma(g_1)^{d'}\Xi^G(g_2)\sigma(g_2)^{-d}$$

\noindent for all $f,f'\in \mathcal{C}(G(F))$ and for all $g_1,g_2\in G(F)$.

\item For all $d\geqslant 0$, there exists a continuous semi-norm $\nu_d$ on $\mathcal{C}(G(F))$ such that

$$\lvert K_{f,f'}^A(x,x)\rvert \leqslant \nu_d(f)\nu_d(f')\Xi^G(x)^2 \sigma_{A_G\backslash G}(x)^{-d}$$

\noindent for all $f\in \mathcal{C}_{\scusp}(G(F))$, all $f'\in \mathcal{C}(G(F))$ and all $x\in A_G(F)\backslash G(F)$.

\item Let $f,f'\in \mathcal{C}(G(F))$ with $f$ strongly cuspidal. Then, there exists $c>0$ such that for all $d\geqslant 0$ there exists $d'\geqslant 0$ such that

$$\lvert K_{f,f'}^A(g,hg)\rvert\ll \Xi^G(g)^2\sigma_{A_G\backslash G}(g)^{-d} e^{c\sigma(h)}\sigma(h)^{d'}$$

\noindent for all $h,g\in G(F)$.

\end{enumerate}

\vspace{2mm}

\noindent By the point (ii), the function $x\in A_G(F)\backslash G(F)\mapsto K^A_{f,f'}(x,x)$ is integrable as soon as $f$ is strongly cuspidal. We set

$$\displaystyle J^A(f,f')=\int_{A_G(F)\backslash G(F)} K^A_{f,f'}(x,x)dx$$

\noindent for all $f\in \mathcal{C}_{\scusp}(G(F))$ and all $f'\in \mathcal{C}(G(F))$.

\begin{enumerate}[(i)]
\setcounter{enumi}{3}

\item We have the geometric expansion

$$\displaystyle J^A(f,f')=\int_{\Gamma(G)} D^G(x)^{1/2}\theta_f(x)J_G(x,f')dx$$

\noindent for all $f\in \mathcal{C}_{\scusp}(G(F))$ and all $f'\in \mathcal{C}(G(F))$, the integral above being absolutely convergent.

\item We have the spectral expansion

$$\displaystyle J^A(f,f')=\int_{\mathcal{X}(G)} D(\pi)\widehat{\theta}_f(\pi) \theta_{\overline{\pi}}(f') d\pi$$

\noindent for all $f\in \mathcal{C}_{\scusp}(G(F))$ and all $f'\in \mathcal{C}(G(F))$, the integral above being absolutely convergent.
\end{enumerate}
\end{theo}

\vspace{3mm}

\begin{theo}\label{theorem 5.5.2}
\begin{enumerate}[(i)]

\item For all $N\geqslant 0$, there exists a continuous semi-norm $\nu_N$ on $\mathcal{S}(\mathfrak{g}(F))$ such that

$$\lvert K_{f,f'}^A(x,x)\rvert \leqslant \nu_N(f)\nu_N(f')\lVert x\rVert_{A_G\backslash G}^{-N}$$

\noindent for all $f\in \mathcal{S}_{\scusp}(\mathfrak{g}(F))$, all $f'\in \mathcal{S}(\mathfrak{g}(F))$ and all $x\in A_G(F)\backslash G(F)$.

\end{enumerate}

\vspace{2mm}

\noindent In particular the function $x\in A_G(F)\backslash G(F)\mapsto K^A_{f,f'}(x,x)$ is integrable as soon as $f$ is strongly cuspidal. We set

$$\displaystyle J^A(f,f')=\int_{A_G(F)\backslash G(F)} K^A_{f,f'}(x,x)dx$$

\noindent for all $f\in \mathcal{S}_{\scusp}(\mathfrak{g}(F))$ and all $f'\in \mathcal{S}(\mathfrak{g}(F))$.

\begin{enumerate}[(i)]
\setcounter{enumi}{1}

\item We have the ``geometric" expansion

$$\displaystyle J^A(f,f')=\int_{\Gamma(\mathfrak{g})} D^G(X)^{1/2}\theta_f(X)J_G(X,f')dX$$

\noindent for all $f\in \mathcal{S}_{\scusp}(\mathfrak{g}(F))$ and all $f'\in \mathcal{S}(\mathfrak{g}(F))$, the integral above being absolutely convergent.

\item We have the ``spectral" expansion

$$\displaystyle J^A(f,f')=\int_{\Gamma(\mathfrak{g})} D^G(X)^{1/2}\theta_{\widehat{f}}(X) J_G(-X,\widehat{f'}) dX$$

\noindent for all $f\in \mathcal{S}_{\scusp}(\mathfrak{g}(F))$ and all $f'\in \mathcal{S}(\mathfrak{g}(F))$, the integral above being absolutely convergent.
\end{enumerate}
\end{theo}

\subsection{Strongly cuspidal functions and quasi-characters}\label{section 5.6}

\begin{prop}\label{proposition 5.6.1}
\begin{enumerate}[(i)]
\item For all $f\in \mathcal{S}_{\scusp}(\mathfrak{g}(F))$, the function $\theta_f$ is a Schwartz quasi-character and we have $\widehat{\theta}_f=\theta_{\widehat{f}}$. Moreover, if $G$ admits an elliptic maximal torus, then the linear map

$$\mathcal{S}_{\scusp}(\mathfrak{g}(F))\to SQC(\mathfrak{g}(F))$$
$$f\mapsto \theta_f$$

\noindent has dense image and for every completely $G(F)$-invariant open subset $\omega\subseteq \mathfrak{g}(F)$ which is relatively compact modulo conjugation, the linear map

$$\mathcal{S}_{\scusp}(\omega)\to QC_c(\omega)$$
$$f\mapsto \theta_f$$

\noindent also has dense image.

\item Let $f\in \mathcal{C}_{\scusp}(G(F))$. Then, the function $\theta_f$ is a quasi-character on $G(F)$ and we have an equality of quasi-characters

$$\displaystyle\theta_f=\int_{\mathcal{X}(G)}D(\pi) \widehat{\theta}_f(\pi)\overline{\theta_{\pi}}d\pi$$

\noindent where the integral above is absolutely convergent in $QC(G(F))$.
\end{enumerate}
\end{prop}

\vspace{2mm}

\noindent\ul{Proof}:

\begin{enumerate}[(i)]
\item If $F$ is $p$-adic, all of these is contained in \cite{Wa1} (Note that for $p$-adic group there always exists a maximal elliptic torus). Let us assume now that $F=\mathbb{R}$. Let $f\in \mathcal{S}_{\scusp}(\mathfrak{g}(F))$. By the estimates of Lemma \ref{lemma 5.2.2}(i), the function $\theta_f$ satisfies the assumption of Proposition \ref{proposition 4.2.1} 1.(i). Hence, by this proposition, there exists a quasi-character $\widehat{\theta}_f$ such that $\widehat{T}_{\theta_f}=T_{\widehat{\theta}_f}$. By Theorem \ref{theorem 5.5.2}(ii) and (iii) and the Weyl integration formula, we have

$$\displaystyle \int_{\mathfrak{g}(F)} \theta_f(X)\widehat{f'}(X)dX=\int_{\mathfrak{g}(F)}\theta_{\widehat{f}}(X) f'(X)dX$$

\noindent for all $f'\in C_c^\infty(\mathfrak{g}(F))$. It follows that $\widehat{\theta}_f=\theta_{\widehat{f}}$. Applying this to the inverse Fourier transform of $f$, we see that $\theta_f$ is a quasi-character. In particular for all $u\in I(\mathfrak{g})$, we have $\partial(u)T_{\theta_f}=T_{\partial(u)\theta_f}$. For all $u\in I(\mathfrak{g})$, $\partial(u)f$ is also strongly cuspidal and by Lemma \ref{lemma 5.2.2}(i) we have $\partial(u)\theta_f=\theta_{\partial(u)f}$. Hence, applying the estimates of Lemma \ref{lemma 5.2.2}(i) to the functions $\theta_{\partial(u)f}$, $u\in I(\mathfrak{g})$, we see that $\theta_f$ satisfies the assumption of Proposition \ref{proposition 4.2.1} 1.(ii). It follows from this proposition that $\theta_{f}$ is a Schwartz quasi-character.

\vspace{2mm}

\noindent Let us now assume that $G$ admits an elliptic maximal torus, hence $\mathfrak{g}_{\reg}(F)_{\elli}\neq\emptyset$. We first show that the linear map

\begin{align}\label{eq 5.6.1}
\mathcal{S}_{\scusp}(\mathfrak{g}(F))\to SQC(\mathfrak{g}(F))
\end{align}
$$f\mapsto \theta_f$$

\noindent has dense image. We start by proving the following

\vspace{3mm}

\begin{num}
\item\label{eq 5.6.2} For all $X\in \mathfrak{g}_{\reg}(F)$, there exists $f\in \mathcal{S}_{\scusp}(\mathfrak{g}(F))$ such that $\theta_f(X)\neq 0$.
\end{num}

\vspace{3mm}

\noindent Let $X\in \mathfrak{g}_{\reg}(F)$. Every function $f\in \mathcal{S}(\mathfrak{g}(F))$ which is supported in $\mathfrak{g}_{\reg}(F)_{\elli}$ is strongly cuspidal. Let $f$ be such a function. Since $\widehat{\theta}_f=\theta_{\widehat{f}}$, by Lemma \ref{lemma 4.2.2}(iii), we have

\begin{align}\label{eq 5.6.3}
\displaystyle \theta_{\widehat{f}}(X)=\int_{\Gamma_{\elli}(\mathfrak{g})}D^G(Y)^{1/2}\theta_{f}(Y)\widehat{j}(Y,X)dY
\end{align}

\noindent By \ref{eq 1.7.4}, there exist $Y_0\in \mathfrak{g}_{\reg}(F)_{\elli}$ such that $\widehat{j}(Y_0,X)\neq 0$. Now, the term $\theta_f(Y_0)$ is just the orbital integral of $f$ at $Y_0$ and it is not hard to see that we may choose $f$ such that this orbital integral is nonzero (just take $f\in C_c^\infty(\mathfrak{g}(F)_{\elli})$ positive and such that $f(Y_0)\neq 0$). Up to multiplying $f$ by a well chosen invariant function $\varphi\in C^\infty(\mathfrak{g}(F))^G$ that is positive, equals $1$ near $Y_0$ and is supported in a small compact modulo conjugation invariant neighborhood of $Y_0$, we see that we may arrange the right hand side of \ref{eq 5.6.3} to be nonzero. This proves \ref{eq 5.6.2}.

\vspace{2mm}

\noindent We now prove the following

\vspace{3mm}

\begin{num}
\item\label{eq 5.6.4} For all $\theta\in SQC(\mathfrak{g}(F))$ and every integer $N\geqslant 1$, there exists a constant $c_N>0$ such that for every invariant and compact modulo conjugation subset $L\subseteq \mathfrak{g}_{\reg}(F)$ there exists a function $f\in \mathcal{S}_{\scusp}(\mathfrak{g}(F))$ such that
$$D^G(X)^{1/2}\lvert \theta(X)-\theta_f(X)\rvert\leqslant c_N\lVert X\rVert_{\Gamma(\mathfrak{g})}^{-N} \mathbf{1}_{L^c}(X)$$
for all $X\in \mathfrak{g}_{\reg}(F)$ (where $\mathbf{1}_{L^c}$ denotes the characteristic function of $\mathfrak{g}_{\reg}(F)-L$).
\end{num}

\vspace{3mm}

\noindent From \ref{eq 5.6.2} and the existence of smooth invariant partition of unity (Proposition \ref{proposition 3.1.1}(ii)), we easily deduce that for all $\theta\in QC_c(\mathfrak{g}(F))$ whose support is contained in $\mathfrak{g}_{\reg}(F)$ there exists $f\in \mathcal{S}_{\scusp}(\mathfrak{g}(F))$ such that $\theta=\theta_f$. Now let $\theta\in SQC(\mathfrak{g}(F))$ and $L\subseteq \mathfrak{g}_{\reg}(F)$ be an invariant compact modulo conjugation subset. Choose $\varphi\in C^\infty(\mathfrak{g}(F))^G$ such that $0\leqslant \varphi\leqslant 1$, $\varphi=1$ on $L$ and the support of $\varphi$ is contained in $\mathfrak{g}_{\reg}(F)$ and compact modulo conjugation. Then, by what we just saw, there exists $f\in \mathcal{S}_{\scusp}(\mathfrak{g}(F))$ such that $\theta_f=\varphi\theta$. Since 

$$\left\lvert\theta(X)-\theta_f(X)\right\rvert=\left\lvert(1-\varphi(X))\theta(X)\right\rvert\leqslant \lvert \theta(X)\rvert$$

\noindent for all $X\in \mathfrak{g}_{\reg}(F)$, we have

$$D^G(X)^{1/2}\lvert \theta(X)-\theta_f(X)\rvert\leqslant c_N\lVert X\rVert_{\Gamma(\mathfrak{g})}^{-N} \mathbf{1}_{L^c}(X)$$

\noindent for all $N\geqslant 1$ and all $X\in \mathfrak{g}_{\reg}(F)$, where

$$c_N=\sup_{X\in \mathfrak{g}_{\reg}(F)} \lVert X\rVert_{\Gamma(\mathfrak{g})}^N D^G(X)^{1/2}\lvert \theta(X)\rvert$$

\noindent This ends the proof of \ref{eq 5.6.4}.

\vspace{2mm}

\noindent With notation of Section \ref{section 4.2}, we set $q_u=q_{u,0}$ for all $u\in I(\mathfrak{g})$ (these are continuous semi-norms on $SQC(\mathfrak{g}(F))$). We now deduce from \ref{eq 5.6.4} the following

\vspace{3mm}

\begin{num}
\item \label{eq 5.6.5} For all $\theta\in SQC(\mathfrak{g}(F))$, every finite family $\{u_1,\ldots,u_k\}\subset I(\mathfrak{g})$ and all $\epsilon>0$, there exists $f\in \mathcal{S}_{\scusp}(\mathfrak{g}(F))$ such that
$$q_{u_i}(\theta-\theta_f)\leqslant \epsilon$$
for all $1\leqslant i\leqslant k$.
\end{num}

\vspace{3mm}

\noindent Let $\theta\in SQC(\mathfrak{g}(F))$ and let us fix a finite family $\{u_1,\ldots,u_k\}\subset I(\mathfrak{g})$. By Lemma \ref{lemma 4.2.2}(iii), for all $f\in \mathcal{S}_{\scusp}(\mathfrak{g}(F))$, we have

$$\displaystyle\widehat{\theta}(X)-\theta_{\widehat{f}}(X)=\int_{\Gamma(\mathfrak{g})} D^G(Y)^{1/2}\left(\theta(Y)-\theta_f(Y)\right) \widehat{j}(Y,X) dY$$

\noindent for all $X\in \mathfrak{g}_{\reg}(F)$. For all $1\leqslant i\leqslant k$, let $p_i\in I(\mathfrak{g}^*)$ be such that $u_{p_i}=u_i$. Applying the above equality to $p_{i}\theta$ and $p_{i}f$ for all $1\leqslant i\leqslant k$, we obtain

$$\displaystyle\partial(u_i)\left(\widehat{\theta}-\theta_{\widehat{f}}\right)(X)=\int_{\Gamma(\mathfrak{g})} D^G(Y)^{1/2}p_{u_i}(Y)\left(\theta(Y)-\theta_f(Y)\right) \widehat{j}(Y,X) dY$$

\noindent for all $X\in \mathfrak{g}_{\reg}(F)$, all $f\in \mathcal{S}_{\scusp}(\mathfrak{g}(F))$ and all $1\leqslant i\leqslant k$. Hence, by \ref{eq 1.7.3}, we get

$$\displaystyle q_{u_i}(\widehat{\theta}-\theta_{\widehat{f}})\ll \int_{\Gamma(\mathfrak{g})} \lvert p_{u_i}(Y)\rvert D^G(Y)^{1/2} \left\lvert \theta(Y)-\theta_f(Y)\right\rvert dY$$

\noindent for all $f\in \mathcal{S}_{\scusp}(\mathfrak{g}(F))$ and all $1\leqslant i\leqslant k$. Of course, there exists $N_0\geqslant 1$ such that $\lvert p_{u_i}(Y)\rvert\ll \lVert Y\rVert_{\Gamma(\mathfrak{g})}^{N_0}$ for all $Y\in \mathfrak{g}(F)$ and $1\leqslant i\leqslant k$. Hence, it follows from \ref{eq 5.6.4} that for all $N\geqslant 1$ we have an inequality 

\begin{align}\label{eq 5.6.6}
\displaystyle \inf_{f\in \mathcal{S}_{\scusp}(\mathfrak{g}(F))} \sup_{i=1,\ldots,k} q_{u_i}(\widehat{\theta}-\theta_{\widehat{f}})\ll \inf_L\int_{\Gamma(\mathfrak{g})} \lVert Y\rVert_{\Gamma(\mathfrak{g})}^{-N} \mathbf{1}_{L^c}(Y)dY
\end{align}

\noindent where $L$ runs through the invariant and compact modulo conjugation subsets of $\mathfrak{g}_{\reg}(F)$. Choose $N\geqslant 1$ such that the function $Y\mapsto \lVert Y\rVert_{\Gamma(\mathfrak{g})}^{-N}$ is integrable on $\Gamma(\mathfrak{g})$. Then it is not hard to see that

$$\displaystyle \inf_{L} \int_{\Gamma(\mathfrak{g})} \lVert Y\rVert_{\Gamma(\mathfrak{g})}^{-N} \mathbf{1}_{L^c}(Y)dY=0$$

\noindent Hence, replacing $\theta$ by $\widehat{\theta}$, \ref{eq 5.6.5} follows from the inequality \ref{eq 5.6.6} above.

\vspace{2mm}

\noindent We are now in position to prove that the image of the linear map \ref{eq 5.6.1} is dense. By Lemma \ref{lemma 4.2.2}(v), it suffices to prove that the intersection

\begin{align}\label{eq 5.6.7}
QC_c(\mathfrak{g}(F))\cap \{\theta_f;\;\; f\in \mathcal{S}_{\scusp}(\mathfrak{g}(F))\}
\end{align}

\noindent is dense in $QC_c(\mathfrak{g}(F))$ (for its own topology). Let $\theta\in QC_c(\mathfrak{g}(F))$. Then, it follows directly from \ref{eq 5.6.5} that we may find a sequence $(f_n)_{n\geqslant 1}$ in $\mathcal{S}_{\scusp}(\mathfrak{g}(F))$ such that

$$\lim\limits_{n\to \infty}\theta_{f_n}=\theta$$

\noindent in $QC(\mathfrak{g}(F))$. Let $\varphi\in C^\infty(\mathfrak{g}(F))^G$ be compactly supported modulo conjugation and such that $\varphi=1$ on $\Supp(\theta)$. Then, by Lemma \ref{lemma 4.2.2}(iv) and the closed graph theorem, we have

$$\lim\limits_{n\to \infty}\varphi\theta_{f_n}=\varphi\theta=\theta$$

\noindent in $QC_c(\mathfrak{g}(F))$. But $\varphi\theta_{f_n}=\theta_{\varphi f_n}$ and $\varphi f_n\in \mathcal{S}_{\scusp}(\mathfrak{g}(F))$ for all $n\geqslant 1$ (by Proposition \ref{proposition 3.1.1}(iv)). Hence, the subspace \ref{eq 5.6.7} is indeed dense in $QC_c(\mathfrak{g}(F))$ and this ends the proof that the linear map \ref{eq 5.6.1} has dense image. Let $\omega\subseteq \mathfrak{g}(F)$ be a completely $G(F)$-invariant open subset which is relatively compact modulo conjugation. The argument we just used actually also show that the linear map

$$\mathcal{S}_{\scusp}(\omega)\to QC_c(\omega)$$
$$f\mapsto \theta_f$$

\noindent has dense image since for $\theta\in QC_c(\omega)$, we can choose $\varphi$ as before which is supported in $\omega$, hence the functions $\varphi f_n$, for $n\geqslant 1$, will belong to $\mathcal{S}_{\scusp}(\omega)$. This ends the proof of (i).

\item Let $f\in \mathcal{C}_{\scusp}(G(F))$. By Theorem \ref{theorem 5.5.1}(iv) and (v) and the Weyl integration formula, we have

\[\begin{aligned}
\displaystyle \int_{G(F)} \theta_f(x)f'(x)dx & =\int_{\mathcal{X}(G)}D(\pi)\widehat{\theta}_f(\pi)\theta_{\overline{\pi}}(f')d\pi \\
 & =\int_{\mathcal{X}(G)}D(\pi)\widehat{\theta}_f(\pi)\int_{G(F)} \overline{\theta_{\pi}(x)}f'(x)dxd\pi
\end{aligned}\]

\noindent for all $f'\in C_c^\infty(G(F))$. By Proposition \ref{proposition 4.8.1}(ii), Lemma \ref{lemma 5.4.2} and \ref{eq 1.8.2} the above double integral is absolutely convergent. It follows that

$$\displaystyle \int_{G(F)} \theta_f(x)f'(x)dx=\int_{G(F)}\int_{\mathcal{X}(G)}D(\pi)\widehat{\theta}_f(\pi) \overline{\theta_{\pi}(x)}d\pi f'(x)dx$$

\noindent for all $f'\in C_c^\infty(G(F))$. Hence, we have

$$\displaystyle \theta_f(x)=\int_{\mathcal{X}(G)}D(\pi)\widehat{\theta}_f(\pi) \overline{\theta_{\pi}(x)}d\pi$$

\noindent for almost all $x\in G_{\reg}(F)$. Consequently, to prove the point (ii), it suffices to show that the integral

$$\displaystyle \int_{\mathcal{X}(G)} D(\pi)\widehat{\theta}_f(\pi)\overline{\theta_{\pi}}d\pi$$

\noindent is absolutely convergent in $QC(G(F))$. But this follows easily from Proposition \ref{proposition 4.8.1}(ii) combined with Lemma \ref{lemma 5.4.2}. $\blacksquare$

\end{enumerate}

\subsection{Lifts of strongly cuspidal functions}\label{section 5.7}

\begin{prop}\label{proposition 5.7.1}
Let $x\in G(F)_{\elli}$ be elliptic and let $\Omega_x\subseteq G_x(F)$ be a $G$-good open neighborhood of $x$ which is relatively compact modulo conjugation. Set $\Omega=\Omega_x^G$. Then, there exists a linear map

$$\mathcal{S}_{\scusp}(\Omega_x)\to \mathcal{S}_{\scusp}(\Omega)$$
$$f\mapsto \widetilde{f}$$

\noindent such that

\begin{enumerate}[(i)]
\item For all $f\in \mathcal{S}_{\scusp}(\Omega_x)$, we have

$$\displaystyle (\theta_{\widetilde{f}})_{x,\Omega_x}=\sum_{z\in Z_G(x)(F)/G_x(F)} {}^z\theta_f$$

\item There exists a function $\alpha\in C_c^\infty(Z_G(x)(F)\backslash G(F))$ satisfying

$$\displaystyle \int_{Z_G(x)(F)\backslash G(F)} \alpha(g)dg=1$$

\noindent and such that for all $f\in \mathcal{S}_{\scusp}(\Omega_x)$ and all $g\in G(F)$, there exists $z\in Z_G(x)(F)$ such that

$$({}^{zg}\widetilde{f})_{x,\Omega_x}=\alpha(g)f$$
\end{enumerate}
\end{prop}

\vspace{2mm}

\noindent\ul{Proof}: Let us denote by $\pi: G(F)\to Z_G(x)(F)\backslash G(F)$ the natural projection. It is an $F$-analytic locally trivial fibration. Let us fix an open subset $\mathcal{U}\subseteq Z_G(x)(F)\backslash G(F)$ and an $F$-analytic section

$$s:\mathcal{U}\to G(F)$$

\noindent Since $\Omega_x$ is a $G$-good open subset, the map

$$\beta:\Omega_x\times \mathcal{U}\to G(F)$$
$$(y,g)\mapsto \beta(y,g)=s(g)^{-1}ys(g)$$

\noindent is an open embedding of $F$-analytic spaces. For all $f\in \mathcal{S}(\Omega_x)$ and all $\varphi\in C_c^\infty(\mathcal{U})$, we define a function $f_\varphi$ on $G(F)$ by

$$\displaystyle f_\varphi(\gamma)=\left\{
    \begin{array}{ll}
        f(y)\varphi(g) & \mbox{ if } \gamma=s(g)^{-1}ys(g) \mbox{ for some } g\in \mathcal{U} \mbox{ and some } y\in \Omega_x \\
        0 & \mbox{ otherwise.}
    \end{array}
\right.
$$

\noindent Let us now prove the following

\vspace{3mm}

\begin{num}
\item\label{eq 5.7.1} For all $\varphi\in C_c^\infty(\mathcal{U})$ and all $f\in \mathcal{S}(\Omega_x)$, the function $f_\varphi$ belongs to $\mathcal{S}(\Omega)$.
\end{num}

\vspace{3mm}

\noindent This is almost straightforward in the $p$-adic case. Assume that $F=\mathbb{R}$. Let $\varphi\in C_c^\infty(\mathcal{U})$ and $f\in \mathcal{S}(\Omega_x)$. The only thing which is not obvious is to prove that for all $u\in \mathcal{U}(\mathfrak{g})$ the function $L(u)f_\varphi$ is rapidly decreasing. This follows at once from the two following facts

\vspace{3mm}

\begin{num}
\item\label{eq 5.7.2} The function $f_\varphi$ is rapidly decreasing i.e., for all $N\geqslant 1$ we have
$$\left\lvert f_\varphi(\gamma)\right\rvert\ll \lVert \gamma\rVert^{-N}$$
for all $\gamma\in G(F)$.

\item\label{eq 5.7.3} For all $X\in \mathfrak{g}(F)$, there exist an integer $k\geqslant 1$, functions $\varphi_1,\ldots,\varphi_k\in C_c^\infty(\mathcal{U})$ and Schwartz functions $f_1,\ldots,f_k\in\mathcal{S}(\Omega_x)$ such that
$$L(X)f_\varphi=(f_1)_{\varphi_1}+\ldots+(f_k)_{\varphi_k}$$
\end{num}

\vspace{3mm}

\noindent The claim \ref{eq 5.7.2} follows from the fact that there exist constants $c_1,c_2>0$ and integers $N_1,N_2\geqslant 1$ such that

$$\displaystyle c_1^{-1}\lVert y\rVert^{1/N_1}\leqslant \left\lVert s(g)^{-1}y s(g)\right\rVert\leqslant c_2\lVert y\rVert^{N_2}$$

\noindent for all $g\in \Supp(\varphi)$ and all $y\in \Omega_x$ together with the fact that $f$ is itself rapidly decreasing. Let us focus on \ref{eq 5.7.3} now. Fix $X\in \mathfrak{g}(F)$. Using a partition of unity if necessary, we may assume that $\mathcal{U}$ is parallelizable (i.e., its tangent bundle is trivial). Assume this is so and let us fix a trivialization of the tangent bundle of $\mathcal{U}$

$$T\mathcal{U}\simeq \mathcal{U}\times V$$

\noindent where $V$ is some finite dimensional $\mathbb{R}$-vector space. Let us also fix trivializations $TG(F)\simeq G(F)\times \mathfrak{g}(F)$ and $TG_x(F)\simeq G_x(F)\times \mathfrak{g}_x(F)$ using right translations. For all $\gamma\in G(F)$, we have

\begin{align}\label{eq 5.7.4}
\displaystyle \left( L(X)f_\varphi\right) & (\gamma)= \\
 \nonumber & \left\{
    \begin{array}{ll}
        \left[\partial\left(d\beta_{(y,g)}^{-1}(X)\right)(f\otimes \varphi)\right](y,g) & \mbox{ if } \gamma=s(g)^{-1}ys(g) \mbox{ for some } g\in \mathcal{U},\; y\in \Omega_x \\
        0 & \mbox{ otherwise.}
    \end{array}
\right.
\end{align}

\noindent where $d\beta_{(y,g)}:\mathfrak{g}_X(F)\oplus V\simeq \mathfrak{g}(F)$ denotes the differential of $\beta$ at $(y,g)$. A painless computation shows that

$$\displaystyle d\beta_{(y,g)}(Y,Z)=\Ad(s(g))^{-1}\left[ Y+\left(\Ad(y)-1\right)ds_g(Z)\right]$$

\noindent for all $(y,g)\in \Omega_x\times \mathcal{U}$ and all $(Y,Z)\in \mathfrak{g}_x(F)\oplus V$, where $ds_g: V\to \mathfrak{g}(F)$ denotes the differential at $g$ of the section $s$. It is obvious from this description that there exist an integer $r\geqslant 1$, polynomials $P_1,\ldots,P_r\in \mathbb{R}[G_x]$, smooth functions $\psi_1,\ldots,\psi_r\in C^\infty(\mathcal{U})$ and vectors $X_1,\ldots,X_r\in \mathfrak{g}_x(F)\oplus V$ such that

\begin{align}\label{eq 5.7.5}
\displaystyle d\beta_{(y,g)}^{-1}(X)=\eta_x^{alg}(y)^{-1}\sum_{i=1}^r P_i(y)\psi_i(g) X_i
\end{align}

\noindent for all $(y,g)\in \Omega_x\times \mathcal{U}$, where we have set $\eta_x^{alg}(y)=\det\left(1-\Ad(y)\right)_{\mid \mathfrak{g}/\mathfrak{g}_x}$. Writing $X_i=Y_i+Z_i$ where $Y_i\in \mathfrak{g}_x(F)$ and $Z_i\in V$ for all $i=1,\ldots,r$, we get from \ref{eq 5.7.4} and \ref{eq 5.7.5} that

$$\displaystyle \left(L(X)f_\varphi\right)=\sum_{i=1}^r \left(f'_i\right)_{\varphi'_i}+\left(f''_i\right)_{\varphi''_i}$$

\noindent where $f'_i=(\eta_x^{alg})^{-1} P_i \left(L(Y_i)f\right)$, $f_i''=(\eta_x^{alg})^{-1}P_i f$, $\varphi_i'=\psi_i\varphi$ and $\varphi''_i=\psi_i \left(\partial(Z_i)\varphi\right)$ for $i=1,\ldots,r$. The claim \ref{eq 5.7.3} follows once we remark that multiplication by $(\eta_x^{alg})^{-1}$ preserves $\mathcal{S}(\Omega_x)$ (this is a consequence of Proposition \ref{proposition 3.1.1}(iv), here we use the fact that $\Omega_x$ is relatively compact modulo conjugation so that functions in $\mathcal{S}(\Omega_x)$ are compactly supported modulo conjugation).

\vspace{2mm}

\noindent We now construct a linear form as in the proposition. The construction is as follows. Choose a function $\alpha\in C_c^\infty(\mathcal{U})$ such that

$$\displaystyle \int_{Z_G(x)(F)\backslash G(F)} \alpha(g) dg=1$$

\noindent and set $\widetilde{f}=\left(\left(\eta_x^G\right)^{-1/2}f\right)_\alpha$ for all $f\in \mathcal{S}_{\scusp}(\Omega_x)$. The second point of the proposition is obvious from this definition. We need to check the two following facts

\vspace{3mm}

\begin{num}
\item\label{eq 5.7.6} For all $f\in \mathcal{S}_{\scusp}(\Omega_x)$, the function $\widetilde{f}$ is strongly cuspidal.

\item\label{eq 5.7.7} For all $f\in \mathcal{S}_{\scusp}(\Omega_x)$, we have
$$\displaystyle (\theta_{\widetilde{f}})_{x,\Omega_x}=\sum_{z\in Z_G(x)(F)/G_x(F)} {}^z\theta_f$$
\end{num}

\vspace{3mm}

\noindent Let us prove \ref{eq 5.7.6}. Let $f\in \mathcal{S}_{\scusp}(\Omega_x)$ and let $P=MU$ be a proper parabolic subgroup of $G$. Let $m\in M(F)\cap G_{\reg}(F)$. We want to show that

$$\displaystyle \int_{U(F)} \widetilde{f}(u^{-1}mu)du=0$$

\noindent If the conjugacy class of $m$ does not meet $\Omega_x$ then the function $u\in U(F)\mapsto \widetilde{f}(u^{-1}mu)$ is identically zero. Assume now that $m=g^{-1}m_xg$ for some $g\in G(F)$ and some $m_x\in \Omega_x$. Set $M'=gMg^{-1}$, $U'=gUg^{-1}$ and $P'=gPg^{-1}=M'U'$. Then, $x$ belongs to $M'(F)$ and $P'_x=M'_xU'_x$ is a parabolic subgroup of $G_x$ which is proper (since $x$ is elliptic). We may now write

\[\begin{aligned}
\displaystyle \int_{U(F)} \widetilde{f}(u^{-1}mu)du & =\int_{U'(F)} {}^g\widetilde{f}({u'}^{-1}m_xu')du' \\
& =\int_{U'_x(F)\backslash U'(F)} \int_{U'_x(F)} ({}^{u'g}\widetilde{f})_{x,\Omega_x}({u'}_x^{-1}m_xu'_x)du'_xdu'
\end{aligned}\]

\noindent By (ii), for all $u'\in U'(F)$ there exists $z\in Z_G(x)(F)$ such that the function $({}^{u'g}\widetilde{f})_{x,\Omega_x}$ is a scalar multiple of ${}^zf$. Since $f$ is strongly cuspidal, it follows that the inner integral above is zero. This ends the proof of \ref{eq 5.7.6}.

\vspace{2mm}

\noindent We now prove \ref{eq 5.7.7}. Let $y\in \Omega_{x,\reg}$. Set $M(y)=Z_G(A_{G_y})$ and $M(y)_x=Z_{G_x}(A_{G_y})$. Returning to the definitions, we need to show that

\begin{align}\label{eq 5.7.8}
\displaystyle J_{M(y)}(y,\widetilde{f})=\sum_{z\in G_x(F)\backslash Z_G(x)(F)} J_{M(y)_x}(y,{}^zf)
\end{align}

\noindent By definition, we have

\begin{align}\label{eq 5.7.9}
\displaystyle J_{M(y)}(y,\widetilde{f}) & =D^G(y)^{1/2}\int_{G_y(F)\backslash G(F)} \widetilde{f}(g^{-1}yg)v_{M(y)}(g) dg \\
\nonumber & =\int_{Z_G(x)(F)\backslash G(F)} \sum_{z\in G_x(F)\backslash Z_G(x)(F)} D^{G_x}(y)^{1/2} \\
\nonumber & \int_{G_y(F)\backslash G_x(F)} \left({}^g \widetilde{f}\right)_{x,\Omega_x}(z^{-1}g_x^{-1}yg_xz)v_{M(y)}(g_xzg)dg_xdg
\end{align}

\noindent Let $g\in G(F)$. By (ii), up to translating $g$ by an element of $Z_G(x)(F)$, we have $\left({}^g\widetilde{f}\right)_{x,\Omega_x}=\alpha(g)f$. Hence, the inner term of the last expression above becomes

\begin{align}\label{eq 5.7.10}
\displaystyle \alpha(g) \sum_{z\in G_x(F)\backslash Z_G(x)(F)} D^{G_x}(y)^{1/2} \int_{G_y(F)\backslash G_x(F)} {}^z f(g_x^{-1}yg_x) v_{M(y)}(g_xzg)dg_x
\end{align}

\noindent By Lemma 3.3 of \cite{Wa1}, we have the descent formula

$$\displaystyle v_{M(y)}(g_x\gamma)=\sum_{L\in \mathcal{L}(M(y)_x)}\sum_{Q\in \mathcal{P}(\boldsymbol{L})} v_{M(y)_x}^{Q_x}(g_x) u_Q\left(H_Q(g_x\gamma)-H_{Q_x}(g_x)\right)$$

\noindent for all $g_x\in G_x(F)$ and all $\gamma\in G(F)$, where $\mathcal{L}(M(y)_x)$ denotes the set of Levi subgroups of $G_x$ containing $M(y)_x$ and for $L\in \mathcal{L}(M(y)_x)$, $\boldsymbol{L}$ denotes the centralizer of $A_L$ in $G$ (a Levi subgroup of $G$), the other terms appearing in the formula above have been defined in Section \ref{section 1.10} (they depend on the choice of two maximal compact subgroup $K$ and $K_x$ of $G(F)$ and $G_x(F)$ which are special in the $p$-adic case). We may thus decompose the expression \ref{eq 5.7.10} further as a sum over $L\in \mathcal{L}(M(y)_x)$, $Q\in \mathcal{P}(\boldsymbol{L})$ and $z\in G_x(F)\backslash Z_G(x)(F)$ of

$$\displaystyle \alpha(g) D^{G_x}(y)^{1/2}\int_{G_y(F)\backslash G_x(F)} {}^z f(g_x^{-1}yg_x) v_{M(y)_x}^{Q_x}(g_x) u_Q\left(H_Q(g_x zg)-H_{Q_x}(g_x)\right)dg_x$$

\noindent Since the function $g_x\mapsto v_{M(y)_x}^{Q_x}(g_x) u_Q\left(H_Q(g_x zg)-H_{Q_x}(g_x)\right)$ is invariant by left translation by $U_{Q_x}(F)$ (the unipotent radical of $Q_x(F)$), $y\in Q_x(F)$ and ${}^z f$ is strongly cuspidal, this last term is zero unless $Q_x=G_x$. As $x$ is elliptic in $G(F)$, this last condition is equivalent to $Q=G$ and $L=G_x$. In this case the expression above reduces to $\alpha(g) J_{M(y)_x}(y,{}^z f)$, hence \ref{eq 5.7.10} is equal to

$$\displaystyle \alpha(g)\sum_{z\in G_x(F)\backslash Z_G(x)(F)} J_{M(y)_x}(y,{}^z f)$$

\noindent Going back to \ref{eq 5.7.9} and recalling that we choose $\alpha$ so that $\int_{Z_G(x)(F)\backslash G(F)}\alpha(g)=1$, we immediately get \ref{eq 5.7.8}. This ends the proof of \ref{eq 5.7.7} and of the proposition. $\blacksquare$

\begin{cor}\label{corollary 5.7.1}
Assume that $G$ admits an elliptic maximal torus and that the center of $G(F)$ is compact. Then

\begin{enumerate}[(i)]
\item There exists $\Omega\subseteq G(F)$ a completely $G(F)$-invariant open subset which is relatively compact modulo conjugation and contains $G(F)_{\elli}$ such that the linear map

$$f\in \mathcal{S}_{\scusp}(\Omega)\mapsto \theta_f\in QC_c(\Omega)$$

\noindent has dense image (in particular in the $p$-adic case this map is surjective).

\item If $F$ is $p$-adic, then for all $\theta\in QC(G(F))$, there exists a compact subset $\Omega_\theta\subseteq \mathcal{X}_{\elli}(G)$ such that

$$\displaystyle \int_{\Gamma_{\elli}(G)}D^G(x)\theta(x)\theta_\pi(x)dx=0$$

\noindent for all $\pi\in \mathcal{X}_{\elli}(G)-\Omega_\theta$, the integral above being absolutely convergent.

\item If $F=\mathbb{R}$, then for all $k\geqslant 0$ there exists a continuous semi-norm $\nu_k$ on $QC(G(F))$ such that

$$\displaystyle \left\lvert \int_{\Gamma_{\elli}(G)} D^G(x)\theta(x)\theta_\pi(x)dx\right\rvert\leqslant \nu_k(\theta)N^G(\pi)^{-k}$$

\noindent for all $\pi\in \mathcal{X}_{\elli}(G)$ and all $\theta\in QC(G(F))$, the integral above being absolutely convergent.

\item For all $\pi\in \mathcal{X}_{\elli}(G)$ there exists $f\in \mathcal{C}_{\scusp}(G(F))$ such that for all $\pi'\in \mathcal{X}_{\elli}(G)$ we have

$$\widehat{\theta}_f(\pi')\neq 0\Leftrightarrow \pi'=\pi$$

\end{enumerate}
\end{cor}

\vspace{2mm}

\noindent\ul{Proof}:

\begin{enumerate}[(i)]
\item For all $x\in G(F)_{\elli}$, choose $\Omega_x\subseteq G_x(F)$ a $G$-good open neighborhood of $x$ which is relatively compact modulo conjugation and such that there exists $\omega_x\subseteq \mathfrak{g}_x(F)$ a $G_x$-excellent open subset with $\Omega_x=x\exp(\omega_x)$. Since $G(F)_{\elli}$ is compact modulo conjugation, we may find $x_1,\ldots,x_k\in G(F)_{\elli}$ such that the family $(\Omega_{x_i}^G)_{1\leqslant i\leqslant k}$ covers $G(F)_{\elli}$. Let us set

$$\Omega= \Omega_{x_1}^G\cup\ldots\cup \Omega_{x_k}^G$$

\noindent By the existence of smooth invariant partition of unity (Proposition \ref{proposition 3.1.1}(ii)), we see that the natural continuous linear map

$$\displaystyle \bigoplus_{i=1}^k QC_c(\Omega_{x_i}^G)\to QC_c(\Omega)$$

\noindent is surjective. Hence, it is sufficient to show that for all $x\in G(F)_{\elli}$, the linear map

$$f\in \mathcal{S}_{\scusp}(\Omega_{x}^G)\mapsto \theta_f\in QC_c(\Omega_{x}^G)$$

\noindent has dense image. Let $x\in G(F)_{\elli}$. By Proposition \ref{proposition 4.4.1}(iii) and the previous proposition, it suffices to prove that the linear map

$$f\in \mathcal{S}_{\scusp}(\Omega_x)\mapsto \theta_f\in QC_c(\Omega_x)$$

\noindent has dense image. We have the following commutative diagram

$$\xymatrix{
\mathcal{S}_{\scusp}(\Omega_x) \ar[d] \ar[r] &  QC_c(\Omega_x) \ar[d] \\
\mathcal{S}_{\scusp}(\omega_x) \ar[r] & QC_c(\omega_x)}$$

\noindent where the two vertical arrows are given by $f\mapsto (R(x)f)_{\omega_x}$ and $\theta\mapsto (R(x)\theta)_{\omega_x}$ respectively and the two horizontal arrows are both given by $f\mapsto \theta_f$. By Proposition \ref{proposition 5.6.1}(i), the bottom map has dense image.On the other hand, by Lemma \ref{lemma 3.3.1} and Proposition \ref{proposition 4.4.1}(i) the two vertical maps are topological isomorphisms. Hence the top map also has dense image. This ends the proof of (i).

\item The integral is absolutely convergent since for every quasi-character $\theta$ the function $(D^G)^{1/2}\theta$ is locally bounded. We may of course assume that $\theta\in QC_c(\Omega)$, where $\Omega$ is as in (i) so that we may find $f\in \mathcal{S}_{\scusp}(\Omega)$ such that $\theta=\theta_f$. But then, by the Weil integration formula, we have

$$\displaystyle \int_{\Gamma_{\elli}(G)}D^G(x)\theta(x)\theta_\pi(x)dx=\theta_\pi(f)$$

\noindent for all $\pi\in \mathcal{X}_{\elli}(G)$ and the result follows from \ref{eq 2.6.1}.

\item The integral is absolutely convergent for the same reason as before. Let $\Omega$ be as in (i) and let $\varphi\in C^\infty(G(F))^G$ be an invariant function that is supported in $\Omega$ and equals $1$ in some neighborhood of $G(F)_{\elli}$. By Proposition \ref{proposition 4.4.1}(iv) and the closed graph theorem, the linear map

$$\theta\in QC(G(F))\mapsto \varphi\theta\in QC_c(\Omega)$$

\noindent is continuous. Hence, we only need to prove the estimate for $\theta\in QC_c(\Omega)$. Since for all $\theta\in QC_c(\Omega)$ the function $(D^G)^{1/2}\theta$ is locally bounded by a continuous semi-norm, for all $\pi\in \underline{\mathcal{X}}_{\elli}(G)$, the linear map

$$\displaystyle \theta\in QC_c(\Omega)\mapsto \int_{\Gamma_{\elli}(G)} D^G(x)^{1/2} \theta(x)\theta_\pi(x)dx$$

\noindent is continuous. Hence, by (i), we only need to prove the estimates for $\theta=\theta_f$ with $f\in \mathcal{S}_{\scusp}(\Omega)$. By the Weyl integration formula and Lemma \ref{lemma 5.2.2}(ii), for all $z\in \mathcal{Z}(\mathfrak{g})$ and all $f\in \mathcal{S}_{\scusp}(\Omega)$, we have

\[\begin{aligned}
\displaystyle \chi_\pi(z)\int_{\Gamma_{\elli}(G)} D^G(x) \theta_f(x)\theta_\pi(x)dx & =\int_{\Gamma_{\elli}(G)} D^G(x) \theta_f(x)(z\theta)_\pi(x)dx \\
 & =\int_{\Gamma_{\elli}(G)} D^G(x) (z^*\theta_f)(x)\theta_\pi(x)dx
\end{aligned}\]

\noindent There exists $z_1,\ldots, z_k\in \mathcal{Z}(\mathfrak{g})$ such that

$$\displaystyle \lvert \chi_\pi(z_1)\rvert+\ldots+\lvert \chi_\pi(z_k)\rvert\geqslant N^G(\pi)$$

\noindent for all $\pi\in \mathcal{X}_{\elli}(G)$ so that the estimates now follows from Proposition \ref{proposition 4.8.1}(ii).

\item Let $\pi\in \mathcal{X}_{\elli}(G)$ and identify it with one of its preimage in $\underline{\mathcal{X}}_{\elli}(G)$. By (i) there exists a sequence $(f_n)_{n\geqslant 1}$ of functions in $\mathcal{C}_{\scusp}(G(F))$ such that the sequence of functions $\left((D^G)^{1/2}\theta_{f_n}\right)_{n\geqslant 1}$ converge uniformly to $(D^G)^{1/2}\overline{\theta}_\pi$ on $G_{\reg}(F)_{\elli}$. By Proposition \ref{proposition 5.6.1}(ii), we have

\begin{align}\label{eq 5.7.11}
\displaystyle \theta_{f_n}=\sum_{\pi'\in \mathcal{X}_{\elli}(G)} D(\pi')\widehat{\theta}_{f_n}(\pi') \overline{\theta_{\pi'}}+\int_{\mathcal{X}_{ind}(G)} D(\pi')\widehat{\theta}_{f_n}(\pi') \overline{\theta_{\pi'}}d\pi'
\end{align}

\noindent for all $n\geqslant 1$. By the orthogonality relations of Arthur (Corollary 6.2 of \cite{A4}), we have that

$$\displaystyle \int_{\Gamma_{\elli}(G)}D^G(x)\theta_{\pi_0}(x)\overline{\theta_{\pi_1}}(x)dx=\left\{
    \begin{array}{ll}
        D(\pi_0)^{-1} & \mbox{ if } \pi_0=\pi_1 \\
        0 & \mbox{ otherwise.}
    \end{array}
\right.
$$

\noindent for all $\pi_0,\pi_1\in \mathcal{X}_{\elli}(G)$. Moreover, all $\pi'\in \mathcal{X}_{ind}(G)$ is properly induced and so has a character that vanish identically on $G_{\reg}(F)_{\elli}$. Hence, we deduce from \ref{eq 5.7.11} that

\begin{align}\label{eq 5.7.12}
\lim\limits_{n\to \infty} \widehat{\theta}_{f_n}(\pi')=\left\{
    \begin{array}{ll}
        D(\pi)^{-1} & \mbox{ if } \pi'=\pi \\
        0 & \mbox{ otherwise.}
    \end{array}
\right.
\end{align}

\noindent Consider the map

\begin{align}\label{eq 5.7.13}
\mathcal{X}_{\elli}(G)\to \mathcal{X}_{\tempe}(G)
\end{align}
$$\pi'\mapsto \Pi'$$

\noindent where for all $\pi'\in \mathcal{X}_{\elli}(G)$, $\Pi'$ denotes the unique element of $\mathcal{X}_{\tempe}(G)$ such that $\pi'$ is built up from the irreducible subrepresentations of $\Pi'$. Let $\Pi\in \mathcal{X}_{\tempe}(G)$ be the image of $\pi$. Then there exists a compact neighborhood $\mathcal{U}$ of $\Pi$ in $\mathcal{X}_{\tempe}(G)$ such that the inverse image of $\mathcal{U}$ by \ref{eq 5.7.13} is a finite set $\{\pi_0=\pi,\pi_1,\ldots,\pi_n\}\subseteq \mathcal{X}_{\elli}(G)$. By \ref{eq 5.7.12} we can certainly find a function $f'\in \mathcal{C}_{\scusp}(G(F))$ that is a finite linear combination of the $f_n$, $n\geqslant 1$, such that

$$\widehat{\theta}_{f'}(\pi_i)\neq 0\Leftrightarrow \pi_i=\pi$$

\noindent for all $i=0,\ldots,n$. Let $\varphi\in C_c^\infty(\mathcal{X}_{\tempe}(G))$ be such that $\varphi(\Pi)=1$ and $\Supp(\varphi)\subseteq \mathcal{U}$. By Lemma \ref{lemma 5.3.1}(i), there exists a unique function $f\in \mathcal{C}_{\scusp}(G(F))$ such that

$$\Pi'(f)=\varphi(\Pi')\Pi'(f')$$

\noindent for all $\Pi'\in \mathcal{X}_{\tempe}(G)$. This function obviously has the desired property. $\blacksquare$
\end{enumerate}

\section{The Gan-Gross-Prasad triples}\label{section 6}

In this chapter, we introduce the main characters of this paper that we call {\em GGP triples}. These are certain triples $(G,H,\xi)$ where $G$ is a connected reductive group, $H$ a closed subgroup and $\xi$ a unitary character of $H(F)$. These GGP triples are themselves associated to pairs $(W,V)$ of hermitian spaces with $W\subset V$ and $W^\perp$ satisfying a certain condition. We call such pairs {\em admissible}. The precise definitions of admissible pairs, GGP triples and related objects are given in Section \ref{section 6.2}. Section \ref{section 6.1} contains some purely group-theoretic background on unitary groups. Given a GGP triple $(G,H,\xi)$, we associate in Section \ref{section 6.3}, following Gan, Gross and Prasad, a multiplicity function $\pi\mapsto m(\pi)$ on the set of admissible irreducible representations $\pi$ of $G(F)$. Fundamental results of Aizenbud-Gourevitch-Rallis-Schiffmann and Jiang-Sun-Zhu, that we recall, state that this multiplicity is always less or equal to $1$. In Section \ref{section 6.4}, we show that $H\backslash G$ is a $F$-spherical variety which means that there exists a minimal parabolic subgroup $P_{\mini}$ with $HP_{\mini}$ open in $G$. More generally, we call a parabolic subgroup $P$ {\em good} if $HP$ is open in $G$ and we list some properties of these subgroups. Sections \ref{section 6.5}, \ref{section 6.7} and \ref{section 6.8} are devoted to the proof of certain estimates that will be needed in Chapters \ref{section 8}, \ref{section 7} and \ref{section 9}. Also crucial for these subsequent chapters, as well as for the proof of the estimates, is the existence of a certain `weak Cartan decomposition' for the homogeneous variety $H\backslash G$ which is the subject of Section \ref{section 6.6}. Such a decomposition is known in some generality and in particular for split $p$-adic spherical varieties (\cite{SV} Lemma 5.3.1), symmetric spaces (\cite{BO}, \cite{DS}) and real spherical varieties \cite{KKSS}. Strictly speaking, these references do not cover the case at hand when the field $F$ is $p$-adic but, fortunately, in this particular situation the author has already established the existence of a weak Cartan decomposition `by hand' in a previous paper \cite{Beu1}.

\subsection{Hermitian spaces and unitary groups}\label{section 6.1}

\noindent We will henceforth fix a quadratic extension $\gls{E}$ of $F$. We will denote by $x\mapsto \gls{xbar}$ the nontrivial $F$-automorphism of $E$ and by $\gls{NEF}$ and $\gls{TrEF}$ the norm and the trace of this extension respectively. We also fix a nonzero element $\gls{eta}\in E$ with zero trace. We will set $\gls{Ebar}=E\otimes_F \overline{F}$. This is an \'etale 2-extension of $\overline{F}$ and as such is isomorphic to $\overline{F}\times \overline{F}$ but we won't fix such an isomorphism.

\vspace{2mm}

\noindent By an hermitian space we will mean a finite dimensional $E$-vector space equipped with a non-degenerate hermitian form $h$ which is linear in the second argument. For $V$ an hermitian space, we will denote by $\gls{U(V)}$ the corresponding unitary group and by $\gls{u(V)}$ the Lie algebra of $U(V)$. Set $V_{\overline{F}}=V\otimes_F \overline{F}$. Then $h$ has a natural extension, still denoted by $h$, to an $\overline{E}$-sesquilinear form on $V_{\overline{F}}$. Then $U(V)$ is the group of $\overline{E}$-linear automorphisms of $V_{\overline{F}}$ preserving the form $h$. We will identify $\mathfrak{u}(V)$ to the subspace of antihermitian, with respect to $h$, elements in $\End_{\overline{E}}(V_{\overline{F}})$. For $v,v'\in V_{\overline{F}}$, we will denote by $\gls{c(v,v')}$ the element of $\mathfrak{u}(V)$ defined by

$$c(v,v')(v'')=h(v,v'')v'-h(v',v'')v$$

\noindent The set $\{c(v,v');\; v,v'\in V\}$ generates $\mathfrak{u}(V)(F)$ as a $F$-vector space.

\vspace{2mm}

\noindent Every parabolic subgroup $P$ of $U(V)$ is the stabilizer of a flag of totally isotropic subspaces

$$Z_1\subset Z_2\subset\ldots\subset Z_k$$

\noindent If $M$ is a Levi component of $P$, then there exists totally isotropic subspaces $Z'_{\pm i}$, $1\leqslant i\leqslant k$, such that $Z_i=Z_{i-1}\oplus Z'_i$ for $1\leqslant i\leqslant k$ and $Z_i^\perp=Z_{i+1}^\perp\oplus Z'_{-i-1}$ for $0\leqslant i\leqslant k-1$, where we have set $Z_0=0$, such that $M$ is the stabilizer in $U(V)$ of the subspaces $Z'_{\pm i}$, $1\leqslant i\leqslant k$. Then, for all $1\leqslant i\leqslant k$, the form $h$ induces a perfect conjugate-duality between $Z'_i$ and $Z'_{-i}$. If $\widetilde{V}$ denotes the orthogonal complement of $\displaystyle \bigoplus_{i=1}^k (Z'_i\oplus Z'_{-i})$, we have a natural isomorphism

$$\displaystyle M\simeq GL_E(Z'_1)\times\ldots\times GL_E(Z'_k)\times U(\widetilde{V})$$

\noindent where for $1\leqslant i\leqslant k$, $GL_E(Z'_i)$ denotes the restriction of scalars from $E$ to $F$ of the general linear group of $Z'_i$.

\vspace{2mm}

\noindent We are now going to describe the regular nilpotent orbits in $\mathfrak{u}(V)(F)$. If $U(V)$ is not quasi-split then there are no such orbit. Assume that $U(V)$ is quasi-split. If $\dim(V)$ is odd or zero, then there is only one regular nilpotent orbits. Assume moreover that $\dim(V)>0$ is even. Then, there are exactly two regular nilpotent orbits. Since $U(V)$ is quasi-split, there exists a basis $(z_i)_{i=\pm 1,\ldots,\pm k}$ such that $h(z_i,z_j)=\delta_{i,-j}$ for all $i,j\in\{\pm 1,\ldots,\pm k\}$ (where $\delta_{i,-j}$ denotes the Kronecker symbol). Let $B$ be the stabilizer in $U(V)$ of the flag

$$\displaystyle \langle z_k\rangle\subset \langle z_k,z_{k-1}\rangle\subset\ldots\subset \langle z_k,\ldots,z_1\rangle$$

\noindent Then, $B$ is a Borel subgroup of $U(V)$. Denote by $U$ its unipotent radical. For all $\mu\in E$ with trace zero, define an element $X(\mu)\in\mathfrak{u}(F)$ by the assignments

$$\displaystyle X(\mu)z_k=0,\; X(\mu)z_i=z_{i+1} \mbox{ for } 1\leqslant i\leqslant k-1,\; X(\mu)z_{-1}=\mu z_1,\; X(\mu)z_{-i}=-z_{1-i} \mbox{ for } 2\leqslant i\leqslant k$$

\noindent Then, for all $\mu\in E^\times$ with $\Tra_{E/F}(\mu)=0$, $X(\mu)$ is regular nilpotent. Moreover the orbits of $X(\mu)$ and $X(\mu')$ coincide if and only if $\mu N_{E/F}(E^\times)=\mu' N_{E/F}(E^\times)$. It follows that for all $\lambda\in F^\times\setminus N_{E/F}(E^\times)$, the elements $X(\eta)$ and $X(\lambda\eta)$ are representatives of the two regular nilpotent conjugacy classes. Notice that in particular multiplication by any element of $F^\times \setminus N_{E/F}(E^\times)$ permutes the two regular nilpotent orbits in $\mathfrak{u}(V)(F)$.

\subsection{Definition of GGP triples}\label{section 6.2}

\noindent Let $\gls{(W,V)}$ be a pair of hermitian spaces. We will call $(W,V)$ an {\em admissible pair} if there exists an hermitian space $\gls{Z}$ satisfying

\vspace{2mm}

\begin{itemize}
\renewcommand{\labelitemi}{$\bullet$}

\item $V\simeq W\oplus^\perp Z$;

\item $Z$ is odd dimensional and $U(Z)$ is quasi-split.
\end{itemize}

\vspace{2mm}

\noindent This last condition admits the following more explicit translation: it means that there exist $\gls{nu}\in F^\times$ and a basis $(z_{-r},\ldots,z_{-1},z_0,z_1,\ldots,z_r)$ of $Z$ such that

\begin{align}\label{eq 6.2.1}
h(z_i,z_j)=\nu\delta_{i,-j}
\end{align}

\noindent for all $i,j\in\{0,\pm 1,\ldots,\pm r\}$.

\vspace{2mm}

\noindent Let $(W,V)$ be an admissible pair. Set $\gls{G}=U(W)\times U(V)$. We are going to associate to $(W,V)$ a triple $(G,H,\xi)$ where $H$ is an algebraic subgroup of $G$ and $\xi:H(F)\to\mathbb{C}^\times$ is a continuous character of $H(F)$ and this triple will be unique up to $G(F)$-conjugacy. Fix an embedding $W\subseteq V$ and set $Z=W^\perp$. We also fix $\nu\in F^\times$ and a basis $\gls{Zbasis}$ (where $\dim(Z)=2\gls{r}+1$) of $Z$ satisfying \ref{eq 6.2.1}. Denote by $\gls{PV}$ the stabilizer in $U(V)$ of the following flag of totally isotropic subspaces of $V$

$$\langle z_r\rangle\subset \ldots\subset \langle z_r,\ldots,z_1\rangle$$

\noindent Then, $P_V$ is a parabolic subgroup of $U(V)$. We will denote by $\gls{N}$ its unipotent radical. Let $\gls{MV}$ the stabilizer in $P_V$ of the lines $\langle z_i\rangle$ for $i=\pm 1,\ldots,\pm r$. It is a Levi component of $P_V$. Set $\gls{P}=U(W)\times P_V$. Then, $P$ is a parabolic subgroup of $G$ with unipotent radical $N$ and $\gls{M}=U(W)\times M_V$ is a Levi component of it. Identifying $U(W)$ with its image by the diagonal embedding $U(W)\hookrightarrow G$, we have $U(W)\subseteq M$. In particular, conjugation by $U(W)$ preserves $N$ and we set

$$\displaystyle \gls{H}=U(W)\ltimes N$$

\noindent It only remains to define the character $\xi$. Let us define a morphism $\gls{lambda}\colon N\to \mathbb{G}_a$ by

$$\displaystyle \lambda(n)=\Tra_{E/F}\left(\sum_{i=0}^{r-1} h(z_{-i-1},nz_i)\right),\;\;\; n\in N$$

\noindent It is easy to check that $\lambda$ is $U(W)$-invariant, hence it admits a unique extension, still denoted by $\lambda$, to a morphism $H\to \mathbb{G}_a$ which is trivial on $U(W)$. We denote by $\gls{lambdaF}:H(F)\to F$ the morphism induced on the groups of $F$-points. Recall that we have fixed a non-trivial continuous additive character $\psi$ of $F$. We set

$$\displaystyle \gls{xi}(h)=\psi(\lambda_F(h))$$

\noindent for all $h\in H(F)$. This ends the definition of the triple $(G,H,\xi)$. We easily check that this definition depends on the various choices only up to $G(F)$-conjugacy. We will call a triple obtained in this way (i.e., from an admissible pair $(W,V)$) a {\em Gan-Gross-Prasad triple} or {\em GGP triple} for short.

\vspace{2mm}

\noindent From now on and until the end of Chapter \ref{section 11}, we fix an admissible pair of hermitian spaces $(W,V)$. We also fix data and notation as above, that is: an embedding $W\subseteq V$, $Z=W^\perp$, an element $\nu\in F^\times$ and a basis $(z_i)_{i=0,\pm 1,\ldots,\pm r}$ of $Z$ satisfying \ref{eq 6.2.1}, the parabolic subgroup $P=MN$, the algebraic character $\lambda$ and the character $\xi$. We will denote by $(G,H,\xi)$ the GGP triple constructed as above. We will also use the following additional notation

\vspace{2mm}

\begin{itemize}

\renewcommand{\labelitemi}{$\bullet$}

\item $\gls{d}=dim(V)$ and $\gls{m}=dim(W)$

\item $\gls{Z+}=\langle z_r,\ldots,z_1\rangle$, and $\gls{Z-}=\langle z_{-1},\ldots,z_{-r}\rangle$;

\item $\gls{D}=Ez_0$ and $\gls{V_0}=W\oplus D$;

\item $\gls{H_0}=U(W)$ and $\gls{G_0}=U(W)\times U(V_0)$. We consider $H_0$ as a subgroup of $G_0$ via the diagonal embedding $H_0\hookrightarrow G_0$. The triple $(G_0,H_0,1)$ is the GGP triple associated to the admissible pair $(W,V_0)$;

\item $\gls{T}$ the subtorus of $U(V)$ preserving the lines $\langle z_i\rangle$, for $i=1,\ldots,r$ and $i=-1,\ldots,-r$ and acting trivially on $V_0$. We have $M=T\times G_0$;

\item $\gls{A}$ the split part of the torus $T$, it is also the split part of the center of $M$;

\item $\xi$ the character of $\mathfrak{h}(F)$, where $\mathfrak{h}=\Lie(H)$, which is trivial on $\mathfrak{u}(W)(F)$ and equal to $\xi\circ \exp$ on $\mathfrak{n}(F)$.

\item $\gls{B}$ is the following non-degenerate $G$-invariant bilinear form on $\mathfrak{g}$:

$$\displaystyle B\left((X_W,X_V),(X'_W,X'_V)\right)=\frac{1}{2}\left( \Tra_{E/F}\left(\Tr(X_WX'_W)\right)+ \Tra_{E/F}\left(\Tr(X_VX'_V)\right)\right)$$

\noindent We will use $B(.,.)$ to normalize the Haar measures on both $\mathfrak{g}(F)$ and $G(F)$ as explained in Section \ref{section 1.6}. We also fix Haar measures on all algebraic subgroups of $G(F)$ and their Lie algebras as explained in Section \ref{section 1.6}.

\end{itemize}

\vspace{2mm}

\noindent Note that when $r=0$ (that is when $Z=D$ is a line), we have $G=G_0$, $H=H_0$ and $\xi=1$. If this is the case, we will say that we are in {\em the codimension one case}. 

\vspace{2mm}

\noindent We will need the following (cf.\ Section \ref{section 1.2} for the definition of {\em norm descent property}):

\vspace{2mm}

\begin{lem}\label{lemma 6.2.1}
\begin{enumerate}[(i)]
\item The map $G\to H\backslash G$ has the norm descent property.

\item The orbit under $M$-conjugacy of $\lambda$ in $\left(\mathfrak{n}/[\mathfrak{n},\mathfrak{n}]\right)^*$ is a Zariski open subset.
\end{enumerate}
\end{lem}

\vspace{2mm}

\noindent\ul{Proof}:

\begin{enumerate}[(i)]
\item  We have a natural identification $H\backslash G=N\backslash U(V)$, so that it is sufficient to prove that $U(V)\to N\backslash U(V)$ has the norm descent property. Since this map is $U(V)$-equivariant for the obvious transitive actions, we only need to show that it admits a section over a nonempty Zariski-open subset. If we denote by $\overline{P}_V=M_V\overline{N}$ the parabolic subgroup opposite to $P_V$ with respect to $M_V$, the multiplication map $N\times M_V\times \overline{N}\to U(V)$ is an open immersion. The image of that open subset is open in $N\backslash U(V)$ and the restriction of the projection $U(V)\to N\backslash U(V)$ to that open set is $N\times M_V\times \overline{N}\to M_V\times \overline{N}$. This map obviously has a section.

\item If $r=0$, i.e., if we are in the codimension one case, we have $\mathfrak{n}=0$ and the result is trivial. Assume now that $r\geqslant 1$. It suffices to show that the dimension of the orbit $M\cdot\lambda$ is equal to the dimension of $\mathfrak{n}/[\mathfrak{n},\mathfrak{n}]$. We easily compute

$$\displaystyle \dim\left(\mathfrak{n}/[\mathfrak{n},\mathfrak{n}]\right)=2(m+r)$$

\noindent and

$$\displaystyle \dim(M)=m^2+2r+(m+1)^2$$

\noindent The stabilizer $M_\lambda$ of $\lambda$ is easily seen to be $M_\lambda=Z(G)\left( U(W)\times U(W)\right)$ (where $Z(G)$ denotes the center of $G$). Hence, we have

$$\displaystyle \dim(M_\lambda)=1+2m^2$$

\noindent and the dimension of the orbit $M\cdot \lambda$ is

$$\displaystyle \dim(M.\lambda)=\dim(M)-\dim(M_\lambda)=m^2+2r+(m+1)^2-1-2m^2=2(r+m)$$

\noindent which is the same as $\dim\left(\mathfrak{n}/[\mathfrak{n},\mathfrak{n}]\right)$. $\blacksquare$
\end{enumerate}

\subsection{The multiplicity $m(\pi)$}\label{section 6.3}

\noindent For $\pi\in \Temp(G)$, let us denote by $\gls{HomHpixi}$ the space of all continuous linear forms $\ell:\pi^\infty\to \mathbb{C}$ such that

$$\ell(\pi(h)e)=\xi(h)\ell(e)$$

\noindent for all $e\in \pi^\infty$ and for all $h\in H(F)$. We define the {\em multiplicity} $\gls{mpi}$ to be the dimension of that space of linear forms, that is

$$m(\pi)=\dim\; \Hom_H(\pi,\xi),\;\;\; \pi\in \Temp(G)$$

\noindent We have the following multiplicity one result which is Theorem A of \cite{JSZ} in the Archimedean case (for $r=0$, it is Theorem B of \cite{SZ}) and follows from the combination of Theorem (1') of \cite{AGRS} (which treat the case $r=0$) and Theorem 15.1 of \cite{GGP} (showing how to extend the result to general $r$) in the $p$-adic case.

\vspace{3mm}

\begin{theo}\label{theorem 6.3.1}
We have

$$m(\pi)\leqslant 1$$

\noindent for all $\pi\in \Temp(G)$.
\end{theo}

\vspace{3mm}

\noindent Note that we have

\begin{align}\label{eq 6.3.1}
m(\pi)=m(\overline{\pi})
\end{align}

\noindent for all $\pi\in \Temp(G)$. Indeed, the conjugation map $\ell\mapsto \overline{\ell}$ induces an isomorphism

$$\Hom_H(\pi,\xi)\simeq \Hom(\overline{\pi},\overline{\xi})$$

\noindent and as we easily check, there exists an element $a\in A(F)$ such that $\xi(aha^{-1})=\overline{\xi(h)}$ for all $h\in H(F)$, hence the linear map $\ell\mapsto \ell\circ\pi(a)$ induces an isomorphism

$$\Hom(\overline{\pi},\overline{\xi})\simeq \Hom_H(\overline{\pi},\xi)$$

\noindent and \ref{eq 6.3.1} follows.

\subsection{$H\backslash G$ is a spherical variety, good parabolic subgroups}\label{section 6.4}

\noindent We will say that a parabolic subgroup $\overline{Q}$ of $G$ is {\em good} if $H\overline{Q}$ is Zariski-open in $G$. This condition is equivalent to $H(F)\overline{Q}(F)$ being open in $G(F)$.

\begin{prop}\label{proposition 6.4.1}
\begin{enumerate}[(i)]
\item There exist minimal parabolic subgroups of $G$ that are good and they are all conjugate under $H(F)$. Moreover, if $\overline{P}_{\mini}=M_{\mini}\overline{U}_{\mini}$ is a good minimal parabolic subgroup we have $H\cap \overline{U}_{\mini}=\{1\}$ and the complement of $H(F)\overline{P}_{\mini}(F)$ in $G(F)$ has null measure;

\item A parabolic subgroup $\overline{Q}$ of $G$ is good if and only if it contains a good minimal parabolic subgroup;

\item Let $\overline{P}_{\mini}=M_{\mini}\overline{U}_{\mini}$ be a good minimal parabolic subgroup and let $A_{\mini}=A_{M_{\mini}}$ be the maximal split central subtorus of $M_{\mini}$. Set

$$A_{\mini}^+=\{a\in A_{\mini}(F); \lvert \alpha(a)\rvert\geqslant 1 \; \forall \alpha\in R(A_{\mini},\overline{P}_{\mini})\}$$

\noindent Then, we have inequalities

\begin{num}

\item\label{eq 6.4.a} $\sigma(h)+\sigma(a)\ll \sigma(ha)$ for all $a\in A_{\mini}^+$, $h\in H(F)$.

\item\label{eq 6.4.b} $\sigma(h)\ll\sigma(a^{-1}ha)$ for all $a\in A_{\mini}^+$, $h\in H(F)$.
\end{num}

\end{enumerate}

\end{prop}

\vspace{3mm}

\noindent\ul{Proof}:

\begin{enumerate}[(i)]

\item Set $w_0=z_0$ and choose a family $(w_1,\ldots,w_\ell)$ of mutually orthogonal vectors in $W$ which is maximal subject to the condition

$$h(w_i)=(-1)^i\nu,\; i=1,\ldots,\ell$$

\noindent Let $\lceil{\frac{\ell}{2}}\rceil$ (resp.\ $\lfloor{\frac{\ell}{2}}\rfloor$) be the smallest (resp.\ the largest) integer which is not less (resp.\ not greater) than $\frac{\ell}{2}$. We define $u_i$, for $i=1,\ldots,\lceil{\frac{\ell}{2}}\rceil$ by

$$u_i=w_{2i-2}+w_{2i-1}$$

\noindent and $u'_i$, for $i=1,\ldots,\lfloor{\frac{\ell}{2}}\rfloor$, by

$$u'_i=w_{2i-1}+w_{2i}$$

\noindent Then, the subspaces

$$Z_{V_0}=\langle u_1,\ldots,u_{\lceil{\frac{\ell}{2}}\rceil}\rangle$$
$$Z_{W}=\langle u'_1,\ldots,u'_{\lfloor{\frac{\ell}{2}}\rfloor}\rangle$$

\noindent are maximal isotropic subspaces of $V_0$ and $W$ respectively. Let $\overline{P}_{V_0}$ and $\overline{P}_W$ be the stabilizers in $U(V_0)$ and $U(W)$ of the totally isotropic flags

$$\langle u_1\rangle \subseteq \langle u_1,u_2\rangle \subseteq\ldots\subseteq \langle u_1,\ldots,u_{\lceil{\frac{\ell}{2}}\rceil}\rangle$$

\noindent and

$$\langle u'_1\rangle\subseteq \langle u'_1,u'_2\rangle\subseteq\ldots\subseteq \langle u'_1,\ldots,u'_{\lfloor{\frac{\ell}{2}}\rfloor}\rangle$$

\noindent respectively. Then $\overline{P}_{V_0}$ and $\overline{P}_W$ are minimal parabolic subgroups of respectively $U(V_0)$ and $U(W)$. Set

$$\overline{P}_0=\overline{P}_{W}\times \overline{P}_{V_0}$$

\noindent It is a minimal parabolic subgroup of $G_0$. Let $W_{an}$ be the orthogonal complement in $W$ of $\langle w_1,\ldots,w_\ell\rangle$. We claim the following

\vspace{3mm}

\begin{num}
\item\label{eq 6.4.1} We have $H_0\cap\overline{P}_0=U(W_{an})$ and $H_0\overline{P}_0$ is Zariski-open in $G_0$ (i.e., $\overline{P}_0$ is a good parabolic subgroup of $G_0$).
\end{num}

\vspace{3mm}

\noindent The second claim follows from the first one by dimension consideration. We prove the first claim. Let $h_0\in H_0\cap \overline{P}_0$. Consider the action of $h_0$ on $V_0$. Since $h_0$ belongs to $H_0$, $h_0$ must stabilize $w_0=z_0$. On the other hand, since $h_0$ belongs to $\overline{P}_0$, $h_0$ must stabilize the line $\langle w_0+w_1\rangle$. Because $w_0$ is orthogonal to $w_1$, it follows that $h_0$ stabilizes $w_1$. We show similarly that $h_0$ stabilizes $w_2,\ldots,w_\ell$, hence $h_0\in U(W_{an})$. This ends the proof of \ref{eq 6.4.1}.

\vspace{2mm}

\noindent Let $\overline{P}=M\overline{N}$ be the parabolic subgroup opposite to $P$ with respect to $M$ and set

$$\overline{P}_{\mini}=\overline{P}_0T\overline{N}$$

\noindent it is a minimal parabolic subgroup of $G$. We deduce easily from \ref{eq 6.4.1} the following

\vspace{3mm}

\begin{num}
\item\label{eq 6.4.2} $\overline{P}_{\mini}$ is a good parabolic subgroup and we have $\overline{P}_{\mini}\cap H=U(W_{an})$.
\end{num}

\vspace{3mm}

\noindent This already proves that there exists minimal parabolic subgroup that are good. Let $\overline{P}_{\mini}'$ be another good minimal parabolic subgroup and let us show that $\overline{P}_{\mini}$ and $\overline{P}'_{\mini}$ are conjugate under $H(F)$. Let $g\in G(F)$ such that $\overline{P}_{\mini}'=g\overline{P}_{\mini}g^{-1}$. Set $\mathcal{U}=H\overline{P}_{\mini}$ and $\mathcal{Z}=G-\mathcal{U}$. Then, $\mathcal{Z}$ is a proper Zariski-closed subset of $G$ which is obviously $H\times \overline{P}_{\mini}$-invariant (for the left and right multiplication respectively). If $g\in \mathcal{Z}$, then we would have

$$H\overline{P}_{\mini}'=Hg\overline{P}_{\mini}g^{-1}\subseteq \mathcal{Z}g^{-1}$$

\noindent which is impossible since $\overline{P}_{\mini}'$ is a good parabolic subgroup. Hence, we have $g\in \mathcal{U}\cap G(F)=\mathcal{U}(F)$. If we can prove that $g\in H(F)\overline{P}_{\mini}(F)$, then we will be done. Hence, it suffices to show that

\begin{align}\label{eq 6.4.3}
\mathcal{U}(F)=H(F)\overline{P}_{\mini}(F)
\end{align}

\noindent by a standard argument, this follows from

\vspace{3mm}

\begin{num}
\item\label{eq 6.4.4} The map $H^1(F,H\cap \overline{P}_{\mini})\to H^1(F,H)$ is injective.
\end{num}

\vspace{3mm}

\noindent By \ref{eq 6.4.2}, we have $H^1(F,H\cap\overline{P}_{\mini})=H^1(F,U(W_{an}))$. Since $H=U(W)\ltimes N$ with $N$ unipotent, we also have $H^1(F,H)=H^1(F,U(W))$. The two sets $H^1(F,U(W_{an}))$ and $H^1(F,U(W))$ classify the (isomorphism classes of) hermitian spaces of the same dimension as $W_{an}$ and $W$ respectively. Moreover, the map $H^1(F,U(W_{an}))\to H^1(F,U(W))$ we are considering sends $W_{an}'$ to $W'_{an}\oplus W_{an}^\perp$, where $W_{an}^\perp$ denotes the orthogonal complement of $W_{an}$ in $W$. By Witt's theorem, this map is injective. This proves \ref{eq 6.4.4} and ends the proof that all good minimal parabolic subgroups are conjugate under $H(F)$.

\vspace{2mm}

\noindent It only remains to show the last part of (i) that is: $H\cap \overline{U}_{\mini}=\{1\}$ and the complement of $H(F)\overline{P}_{\mini}(F)$ in $G(F)$ has null measure for every good minimal parabolic subgroup $\overline{P}_{\mini}=M_{\mini}\overline{U}_{\mini}$. Since we already proved that all good minimal parabolic subgroups are $H(F)$-conjugate, we only need to consider one of them. Choosing for $\overline{P}_{\mini}$ the parabolic subgroup that we constructed above, the result follows directly from \ref{eq 6.4.2} and \ref{eq 6.4.3}.

\item Let $\overline{Q}$ be a good parabolic subgroup and choose $P_{\mini}\subseteq \overline{Q}$ a minimal parabolic subgroup. Set

$$\displaystyle \mathcal{G}:=\{g\in G;g^{-1}P_{\mini}g\mbox { is good}\}$$

\noindent It is a Zariski-open subset of $G$ since it is the inverse image of the Zariski-open subset $\{\mathcal{V}\in Gr_n(\mathfrak{g}); \mathcal{V}+\mathfrak{h}=\mathfrak{g}\}$ of the Grassmannian variety $Gr_n(\mathfrak{g})$, where $n=\dim(P_{\mini})$, by the regular map $g\in G\mapsto g^{-1}\mathfrak{p}_{\mini}g \in Gr_n(\mathfrak{g})$. Moreover, it is non-empty (since by (i) there exists good minimal parabolic subgroups). Since $\overline{Q}$ is good, the intersection $\overline{Q}H\cap \mathcal{G}$ is non-empty too. Hence, we may find $\overline{q}_0\in \overline{Q}$ such that $\overline{q}_0^{-1}P_{\mini}\overline{q}_0$ is a good parabolic subgroup. This parabolic subgroup is contained in $\overline{Q}$ but it may not be defined over $F$. Define

$$\mathcal{Q}:=\{\overline{q}\in \overline{Q}; \overline{q}^{-1}P_{\mini}\overline{q} \mbox{ is good}\}$$

\noindent Then again $\mathcal{Q}$ is a Zariski-open subset of $\overline{Q}$ and we just proved that it is non-empty. Since $\overline{Q}(F)$ is Zariski-dense in $\overline{Q}$, the set $\mathcal{Q}(F)$ is non-empty. Then, for all $\overline{q}\in \mathcal{Q}(F)$ the parabolic subgroup $\overline{q}^{-1}P_{\mini}\overline{q}$ has all the desired properties.

\item First we prove that both \ref{eq 6.4.a} and \ref{eq 6.4.b} don't depend on the particular pair $(\overline{P}_{\mini},M_{\mini})$ chosen. Let $(\overline{P}_{\mini}',M'_{\mini})$ be a similar pair, that is : $\overline{P}'_{\mini}$ is a good parabolic subgroup and $M'_{\mini}$ is a Levi component of it. Then, by (i), there exists $h\in H(F)$ such that $\overline{P}'_{\mini}=h\overline{P}_{\mini}h^{-1}$ and obviously the inequalities \ref{eq 6.4.a} and \ref{eq 6.4.b} are true for the pair $(\overline{P}_{\mini},M_{\mini})$ if and only if they are true for the pair $(h\overline{P}_{\mini}h^{-1},hM_{\mini}h^{-1})=(\overline{P}'_{\mini},hM_{\mini}h^{-1})$. Hence, we may assume without loss of generality that $\overline{P}_{\mini}=\overline{P'}_{\mini}$. Then, there exists $\overline{u}\in \overline{U}_{\mini}(F)$ such that $M'_{\mini}=\overline{u}M_{\mini}\overline{u}^{-1}$ and we have ${A'}_{\mini}^+=\overline{u}A_{\mini}^+\overline{u}^{-1}$. By definition of $A_{\mini}^+$, the sets $\{a^{-1}\overline{u}a;\; a\in A_{\mini}^+\}$ and $\{a^{-1}\overline{u}^{-1}a;\; a\in A_{\mini}^+\}$ are bounded. It follows that

$$\sigma\left(h\overline{u}a\overline{u}^{-1}\right)\sim \sigma(ha)$$

$$\sigma\left(\overline{u}a^{-1}\overline{u}^{-1}h\overline{u}a\overline{u}^{-1}\right)\sim \sigma(a^{-1}ha)$$

\noindent for all $a\in A_{\mini}^+$ and all $h\in H(F)$. We easily deduce that the inequalities \ref{eq 6.4.a} and \ref{eq 6.4.b} are satisfied for the pair $(\overline{P}_{\mini},M_{\mini})$ if and only if they are satisfied for the pair $(\overline{P}'_{\mini},M'_{\mini})$.

\vspace{2mm}

\noindent We now reduce the proof of \ref{eq 6.4.a} and \ref{eq 6.4.b} to the codimension one case. Let $\overline{P}_0=M_0\overline{U}_0$ be a good minimal parabolic subgroup of $G_0$. Let $A_0=A_{M_0}$ be the split part of the center of $M_0$ and let

$$\displaystyle A_0^+=\{a_0\in A_0(F);\; \lvert \alpha(a)\rvert \geqslant 1\; \forall \alpha\in R(A_0,\overline{P}_0)\}$$

\noindent Set $\overline{P}_{\mini}=\overline{P}_0T\overline{N}$ and $M_{\mini}=M_0T$. Then, $\overline{P}_{\mini}$ is a good parabolic subgroup of $G$, $M_{\mini}$ is a Levi component of it and $A_{\mini}^+\subseteq A(F)A_0^+$. We have

$$\sigma(nh_0aa_0)\gg\sigma(n)+\sigma(a)+\sigma(h_0a_0)$$

\noindent for all $h=nh_0\in H(F)=N(F)H_0(F)$ and all $(a,a_0)\in A(F)\times A_0^+$. Since, $\sigma(aa_0)\sim \sigma(a)+\sigma(a_0)$ and $\sigma(nh_0)\sim \sigma(n)+\sigma(h_0)$ for all $(a,a_0)\in A(F)\times A_0^+$ and all $(n,h_0)\in N(F)\times H_0(F)$, the inequality \ref{eq 6.4.a} will follow from

\vspace{3mm}

\begin{num}
\item\label{eq 6.4.5} $\sigma(h_0a_0)\gg\sigma(h_0)+\sigma(a_0)$, for all $a_0\in A_0^+$ and all $h_0\in H_0(F)$.
\end{num}

\vspace{3mm}

\noindent On the other hand, we have $\sigma(a^{-1}na)\gg \sigma(n)$ for all $a\in A_{\mini}^+$ and all $n\in N(F)$. Hence,

$$\sigma(a^{-1}nh_0a)\gg\sigma(n)+\sigma(a^{-1}h_0a)=\sigma(n)+\sigma(a_0^{-1}h_0a_0)$$

\noindent for all $(n,h_0)\in N(F)\times H_0(F)$ and all $a\in A_{\mini}^+$, where $a_0$ denote the unique element of $A_0^+$ such that $aa_0^{-1}\in A(F)$. Hence, the point \ref{eq 6.4.b} will follow from

\vspace{3mm}

\begin{num}
\item\label{eq 6.4.6} $\sigma(a_0^{-1}h_0a_0)\gg \sigma(h_0)$, for all $a_0\in A_0^+$ and all $h_0\in H_0(F)$.
\end{num}

\vspace{3mm}

\noindent Of course, to prove \ref{eq 6.4.5} and \ref{eq 6.4.6} we may work with any pair $(\overline{P}_0,M_0)$ that we want. Introduce a sequence $(w_0,\ldots,w_\ell)$ and a parabolic subgroup $\overline{P}_0=\overline{P}_W\times \overline{P}_{V_0}$ of $G_0$ as in (i). By \ref{eq 6.4.1}, $\overline{P}_0$ is a good parabolic subgroup of $G_0$. Let $M_{V_0}$ be the Levi component of $P_{V_0}$ that preserves the lines

$$\displaystyle \langle u_1\rangle,\ldots,\langle u_{\lceil \frac{\ell}{2} \rceil}\rangle \mbox{ and } \langle u_{-1}\rangle,\ldots, \langle u_{-\lceil \frac{\ell}{2}\rceil}\rangle$$

\noindent where we have set $u_{-i}=w_{2i-2}-w_{2i-1}$ for $i=1,\ldots,\lceil\frac{\ell}{2}\rceil$, and let $M_W$ be the Levi component of $\overline{P}_W$ that preserves the lines

$$\displaystyle \langle u'_1\rangle,\ldots,\langle u'_{\lfloor \frac{\ell}{2} \rfloor}\rangle \mbox{ and } \langle u'_{-1}\rangle,\ldots, \langle u'_{-\lfloor \frac{\ell}{2}\rfloor}\rangle$$

\noindent where we have set $u'_{-i}=w_{2i-1}-w_{2i}$ for $i=1,\ldots,\lfloor\frac{\ell}{2}\rfloor$. Set

$$M_0=M_W\times M_{V_0}$$

\noindent It is a Levi component of $\overline{P}_0$. We are going to prove \ref{eq 6.4.5} and \ref{eq 6.4.6} for the particular pair $(\overline{P}_0,M_0)$. We have a decomposition

$$A_0^+=A_W^+\times A_{V_0}^+$$

\noindent where $A_W^+$ and $A_{V_0}^+$ are defined in the obvious way. For all $a_{V_0}\in A_{V_0}^+$ (resp.\ $a_W\in A_W^+$) let us denote by $a_{V_0}^1,\ldots,a_{V_0}^{\lceil\frac{\ell}{2}\rceil}$ (resp.\ $a_W^1,\ldots,a_W^{\lfloor \frac{\ell}{2}\rfloor}$) the eigenvalues of $a_{V_0}$ (resp.\ $a_W$) acting on $u_1,\ldots,u_{\lceil\frac{\ell}{2}\rceil}$ (resp.\ on $u'_1,\ldots,u'_{\lfloor\frac{\ell}{2}\rfloor}$). Then, we have

$$\lvert a_{V_0}^1\rvert\geqslant\ldots\geqslant \lvert a_{V_0}^{\lceil \frac{\ell}{2}\rceil}\rvert\geqslant 1$$

$$\lvert a_W^1\rvert\geqslant\ldots\geqslant \lvert a_W^{\lfloor\frac{\ell}{2}\rfloor}\rvert\geqslant 1$$

\noindent for all $a_{V_0}\in A_{V_0}^+$ and all $a_W\in A_W^+$.

\vspace{2mm}

\noindent Of course we have

$$\sigma(h_0)+\sigma(a_0)\ll \sigma(h_0a_0)+\sigma(h_0)$$
$$\sigma(h_0)\ll \sigma_{U(V_0)}(h_0a_{V_0})+\sigma_{U(V_0)}(a_{V_0})$$
$$\sigma_{U(V_0)}(h_0a_{V_0})\ll \sigma(h_0a_0)$$

\noindent for all $a_0=(a_W,a_{V_0})\in A_0^+=A_W^+\times A_{V_0}^+$ and all $h_0\in H_0(F)$. Hence \ref{eq 6.4.5} will follow from

\vspace{3mm}

\begin{num}
\item\label{eq 6.4.7} $\sigma_{U(V_0)}(a_{V_0})\ll \sigma_{U(V_0)}(h_0a_{V_0})$, for all $a_{V_0}\in A_{V_0}^+$ and all $h_0\in H_0(F)$.
\end{num}

\vspace{3mm}

\noindent We have

\begin{align}\label{eq 6.4.8}
\sigma(a_{V_0})\sim \log\left(1+\lvert a_{V_0}^1\rvert\right)
\end{align}

\noindent for all $a_{V_0}\in A_{V_0}^+$. Moreover, for all $a_{V_0}\in A_{V_0}^+$ and all $h_0\in H_0(F)$ we have

$$h(h_0a_{V_0}u_1,w_0)=a_{V_0}^1h(h_0u_1,w_0)=a_{V_0}^1\left(h(w_0,w_0)+h(h_0w_1,w_0)\right)=a_{V_0}^1\nu$$

\noindent Since $\sigma_{U(V_0)}(g)\gg \log\left(1+\lvert h(gu_1,w_0)\rvert\right)$ for all $g\in U(W)$, \ref{eq 6.4.7} follows.

\vspace{2mm}

\noindent We now concentrate on the proof of \ref{eq 6.4.6}. Obviously, we only need to prove the following

\vspace{3mm}

\begin{num}
\item\label{eq 6.4.9} For all $v,v'\in V_0$, we have an inequality
$$\log\left(2+\lvert h(h_0v,v')\rvert\right)\ll \sigma(a_0^{-1}h_0a_0)$$
for all $a_0\in A_0^+$ and all $h_0\in H_0(F)$.
\end{num}

\vspace{3mm}

\noindent By sesquilinearity and since $\lvert h(h_0v,v')\rvert=\lvert h(h_0^{-1}v',v)\rvert$, it suffices to prove \ref{eq 6.4.9} in the following cases

\vspace{2mm}

\begin{itemize}
\item $v=w_i$ and $v'\in \langle w_i,\ldots,w_\ell\rangle\oplus W_{an}$ for $0\leqslant i \leqslant \ell$;
\item $v,v'\in W_{an}$
\end{itemize}

\vspace{2mm}

\noindent (recall that $W_{an}$ denote the orthogonal complement of $\langle w_0,\ldots,w_\ell\rangle$ in $V_0$). The proof of \ref{eq 6.4.9} in the second case is easy since we have $h(a_{V_0}^{-1}h_0a_{V_0}v,v')=h(h_0v,v')$ for all $a_{V_0}\in A_{V_0}(F)$, all $h_0\in H_0(F)$ and all $v,v'\in W_{an}$. Let us do the first case. The proof is by induction on $i$. For $i=0$, the result is obvious since $h_0w_0=w_0$ for all $h_0\in H_0(F)$. Let $1\leqslant i\leqslant \ell$ and assume that \ref{eq 6.4.9} is satisfied for $v=w_{i-1}$ and all $v'\in \langle w_{i-1},\ldots,w_\ell\rangle \oplus W_{an}$. If $i$ is odd, then the subspace $\langle w_{i-1},\ldots,w_\ell\rangle\oplus W_{an}$ is preserved by $A_{V_0}(F)$. Obviously, we only need to prove \ref{eq 6.4.9} for $v'$ an eigenvector for the action of $A_{V_0}(F)$ on that subspace. For all $a_{V_0}\in A_{V_0}^+$, the eigenvalue of $a_{V_0}$ on $v'$ have an absolute value which is greater or equal to $\lvert a_{V_0}^{(i-1)/2}\rvert^{-1}$. Hence, we have

\[\begin{aligned}
\sigma(a_0^{-1}h_0a_0) & \gg \log\left(2+\lvert h(a_{V_0}^{-1}h_0a_{V_0}u_{(i-1)/2},v')\rvert\right) \\
 & \gg \log\left(2+\lvert a_{V_0}^{(i-1)/2}\rvert \lvert h(h_0u_{(i-1)/2}, a_{V_0}v')\rvert\right) \\
 & \gg \log\left(2+\lvert h(h_0u_{(i-1)/2},v')\rvert\right)
\end{aligned}\]

\noindent for all $a_0=(a_W,a_{V_0})\in A_0^+=A_W^+\times A_{V_0}^+$ and all $h_0\in H_0(F)$. On the other hand we have $w_i=u_{(i-1)/2}-w_{i-1}$, so that

\[\begin{aligned}
\log\left(2+\lvert h(h_0w_i,v')\rvert\right)\ll \log\left(2+\lvert h(h_0u_{(i-1)/2},v')\rvert\right)+\log\left(2+\lvert h(h_0w_{i-1},v')\rvert\right)
\end{aligned}\]

\noindent for all $h_0\in H_0(F)$. Combining the two previous inequalities and the induction hypothesis we get the desired inequality. If $i$ is even, the proof is similar using the action on $W$ rather than on $V_0$. $\blacksquare$
\end{enumerate}

\subsection{Some estimates}\label{section 6.5}

\begin{lem}\label{lemma 6.5.1}
\begin{enumerate}[(i)]
\item There exists $\epsilon>0$ such that the integral

$$\displaystyle \int_{H_0(F)} \Xi^{G_0}(h_0)e^{\epsilon \sigma(h_0)} dh_0$$

\noindent is absolutely convergent.

\item There exists $d>0$ such that the integral

$$\displaystyle \int_{H(F)} \Xi^G(h)\sigma(h)^{-d} dh$$

\noindent is absolutely convergent.

\item For all $\delta>0$ there exists $\epsilon>0$ such that the integral

$$\displaystyle \int_{H(F)} \Xi^G(h)e^{\epsilon\sigma(h)} \left(1+\lvert \lambda(h)\rvert\right)^{-\delta}dh$$

\noindent is absolutely convergent (where $\lambda:H\to \mathbb{G}_a$ is the homomorphism defined in Section \ref{section 6.2}).

\end{enumerate}

\vspace{3mm}

\noindent Let $\overline{P}_{\mini}=M_{\mini}\overline{U}_{\mini}$ be a good minimal parabolic subgroup of $G$. We have the following

\vspace{3mm}

\begin{enumerate}[(i)]
\setcounter{enumi}{3}

\item For all $\delta>0$ there exists $\epsilon>0$ such that the integral

$$\displaystyle I^1_{\epsilon,\delta}(m_{\mini})=\int_{H(F)} \Xi^G(hm_{\mini})e^{\epsilon\sigma(h)} \left(1+\lvert \lambda(h)\rvert\right)^{-\delta}dh$$

\noindent is absolutely convergent for all $m_{\mini}\in M_{\mini}(F)$ and there exists $d>0$ such that

$$I^1_{\epsilon,\delta}(m_{\mini})\ll \delta_{\overline{P}_{\mini}}(m_{\mini})^{-1/2}\sigma(m_{\mini})^d$$

\noindent for all $m_{\mini}\in M_{\mini}(F)$.

\item Assume moreover that $A$ is contained in $A_{M_{\mini}}$. Then, for all $\delta>0$ there exists $\epsilon>0$ such that the integral

$$\displaystyle I^2_{\epsilon,\delta}(m_{\mini})=\int_{H(F)}\int_{H(F)} \Xi^G(hm_{\mini})\Xi^G(h'hm_{\mini})e^{\epsilon\sigma(h)}e^{\epsilon\sigma(h')} \left(1+\lvert \lambda(h')\rvert\right)^{-\delta}dh'dh$$

\noindent is absolutely convergent for all $m_{\mini}\in M_{\mini}(F)$ and there exists $d>0$ such that

$$I^2_{\epsilon,\delta}(m)\ll \delta_{\overline{P}_{\mini}}(m_{\mini})^{-1}\sigma(m_{\mini})^d$$

\noindent for all $m_{\mini}\in M_{\mini}(F)$.

\end{enumerate}
\end{lem}

\vspace{2mm}

\noindent\ul{Proof}: 

\begin{enumerate}[(i)]

\item This follows from the following fact

\vspace{3mm}

\begin{num}
\item\label{eq 6.5.1} There exists $\epsilon'>0$ such that
$$\Xi^{G_0}(h_0)\ll \Xi^{H_0}(h_0)^2e^{-\epsilon' \sigma(h_0)}$$
for all $h_0\in H_0(F)$.
\end{num}

\vspace{3mm}

If $F$ is $p$-adic, this is proved in \cite{Beu1} (Lemme 12.0.5). The proof works equally well in the real case.

\item Let $d>0$. By Proposition \ref{proposition 1.5.1}(iv), if $d$ is sufficiently large, we have

\[\begin{aligned}
\displaystyle \int_{H(F)}\Xi^G(h) \sigma(h)^{-d}dh & =\int_{H_0(F)}\int_{N(F)} \Xi^G(h_0n)\sigma(h_0n)^{-d}dndh_0 \\
 & \ll \int_{H_0(F)} \Xi^{G_0}(h_0)dh_0
\end{aligned}\]

\noindent (Note that $\delta_P(h_0)=1$ and $\Xi^M(h_0)=\Xi^{G_0}(h_0)$ for all $h_0\in H_0(F)$) and this last integral is absolutely convergent by (i).

\item By (i) and since $\sigma(h_0n)\ll \sigma(h_0)+\sigma(n)$ for all $h_0\in H_0(F)$ and all $n\in N(F)$, it suffices to establish

\vspace{3mm}

\begin{num}
\item\label{eq 6.5.2} For all $\delta>0$ and all $\epsilon_0>0$, there exists $\epsilon>0$ such that the integral
$$I^0_{\epsilon,\delta}(h_0)=\int_{N(F)}\Xi^G(nh_0)e^{\epsilon\sigma(n)} \left(1+\lvert \lambda(n)\rvert\right)^{-\delta}dn$$
is absolutely convergent for all $h_0\in H_0(F)$ and satisfies the inequality
$$I^0_{\epsilon,\delta}(h_0)\ll \Xi^{G_0}(h_0)e^{\epsilon_0\sigma(h_0)}$$
for all $h_0\in H_0(F)$.
\end{num}

\vspace{3mm}

Let $\delta>0$, $\epsilon_0>0$ and $\epsilon>0$. We want to prove that \ref{eq 6.5.2} holds if $\epsilon$ is sufficiently small (compared to $\delta$ and $\epsilon_0$). We shall introduce an auxiliary parameter $b>0$ that we will precise later. For all $h_0\in H_0(F)$, we have $I^0_{\epsilon,\delta}(h_0)=I^0_{\epsilon,\delta,\leqslant b}(h_0)+I^0_{\epsilon,\delta,>b}(h_0)$ where

$$\displaystyle I^0_{\epsilon,\delta,\leqslant b}(h_0)=\int_{N(F)} \mathbf{1}_{\sigma\leqslant b}(n)\Xi^G(nh_0) e^{\epsilon \sigma(n)} \left(1+\lvert \lambda(n)\rvert\right)^{-\delta} dn$$

$$\displaystyle I^0_{\epsilon,\delta,> b}(h_0)=\int_{N(F)} \mathbf{1}_{\sigma> b}(n)\Xi^G(nh_0) e^{\epsilon \sigma(n)} \left(1+\lvert \lambda(n)\rvert\right)^{-\delta} dn$$

\noindent For all $d>0$, we have

$$\displaystyle I^0_{\epsilon,\delta,\leqslant b}(h_0)\leqslant e^{\epsilon b}b^{d}\int_{N(F)} \Xi^G(nh_0) \sigma(n)^{-d}dn$$

\noindent for all $h_0\in H_0(F)$ and all $b>0$. By Proposition \ref{proposition 1.5.1}(iv), we may choose $d>0$ such that the last integral above is essentially bounded by $\delta_P(h_0)^{1/2}\Xi^M(h_0)=\Xi^{G_0}(h_0)$ for all $h_0\in H_0(F)$. We henceforth fix such a $d>0$. Hence, we have

\begin{align}\label{eq 6.5.3}
\displaystyle I^0_{\epsilon,\delta,\leqslant b}(h_0)\ll e^{\epsilon b}b^{d}\Xi^{G_0}(h_0)
\end{align}

\noindent for all $h_0\in H_0(F)$ and all $b>0$.

\vspace{2mm}

\noindent There exists $\alpha>0$ such that $\Xi^G(g_1g_2)\ll e^{\alpha\sigma(g_2)} \Xi^G(g_1)$ for all $g_1,g_2\in G(F)$. It follows that

\begin{align}\label{eq 6.5.4}
\displaystyle I^0_{\epsilon,\delta,> b}(h_0)\ll e^{\alpha\sigma(h_0)-\sqrt{\epsilon}b}\int_{N(F)} \Xi^G(n) e^{(\epsilon+\sqrt{\epsilon}) \sigma(n)} \left(1+\lvert \lambda(n)\rvert\right)^{-\delta} dn
\end{align}

\noindent for all $h_0\in H_0(F)$ and all $b>0$. Assume one moment that the last integral above is convergent if $\epsilon$ is sufficiently small. Taking $\epsilon$ that sufficiently small and combining \ref{eq 6.5.3} with \ref{eq 6.5.4}, we get

$$\displaystyle I^0_{\epsilon,\delta}(h_0)\ll e^{\epsilon b}b^{d} \Xi^{G_0}(h_0)+e^{\alpha\sigma(h_0)-\sqrt{\epsilon}b}$$

\noindent for all $h_0\in H_0(F)$ and for all $b>0$. There exists $\beta>0$ such that $e^{-\beta \sigma(h_0)}\ll  \Xi^{G_0}(h_0)$ for all $h_0\in H_0(F)$. Plugging $b=\frac{\alpha+\beta}{\sqrt{\epsilon}}\sigma(h_0)$ in the last inequality, we obtain

$$\displaystyle I^0_{\epsilon,\delta}(h_0)\ll e^{\sqrt{\epsilon}(\alpha+\beta+1)\sigma(h_0)} \Xi^{G_0}(h_0)$$

\noindent for all $h_0\in H_0(F)$. Hence, for $\epsilon\leqslant \epsilon_0^2 (\alpha+\beta+1)^{-2}$, \ref{eq 6.5.2} indeed holds.

\vspace{2mm}

\noindent It remains to prove the convergence of the integral on the right hand side of \ref{eq 6.5.4} for $\epsilon$ sufficiently small. If $P$ is a minimal parabolic subgroup of $G$ then it follows from Corollary \ref{corollary B.3.1} (since in this case $\lambda$ is a generic additive character of $N$). Assume this is not the case. Then, we can find two isotropic vectors $z_{0,+},z_{0,-}\in V_0$ such that $z_0=z_{0,+}-z_{0,-}$. We have a decomposition $\lambda=\lambda_+-\lambda_-$ where

$$\displaystyle \lambda_+(n)=\Tra_{E/F}\left(\sum_{i=1}^{r-1} h(z_{-i-1},nz_i)+h(z_{-1},nz_{0,+})\right),\;\;\; n\in N$$

$$\displaystyle \lambda_-(n)=\Tra_{E/F}\left(h(z_{-1},nz_{0,-})\right),\;\;\; n\in N$$

\noindent Note that the additive character $\lambda_+$ is the restriction to $N$ of a generic additive character of the unipotent radical of a minimal parabolic subgroup contained in $P$. Hence, Corollary \ref{corollary B.3.1} applies to $\lambda_+$. Choose a one-parameter subgroup $a:\mathbb{G}_m\to M$ such that $\lambda_+(a(t)na(t)^{-1})=t\lambda_+(n)$ and $\lambda_-(a(t)na(t)^{-1})=t^{-1}\lambda_-(n)$ for all $t\in \mathbb{G}_m$ and all $n\in N$ (such a one-parameter subgroup is easily seen to exist). Let $\mathcal{U}\subset F^\times$ be a compact neighborhood of $1$. Then, for all $\epsilon>0$, we have

\[\begin{aligned}
\displaystyle \int_{N(F)} \Xi^G(n) e^{\epsilon \sigma(n)} \left(1+\lvert \lambda(n)\rvert\right)^{-\delta}  dn & \ll \int_{N(F)} \Xi^G(n) e^{\epsilon \sigma(n)} \left(1+\lvert \lambda(a(t)na(t)^{-1})\rvert\right)^{-\delta} dn \\
 & =\int_{N(F)} \Xi^G(n) e^{\epsilon \sigma(n)} \left(1+\lvert t\lambda_+(n)-t^{-1}\lambda_-(n)\rvert\right)^{-\delta} dn
\end{aligned}\]

\noindent for all $t\in \mathcal{U}$. Integrating this last inequality over $\mathcal{U}$, we get that for all $\epsilon>0$, we have

$$\displaystyle \int_{N(F)} \Xi^G(n) e^{\epsilon \sigma(n)} \left(1+\lvert \lambda(n)\rvert\right)^{-\delta}  dn  \ll \int_{N(F)} \Xi^G(n) e^{\epsilon \sigma(n)} \int_{\mathcal{U}}\left(1+\lvert t\lambda_+(n)-t^{-1}\lambda_-(n)\rvert\right)^{-\delta}dt dn$$

\noindent By Lemma \ref{lemma B.1.1}, there exists $\delta'>0$ depending only on $\delta>0$ such that the last expression above is essentially bounded by

$$\displaystyle \int_{N(F)} \Xi^G(n) e^{\epsilon \sigma(n)} \left(1+\lvert \lambda_+(n)\rvert\right)^{-\delta'}dn$$

\noindent Now by Corollary \ref{corollary B.3.1}, this last integral is convergent if $\epsilon$ is sufficiently small.

\item Let $\delta>0$ and $\epsilon>0$. We want to show that (iv) holds if $\epsilon$ is sufficiently small (compared to $\delta$). Since $\Xi^G(g^{-1})\sim \Xi^G(g)$, $\sigma(g^{-1})\sim \sigma(g)$ and $\lambda(h^{-1})=-\lambda(h)$ for all $g\in G(F)$ and all $h\in H(F)$, it is equivalent to show the following

\vspace{3mm}

\begin{num}
\item\label{eq 6.5.5} If $\epsilon$ is sufficiently small the integral
$$\displaystyle J^1_{\epsilon,\delta}(m_{\mini})=\int_{H(F)} \Xi^G(m_{\mini}h)e^{\epsilon\sigma(h)} \left(1+\lvert \lambda(h)\rvert\right)^{-\delta} dh$$
is absolutely convergent for all $m_{\mini}\in M_{\mini}(F)$ and there exists $d>0$ such that
$$J^1_{\epsilon,\delta}(m_{\mini})\ll \delta_{\overline{P}_{\mini}}(m_{\mini})^{1/2}\sigma(m_{\mini})^d$$
for all $m_{\mini}\in M_{\mini}(F)$.
\end{num}

\vspace{3mm}

Let $K$ be a maximal compact subgroup of $G(F)$ that is special in the $p$-adic case. Fix a map $m_{\overline{P}_{\mini}}:G(F)\to M_{\mini}(F)$ such that $g\in m_{\overline{P}_{\mini}}(g)\overline{U}_{\mini}(F)K$ for all $g\in G(F)$. By Proposition \ref{proposition 1.5.1}(ii), there exists $d>0$ such that we have

$$\displaystyle J^1_{\epsilon,\delta}(m_{\mini})\ll \delta_{\overline{P}_{\mini}}(m_{\mini})^{1/2}\sigma(m_{\mini})^d \int_{H(F)} \delta_{\overline{P}_{\mini}}(m_{\overline{P}_{\mini}}(h))^{1/2} \sigma(h)^d e^{\epsilon\sigma(h)} \left(1+\lvert \lambda(h)\rvert\right)^{-\delta} dh$$

\noindent for all $m_{\mini}\in M_{\mini}(F)$. Of course, for any $\epsilon'>\epsilon$ we have $\sigma(h)^de^{\epsilon\sigma(h)}\ll e^{\epsilon'\sigma(h)}$, for all $h\in H(F)$. Hence, we only need to prove that for $\epsilon$ sufficiently small the integral

\begin{align}\label{eq 6.5.6}
\displaystyle  \int_{H(F)} \delta_{\overline{P}_{\mini}}(m_{\overline{P}_{\mini}}(h))^{1/2} e^{\epsilon\sigma(h)} \left(1+\lvert \lambda(h)\rvert\right)^{-\delta} dh
\end{align}

\noindent is absolutely convergent. Since $\overline{P}_{\mini}$ is a good parabolic subgroup, we may find compact neighborhood of the identity $\mathcal{U}_K\subset K$, $\mathcal{U}_H\subset H(F)$ and $\mathcal{U}_{\overline{P}}\subset \overline{P}_{\mini}(F)$ such that $\mathcal{U}_K\subset \mathcal{U}_{\overline{P}}\mathcal{U}_H$. We have inequalities

$$e^{\epsilon\sigma(k_Hh)}\ll e^{\epsilon \sigma(h)} \mbox{ and } \left(1+\lvert \lambda(k_Hh)\rvert\right)^{-\delta}\ll \left(1+\lvert \lambda(h)\rvert\right)^{-\delta}$$

\noindent for all $h\in H(F)$ and for all $k_H\in \mathcal{U}_H$. Hence, we have

\[\begin{aligned}
\displaystyle & \int_{H(F)} \delta_{\overline{P}_{\mini}}(m_{\overline{P}_{\mini}}(h))^{1/2} e^{\epsilon\sigma(h)} \left(1+\lvert \lambda(h)\rvert\right)^{-\delta} dh \\
 & \ll \delta_{\overline{P}_{\mini}}(k_{\overline{P}})^{1/2}  \int_{H(F)} \delta_{\overline{P}_{\mini}}(m_{\overline{P}_{\mini}}(k_Hh))^{1/2} e^{\epsilon\sigma(h)} \left(1+\lvert \lambda(h)\rvert\right)^{-\delta} dh \\
 & =\int_{H(F)} \delta_{\overline{P}_{\mini}}(m_{\overline{P}_{\mini}}(k_{\overline{P}}k_Hh))^{1/2} e^{\epsilon\sigma(h)} \left(1+\lvert \lambda(h)\rvert\right)^{-\delta} dh
\end{aligned}\]

\noindent for all $k_H\in \mathcal{U}_H$ and all $k_{\overline{P}}\in \mathcal{U}_{\overline{P}}$. It follows that

\[\begin{aligned}
\displaystyle & \int_{H(F)} \delta_{\overline{P}_{\mini}}(m_{\overline{P}_{\mini}}(h))^{1/2} e^{\epsilon\sigma(h)} \left(1+\lvert \lambda(h)\rvert\right)^{-\delta} dh \\
 & \ll \int_{H(F)} \int_{\mathcal{U}_K} \delta_{\overline{P}_{\mini}}(m_{\overline{P}_{\mini}}(kh))^{1/2}dk e^{\epsilon\sigma(h)} \left(1+\lvert \lambda(h)\rvert\right)^{-\delta} dh \\
 & \ll \int_{H(F)} \int_{K} \delta_{\overline{P}_{\mini}}(m_{\overline{P}_{\mini}}(kh))^{1/2}dk e^{\epsilon\sigma(h)} \left(1+\lvert \lambda(h)\rvert\right)^{-\delta} dh
\end{aligned}\]

\noindent By Proposition \ref{proposition 1.5.1}(iii), the inner integral above is equal to $\Xi^G(h)$ (for a suitable normalization) and the convergence of \ref{eq 6.5.6} for $\epsilon$ sufficiently small now follows from (iii).

\item Let $\delta>0$ and $\epsilon>0$. We want to prove that (v) holds if $\epsilon$ is sufficiently small (compared to $\delta$). After the variable change $h'\mapsto h'h^{-1}$, we are left with proving that for $\epsilon>0$ sufficiently small the integral

$$\displaystyle I^3_{\epsilon,\delta}(m_{\mini})=\int_{H(F)}\int_{H(F)} \Xi^G(hm_{\mini})\Xi^G(h'm_{\mini}) e^{\epsilon\sigma(h)}e^{\epsilon\sigma(h')} \left(1+\lvert\lambda(h')-\lambda(h)\rvert\right)^{-\delta} dh'dh$$

\noindent is absolutely convergent for all $m_{\mini}\in M_{\mini}(F)$ and that there exists $d>0$ such that

\begin{align}\label{eq 6.5.7}
I^3_{\epsilon,\delta}(m_{\mini})\ll \delta_{\overline{P}_{\mini}}(m_{\mini})^{-1}\sigma(m_{\mini})^d
\end{align}

\noindent for all $m_{\mini}\in M_{\mini}(F)$. Let $a:\mathbb{G}_m\to A$ be a one-parameter subgroup such that $\lambda(a(t)ha(t)^{-1})=t\lambda(h)$ for all $h\in H$ and all $t\in \mathbb{G}_m$. Let $\mathcal{U}\subset F^\times$ be a compact neighborhood of $1$. Since $A$ is in the center of $M_{\mini}$, we have the inequality

\[\begin{aligned}
\displaystyle & I^3_{\epsilon,\delta}(m_{\mini}) \\
 & \ll \int_{H(F)}\int_{H(F)} \Xi^G(hm_{\mini}) \Xi^G(h'm_{\mini}) e^{\epsilon\sigma(h)}e^{\epsilon\sigma(h')} \int_{\mathcal{U}} \left(1+\lvert \lambda(a(t)h'a(t)^{-1})-\lambda(h)\rvert\right)^{-\delta} dtdh'dh \\
 & = \int_{H(F)}\int_{H(F)} \Xi^G(hm_{\mini}) \Xi^G(h'm_{\mini}) e^{\epsilon\sigma(h)}e^{\epsilon\sigma(h')} \int_{\mathcal{U}} \left(1+\lvert t\lambda(h')-\lambda(h)\rvert\right)^{-\delta} dtdh'dh
\end{aligned}\]

\noindent for all $m_{\mini}\in M_{\mini}(F)$. By Lemma \ref{lemma B.1.1}, there exists $\delta'>0$ depending only on $\delta$ such that the last integral above is essentially bounded by

$$\displaystyle \int_{H(F)}\int_{H(F)} \Xi^G(hm_{\mini})\Xi^G(h'm_{\mini}) e^{\epsilon\sigma(h)}e^{\epsilon\sigma(h')}\left(1+\lvert \lambda(h')\rvert\right)^{-\delta'}\left(1+\lvert \lambda(h)\rvert\right)^{-\delta'} dh'dh$$

\noindent for all $m_{\mini}\in M_{\mini}(F)$. This last integral is equal to $I^1_{\epsilon,\delta'}(m_{\mini})^2$. Hence, the inequality \ref{eq 6.5.7} for $\epsilon$ sufficiently small now follows from (iv). $\blacksquare$

\end{enumerate}

\subsection{Relative weak Cartan decompositions}\label{section 6.6}

\subsubsection{Relative weak Cartan decomposition for $G_0$}\label{section 6.6.1}

\noindent Recall that in Section \ref{section 6.2}, we have defined two subgroups $G_0$ and $H_0$ of $G$ and $H$ respectively. The triple $(G_0,H_0,1)$ is a GGP triple which fall into the ``codimension one case". Of course, Proposition \ref{proposition 6.4.1} applies as well to this case. In particular, $G_0$ admits good minimal parabolic subgroups. Let $\gls{P0bar}=M_0\overline{U}_0$ be such a minimal parabolic subgroup of $G_0$ and denote by $\gls{A_0}=A_{M_0}$ the maximal central split subtorus of $M_0$. Set

$$\displaystyle \gls{A0+}=\{a\in A_0(F);\; \lvert \alpha(a)\rvert\geqslant 1\; \forall \alpha\in R(A_0,\overline{P}_0)\}$$

\begin{prop}\label{proposition 6.6.1.1}
There exists a compact subset $\mathcal{K}_0\subseteq G_0(F)$ such that

$$G_0(F)=H_0(F)A_0^+\mathcal{K}_0$$
\end{prop}

\vspace{2mm}

\noindent\ul{Proof}: First, we prove that the result doesn't depend on the particular pair $(\overline{P}_0,M_0)$ that has been chosen. Let $(\overline{P}_0',M_0')$ be another such pair i.e., $\overline{P}_0'$ is a good minimal parabolic subgroup of $G_0$ and $M_0'$ is a Levi component of it. By Proposition \ref{proposition 6.4.1}(i), there exists $h\in H(F)$ such that $\overline{P}'_0=h\overline{P}_0h^{-1}$. Obviously, the result of the proposition for $(\overline{P}_0,M_0)$ implies the same result for the pair $(h\overline{P}_0h^{-1},hM_0h^{-1})=(\overline{P}'_0,hM_0h^{-1})$. Moreover, there exists a $\overline{p}_0'\in \overline{P}_0'(F)$ such that $hM_0h^{-1}=\overline{p}_0'M'_0{\overline{p}'}_0^{-1}$. The result for the pair $(\overline{P}'_0,hM_0h^{-1})$ now implies the result for $(\overline{P}_0',M_0')$ because, by definition of $A_0'^+$, the set

$$\displaystyle \{a_0'^{-1}\overline{p}_0'a_0',\; a_0'\in A_0'^+\}$$

\noindent is bounded. Thus, it suffices to prove that the proposition holds for one particular pair $(\overline{P}_0,M_0)$. In the $p$-adic case, this follows from Proposition 11.0.1 of \cite{Beu1}. We could argue that in the real case the same proof works. Instead, we prefer to rely on the main result of \cite{KKSS}. Fix a good minimal parabolic subgroup $\overline{P}_0\subseteq G_0$. By Proposition \ref{proposition 6.4.1}(i), there exists a Levi component $M_0$ of $\overline{P}_0$ such that $H_0\cap \overline{P}_0\subseteq M_0$. By Theorem 5.13 of \cite{KKSS}, there exists a compact subset $\mathcal{K}_0\subseteq G_0(\mathbb{R})$ such that

$$G_0(\mathbb{R})=H_0(\mathbb{R})F''A_Z^-\mathcal{K}_0$$

\noindent where $A_Z^{-}$ is a certain submonoid of $A_0(\mathbb{R})$ (the exponential of the so-called ``compression cone" associated to the real spherical variety $Z=H_0(\mathbb{R})\backslash G_0(\mathbb{R})$, cf.\ Section 5.1 of \cite{KKSS}) and $F''$ is a subset of $N_{G_0(\mathbb{R})}(H_0)F$, $F$ being any set of representatives for the open $H_0(\mathbb{R})\times \overline{P}_0(\mathbb{R})$ double cosets in $G_0(\mathbb{R})$. By Proposition \ref{proposition 6.4.1}(i), we can take $F=\{1\}$. Moreover, we easily check that $N_{G_0(\mathbb{R})}(H_0)=H_0(\mathbb{R})Z_{G_0}(\mathbb{R})$. As $Z_{G_0}(\mathbb{R})$ is compact, up to multiplying $\mathcal{K}_0$ by $Z_{G_0}(\mathbb{R})$, we may also assume that $F''=\{1\}$. To end the proof of the proposition, it only remains to see that $A_Z^-\subseteq A_0^+$ (note that our convention for the positive chamber is the opposite to that of \cite{KKSS}, this is because we are denoting our good parabolic subgroup by $\overline{P}_0$ and not by $P_0$). But this follows from the fact that the real spherical variety $Z=H_0(\mathbb{R})\backslash G_0(\mathbb{R})$ is wavefront (cf.\ Definition 6.1 of \cite{KKSS} noting that here $\mathfrak{a}_H=0$, the notion of wavefront spherical variety has been first introduced in \cite{SV}). To see this, we can proceed as follows. Consider the complex homogeneous space $Z_{\mathbb{C}}=H_0(\mathbb{C})\backslash G_0(\mathbb{C})\simeq GL_{d-1}(\mathbb{C})\backslash \left(GL_{d-1}(\mathbb{C})\times GL_d(\mathbb{C})\right)$. It is spherical (it follows for example from Proposition \ref{proposition 6.4.1}(i) applied to GGP triples of codimension one with $G$ quasi-split) and wavefront by Remark 6.2 of \cite{KKSS}. On the other hand, it is easy to see from the characterization of the compression cone given in Lemma 5.9 of \cite{KKSS} that a sufficient condition for a real spherical variety to be wavefront is that its complexification is spherical and wavefront. Thus $Z$ is wavefront and this ends the proof of the proposition in the real case. $\blacksquare$

\subsubsection{Relative weak Cartan decomposition for $G$}\label{section 6.6.2}

\noindent Let the quadruple $\left(\overline{P}_0,M_0,A_0,A_0^+\right)$ be as in the previous section. Denote by $\overline{P}=M\overline{N}$ the parabolic subgroup opposite to $P$ with respect to $M$ and define the following subgroups of $G$:

$$A_{\mini}=A_0A\subseteq M_{\mini}=M_0T\subseteq \overline{P}_{\mini}=\overline{P}_0T\overline{N}$$

\noindent Then, $\overline{P}_{\mini}$ is a parabolic subgroup, $M_{\mini}$ is a Levi component of it and $A_{\mini}$ is the maximal split central subtorus of $M_{\mini}$. Moreover, it is easy to see that $\overline{P}_{\mini}$ is a good parabolic subgroup of $G$. Set

$$A_{\mini}^+=\{a\in A_{\mini}(F);\; \lvert \alpha(a)\rvert\geqslant 1 \; \forall \alpha\in R(A_{\mini},\overline{P}_{\mini})\}$$

\noindent We will denote by $P_{\mini}$ the parabolic subgroup opposite to $\overline{P}_{\mini}$ with respect to $M_{\mini}$. We have $P_{\mini}\subseteq P$. Let $\Delta$ be the set of simple roots of $A_{\mini}$ in $P_{\mini}$ and $\gls{DeltaP}=\Delta\cap R(A_{\mini},N)$ be the subset of simple roots appearing in $\mathfrak{n}=\Lie(N)$. For $\alpha\in \Delta_P$, we will denote by $\mathfrak{n}_\alpha$ the corresponding root subspace. Recall also that we have defined in Section \ref{section 6.2}, a character $\xi$ of $\mathfrak{n}(F)$.

\begin{lem}\label{lemma 6.6.2.1}
We have the following

\begin{enumerate}[(i)]

\item $$A_{\mini}^+=\{a\in A_0^+A(F);\; \lvert \alpha(a)\rvert\leqslant 1\; \forall \alpha\in \Delta_P\}$$

\item There exists a compact subset $\mathcal{K}_G\subseteq G(F)$ such that

$$G(F)=H(F)A_0^+A(F)\mathcal{K}_G$$

\item For all $\alpha\in \Delta_P$, the restriction of $\xi$ to $\mathfrak{n}_\alpha(F)$ is nontrivial.

\end{enumerate}
\end{lem}

\vspace{2mm}

\noindent\ul{Proof}: (i) is obvious, so we only provide a proof of (ii) and (iii).

\begin{enumerate}[(i)]
\setcounter{enumi}{1}

\item Let $K$ be a maximal compact subgroup of $G(F)$ which is special in the $p$-adic case. Then we have the Iwasawa decomposition

\begin{align}\label{eq 6.6.1}
G(F)=P(F)K=N(F)G_0(F)T(F)K
\end{align}

Since $A=A_T$ is the maximal split subtorus of $T$, there exists a compact subset $\mathcal{K}_T\subseteq T(F)$ such that

\begin{align}\label{eq 6.6.2}
T(F)=A(F)\mathcal{K}_T
\end{align}

\noindent Also by Proposition \ref{proposition 6.6.1.1}, we know there exists a compact subset $\mathcal{K}_0\subseteq G_0(F)$ such that

\begin{align}\label{eq 6.6.3}
G_0(F)=H_0(F)A_0^+\mathcal{K}_0
\end{align}

\noindent Combining \ref{eq 6.6.1}, \ref{eq 6.6.2} and \ref{eq 6.6.3}, and since $A$ and $G_0$ centralize each other, we get

$$G(F)=H(F)A_0^+A(F)\mathcal{K}_G$$

\noindent where $\mathcal{K}_G=\mathcal{K}_0\mathcal{K}_TK$.

\item Let $\alpha\in \Delta_P$ and assume, by way of contradiction, that $\xi$ is trivial when restricted to $\mathfrak{n}_\alpha(F)$. Recall that $\xi$ is the composition $\xi=\psi\circ \lambda_F$ where $\lambda$ is an algebraic additive character $\mathfrak{n}\to \mathbb{G}_a$. Since $\mathfrak{n}_\alpha$ is a linear subspace of $\mathfrak{n}$, the condition that $\xi$ is trivial on $\mathfrak{n}_\alpha(F)$ amounts to saying that $\lambda$ is trivial on $\mathfrak{n}_\alpha$. This has the advantage of reducing everything to a statement over $\overline{F}$. Since $\lambda$ is invariant by $H_0$-conjugation and $\mathfrak{n}_\alpha$ is invariant by both $T$-conjugation and $\overline{P}_0$-conjugation, it follows that $\lambda$ is trivial on $m\mathfrak{n}_\alpha m^{-1}$ for all $m\in H_0\overline{P}_0T$. But $\overline{P}_0$ being a good parabolic subgroup of $G_0$, $H_0\overline{P}_0T$ is Zariski-dense in $M=G_0T$. Hence, $\lambda$ is trivial on $m\mathfrak{n}_\alpha m^{-1}$ for all $m\in M$. This is a contradiction in view of Lemma \ref{lemma 6.2.1}(ii) (since $\mathfrak{n}_\alpha$ is not included in $[\mathfrak{n},\mathfrak{n}]$). $\blacksquare$

\end{enumerate}

\subsection{The function $\Xi^{H\backslash G}$}\label{section 6.7}

\noindent Let $C\subseteq G(F)$ be a compact subset with nonempty interior. We define a function $\Xi^{H\backslash G}_C$ on $H(F)\backslash G(F)$ by

$$\displaystyle \Xi^{H\backslash G}_C(x)=\vol_{H\backslash G}(xC)^{-1/2}$$

\noindent for all $x\in H(F)\backslash G(F)$. It is not hard to see that if $C'\subseteq G(F)$ is another compact subset with nonempty interior, we have

$$\displaystyle \Xi^{H\backslash G}_C(x)\sim \Xi^{H\backslash G}_{C'}(x)$$

\noindent for all $x\in H(F)\backslash G(F)$. From now on, we will assume implicitly fixed a compact subset with nonempty interior $C\subseteq G(F)$ and we will set

$$\displaystyle \gls{XiHG}(x)=\Xi^{H\backslash G}_C(x)$$

\noindent for all $x\in H(F)\backslash G(F)$. The precise choice of $C$ won't matter because the function $\Xi^{H\backslash G}$ will only be used for the purpose of estimates.

\begin{prop}\label{proposition 6.7.1}
\begin{enumerate}[(i)]

\item For every compact subset $\mathcal{K}\subseteq G(F)$, we have the following equivalences of functions

\vspace{2mm}

(a)\hfill $\Xi^{H\backslash G}(xk)\sim\Xi^{H\backslash G}(x)$ \hfill.

\vspace{2mm}

(b)\hfill $\sigma_{H\backslash G}(xk)\sim \sigma_{H\backslash G}(x)$ \hfill.

\vspace{2mm}

\noindent for all $x\in H(F)\backslash G(F)$ and all $k\in \mathcal{K}$.

\item Let $\overline{P}_0=M_0\overline{U}_0\subseteq G_0$ be a good minimal parabolic subgroup of $G_0$ and $A_0=A_{M_0}$ be the split part of the center of $M_0$. Set

$$A_0^+=\{a_0\in A_0(F);\; \lvert \alpha(a_0)\rvert\geqslant 1 \; \forall \alpha\in R(A_0,\overline{P}_0)\}$$

\noindent then there exists a positive constant $d>0$ such that

\vspace{2mm}

(a)\hfill $\Xi^{G_0}(a_0)\delta_P(a)^{1/2}\sigma(a_0)^{-d}\ll\Xi^{H\backslash G}(aa_0)\ll \Xi^{G_0}(a_0)\delta_P(a)^{1/2}$ \hfill.

\vspace{2mm}

(b)\hfill $\sigma_{H\backslash G}(aa_0)\sim \sigma_G(aa_0)$ \hfill.

\vspace{2mm}

\noindent for all $a_0\in A_0^+$ and all $a\in A(F)$.

\item There exists $d>0$ such that the integral

$$\displaystyle \int_{H(F)\backslash G(F)} \Xi^{H\backslash G}(x)^2 \sigma_{H\backslash G}(x)^{-d} dx$$

\noindent is absolutely convergent.

\item For all $d>0$, there exists $d'>0$ such that

$$\displaystyle \int_{H(F)\backslash G(F)} \mathbf{1}_{\sigma_{H\backslash G}\leqslant c}(x)\Xi^{H\backslash G}(x)^2\sigma_{H\backslash G}(x)^d dx\ll c^{d'}$$

\noindent for all $c\geqslant 1$.

\item There exist $d>0$ and $d'>0$ such that

$$\displaystyle \int_{H(F)} \Xi^G(x^{-1}hx) \sigma_G(x^{-1}hx)^{-d} dh\ll\Xi^{H\backslash G}(x)^2\sigma_{H\backslash G}(x)^{d'}$$

\noindent for all $x\in H(F)\backslash G(F)$.

\item For all $d>0$, there exists $d'>0$ such that

$$\displaystyle \int_{H(F)} \Xi^G(hx)\sigma(hx)^{-d'}dh\ll \Xi^{H\backslash G}(x)\sigma_{H\backslash G}(x)^{-d}$$

\noindent for all $x\in H(F)\backslash G(F)$.

\item Let $\delta>0$ and $d>0$. Then, the integral

$$\displaystyle I_{\delta,d}(c,x)=\int_{H(F)}\int_{H(F)} \mathbf{1}_{\sigma\geqslant c}(h')\Xi^G(hx)\Xi^G(h'hx)\sigma(hx)^d\sigma(h'hx)^d \left(1+\lvert\lambda(h') \rvert\right)^{-\delta} dh'dh$$

\noindent is absolutely convergent for all $x\in H(F)\backslash G(F)$ and all $c\geqslant 1$. Moreover, there exist $\epsilon>0$ and $d'>0$ such that

$$I_{\delta,d}(c,x)\ll \Xi^{H\backslash G}(x)^2\sigma_{H\backslash G}(x)^{d'} e^{-\epsilon c}$$

\noindent for all $x\in H(F)\backslash G(F)$ and all $c\geqslant 1$.

\end{enumerate}
\end{prop}

\vspace{2mm}

\noindent\ul{Proof}: Since $G\to H\backslash G$ has the norm descent property (Lemma \ref{lemma 6.2.1}(i)), we may assume (and we will) that

\begin{align}\label{eq 6.7.1}
\sigma_{H\backslash G}(x)=\inf_{h\in H(F)} \sigma_G(hx)
\end{align}

\noindent for all $x\in H(F)\backslash G(F)$.

\begin{enumerate}[(i)]

\item is easy and left to the reader.

\item

\begin{enumerate}[(a)]

\item Let $\overline{P}=M\overline{N}$ be the parabolic subgroup opposite to $P$ with respect to $M$. Fix compact subsets with nonempty interior

\begin{center}
$C_{\overline{N}}\subseteq \overline{N}(F)$, $C_T\subseteq T(F)$, $C_0\subseteq G_0(F)$ and $C_N\subseteq N(F)$
\end{center}

\noindent Then $C=C_NC_TC_0C_{\overline{N}}$ is a compact subset of $G(F)$ with nonempty interior. We have

\begin{center}
$\Xi^{H\backslash G}(g)\sim \vol_{H\backslash G}\left(H(F)gC\right)^{-1/2}$, for all $g\in G(F)$
\end{center}

\noindent and there exists a $d>0$ such that

\begin{center}
$\Xi^{G_0}(g_0)\sigma(g_0)^{-d}\ll \vol_{G_0}\left(C_0g_0C_0\right)^{-1/2}\ll \Xi^{G_0}(g_0)$, for all $g_0\in G_0(F)$.
\end{center}

\noindent So (ii)(a) is equivalent to

$$\delta_P(a)^{-1}\vol_{G_0}\left(C_0a_0C_0\right)\sim \vol_{H\backslash G}\left(H(F)aa_0C\right)$$

\noindent for all $a_0\in A_0^+$ and all $a\in A(F)$. We have

$$H(F)aa_0C=H(F)aa_0C_{\overline{P}}$$

\noindent for all $a_0\in A_0(F)$ and $a\in A(F)$, where $C_{\overline{P}}=C_TC_0C_{\overline{N}}$. Hence, we need to prove that

\vspace{3mm}

\begin{num}
\item\label{eq 6.7.2} $\delta_P(a)^{-1}\vol_{G_0}\left(C_0a_0C_0\right)\sim \vol_{H\backslash G}\left(H(F)aa_0C_{\overline{P}}\right)$, for all $(a_0,a)\in A_0^+\times A(F)$.
\end{num}

\vspace{3mm}

\noindent Let $C_{H_0}\subseteq H_0(F)$ be a compact subset with nonempty interior and set $C_H=C_NC_{H_0}$. It is a compact subset of $H(F)$ with nonempty interior. We claim that

\vspace{3mm}

\begin{num}
\item\label{eq 6.7.3} $\vol_{H\backslash G}\left(H(F)aa_0C_{\overline{P}}\right)\sim \vol_G\left(C_Haa_0C_{\overline{P}}\right)$, for all $(a_0,a)\in A_0^+\times A(F)$.
\end{num}

\vspace{3mm}

\noindent We have

$$\displaystyle \vol_G\left(C_Haa_0C_{\overline{P}}\right)=\int_{H(F)\backslash G(F)}\int_{H(F)} \mathbf{1}_{C_Haa_0C_{\overline{P}}}(hx)dhdx$$

\noindent for all $(a_0,a)\in A_0^+\times A(F)$. The inner integral above is nonzero only if $x\in H(F)aa_0C_{\overline{P}}$ and is then equal to
 
$$\displaystyle \vol_H\left(H(F)\cap C_Haa_0C_{\overline{P}}x^{-1}\right)=\vol_H\left(C_H\left(H(F)\cap aa_0C_{\overline{P}}x^{-1}\right)\right)$$

\noindent Hence, to get \ref{eq 6.7.3}, it suffices to show that

$$\vol_H\left(C_H\left(H(F)\cap aa_0C_{\overline{P}}x^{-1}\right)\right)\sim 1$$

\noindent for all $(a_0,a)\in A_0^+\times A(F)$ and all $x\in aa_0C_{\overline{P}}$. For such an $x$, we have $C_H\subseteq C_H\left(H(F)\cap aa_0C_{\overline{P}}x^{-1}\right)$, so that we easily get the inequality

$$\vol_H\left(C_H\left(H(F)\cap aa_0C_{\overline{P}}x^{-1}\right)\right)\gg 1$$

\noindent for all $(a_0,a)\in A_0^+\times A(F)$ and all $x\in aa_0C_{\overline{P}}$. Let $C'_{\overline{P}}=C_{\overline{P}}.C_{\overline{P}}^{-1}$. To get the reverse inequality, it suffices to show that the subsets $H(F)\cap aa_0C'_{\overline{P}}(aa_0)^{-1}$ remain uniformly bounded as $(a,a_0)$ runs through $A_0^+\times A(F)$. Since $\overline{P}\cap H=H_0$, we have

$$H(F)\cap aa_0C'_{\overline{P}}(aa_0)^{-1}=H_0(F)\cap a_0C'_0a_0^{-1}$$

\noindent for all $(a_0,a)\in A_0^+\times A(F)$, where $C'_0=C'_{\overline{P}}\cap G_0(F)$. Now, the subsets $H_0(F)\cap a_0C_0'a_0^{-1}$, $a_0\in A_0^+$, are uniformly bounded by Proposition \ref{proposition 6.4.1}(iii). This ends the proof of \ref{eq 6.7.3}.

\vspace{3mm}

By \ref{eq 6.7.3}, \ref{eq 6.7.2} is now equivalent to

\vspace{3mm}

\begin{num}
\item\label{eq 6.7.4} $\delta_P(a)^{-1}\vol_{G_0}\left(C_0a_0C_0\right)\sim \vol_G\left(C_Haa_0C_{\overline{P}}\right)$, for all $(a_0,a)\in A_0^+\times A(F)$.
\end{num}

\vspace{3mm}

\noindent Recall that we have $C_H=C_NC_{H_0}$ and $C_{\overline{P}}=C_TC_0C_{\overline{N}}$. Hence

$$C_Haa_0C_{\overline{P}}=C_N\left(aC_T\right)\left(C_{H_0}a_0C_0\right) C_{\overline{N}}$$

\noindent for all $(a_0,a)\in A_0(F)\times A(F)$. For suitable choices of Haar measures, we have the decomposition $dg=\delta_P(t)^{-1}dndtdg_0d\overline{n}$ where $dn$, $dt$, $dg_0$ and $d\overline{n}$ are Haar measures on respectively $N(F)$, $T(F)$, $G_0(F)$ and $\overline{N}(F)$. From these, it follows easily that

$$\displaystyle \vol_G\left(C_Haa_0C_{\overline{P}}\right)\sim \delta_P(a)^{-1}\vol_{G_0}\left(C_{H_0}a_0C_0\right)$$

\noindent for all $(a_0,a)\in A_0(F)\times A(F)$. Hence, the last thing to show to get \ref{eq 6.7.4} is

\vspace{3mm}

\begin{num}
\item\label{eq 6.7.5} $\vol_{G_0}\left(C_0a_0C_0\right)\sim \vol_{G_0}\left(C_{H_0}a_0C_0\right)$, for all $a_0\in A_0^+$.
\end{num}

\vspace{3mm}

\noindent The inequality $\vol_{G_0}\left(C_{H_0}a_0C_0\right)\ll \vol_{G_0}\left(C_0a_0C_0\right)$ is obvious. So, we only need to prove the reverse one. The choice of $C_0$ doesn't matter. Since $H_0(F)\overline{P}_0(F)$ is open in $G_0(F)$, we may assume that $C_0=C_{H_0}C_{\overline{P}_0}$ where $C_{\overline{P}_0}\subseteq \overline{P}_0(F)$ is a compact subset with nonempty interior. By definition of $A_0^+$, the subsets $a_0^{-1}C_{\overline{P}_0}a_0$ remain uniformly bounded as $a_0$ runs through $A_0^+$. Hence, there exists a compact subset $C_0'\subseteq G_0(F)$ such that

$$a_0^{-1}C_{\overline{P}_0}a_0C_0\subseteq C_0'$$

\noindent for all $a_0\in A_0^+$. From this, we get

$$\displaystyle \vol_{G_0}\left(C_0a_0C_0\right)\leqslant \vol_{G_0}\left(C_{H_0}a_0C'_0\right)\ll \vol_{G_0}\left(C_{H_0}a_0 C_0\right)$$

\noindent for all $a_0\in A_0^+$. This ends the proof of \ref{eq 6.7.5} and hence the proof of (ii)(a).

\item Obviously, we have the inequality $\sigma_{H\backslash G}(g)\ll \sigma_G(g)$, for all $g\in G(F)$. So, we only need to show that

$$\sigma_G(a_0a)\ll\sigma_{H\backslash G}(a_0a)$$

\noindent for all $(a_0,a)\in A_0^+\times A(F)$. Because of \ref{eq 6.7.1}, this is equivalent to the inequality

\vspace{3mm}

\begin{num}
\item\label{eq 6.7.6} $\sigma_G(a_0a)\ll\sigma_G(ha_0a)$, for all $(a_0,a)\in A_0^+\times A(F)$ and all $h\in H(F)$.
\end{num}

\vspace{3mm}

\noindent Every $h\in H(F)$ may be written $h=nh_0$ where $n\in N(F)$, $h_0\in H_0(F)$, and we have

$$\sigma_G(ng_0t)\gg \sigma_{G_0}(g_0)+\sigma_G(t)$$

\noindent for all $n\in N(F)$, $g_0\in G_0(F)$ and $t\in T(F)$. Hence, we have

$$\sigma_G(nh_0a_0a)\gg \sigma_{G_0}(h_0a_0)+\sigma_G(a)$$

\noindent for all $(a_0,a)\in A_0^+\times A(F)$, all $n\in N(F)$ and all $h_0\in H_0(F)$. Since, $\sigma_G(a_0a)\sim \sigma_{G_0}(a_0)+\sigma_G(a)$ for all $(a_0,a)\in A_0(F)\times A(F)$, to get \ref{eq 6.7.6} it suffices to show that

$$\sigma_{G_0}(a_0)\ll \sigma_{G_0}(h_0a_0)$$

\noindent for all $a_0\in A_0^+$ and $h_0\in H_0(F)$. But, this inequality is a straightforward consequence of Proposition \ref{proposition 6.4.1}(iii).

\end{enumerate}

\item Let $C\subseteq G(F)$ be a compact subset with non-empty interior. Let us first show that (iii) follows from the following fact

\vspace{3mm}

\begin{num}
\item\label{eq 6.7.7} There exists a positive integer $N>0$ such that for all $R\geqslant 1$, the subset $B(R)=\{x\in H(F)\backslash G(F);\; \sigma_{H\backslash G}(x)< R\}$ may be covered by less than $(1+R)^N$ subsets of the form $xC$, $x\in H(F)\backslash G(F)$.
\end{num}

\vspace{3mm}

\noindent (we then say that $H(F)\backslash G(F)$ has polynomial growth following \cite{Ber1}). Indeed, set

$$\displaystyle \lambda(R,d)=\int_{B(R+1)\backslash B(R)} \Xi^{H\backslash G}(x)^2 \sigma_{H\backslash G}(x)^{-d} dx$$

\noindent for all $d>0$ and $R\geqslant 1$. Then we have

\begin{align}\label{eq 6.7.8}
\displaystyle \int_{H(F)\backslash G(F)} \Xi^{H\backslash G}(x)^2 \sigma_{H\backslash G}(x)^{-d}dx=\sum_{R=1}^{\infty} \lambda (R,d)
\end{align}

\noindent for all $d>0$. By \ref{eq 6.7.7}, for all $R\geqslant 1$, $B(R+1)\backslash B(R)$ may be covered by subsets $x_1C,\ldots,x_{k_R}C$ where $k_R\leqslant (R+2)^N$. Hence,

\begin{align}\label{eq 6.7.9}
\displaystyle  \lambda(R,d)\leqslant \sum_{i=1}^{k_R} \int_{x_iC} \Xi^{H\backslash G}(x)^2 \sigma_{H\backslash G}(x)^{-d} dx
\end{align}

\noindent for all $d>0$ and all $R\geqslant 1$. By (i).(a) and (i).(b) and the definition of $\Xi^{H\backslash G}$, we have

\[\begin{aligned}
\displaystyle \int_{yC} \Xi^{H\backslash G}(x)^2\sigma_{H\backslash G}(x)^{-d} dx & \ll \vol_{H\backslash G}\left(yC\right) \Xi^{H\backslash G}(y)^2 \sigma_{H\backslash G}(y)^{-d} \\
 & \ll \sigma_{H\backslash G}(y)^{-d}
\end{aligned}\]

\noindent for all $y\in H(F)\backslash G(F)$. Consequently, by \ref{eq 6.7.9}, we get

\begin{align}\label{eq 6.7.10}
\displaystyle \lambda(R,d)\ll \sum_{i=1}^{k_R} \sigma_{H\backslash G}(x_i)^{-d}
\end{align}

\noindent for all $d>0$ and all $R\geqslant 1$. We may of course assume that $x_iC\cap \left(B(R+1)\backslash B(R)\right)\neq \emptyset$ for all $R\geqslant 1$ and all $1\leqslant i\leqslant k_R$. Then by (i).(b), we have

$$\sigma_{H\backslash G}(x_i)^{-1}\ll R^{-1}$$

\noindent for all $R\geqslant 1$ and all $1\leqslant i\leqslant k_R$. Combining this with \ref{eq 6.7.10}, we finally get

$$\lambda(R,d)\ll R^{-d}k_R\leqslant (R+2)^{N}R^{-d}$$

\noindent for all $d>0$ and $R\geqslant 1$. Hence, for $d>N+1$, \ref{eq 6.7.8} is absolutely convergent. This ends the proof that \ref{eq 6.7.7} implies (iii).

\vspace{2mm}

Let us now prove \ref{eq 6.7.7}. By Lemma \ref{lemma 6.6.2.1}(ii), there exists a compact subset $\mathcal{K}\subseteq G(F)$ such that

\begin{align}\label{eq 6.7.11}
G(F)=H(F)A_0^+A(F)\mathcal{K}
\end{align}

\noindent Thus by (i).(b) and (ii).(b), we see that there exists a constant $c_0>0$ such that

$$\displaystyle B(R)\subseteq H(F)\{a_0a;\; a_0\in A_0^+\; a\in A(F)\; \sigma_G(aa_0)\leqslant c_0R\}\mathcal{K}$$

\noindent for all $R\geqslant 1$. Set $A_{\mini}=A_0A$. Using the above, it is easy to see that \ref{eq 6.7.7} is a consequence of the following fact which is not hard to prove and left to the reader.

\vspace{3mm}

\begin{num}
\item Let $C_{\mini}\subseteq A_{\mini}(F)$ be a compact subset with nonempty interior. Then, there exists a positive integer $N>0$ such that for all $R\geqslant 1$, the subset $\{a\in A_{\mini}(F); \sigma_G(a)<R\}$ may be covered by less than $(1+R)^N$ subsets of the form $aC_{\mini}$, $a\in A_{\mini}(F)$.
\end{num}

\vspace{3mm}

\item By similar arguments, this also follows from \ref{eq 6.7.7}.

\item By (i), (ii) and \ref{eq 6.7.11} together with the fact that $\Xi^G(kgk^{-1})\sim \Xi^G(g)$ and $\sigma_G(kgk^{-1})\sim \sigma_G(g)$ for all $k\in \mathcal{K}$ and $g\in G(F)$, we only need to show the existence of $d>0$ such that

\begin{align}\label{eq 6.7.12}
\displaystyle \int_{H(F)} \Xi^G\left(a^{-1}a_0^{-1} ha_0a\right)\sigma_G\left(a^{-1}a_0^{-1}ha_0a\right)^{-d} dh\ll \Xi^{G_0}(a_0)^2 \delta_P(a)
\end{align}

\noindent for all $(a_0,a)\in A_0^+\times A(F)$. We have

\[\begin{aligned}
\displaystyle & \int_{H(F)}  \Xi^G\left(a^{-1}a_0^{-1} ha_0a\right)\sigma_G\left(a^{-1}a_0^{-1}ha_0a\right)^{-d} dh \\
 & =\int_{H_0(F)}\int_{N(F)}  \Xi^G\left(a^{-1}a_0^{-1} h_0na_0a\right)\sigma_G\left(a^{-1}a_0^{-1}h_0na_0a\right)^{-d} dndh_0
\end{aligned}\]

\noindent for all $(a_0,a)\in A_0(F)\times A(F)$ and all $d>0$. After the variable change $n\mapsto a_0ana^{-1}a_0$, the last integral above becomes

$$\displaystyle \delta_P(a) \int_{H_0(F)}\int_{N(F)}  \Xi^G\left(a_0^{-1} h_0a_0n\right)\sigma_G\left(a_0^{-1}h_0a_0n\right)^{-d} dndh_0$$

\noindent By Proposition \ref{proposition 1.5.1}(iv), for $d>0$ sufficiently large, we have

$$\displaystyle \int_{N(F)}  \Xi^G\left(a_0^{-1} h_0a_0n\right)\sigma_G\left(a_0^{-1}h_0a_0n\right)^{-d} dn\ll \Xi^{G_0}\left(a_0^{-1}h_0a_0\right)$$

\noindent for all $a_0\in A_0(F)$ and all $h_0\in H_0(F)$. Hence, for $d>0$ sufficiently large we have

$$\displaystyle  \int_{H(F)}  \Xi^G\left(a^{-1}a_0^{-1} ha_0a\right)\sigma_G\left(a^{-1}a_0^{-1}ha_0a\right)^{-d} dh\ll \delta_P(a)^{-1}\int_{H_0(F)} \Xi^{G_0}\left(a_0^{-1}h_0a_0\right) dh_0$$

\noindent for all $(a_0,a)\in A_0(F)\times A(F)$. Thus to get \ref{eq 6.7.12}, it is sufficient to show that

\begin{align}\label{eq 6.7.13}
\displaystyle \int_{H_0(F)} \Xi^{G_0}\left(a_0^{-1}h_0a_0\right)dh_0\ll \Xi^{G_0}(a_0)^2
\end{align}

\noindent for all $a_0\in A_0^+$. Let $\mathcal{U}_{H_0}\subset H_0(F)$ and $\mathcal{U}_{\overline{P}_0}\subset \overline{P}_0(F)$ be compact neighborhood of the identity. Since the subsets $a_0^{-1}\mathcal{U}_{\overline{P}_0}a_0$ remain uniformly bounded as $a_0$ runs through $A_0^+$, we have

$$\displaystyle \int_{H_0(F)} \Xi^{G_0}\left(a_0^{-1}h_0a_0\right)dh_0\ll \int_{H_0(F)} \Xi^{G_0}\left(a_0^{-1}k_{\overline{P}_0}^1k^1_{H_0}h_0k^2_{H_0}k^2_{\overline{P}_0}a_0\right)dh_0$$

\noindent for all $a_0\in A_0^+$, all $k^1_{H_0},k^2_{H_0}\in \mathcal{U}_{H_0}$ and all $k^1_{\overline{P}_0},k^2_{\overline{P}_0}\in \mathcal{U}_{\overline{P}_0}$. Let $K_0$ be a maximal compact subgroup of $G_0(F)$. Since $\overline{P}_0$ is a good parabolic subgroup of $G_0$, there exists a compact neighborhood of the identity $\mathcal{U}_{K_0}\subset K_0$ such that $\mathcal{U}_{K_0}\subset \mathcal{U}_{\overline{P}_0}\mathcal{U}_{H_0}\cap \mathcal{U}_{H_0}\mathcal{U}_{\overline{P}_0}$. From the last inequality above, we deduce

\[\begin{aligned}
\displaystyle \int_{H_0(F)} \Xi^{G_0}\left(a_0^{-1}h_0a_0\right)dh_0 & \ll \int_{H_0(F)} \int_{\mathcal{U}_{K_0}^2} \Xi^{G_0}\left(a_0^{-1}k^1h_0k^2a_0\right)dk_1dk_2dh_0 \\
& \ll \int_{H_0(F)} \int_{K_0^2} \Xi^{G_0}\left(a_0^{-1}k^1h_0k^2a_0\right)dk_1dk_2dh_0
\end{aligned}\]

\noindent for all $a_0\in A_0^+$. By the ``doubling principle" (Proposition \ref{proposition 1.5.1}(vi)), this last integral is essentially bounded by

$$\displaystyle \Xi^{G_0}(a_0)^2\int_{H_0(F)} \Xi^{G_0}(h_0) dh_0$$

\noindent for all $a_0\in A_0^+$. By Lemma \ref{lemma 6.5.1}(i), the last integral above is convergent. This proves \ref{eq 6.7.13} and ends the proof of (v).

\item By (i), (ii), \ref{eq 6.7.1} and \ref{eq 6.7.11}, it suffices to show the following

\vspace{3mm}

\begin{num}
\item There exist $d>0$ such that
$$\displaystyle \int_{H(F)}\Xi^G(haa_0)\sigma(haa_0)^{-d}dh\ll \delta_P(a)^{1/2}\Xi^{G_0}(a_0)$$
for all $(a_0,a)\in A_0^+\times A(F)$.
\end{num}

\vspace{3mm}

Using again Proposition \ref{proposition 1.5.1}(iv), this will follow from the following inequality

\vspace{3mm}

\begin{num}
\item $\displaystyle \int_{H_0(F)}\Xi^{G_0}(h_0a_0)dh_0\ll \Xi^{G_0}(a_0)$, for all $a_0\in A_0^+$.
\end{num}

\vspace{3mm}

\noindent  To obtain this last inequality, we can argue as in the end of the proof of (v), using the ``doubling principle" (Proposition \ref{proposition 1.5.1}(vi)) and the fact that $\overline{P}_0$ is a good parabolic subgroup of $G_0$ to reduce it to the convergence of the integral

$$\displaystyle \int_{H_0(F)} \Xi^{G_0}(h_0)dh_0$$

\noindent which is a consequence of Lemma \ref{lemma 6.5.1}(i).

\item By (i), (ii) and \ref{eq 6.7.11} and since for all $d>0$ and all $\epsilon>0$ we have $\mathbf{1}_{\sigma\geqslant c}(h)\sigma(h)^d\ll e^{\epsilon\sigma(h)}e^{-\epsilon c/2}$ for all $h\in H(F)$, it suffices to show the following

\vspace{3mm}

\begin{num}
\item\label{eq 6.7.14} For all $\delta>0$, there exist $d>0$ and $\epsilon>0$ such that
\[\begin{aligned}
\displaystyle \int_{H(F)}\int_{H(F)} \Xi^G(haa_0)\Xi^G(h'haa_0) e^{\epsilon \sigma(h')}e^{\epsilon\sigma(h)} & \left(1+\lvert \lambda(h')\rvert\right)^{-\delta} dh'dh \\
 & \ll \delta_P(a) \Xi^{G_0}(a_0)^2\sigma(aa_0)^d
\end{aligned}\]
for all $a_0\in A_0^+$ and all $a\in A(F)$.
\end{num} 

\vspace{3mm}

Let $\delta>0$. Let $\overline{P}=M\overline{N}$ be the parabolic subgroup opposite to $P$ with respect to $M$ and set $\overline{P}_{\mini}=\overline{P}_0T\overline{N}$, $M_{\mini}=M_0T$. Then $\overline{P}_{\mini}$ is a good parabolic subgroup of $G$ and $M_{\mini}$ is a Levi component of it that contains $A$. Hence, by Lemma \ref{lemma 6.5.1}(v), there exists $\epsilon>0$ and $d>0$ such that

$$\displaystyle \int_{H(F)}\int_{H(F)} \Xi^G(haa_0)\Xi^G(h'haa_0)e^{\epsilon \sigma(h')}e^{\epsilon\sigma(h)} \left(1+\lvert \lambda(h')\rvert\right)^{-\delta} dh'dh \ll \delta_{\overline{P}_{\mini}}(aa_0)^{-1}\sigma(aa_0)^d$$

\noindent for all $a_0\in A_0^+$ and all $a\in A(F)$. We have $\delta_{\overline{P}_{\mini}}(aa_0)^{-1}=\delta_P(a)\delta_{\overline{P}_0}(a_0)$ and by Proposition \ref{proposition 1.5.1}(i), we have $\delta_{\overline{P}_0}(a_0)\ll \Xi^{G_0}(a_0)^2$ for all $a_0\in A_0^+$. It follows that the inequality \ref{eq 6.7.14} is satisfied for such an $\epsilon>0$ and such a $d>0$. $\blacksquare$

\end{enumerate}

\subsection{Parabolic degenerations}\label{section 6.8}

Let $\overline{Q}=LU_{\overline{Q}}$ be a good parabolic subgroup of $G$ (recall that it means that $H\overline{Q}$ is Zariski open in $G$ see \S \ref{section 6.4}). Let $\overline{P}_{\mini}=M_{\mini}\overline{U}_{\mini}\subseteq \overline{Q}$ be a good minimal parabolic subgroup of $G$ (Proposition \ref{proposition 6.4.1}(ii)) with the Levi component chosen so that $M_{\mini}\subset L$. Let $A_{\mini}=A_{M_{\mini}}$ be the maximal central split torus of $M_{\mini}$ and set
$$A_{\mini}^+=\{a\in A_{\mini}(F);\;\lvert \alpha(a)\rvert \geqslant 1\;\forall\alpha\in R(A_{\mini},\overline{P}_{\mini})\}$$
Let $H_{\overline{Q}}=H\cap \overline{Q}$ and $H_L$ be the image of $H_{\overline{Q}}$ by the natural pojection $\overline{Q}\twoheadrightarrow L$. Let $Q=LU_Q$ be the parabolic subgroup opposite to $\overline{Q}$ with respect to $L$. We define $H^Q=H_L\ltimes U_Q$.

\begin{prop}\label{proposition 6.8.1}
\begin{enumerate}[(i)]

\item\label{proposition 6.8.1(i)} $H_{\overline{Q}}\cap U_{\overline{Q}}=\{1\}$ so that the natural projection $H_{\overline{Q}}\to H_L$ is an isomorphism;

\item\label{proposition 6.8.1(ii)} $\delta_{\overline{Q}}(h_{\overline{Q}})=\delta_{H_{\overline{Q}}}(h_{\overline{Q}})$ and $\delta_{\overline{Q}}(h_{L})=\delta_{H_{L}}(h_{L})$ for all $h_{\overline{Q}}\in H_{\overline{Q}}(F)$ and all $h_L\in H_L(F)$. In particular, the group $H^Q(F)$ is unimodular.
\end{enumerate}
Fix a left Haar measure $d_Lh_L$ on $H_L(F)$ and a Haar measure $dh^Q$ on $H^Q(F)$.
\begin{enumerate}[(i)]
\setcounter{enumi}{2}
\item\label{proposition 6.8.1(iii)} There exists $d>0$ such that the integral
$$\displaystyle \int_{H_{L}(F)} \Xi^{L}(h_{L}) \sigma(h_{L})^{-d} \delta_{H_{L}}(h_{L})^{1/2} d_L h_{L}$$
converges. Moreover, in the codimension one case (that is when $G=G_0$ and $H=H_0$), the integral
$$\displaystyle \int_{H_{L}(F)} \Xi^{L}(h_{L}) \sigma(h_{L})^{d} \delta_{H_L}(h_{L})^{1/2} d_L h_{L}$$
is convergent for all $d>0$.

\item\label{proposition 6.8.1(iv)} There exists $d>0$ such that the integral
$$\displaystyle \int_{H^Q(F)} \Xi^G(h^Q) \sigma(h^Q)^{-d} dh^Q$$
converges.

\item\label{proposition 6.8.1(v)} We have $\sigma(h^Q)\ll \sigma(a^{-1}h^Qa)$ for all $a\in A_{\mini}^+$ and all $h^Q\in H^Q(F)$.

\item\label{proposition 6.8.1(vi)} There exist $d>0$ and $d'>0$ such that
$$\displaystyle \int_{H^Q(F)} \Xi^{G}(a^{-1}h^Qa) \sigma(a^{-1}h^Qa)^{-d} dh^Q\ll \Xi^{H\backslash G}(a)^2 \sigma_{H\backslash G}(a)^{d'}$$
for all $a\in A_{\mini}^+$.
\end{enumerate}
\end{prop}

\vspace{2mm}

\noindent\ul{Proof}:

\begin{enumerate}[(i)]

\item This follows directly from Proposition \ref{proposition 6.4.1}(i).

\item For $h_{\overline{Q}}\in H_{\overline{Q}}(F)$ which maps to $h_L\in H_L(F)$ via the isomorphism $H_{\overline{Q}}\simeq H_L$, we have $\delta_{\overline{Q}}(h_{\overline{Q}})=\delta_{\overline{Q}}(h_{L})$ and  $\delta_{H_{\overline{Q}}}(h_{\overline{Q}})=\delta_{H_{L}}(h_{L})$. Thus, it suffices to show that $\delta_{\overline{Q}}(h_{\overline{Q}})=\delta_{H_{\overline{Q}}}(h_{\overline{Q}})$ or equivalently

\begin{align}\label{eq 6.8.1}
\det\left( \Ad(h_{\overline{Q}})_{\mid \overline{\mathfrak{q}}/\mathfrak{h}_{\overline{Q}}}\right)=1
\end{align}

\noindent for all $h_{\overline{Q}}\in H_{\overline{Q}}(F)$. We have $\overline{\mathfrak{q}}+\mathfrak{h}=\mathfrak{g}$ (because $\overline{Q}$ is a good parabolic subgroup) and $\mathfrak{h}_{\overline{Q}}=\mathfrak{h}\cap\overline{\mathfrak{q}}$, hence the inclusion $\overline{\mathfrak{q}}\subseteq \mathfrak{g}$ induces an isomorphism $\overline{\mathfrak{q}}/\mathfrak{h}_{\overline{Q}}\simeq \mathfrak{g}/\mathfrak{h}$, from which it follows that

$$\det\left( \Ad(h_{\overline{Q}})_{\mid \overline{\mathfrak{q}}/\mathfrak{h}_{\overline{Q}}}\right)=\det\left( \Ad(h_{\overline{Q}})_{\mid \mathfrak{g}/\mathfrak{h}}\right)=\det\left( \Ad(h_{\overline{Q}})_{\mid \mathfrak{g}}\right) \det\left( \Ad(h_{\overline{Q}})_{\mid \mathfrak{h}}\right)^{-1}$$

\noindent for all $h_{\overline{Q}}\in H_{\overline{Q}}(F)$. But since $G$ and $H$ are unimodular groups, we have $\det\left( \Ad(h_{\overline{Q}})_{\mid \mathfrak{g}}\right)=\det\left( \Ad(h_{\overline{Q}})_{\mid \mathfrak{h}}\right)=1$ for all $h_{\overline{Q}}\in H_{\overline{Q}}$ and \ref{eq 6.8.1} follows.

\item Let $K$ be a maximal compact subgroup of $G(F)$ which is special in good position with respect to $L$ in the $p$-adic case. Set $K_L=K\cap L(F)$ (a maximal compact subgroup of $L(F)$ which is special in the $p$-adic case), $\tau=i_{\overline{P}_{\mini}\cap L}^L(1)$ and $\pi=i_{\overline{Q}}^G(\tau)$. We will denote by $(.,.)$ and $(.,.)_\tau$ invariant scalar products on $\pi$ and $\tau$ respectively. Let $e_K\in \pi^\infty$ and $e_{K_L}\in \tau^\infty$ be the unique $K$-fixed and $K_L$-fixed vectors respectively. Note that we have $e_K(k)=e_{K_L}$ for all $k\in K$. We may assume that the functions $\Xi^G$ and $\Xi^L$ are given by

\begin{align}\label{eq 6.8.2}
\Xi^G(g)=(\pi(g)e_K,e_K)\;\;,\;\; g\in G(F)
\end{align}

\begin{align}\label{eq 6.8.3}
\Xi^L(l)=(\tau(l)e_{K_L},e_{K_L})_\tau \;\;,\;\; l\in L(F)
\end{align}

\noindent (Note that by the process of induction by stages, we have a natural isomorphism $\pi\simeq i_{\overline{P}_{\mini}}^G(1)$). If we choose Haar measures suitably, \ref{eq 6.8.2} gives

$$\displaystyle \Xi^G(g)=\int_{\overline{Q}(F)\backslash G(F)} (e_K(g'g),e_K(g'))_\tau dg'$$

\noindent for all $g\in G(F)$. Since $\overline{Q}$ is a good parabolic subgroup, by Proposition \ref{proposition 6.4.1}(i) (and since $g\mapsto g^{-1}$ is an automorphism of $G$) the subset $H_{\overline{Q}}(F)\backslash H(F)\subseteq \overline{Q}(F)\backslash G(F)$ has a complement which is negligible. Hence, by (\ref{proposition 6.8.1(ii)}), if we choose Haar measures compatibly, we have

$$\displaystyle \int_{\overline{Q}(F)\backslash G(F)}\varphi(g) dg=\int_{H_{\overline{Q}}(F)\backslash H(F)} \varphi(h) dh$$

\noindent for all $\varphi\in L^1(\overline{Q}(F)\backslash G(F),\delta_{\overline{Q}})$. In particular we get

\begin{align}\label{eq 6.8.4}
\displaystyle \Xi^G(g)=\int_{H_{\overline{Q}}(F)\backslash H(F)} (e_K(hg),e_K(h))_\tau dh
\end{align}

\noindent for all $g\in G(F)$.

\vspace{2mm}

By Lemma \ref{lemma 6.5.1}(ii), there exists $d>0$ such that the integral

$$\displaystyle \int_{H(F)} \Xi^G(h)\sigma(h)^{-d}dh$$

\noindent converges. Choose such a $d>0$. Then by \ref{eq 6.8.4}, we have

$$\displaystyle \int_{H(F)} \Xi^G(h)\sigma(h)^{-d}dh=\int_{H(F)}\int_{H_{\overline{Q}}(F)\backslash H(F)} (e_K(h'h),e_K(h'))_\tau dh' \sigma(h)^{-d}dh$$

\noindent Note that this last double integral is absolutely convergent. Indeed, since the integral is convergent in that order, it suffices to check that $(e_K(h'h),e_K(h'))_\tau$ is positive for all $h', h\in H(F)$. But by the Iwasawa decomposition and \ref{eq 6.8.3} the terms $(e_K(h'h),e_K(h'))_\tau$ are values of $\Xi^L$ hence positive. Switching the two integrals, making the variable change $h\mapsto h'^{-1}h$ and decomposing the integral over $H(F)$ as a double integral over $H_{\overline{Q}}(F)\backslash H(F)$ and $H_{\overline{Q}}(F)\simeq H_L(F)$, by (\ref{proposition 6.8.1(ii)}) we get that the expression

$$\displaystyle \int_{\left(H_{\overline{Q}}(F)\backslash H(F)\right)^2} \int_{H_{L}(F)} \left(\tau(h_{L})e_K(h),e_K(h')\right)_\tau \sigma(h'^{-1}h_{L}h)^{-d}\delta_{H_{L}}(h_{L})^{1/2}d_Lh_{L}dhdh'$$

\noindent is absolutely convergent. By Fubini, it follows that there exist $h,h'\in H(F)$ such that the inner integral

$$\displaystyle \int_{H_{L}(F)} \left(\tau(h_{L})e_K(h),e_K(h')\right)_\tau \sigma(h'^{-1}h_{L}h)^{-d}\delta_{H_{L}}(h_{L})^{1/2}d_Lh_{L}$$

\noindent is absolutely convergent. Fix such $h,h'\in H(F)$. By the Iwasawa decomposition we may write $h=luk$ and $h'=l'u'k'$ with $l,l'\in L(F)$, $u,u'\in U_{\overline{Q}}(F)$ and $k,k'\in K$. Then, by \ref{eq 6.8.3} we have $\left(\tau(h_{L})e_K(h),e_K(h')\right)_\tau=\delta_{\overline{Q}}(l'l)^{1/2}\Xi^L(l'^{-1}h_{L}l)$ for all $h_{L}\in H_{L}(F)$. Since $ \Xi^L(h_{L}) \ll \Xi^L(l'^{-1}h_{L}l)$ and $\sigma(h'^{-1}h_{L}h)\ll \sigma(h_{L})$ for all $h_{L}\in H_{L}(F)$, it follows that the integral

$$\displaystyle \int_{H_{\overline{Q}}(F)} \Xi^L(h_{\overline{Q}})\sigma(h_{\overline{Q}})^{-d} \delta_{H_{\overline{Q}}}(h_{\overline{Q}})^{1/2}d_L h_{\overline{Q}}$$

\noindent is also absolutely convergent. This proves the first part of (\ref{proposition 6.8.1(iii)}). The second part follows from the same arguments using Lemma \ref{lemma 6.5.1}(i) instead of Lemma \ref{lemma 6.5.1}(ii).

\item This follows from (\ref{proposition 6.8.1(ii)}), (\ref{proposition 6.8.1(iii)}) and Proposition \ref{proposition 1.5.1}(iv).

\item Every $h^Q\in H^Q(F)$ can be written $h^Q=h_Lu_Q$ where $h_L\in H_L(F)\subset L(F)$ and $u_Q\in U_Q(F)$. Moreover, we have $\sigma(lu_Q)\sim \sigma(l)+\sigma(u_Q)$ and $\sigma(u_Q)\ll\sigma(a^{-1}u_Qa)$ for all $l\in L(F)$, $u_Q\in U_Q(F)$ and $a\in A_{\mini}^+$. Thus, it suffices to show that

\begin{align}\label{eq 6.8.a}
\displaystyle \sigma(h_L)\ll \sigma(a^{-1}h_La)
\end{align}

\noindent for all $h_L\in H_L(F)$ and $a\in A_{\mini}^+$.

Let $l:\overline{Q}\to L$ be the unique regular map such that $l(\overline{q})^{-1}\overline{q}\in U_{\overline{Q}}$ for all $\overline{q}\in \overline{Q}$. By (\ref{proposition 6.8.1(i)}) the map $h_{\overline{Q}}\mapsto l(h_{\overline{Q}})$ induces an isomorphism $H_{\overline{Q}}\simeq H_L$.In particular, we have
\begin{align}\label{eq 6.8.b}
\sigma(h_{\overline{Q}})\sim \sigma(l(h_{\overline{Q}}))
\end{align}
for all $h_{\overline{Q}}\in H_{\overline{Q}}(F)$. Let $A_L$ be the maximal split central torus of $L$ and set $A_L^+=A_L\cap A_{\mini}^+$. By definition of $A_{\mini}^+$, it is not hard to see that there exists a map $\overline{q}\in \overline{Q}(F)\mapsto a_L(\overline{q})\in A_L^+$ such that 
\begin{align}\label{eq 6.8.c}
\sigma(a_L(\overline{q})^{-1}a^{-1}\overline{q}aa_L(\overline{q}))\ll \sigma(a^{-1}l(\overline{q})a)
\end{align}
for all $\overline{q}\in \overline{Q}$ and all $a\in A_{\mini}^+$. Moreover, since $H_{\overline{Q}}\subset H$, by Proposition \ref{proposition 6.4.1}(iii), we have
\begin{align}\label{eq 6.8.d}
\sigma(h_{\overline{Q}})\ll \sigma(a^{-1}h_{\overline{Q}}a)
\end{align}
for all $h_{\overline{Q}}\in H_{\overline{Q}}(F)$ and all $a\in A_{\mini}^+$. From \ref{eq 6.8.b}, \ref{eq 6.8.c} and \ref{eq 6.8.d}, it follows that
$$\displaystyle \sigma(h_L)\sim \sigma(h_{\overline{Q}})\ll \sigma(a_L(h_{\overline{Q}})^{-1}a^{-1}h_{\overline{Q}}aa_L(h_{\overline{Q}}))\ll \sigma(a^{-1}h_La)$$
for all $h_{\overline{Q}}\in H_{\overline{Q}}(F)$ and all $a\in A_{\mini}^+$ where we have set $h_L=l(h_{\overline{Q}})\in H_L(F)$. This shows \ref{eq 6.8.a} and ends the proof of (\ref{proposition 6.8.1(v)}).

\item By (\ref{proposition 6.8.1(v)}), Proposition \ref{proposition 6.7.1}(ii) and Proposition \ref{proposition 1.5.1}(i), it suffices to show the existence of $d>0$, such that

\begin{align}\label{eq 6.8.5}
\displaystyle \int_{H^Q(F)} \Xi^{G}(a^{-1}h^Qa) \sigma(h^Q)^{-d} dh^Q\ll \Xi^{G}(a)^2
\end{align}

\noindent for all $a\in A_{\mini}^+$. As $\overline{P}_{\mini}$ is a good parabolic subgroup, it easily follows that $\overline{P}_{\mini}H^Q$ is Zariski open in $G$. Using this, and the doubling principle (Proposition \ref{proposition 1.5.1}(vi)), we show as in the proof of Proposition \ref{proposition 6.7.1}(v) that 

\[\begin{aligned}
\displaystyle \int_{H^Q(F)} \Xi^{G}(a^{-1}h^Qa) \sigma(h^Q)^{-d} dh^Q \ll \Xi^{G}(a)^2 \int_{H^Q(F)} \Xi^{G}(h^Q) \sigma(h^Q)^{-d} dh^Q
\end{aligned}\]
for all $a\in A_{\mini}^+$. This proves \ref{eq 6.8.5} since the right hand side is absolutely convergent for $d$ sufficiently large by (\ref{proposition 6.8.1(iv)}). $\blacksquare$
\end{enumerate}

\section{Explicit tempered intertwinings}\label{section 8}

We keep the notation of the previous chapter. Given a tempered representation $\pi$ of $G(F)$, the present chapter is devoted to the study of a certain explicit $(H,\xi)\times (H,\xi)$-equivariant sesquilinear form $\mathcal{L}_\pi$ on (the space of) $\pi$, the main result being that $\mathcal{L}_\pi$ is nonzero if and only if the multiplicity $m(\pi)$ is nonzero (Theorem \ref{theorem 8.2.1}). This will be a crucial ingredient in the proof of the spectral side of our local trace formula (Theorem \ref{theorem 9.1.1}). The sesquilinear form $\mathcal{L}_\pi$ is introduced in \ref{section 8.2}. It is essentially defined by integrating matrix coefficients of $\pi$ against the character $\xi$ of $H(F)$. Unfortunately, this integral does not converge for all tempered representations unless we are in the codimension one case (i.e. when $\xi=1$) but the oscillatory nature of the integral allows to regularize it in some canonical way. This is the content of Section \ref{section 8.1} (a similar regularization has actually been constructed, in greater generality, by Sakellaridis-Venkatesh, see \cite{SV} Corollary 6.3.3). In Section \ref{section 8.3}, we prove some a priori estimates for $(H,\xi)$-equivariant linear forms on $\pi$. Then, we discuss a certain relation between the (non-vanishing of) sesquilinear form $\mathcal{L}_\pi$ and parabolic induction in Section \ref{section 8.4}. The proof of the main theorem is given in Section \ref{section 8.5} and in Section \ref{section 8.6} we draw some consequences of this result.

\subsection{The $\xi$-integral}\label{section 8.1}

\noindent For all $f\in \mathcal{C}(G(F))$, the integral

$$\displaystyle \int_{H(F)} f(h)\xi(h)dh$$

\noindent is absolutely convergent by Lemma \ref{lemma 6.5.1}(ii). Moreover, by Lemma \ref{lemma 6.5.1}(ii) again, this defines a continuous linear form on $\mathcal{C}(G(F))$. Recall that $\mathcal{C}(G(F))$ is a dense subspace of the weak Harish-Chandra Schwartz space $\mathcal{C}^w(G(F))$ (by \ref{eq 1.5.1}).

\begin{prop}\label{proposition 8.1.1}
The linear form

$$\displaystyle f\in\mathcal{C}(G(F))\mapsto \int_{H(F)} f(h)\xi(h)dh$$

\noindent extends continuously to $\mathcal{C}^w(G(F))$.
\end{prop}

\vspace{2mm}

\noindent\ul{Proof}: Let us fix a one-parameter subgroup $a:\mathbb{G}_m\to A$ such that $\lambda(a(t)ha(t)^{-1})=t\lambda(h)$ for all $t\in\mathbb{G}_m$ and all $h\in H$ (Recall that $\lambda:H\to \mathbb{G}_a$ is the algebraic character such that $\xi=\psi\circ \lambda_F$), such a one-parameter subgroup is easy to construct. We shall now divide the proof according to whether $F$ is $p$-adic or real.

\vspace{2mm}

\begin{itemize}
\renewcommand{\labelitemi}{$\bullet$}

\item If $F$ is a $p$-adic field, then we may fix a compact-open subgroup $K\subseteq G(F)$ and prove that the linear form

$$\displaystyle f\in\mathcal{C}_K(G(F))\mapsto \int_{H(F)} f(h)\xi(h)dh$$

\noindent extends continuously to $\mathcal{C}_K^w(G(F))$. Set $K_a=a^{-1}\left(K\cap A(F)\right)\subseteq F^\times$. Then for all $f\in \mathcal{C}_K(G(F))$, we have

\begin{align}\label{eq 8.1.1}
\displaystyle \int_{H(F)} f(h)\xi(h)dh & =meas(K_a)^{-1}\int_{K_a}\int_{H(F)} f(a(t)^{-1}ha(t))\xi(h)dhd^\times t \\
\nonumber & =meas(K_a)^{-1}\int_{H(F)}f(h)\int_{K_a} \xi(a(t)ha(t)^{-1})d^\times tdh \\
\nonumber & =meas(K_a)^{-1}\int_{H(F)} f(h)\int_{K_a} \psi\left(t\lambda(h)\right)\lvert t\rvert^{-1}dtdh
\end{align}

\noindent The function $x\in F\mapsto \int_{K_a}\psi(tx)\lvert t\rvert^{-1}dt$ is the Fourier transform of the function $\lvert .\rvert^{-1}\mathbf{1}_{K_a}\in C_c^\infty(F)$ hence it belongs to $C_c^\infty(F)$. Now by Lemma \ref{lemma 6.5.1}(iii), the last integral of \ref{eq 8.1.1} is absolutely convergent for all $f\in\mathcal{C}_K^w(G(F))$ and defines a continuous linear form on that space. This is the extension we were looking for.

\item Now assume that $F=\mathbb{R}$. Let us denote by $\Ad$ the adjoint action of $G(F)$ on $\mathcal{C}^w(G(F))$ i.e., one has

$$\displaystyle \left(\Ad(g)f\right)(x)=f(g^{-1}xg),\;\; f\in \mathcal{C}^w(G(F)),\;\; g,x\in G(F)$$

\noindent Set $\Ad_a=\Ad\circ a$. Then $\Ad_a$ is a smooth representation of $F^\times$ on $\mathcal{C}^w(G(F))$ and hence induces an action, also denoted by $\Ad_a$, of $\mathcal{U}\left(\Lie(F^\times)\right)$ on $\mathcal{C}^w(G(F))$. Set $\Delta=1-\left(t\frac{d }{dt}\right)^2\in \mathcal{U}\left(\Lie(F^\times)\right)$. By elliptic regularity (\ref{eq 2.1.2}), for every integer $m\geqslant 1$, there exist functions $\varphi_1\in C_c^{2m-2}(F^\times)$ and $\varphi_2\in C_c^\infty(F^\times)$ such that

$$\displaystyle \varphi_1\ast \Delta^m+\varphi_2=\delta_1$$

\noindent Hence, we have the equality

$$\displaystyle \Ad_a(\varphi_1)\Ad_a(\Delta^m)+\Ad_a(\varphi_2)=Id$$

\noindent It follows that for all $f\in \mathcal{C}(G(F))$, we have

\begin{align}\label{eq 8.1.2}
\displaystyle \int_{H(F)} f(h) \xi(h)dh & =\int_{H(F)} \left(\Ad_a(\varphi_1)\Ad_a(\Delta^m)f\right) \xi(h)dh+ \int_{H(F)} \left(\Ad_a(\varphi_2)f\right)(h) \xi(h) dh \\
\nonumber & =\int_{H(F)} \left(\Ad_a(\Delta^m)f\right)(h) \int_{F^\times} \varphi_1(t)\xi(a(t)ha(t)^{-1}) \delta_P(a(t))d^\times tdh \\
\nonumber & +\int_{H(F)} f(h) \int_{F^\times} \varphi_2(t)\xi(a(t)ha(t)^{-1}) \delta_P(a(t))d^\times tdh \\
\nonumber & =\int_{H(F)} \left(\Ad_a(\Delta^m)f\right)(h) \int_{F} \varphi_1(t) \delta_P(a(t))\lvert t\rvert^{-1} \psi(t\lambda(h))dtdh \\
\nonumber & +\int_{H(F)} f(h) \int_{F} \varphi_2(t)\delta_P(a(t))\lvert t\rvert^{-1} \psi(t\lambda(h))dtdh
\end{align}

\noindent Consider the functions $f_i:x\in F\mapsto \int_{F} \varphi_i(t)\delta_P(a(t))\lvert t\rvert^{-1} \psi(tx)dt$, $i=1,2$. These are the Fourier transforms of the functions $t\mapsto \varphi_i(t)\delta_H(a(t))\lvert t\rvert^{-1}$, $i=1,2$, which both belong to $C_c^{2m-2}(F)$. Hence, $f_1$ and $f_2$ are both essentially bounded by $\left(1+\lvert x\rvert\right)^{-2m+2}$. Now, by Lemma \ref{lemma 6.5.1}(iii), if $m\geqslant 2$ the two integrals in the last term of \ref{eq 8.1.2} are absolutely convergent for all $f\in\mathcal{C}^w(G(F))$ and define on that space continuous linear forms. The extension we were looking for is just the sum of these two integrals. $\blacksquare$
\end{itemize}

\vspace{5mm}

\noindent The continuous linear form on $\mathcal{C}^w(G(F))$ whose existence is proved by the proposition above will be called the {\em $\xi$-integral} on $H(F)$ and will be denoted by

$$\displaystyle f\in\mathcal{C}^w(G(F))\mapsto \int_{H(F)}^* f(h)\xi(h)dh$$

\noindent or

$$\displaystyle f\in\mathcal{C}^w(G(F))\mapsto \gls{PHxi}(f)$$

\vspace{2mm}

\noindent We now note the following properties of the $\xi$-integral:

\begin{lem}\label{lemma 8.1.1}

\begin{enumerate}[(i)]
\item For all $f\in\mathcal{C}^w(G(F))$ and all $h_0,h_1\in H(F)$, we have

$$\displaystyle \mathcal{P}_{H,\xi}(L(h_0)R(h_1)f)=\xi(h_0)\xi(h_1)^{-1}\mathcal{P}_{H,\xi}(f)$$

\item Let $a:\mathbb{G}_m\to A$ be a one-parameter subgroup such that $\lambda(a(t)ha(t)^{-1})=t\lambda(h)$ for all $t\in \mathbb{G}_m$ and all $h\in H$. Denote by $\Ad_a$ the representation of $F^\times$ on $\mathcal{C}^w(G(F))$ given by $\Ad_a(t)=L(a(t))R(a(t))$ for all $t\in F^\times$. Let $\varphi\in C_c^\infty(F^\times)$. Set $\varphi'(t)=\lvert t\rvert^{-1}\delta_P(a(t))\varphi(t)$ for all $t\in F^\times$ and denote by $\widehat{\varphi'}$ its Fourier transform, that is

$$\displaystyle \widehat{\varphi'}(x)=\int_F \varphi'(t) \psi(tx)dt,\;\;\; x\in F$$

\noindent Then, we have

$$\displaystyle \mathcal{P}_{H,\xi} (\Ad_a(\varphi)f)=\int_{H(F)} f(h) \widehat{\varphi'}(\lambda(h)) dh$$

\noindent for all $f\in\mathcal{C}^w(G(F))$, where the second integral is absolutely convergent.
\end{enumerate}
\end{lem}

\vspace{2mm}

\noindent\ul{Proof}: In both (i) and (ii), both sides of the equality to be proved are continuous in $f\in\mathcal{C}^w(G(F))$ (for (ii) this follows from Lemma \ref{lemma 6.5.1}(iii)). Hence it is sufficient to check the relations for $f\in\mathcal{C}(G(F))$ where by straightforward variable changes we can pass from the left hand side to the right hand side. $\blacksquare$

\subsection{Definition of $\mathcal{L}_\pi$}\label{section 8.2}

\noindent Let $\pi$ be a tempered representation of $G(F)$. For all $T\in \End(\pi)^\infty$, the function

$$g\in G(F)\mapsto \Tr\left(\pi(g^{-1})T\right)$$

\noindent belongs to the weak Harish-Chandra Schwartz space $\mathcal{C}^w(G(F))$ by \ref{eq 2.2.4}. We can thus define a linear form $\gls{Lpi}\colon \End(\pi)^\infty\to \mathbb{C}$ by setting

$$\displaystyle \mathcal{L}_{\pi}(T)=\int_{H(F)}^* \Tr\left(\pi(h^{-1})T\right) \xi(h) dh,\;\;\; T\in \End(\pi)^\infty$$

\noindent By Lemma \ref{lemma 8.1.1}(i), we have

$$\displaystyle \mathcal{L}_{\pi}(\pi(h)T\pi(h'))=\xi(h)\xi(h')\mathcal{L}_\pi(T)$$

\noindent for all $h,h'\in H(F)$ and $T\in \End(\pi)^\infty$. By \ref{eq 2.2.5}, the map which associates to $T\in \End(\pi)^\infty$ the function

$$\displaystyle g\mapsto \Tr\left(\pi(g^{-1})T\right)$$

\noindent in $\mathcal{C}^w(G(F))$ is continuous. Since the $\xi$-integral is a continuous linear form on $\mathcal{C}^w(G(F))$, it follows that the linear form $\mathcal{L}_{\pi}$ is continuous. \\

\noindent Recall that we have a continuous $G(F)\times G(F)$-equivariant embedding with dense image $\pi^\infty\otimes \overline{\pi^\infty}\hookrightarrow \End(\pi)^\infty$, $e\otimes e'\mapsto T_{e,e'}$ (which is an isomorphism in the $p$-adic case). In any case, $\End(\pi)^\infty$ is naturally isomorphic to the completed projective tensor product $\pi^\infty\widehat{\otimes}_p \overline{\pi^\infty}$. Thus we may identify $\mathcal{L}_\pi$ with the continuous sesquilinear form on $\pi^\infty$ given by

$$\mathcal{L}_\pi(e,e'):=\mathcal{L}_\pi(T_{e,e'})$$

\noindent for all $e,e'\in \pi^\infty$. Expanding definitions, we have

$$\displaystyle \mathcal{L}_{\pi}(e,e')=\int_{H(F)}^* (e,\pi(h)e')\xi(h) dh$$

\noindent for all $e,e'\in \pi^\infty$. Fixing $e'\in \pi^\infty$, we see that the map $e\in\pi^\infty\mapsto \mathcal{L}_\pi(e,e')$ belongs to $\Hom_H(\pi^\infty,\xi)$. By the density of $\pi^\infty\otimes \overline{\pi}^\infty$ in $\End(\pi)^\infty$, it follows that

$$\mathcal{L}_\pi\neq 0\Rightarrow m(\pi)\neq 0$$

\noindent The purpose of this chapter is to prove the converse direction. Namely, we will show

\vspace{5mm}

\begin{theo}\label{theorem 8.2.1}
For all $\pi\in \Temp(G)$, we have

$$\mathcal{L}_\pi\neq 0\Leftrightarrow m(\pi)\neq 0$$
\end{theo}

\vspace{5mm}

\noindent As said in the introduction, this result has already been proved in \cite{Beu1} (Theorem 14.3.1) in the $p$-adic case following closely the proof of the analogous result for special orthogonal groups given by Y. Sakellaridis and A. Venkatesh (Theorem 6.4.1 of \cite{SV}). The proof, that goes through the 3 following sections, is closer to the original treatment of Waldspurger (Proposition 5.7 of \cite{Wa4}). 

\vspace{3mm}

\noindent To end this section, we will content ourself with giving some of the basic properties of $\mathcal{L}_\pi$. First, since $\mathcal{L}_\pi$ is a continuous sesquilinear form on $\pi^\infty$, it defines a continuous linear map

$$\displaystyle \gls{Lpid}: \pi^\infty\to \overline{\pi^{-\infty}}$$
$$\displaystyle e\mapsto \mathcal{L}_\pi(e,.)$$

\noindent where we recall that $\overline{\pi^{-\infty}}$ denotes the topological conjugate-dual of $\pi^\infty$ endowed with the strong topology. This operator $L_\pi$ has its image included in $\overline{\pi^{-\infty}}^{H,\xi}=\Hom_H(\overline{\pi^{\infty}},\xi)$. By Theorem \ref{theorem 6.3.1}, this subspace is finite-dimensional and even of dimension less or equal to $1$ if $\pi$ is irreducible. Let $T\in \End(\pi)^\infty$. Recall that it extends uniquely to a continuous operator $T:\overline{\pi^{-\infty}}\to \pi^\infty$. Thus, we may form the two compositions

$$TL_\pi:\pi^\infty\to \pi^\infty$$

$$L_\pi T:\overline{\pi^{-\infty}}\to \overline{\pi^{-\infty}}$$

\noindent which are both finite-rank operators. In particular, their traces are well-defined and we have, almost by definition,

\begin{align}\label{eq 8.2.1}
\Tr(TL_\pi)=\Tr(L_\pi T)=\mathcal{L}_\pi(T)
\end{align}

\begin{lem}\label{lemma 8.2.1}
We have the following

\begin{enumerate}[(i)]

\item The maps

$$\pi\in \mathcal{X}_{\tempe}(G)\mapsto L_\pi\in \Hom(\pi^\infty,\overline{\pi^{-\infty}})$$

$$\pi\in \mathcal{X}_{\tempe}(G)\mapsto \mathcal{L}_\pi\in \End(\pi)^{-\infty}$$

\noindent are smooth in the following sense: For every parabolic subgroup $Q=LU_Q$ of $G$, for all $\sigma\in \Pi_2(L)$ and for every maximal compact subgroup $K$ of $G(F)$, which is special in the $p$-adic case, the maps

$$\lambda\in i\mathcal{A}_L^*\mapsto \mathcal{L}_{\pi_\lambda}\in \End(\pi_\lambda)^{-\infty}\simeq \End(\pi_K)^{-\infty}$$

$$\lambda\in i\mathcal{A}_L^*\mapsto L_{\pi_\lambda}\in \Hom(\pi_\lambda^\infty,\overline{\pi_\lambda^{-\infty}})\simeq \Hom(\pi_K^\infty,\overline{\pi_K^{-\infty}})$$

\noindent are smooth, where we have set $\pi_\lambda=i_Q^G(\sigma_\lambda)$ and $\pi_K=i_{Q\cap K}^K(\sigma)$.

\item Let $\pi$ be in $\Temp(G)$ or $\mathcal{X}_{\tempe}(G)$. Then for all $S,T\in \End(\pi)^\infty$, we have $SL_\pi\in \End(\pi)^\infty$ and

$$\mathcal{L}_\pi(S)\mathcal{L}_\pi(T)=\mathcal{L}_\pi(SL_\pi T)$$

\item Let $S,T\in \mathcal{C}(\mathcal{X}_{\tempe}(G),\mathcal{E}(G))$. Then, the section $\pi\in \Temp(G)\mapsto S_\pi L_\pi T_\pi\in \End(\pi)^\infty$ belongs to $C^\infty(\mathcal{X}_{\tempe}(G),\mathcal{E}(G))$.

\item Let $f\in \mathcal{C}(G(F))$ and assume that its Fourier transform $\pi\in \mathcal{X}_{\tempe}(G)\mapsto \pi(f)$ is compactly supported (this condition is automatically satisfied when $F$ is $p$-adic). Then, we have the equality

$$\displaystyle \int_{H(F)} f(h)\xi(h)dh=\int_{\mathcal{X}_{\tempe}(G)} \mathcal{L}_\pi(\pi(f)) \mu(\pi) d\pi$$

\noindent both integrals being absolutely convergent.

\item Let $f,f'\in \mathcal{C}(G(F))$ and assume that the Fourier transform of $f$ is compactly supported. Then we have the equality

$$\displaystyle \int_{\mathcal{X}_{\tempe}(G)} \mathcal{L}_\pi(\pi(f))\overline{\mathcal{L}_{\pi}(\pi(\overline{f'}))} \mu(\pi) d\pi=\int_{H(F)} \int_{H(F)} \int_{G(F)} f(hgh')f'(g)dg\xi(h')dh'\xi(h) dh$$

\noindent where the first integral is absolutely convergent and the second integral is convergent in that order but not necessarily as a triple integral.

\end{enumerate}
\end{lem}

\vspace{3mm}

\noindent\ul{Proof}:

\begin{enumerate}[(i)]

\item Let $Q=LU_Q$, $\sigma\in \Pi_2(L)$ and $K$ be as in the statement. Recall that our convention is to equip all the spaces that appear in the statement with the strong topology.

\vspace{2mm}

\noindent We have $\End(\pi_K)^\infty\simeq \pi_K^\infty\widehat{\otimes}_p \overline{\pi_K^\infty}$. Hence, the space $\End(\pi_K)^{-\infty}$ may be identified with the space of continuous sesquilinear forms on $\pi_K^\infty$ and we get a natural continuous embedding

$$\End(\pi_K)^{-\infty}\hookrightarrow \Hom(\pi_K^\infty,\overline{\pi_K^{-\infty}})$$

\noindent Of course, the image of $\mathcal{L}_{\pi_\lambda}$ by this map is $L_{\pi_\lambda}$, for all $\lambda\in i\mathcal{A}_L^*$. Consequently, it suffices to prove the smoothness of the map $\lambda\mapsto \mathcal{L}_{\pi_\lambda}$. By Proposition \ref{proposition A.3.1}(iv), this is equivalent to the smoothness of

$$\lambda\mapsto \mathcal{L}_{\pi_\lambda}(T)$$

\noindent for all $T\in \End(\pi_K)^\infty$. Because the $\xi$-integral is a continuous linear form on $\mathcal{C}^w(G(F))$, the smoothness of this last map follows from Lemma \ref{lemma 2.3.1}(ii).

\item The two inclusions $\End(\pi)^\infty\subset \Hom(\overline{\pi^{-\infty}},\pi)$ and $\End(\pi)^\infty\subset \Hom(\pi,\pi^\infty)$ are continuous. It follows that the bilinear map

$$\End(\pi)^\infty\times \End(\pi)^\infty\to \End(\pi)$$
$$(S,T)\mapsto SL_\pi T$$

\noindent is separately continuous. For all $S,T\in \End(\pi)^\infty$, the maps $g\in G(F)\mapsto \pi(g)S\in \End(\pi)^\infty$ and $g\in G(F)\mapsto T\pi(g)\in \End(\pi)^\infty$ are smooth. Hence, by Proposition \ref{proposition A.3.1}(v) in the real case, we have $SL_\pi T\in \End(\pi)^\infty$ for all $S,T\in \End(\pi)^\infty$. We now prove the equality

$$\mathcal{L}_\pi(S)\mathcal{L}_\pi(T)=\mathcal{L}_\pi(SL_\pi T)$$

\noindent for all $S,T\in \End(\pi)^\infty$. Assume first that $\pi\in \Temp(G)$. Then, this follows directly from \ref{eq 8.2.1} since the operators $L_\pi S,L_\pi T:\overline{\pi^{-\infty}}\to \overline{\pi^{-\infty}}$ have their images contained in the same line (which is $\Hom_H(\overline{\pi}^\infty,\xi)$). Assume now that $\pi\in \mathcal{X}_{\tempe}(G)$. We may then find a parabolic subgroup $Q=LU_Q$ of $G$ and a square-integrable representation $\sigma\in \Pi_2(L)$ such that $\pi=i_P^G(\sigma)$. Let $K$ be a maximal compact subgroup of $G(F)$ which is special in the $p$-adic case and set $\pi_K=i_{Q\cap K}^K(\sigma)$ and $\pi_\lambda=i_P^G(\sigma_\lambda)$ for all $\lambda\in i\mathcal{A}_L^*$. Then, we have isomorphisms $\End(\pi_\lambda)^\infty\simeq \End(\pi_K)^\infty$ for all $\lambda\in i\mathcal{A}_L^*$. Let $S,T\in \End(\pi)^\infty$ and identify them to their images in $\End(\pi_K)^\infty$ by the previous isomorphism. For $\lambda$ in a dense subset of $i\mathcal{A}_L^*$, the representation $\pi_\lambda$ is irreducible. Hence, by what we just saw, for every such $\lambda\in i\mathcal{A}_L^*$ we have

$$\mathcal{L}_{\pi_\lambda}(S)\mathcal{L}_{\pi_\lambda}(T)=\mathcal{L}_{\pi_\lambda}(SL_{\pi_\lambda} T)$$

\noindent By (i), the left hand side of the above equality is continuous in $\lambda\in i\mathcal{A}_L^*$. To deduce the equality at $\lambda=0$ (what we want), it thus suffices to show that the function

$$\lambda\in i\mathcal{A}_L^*\mapsto \mathcal{L}_{\pi_\lambda}(SL_{\pi_\lambda} T)$$

\noindent is continuous. We are even going to prove that it is a smooth function. By (i) and Proposition \ref{proposition A.3.1}(v), it suffices to show that for all $\lambda\in i\mathcal{A}_L^*$, the trilinear map

\begin{align}\label{eq 8.2.2}
\End(\pi_\lambda)^\infty\times \Hom(\pi_\lambda^\infty,\overline{\pi_{\lambda}^{-\infty}})\times \End(\pi_\lambda)^\infty\to \End(\pi_\lambda)^\infty
\end{align}
$$(S,L,T)\mapsto SLT$$

\noindent is separately continuous. As the inclusions $\End(\pi_\lambda)^\infty\subset \Hom(\overline{\pi_\lambda^{-\infty}},\pi_\lambda)$ and $\End(\pi_\lambda)^\infty\subset \Hom(\pi_\lambda,\pi_\lambda^\infty)$ are continuous, the trilinear map

$$\End(\pi_\lambda)^\infty\times \Hom(\pi_\lambda^\infty,\overline{\pi_{\lambda}^{-\infty}})\times \End(\pi_\lambda)^\infty\to \End(\pi_\lambda)$$
$$(S,L,T)\mapsto SLT$$

\noindent is separately continuous fo all $\lambda\in i\mathcal{A}_L^*$. By definition of the topology on $\End(\pi_\lambda)^\infty$, this immediately implies that \ref{eq 8.2.2} is separately continuous for all $\lambda\in i\mathcal{A}_L^*$. This ends the proof of (ii)

\item This is also a direct consequence of (i) and of the fact that the trilinear map \ref{eq 8.2.2} is separately continuous.

\item Let $f\in \mathcal{C}(G(F))$. The left hand side of $(iv)$ is absolutely convergent by Lemma \ref{lemma 6.5.1}(ii). By Lemma \ref{lemma 2.3.1}(ii), the map

$$\pi\in \mathcal{X}_{\tempe}(G)\mapsto \varphi(f,\pi)\in \mathcal{C}^w(G(F))$$

\noindent where $\varphi(f,\pi)(g)=\Tr(\pi(g^{-1})\pi(f))$, is continuous. By the hypothesis made on $f$, this map is also compactly supported. It follows that the function $\pi\in \mathcal{X}_{\tempe}(G)\mapsto \mu(\pi)\varphi(f,\pi)\in \mathcal{C}^w(G(F))$ is absolutely integrable. Hence, the function

$$\pi\in \mathcal{X}_{\tempe}(G)\mapsto \mu(\pi)\mathcal{L}_\pi(\pi(f))=\mu(\pi)\mathcal{P}_{H,\xi}(\varphi(f,\pi))$$

\noindent where $\mathcal{P}_{H,\xi}:\mathcal{C}^w(G(F))\to \mathbb{C}$ denotes the $\xi$-integral, is also absolutely integrable, proving the convergence of the right hand side of $(iv)$. We also have the equality

$$\displaystyle f=\int_{\mathcal{X}_{\tempe}(G)}\varphi(f,\pi)\mu(\pi)d\pi$$

\noindent in $\mathcal{C}^w(G(F))$ (or its completion),indeed by the Harish-Chandra Plancherel formula both sides are equal after applying the evaluation map at $g$ for all $g\in G(F)$. It follows that

$$\displaystyle \mathcal{P}_{H,\xi}(f)=\int_{\mathcal{X}_{\tempe}(G)} \mathcal{P}_{H,\xi}(\varphi(\pi,f))\mu(\pi)d\pi$$

\noindent which is exactly the content of $(iv)$.

\item The right hand side of (v) may be rewritten as

\begin{align}\label{eq 8.2.3}
\displaystyle \int_{H(F)} \int_{H(F)} ({f'}^\vee\ast L(h^{-1})f)(h') \xi(h')dh'\xi(h)dh
\end{align}

\noindent where ${f'}^\vee(g)=f'(g^{-1})$. The Fourier transform of ${f'}^\vee\ast L(h^{-1})f$ is given by

$$\pi\in \mathcal{X}_{\tempe}(G)\mapsto \pi({f'}^\vee\ast L(h^{-1})f)=\pi({f'}^\vee)\pi(h^{-1})\pi(f)$$

\noindent In particular, it is compactly supported. Applying (iv) to ${f'}^\vee\ast L(h^{-1})f$, we deduce that the integral

$$\displaystyle \int_{H(F)} ({f'}^\vee\ast L(h^{-1})f)(h') \xi(h')dh'$$

\noindent is absolutely convergent and is equal to

$$\displaystyle \int_{\mathcal{X}_{\tempe}(G)} \mathcal{L}_\pi(\pi({f'}^\vee)\pi(h^{-1})\pi(f)) \mu(\pi) d\pi$$

\noindent By \ref{eq 8.2.1}, this last integral is equal to 

$$\displaystyle \int_{\mathcal{X}_{\tempe}(G)} \Tr(\pi(h^{-1})\pi(f)L_\pi\pi({f'}^\vee)) \mu(\pi) d\pi$$

\noindent By (iii), the section $\pi\in \mathcal{X}_{\tempe}(G)\mapsto \pi(f)L_\pi\pi({f'}^\vee)\in \End(\pi)^\infty$ is smooth. Moreover, it is compactly supported and so it belongs to $\mathcal{C}(\mathcal{X}_{\tempe}(G),\mathcal{E}(G))$. By the matricial Paley-Wiener theorem (Theorem \ref{theorem 2.6.1}), it is thus the Fourier transform of a Harish-Chandra Schwartz function. Applying (iv) to this function, we see that the exterior integral of \ref{eq 8.2.3} is absolutely convergent and that the whole expression is equal to the absolutely convergent integral

$$\displaystyle \int_{\mathcal{X}_{\tempe}(G)}\mathcal{L}_{\pi}(\pi(f)L_\pi\pi({f'}^\vee)) \mu(\pi) d\pi$$

\noindent which by (ii) is equal to

$$\displaystyle \int_{\mathcal{X}_{\tempe}(G)}\mathcal{L}_{\pi}(\pi(f))\mathcal{L}_\pi(\pi(f'^\vee)) \mu(\pi) d\pi$$

\noindent The point (v) now follows from this and the easily checked equality

$$\mathcal{L}_\pi(\pi(f'^\vee))=\overline{\mathcal{L}_{\pi}(\pi(\overline{f'}))}, \;\;\; \pi\in\mathcal{X}_{\tempe}(G)$$

$\blacksquare$

\end{enumerate}

\subsection{Asymptotics of tempered intertwinings}\label{section 8.3}

\begin{lem}\label{lemma 8.3.1}
\begin{enumerate}[(i)]
\item Let $\pi$ be a tempered representation of $G(F)$ and $\ell\in \Hom_H(\pi^\infty,\xi)$ be a continuous $(H,\xi)$-equivariant linear form. Then, there exist $d>0$ and a continuous semi-norm $\nu_d$ on $\pi^\infty$ such that

$$\lvert \ell(\pi(x)e)\rvert\leqslant \nu_d(e)\Xi^{H\backslash G}(x)\sigma_{H\backslash G}(x)^d$$

\noindent for all $e\in \pi^\infty$ and all $x\in H(F)\backslash G(F)$.

\item For all $d>0$, there exists $d'>0$ and a continuous semi-norm $\nu_{d,d'}$ on $\mathcal{C}^w_d(G(F))$ such that

$$\left\lvert \mathcal{P}_{H,\xi}(R(x)L(y)\varphi)\right\rvert\leqslant \nu_{d,d'}(\varphi) \Xi^{H\backslash G}(x)\Xi^{H\backslash G}(y) \sigma_{H\backslash G}(x)^{d'} \sigma_{H\backslash G}(y)^{d'}$$

\noindent for all $\varphi\in \mathcal{C}_d^w(G(F))$ and all $x,y\in H(F)\backslash G(F)$.
\end{enumerate}
\end{lem}

\vspace{2mm}

\noindent\ul{Proof}: We will use the notation of Section \ref{section 6.6.2}. Namely, we have

\vspace{2mm}

\begin{itemize}
\renewcommand{\labelitemi}{$\bullet$}
\item $\overline{P}_0=M_0\overline{N}_0$ is a good parabolic subgroup of $G_0$, $A_0$ the split component of $M_0$;

\item $A_0^+=\{a\in A_0(F); \lvert \alpha(a)\rvert\geqslant 1 \;\forall \alpha\in R(A_0,\overline{P}_0)\}$;

\item $\overline{P}=M\overline{N}$ is the parabolic subgroup opposite to $P$ with respect to $M$;

\item $\overline{P}_{\mini}=\overline{P}_0T\overline{N}$, $M_{\mini}=M_0T$ and $A_{\mini}=A_0A$;

\item $P_{\mini}$ the parabolic subgroup opposite to $\overline{P}_{\mini}$ with respect to $M_{\mini}$;

\item $\Delta$ the set of simple roots of $A_{\mini}$ in $P_{\mini}$ and $\Delta_P=\Delta\cap R(A_{\mini},P)$.
\end{itemize}

\vspace{2mm}

\noindent Then, $\overline{P}_{\mini}$ is a good parabolic subgroup of $G$ and by Lemma \ref{lemma 6.6.2.1}(ii), there exists a compact subset $\mathcal{K}\subseteq G(F)$ such that

\begin{align}\label{eq 8.3.1}
G(F)=H(F)A_0^+A(F)\mathcal{K}
\end{align}

\begin{enumerate}[(i)]

\item By \ref{eq 8.3.1}, Proposition \ref{proposition 6.7.1}(i) and (ii) and since the family $\left(\pi(k)\right)_{k\in \mathcal{K}}$ is equicontinuous on $\pi^\infty$, it is sufficient to prove the following

\vspace{3mm}

\begin{num}
\item\label{eq 8.3.2} There exists a continuous semi-norm $\nu$ on $\pi^\infty$ such that
$$\lvert \ell(\pi(a)e)\rvert \leqslant \Xi^G(a)\nu(e)$$
for all $e\in \pi^\infty$ and all $a\in A_0^+A(F)$.
\end{num}

\vspace{3mm}

We divide the proof of \ref{eq 8.3.2} according to whether $F$ is $p$-adic or real.

\vspace{3mm}

\begin{itemize}
\renewcommand{\labelitemi}{$\bullet$}

\item Assume first that $F$ is $p$-adic. Since the topology on $\pi^\infty$ is the finest locally convex topology, we only need to show that for all $e\in \pi^\infty$, we have

\begin{align}\label{eq 8.3.3}
\lvert \ell(\pi(a)e)\rvert\ll \Xi^G(a)
\end{align}

\noindent for all $a\in A_0^+A(F)$. Let $e\in \pi^\infty$ and let $K$ be a compact-open subgroup of $G(F)$ such that $e\in (\pi^\infty)^{K}$. First, we show

\vspace{3mm}

\begin{num}
\item\label{eq 8.3.4} There exists $c=c_{K}\geqslant 1$ such that for all $a\in A_{\mini}(F)$, if there exists $\alpha\in \Delta_P$ such that $\lvert \alpha(a)\rvert\geqslant c$, then
$$\ell(\pi(a)e)=0$$
\end{num}

\vspace{3mm}

\noindent Since $\Delta_P$ is finite, it is sufficient to fix $\alpha\in \Delta_P$ and prove that for $a\in A_{\mini}(F)$, if $\lvert \alpha(a)\rvert$ is big enough then $\ell(\pi(a)e)=0$. So, let $\alpha\in \Delta_P$ and $a\in A_{\mini}(F)$. By Lemma \ref{lemma 6.6.2.1}(iii), there exists $X\in \mathfrak{n}_{\alpha}(F)$ (where $\mathfrak{n}_\alpha$ is the eigensubspace of $\mathfrak{n}$ corresponding to $\alpha$) such that $\xi(e^X)\neq 1$. Now, if $\lvert \alpha(a)\rvert$ is big enough, we will have $a^{-1}e^Xa\in K$ and so

$$\xi(e^X)\ell(\pi(a)e)=\ell(\pi(e^X)\pi(a)e)=\ell(\pi(a)e)$$

\noindent From which it follows that $\ell(\pi(a)e)=0$. This ends the proof of \ref{eq 8.3.4}.

\vspace{2mm}

\noindent Let $c\geqslant 1$ be as in \ref{eq 8.3.4} and set

$$A_{\mini}^+(c)=\{a\in A_{\mini}(F); \lvert \alpha(a)\rvert\leqslant c\; \forall \alpha\in \Delta\}$$

\noindent Note that we have $\lvert \alpha(a)\rvert\leqslant 1$ for all $a\in A_0^+A(F)$ and all $\alpha\in \Delta\backslash \Delta_P$. Hence, by \ref{eq 8.3.4}, it is sufficient to prove the inequality \ref{eq 8.3.3} for $a\in A_{\mini}^+(c)$. By definition of $A_{\mini}^+(c)$, there exists a compact-open subgroup $K'_{\overline{P}_{\mini}}$ of $\overline{P}_{\mini}(F)$ such that

$$K'_{\overline{P}_{\mini}}\subseteq aKa^{-1}$$

\noindent for all $a\in A_{\mini}^+(c)$. Also, let $K'_H$ be a compact-open subgroup of $H(F)$ on which $\xi$ is trivial. Since $\overline{P}_{\mini}$ is a good parabolic subgroup of $G$, we may find a compact-open subgroup $K'$ of $G(F)$ such that $K'\subseteq K'_HK'_{\overline{P}_{\mini}}$. Let $k'=k'_{H}k'_{\overline{P}_{\mini}}\in K'$, where $k'_H\in K'_H$ and $k'_{\overline{P}_{\mini}}\in K'_{\overline{P}_{\mini}}$, then we have

$$\ell(\pi(k')\pi(a)e)=\ell\left(\pi(k'_H)\pi(a)\pi(a^{-1}k'_{\overline{P}_{\mini}}a)e\right)=\xi(k'_H)\ell(\pi(a)e)=\ell(\pi(a)e)$$

\noindent for all $a\in A^+_{\mini}(c)$. It follows that

$$\ell(\pi(a)e)=\ell(\pi(e_{K'})\pi(a)e)$$

\noindent for all $a\in A_{\mini}^+(c)$. We have $\ell\circ \pi(e_{K'})\in \overline{\pi^\infty}$ and the inequality \ref{eq 8.3.3} for $a\in A_{\mini}^+(c)$ now follows from \ref{eq 2.2.3}.

\item Assume now that $F=\mathbb{R}$. Let us set for all $I\subseteq \Delta$,

$$A_{\mini}^+(I)=\{a\in A_{\mini}(F); \lvert \alpha(a)\rvert\leqslant 1 \; \forall \alpha\in \Delta\backslash I \mbox{ and } \lvert \alpha(a)\rvert>1\;\forall \alpha\in I\}$$

\noindent Since $\lvert \alpha(a)\rvert\leqslant 1$ for all $a\in A_0^+A(F)$ and all $\alpha\in \Delta\backslash \Delta_P$, we have

\begin{align}\label{eq 8.3.5}
\displaystyle A_0^+A(F)=\bigsqcup_{I\subseteq \Delta_P} A_{\mini}^+(I)
\end{align}

\noindent Hence, we may fix $I\subseteq \Delta_P$ and prove the inequality \ref{eq 8.3.2} restricted to $a\in A_{\mini}^+(I)$. Let $X_1,\ldots,X_p$ be a basis of $\overline{\mathfrak{p}}_{\mini}(F)$ and set

$$\Delta_{\mini}=1-\left(X_1^2+\ldots+X_p^2\right)\in \mathcal{U}(\overline{\mathfrak{p}}_{\mini})$$

\noindent Let $k\geqslant \dim(\overline{P}_{\mini})+1$ be an integer. We will need the following

\vspace{3mm}

\begin{num}
\item\label{eq 8.3.6} There exists $u=u_{I,k}\in \mathcal{U}(\mathfrak{n})$ such that the two maps
$$a\in A_{\mini}^+(I)\mapsto a^{-1}(\Delta_{\mini}^ku)a\in \mathcal{U}(\mathfrak{g})$$
$$a\in A_{\mini}^+(I)\mapsto a^{-1}ua\in \mathcal{U}(\mathfrak{g})$$
have bounded images and $d\xi(u)=1$.
\end{num}

\vspace{3mm}

\noindent This follows rather easily from Lemma \ref{lemma 6.6.2.1}(iii). We henceforth fix such an element $u=u_{I,k}\in\mathcal{U}(\mathfrak{n})$. By elliptic regularity \ref{eq 2.1.2}, there exist two functions $\varphi_1\in C_c^{k_1}(\overline{P}_{\mini}(F))$ and $\varphi_2\in C_c^\infty(\overline{P}_{\mini}(F))$, where $k_1=2k-\dim(\overline{P}_{\mini})-1$, such that

$$\pi(\varphi_1)\pi(\Delta_{\mini}^k)+\pi(\varphi_2)=Id_{\pi^\infty}$$

\noindent Choose a function $\varphi_H\in C_c^\infty(H(F))$ such that $\int_{H(F)}\varphi_H(h)\xi(h)dh=1$. Then, for all $e\in \pi^\infty$ and all $a\in A_{\mini}^+(I)$, we have

\[\begin{aligned}
\ell(\pi(a)e) & =d\xi(u)\ell(\pi(a)e)=\ell(\pi(u)\pi(a)e) \\
 & =\ell\left(\pi(\varphi_1)\pi(\Delta_{\mini}^ku)\pi(a)e\right)+\ell\left(\pi(\varphi_2)\pi(u)\pi(a)e\right) \\
 & =\ell\left(\pi(\varphi_1)\pi(a)\pi(a^{-1}(\Delta_{\mini}^ku)a)e\right) +\ell\left(\pi(\varphi_2)\pi(a)\pi(a^{-1}ua)e\right) \\
 & =\ell\left(\pi(\varphi_H\ast\varphi_1)\pi(a)\pi(a^{-1}(\Delta_{\mini}^ku)a)e\right) +\ell\left(\pi(\varphi_H\ast\varphi_2)\pi(a)\pi(a^{-1}ua)e\right)
\end{aligned}\]

\noindent Note that the functions $\varphi_H\ast\varphi_1$ and $\varphi_H\ast\varphi_2$ both belong to $C_c^{k_1}(G(F))$. The inequality \ref{eq 8.3.2} for $a\in A_{\mini}^+(I)$ is now a consequence of \ref{eq 8.3.6} and the equality above if we choose $k$ sufficiently large by the following fact 

\vspace{3mm}

\begin{num}
\item\label{eq 8.3.7} There exists an integer $k'_1\geqslant 1$ such that for all $\varphi\in C_c^{k'_1}(G(F))$, there exists a continuous semi-norm $\nu_\varphi$ on $\pi^\infty$ such that
$$\lvert\ell(\pi(\varphi)\pi(g)e)\rvert\leqslant\nu_\varphi(e)\Xi^G(g)$$
for all $e\in \pi^\infty$ and for all $g\in G(F)$.
\end{num}

\vspace{3mm}

\noindent (This is an easy consequence of the fact that $\ell$ is a continuous linear form and of \ref{eq 2.2.6}).
\end{itemize}

\item Let $d>0$. Again by \ref{eq 8.3.1} and Proposition \ref{proposition 6.7.1}(i) and (ii), we only need to prove the following

\vspace{3mm}

\begin{num}
\item\label{eq 8.3.8} There exists a continuous semi-norm $\nu_d$ on $\mathcal{C}^w_d(G(F))$ such that
$$\left\lvert \mathcal{P}_{H,\xi}(R(a_1)L(a_2)\varphi)\right\rvert\leqslant \nu_d(\varphi)\Xi^G(a_1)\Xi^G(a_2)\sigma(a_1)^d \sigma(a_2)^d$$
for all $\varphi\in\mathcal{C}_d^w(G(F))$ and all $a_1,a_2\in A_0^+A(F)$.
\end{num}

\vspace{3mm}

\noindent The proof of this fact is very close to the proof of \ref{eq 8.3.2}. We shall only sketch it, distinguishing again between the case where $F$ is $p$-adic and the case where $F$ is real.

\vspace{3mm}

\begin{itemize}
\renewcommand{\labelitemi}{$\bullet$}
\item Assume first that $F$ is $p$-adic. Then, we may fix a compact-open subgroup $K\subset G(F)$ and we are reduced to proving the following

\vspace{3mm}

\begin{num}
\item\label{eq 8.3.9} There exists a continuous semi-norm $\nu_{K,d}$ on $\mathcal{C}^w_{K,d}(G(F))$ such that
$$\left\lvert \mathcal{P}_{H,\xi}(R(a_1)L(a_2)\varphi)\right\rvert\leqslant \nu_{K,d}(\varphi) \Xi^G(a_1)\Xi^G(a_2)\sigma(a_1)^d\sigma(a_2)^d$$
for all $\varphi\in \mathcal{C}_{K,d}^w(G(F))$ and all $a_1,a_2\in A_0^+A(F)$.
\end{num}

\vspace{3mm}

\noindent We prove as in the proof of (i) that there exists a constant $c=c_K\geqslant 1$ such that

$$\mathcal{P}_{H,\xi}(R(a_1)L(a_2)\varphi)=0$$

\noindent for all $\varphi\in \mathcal{C}_{K,d}^w(G(F))$ as soon as $a_1\in A_0^+A(F)-A_{\mini}^+(c)$ or $a_2\in A_0^+A(F)-A_{\mini}^+(c)$. We also prove as in the proof of (i) that there exists a compact-open subgroup $K'\subset G(F)$ such that

$$\mathcal{P}_{H,\xi}(R(a_1)L(a_2)\varphi)=\mathcal{P}_{H,\xi}(R(e_{K'})L(e_{K'})R(a_1)L(a_2)\varphi)$$

\noindent for all $\varphi\in \mathcal{C}_{K,d}^w(G(F))$ and all $a_1,a_2\in A_{\mini}^+(c)$. The inequality \ref{eq 8.3.9} now follows from Lemma \ref{lemma 1.5.1}(i)(a).

\item Assume now that $F=\mathbb{R}$. Then, by \ref{eq 8.3.5}, we may fix $I,J\subseteq \Delta_P$ and just prove

\vspace{3mm}

\begin{num}
\item\label{eq 8.3.10} There exists a continuous semi-norm $\nu_{I,J,d}$ on $\mathcal{C}_d^w(G(F))$ such that
$$\left\lvert \mathcal{P}_{H,\xi}(R(a_1)L(a_2)\varphi)\right\rvert\leqslant \nu_{I,J,d}(\varphi) \Xi^G(a_1)\Xi^G(a_2)\sigma(a_1)^d\sigma(a_2)^d$$
for all $\varphi\in \mathcal{C}_d^w(G(F))$, all $a_1\in A_{\mini}^+(I)$ and all $a_2\in A_{\mini}^+(J)$.
\end{num}

\vspace{3mm}

\noindent Let $k\geqslant \dim(\overline{P}_{\mini})+1$ be an integer and choose elements $u_I=u_{I,k},u_J=u_{J,k}\in\mathcal{U}(\mathfrak{n})$ as in \ref{eq 8.3.6}. Then, as in the proof of (i), we show that there exists functions $\varphi_1,\varphi_2,\varphi_3,\varphi_4\in C_c^{k_1}(G(F))$, where $k_1=2k-\dim(\overline{P}_{\mini})-1$, such that

\[\begin{aligned}
\displaystyle \mathcal{P}_{H,\xi}(R(a_1)L(a_2)\varphi) & = \mathcal{P}_{H,\xi}\left[R(\varphi_1)L(\varphi_3)R(a_1)L(a_2)R(a_1^{-1}\Delta_{\mini}^ku_Ia_1)L(a_2^{-1}\Delta_{\mini}^ku_Ja_2)\varphi\right] \\
 & +\mathcal{P}_{H,\xi}\left[R(\varphi_1)L(\varphi_4)R(a_1)L(a_2)R(a_1^{-1}\Delta_{\mini}^ku_Ia_1)L(a_2^{-1}u_Ja_2)\varphi\right] \\
  & \mathcal{P}_{H,\xi}\left[R(\varphi_2)L(\varphi_3)R(a_1)L(a_2)R(a_1^{-1}u_Ia_1)L(a_2^{-1}\Delta_{\mini}^ku_Ja_2)\varphi\right] \\
   & +\mathcal{P}_{H,\xi}\left[R(\varphi_2)L(\varphi_4)R(a_1)L(a_2)R(a_1^{-1}u_Ia_1)L(a_2^{-1}u_Ja_2)\varphi\right]
\end{aligned}\]

\noindent for all $\varphi\in \mathcal{C}_d^w(G(F))$ and all $(a_1,a_2)\in A_{\mini}^+(I)\times A_{\mini}(J)^+$. The inequality \ref{eq 8.3.10} now follows from Lemma \ref{lemma 1.5.1}(i)(b). $\blacksquare$
\end{itemize}

\end{enumerate}

\subsection{Explicit intertwinings and parabolic induction}\label{section 8.4}

\noindent Let $Q=LU_Q$ be a parabolic subgroup of $G$. Because $G=U(W)\times U(V)$, we have decompositions

\begin{align}\label{eq 8.4.1}
Q=Q_W\times Q_V \mbox{ and } L=L_W\times L_V
\end{align}

\noindent where $Q_W$ and $Q_V$ are parabolic subgroups of $U(W)$ and $U(V)$ respectively and $L_W$, $L_V$ are Levi components of these. By the explicit description of parabolic subgroups of unitary groups (cf.\ Section \ref{section 6.1}), we have

\begin{align}\label{eq 8.4.2}
L_W=GL_E(Z_{1,W})\times\ldots\times GL_E(Z_{a,W})\times U(\widetilde{W})
\end{align}

\begin{align}\label{eq 8.4.3}
L_V=GL_E(Z_{1,V})\times\ldots\times GL_E(Z_{b,V})\times U(\widetilde{V})
\end{align}

\noindent where $Z_{i,W}$, $1\leqslant i\leqslant a$ (respectively $Z_{i,V}$, $1\leqslant i\leqslant b$) are totally isotropic subspaces of $W$ (respectively of $V$) and $\widetilde{W}$ (respectively $\widetilde{V}$) is a non-degenerate subspace of $W$ (respectively of $V$). Let $\widetilde{G}=U(\widetilde{W})\times U(\widetilde{V})$. The pair $(\widetilde{V},\widetilde{W})$ is easily seen to be admissible up to permutation, hence it defines a GGP triple $(\widetilde{G},\widetilde{H},\widetilde{\xi})$ well-defined up to $\widetilde{G}(F)$-conjugation where $\widetilde{G}=U(\widetilde{W})\times U(\widetilde{V})$. For all tempered representations $\widetilde{\sigma}$ of $\widetilde{G}(F)$, we may define as in Section \ref{section 8.2} a continuous linear form $\mathcal{L}_{\widetilde{\sigma}}:\End(\widetilde{\sigma})^\infty\to \mathbb{C}$.

\vspace{2mm}

\noindent Let $\sigma$ be a tempered representation of $L(F)$ which decomposes according to the decompositions \ref{eq 8.4.1}, \ref{eq 8.4.2} and \ref{eq 8.4.3} as a tensor product

\begin{align}\label{eq 8.4.4}
\sigma=\sigma_W\boxtimes \sigma_V
\end{align}

\begin{align}\label{eq 8.4.5}
\sigma_W=\sigma_{1,W}\boxtimes\ldots\boxtimes\sigma_{a,W}\boxtimes\widetilde{\sigma}_{W}
\end{align}

\begin{align}\label{eq 8.4.6}
\sigma_V=\sigma_{1,V}\boxtimes\ldots\boxtimes\sigma_{b,V}\boxtimes\widetilde{\sigma}_{V}
\end{align}

\noindent where $\sigma_{i,W}\in \Temp(GL_E(Z_{i,W}))$ for $1\leqslant i\leqslant a$, $\sigma_{i,V}\in \Temp(GL_E(Z_{i,V}))$ for $1\leqslant i\leqslant b$, $\widetilde{\sigma}_W$ is a tempered representation of $U(\widetilde{W})(F)$ and $\widetilde{\sigma}_{V}$ is a tempered representation of $U(\widetilde{V})(F)$. Let us set $\tilde{\sigma}=\widetilde{\sigma}_{W}\boxtimes\widetilde{\sigma}_{V}$. It is a tempered representation of $\widetilde{G}(F)$. Finally, let us set $\pi=i_Q^G(\sigma)$, $\pi_W=i_{Q_W}^{U(W)}(\sigma_W)$ and $\pi_V=i_{Q_V}^{U(V)}(\sigma_V)$ for the parabolic inductions of $\sigma$, $\sigma_W$ and $\sigma_V$ respectively. We have $\pi=\pi_W\boxtimes \pi_V$.

\vspace{3mm}

\begin{prop}\label{proposition 8.4.1}
With notation as above, we have

$$\mathcal{L}_{\pi}\neq 0\Leftrightarrow \mathcal{L}_{\tilde{\sigma}}\neq 0$$
\end{prop}

\vspace{3mm}

\noindent\ul{Proof}: We will use the following convenient notation. If $X$ is an $E$-vector space of finite dimension, $Q_X=L_XU_X$ is a parabolic subgroup of $GL_E(X)$ with

$$L_X=GL_E(X_1)\times\ldots\times GL_E(X_c)$$

\noindent and we have tempered representations $\sigma_{i,X}$ of $GL_E(X_i)$ for $1\leqslant i\leqslant c$, then we will denote by

$$\sigma_{1,X}\times\ldots\times\sigma_{c,X}$$

\noindent the induced representation $i_{Q_X}^{GL_E(X)}(\sigma_{1,X}\boxtimes\ldots\boxtimes \sigma_{c,X})$. Note that if all the $\sigma_{i,X}$, $1\leqslant i\leqslant c$, are irreducible so is $\sigma_{1,X}\times\ldots\times\sigma_{c,X}$. Similarly, if $X$ is an hermitian space, $Q_X=L_XU_X$ is a parabolic subgroup of $U(X)$ with

$$L_X=GL_E(Z_{1,X})\times\ldots\times GL_E(Z_{d,X})\times U(\widetilde{X})$$

\noindent and we have tempered representations $\sigma_{i,X}$ of $GL_E(Z_{i,X})$ for $1\leqslant i\leqslant d$ and a tempered representation $\widetilde{\sigma}_{X}$ of $U(\widetilde{X})(F)$, then we will denote by

$$\sigma_{1,X}\times\ldots\times\sigma_{d,X}\times\widetilde{\sigma}_{X}$$

\noindent the induced representation $i_{Q_X}^{U(X)}(\sigma_{1,X}\boxtimes\ldots\boxtimes \sigma_{d,X}\boxtimes\widetilde{\sigma}_{X})$. In particular, with these notation, we have

$$\pi_W=\sigma_{1,W}\times\ldots\times\sigma_{a,W}\times\widetilde{\sigma}_{W}$$

\noindent and

$$\pi_V=\sigma_{1,V}\times\ldots\times\sigma_{b,V}\times\widetilde{\sigma}_{V}$$

\noindent By the process of induction by stages, we also have $\pi_W=\sigma'_W\times \widetilde{\sigma}_W$ and $\pi_V=\sigma'_V\times \widetilde{\sigma}_V$, where

$$\sigma'_W=\sigma_{1,W}\times\ldots\times\sigma_{a,W}$$
$$\sigma'_V=\sigma_{1,V}\times\ldots\times\sigma_{b,V}$$

\noindent These are tempered irreducible representations of general linear groups (over $E$). Hence, the statement immediately reduces to the case where $a\leqslant 1$ and $b\leqslant 1$.

\vspace{2mm}

\noindent First, we treat the particular case where $G=G_0$, $H=H_0$ (codimension one case) and $(a,b)=(0,1)$. Then $Q=U(W)\times Q_V$ where $Q_V$ is a maximal proper parabolic subgroup of $G$. Up to conjugating $Q$, we may assume that it is a good parabolic subgroup (cf.\ Section \ref{section 6.4}). Then, to fit with our general notation of Chapter \ref{section 6}, we will change our notation and denote $Q$ by $\overline{Q}$ and $Q_V$ by $\overline{Q}_V$. Set $H_{\overline{Q}}=H\cap \overline{Q}$. Clearly, we have a natural embedding $H_{\overline{Q}}\hookrightarrow L$.

\vspace{2mm}

\noindent We may assume without loss of generality that the invariant scalar product on $\pi$ is given by

$$\displaystyle (e,e')=\int_{\overline{Q}(F)\backslash G(F)} (e(g),e'(g)) dg,\;\;\; e,e'\in \pi$$

\noindent where the scalar product in the integral is the scalar product on $\sigma$. Since we are in the codimension one case, the integral defining $\mathcal{L}_{\pi}$ is absolutely convergent and we have

\begin{align}\label{eq 8.4.7}
\displaystyle \mathcal{L}_{\pi}(e,e')=\int_{H(F)}\int_{\overline{Q}(F)\backslash G(F)} (e(g),e'(gh)) dg\xi(h)dh
\end{align}

\noindent for all $e,e'\in \pi^\infty$. Let us show

\vspace{3mm}

\begin{num}
\item\label{eq 8.4.8} The expression \ref{eq 8.4.7} is absolutely convergent for all $e,e'\in \pi^\infty$.
\end{num}

\vspace{3mm}

\noindent Let $e,e'\in \pi^\infty$ and choose a maximal compact subgroup $K$ of $G(F)$ which is special in the $p$-adic case. Then, for a suitable choice of Haar measure on $K$, we have

\[\begin{aligned}
\displaystyle \int_{\overline{Q}(F)\backslash G(F)} \left\lvert(e(g),e'(gh))\right\rvert dg & =\int_K \left\lvert(e(k),e'(kh))\right\rvert dk \\
 & =\int_K \delta_{\overline{Q}}(l_{\overline{Q}}(kh))^{1/2} \left\lvert\left(e(k),\sigma\left(l_{\overline{Q}}(kh)\right)e'\left(k_{\overline{Q}}(kh)\right)\right)\right\rvert dk
\end{aligned}\]

\noindent for all $h\in H(F)$. Here, as usual $l_{\overline{Q}}:G(F)\to L(F)$ and $k_{\overline{Q}}:G(F)\to K$ are maps such that $l_{\overline{Q}}(g)^{-1}gk_{\overline{Q}}(g)^{-1}\in U_{\overline{Q}}(F)$ (the unipotent radical of $\overline{Q}(F)$) for all $g\in G(F)$. Since $\sigma$ is tempered and the maps $k\in K\mapsto e(k)\in \sigma$, $k\in K\mapsto e'(k)\in \sigma$ have bounded image, it follows that

$$\displaystyle \int_{\overline{Q}(F)\backslash G(F)} \left\lvert(e(g),e'(gh))\right\rvert dg\ll \int_K \delta_{\overline{Q}}(l_{\overline{Q}}(kg))^{1/2} \Xi^L(l_{\overline{Q}}(kh))dk=\Xi^G(h)$$

\noindent for all $h\in H(F)$, where the last equality is Proposition \ref{proposition 1.5.1}(iii). The absolute convergence of \ref{eq 8.4.7} now follows from Lemma \ref{lemma 6.5.1}(i).

\vspace{2mm}

\noindent Since $\overline{Q}$ is a good parabolic subgroup, by Proposition \ref{proposition 6.4.1}(i) the quotient $H_{\overline{Q}}(F)\backslash H(F)$ has negligible complement in $\overline{Q}(F)\backslash G(F)$. Hence, if we choose Haar measures compatibly, we have

$$\displaystyle \int_{\overline{Q}(F)\backslash G(F)}\varphi(g) dg=\int_{H_{\overline{Q}}(F)\backslash H(F)} \varphi(h)dh$$

\noindent for all $\varphi\in L^1\left(\overline{Q}(F)\backslash G(F),\delta_{\overline{Q}}\right)$. Thus, \ref{eq 8.4.7} becomes

$$\displaystyle \mathcal{L}_{\pi}(e,e')=\int_{H(F)}\int_{H_{\overline{Q}}(F)\backslash H(F)} \left(e(h'),e'(h'h)\right) dh'dh$$

\noindent for all $e,e'\in \pi^\infty$, the double integral being absolutely convergent by \ref{eq 8.4.8}. Switching the two integrals, we get

\begin{align}\label{eq 8.4.9}
\displaystyle \mathcal{L}_{\pi}(e,e')=\int_{\left(H_{\overline{Q}}(F)\backslash H(F)\right)^2} \mathcal{L}_{\sigma}(e(h),e'(h'))dhdh'
\end{align}

\noindent for all $e,e'\in \pi^\infty$, where we have set

$$\displaystyle \mathcal{L}_{\sigma}(v,v')=\int_{H_{\overline{Q}}(F)} \left(v,\sigma(h_{\overline{Q}})v'\right)\delta_{H_{\overline{Q}}}(h_{\overline{Q}})^{1/2}d_Lh_{\overline{Q}}$$

\noindent for all $v,v'\in \sigma^\infty$. The presence in the integral above of $\delta_{H_{\overline{Q}}}$ instead of $\delta_{\overline{Q}}$ follows from Proposition \ref{proposition 6.8.1}(ii). Set

$$\displaystyle \mathcal{L}_{\sigma}(T_{\sigma})=\int_{H_{\overline{Q}}(F)} \Tr\left(\sigma(h_{\overline{Q}}^{-1})T_\sigma\right)\delta_{H_{\overline{Q}}}(h_{\overline{Q}})^{1/2}d_Lh_{\overline{Q}}$$

\noindent for all $T_{\sigma}\in \End(\sigma)^\infty$. We have

\vspace{3mm}

\begin{num}
\item\label{eq 8.4.10} The integral defining $\mathcal{L}_{\sigma}$ is absolutely convergent and $\mathcal{L}_{\sigma}$ is a continuous linear form on $\End(\sigma)^\infty$.
\end{num}

\vspace{3mm}

\noindent Indeed, this follows from Proposition \ref{proposition 6.8.1}(iii) as $\sigma$ is tempered. We now prove the following

\begin{align}\label{eq 8.4.11}
\mathcal{L}_{\pi}\neq 0\Leftrightarrow \mathcal{L}_{\sigma}\neq 0
\end{align}

\noindent By \ref{eq 8.4.9} and the density of $\pi^\infty\otimes \overline{\pi^\infty}$ in $\End(\pi)^\infty$, we see that if $\mathcal{L}_{\pi}$ is nonzero then $\mathcal{L}_{\sigma}$ is nonzero. Let us prove the converse. The analytic fibration $H(F)\to H_{\overline{Q}}(F)\backslash H(F)$ is locally trivial. Let $s:\mathcal{U}\to H(F)$ be an analytic section over an open subset $\mathcal{U}$ of $H_{\overline{Q}}(F)\backslash H(F)$. For $\varphi\in C_c^\infty(\mathcal{U},\sigma^\infty)$ a smooth compactly supported function from $\mathcal{U}$ to $\sigma^\infty$, the following assignment

$$e_\varphi(g)=\left\{
\begin{array}{ll}
\delta_{\overline{Q}}(l)^{1/2}\sigma(l)\varphi(h) \mbox{ if } g=lus(h) \mbox{ with } l\in L(F), u\in U_{\overline{Q}}(F), h\in \mathcal{U} \\
0 \mbox{ otherwise}
\end{array}
\right.
$$

\noindent defines an element of $\pi^\infty$. By \ref{eq 8.4.9}, we have

$$\displaystyle \mathcal{L}_{\pi}(e_{\varphi},e_{\varphi'})=\int_{\left(H_{\overline{Q}}(F)\backslash H(F)\right)^2} \mathcal{L}_{\sigma}\left(\varphi(h),\varphi'(h')\right) dhdh'$$

\noindent for all $\varphi,\varphi'\in C_c^\infty(\mathcal{U},\sigma^\infty)$. Now, assume that $\mathcal{L}_{\sigma}$ is nonzero. Because $\sigma^\infty\otimes \overline{\sigma^\infty}$ is dense in $\End(\sigma)^\infty$, there exists $v_0,v'_0\in \sigma^\infty$ such that $\mathcal{L}_{\sigma}(v_0,v'_0)\neq 0$. Setting $\varphi(h)=f(h)v_0$ and $\varphi'(h)=\overline{f(h)}v'_0$ where $f\in C_c^\infty(\mathcal{U})$ is nonzero in the formula above, we get $\mathcal{L}_{\pi}(e_{\varphi},e_{\varphi'})\neq 0$ hence $\mathcal{L}_{\pi}$ doesn't vanish. This ends the proof of \ref{eq 8.4.11}.

\vspace{2mm}

\noindent By \ref{eq 8.4.11}, we are now reduced to proving

\begin{align}\label{eq 8.4.12}
\mathcal{L}_{\sigma}\neq 0\Leftrightarrow \mathcal{L}_{\tilde{\sigma}}\neq 0
\end{align}

\noindent In order to prove \ref{eq 8.4.12}, we will need a precise description of $H_{\overline{Q}}$ and of the embedding $H_{\overline{Q}}\hookrightarrow L$. Since $\overline{Q}_V$ is a maximal proper parabolic subgroup of $U(V)$ it is the stabilizer of a totally isotropic subspace $Z'$ of $V$. The quotient $\overline{Q}\backslash G$ classifies the totally isotropic subspaces of $V$ of the same dimension as $Z'$. The action of $H=U(W)$ on that space has two orbits: one is open and corresponds to the subspaces $Z''$ such that $dim(Z''\cap W)=dim(Z')-1$, the other is closed and corresponds to the subspaces $Z''$ such that $dim(Z''\cap W)=dim(Z')$. Since we are assuming that $\overline{Q}$ is a good parabolic subgroup, $Z'$ is in the open orbit. Hence $Z'_W=Z'\cap W$ is an hyperplane in $Z'$. Moreover, we have

$$L_V=GL_E(Z')\times U(\widetilde{V})$$

\noindent where $\widetilde{V}$ is a non-degenerate subspace of $V$ orthogonal to $Z'$. Since $Z'_W$ is an hyperplane in $Z'$, up to $\overline{Q}(F)$-conjugation we have $\widetilde{V}\subset W$ and so we may as well assume that it is the case. Note that we have a natural identification $H_{\overline{Q}}=U(W)\cap \overline{Q}_V$. The exact sequence

$$1\to U_{\overline{Q}}\to \overline{Q}_V\to GL_E(Z')\times U(\widetilde{V})\to 1$$

\noindent induces an exact sequence

\begin{align}\label{eq 8.4.13}
1\to U^\natural_{\overline{Q}}\to H_{\overline{Q}}\to P_{Z'}\times U(\widetilde{V})\to 1
\end{align}

\noindent where $U_{\overline{Q}}^\natural=U(W)\cap U_{\overline{Q}}$ and $P_{Z'}$ is the mirabolic subgroup of elements $g\in GL_E(Z')$ preserving the hyperplane $Z'_W$ and acting trivially on the quotient $Z'/Z'_W$. On the other hand, we have

$$L=U(W)\times GL_E(Z')\times U(\widetilde{V})$$

\noindent and the embedding $H_{\overline{Q}}\hookrightarrow L$ is the product of the three following maps

\vspace{2mm}

\begin{itemize}
\renewcommand{\labelitemi}{$\bullet$}
\item The natural inclusion $H_{\overline{Q}}\subseteq U(W)$;

\item The projection $H_{\overline{Q}}\twoheadrightarrow P_{Z'}$ followed by the natural inclusion $P_{Z'}\subseteq GL_E(Z')$;

\item The projection $H_{\overline{Q}}\twoheadrightarrow U(\widetilde{V})$.
\end{itemize}

\vspace{2mm}

\noindent Let $\widetilde{D}$ be a line such that $(Z'_W)^\perp\cap W=Z'_W\oplus\left(\widetilde{D}\oplus^\perp \widetilde{V}\right)$. Then, the unipotent group $U_{\overline{Q}}^\natural$ may be described as the subgroup of $U(W)$ stabilizing the subspace $\widetilde{D}\oplus Z'_W$. Fix a basis $z'_1,\ldots,z'_\ell$ of $Z'_W$ and let $B_{Z'}$ be the Borel subgroup of $GL_E(Z')$ fixing the complete flag

$$\langle z'_1\rangle\subsetneq \langle z'_1,z'_2\rangle\subsetneq\ldots\subsetneq \langle z'_1,\ldots,z'_\ell \rangle=Z'_W\subsetneq Z'$$

\noindent and denote by $N_{Z'}$ its unipotent radical. Let $\widetilde{\xi}$ be a generic character of $N_{Z'}(F)$. Let us denote by $\widetilde{N}$ and $\widetilde{H}$ the inverse images in $H_{\overline{Q}}$ of the subgroups $N_{Z'}$ and $N_{Z'}\times U(\widetilde{V})$ of $P_{Z'}\times U(\widetilde{V})$ by the last map of \ref{eq 8.4.13}. Recall that $\widetilde{G}=U(\widetilde{V})\times U(W)$. We have a natural 'diagonal' embedding $\widetilde{H}\hookrightarrow \widetilde{G}$ (that is: the product of the natural projection $\widetilde{H}\twoheadrightarrow U(\widetilde{V})$ and the natural inclusion $\widetilde{H}\hookrightarrow U(W)$) and if we extend $\widetilde{\xi}$ to a character of $\widetilde{H}(F)$ through the projection $\widetilde{H}(F)\twoheadrightarrow N_{Z'}(F)$, then the triple $(\widetilde{G},\widetilde{H},\widetilde{\xi})$ is a GGP triple corresponding to the pair of hermitian spaces $(\widetilde{V},\widetilde{W})$. We can of course use this triple to define $\mathcal{L}_{\widetilde{\sigma}}$. We henceforth assume that it is the case.

\vspace{2mm}

\noindent The representation $\sigma$ decomposes as $\sigma=\sigma_{1,V}\boxtimes \widetilde{\sigma}$ where $\sigma_{1,V}$ is a tempered irreducible representation of $GL_E(Z')$. The subspace $\sigma_{1,V}^\infty\otimes \widetilde{\sigma}^\infty$ is dense in $\sigma^\infty$. Hence $\mathcal{L}_{\sigma}$ is nonzero if and only if there exist vectors $\widetilde{v},\widetilde{v}'\in \widetilde{\sigma}^\infty$ and $w,w'\in \sigma_{1,V}^\infty$ such that

$$\mathcal{L}_{\sigma}(w\otimes \widetilde{v},w'\otimes \widetilde{v}')\neq 0$$

\noindent Let us fix a Whittaker model $\sigma_{1,V}^\infty\hookrightarrow C^\infty\left( N_{Z'}(F)\backslash GL_E(Z'),\widetilde{\xi}^{-1}\right)$ for $\sigma_{1,V}^\infty$ (such a model exists because $\sigma_{1,V}$ is a tempered representation of a general linear group). For $w\in \sigma_{1,V}^\infty$, we will denote by $W_w$ the corresponding Whittaker function. We have the following

\vspace{3mm}

\begin{num}
\item\label{eq 8.4.14} For a suitable choice of Haar measures, we have the equality (recall that $P_{Z'}$ is a subgroup of $U(W)$)
$$\displaystyle \mathcal{L}_{\sigma}(w\otimes \widetilde{v},w'\otimes\widetilde{v}')=\int_{\left( N_{Z'}(F)\backslash P_{Z'}(F)\right)^2} W_w(p)\overline{W_{w'}(p')} \mathcal{L}_{\widetilde{\sigma}}\left(\widetilde{\sigma}(p)\widetilde{v},\widetilde{\sigma}(p')\widetilde{v}'\right) \delta_{H_{\overline{Q}}}(pp')^{-1/2} d_Rpd_Rp'$$
for all $w,w'\in \sigma_{1,V}^\infty$ and all $\widetilde{v},\widetilde{v}'\in \widetilde{\sigma}^\infty$ and where the integral on the right is absolutely convergent.
\end{num}

\vspace{3mm}

\noindent Before proving this equality, let us explain how we can deduce \ref{eq 8.4.12} from it. The implication $\mathcal{L}_{\sigma}\neq 0\Rightarrow \mathcal{L}_{\widetilde{\sigma}}\neq 0$ is now obvious. For the converse direction, we use the fact that the Kirillov model of $\sigma_{1,V}^\infty$, that is the restriction of the Whittaker model to functions on $P_{Z'}(F)$, contains $C_c^\infty\left(N_{Z'}(F)\backslash P_{Z'}(F),\widetilde{\xi}^{-1}\right)$ (cf.\ Theorem 6 of \cite{GK} in the $p$-adic case and Theorem 1 of \cite{Kem} in the real case). Since the analytic fibration $P_{Z'}(F)\to N_{Z'}(F)\backslash P_{Z'}(F)$ is locally trivial, we can now argue in a similar way as in the proof of \ref{eq 8.4.11}.

\vspace{2mm}

\noindent We will now prove \ref{eq 8.4.14}. Let us fix $w,w'\in \sigma_{1,V}^\infty$. Let $H_\natural$ be the inverse image of $U(\widetilde{V})$ in $H_{\overline{Q}}$ via the last map of the exact sequence \ref{eq 8.4.13}. The group $H_\natural(F)$ is unimodular and for suitable choices of Haar measures, we have $d_L h_{\overline{Q}}=dh_\natural d_Lp$ where $dh_\natural$ is a Haar measure on $H_\natural(F)$ and $d_Lp$ is a left Haar measure on $P_{Z'}(F)$. Moreover, the modular character $\delta_{H_{\overline{Q}}}$ is trivial on $H_\natural (F)$. Thus, because of the description of the embedding $H_{\overline{Q}}\hookrightarrow L$, we have

\begin{align}\label{eq 8.4.15}
\displaystyle \mathcal{L}_{\sigma}(w\otimes \widetilde{v},w'\otimes \widetilde{v}')=\int_{P_{Z'}(F)}\int_{H_\natural(F)}  \left(\widetilde{v},\widetilde{\sigma}(ph_\natural)\widetilde{v}'\right) (w,\sigma_{1,V}(p)w')\delta_{H_{\overline{Q}}}(p)^{1/2} dh_\natural d_Lp
\end{align}

\noindent for all $\widetilde{v},\widetilde{v}'\in \widetilde{\sigma}^\infty$, the integral being absolutely convergent by \ref{eq 8.4.10}. Let us define

$$\displaystyle \mathcal{P}_{\widetilde{H},\widetilde{\xi}}(\varphi)=\int_{\widetilde{H}(F)}^* \varphi(\widetilde{h})\widetilde{\xi}(\widetilde{h})d\widetilde{h}$$

$$\displaystyle \mathcal{P}^1_{w,w'}(\varphi)=\int_{\left( N_{Z'}(F)\backslash P_{Z'}(F)\right)^2} W_w(p)\overline{W_{w'}(p')} \mathcal{P}_{\widetilde{H},\widetilde{\xi}}\left(L(p)R(p')\varphi\right) \delta_{H_{\overline{Q}}}(pp')^{-1/2} d_Rpd_Rp'$$

$$\displaystyle \mathcal{P}^2_{w,w'}(\varphi)= \int_{P_{Z'}(F)} \int_{H_\natural(F)} \varphi(ph_\natural) (w,\sigma_{1,V}(p)w')\delta_{H_{\overline{Q}}}(p)^{1/2} dh_\natural d_Lp$$

\noindent for all $\varphi\in \mathcal{C}^w(\widetilde{G}(F))$. The first expression above is the generalized $\widetilde{\xi}$-integral on $\widetilde{H}(F)$ of Section \ref{section 8.1}. Because of \ref{eq 8.4.15} and since $\widetilde{\sigma}$ is tempered, the claim \ref{eq 8.4.14} will follow from

\vspace{3mm}

\begin{num}
\item\label{eq 8.4.16} The integrals defining $\mathcal{P}^1_{w,w'}$ and $\mathcal{P}_{w,w'}^2$ are absolutely convergent, moreover they define continuous linear forms on $\mathcal{C}^w(\widetilde{G}(F))$ and if Haar measures are chosen suitably we have
$$\mathcal{P}^1_{w,w'}(\varphi)=\mathcal{P}_{w,w'}^2(\varphi)$$
for all $\varphi \in \mathcal{C}^w(\widetilde{G}(F))$.
\end{num}

\vspace{3mm}

\noindent Let us start by proving the absolute convergence and continuity of $\mathcal{P}_{w,w'}^2$. It suffices to show that the integral

$$\displaystyle \int_{P_{Z'}(F)} \int_{H_\natural(F)} \Xi^{\widetilde{G}}(ph_\natural)\sigma_{\widetilde{G}}(ph_\natural)^d \left\lvert (w,\sigma_{1,V}(p)w')\right\rvert\delta_{H_{\overline{Q}}}(p)^{1/2} dh_\natural d_Lp$$

\noindent is absolutely convergent for all $d>0$. Since $\sigma_{1,V}$ is a tempered representation this last integral is essentially bounded by

$$\displaystyle \int_{P_{Z'}(F)} \int_{H_\natural(F)} \Xi^{\widetilde{G}}(ph_\natural)\sigma_{\widetilde{G}}(ph_\natural)^d \Xi^{GL_E(Z')}(p)\delta_{H_{\overline{Q}}}(p)^{1/2} dh_\natural d_Lp=\int_{H_{\overline{Q}}(F)} \Xi^L(h_{\overline{Q}}) \sigma(h_{\overline{Q}})^d\delta_{H_{\overline{Q}}}(h_{\overline{Q}})^{1/2}d_Lh_{\overline{Q}}$$

\noindent which by Proposition \ref{proposition 6.8.1}(iii) is an absolutely convergent integral.

\vspace{2mm}

\noindent Let us now show the absolute convergence and continuity of $\mathcal{P}_{w,w'}^1$. By Lemma \ref{lemma 8.3.1}(ii), it is sufficient to show the following

\vspace{3mm}

\begin{num}
\item\label{eq 8.4.17} For all $d>0$ and all $w_0\in \sigma_{1,V}^\infty$ the integral
$$\displaystyle \int_{N_{Z'}(F)\backslash P_{Z'}(F)} \left\lvert W_{w_0}(p)\right\rvert \Xi^{\widetilde{H}\backslash \widetilde{G}}(p) \sigma(p)^d \delta_{H_{\overline{Q}}}(p)^{-1/2}d_Rp$$
is absolutely convergent.
\end{num}

\vspace{3mm}

\noindent Let $d>0$ and $w_0\in \sigma_{1,V}$. Let $\widetilde{T}$ be the maximal subtorus of $GL_E(Z'_W)$ stabilizing the lines $\langle z'_1\rangle,\ldots, \langle z'_\ell\rangle$. Notice that $\widetilde{T}$ plays the role of the torus $T$ for the GGP triple $(\widetilde{G},\widetilde{H},\widetilde{\xi})$ (cf.\ Section \ref{section 6.2}). Also, let us fix a maximal (special in the $p$-adic case) compact subgroup $K_{Z'_W}$ of $GL_E(Z'_W)$. Then, we have $P_{Z'}(F)=N_{Z'}(F)\widetilde{T}(F)K_{Z'_W}$ together with a decomposition of the Haar measure $d_Rp=\delta_{P_{Z'}}(\widetilde{t})\delta_{B_{Z'}}(\widetilde{t})^{-1}dnd\widetilde{t}dk$. Hence,

\[\begin{aligned}
\displaystyle \int_{N_{Z'}(F)\backslash P_{Z'}(F)} \left\lvert W_{w_0}(p)\right\rvert \Xi^{\widetilde{H}\backslash \widetilde{G}}(p) & \sigma(p)^d \delta_{H_{\overline{Q}}}(p)^{-1/2}d_Rp= \\
 & \int_{K_{Z'_W}}\int_{\widetilde{T}(F)} \left\lvert W_{w_0}(\widetilde{t}k)\right\rvert \Xi^{\widetilde{H}\backslash \widetilde{G}}(\widetilde{t}k) \sigma(\widetilde{t}k)^d \delta_{H_{\overline{Q}}}(\widetilde{t})^{-1/2} \delta_{P_{Z'}}(\widetilde{t}) \delta_{B_{Z'}}(\widetilde{t})^{-1}d\widetilde{t} dk
\end{aligned}\]

\noindent By Proposition \ref{proposition 6.7.1}(i)(a) and (ii)(a), there exists a $d_1>0$ such that $\Xi^{\widetilde{H}\backslash \widetilde{G}}(\widetilde{t}k)\ll \delta_{\widetilde{N}}(\widetilde{t})^{1/2}\sigma(\widetilde{t})^{d_1}$ for all $\widetilde{t}\in \widetilde{T}(F)$ and all $k\in K_{Z'_W}$. Also,we have $\sigma(\widetilde{t}k)\ll \sigma(\widetilde{t})$ for all $\widetilde{t}\in \widetilde{T}(F)$ and all $k\in K_{Z'_W}$. We have an isomorphism

$$\widetilde{T}\simeq \left(R_{E/F}\mathbb{G}_m\right)^\ell$$
$$\widetilde{t}\mapsto (\widetilde{t}_i)_{1\leqslant i\leqslant\ell}$$ 

\noindent where $\widetilde{t}_i$, $1\leqslant i\leqslant \ell$, denotes the eigenvalue of $\widetilde{t}$ acting on $z'_i$. Since $W_{w_0}(\widetilde{t}k)=W_{\sigma_{1,V}(k)w_0}(\widetilde{t})$ and the map $k\in K_{Z'_W}\mapsto \sigma_{1,V}(k)w_0$ has bounded image, by Lemma \ref{lemma B.2.1}, it is sufficient to prove the existence of $R>0$ such that the integral

$$\displaystyle \int_{\widetilde{T}(F)} \Xi^{GL_E(Z')}(\widetilde{t})\delta_{\widetilde{N}}(\widetilde{t})^{1/2}\sigma(\widetilde{t})^{d+d_1} \delta_{H_{\overline{Q}}}(\widetilde{t})^{-1/2}\delta_{P_{Z'}}(\widetilde{t}) \delta_{B_{Z'}}(\widetilde{t})^{-1}\prod_{i=1}^{\ell} \max(1,\lvert \widetilde{t}_i\rvert_E)^{-R}d\widetilde{t} $$

\noindent converges absolutely. Here, we have denoted by $\lvert .\rvert_E$ the normalized absolute value of $E$. By Proposition \ref{proposition 1.5.1}(ii), there exists a $d_2>0$ such that $\Xi^{GL_E(Z')}(\widetilde{t})\ll \delta_{B_{Z'}}(\widetilde{t})^{1/2}\sigma(\widetilde{t})^{d_2}$ for all $\widetilde{t}\in \widetilde{T}(F)$. Moreover, we have $\delta_{\widetilde{N}}(\widetilde{t})=\delta_{H_{\overline{Q}}}(\widetilde{t})\delta_{B_{Z'}}(\widetilde{t})\delta_{P_{Z'}}(\widetilde{t})^{-1}$ and

$$\displaystyle \delta_{P_{Z'}}(\widetilde{t})=\prod_{i=1}^{\ell}\lvert \widetilde{t}_i\rvert_E$$

\noindent for all $\widetilde{t}\in\widetilde{T}(F)$. Hence, it suffices to prove that there exists $R>0$ such that the integral

$$\displaystyle \int_{\widetilde{T}(F)}\sigma(\widetilde{t})^{d+d_1+d_2}\prod_{i=1}^{\ell} \lvert \widetilde{t}_i\rvert_E \prod_{i=1}^{\ell}\max(1,\lvert \widetilde{t}_i\rvert_E)^{-R} d\widetilde{t}$$

\noindent converges. We have

$$\displaystyle \sigma(\widetilde{t})\ll\prod_{i=1}^{\ell} \log \left(1+\max ( \lvert\widetilde{t}_i^{-1}\rvert_E, \lvert\widetilde{t}_i\rvert_E)\right)$$

\noindent for all $\widetilde{t}\in \widetilde{T}(F)$ and for suitable additive Haar measures $d\widetilde{t}_i$ on $E$, $1\leqslant i\leqslant \ell$, we have

$$\displaystyle d\widetilde{t}=\prod_{i=1}^{\ell}\lvert\widetilde{t}_i\rvert_E^{-1} d\widetilde{t}_i$$

\noindent So finally, we only need to prove the existence of $R>0$ such that the integral

$$\displaystyle \int_E \log\left(1+\max(\lvert x\rvert_E^{-1},\lvert x\rvert_E)\right)^{d+d_1+d_2} \max(1,\lvert x\rvert_E)^{-R}dx$$

\noindent converges. Any $R>1$ fulfills this condition. This ends the proof of \ref{eq 8.4.17} and hence of the absolute convergence and continuity of $\mathcal{P}_{w,w'}^1$.

\vspace{2mm}

\noindent We are now only left with proving the equality of the functionals $\mathcal{P}^1_{w,w'}$ and $\mathcal{P}_{w,w'}^2$. Because we have shown that these two functionals are continuous linear forms on $\mathcal{C}^w(\widetilde{G}(F))$ and since $\mathcal{C}(\widetilde{G}(F))$ is a dense subspace of $\mathcal{C}^w(\widetilde{G}(F))$ (by \ref{eq 1.5.1}), we only need to prove the equality

$$\mathcal{P}^1_{w,w'}(\varphi)=\mathcal{P}^2_{w,w'}(\varphi)$$

\noindent for all $\varphi\in \mathcal{C}(\widetilde{G}(F))$. Let $\varphi\in \mathcal{C}(\widetilde{G}(F))$. By Theorem 6.2 of \cite{Ber2} in the $p$-adic case and Theorem 10.3 of \cite{Bar} in the real case, if we choose the Haar measures suitably, we have

$$\displaystyle (w_0,w'_0)=\int_{N_{Z'}(F)\backslash P_{Z'}(F)} W_{w_0}(p) \overline{W_{w'_0}(p)} d_Rp$$

\noindent for all $w_0,w'_0\in \sigma_{1,V}^\infty$, the integral being absolutely convergent. Hence, we have

\begin{align}\label{eq 8.4.18}
\displaystyle \mathcal{P}^2_{w,w'}(\varphi) =\int_{P_{Z'}(F)} \int_{H_\natural(F)} \varphi(ph_\natural) \int_{N_{Z'}(F)\backslash P_{Z'}(F)} W_{w}(p') \overline{W_{w'}(p'p)} d_Rp' \delta_{H_{\overline{Q}}}(p)^{1/2}dh_\natural d_Lp
\end{align}

\noindent Let us assume one moment that the above triple integral is absolutely convergent. Then, we would have

\[\begin{aligned}
\displaystyle \mathcal{P}^2_{w,w'}(\varphi) & =\int_{N_{Z'}(F)\backslash P_{Z'}(F)} W_{w}(p') \int_{P_{Z'}(F)} \int_{H_\natural(F)} \varphi(ph_\natural) \overline{W_{w'}(p'p)} \delta_{H_{\overline{Q}}}(p)^{1/2} dh_\natural d_Lp d_Rp' \\
 & =\int_{N_{Z'}(F)\backslash P_{Z'}(F)} W_{w}(p') \int_{P_{Z'}(F)}\int_{H_\natural(F)} \varphi(p'^{-1}ph_\natural) \overline{W_{w'}(p)} \delta_{H_{\overline{Q}}}(p)^{1/2}  dh_\natural d_Lp \delta_{H_{\overline{Q}}}(p')^{-1/2}d_Rp' \\
 & =\int_{N_{Z'}(F)\backslash P_{Z'}(F)} W_{w}(p') \int_{P_{Z'}(F)} \int_{H_\natural(F)} \varphi(p'^{-1}h_\natural p) \overline{W_{w'}(p)} \delta_{H_{\overline{Q}}}(p)^{1/2} \delta_{H_\natural}(p)^{-1} dh_\natural d_Lp \\
 &\omit\hfill $\delta_{H_{\overline{Q}}}(p')^{-1/2}d_Rp'$ \\
 & =\int_{N_{Z'}(F)\backslash P_{Z'}(F)} W_{w}(p') \int_{P_{Z'}(F)} \int_{H_\natural(F)} \varphi(p'^{-1}h_\natural p) \overline{W_{w'}(p)} \delta_{H_{\overline{Q}}}(p)^{-1/2} dh_\natural d_Rp \delta_{H_{\overline{Q}}}(p')^{-1/2}d_Rp'\\
 & =\int_{N_{Z'}(F)\backslash P_{Z'}(F)} W_{w}(p') \int_{N_{Z'}(F)\backslash P_{Z'}(F)} \overline{W_{w'}(p)} \int_{N_{Z'}(F)} \int_{H_\natural(F)} \varphi(p'^{-1} h_\natural n p)\widetilde{\xi}(n)dh_\natural dn  \\
 &\omit\hfill $\delta_{H_{\overline{Q}}}(p)^{-1/2} d_Rp \delta_{H_{\overline{Q}}}(p')^{-1/2}d_Rp'$ \\
 & =\int_{\left(N_{Z'}(F)\backslash P_{Z'}(F)\right)^2} W_w(p')\overline{W_{w'}(p)}\int_{\widetilde{H}(F)} \varphi(p'^{-1}\widetilde{h}p)\widetilde{\xi}(\widetilde{h})d\widetilde{h}\delta_{H_{\overline{Q}}}(pp')^{-1/2} d_Rpd_Rp' \\
 & =\int_{\left(N_{Z'}(F)\backslash P_{Z'}(F)\right)^2} W_w(p')\overline{W_{w'}(p)}\mathcal{P}_{\widetilde{H},\widetilde{\xi}}\left(L(p')R(p)\varphi\right)\delta_{H_{\overline{Q}}}(pp')^{-1/2} d_Rpd_Rp' \\
 & =\mathcal{P}^1_{w,w'}(\varphi)
\end{aligned}\]

\noindent where on the second line we made the variable change $p\mapsto p'^{-1}p$, on the third line we made the variable change $h_\natural\mapsto p^{-1}h_\natural p$ and on the fourth line we used the easily checked equality $\delta_{H_{\overline{Q}}}(p)d_Lp=\delta_{H_\natural}(p)d_Rp$. This proves the equality we wanted. Hence, it only remains to show the absolute convergence of \ref{eq 8.4.18}. It is certainly equivalent to show the absolute convergence of the antepenultimate integral above i.e., we need to show that the integral

$$\displaystyle \int_{\left(N_{Z'}(F)\backslash P_{Z'}(F)\right)^2} \left\lvert W_w(p')\right\rvert \left\lvert W_{w'}(p)\right\rvert\int_{\widetilde{H}(F)} \lvert\varphi(p'^{-1}\widetilde{h}p)\rvert d\widetilde{h}\delta_{H_{\overline{Q}}}(pp')^{-1/2} d_Rpd_Rp'$$

\noindent converges. By Proposition \ref{proposition 6.7.1}(v), there exists $d>0$ such that this integral is essentially bounded by

$$\displaystyle \int_{\left(N_{Z'}(F)\backslash P_{Z'}(F)\right)^2} \left\lvert W_w(p')\right\rvert \left\lvert W_{w'}(p)\right\rvert \Xi^{\widetilde{H}\backslash \widetilde{G}}(p) \Xi^{\widetilde{H}\backslash \widetilde{G}}(p')\sigma(p)^d\sigma(p')^d\delta_{H_{\overline{Q}}}(pp')^{-1/2} d_Rpd_Rp'$$

\noindent But by \ref{eq 8.4.17}, we already know that such an integral converges. This ends the proof of \ref{eq 8.4.16} and hence of the proposition in the particular case we have been considering.

\vspace{5mm}

\noindent We now explain how the other cases reduce to the particular case we just treated. We distinguish three cases

\vspace{5mm}

\begin{itemize}
\renewcommand{\labelitemi}{$\bullet$}

\item Case where $a=1$ and $b=0$. Let us define a new hermitian space $(V',h')$ by $V'=V\oplus Ez'_0$ and $h'(v+\lambda z'_0)=h(v)-N(\lambda)h(z_0)$. We have a natural inclusion of hermitian spaces $V\subset V'$ and the pair $(V,V')$ is admissible and gives rise to a GGP triple $(G',H',\xi')$ where $G'=U(V)\times U(V')$, $H'=U(V)$ and $\xi'=1$. Note that this new GGP triple is of codimension one. Let $Z'_+=Z_+\oplus Ez_{r+1}$ where $z_{r+1}=z_0+z'_0$ (where $Z_+$ is defined as in Section \ref{section 6.2}). Then the group $GL_E(Z'_+)\times U(W)$ is a Levi subgroup of $U(V')$. Let $\sigma'$ be any tempered irreducible representation of $GL_E(Z'_+)$. Then, we have the chain of equivalences

\[\begin{aligned}
\mathcal{L}_\pi=\mathcal{L}_{\pi_W\boxtimes \pi_V}\neq 0 \Leftrightarrow \mathcal{L}_{\pi_V\boxtimes(\pi_W\times \sigma')}\neq 0 & \Leftrightarrow \mathcal{L}_{\pi_V\boxtimes (\widetilde{\sigma}_W\times (\sigma_{1,W}\times \sigma'))}\neq 0 \\
 & \Leftrightarrow \mathcal{L}_{\widetilde{\sigma}}=\mathcal{L}_{\widetilde{\sigma}_W\boxtimes \pi_V}\neq 0
\end{aligned}\]

\noindent where in the first and third equivalences we have used the case already treated (for the triple $(G',H',\xi')$) and in the second equivalence we have used induction by stages.

\item Case where $a=0$ and $b=1$. Note that we already solved this case when the triple $(G,H,\xi)$ is of codimension one. We proceed by induction on $\dim(V)$ (if $\dim(V)=1$ then we are in the codimension one case). Let us introduce two GGP triples $(G',H',\xi')$ and $(G'',H'',\xi'')$ relative to the admissible pairs $(W\oplus Z,V)$ and $(\widetilde{V},W\oplus Z)$ respectively ($Z$ is defined in Section \ref{section 6.2}). Let $\sigma'$ be any irreducible tempered representation of $GL_E(Z_+)$. Then, we have the chain of equivalences

\[\begin{aligned}
\mathcal{L}_\pi=\mathcal{L}_{\pi_W\boxtimes \pi_V}\neq 0\Leftrightarrow \mathcal{L}_{(\pi_W\times \sigma')\boxtimes \pi_V}=\mathcal{L}_{(\pi_W\times \sigma')\boxtimes (\sigma_{1,V}\times\widetilde{\sigma}_V)}\neq 0 & \Leftrightarrow \mathcal{L}_{\widetilde{\sigma}_V\boxtimes (\pi_W\times \sigma')}\neq 0 \\
 & \Leftrightarrow \mathcal{L}_{\widetilde{\sigma}}=\mathcal{L}_{\widetilde{\sigma}_V\boxtimes \pi_W}\neq 0
\end{aligned}\]

\noindent where in the first equivalence we applied the case $(a,b)=(1,0)$ to the triple $(G',H',\xi')$, in the second equivalence we applied the case $(a,b)=(0,1)$ to the same triple (which is of codimension one) and in the third equivalence we applied the induction hypothesis to the triple $(G'',H'',\xi'')$. 

\item Case where $a=b=1$. Then, we have the chain of equivalences

$$\mathcal{L}_\pi=\mathcal{L}_{(\sigma_{1,W}\times \widetilde{\sigma}_W)\boxtimes \pi_V}\neq 0 \Leftrightarrow \mathcal{L}_{\widetilde{\sigma}_W\boxtimes (\sigma_{1,V}\times \widetilde{\sigma}_V)}=\mathcal{L}_{\widetilde{\sigma}_W\boxtimes \pi_V}\neq 0 \Leftrightarrow \mathcal{L}_{\widetilde{\sigma}}=\mathcal{L}_{\widetilde{\sigma}_W\boxtimes \widetilde{\sigma}_V}\neq 0$$

\noindent where in the first equivalence we applied the case $(a,b)=(1,0)$ and in the second equivalence we applied the case $(a,b)=(0,1)$. $\blacksquare$

\end{itemize}

\subsection{Proof of Theorem \ref{theorem 8.2.1}}\label{section 8.5}

\noindent Let $\pi_0\in \Temp(G)$. We already saw in Section \ref{section 8.2} the implication $\mathcal{L}_{\pi_0}\Rightarrow m(\pi_0)\neq 0$. Let us prove the converse. Assume that $m(\pi_0)\neq 0$ and let $\ell\in \Hom_H(\pi_0^\infty,\xi)$ be nonzero. We begin by establishing the following

\vspace{3mm}

\begin{num}
\item\label{eq 8.5.1} For all $e\in \pi_0^\infty$ and all $f\in\mathcal{C}(G(F))$, the integral
$$\displaystyle \int_{G(F)} \ell(\pi_0(g)e) f(g) dg$$
is absolutely convergent.
\end{num}

\vspace{3mm}

\noindent Indeed, this is equivalent to the convergence of

$$\displaystyle \int_{H(F)\backslash G(F)} \lvert \ell(\pi_0(x)e)\rvert \int_{H(F)} \lvert f(hx)\rvert dhdx$$

\noindent By Proposition \ref{proposition 6.7.1}(vi), for all $d>0$, we have an inequality

\begin{align}\label{eq 8.5.2}
\displaystyle \int_{H(F)} \lvert f(hx)\rvert dh\ll \Xi^{H\backslash G}(x) \sigma_{H\backslash G}(x)^{-d}
\end{align}

\noindent for all $x\in H(F)\backslash G(F)$. On the other hand, by Lemma \ref{lemma 8.3.1}(i), there exists $d'>0$ such that we have the inequality

\begin{align}\label{eq 8.5.3}
\lvert \ell(\pi_0(x)e)\rvert\ll \Xi^{H\backslash G}(x)\sigma_{H\backslash G}(x)^{d'}
\end{align}

\noindent for all $x\in H(F)\backslash G(F)$. To conclude, it suffices to combine \ref{eq 8.5.2} and \ref{eq 8.5.3} with Proposition \ref{proposition 6.7.1}(iii).

\vspace{2mm}

\noindent We may compute the integral \ref{eq 8.5.1} in two different ways. First, using the decomposition $\mathcal{C}(G(F))=C_c^\infty(G(F))\ast \mathcal{C}(G(F))$ (cf.\ \ref{eq 2.1.1}), we easily get the equality

\begin{align}\label{eq 8.5.4}
\displaystyle \int_{G(F)} \ell(\pi_0(g)e)f(g)dg=\ell(\pi_0(f)e)
\end{align}

\noindent for all $e\in \pi_0^\infty$ and all $f\in \mathcal{C}(G(F))$. Indeed, by linearity and the aforementioned decomposition, we only need to prove \ref{eq 8.5.4} when $f$ is of the form $f=\varphi\ast f'$ for some $\varphi\in C_c^\infty(G(F))$ and some $f'\in \mathcal{C}(G(F))$. Then, we have

\[\begin{aligned}
\displaystyle \int_{G(F)} \ell(\pi_0(g)e)f(g)dg & =\int_{G(F)}\int_{G(F)} \ell(\pi_0(g)e) \varphi(\gamma)f'(\gamma^{-1}g)d\gamma dg \\
 & =\int_{G(F)}\int_{G(F)} \ell(\pi_0(\gamma g)e) \varphi(\gamma)d\gamma f'(g) dg \\
 & =\int_{G(F)} \ell(\pi_0(\varphi)\pi_0(g)e)f'(g)dg
\end{aligned}\]

\noindent where in the second line we have performed the variable change $g\mapsto \gamma g$. This step requires the switch of the two integrals. This is justified by the fact that the double integral

$$\displaystyle \int_{G(F)}\int_{G(F)} \left\lvert\ell(\pi_0(g)e)\right\rvert \lvert\varphi(\gamma)\rvert \lvert f'(\gamma^{-1}g)\rvert d\gamma dg$$

\noindent is absolutely convergent. But this easily follows from \ref{eq 8.5.1} since $\varphi$ is compactly supported. Now, the vector $\ell\circ \pi_0(\varphi)\in \pi_0^{-\infty}$ actually belongs to $\overline{\pi}_0^\infty$ and by definition of the action of $\mathcal{C}(G(F))$ on $\pi_0^\infty$, we have

\[\begin{aligned}
\displaystyle \int_{G(F)} \ell(\pi_0(\varphi)\pi_0(g)e)f'(g)dg & =\int_{G(F)} f'(g) \left(\pi_0(g)e, \ell\circ \pi_0(\varphi)\right) dg=\left(\pi_0(f')e,\ell\circ \pi_0(\varphi)\right) \\
 & =\ell\left(\pi_0(\varphi)\pi_0(f')e\right)=\ell(\pi_0(f)e)
\end{aligned}\]

\noindent and \ref{eq 8.5.4} follows.

\vspace{3mm}

\noindent On the other hand, we may write

$$\displaystyle \int_{G(F)} \ell(\pi_0(g)e)f(g)dg=\int_{H(F)\backslash G(F)} \ell(\pi_0(x)e)\int_{H(F)} f(hx)\xi(h)dhdx$$

\noindent Hence, by Lemma \ref{lemma 8.2.1}(iv), if $\pi\in \mathcal{X}_{\tempe}(G)\mapsto \pi(f)$ is compactly supported, we have

\begin{align}\label{eq 8.5.5}
\displaystyle \int_{G(F)} \ell(\pi_0(g)e)f(g)dg=\int_{H(F)\backslash G(F)} \ell(\pi_0(x)e)\int_{\mathcal{X}_{\tempe}(G)} \mathcal{L}_\pi(\pi(f)\pi(x^{-1})) \mu(\pi) d\pi dx
\end{align}

\noindent for all $e\in \pi_0^\infty$.

\vspace{2mm}

\noindent Let $T\in C_c^\infty(\mathcal{X}_{\tempe}(G),\mathcal{E}(G))$. Applying \ref{eq 8.5.4} and \ref{eq 8.5.5} to $f=f_T$, we get

\begin{align}\label{eq 8.5.6}
\displaystyle \ell(T_{\pi_0}e)=\int_{H(F)\backslash G(F)} \ell(\pi_0(x)e)\int_{\mathcal{X}_{\tempe}(G)} \mathcal{L}_\pi(T_\pi\pi(x^{-1})) \mu(\pi)d\pi dx
\end{align}

\noindent for all $e\in \pi_0^\infty$.

\vspace{2mm}

\noindent Let $Q=LU_Q$ be a parabolic subgroup of $G$ and $\sigma\in \Pi_2(L)$ such that $\pi_0$ appears as a subrepresentation of $\pi'=i_Q^G(\sigma)$. Set 

$$\mathcal{O}=\{i_Q^G(\sigma_\lambda);\; \lambda\in i\mathcal{A}_L^*\}\subseteq \mathcal{X}_{\tempe}(G)$$

\noindent It is the connected component of $\pi'\in \mathcal{X}_{\tempe}(G)$. Let $e_0\in \pi_0^\infty$ be such that $\ell(e_0)\neq 0$ and let $T_0\in \End(\pi_0)^\infty$ be such that $T_0e_0=e_0$. We may find a section $T^0\in C_c^\infty(\mathcal{X}_{\tempe}(G),\mathcal{E}(G))$, such that

\vspace{2mm}

\begin{itemize}
\renewcommand{\labelitemi}{$\bullet$}
\item $T^0_{\pi_0}=T_0$;

\item $\Supp(T^0)\subseteq \mathcal{O}$.
\end{itemize}

\vspace{2mm}

\noindent Applying \ref{eq 8.5.6} to $e=e_0$ and $T=T^0$, we see that there exists $\lambda\in i\mathcal{A}_L^*$ such that $\mathcal{L}_{\pi'_\lambda}\neq 0$, where $\pi'_\lambda=i_Q^G(\sigma_\lambda)$. Introducing data as in Section \ref{section 8.4}, we may write $\sigma=\sigma_{GL}\boxtimes \widetilde{\sigma}$ where $\sigma_{GL}$ is a tempered representation of a product of general linear groups and $\widetilde{\sigma}$ is a tempered representation of a group $\widetilde{G}(F)$ which is the first component of a GGP triple $(\widetilde{G},\widetilde{H},\widetilde{\xi})$. Twists of $\sigma$ by $i\mathcal{A}_L^*$ leave the component $\widetilde{\sigma}$ unchanged and so, by Proposition \ref{proposition 8.4.1}, we have

$$\mathcal{L}_{\pi'_\lambda}\neq 0\Leftrightarrow \mathcal{L}_{\widetilde{\sigma}}\neq 0\Leftrightarrow \mathcal{L}_{\pi'}\neq 0$$

\noindent Hence, $\mathcal{L}_{\pi'}\neq 0$. Thus, we may find a section $T^1\in \mathcal{C}(\mathcal{X}_{\tempe}(G),\mathcal{E}(G))$ such that $\mathcal{L}_{\pi'}(T^1_{\pi'})\neq 0$. By Lemma \ref{lemma 8.2.1}(i), the function $\pi\in \mathcal{X}_{\tempe}(G)\mapsto \mathcal{L}_\pi(T^1_\pi)$ is smooth. Since $T^0$ is compactly supported, it follows that the section $\pi\in \mathcal{X}_{\tempe}(G)\mapsto T^2_\pi=\mathcal{L}_\pi(T^1_\pi)T^0_\pi$ belongs to $C_c^\infty(\mathcal{X}_{\tempe}(G),\mathcal{E}(G))$. Applying \ref{eq 8.5.6} to $e=e_0$ and $T=T^2$, we get

\begin{align}\label{eq 8.5.7}
\displaystyle \mathcal{L}_{\pi'}(T^1_{\pi'})\ell(e_0)=\int_{H(F)\backslash G(F)} \ell(\pi_0(x)e_0)\int_{\mathcal{X}_{\tempe}(G)} \mathcal{L}_{\pi}(T^1_\pi)\mathcal{L}_{\pi}(T^0_\pi \pi(x^{-1})) \mu(\pi)d\pi dx
\end{align}

\noindent Notice that by the choices of $T^1$ and $e_0$, the left hand side of \ref{eq 8.5.7} is nonzero. On the other hand, by Lemma \ref{lemma 8.2.1}(ii), we have

$$\mathcal{L}_{\pi}(T^1_\pi)\mathcal{L}_{\pi}(T^0_\pi \pi(x^{-1}))=\mathcal{L}_\pi(T^1_\pi L_\pi T^0_\pi\pi(x^{-1}))$$

\noindent for all $\pi\in \mathcal{X}_{\tempe}(G)$ and all $x\in G(F)$. By Lemma \ref{lemma 8.2.1}(iii), the section $\pi\mapsto T^3_\pi=T^1_\pi L_\pi T^0_\pi$ belongs to $C_c^\infty(\mathcal{X}_{\tempe}(G),\mathcal{E}(G))$. Hence, applying \ref{eq 8.5.6} to $e=e_0$ and $T=T^3$, we also get

$$\displaystyle \int_{H(F)\backslash G(F)} \ell(\pi_0(x)e_0)\int_{\mathcal{X}_{\tempe}(G)} \mathcal{L}_{\pi}(T^1_\pi)\mathcal{L}_{\pi}(T^0_\pi \pi(x^{-1})) \mu(\pi)d\pi dx=\ell(T^1_{\pi_0}L_{\pi_0}T^0_{\pi_0}e_0)$$

\noindent By the non-vanishing of the left hand side of \ref{eq 8.5.7}, we deduce that $\ell(T^1_{\pi_0}L_{\pi_0}T^0_{\pi_0}e_0)\neq 0$. Hence, in particular $L_{\pi_0}\neq 0$ which is equivalent to the non-vanishing of $\mathcal{L}_{\pi_0}$ i.e., what we want. This ends the proof of Theorem \ref{theorem 8.2.1}. $\blacksquare$

\subsection{A corollary}\label{section 8.6}

\noindent Let us adopt the notation and hypothesis of Section \ref{section 8.4}. In particular, $Q=LU_Q$ is a parabolic subgroup of $G$ with $L$ decomposing as in \ref{eq 8.4.1}, \ref{eq 8.4.2} and \ref{eq 8.4.3}, $\sigma$ is a tempered representation of $L(F)$ admitting decompositions as in \ref{eq 8.4.4}, \ref{eq 8.4.5} and \ref{eq 8.4.6} and we set $\widetilde{\sigma}=\widetilde{\sigma}_W\boxtimes\widetilde{\sigma}_V$. This is a tempered representation of $\widetilde{G}(F)$ where $\widetilde{G}=U(\widetilde{W})\times U(\widetilde{V})$. Recall that the admissible pair $(\widetilde{W},\widetilde{V})$ (up to permutation) defines a GGP triple $(\widetilde{G},\widetilde{H},\widetilde{\xi})$. Hence, we can define the multiplicity $m(\widetilde{\sigma})$ of $\widetilde{\sigma}$ relative to this GGP triple. We also set, as in Section \ref{section 8.4}, $\pi=i_Q^G(\sigma)$.

\begin{cor}\label{corollary 8.6.1}
\begin{enumerate}[(i)]

\item  Assume that $\sigma$ is irreducible, then we have

$$m(\pi)=m(\widetilde{\sigma})$$

\item Let $\mathcal{K}\subseteq \mathcal{X}_{\tempe}(G)$ be a compact subset. There exists a section $T\in \mathcal{C}(\mathcal{X}_{\tempe}(G),\mathcal{E}(G))$ such that

$$\mathcal{L}_\pi(T_\pi)=m(\pi)$$

\noindent for all $\pi\in \mathcal{K}$ and moreover in this case, the same equality is satisfied for every subrepresentation $\pi$ of some $\pi'\in \mathcal{K}$. 
\end{enumerate}
\end{cor}

\vspace{4mm}

\noindent\ul{Proof}: 

\begin{enumerate}[(i)]

\item By Proposition \ref{proposition 8.4.1} and Theorem \ref{theorem 8.2.1}, we have the chain of equivalences

$$m(\pi)\neq 0\Leftrightarrow \mathcal{L}_\pi\neq 0\Leftrightarrow \mathcal{L}_{\widetilde{\sigma}}\neq 0\Leftrightarrow m(\widetilde{\sigma})\neq 0$$

\noindent Moreover, by Theorem \ref{theorem 6.3.1} the multiplicity $m(\widetilde{\sigma})$ is at most one. Hence it suffices to show that

$$m(\pi)\leqslant 1$$

\noindent Equivalently: there is at most one irreducible subrepresentation of $\pi$ with nonzero multiplicity. Assume this is not the case and let $\pi_1,\pi_2\subset \pi$ be two orthogonal irreducible subrepresentations such that $m(\pi_1)=m(\pi_2)=1$. By Theorem \ref{theorem 8.2.1}, we have $\mathcal{L}_{\pi_1}\neq 0$ and $\mathcal{L}_{\pi_2}\neq 0$. Let $T_1\in \End(\pi_1)^\infty\subset \End(\pi)^\infty$ and $T_2\in \End(\pi_2)^\infty\subset \End(\pi)^\infty$ be such that $\mathcal{L}_\pi(T_1)=\mathcal{L}_{\pi_1}(T_1)\neq 0$ and $\mathcal{L}_\pi(T_2)=\mathcal{L}_{\pi_2}(T_2)\neq 0$. Then, by Lemma \ref{lemma 8.2.1}(ii), we have $\mathcal{L}_\pi(T_1)\mathcal{L}_\pi(T_2)=\mathcal{L}_\pi(T_1L_\pi T_2)$. But obviously $T_1L_\pi T_2=0$ and this yields a contradiction.

\item Using a partition of unity, it clearly suffices to show that for all $\pi'\in \mathcal{X}_{\tempe}(G)$ there exists a section $T\in \mathcal{C}(\mathcal{X}_{\tempe},\mathcal{E}(G))$ such that

$$\mathcal{L}_\pi(T_\pi)=m(\pi)$$

\noindent for $\pi$ in some neighborhood of $\pi'$ in $\mathcal{X}_{\tempe}(G)$. By (i), we know that the function $\pi\in \mathcal{X}_{\tempe}(G)\mapsto m(\pi)$ is locally constant. If $m(\pi')=0$, then there is nothing to prove (just take $T=0$). If $m(\pi')\neq 0$, then by (i) there exists a unique irreducible subrepresentation $\pi_0\subset \pi'$ such that $m(\pi_0)=1$. Then by Theorem \ref{theorem 8.2.1}, we may find $T_0\in \End(\pi_0)^\infty$ such that $\mathcal{L}_{\pi_0}(T_0)\neq 0$. Let $T^0\in \mathcal{C}(\mathcal{X}_{\tempe}(G),\mathcal{E}(G))$ be such that $T^0_{\pi_0}=T_0$. Then we have

$$\mathcal{L}_{\pi'}(T^0_{\pi'})=\mathcal{L}_{\pi_0}(T_0)\neq 0$$

\noindent By Lemma \ref{lemma 8.2.1}(i), the function $\pi\in \mathcal{X}_{\tempe}(G)\mapsto \mathcal{L}_\pi(T^0_\pi)$ is smooth and so we can certainly find a smooth and compactly supported function $\varphi$ on $\mathcal{X}_{\tempe}(G)$ such that $\varphi(\pi)\mathcal{L}_\pi(T^0_\pi)=1$ in some neighborhood of $\pi'$. It then suffices to take $T=\varphi T^0$. $\blacksquare$
\end{enumerate}

\section{The distributions $J$ and $J^{\Lie}$}\label{section 7}

We keep the notation introduced in Chapter \ref{section 6}. The goal of this chapter is to define two functionals $J$ and $J^{\Lie}$ on the spaces of strongly cuspidal functions on the group $G(F)$ and its Lie algebra respectively. In the subsequent Chapters \ref{section 9}, \ref{section 10} and \ref{section 11} we will establish spectral and geometric expansions for these distributions resulting in the local trace formulas alluded to in the introduction. Both the definitions of $J(f)$ and $J^{\Lie}(f)$ involve integrating a certain kernel over the diagonal of $H\backslash G\times H\backslash G$ and the main result of this chapter is that the resulting integrals are absolutely convergent. The convergence of $J(f)$ is proved in Section \ref{section 7.1} using some crucial estimates from Chapter \ref{section 6}. The proof of the convergence of $J^{\Lie}(f)$ is completely similar and thus the result is only stated in Section \ref{section 7.2}.

\subsection{The distribution $J$}\label{section 7.1}

\noindent For all $f\in \mathcal{C}(G(F))$, let us define a function $\gls{Kf}$ on $H(F)\backslash G(F)$ by

$$\displaystyle K(f,x)=\int_{H(F)} f(x^{-1}hx)\xi(h)dh,\;\;\;\; x\in H(F)\backslash G(F)$$

\noindent Notice that by Lemma \ref{lemma 6.5.1}(ii), the above integral is absolutely convergent. The theorem below and Proposition \ref{proposition 6.7.1}(iii) show that the integral

$$\displaystyle \gls{Jf}=\int_{H(F)\backslash G(F)} K(f,x)dx$$

\noindent is absolutely convergent for all $f\in \mathcal{C}_{\scusp}(G(F))$ and defines a continuous linear form

$$\mathcal{C}_{\scusp}(G(F))\to \mathbb{C}$$
$$f\mapsto J(f)$$

\begin{theo}\label{proposition 7.1.1}
\begin{enumerate}[(i)]

\item There exists $d>0$ and a continuous semi-norm $\nu$ on $\mathcal{C}(G(F))$ such that

$$\displaystyle \left\lvert K(f,x)\right\rvert\leqslant \nu(f)\Xi^{H\backslash G}(x)^2\sigma_{H\backslash G}(x)^d$$

\noindent for all $f\in \mathcal{C}(G(F))$ and all $x\in H(F)\backslash G(F)$.

\item For all $d>0$, there exists a continuous semi-norm $\nu_d$ on $\mathcal{C}(G(F))$ such that

$$\lvert K(f,x)\rvert\leqslant \nu_d(f)\Xi^{H\backslash G}(x)^2\sigma_{H\backslash G}(x)^{-d}$$

\noindent for all $x\in H(F)\backslash G(F)$ and all $f\in\mathcal{C}_{\scusp}(G(F))$.
\end{enumerate}
\end{theo}

\vspace{2mm}

\noindent\ul{Proof}: Recall that for all $R>0$, $p_R$ denotes the continuous semi-norm on $\mathcal{C}(G(F))$ given by

$$p_R(f)=\sup_{g\in G(F)}\lvert f(g)\rvert \Xi^G(g)^{-1}\sigma_G(g)^R$$

\noindent for all $f\in \mathcal{C}(G(F))$.

\begin{enumerate}[(i)]

\item Let $d'>0$' Then, we have

$$\displaystyle \left\lvert K(f,x)\right\rvert\leqslant p_{d'}(f)\int_{H(F)} \Xi^G(x^{-1}hx)\sigma(x^{-1}hx)^{-d'} dh$$

\noindent for all $f\in \mathcal{C}(G(F))$ and all $x\in H(F)\backslash G(F)$. By Proposition \ref{proposition 6.7.1}(v), if $d'$ is sufficiently large, there exists $d>0$ such that

$$\displaystyle \int_{H(F)} \Xi^G(x^{-1}hx)\sigma(x^{-1}hx)^{-d'} dh\ll \Xi^{H\backslash G}(x)^2\sigma_{H\backslash G}(x)^d$$

\noindent for all $x\in H(F)\backslash G(F)$. This proves (i) for $\nu$ a scalar multiple of $p_{d'}$ when $d'$ is sufficiently large.

\item We will adopt the notation of Section \ref{section 6.6.2}. That is:

\vspace{3mm}

\begin{itemize}
\renewcommand{\labelitemi}{$\bullet$}

\item $\overline{P}_0$ is a good minimal parabolic subgroup of $G_0$, $M_0\subseteq \overline{P}_0$ a Levi component and $A_0$ the maximal central split subtorus of $M_0$;

\item $\overline{P}$ is the parabolic subgroup opposite to $P$ with respect to $M$ and $\overline{N}$ is its unipotent radical;

\item $\overline{P}_{\mini}=\overline{P}_0T\overline{N}$, $M_{\mini}=M_0T$ and $A_{\mini}=A_0A$ are respectively a good minimal parabolic subgroup of $G$, a Levi component of it and the maximal central split torus of $M_{\mini}$;

\item $A_0^+=\{a\in A_0(F);\; \lvert \alpha(a)\rvert\geqslant 1 \;\forall \alpha\in R(A_0,\overline{P}_0)\}$ and \\
$A_{\mini}^+=\{a\in A_{\mini}(F);\; \lvert \alpha(a)\rvert\geqslant 1 \;\forall \alpha\in R(A_{\mini},\overline{P}_{\mini})\}$;

\item $P_{\mini}$ is the parabolic subgroup opposite to $\overline{P}_{\mini}$ with respect to $M_{\mini}$ (we have $P_{\mini}\subseteq P$);

\item $\Delta$ is the set of simple roots of $A_{\mini}$ in $P_{\mini}$ and $\Delta_P=\Delta\cap R(A_{\mini},N)$ is the subset of simple roots appearing in $\mathfrak{n}=\Lie(N)$.
\end{itemize}

In the $p$-adic case, we fix a compact-open subgroup $K\subseteq G(F)$ and to get uniform notation, we set $\mathcal{C}_K(G(F))=\mathcal{C}(G(F))$ and $\mathcal{C}_{\scusp,K}(G(F))=\mathcal{C}_{\scusp}(G(F))$ when $F=\mathbb{R}$. Clearly, we just need to establish the estimate of the theorem for $f\in \mathcal{C}_{\scusp,K}(G(F))$. The first step is to reduce the range of the inequality to prove to those $x=a\in A_{\mini}^+$. More precisely, we have

\begin{lem}
Point (ii) of the theorem follows from the following estimate:

\begin{num}
\item\label{eq 7.1.4} For all $d>0$, there exists a continuous semi-norm $\nu_{d,K}$ on $\mathcal{C}_K(G(F))$ such that
$$\lvert K(f,a)\rvert\leqslant \nu_{d,K}(f)\Xi^{H\backslash G}(a)^2 \sigma_{H\backslash G}(a)^{-d}$$
for all $a\in A_{\mini}^+$ and all $f\in \mathcal{C}_{\scusp,K}(G(F))$.
\end{num}
\end{lem}

\noindent\ul{Proof}: By Lemma \ref{lemma 6.6.2.1}(ii), there exists a compact subset $\mathcal{K}\subseteq G(F)$ such that

$$G(F)=H(F)A_0^+A(F)\mathcal{K}$$

\noindent Hence, by Proposition \ref{proposition 6.7.1}(i) and the fact that for every semi-norm $\nu$ on $\mathcal{C}(G(F))$, there exists a continuous semi-norm $\nu'$ on $\mathcal{C}(G(F))$ such that $\nu({}^k f)\ll \nu'(f)$ for all $k\in \mathcal{K}$ and all $f\in \mathcal{C}(G(F))$, it certainly suffices to prove the estimate of the theorem for $x=a\in A_0^+A(F)$ i.e. we have a reduction to the following statement:

\vspace{3mm}

\begin{num}
\item\label{eq 7.1.1} For all $d>0$, there exists a continuous semi-norm $\nu_{d,K}$ on $\mathcal{C}_K(G(F))$ such that
$$\displaystyle \lvert K(f,a)\rvert\leqslant \nu_{d,K}(f)\Xi^{H\backslash G}(a)^2 \sigma_{H\backslash G}(a)^{-d}$$
for all $a\in A_0^+A(F)$ and all $f\in \mathcal{C}_{\scusp,K}(G(F))$.
\end{num}

\vspace{3mm}

Set

$$\displaystyle \mathcal{A}_{\mini}^{a +}=\{a\in A_0^+A(F);\; \lvert \alpha(a)\rvert\leqslant \sigma_{H\backslash G}(a)\; \forall \alpha\in\Delta_P\}$$

\noindent where the exponent ``a" stand for ``almost" since by Lemma \ref{lemma 6.6.2.1}(i), elements of $\mathcal{A}_{\mini}^{a +}$ are ``almost" in $A_{\mini}^+$. We now claim:

\vspace{3mm}

\begin{num}
\item\label{eq 7.1.2} In the Archimedean case, for all $d>0$ there exists a continuous semi-norm $\nu_{d}$ on $\mathcal{C}(G(F))$ such that
$$\displaystyle \lvert K(f,a)\rvert\leqslant \nu_d(f) \Xi^{H\backslash G}(a)^2 \sigma_{H\backslash G}(a)^{-d}$$
for all $a\in A_0^+A(F)\backslash \mathcal{A}_{\mini}^{a +}$ and all $f\in \mathcal{C}(G(F))$.
\end{num}

\vspace{3mm}

\begin{num}
\item\label{eq 7.1.2bis} In the non-Archimedean case, there exists a compact subset $C_A\subset A(F)$ such that $$\displaystyle K(f,a)=0$$
for all $a\in A_0^+A(F)\backslash A_{\mini}^+C_A$ and all $f\in \mathcal{C}_K(G(F))$.
\end{num}

\vspace{3mm}

\noindent First we prove \ref{eq 7.1.2}. Let $a\in A_0^+A(F)\backslash \mathcal{A}_{\mini}^{a +}$. Then, there exists $\alpha\in \Delta_P$ such that $\lvert \alpha(a)\rvert>\sigma_{H\backslash G}(a)$. Since $\Delta_P$ is finite, we may as well fix $\alpha\in \Delta_P$ and prove the estimate \ref{eq 7.1.2} only for those $a\in A_0^+A(F)$ such that $\lvert \alpha(a)\rvert>\sigma_{H\backslash G}(a)$. So, let us fix $\alpha\in \Delta_P$. By Lemma \ref{lemma 6.6.2.1}(iii), we know that $\xi$ is nontrivial on $\mathfrak{n}_\alpha(F)$. Choose $X\in \mathfrak{n}_\alpha(F)$ so that $\xi(X)\neq 1$. Denote by $d\xi:\mathfrak{n}(F)\to \mathbb{C}$ the differential of $\xi$ at the origin. Since $\xi(X)\neq 1$ and $\xi$ is a character, we have $d\xi(X)\neq 0$. By integration by part, for all $f\in \mathcal{C}(G(F))$ all $a\in A_{\mini}(F)$ and all positive integers $N$, we have

\[\begin{aligned}
\displaystyle d\xi(X)^NK(f,a) & =\int_{H(F)} {}^a f(h) \left(L(X^N)\xi\right)(h) dh \\
 & =(-1)^N\int_{H(F)}\left(L(X^N){}^a f\right)(h)\xi(h) dh\\
 & =(-1)^N\int_{H(F)} {}^a\left(L(a^{-1}X^Na)f\right)(h)\xi(h)dh \\
 & =(-1)^N \alpha(a)^{-N}\int_{H(F)} \left(L(X^N)f\right)(a^{-1}ha)\xi(h) dh\\
 & =(-1)^N\alpha(a)^{-N} K(L(X^N)f,a)
\end{aligned}\]

\noindent This implies in particular that for every positive integer $N$, we have

$$\displaystyle \lvert K(f,a)\rvert \ll\sigma_{H\backslash G}(a)^{-N} \lvert K(L(X^N)f,a)\rvert$$

\noindent for all $f\in \mathcal{C}(G(F))$ and all $a\in A_{\mini}(F)$ such that $\lvert \alpha(a)\rvert>\sigma_{H\backslash G}(a)$. Together with (i), this implies the desired inequality.

\vspace{2mm}

\noindent We now prove \ref{eq 7.1.2bis}. Let $C>0$. Certainly, there exists a compact subset $C_A\subset A(F)$ such that 
$A_{\mini}^+C_A$ contains the set of $a\in A_{\mini}^+A(F)$ such that $\lvert \alpha(a)\rvert\leqslant C$ for all $\alpha\in \Delta_P$. Hence, fixing $\alpha\in \Delta_P$, it suffices to show that for $C$ sufficiently large we have $K(f,a)=0$ for all $f\in \mathcal{C}_K(G(F))$ and all $a\in A_{\mini}(F)$ with $\lvert \alpha(a)\rvert>C$. Set $K_N=K\cap N(F)$, $L_N=\log(K_N)$ and $L_{\alpha}=L_N\cap \mathfrak{n}_\alpha(F)$. Then $L_\alpha$ is a lattice of $\mathfrak{n}_\alpha(F)$ and there exists a constant $C>0$ such that $\lambda^{-1}X\in L_\alpha$ for all $\lambda\in F^\times$ satisfying $\lvert \lambda\rvert>C$. Hence, for all $f\in \mathcal{C}_K(G(F))$ and all $a\in A_{\mini}(F)$ with $\lvert \alpha(a)\rvert>C$, we have

\[\begin{aligned}
\displaystyle \xi(X)K(f,a) & =\int_{H(F)} f(a^{-1}ha)\xi(he^X)dh=\int_{H(F)} f\left(a^{-1}he^{-X}a\right)\xi(h) dh\\
 & =\int_{H(F)} f\left(a^{-1}hae^{-X/\alpha(a)}\right)\xi(h) dh=\int_{H(F)} f\left(a^{-1}ha\right)\xi(h) dh=K(f,a)
\end{aligned}\]

\noindent Since $\xi(X)\neq 1$, this implies $K(f,a)=0$ and \ref{eq 7.1.2bis} follows.

\vspace{4mm}

By \ref{eq 7.1.2} and \ref{eq 7.1.2bis}, it is thus sufficient to establish the estimate \ref{eq 7.1.1} only for those $x=a\in \mathcal{A}_{\mini}^{a +}$ in the Archimedean case and only for those $x=a\in A_{\mini}^+C_A$ for a certain compact subset $C_A\subset A(F)$ in the non-Archimedean case. The fact that we can further reduce to the case where $x=a\in A_{\mini}^+$, up to changing the level in the $p$-adic case, is now a consequence of the next lemma, whose straightforward proof is left to the reader.

\begin{lem}
\begin{enumerate}[(i)]
\item For all $a_+,a_-\in A_{\mini}(F)$ and $f\in \mathcal{C}_{\scusp}(G(F))$, we have $K(f,a_+a_-)=K({}^{a_-}f,a_+)$ where ${}^{a_-}f\in \mathcal{C}_{\scusp}(G(F))$ is defined by $({}^{a_-}f)(g)=f(a_-^{-1}ga_-)$ for all $g\in G(F)$.

\item In the Archimedean case, any $a\in\mathcal{A}_{\mini}^{a +}$ can be written $a=a_+a_-$ where $a_+\in A_{\mini}^+$ and $a_-\in A(F)$ satisfy an inequality
$$\displaystyle \sigma(a_-)\ll \log\left(1+\sigma_{H\backslash G}(a_+)\right)$$
Moreover, fixing such a decomposition for all $a\in\mathcal{A}_{\mini}^{a +}$, there exist $c>0$ such that for every continuous semi-norm $\nu$ on $\mathcal{C}(G(F))$ there exists a continuous semi-norm $\nu'$ on $\mathcal{C}(G(F))$ satisfying
$$\nu({}^{a_-} f)\leqslant \nu'(f)\sigma_{H\backslash G}(a_+)^c,\;\;\; \sigma_{H\backslash G}(a)\ll\sigma_{H\backslash G}(a_+) \mbox{ and } \Xi^{H\backslash G}(a_+)\ll \Xi^{H\backslash G}(a)\sigma_{H\backslash G}(a_+)^c$$
for all $a\in\mathcal{A}_{\mini}^{a +}$ and all $f\in \mathcal{C}(G(F))$.

\item In the non-Archimedean case, let $C_A\subset A(F)$ be a compact. Then, any $a\in A_{\mini}^+C_A$ can be written $a=a_+a_-$ where $a_+\in A_{\mini}^+$ and $a_-\in C_A$. Moreover, fixing such a decomposition for all $a\in A_{\mini}^+C_A$, there exists a compact-open subgroup $K'\subset G(F)$ and for all continuous semi-norm $\nu$ on $\mathcal{C}_K(G(F))$, there exists a continuous semi-norm $\nu'$ on $\mathcal{C}_{K'}(G(F))$ such that
$$\displaystyle {}^{a_-}f\in \mathcal{C}_{K'}(G(F)),\;\;\; \nu({}^{a_-} f)\leqslant \nu'(f),\;\;\; \sigma_{H\backslash G}(a)\ll\sigma_{H\backslash G}(a_+) \mbox{ and } \Xi^{H\backslash G}(a_+)\ll \Xi^{H\backslash G}(a)$$
for all $a\in A_{\mini}^{+}$ and all $f\in \mathcal{C}_K(G(F))$.
\end{enumerate}
\end{lem}

$\blacksquare$

\vspace{2mm}

We are thus left with proving \ref{eq 7.1.4}. In order not to have to keep track of semi-norms, we now make the following useful remark: by (i), for all $a\in A_{\mini}^+$ the linear form $f\in \mathcal{C}(G(F))\mapsto K(f,a)$ is continuous and therefore, by the uniform boundedness principle, in order to prove \ref{eq 7.1.4} we just need to establish that for any fixed $f\in \mathcal{C}_{K,\scusp}(G(F))$ and $d>0$, we have

\begin{align}\label{eq 7.1.0}
\displaystyle \lvert K(f,a)\rvert\ll \Xi^{H\backslash G}(a)^2 \sigma_{H\backslash G}(a)^{-d}
\end{align}

\noindent for all $a\in A_{\mini}^+$.

\vspace{3mm}

For every maximal proper parabolic subgroup $\overline{Q}$ of $G$ containing $\overline{P}_{\mini}$ with unipotent radical $U_{\overline{Q}}$ and any $\delta>0$, we set

$$A_{\mini}^{\overline{Q},+}(\delta)=\{a\in A_{\mini}^+;\; \lvert \alpha(a)\rvert\geqslant e^{\delta \sigma(a)}\; \forall \alpha\in R(A_{\mini},U_{\overline{Q}})\}$$

\noindent Obviously, we can choose $\delta>0$ such that the complement of

$$\displaystyle \bigcup_{\overline{Q}} A_{\mini}^{\overline{Q},+}(\delta),$$

\noindent the union being taken over all maximal proper parabolic subgroups $\overline{Q}\supseteq \overline{P}_{\mini}$, in $A_{\mini}^+$ is relatively compact. We fix such a $\delta>0$ henceforth and let $\overline{Q}$ be a maximal proper parabolic subgroup containing $\overline{P}_{\mini}$. By the previous decomposition, we only need to prove the estimate \ref{eq 7.1.0} for $a\in A_{\mini}^{\overline{Q},+}(\delta)$.

Let $L$ be the unique Levi component of $\overline{Q}$ containing $A_{\mini}$. Set $H_{\overline{Q}}=H\cap \overline{Q}$ and let $H_L$ be the image of $H_{\overline{Q}}$ by the projection $\overline{Q}\twoheadrightarrow L$. Since $\overline{Q}$ is a good parabolic subgroup, by Proposition \ref{proposition 6.8.1}(i) this projection induces an isomorphism $H_{\overline{Q}}\simeq H_L$. We define a character $\xi_L$ on $H_L(F)$ by setting $\xi_L(h_L)=\xi(h_{\overline{Q}})$ for all $h_{\overline{Q}}\in H_{\overline{Q}}(F)$ with image $h_L\in H_L(F)$. Let $Q=LU_Q$ be the parabolic subgroup opposite to $\overline{Q}$ with respect to $L$ and set $H^Q=H_L\ltimes U_Q$. By Proposition \ref{proposition 6.8.1}(ii), $H^Q(F)$ is a unimodular group and we fix a Haar measure $dh^Q$ on it. Let $\xi^Q$ be the character of $H^Q(F)$ defined by $\xi^Q(h_Lu_Q)=\xi_L(h_L)$ for all $h_L\in H_L(F)$ and $u_Q\in U_Q(F)$. For all $f\in \mathcal{C}(G(F))$ and $a\in A_{\mini}^+$, we define
$$\displaystyle K^Q(f,a)=\int_{H^Q(F)} f(a^{-1}h^Qa)\xi^Q(h^Q)dh^Q$$
This expression is absolutely convergent by Proposition \ref{proposition 6.8.1}(iv). Moreover, as $U_Q\subset H^Q\subset Q$, it is clear that for every strongly cuspidal function $f\in \mathcal{C}_{\scusp}(G(F))$ we have $K^Q(f,a)=0$ for all $a\in A_{\mini}^+$. Hence, the following proposition implies \ref{eq 7.1.0} and thus will end the proof of the theorem.

\begin{prop}
There exists a constant $c>0$ (depending on the choices of Haar measures) so that for any $f\in \mathcal{C}(G(F))$ and any $d>0$ we have
$$\displaystyle \left\lvert K(f,a)-cK^Q(f,a)\right\rvert\ll \Xi^{H\backslash G}(a)^2\sigma_{H\backslash G}(a)^{-d}$$
for all $a\in A^{\overline{Q},+}_{\mini}(\delta)$.
\end{prop}

\noindent\ul{Proof}: Fix $f\in \mathcal{C}(G(F))$ and $d>0$. Let $l:\overline{Q}\to L$ be the unique regular map such that $\overline{q}l(\overline{q})^{-1}\in U_{\overline{Q}}$ for all $\overline{q}\in \overline{Q}$. Then $h_{\overline{Q}}\in H_{\overline{Q}}\mapsto h_L:=l(h_{\overline{Q}})$ is precisely the aforementioned isomorphism $H_{\overline{Q}}\simeq H_L$. To simplify some arguments, we will assume, as we may, that $\sigma(l(h_{\overline{Q}}))=\sigma(h_{\overline{Q}})$ for all $h_{\overline{Q}}\in H_{\overline{Q}}(F)$. This in particular implies that for any $C>0$ we have
$$\displaystyle l(H_{\overline{Q}}[<C])=H_L[<C]$$
We also fix left Haar measures $d_Lh_{\overline{Q}}$, $d_Lh_L$ on $H_{\overline{Q}}(F)$ and $H_L(F)$ respectively which correspond via the isomorphism $H_{\overline{Q}}\simeq H_L$ and we equip $U_Q(F)$ with the unique Haar measure so that
$$\displaystyle \int_{H^Q(F)}\varphi(h^Q)dh^Q=\int_{H_L(F)}\int_{U_Q(F)} \varphi(h_Lu_Q)du_Qd_Lh_L$$
for all $\varphi\in L^1(H^Q(F))$.

Let $U_{\mini}$ be the unipotent radical of $P_{\mini}$ and set

$$\accentset{\circ}{H}=H(F)\cap \overline{P}_{\mini}(F)U_{\mini}(F)$$

\noindent Then $\accentset{\circ}{H}$ is an open subset of $H(F)$ containing the identity. Let

$$u:\accentset{\circ}{H}\to U_{\mini}(F)$$

\noindent be the $F$-analytic map sending $h\in \accentset{\circ}{H}$ to the unique element $u(h)\in U_{\mini}(F)$ such that $hu(h)^{-1}\in \overline{P}_{\mini}(F)$. We have

\vspace{3mm}

\begin{num}
\item\label{eq 7.1.6} The map $u$ is submersive at the identity.
\end{num}

\vspace{3mm}

\noindent Indeed, the differential of $u$ at $1$ is given by $d_1 u(X)=p_{\mathfrak{u}_{\mini}}(X)$ for all $X\in \mathfrak{h}(F)$, where $p_{\mathfrak{u}_{\mini}}$ denotes the linear projection of $\mathfrak{g}$ onto $\mathfrak{u}_{\mini}$ relative to the decomposition $\mathfrak{g}=\overline{\mathfrak{p}}_{\mini}\oplus \mathfrak{u}_{\mini}$, and $p_{\mathfrak{u}_{\mini}}(\mathfrak{h})=\mathfrak{u}_{\mini}$ since $\overline{P}_{\mini}$ is a good parabolic subgroup.

\vspace{2mm}

\noindent Because of \ref{eq 7.1.6}, we can find a relatively compact open neighborhood $\mathcal{U}_{\mini}$ of $1$ in $U_{\mini}(F)$ and an $F$-analytic section

$$h:\mathcal{U}_{\mini}\to \accentset{\circ}{H}$$
$$u\mapsto h(u)$$

\noindent to the map $u(.)$ over $\mathcal{U}_{\mini}$ such that $h(1)=1$. Set $\mathcal{U}_Q=\mathcal{U}_{\mini}\cap U_Q(F)$ and $\mathcal{H}=H_{\overline{Q}}(F)h\left(\mathcal{U}_Q\right)$. We will need the following fact:

\vspace{3mm}

\begin{num}
\item\label{eq 7.1.7} The map $\iota: H_{\overline{Q}}(F)\times \mathcal{U}_Q\to H(F)$, $(h_{\overline{Q}},u_Q)\mapsto h_{\overline{Q}}h(u_Q)$, is an $F$-analytic open embedding with image $\mathcal{H}$ and there exists a smooth function $j\in C^\infty(\mathcal{U}_Q)$ such that
$$\displaystyle \int_{\mathcal{H}}\varphi(h)dh=\int_{H_{\overline{Q}}(F)}\int_{\mathcal{U}_Q} \varphi(h_{\overline{Q}}h(u_Q)) j(u_Q)du_Qd_Lh_{\overline{Q}}$$
for all $\varphi \in L^1(\mathcal{H})$.
\end{num}

\vspace{3mm}

\noindent Indeed, we have the following Cartesian diagram
$$\xymatrix{
H_{\overline{Q}}(F)\times \mathcal{U}_Q \ar[d] \ar[r]^{\iota} & H(F)\ar[d]\\
\mathcal{U}_Q  \ar[r] & \overline{Q}(F)\backslash G(F)}$$
where the left arrow is the natural projection $H_{\overline{Q}}(F)\times \mathcal{U}_Q\twoheadrightarrow \mathcal{U}_Q$ and the right and bottom arrows are the restrictions to $H(F)$ and $\mathcal{U}_Q$ respectively of the natural projection $G(F)\twoheadrightarrow \overline{Q}(F)\backslash G(F)$. Since the bottom arrow is an open embedding, so is $\iota$. Let $j$ be the absolute value of the Jacobian of $\iota$. This is a smooth function on $H_{\overline{Q}}(F)\times \mathcal{U}_Q$ which is left invariant by $H_{\overline{Q}}(F)$ as $\iota$ is clearly $H_{\overline{Q}}(F)$-equivariant on the left. The claim \ref{eq 7.1.7} follows.

\vspace{2mm}

\noindent Fix $\epsilon>0$ that we will assume sufficiently small in what follows. By Proposition \ref{proposition 1.3.1}(ii), for $\epsilon$ small enough we have
$$aU_Q\left[<\epsilon \sigma(a)\right]a^{-1}\subseteq \mathcal{U}_Q$$
for all $a\in A_{\mini}^{\overline{Q},+}(\delta)$. Then, we set
$$\displaystyle H^{<\epsilon,a}=H_{\overline{Q}}\left[<\epsilon \sigma(a)\right]\; h\left(a U_Q\left[<\epsilon \sigma(a)\right]a^{-1}\right),$$
$$\displaystyle H^{Q,<\epsilon,a}=H_{L}\left[<\epsilon \sigma(a)\right]a U_Q\left[<\epsilon \sigma(a)\right]a^{-1}$$
and we introduce the following expressions
$$\displaystyle K^{<\epsilon}(f,a)=\int_{H^{<\epsilon,a}} f(a^{-1}ha)\xi(h)dh$$
$$\displaystyle K^{Q,<\epsilon}(f,a)=\int_{H^{Q,<\epsilon,a}} f(a^{-1}h^Qa)\xi^Q(h^Q)dh^Q$$
for all $a\in A_{\mini}^{\overline{Q},+}(\delta)$. Set $c=j(1)$ (where the function $j(.)$ is the one appearing in \ref{eq 7.1.7}) and let $0<\delta_0<\delta$. As $\sigma_{H\backslash G}(a)\ll \sigma(a)$ for all $a\in A_{\mini}(F)$, the proposition is obviously a consequence of the following estimates:

\begin{align}\label{eq 7.1.9}
\displaystyle \left\lvert K(f,a)-K^{<\epsilon}(f,a)\right\rvert\ll \Xi^{H\backslash G}(a)^2\sigma_{H\backslash G}(a)^{-d}
\end{align}

\begin{align}\label{eq 7.1.9bis}
\displaystyle \left\lvert K^Q(f,a)-K^{Q,<\epsilon}(f,a)\right\rvert\ll \Xi^{H\backslash G}(a)^2\sigma_{H\backslash G}(a)^{-d}
\end{align}

\begin{align}\label{eq 7.1.10}
\displaystyle \left\lvert K^{<\epsilon}(f,a)-cK^{Q,<\epsilon}(f,a)\right\rvert\ll \Xi^{H\backslash G}(a)^2e^{-\delta_0\sigma(a)}
\end{align}

\noindent for all $a\in A_{\mini}^{\overline{Q},+}(\delta)$.

\vspace{2mm}

\noindent We start by proving \ref{eq 7.1.9} and \ref{eq 7.1.9bis}. For this we need the following

\vspace{3mm}

\begin{num}
\item\label{eq 7.1.11} We have
$$\displaystyle \sigma(a)\ll \sigma\left(a^{-1}h^Qa\right) \mbox{ and } \sigma(a)\ll \sigma\left(a^{-1}ha\right)$$
for all $a\in A_{\mini}^{\overline{Q},+}(\delta)$, all $h^Q\in H^Q(F)\setminus H^{Q,<\epsilon,a}$ and all $h\in H(F)\setminus H^{<\epsilon,a}$.
\end{num}

\vspace{3mm}

\noindent The first inequality follows from Proposition \ref{proposition 1.3.1}(i) and Proposition \ref{proposition 6.8.1}(v). For the second inequality, combining Proposition \ref{proposition 1.3.1}(i) with Proposition \ref{proposition 6.4.1}(iii), we see that it suffices to show the existence of $\epsilon'>0$ such that

\begin{align}\label{eq 7.1.12}
\displaystyle H(F)\cap \left(\overline{Q}\left[<\epsilon'\sigma(a)\right]\; aU_Q\left[<\epsilon' \sigma(a)\right]a^{-1}\right)\subseteq H^{<\epsilon,a}
\end{align}

\noindent for all $a\in A_{\mini}^{\overline{Q},+}(\delta)$. Fix $\epsilon'>0$ and let us show that the above inclusion is satisfied if $\epsilon'$ is sufficiently small for any $a\in A_{\mini}^{\overline{Q},+}(\delta)$. If $\sigma(a)\leqslant \epsilon'^{-1}$, then the left hand side of \ref{eq 7.1.12} is empty and there is nothing to prove. Hence, we may assume that $\sigma(a)>\epsilon'^{-1}$. Let $\displaystyle h\in H(F)\cap \left(\overline{Q}\left[<\epsilon'\sigma(a)\right]\; aU_Q\left[<\epsilon' \sigma(a)\right]a^{-1}\right)$ and assume $\epsilon'<\epsilon$. Then, there exists $\displaystyle u\in aU_Q\left[<\epsilon' \sigma(a)\right]a^{-1}$ such that $hh(u)^{-1}\in H_{\overline{Q}}(F)$ and we only need to show that for $\epsilon'$ sufficiently small we have $\sigma(hh(u)^{-1})<\epsilon \sigma(a)$. By Lemma \ref{proposition 1.3.1}(ii), for $\epsilon'$ sufficiently small, the set $\displaystyle aU_Q\left[<\epsilon' \sigma(a)\right]a^{-1}$ remains in a fixed compact subset of $\mathcal{U}_{\mini}$ independent of $a\in A_{\mini}^{\overline{Q},+}(\delta)$. This immediately implies $\sigma\left(hh(u)^{-1}\right)\ll \sigma(h)+\sigma(h(u))\ll \epsilon'\sigma(a)$ (as we are assuming $\sigma(a)>\epsilon'^{-1}$) where the implicit constant depends only on this fixed compact set. This proves \ref{eq 7.1.12} for $\epsilon'$ sufficiently small and ends the proof of \ref{eq 7.1.11}.

\vspace{2mm}

\noindent We are now in position to prove \ref{eq 7.1.9} and \ref{eq 7.1.9bis}. Indeed, by \ref{eq 7.1.11}, for all $d'>0$ we have

$$\displaystyle \left\lvert K(f,a)- K^{<\epsilon}(f,a)\right\rvert\ll \sigma(a)^{-d'/2} \int_{H(F)} \Xi^G(a^{-1}ha)\sigma(a^{-1}ha)^{-d'/2} dh$$
and
$$\displaystyle \left\lvert K^Q(f,a)- K^{Q,<\epsilon}(f,a)\right\rvert\ll \sigma(a)^{-d'/2} \int_{H^Q(F)} \Xi^G(a^{-1}h^Qa)\sigma(a^{-1}h^Qa)^{-d'/2} dh$$

\noindent for all $a\in A_{\mini}^{\overline{Q},+}(\delta)$. By Proposition \ref{proposition 6.7.1}(v), Proposition \ref{proposition 6.8.1}(vi) and the inequality $\sigma_{H\backslash G}(g)\ll \sigma(g)$ (for all $g\in G(F)$), for $d'$ sufficiently large the two last expressions above are essentially bounded by
$$\displaystyle \sigma_{H\backslash G}(a)^{-d} \Xi^{H\backslash G}(a)^2$$

\noindent for all $a\in A_{\mini}^+$ and this proves \ref{eq 7.1.9} and \ref{eq 7.1.9bis}.

\vspace{2mm}

We now go on to the proof of \ref{eq 7.1.10}. By \ref{eq 7.1.7}, we have

\begin{align}\label{eq 7.1.13}
\displaystyle  K^{<\epsilon}(f,a)= \int_{H_{\overline{Q}}\left[<\epsilon \sigma(a)\right]} \int_{aU_Q\left[<\epsilon \sigma(a) \right]a^{-1}} f\left(a^{-1} h_{\overline{Q}} h(u_Q) a\right) \xi\left(h_{\overline{Q}} h(u_Q)\right) j(u_Q) du_Q dh_{\overline{Q}}
\end{align}

\noindent and

\begin{align}\label{eq 7.1.13bis}
\displaystyle  K^{Q,<\epsilon}(f,a)= \int_{H_L\left[<\epsilon \sigma(a)\right]} \int_{aU_Q\left[<\epsilon \sigma(a) \right]a^{-1}} f\left(a^{-1} h_L u_Q a\right) \xi_L(h_{L}) du_Q dh_{L}
\end{align}

\noindent for all $a\in A_{\mini}^{\overline{Q},+}(\delta)$. Let $\delta_0<\delta'<\delta$ and $d'>0$ and assume for one moment the following estimate:

\vspace{3mm}

\begin{num}
\item\label{eq 7.1.14}  If $\epsilon$ is sufficiently small, we have
\[\begin{aligned}
\displaystyle \big\lvert f\left(a^{-1} h_{\overline{Q}} h(u_Q) a\right) \xi\left(h_{\overline{Q}} h(u_Q)\right) j(u_Q)-c & f\left(a^{-1} h^Q a\right) \xi(h_{\overline{Q}})\big\rvert \\
 & \ll \Xi^G\left(a^{-1}h^Qa\right)\sigma\left(a^{-1}h^Qa\right)^{-d'} e^{-\delta'\sigma(a)}
\end{aligned}\]
for all $a\in A_{\mini}^{\overline{Q},+}(\delta)$, all $u_Q\in aU_Q\left[<\epsilon\sigma(a)\right]a^{-1}$ and all $h_{\overline{Q}}\in H_{\overline{Q}}\left[<\epsilon \sigma(a)\right]$ where we have set $h^Q=l(h_{\overline{Q}}) u_Q$ (and where we recall that $c=j(1)$).
\end{num}

\vspace{3mm}

\noindent Then, as $h_{\overline{Q}}\mapsto l(h_{\overline{Q}})$ is an isomorphism $H_{\overline{Q}}(F)\simeq H_L(F)$ preserving the measures, sending $\xi_{\mid H_{\overline{Q}}(F)}$ to $\xi_L$ and $H_{\overline{Q}}[<C]$ to $H_L[<C]$ for any $C>0$, by \ref{eq 7.1.13} and \ref{eq 7.1.13bis} we get

\[\begin{aligned}
\displaystyle & \left\lvert K^{<\epsilon}(f,a)-cK^{Q,<\epsilon}(f,a)\right\rvert\leqslant \\
& \int_{H_{\overline{Q}}\left[<\epsilon \sigma(a)\right]} \int_{aU_Q\left[<\epsilon \sigma(a) \right]a^{-1}} \left\lvert f\left(a^{-1} h_{\overline{Q}} h(u_Q) a\right) \xi\left(h_{\overline{Q}} h(u_Q)\right) j(u_Q)-cf\left(a^{-1} l(h_{\overline{Q}}) u_Q a\right) \xi(h_{\overline{Q}})\right\rvert du_Q dh_{\overline{Q}} \\
 & \ll e^{-\delta'\sigma(a)}\int_{H^{Q,<\epsilon,a}} \Xi^G\left(a^{-1}h^Qa\right)\sigma\left(a^{-1}h^Qa\right)^{-d'}dh^Q\leqslant e^{-\delta'\sigma(a)}\int_{H^Q(F)} \Xi^G\left(a^{-1}h^Qa\right)\sigma\left(a^{-1}h^Qa\right)^{-d'}dh^Q
\end{aligned}\]

\noindent for all $a\in A_{\mini}^{\overline{Q},+}(\delta)$. By Proposition \ref{proposition 6.8.1}(vi), for $d'$ sufficiently large the last expression above is essentially bounded by $e^{-\delta'\sigma(a)}\Xi^{H\backslash G}(a)^2\sigma_{H\backslash G}(a)^{d_0}$ for some $d_0>0$. As $\sigma_{H\backslash G}(g)\ll \sigma(g)$ for all $g\in G$, we have $e^{-\delta'\sigma(a)}\sigma_{H\backslash G}(a)^{d_0}\ll e^{-\delta_0\sigma(a)}$ for all $a\in A_{\mini}^{\overline{Q},+}(\delta)$ and the estimate \ref{eq 7.1.10} follows.

\vspace{2mm}

Thus, it only remains to establish \ref{eq 7.1.14}. As $f$ is a Harish-Chandra Schwartz function, this estimate is itself a consequence of the two following ones:

\vspace{3mm}

\begin{num}
\item\label{eq 7.1.15} If $\epsilon$ is sufficiently small, we have
$$\displaystyle \left\lvert \xi(h(u_Q))j(u_Q)-j(1)\right\rvert\ll e^{-\delta' \sigma(a)}$$
for all $a\in A_{\mini}^{\overline{Q},+}(\delta)$ and all $u_Q\in aU_Q\left[<\epsilon\sigma(a)\right]a^{-1}$.
\end{num}

\vspace{3mm}

\begin{num}
\item\label{eq 7.1.16} If $\epsilon$ is sufficiently small, we have
$$\displaystyle \big\lvert f\left(a^{-1} h_{\overline{Q}} h(u_Q) a\right)-f\left(a^{-1} h^Q a\right)\big\rvert \ll \Xi^G\left(a^{-1}h^Qa\right)\sigma\left(a^{-1}h^Qa\right)^{-d'} e^{-\delta'\sigma(a)}$$
for all $a\in A_{\mini}^{\overline{Q},+}(\delta)$, all $u_Q\in aU_Q\left[<\epsilon\sigma(a)\right]a^{-1}$ and all $h_{\overline{Q}}\in H_{\overline{Q}}\left[<\epsilon \sigma(a)\right]$ where as before we have set $h^Q=l(h_{\overline{Q}}) u_Q$.
\end{num}

\vspace{3mm}

Actually, \ref{eq 7.1.15} is an easy consequence of Lemma \ref{proposition 1.3.1}(ii) and the fact that the function $u_Q\mapsto \xi(h(u_Q))j(u_Q)$ is smooth in a neighborhood of $1$ (so that in particular in the $p$-adic case, the left hand side of \ref{eq 7.1.15} is identically $0$ for $\epsilon$ small enough). We now prove \ref{eq 7.1.16}. By Lemma \ref{lemma 1.5.1bis} in the Archimedean case and the smoothness of $f$ in the non-Archimedean case, it is sufficient to show:

\vspace{3mm}

\begin{num}
\item\label{eq 7.1.17} If $\epsilon$ is sufficiently small, for all $a\in A_{\mini}^{\overline{Q},+}(\delta)$, all $u_Q\in aU_Q\left[<\epsilon\sigma(a)\right]a^{-1}$ and $h_{\overline{Q}}\in H_{\overline{Q}}\left[<\epsilon \sigma(a)\right]$, there exists $X,Y\in B(0,e^{-\delta'\sigma(a)})$ such that
$$\displaystyle a^{-1}h_{\overline{Q}}h(u_Q)a=e^X a^{-1}l(h_{\overline{Q}})u_Qa e^Y$$
\end{num}

\vspace{3mm}

\noindent Let $a\in A_{\mini}^{\overline{Q},+}(\delta)$, $u_Q\in aU_Q\left[<\epsilon\sigma(a)\right]a^{-1}$ and $h_{\overline{Q}}\in H_{\overline{Q}}\left[<\epsilon \sigma(a)\right]$. Then we have
$$\displaystyle a^{-1}h_{\overline{Q}}h(u_Q)a=\gamma_1 a^{-1}l(h_{\overline{Q}})u_Qa \gamma_2$$
where $\gamma_1=a^{-1}h_{\overline{Q}}l(h_{\overline{Q}})^{-1}a$ and $\gamma_2=a^{-1}u_Q^{-1}h(u_Q)a$.

By Lemma \ref{proposition 1.3.1}(ii), there exists $\epsilon'>0$ such that $a^{-1}U_{\overline{Q}}[<\epsilon' \sigma(a)]a\subset \exp(B(0,e^{-\delta'\sigma(a)}))$. Moreover, for $\epsilon$ sufficiently small we have $h_{\overline{Q}}l(h_{\overline{Q}})^{-1}\in U_{\overline{Q}}[<\epsilon' \sigma(a)]$ and thus $\gamma_1\in \exp(B(0,e^{-\delta'\sigma(a)}))$.

Choose $\delta'<\delta''<\delta$. By definition of $h(.)$, the map $u_Q\mapsto h(u_Q)u_Q^{-1}$ is an $F$-analytic map sending $1$ to itself and taking values into $\overline{P}_{\mini}(F)$. By Lemma \ref{proposition 1.3.1}(ii) and since $A_{\mini}^{\overline{Q},+}(\delta)\subset A_{\mini}^+$, it follows that for $\epsilon$ sufficiently small we have $a^{-1}h(u_Q)u_Q^{-1}a\in \exp(B(0,e^{-\delta''\sigma(a)}))$. Moreover, we have
$$\displaystyle \gamma_2=(a^{-1}u_Qa)^{-1}(a^{-1}h(u_Q)u_Q^{-1}a)(a^{-1}u_Qa)$$
where $a^{-1}u_Qa\in U_Q[<\epsilon \sigma(a)]$ and there exists $\alpha>0$ such that
$$\displaystyle \lvert \Ad(g^{-1})X\rvert_{\mathfrak{g}}\leqslant  e^{\alpha\sigma(g)}\lvert X\rvert_{\mathfrak{g}}$$
for all $g\in G(F)$ and all $X\in \mathfrak{g}(F)$. From this it follows that for $\epsilon$ sufficiently small we have $\gamma_2\in \exp(B(0,e^{-\delta'\sigma(a)}))$. This ends the proof of \ref{eq 7.1.17} and thus of the proposition $\blacksquare$
\end{enumerate}

\subsection{The distribution $J^{\Lie}$}\label{section 7.2}

\noindent For all $f\in \mathcal{S}(\mathfrak{g}(F))$, let us define a function $\gls{KLief}$ on $H(F)\backslash G(F)$ by

$$\displaystyle K^{\Lie}(f,x)=\int_{\mathfrak{h}(F)} f(x^{-1}Xx)\xi(X)dX,\;\;\;\; x\in H(F)\backslash G(F)$$

\noindent the above integral being absolutely convergent. The theorem below, whose proof is similar to the proof of Theorem \ref{proposition 7.1.1}, shows that the integral

$$\displaystyle \gls{JLief}=\int_{H(F)\backslash G(F)} K^{\Lie}(f,x)dx$$

\noindent is absolutely convergent for all $f\in \mathcal{S}_{\scusp}(\mathfrak{g}(F))$ and define a continuous linear form

$$\mathcal{S}_{\scusp}(\mathfrak{g}(F))\to \mathbb{C}$$
$$f\mapsto J^{\Lie}(f)$$

\begin{theo}\label{proposition 7.2.1}
\begin{enumerate}[(i)]

\item There exists $c>0$ and a continuous semi-norm $\nu$ on $\mathcal{S}(\mathfrak{g}(F))$ such that

$$\lvert K^{\Lie}(f,x)\rvert\leqslant \nu(f)e^{c\sigma_{H\backslash G}(x)}$$

\noindent for all $x\in H(F)\backslash G(F)$ and all $f\in \mathcal{S}(\mathfrak{g}(F))$.

\item For all $c>0$, there exists a continuous semi-norm $\nu_c$ on $\mathcal{S}(\mathfrak{g}(F))$ such that

$$\lvert K^{\Lie}(f,x)\rvert\leqslant \nu_c(f)e^{-c\sigma_{H\backslash G}(x)}$$

\noindent for all $x\in H(F)\backslash G(F)$ and all $f\in\mathcal{S}_{\scusp}(\mathfrak{g}(F))$.
\end{enumerate}
\end{theo}

\section{Spectral expansion}\label{section 9}

The goal of this chapter is to give a spectral expansion for the distribution $J$ introduced in the previous chapter. The result is stated in Section \ref{section 9.1} and the proof goes through Sections \ref{section 9.2} and \ref{section 9.3}.

\subsection{The theorem}\label{section 9.1}

\noindent Let us set

$$\displaystyle \gls{Jspecf}=\int_{\mathcal{X}(G)} D(\pi)\widehat{\theta}_f(\pi)m(\overline{\pi})d\pi$$

\noindent for all $f\in \mathcal{C}_{\scusp}(G(F))$. Note that by Lemma \ref{lemma 5.4.2} and \ref{eq 2.7.2}, this integral is absolutely convergent. The purpose of this chapter is to show the following

\vspace{2mm}

\begin{theo}\label{theorem 9.1.1}
For all $f\in\mathcal{C}_{\scusp}(G(F))$, we have

$$\displaystyle J(f)=J_{\spec}(f)$$
\end{theo}

\vspace{2mm}

\noindent By Lemma \ref{lemma 5.4.2} and Theorem \ref{proposition 7.1.1}, both sides of the equality of the theorem are continuous in $f\in \mathcal{C}_{\scusp}(G(F))$. Hence, by Lemma \ref{lemma 5.3.1}(ii) it is sufficient to establish the equality for functions $f\in \mathcal{C}_{\scusp}(G(F))$ which have a compactly supported Fourier transform. We fix until the end of Section \ref{section 9.3} a function $f\in \mathcal{C}_{\scusp}(G(F))$ having a compactly supported Fourier transform.

\subsection{Study of an auxiliary distribution}\label{section 9.2}

\noindent Let us introduce, for all $f'\in \mathcal{C}(G(F))$, the following integrals

$$\displaystyle K^A_{f,f'}(g_1,g_2)=\int_{G(F)} f(g_1^{-1}gg_2) f'(g)dg,\;\;\; g_1,g_2\in G(F)$$

$$\displaystyle K^1_{f,f'}(g,x)=\int_{H(F)} K^A_{f,f'}(g,hx) \xi(h) dh,\;\;\; g,x\in G(F)$$

$$\displaystyle K^2_{f,f'}(x,y)=\int_{H(F)} K^1_{f,f'}(h^{-1}x,y) \xi(h) dh,\;\;\; x,y\in G(F)$$

$$\displaystyle \gls{Jauxff'}=\int_{H(F)\backslash G(F)} K^2_{f,f'}(x,x) dx$$

\begin{prop}\label{proposition 9.2.1}
\begin{enumerate}[(i)]

\item The integral defining $K^A_{f,f'}(g_1,g_2)$ is absolutely convergent. For all $g_1\in G(F)$ the map

$$g_2\in G(F)\mapsto K^A_{f,f'}(g_1,g_2)$$

\noindent belongs to $\mathcal{C}(G(F))$. Moreover, for all $d>0$ there exists $d'>0$ such that for every continuous semi-norm $\nu$ on $\mathcal{C}_{d'}^w(G(F))$, there exists a continuous semi-norm $\mu$ on $\mathcal{C}(G(F))$ satisfying

$$\displaystyle \nu\left(K^A_{f,f'}(g,.)\right)\leqslant \mu(f') \Xi^G(g)\sigma(g)^{-d}$$

\noindent for all $f'\in\mathcal{C}(G(F))$ and all $g\in G(F)$.

\item The integral defining $K^1_{f,f'}(g,x)$ is absolutely convergent. Moreover, for all $d>0$, there exist $d'>0$ and a continuous semi-norm $\nu_{d,d'}$ on $\mathcal{C}(G(F))$ such that

$$\left\lvert K_{f,f'}^1(g,x)\right\rvert \leqslant \nu_{d,d'}(f')\Xi^G(g)\sigma(g)^{-d}\Xi^{H\backslash G}(x)\sigma_{H\backslash G}(x)^{d'}$$

\noindent for all $f'\in \mathcal{C}(G(F))$ and all $g,x\in G(F)$.

\item The integral defining $K^2_{f,f'}(x,y)$ is absolutely convergent. Moreover, we have

$$\displaystyle K^2_{f,f'}(x,y)=\int_{\mathcal{X}_{\tempe}(G)} \mathcal{L}_\pi(\pi(x)\pi(f)\pi(y^{-1}))\overline{\mathcal{L}_{\pi}(\pi(\overline{f'}))} \mu(\pi) d\pi$$

\noindent for all $f'\in \mathcal{C}(G(F))$ and all $x,y\in G(F)$, the integral above being absolutely convergent.

\item The integral defining $J_{\auxi}(f,f')$ is absolutely convergent. More precisely, for every $d>0$ there exists a continuous semi-norm $\nu_d$ on $\mathcal{C}(G(F))$ such that

$$\left\lvert K^2_{f.f'}(x,x)\right\rvert\leqslant \nu_d(f')\Xi^{H\backslash G}(x)^2\sigma_{H\backslash G}(x)^{-d}$$\

\noindent for all $f'\in \mathcal{C}(G(F))$ and all $x\in H(F)\backslash G(F)$. In particular, the linear form

$$f'\in \mathcal{C}(G(F))\mapsto J_{\auxi}(f,f')$$

\noindent is continuous.
\end{enumerate}
\end{prop}

\vspace{4mm}

\noindent\ul{Proof}: The point (i) follows from Theorem \ref{theorem 5.5.1}(i). The point (ii) follows from (i), Lemma \ref{lemma 6.5.1}(ii) and Lemma \ref{lemma 8.3.1}(ii). The absolute convergence of the integral defining $K^2_{f,f'}(x,y)$ follows from (ii) and Lemma \ref{lemma 6.5.1}(ii). The spectral formula for $K^2_{f,f'}(x,y)$ is a direct application of Lemma \ref{lemma 8.2.1}(v).  We are thus only left with proving the estimate (iv). For $f'\in \mathcal{C}(G(F))$ the section

$$T(f'): \pi\in \mathcal{X}_{\tempe}(G)\mapsto \overline{\mathcal{L}_\pi(\pi(\overline{f'}))}\pi(f)\in \End(\pi)^\infty$$

\noindent is smooth by Lemma \ref{lemma 8.2.1}(i) and is compactly supported by the hypothesis on $f$. Hence it belongs to $\mathcal{C}(\mathcal{X}_{\tempe}(G),\mathcal{E}(G))$ and by the matricial Paley-Wiener theorem (Theorem \ref{theorem 2.6.1}) there exists a unique function $\varphi_{f'}\in \mathcal{C}(G(F))$ such that $\pi(\varphi_{f'})=\mathcal{L}_\pi(\pi(f'))\pi(f)$ for all $\pi\in \mathcal{X}_{\tempe}(G)$. By Lemma \ref{lemma 5.3.1}(i), the function $\varphi_{f'}$ is strongly cuspidal for all $f'\in \mathcal{C}(G(F))$. By the formula (iii) for $K^2_{f,f'}$ and Lemma \ref{lemma 8.2.1}(iv), we have

$$K^2_{f,f'}(x,x)=K(\varphi_{f'},x)$$

\noindent for all $f'\in \mathcal{C}(G(F))$ and for all $x\in H(F)\backslash G(F)$. Hence by Theorem \ref{proposition 7.1.1}, for all $d>0$, there exists a continuous semi-norm $\mu_d$ on $\mathcal{C}(G(F))$ such that

$$\lvert K^2_{f,f'}(x,x)\rvert \leqslant \mu_d(\varphi_{f'})\Xi^{H\backslash G}(x)^2\sigma_{H\backslash G}(x)^{-d}$$

\noindent for all $f'\in \mathcal{C}(G(F))$ and for all $x\in H(F)\backslash G(F)$. To conclude, it is thus sufficient to show that the linear map

$$\mathcal{C}(G(F))\to \mathcal{C}(G(F))$$
$$f'\mapsto \varphi_{f'}$$

\noindent is continuous. By Theorem \ref{theorem 2.6.1}, it suffices to show that the linear map

$$\displaystyle f'\in \mathcal{C}(G(F))\mapsto \left(\pi\in \mathcal{X}_{\tempe}(G)\mapsto \overline{\mathcal{L}_\pi(\pi(\overline{f'}))}\right)\in C^\infty(\mathcal{X}_{\tempe}(G))$$

\noindent is continuous, where the topology on the target space is the obvious one. This follows easily from Lemma \ref{lemma 8.2.1}(i) and Theorem \ref{theorem 2.6.1}. $\blacksquare$

\vspace{5mm}

\begin{prop}\label{proposition 9.2.2}
We have the equality

$$\displaystyle J_{\auxi}(f,f')=\int_{\mathcal{X}(G)} D(\pi) \widehat{\theta}_f(\pi) \overline{\mathcal{L}_{\pi}(\pi(\overline{f'}))}d\pi$$

\noindent for all $f'\in \mathcal{C}(G(F))$.
\end{prop}

\vspace{5mm}

\noindent\ul{Proof}: Let $a:\mathbb{G}_m\to A$ be a one parameter subgroup such that $\lambda(a(t)ha(t)^{-1})=t\lambda(h)$ for all $t\in \mathbb{G}_m$ and all $h\in H$ (recall that we denote by $\lambda:H\to \mathbb{G}_a$ the additive character such that $\xi=\psi\circ\lambda_F$). We denote as usual by $R$, $L$ and $\Ad$ the action by right translation, left translation and conjugation of $G(F)$ on functions on $G(F)$. We will set $R_a=R\circ a_F$, $L_a=L\circ a_F$ and $\Ad_a=\Ad\circ a_F$. These provide smooth representations of $F^\times$ on $\mathcal{C}(G(F))$. Let $f'\in \mathcal{C}(G(F))$. We want to prove the formula of the proposition for this function $f'$. By Dixmier-Malliavin in the real case, we may write $f'$ as a finite sum $f'=\sum_{i=1}^k \Ad_a(\varphi_i)(f_i'')$ where $\varphi_i\in C_c^\infty(F^\times)$ and $f''_i\in \mathcal{C}(G(F))$ for $1\leqslant i\leqslant k$. By linearity, we may assume that this sum has only one element, that is $f'=\Ad_a(\varphi)(f'')$ with $\varphi\in C_c^\infty(G(F))$ and $f''\in \mathcal{C}(G(F))$. By continuity of the linear form $J_{\auxi}$, we have

$$\displaystyle J_{\auxi}(f,f')=\int_{F^\times} \varphi(t)J_{\auxi}(f,\Ad_a(t)f'') d^\times t$$

\noindent Returning to the definition of $J_{\auxi}$, we get

$$\displaystyle J_{\auxi}(f,f')=\int_{F^\times} \int_{H(F)\backslash G(F)} \varphi(t) K^2_{f,\Ad_a(t)f''}(x,x)dxd^\times t$$

\noindent By Proposition \ref{proposition 9.2.1}(iv), this double integral is absolutely convergent. Doing the variable change $x\mapsto a(t)x$ and switching the two integrals, we obtain

\begin{align}\label{eq 9.2.1}
\displaystyle J_{\auxi}(f,f')=\int_{H(F)\backslash G(F)}\int_{F^\times} \varphi(t)\delta_H(a(t))^{-1} K^2_{f,\Ad_a(t)f''}(a(t)x,a(t)x) d^\times t dx
\end{align}

\noindent By definition, the inner integral above is equal to

$$\displaystyle \int_{F^\times} \varphi(t)\delta_H(a(t))^{-1} \int_{H(F)} K^1_{f,\Ad_a(t)f''}(ha(t)x,a(t)x)\xi(h)^{-1} dh d^\times t$$

\noindent By Proposition \ref{proposition 9.2.1}(ii), this double integral is also absolutely convergent. Doing the variable change $h\mapsto a(t)ha(t)^{-1}$, switching the two integrals and noticing that $K^1_{f,\Ad_a(t)f''}(a(t)hx,a(t)x)=K^1_{f,R_a(t)f''}(hx,a(t)x)$, we obtain the equality

\begin{align}\label{eq 9.2.2}
\displaystyle \int_{F^\times} \varphi(t)\delta_H(a(t))^{-1} & K^2_{f,\Ad_a(t)f''}(a(t)x,a(t)x) d^\times t \\
\nonumber & =\int_{H(F)}\int_{F^\times} \varphi(t)K^1_{f,R_a(t)f''}(hx,a(t)x) \psi(-t\lambda(h))d^\times t dh
\end{align}

\noindent By definition, the inner integral of the last expression above is equal to

$$\displaystyle \int_{F^\times} \varphi(t)\int_{H(F)} K^A_{f,R_a(t)f''}(hx,h'a(t)x)\xi(h')dh'\psi(-t\lambda(h))d^\times t$$

\noindent After the variable change $h'\mapsto a(t)h'ha(t)^{-1}$, this becomes

$$\displaystyle \int_{F^\times}\varphi(t)\int_{H(F)} K^A_{f,f''}(hx,h'hx)\delta_H(a(t)) \psi\left(t\lambda(h')\right) d^\times t dh'$$

\noindent By Proposition \ref{proposition 9.2.1}(i), this double integral is absolutely convergent. Switching the two integrals, we obtain

\begin{align}\label{eq 9.2.3}
\displaystyle \int_{F^\times} \varphi(t) & K^1_{f,R_a(t)f''}(hx,a(t)x) \psi(-t\lambda(h))d^\times t \\
\nonumber & =\int_{H(F)}\int_{F^\times} K^A_{f,f''}(hx,h'hx)\varphi(t)\delta_H(a(t)) \psi\left(t\lambda(h')\right) d^\times t dh'
\end{align}

\noindent We have $d^\times t=\lvert t\rvert^{-1} dt$ where $dt$ is an additive Haar measure on $F$. Let us set $\varphi'(t)=\varphi(t)\delta_H(a(t))\lvert t\rvert^{-1}$ and

$$\displaystyle \widehat{\varphi'}(x)=\int_F \varphi'(t) \psi(tx) dt,\;\;\;x\in F$$

\noindent for its Fourier transform. By \ref{eq 9.2.1}, \ref{eq 9.2.2} and \ref{eq 9.2.3}, we have

\begin{align}\label{eq 9.2.4}
\displaystyle J_{\auxi}(f,f')=\int_{H(F)\backslash G(F)} \int_{H(F)}\int_{H(F)} K^A_{f,f''}(hx,h'hx) \widehat{\varphi'}(\lambda(h')) dh'dhdx
\end{align}

\noindent For $N>0$ and $M>0$, let us denote by $\alpha_N:H(F)\backslash G(F)\to \{0, 1\}$ and $\beta_M:G(F)\to \{0, 1\}$ the characteristic functions of the sets $\{x\in H(F)\backslash G(F);\sigma_{H\backslash G}(x)\leqslant N\}$ and $\{g\in G(F); \sigma(g)\leqslant M\}$ respectively. For all $N\geqslant 1$ and $C>0$, we set

$$\displaystyle J_{\auxi,N}(f,f'):=\int_{H(F)\backslash G(F)} \alpha_N(x)\int_{H(F)}\int_{H(F)} K^A_{f,f''}(hx,h'hx) \widehat{\varphi'}(\lambda(h')) dh'dhdx$$

\[\begin{aligned}
\displaystyle J_{\auxi,N,C}(f,f'):=\int_{H(F)\backslash G(F)} \alpha_N(x) \int_{H(F)}\int_{H(F)} \beta_{C\log(N)}(h') K^A_{f,f''}(hx,h'hx) \widehat{\varphi'}(\lambda(h')) dh'dhdx
\end{aligned}\]

\noindent By \ref{eq 9.2.4}, we have

\begin{align}\label{eq 9.2.5}
J_{\auxi}(f,f')=\lim\limits_{N\to\infty} J_{\auxi,N}(f,f')
\end{align}

\noindent Moreover, we have the following

\vspace{3mm}

\begin{num}
\item\label{eq 9.2.6} The triple integrals defining $J_{\auxi,N}(f,f')$ and $J_{\auxi,N,C}(f,f')$ are absolutely convergent and there exists $C>0$ such that
$$\left\lvert J_{\auxi,N}(f,f')-J_{\auxi,N,C}(f,f')\right\rvert \ll N^{-1}$$
for all $N\geqslant 1$.
\end{num}
 
\vspace{3mm}

\noindent Indeed, since $\widehat{\varphi'}\in \mathcal{S}(F)$, we have $\lvert \widehat{\varphi'}(\lambda)\rvert\ll (1+\lvert \lambda\rvert)^{-1}$ for all $\lambda\in F$. Hence, by Theorem \ref{theorem 5.5.1}(i), there exists $d>0$ such that

\[\begin{aligned}
\displaystyle \left\lvert J_{\auxi,N}(f,f')\right\rvert \ll \int_{H(F)\backslash G(F)} \alpha_N(x)\int_{H(F)}\int_{H(F)} & \Xi^G(hx)\Xi^G(h'hx)\sigma(hx)^d \\
 & \sigma(h'hx)^d \left(1+\lvert \lambda(h')\rvert\right)^{-1}dh'dhdx
\end{aligned}\]

\[\begin{aligned}
\displaystyle \left\lvert J_{\auxi,N,C}(f,f')\right\rvert\ll \int_{H(F)\backslash G(F)} \alpha_N(x)\int_{H(F)}\int_{H(F)} & \beta_{C \log(N)}(h') \Xi^G(hx)\Xi^G(h'hx) \\
 & \sigma(hx)^d\sigma(h'hx)^d \left(1+\lvert \lambda(h')\rvert\right)^{-1}dh'dhdx
\end{aligned}\]

\noindent and

\[\begin{aligned}
\displaystyle \left\lvert J_{\auxi,N}(f,f')-J_{\auxi,N,C}(f,f')\right\rvert \ll \int_{H(F)\backslash G(F)}\alpha_N(x) & \int_{H(F)}\int_{H(F)}\mathbf{1}_{\sigma\geqslant C\log(N)}(h')\Xi^G(hx)\Xi^G(h'hx) \\
 & \sigma(hx)^d\sigma(h'hx)^d \left(1+\lvert \lambda(h')\rvert\right)^{-1}dh'dhdx
\end{aligned}\]

\noindent for all $N\geqslant 1$ and all $C>0$. By Proposition \ref{proposition 6.7.1}(vii) (applied to $c=1$), there exists $d'>0$ such that the first two integrals above are essentially bounded by

$$\displaystyle \int_{H(F)\backslash G(F)} \alpha_N(x)\Xi^{H\backslash G}(x)^2\sigma_{H\backslash G}(x)^{d'}dx$$

\noindent which is of course an absolutely convergent integral (the integrand is bounded and compactly supported). On the other hand, by Proposition \ref{proposition 6.7.1}(vii), there exist $\epsilon>0$ and $d'>0$ such that the third integral above is essentially bounded by

$$\displaystyle e^{-\epsilon C\log(N)}\int_{H(F)\backslash G(F)} \alpha_N(x)\Xi^{H\backslash G}(x)^2\sigma_{H\backslash G}(x)^{d'}dx$$

\noindent for all $N\geqslant 1$ and all $C>0$. By Proposition \ref{proposition 6.7.1}(iv), there exists $d''>0$ such that this last term is itself essentially bounded by $e^{-\epsilon C\log(N)}N^{d''}$, for all $N\geqslant 1$ and all $C>0$. Choosing $C$ to be bigger than $(d''+1)/\epsilon$ gives the estimate of \ref{eq 9.2.6}.

\vspace{2mm}

\noindent Let us fix $C>0$ which satisfies \ref{eq 9.2.6}. By \ref{eq 9.2.5}, it follows that

\begin{align}\label{eq 9.2.7}
\displaystyle J_{\auxi}(f,f')=\lim\limits_{N\to \infty} J_{\auxi,N,C}(f,f')
\end{align}

\noindent Since the triple integral defining $J_{\auxi,N,C}(f,f')$ is absolutely convergent, we may write

$$J_{\auxi,N,C}(f,f')=\int_{H(F)} \beta_{C\log(N)}(h)\widehat{\varphi'}(\lambda(h)) \int_{G(F)} \alpha_N(g) K^A_{f,f''}(g,hg) dg dh$$

\noindent We now prove the following estimate

\begin{align}\label{eq 9.2.8}
\displaystyle \left\lvert J_{\auxi,N,C}(f,f')-\int_{H(F)} \beta_{C\log(N)}(h)\widehat{\varphi'}(\lambda(h)) \int_{G(F)} K^A_{f,f''}(g,hg) dg dh \right\rvert \ll N^{-1}
\end{align}

\noindent for all $N\geqslant 1$.

\vspace{2mm}

\noindent Indeed, since $f$ is strongly cuspidal, by Theorem \ref{theorem 5.5.1}(iii), there exists $c_1>0$ such that for all $d>0$, there exists $d'>0$ such that

$$\left\lvert K^A_{f,f''}(g,hg)\right\rvert\ll \Xi^G(g)^2 \sigma(g)^{-d} e^{c_1 \sigma(h)}\sigma(h)^{d'}$$

\noindent for all $g\in G(F)$ and all $h\in H(F)$. Fix such a $c_1>0$. Also, choose $d_0>0$ such that the function $g\mapsto \Xi^G(g)^2\sigma(g)^{-d_0}$ is integrable over $G(F)$ (Proposition \ref{proposition 1.5.1}(v)). Then, by the above inequality, for all $d>d_0$ there exists $d'>0$ such that the left hand side of \ref{eq 9.2.8} is essentially bounded by

$$\displaystyle N^{c_1C-d+d_0}\log(N)^{d'}\int_{H(F)} \beta_{C\log(N)}(h) dh$$

\noindent for all $N\geqslant 1$ (here the factor $N^{d_0-d}$ comes from the fact that $\sigma(g)^{-1}\ll N^{-1}$ for all $g\in G(F)$ such that $\alpha_N(g)=0$). Now, by Lemma \ref{lemma B.1.3}, there exists $c_2>0$ such that

$$\displaystyle \int_{H(F)} \beta_{C\log(N)}(h) dh\ll N^{c_2}$$

\noindent for all $N\geqslant 1$. Hence, it suffices to choose $d>c_1C+d_0+c_2+1$ to get the estimate \ref{eq 9.2.8}.

\vspace{2mm}

\noindent From \ref{eq 9.2.7} and \ref{eq 9.2.8}, we deduce that

\begin{align}\label{eq 9.2.9}
\displaystyle J_{\auxi}(f,f')=\lim\limits_{N\to \infty}\int_{H(F)} \beta_{C\log(N)}(h)\widehat{\varphi'}(\lambda(h)) \int_{G(F)} K^A_{f,f''}(g,hg) dg dh
\end{align}

\noindent Arthur's local trace formula allows us to express the inner integral above in spectral terms. Indeed, since $f$ is strongly cuspidal, by Theorem \ref{theorem 5.5.1}(v), we have the equality

\begin{align}\label{eq 9.2.10}
\displaystyle \int_{G(F)} K^A_{f,f''}(g,hg) dg=\int_{\mathcal{X}(G)}D(\pi)\widehat{\theta}_f(\pi) \theta_{\overline{\pi}}(R(h^{-1})f'') d\pi
\end{align}

\noindent for all $h\in H(F)$. By \ref{eq 2.2.5} and \ref{eq 2.2.7}, it is easy to see that

$$\displaystyle \lvert \theta_{\overline{\pi}}(R(h^{-1})f'')\rvert\ll \Xi^G(h)$$

\noindent for all $\pi\in \mathcal{X}(G)$ and all $h\in H(F)$. Hence, by \ref{eq 2.7.2} and Lemma \ref{lemma 5.4.2}, we have an inequality

$$\displaystyle \int_{\mathcal{X}(G)} D(\pi) \left\lvert \widehat{\theta}_f(\pi) \theta_{\overline{\pi}}(R(h^{-1})f'')\right\rvert d\pi\ll \Xi^G(h)$$

\noindent for all $h\in H(F)$. Hence, by Lemma \ref{lemma 6.5.1}(iii), the double integral

$$\displaystyle \int_{H(F)}\widehat{\varphi'}(\lambda(h)) \int_{\mathcal{X}(G)} D(\pi)\widehat{\theta}_f(\pi) \theta_{\overline{\pi}}(R(h^{-1})f'')d\pi dh$$

\noindent is absolutely convergent and by \ref{eq 9.2.9} and \ref{eq 9.2.10} this double integral is equal to $J_{\auxi}(f,f')$. Switching the two integrals and applying Lemma \ref{lemma 8.1.1}(ii), we get

\[\begin{aligned}
\displaystyle J_{\auxi}(f,f') & = \int_{\mathcal{X}(G)} D(\pi) \widehat{\theta}_f(\pi)\overline{\mathcal{L}_{\pi}\left(\pi(\overline{\Ad_a(\varphi)f''})\right)}d\pi \\
 & = \int_{\mathcal{X}(G)} D(\pi) \widehat{\theta}_f(\pi)\overline{\mathcal{L}_{\pi}(\pi(\overline{f'}))}d\pi
\end{aligned}\]

\noindent which is the equality we were looking for and this ends the proof of the proposition. $\blacksquare$

\subsection{End of the proof of Theorem \ref{theorem 9.1.1}}\label{section 9.3}

\noindent Recall that we have fixed a function $f\in \mathcal{C}_{\scusp}(G(F))$ having a compactly supported Fourier transform. By Lemma \ref{lemma 8.2.1}(iv), we have

\begin{align}\label{eq 9.3.1}
\displaystyle K(f,x)=\int_{\mathcal{X}_{\tempe}(G)} \mathcal{L}_{\pi}(\pi(x)\pi(f)\pi(x^{-1})) \mu(\pi) d\pi
\end{align}

\noindent for all $x\in H(F)\backslash G(F)$. By Corollary \ref{corollary 8.6.1}(ii), there exists a function $f'\in \mathcal{C}(G(F))$ such that

\begin{align}\label{eq 9.3.2}
\mathcal{L}_{\pi}(\pi(\overline{f'}))=m(\pi)
\end{align}

\noindent for all $\pi\in \mathcal{X}_{\tempe}(G)$ such that $\pi(f)\neq 0$. Also, by Theorem \ref{theorem 8.2.1} and Corollary \ref{corollary 8.6.1}(i), for all $\pi\in \mathcal{X}_{\tempe}(G)$, we have

$$\mathcal{L}_\pi\neq 0\Leftrightarrow m(\pi)=1$$

\noindent Hence, by \ref{eq 9.3.1}, we have the equality

$$\displaystyle K(f,x)=\int_{\mathcal{X}_{\tempe}(G)} \mathcal{L}_{\pi}(\pi(x)\pi(f)\pi(x^{-1})) \overline{\mathcal{L}_{\pi}(\pi(\overline{f'}))}\mu(\pi) d\pi$$

\noindent and by Proposition \ref{proposition 9.2.1}(iii), it follows that

$$K(f,x)=K^2_{f,f'}(x,x)$$

\noindent for all $x\in H(F)\backslash G(F)$. Consequently, we have the equality

$$J(f)=J_{\auxi}(f,f')$$

\noindent Applying Proposition \ref{proposition 9.2.2}, we deduce that

$$\displaystyle J(f)=\int_{\mathcal{X}(G)} D(\pi)\widehat{\theta}_f(\pi)\overline{\mathcal{L}_{\pi}(\pi(\overline{f'}))}d\pi$$

\noindent Let $\pi\in\mathcal{X}(G)$ be such that $\widehat{\theta}_f(\pi)\neq 0$ and let $\pi'$ be the unique representation in $\mathcal{X}_{\tempe}(G)$ such that $\pi$ is a linear combination of subrepresentations of $\pi'$. Then, we have $\pi'(f)\neq 0$. Hence, by \ref{eq 9.3.2} and Corollary \ref{corollary 8.6.1}(ii), we have $\overline{\mathcal{L}_{\pi}(\pi(\overline{f'}))}=\overline{m(\pi)}=m(\overline{\pi})$. It follows that

$$\displaystyle J(f)=\int_{\mathcal{X}(G)} D(\pi)\widehat{\theta}_f(\pi)m(\overline{\pi})d\pi$$

\noindent and this ends the proof of Theorem \ref{theorem 9.1.1}. $\blacksquare$

\section{The spectral expansion of $J^{\Lie}$}\label{section 10}

In this chapter, we establish a `spectral' expansion for the distribution $J^{\Lie}$ that was introduced in Chapter \ref{section 7}. More precisely, we will express $J^{\Lie}(f)$ in terms of (weighted) orbital integrals of the Fourier transform of the test function $f$. Even the statement of this expansion (Theorem \ref{theorem 10.8.1}) needs some preparation which is the content of Sections \ref{section 10.1} to \ref{section 10.7}. To be a little bit more specific, by Fourier inversion we can rewrite the kernel function involved into the definition of $J^{\Lie}(f)$ as an integral over a certain affine subspace $\Sigma(F)$ of $\mathfrak{g}(F)$. Then, Sections \ref{section 10.1} to \ref{section 10.7} are devoted to a thorough study of this affine space $\Sigma$ and in particular of the adjoint action of $H$ on it. Doing so, we isolate a certain $H$-invariant Zariski open subset $\Sigma'$ of $\Sigma$ with particularly nice properties including freeness of the adjoint action of $H$ (Proposition \ref{proposition 10.5.1}) and a certain genericity property for the Borel subalgebras intersecting $\Sigma'$ (Proposition \ref{proposition 10.6.1}). Once these preparations are in place, we can state in Section \ref{section 10.8} the main result of this chapter (Theorem \ref{theorem 10.8.1}) the proof of which goes though Sections \ref{section 10.9} to \ref{section 10.11}. The basic scheme of the proof is inspired from \cite{Wa1} Section 9. We first introduce some truncation in the original expression defining $J^{\Lie}(f)$ (there are lot of freedom in the choice of this truncation, see Section \ref{section 10.9}). After this, we end up with an expression (depending on some integer $N$) $J_N^{\Lie}(f)$ which naturally decomposes as a sum (or rather an integral) of orbital integrals weighted by certain weights which depend on the truncation and in Section \ref{section 10.10} we show that we can replace these weights by other more rigid ones which have been studied by Arthur. Then using some computations made by Arthur in the course of establishing his local trace formula \cite{A1}, we readily finish the proof in Section \ref{section 10.11}. 

\subsection{The affine subspace $\Sigma$}\label{section 10.1}

\noindent Recall that in Section \ref{section 6.2} we have defined a parabolic subgroup $P=MN$ of $G=U(W)\times U(V)$ as a product $U(W)\times P_V$ where $P_V$ is a certain parabolic subgroup of $U(V)$. Let $\gls{Pbar}=M\overline{N}$ be the parabolic subgroup opposite to $P$ with respect to $M$. Then the unipotent radicals $N$ and $\gls{Nbar}$ can also be seen as subgroups of $U(V)$ and we will identify them as such in what follows.

\noindent Recall also that in Section \ref{section 6.2} we have defined a character $\xi$ on $N(F)$ which extends to a character of $H(F)=U(W)(F)\ltimes N(F)$ trivial on $U(W)(F)$. Using the $G$-invariant bilinear pairing $B$ on $\mathfrak{g}$ defined in the same section, there exist a unique element $\gls{Xi}\in \overline{\mathfrak{n}}(F)$ such that

$$\xi(X)=\psi(B(\Xi,X))$$

\noindent for all $X\in \mathfrak{n}(F)$.

\noindent We have the following explicit description of $\Xi$ (seen as an element of $\mathfrak{u}(V)$):

\vspace{2mm}

\begin{num}
\item\label{eq 10.1.1} $\Xi z_i=z_{i-1}$, for $1\leqslant i\leqslant r$, $\Xi z_{-i}=-z_{-i-1}$, for $0\leqslant i\leqslant r-1$, $\Xi z_{-r}=0$ and $\Xi(W)=0$.
\end{num}

\vspace{2mm}

\noindent Set $\gls{Sigma}=\Xi+ \mathfrak{h}^\perp$ where $\gls{hperp}$ is the orthogonal of $\mathfrak{h}$ in $\mathfrak{g}$ for $B(.,.)$. 

\vspace{2mm}

\noindent Recall that we have fixed a Haar measure $dX$ on $\mathfrak{h}(F)$. In this whole chapter we will denote this Haar measure by $d\mu_\mathfrak{h}(X)$. In Section \ref{section 1.6}, we explained how to associate to $d\mu_\mathfrak{h}(X)$, using $B(.,.)$, a Haar measure $d\mu_\mathfrak{h}^\perp$ on $\mathfrak{h}^\perp(F)$. Let us denote by $d\mu_\Sigma$ the translate of this measure to $\Sigma(F)$. Then, by \ref{eq 1.6.1}, we have the following equality

\begin{align}\label{eq 10.1.2}
\displaystyle \int_{\mathfrak{h}(F)} f(X)\xi(X)d\mu_{\mathfrak{h}}(X)=\int_{\Sigma(F)} \widehat{f}(Y) d\mu_{\Sigma}(Y)
\end{align}

\noindent for all $f\in \mathcal{S}(\mathfrak{g}(F))$.

\subsection{Conjugation by $N$}\label{section 10.2}

\noindent We have the following explicit description of $\mathfrak{h}^\perp$: an element $X=(X_W,X_V)\in \mathfrak{g}=\mathfrak{u}(W)\oplus \mathfrak{u}(V)$ is in $\mathfrak{h}^\perp$ if and only if we have a decomposition

$$X_V=-X_W+c(z_0,w)+\lambda c(z_0,\eta z_0)+A+N$$

\noindent for some $w\in W_{\overline{F}}$, $\lambda\in \overline{F}$, $A\in \mathfrak{a}$ and $N\in \mathfrak{n}$ (recall that $\eta\in E$ is a nonzero element with trace zero and cf.\ Section \ref{section 6.1} for the notation $c(v,v')$). Thus for every element $X=(X_V,X_W)$ of $\Sigma$ we have a decomposition

\begin{align}\label{eq 10.2.1}
X_V=\Xi-X_W+c(z_0,w)+\lambda c(z_0,\eta z_0)+A+N
\end{align}

\noindent where $w$, $\lambda$, $A$ and $N$ are as above. Let us define the following affine subspaces of $\mathfrak{g}$:

\vspace{2mm}

\begin{itemize}
\renewcommand{\labelitemi}{$\bullet$}
\item $\mathfrak{u}(W)^-=\{(X_W,-X_W); \;\; X_W\in \mathfrak{u}(W)\}$;

\item $\Lambda_0$ is the subspace of $\mathfrak{u}(V)\subset \mathfrak{g}$ generated by the $c(z_i,\eta z_i)$ for $i=0,\ldots,r$, the $c(z_i,z_{i+1})$ for $i=0,\ldots,r-1$, and the $c(z_r,w)$ for $w\in W$;

\item $\gls{Lambda}=\Xi+\big( \mathfrak{u}(W)^-\oplus \Lambda_0\big)$.
\end{itemize}

\vspace{3mm}

\begin{prop}\label{proposition 10.2.1}
Conjugation by $N$ preserves $\Sigma$ and induces an isomorphism of algebraic varieties:

$$N\times \Lambda\to \Sigma$$
$$(n,X)\mapsto nXn^{-1}$$
\end{prop}

\vspace{3mm}

\noindent\ul{Proof}: First we show that the map

\begin{align}\label{eq 10.2.2}
& N\times \Lambda\to \Sigma \\
\nonumber & (n,X)\mapsto nXn^{-1}
\end{align}

\noindent is injective. This amounts to proving that for all $n\in N$ and all $X\in \Lambda$ if $nXn^{-1}\in \Lambda$ then $n=1$. So let $n\in N$ and $X=(X_W,X_V)\in \Lambda$ be such that $nXn^{-1}\in \Lambda$. By definition of $\Lambda$, we may write $X_V$ and $nX_Vn^{-1}$ as

\begin{align}\label{eq 10.2.3}
\displaystyle X_V=\Xi-X_W+c(z_r,w)+\sum_{i=0}^r \lambda_i c(z_i,\eta z_i)+ \sum_{i=0}^{r-1} \mu_i c(z_i,z_{i+1})
\end{align}

\begin{align}\label{eq 10.2.4}
\displaystyle nX_Vn^{-1}=\Xi-X_W+c(z_r,w')+\sum_{i=0}^r \lambda'_i c(z_i,\eta z_i)+ \sum_{i=0}^{r-1} \mu'_i c(z_i,z_{i+1})
\end{align}

\noindent where $w,w'\in W_{\overline{F}}$, $\lambda_i,\lambda_i'\in \overline{F}$, $0\leqslant i\leqslant r$, and $\mu_i,\mu'_i\in \overline{F}$, $0\leqslant i\leqslant r-1$. Let us prove first that

\begin{align}\label{eq 10.2.5}
nz_i=z_i \mbox{ for all } 0\leqslant i\leqslant r.
\end{align}

\noindent The proof is by descending induction. The result is trivial for $i=r$ by definition of $N$. Assume that the equality \ref{eq 10.2.5} is true for some $1\leqslant i\leqslant r$. Then, from \ref{eq 10.2.3} we easily deduce that
$$(nX_Vn^{-1})z_i=nX_Vz_i=n\Xi z_i=n z_{i-1}$$
and from \ref{eq 10.2.4}, it is not hard to see that
$$(nX_Vn^{-1})z_i=\Xi z_i=z_{i-1}$$
Thus, the equality \ref{eq 10.2.5} is satisfied with $i-1$ instead of $i$ and this ends the proof of \ref{eq 10.2.5}.

\vspace{2mm}

\noindent We now prove the following

\begin{align}\label{eq 10.2.6}
nz_{-i}=z_{-i} \mbox{ for all } 1\leqslant i\leqslant r.
\end{align}

\noindent We prove this by strong induction on $i$. First, we treat the case $i=1$. By \ref{eq 10.2.3} and \ref{eq 10.2.5}, we have
$$(nX_Vn^{-1})z_0=nX_Vz_0=n(-z_{-1}+2\lambda_0\eta\nu z_0+\mu_0\nu z_1)=-nz_{-1}+2\lambda_0\eta\nu z_0+\mu_0\nu z_1$$
On the other hand, by \ref{eq 10.2.4}, we have
$$(nX_Vn^{-1})z_0=-z_{-1}+2\lambda_0'\eta\nu z_0+\mu_0'\nu z_1$$
It follows that
$$nz_{-1}-z_{-1}=2(\lambda_0-\lambda_0')\eta\nu z_0+(\mu_0-\mu_0')\nu z_1$$
But, since $nz_0=z_0$ and $n\in U(V)$ we have $h(nz_{-1},z_0)=h(z_{-1},z_0)=0$ and $h(nz_{-1},z_{-1})\in \overline{F}\eta$. From this we deduce that $\lambda_0=\lambda_0'$ and $\mu_0=\mu'_0$ so that indeed $nz_{-1}=z_{-1}$. Let $1\leqslant j\leqslant r-1$ and assume now that \ref{eq 10.2.6} is true for all $1\leqslant i\leqslant j$. By \ref{eq 10.2.3} and \ref{eq 10.2.5}, we have
\[\begin{aligned}
(nX_Vn^{-1})z_{-j} & =nX_Vz_{-j}=n(-z_{-j-1}+2\lambda_j\eta\nu z_j-\mu_{j-1}\nu z_{j-1}+\mu_j\nu z_{j+1}) \\
 & =-nz_{-j-1}+2\lambda_j\eta\nu z_j-\mu_{j-1}\nu z_{j-1}+\mu_j\nu z_{j+1}
\end{aligned}\]

\noindent On the other hand, by \ref{eq 10.2.4}, we have

$$(nX_Vn^{-1})z_{-j}=-z_{-j-1}+2\lambda_j'\eta\nu z_j-\mu'_{j-1}\nu z_{j-1}+\mu_j'\nu z_{j+1}$$

\noindent It follows that

$$nz_{-j-1}-z_{-j-1}=2(\lambda_j-\lambda_j')\eta\nu z_j+(\mu'_{j-1}-\mu_{j-1})\nu z_{j-1}+(\mu_j-\mu_j')\nu z_{j+1}$$

\noindent Since $nz_{-j}=z_{-j}$, $nz_{-j+1}=z_{-j+1}$ (by the induction hypothesis) and $n\in U(V)$, we have $h(nz_{-j-1},z_{-j})=h(z_{-j-1},z_{-j})=0$, $h(nz_{-j-1},z_{-j+1})=h(z_{-j-1},z_{-j+1})=0$ and \\
\noindent $h(nz_{-j-1},z_{j+1})\in \overline{F}\eta$. From this we deduce that $\lambda_j=\lambda_j'$, $\mu'_{j-1}=\mu_{j-1}$ and $\mu_j=\mu'_j$ so that indeed $nz_{-j-1}=z_{-j-1}$. This ends the proof of \ref{eq 10.2.6}.

\vspace{2mm}

\noindent From \ref{eq 10.2.5} and \ref{eq 10.2.6} and since $n\in N$, we may now deduce that $n=1$. This ends the proof that the map \ref{eq 10.2.2} is injective. We easily compute

$$\dim(N\times \Lambda)=\dim(N)+\dim(\Lambda)=(2r^2+r+2mr)+(m^2+2m+2r+1)=2r^2+3r+2mr+(m+1)^2$$

\[\begin{aligned}
\dim(\Sigma)=\dim(\mathfrak{h}^\perp)=\dim(G)-\dim(H) & =(m+2r+1)^2+m^2-m^2-\dim(N) \\
 & =2r^2+3r+2mr+(m+1)^2
\end{aligned}\]

\noindent where $m=\dim(W)$. Hence, we have $\dim(\Sigma)=\dim(N\times \Lambda)$. Since we are in characteristic zero, it follows that the regular map \ref{eq 10.2.2} induces an isomorphism between $N\times\Lambda$ and a Zariski open subset of $\Sigma$. But, obviously $N\times\Lambda$ and $\Sigma$ are both affine spaces so that the only Zariski open subset of $\Sigma$ that can be isomorphic to $N\times \Lambda$ is $\Sigma$ itself. It follows that the regular map \ref{eq 10.2.2} is an isomorphism. $\blacksquare$

\subsection{Characteristic polynomial}\label{section 10.3}

\noindent Let $X=(X_W,X_V)\in \Lambda$. By definition of $\Lambda$, we may write

\begin{align}\label{eq 10.3.1}
\displaystyle X_V=\Xi-X_W+c(z_r,w)+\sum_{i=0}^r \lambda_i c(z_i,\eta z_i)+ \sum_{i=0}^{r-1} \mu_i c(z_i,z_{i+1})
\end{align}

\noindent where $w\in W_{\overline{F}}$, $\lambda_i\in \overline{F}$ and $\mu_i\in \overline{F}$. Denote by $P_{X_V}$ and $P_{-X_W}$ the characteristic polynomials of $X_V$ and $-X_W$ acting on $V_{\overline{F}}$ and $W_{\overline{F}}$ respectively (these are elements of $\overline{E}[T]$). Let $D$ be the $\overline{E}$-linear endomorphism of $\overline{E}[T]$ given by $D(T^{i+1})=T^i$, for $i\geqslant 0$ and $D(1)=0$.

\vspace{2mm}

\begin{prop}\label{proposition 10.3.1}
We have the following equality

\[\begin{aligned}
\displaystyle P_{X_V}(T)= & \sum_{j=0}^{m-1} (-1)^r h(w,X_W^j w) D^{j+1}\big(P_{-X_W}(T)\big) + \\
 & P_{-X_W}(T) \bigg(T^{2r+1}+ \sum_{j=0}^r (-1)^{j+1}2\lambda_j \eta T^{2r-2j}+\sum_{j=0}^{r-1} (-1)^{j+1} 2\mu_j T^{2r-1-2j}\bigg)
\end{aligned}\]

\noindent (Recall that $m=\dim(W)$).
\end{prop}

\vspace{2mm}

\noindent\ul{Proof}: This can be proved by induction on $r$. The computation, fastidious but direct, is left to the reader. $\blacksquare$

\vspace{2mm}

\begin{cor}\label{corollary 10.3.1}
The following $U(W)$-invariant polynomial functions on $\Lambda$

$$X=(X_W,X_V)\mapsto \eta^j h(w,X_W^jw)\in \overline{F}, \;\;\; j=0,\ldots,m-1$$

\noindent (where we have written $X_V$ as in \ref{eq 10.3.1}) extend to $G$-invariant polynomial functions on $\mathfrak{g}$ defined over $F$.
\end{cor}

\vspace{3mm}

\noindent In particular, the polynomial function

$$X=(X_W,X_V)\mapsto \det\left(h\left(X_W^i w, X_W^j w\right)\right)_{0\leqslant i,j\leqslant m-1}\in \overline {F}$$

\noindent extends to a $G$-invariant polynomial function on $\mathfrak{g}$ defined over $F$. Let us denote by $Q_0$ such an extension. Set $d^G(X):=\det(1-Ad(X))_{\mid \mathfrak{g}/\mathfrak{g}_X}$ for all $X\in \mathfrak{g}_{\reg}$. Then $d^G$ extends uniquely to a polynomial $d^G\in F[\mathfrak{g}]^G$. Set $\gls{Q}=Q_0d^G\in F[\mathfrak{g}]^G$ and let $\gls{Lambda'}$ and $\gls{Sigma'}$ be the non-vanishing loci of $Q$ in $\Lambda$ and $\Sigma$ respectively. Notice that we have $\Lambda'\subset \Lambda_{\reg}$ and $\Sigma'\subset \Sigma_{\reg}$ (since $d^G$ divides $Q$).

\subsection{Characterization of $\Sigma'$}\label{section 10.4}

\begin{prop}\label{proposition 10.4.1}
$\Sigma'$ is precisely the set of $X=(X_W,X_V)\in \Sigma_{\reg}$ such that the family

$$z_r,X_V z_r,\ldots, X_V^{d-1}z_r$$

\noindent generates $V_{\overline{F}}$ as a $\overline{E}$-module (Recall that $d=\dim(V)$).
\end{prop}

\vspace{2mm}

\noindent\ul{Proof}: Let $X=(X_W,X_V)\in \Sigma$. It suffices to prove that the family

$$z_r,X_V z_r,\ldots, X_V^{d-1}z_r$$

\noindent generates $V_{\overline{F}}$ as a $\overline{E}$-module if and only if $Q_0(X)\neq 0$. By Proposition \ref{proposition 10.2.1} $X$ is $N$-conjugate to an element of $\Lambda$. Since $Q_0$ is $G$-invariant and $nz_r=z_r$ for all $n\in N$, we may as well assume that $X\in \Lambda$. We assume that it is the case in what follows.

\vspace{2mm}

\noindent By the decomposition \ref{eq 10.3.1}, we see that

$$X_Vz_i=\Xi z_i=z_{i-1}$$

\noindent for all $i=1,\ldots,r$. It follows that

$$\displaystyle \left(z_r,X_Vz_r,\ldots,X_V^r z_r\right)=\left(z_,z_{r-1},\ldots,z_0\right)$$

\noindent Next, again by the decomposition \ref{eq 10.3.1}, it is easy to see that

$$\displaystyle X_Vz_{-i}\equiv \Xi z_{-i}=-z_{-i-1}\;\mod \left\langle z_0,\ldots,z_r\right\rangle$$

\noindent for all $0\leqslant i\leqslant r-1$. Hence, we have

$$\displaystyle \left\langle z_r,\ldots,X_V^{2r}z_r\right\rangle=\left\langle z_r,\ldots,z_1,z_0,z_{-1},\ldots,z_{-r}\right\rangle=W_{\overline{F}}^\perp$$

\noindent and

$$X_V^{2r}z_r\equiv (-1)^r z_{-r}\;\mod W_{\overline{F}}^\perp$$

\noindent It follows that the family

$$z_r,X_V z_r,\ldots, X_V^{d-1}z_r$$

\noindent generates $V_{\overline{F}}$ as a $\overline{E}$-module if and only if the image of the family

$$X_Vz_{-r},\ldots,X_V^mz_{-r}$$

\noindent in $V_{\overline{F}}/W_{\overline{F}}^\perp\simeq W_{\overline{F}}$ is a basis of $W_{\overline{F}}$. But, using again the decomposition \ref{eq 10.3.1} we see that

$$\left(X_Vz_{-r},\ldots,X_V^mz_{-r}\right)\equiv \nu\left(w,-X_Ww,\ldots,\left(-X_W\right)^{m-1}w\right)\;\mod W_{\overline{F}}^\perp$$

\noindent Hence, the family

$$\left( X_Vz_{-r},\ldots,X_V^mz_{-r}\right) \;\mod W_{\overline{F}}^\perp$$

\noindent generates $W_{\overline{F}}$ as an $\overline{E}$-module if and only if the determinant

$$Q_0(X)=\det\left(h\left(X_W^iw,X_W^jw\right)\right)_{0\leqslant i,j\leqslant m-1}\in \overline{F}$$

\noindent is non-zero and this ends the proof of the proposition. $\blacksquare$

\subsection{Conjugacy classes in $\Sigma'$}\label{section 10.5}

\begin{prop}\label{proposition 10.5.1}
The action by conjugation of $H$ on $\Sigma'$ is free and moreover two elements of $\Sigma'$ are $G$-conjugate if and only if they are $H$-conjugate.
\end{prop}

\vspace{2mm}

\noindent\ul{Proof}: Recall that by definition, $H$ acts freely on $\Sigma'$ if the map

$$H\times \Sigma'\to \Sigma'\times \Sigma'$$
$$(h,X)\mapsto (X,hXh^{-1})$$

\noindent is a closed immersion. Because of Proposition \ref{proposition 10.2.1}, this is equivalent to proving that

\begin{align}\label{eq 10.5.1}
& U(W)\times \Lambda'\to \Lambda'\times \Lambda' \\
\nonumber & (h,X)\mapsto (X,hXh^{-1})
\end{align}

\noindent is a closed immersion. For $X=(X_W,X_V)\in \mathfrak{g}$, we define the characteristic polynomial of $X$ to be the pair $P_X=\left(P_{X_W},P_{X_V}\right)$. Let $\mathcal{Y}\subset \Lambda'\times \Lambda'$ be the closed subset of pairs $(X,X')$ such that $P_X=P_{X'}$. We claim the following

\vspace{3mm}

\begin{num}
\item\label{eq 10.5.2} The map \ref{eq 10.5.1} is a closed immersion whose image is $\mathcal{Y}$.
\end{num}

\vspace{3mm}

\noindent This will prove the two points of the proposition (if two elements of $\mathfrak{g}$ are $G$-conjugate, they share the same characteristic polynomial). First, of course, the image of \ref{eq 10.5.1} is contained in $\mathcal{Y}$. Let $(X,X')\in \mathcal{Y}$. We may write

$$X_V=\Xi-X_W+c(z_r,w)+\sum_{i=0}^r \lambda_i c(z_i,\eta z_i)+ \sum_{i=0}^{r-1} \mu_i c(z_i,z_{i+1})$$

$$X'_V=\Xi-X'_W+c(z_r,w')+\sum_{i=0}^r \lambda'_i c(z_i,\eta z_i)+ \sum_{i=0}^{r-1} \mu'_i c(z_i,z_{i+1})$$

\noindent where $X=(X_V,X_W)$, $X'=(X'_V,X'_W)$, $w,w'\in W_{\overline{F}}$, $\lambda_i,\lambda'_i\in \overline{F}$ and $\mu_i,\mu'_i\in \overline{F}$. By Proposition \ref{proposition 10.3.1}, we have $\lambda_i=\lambda_i'$ for $i=0,\ldots,r$, $\mu_i=\mu'_i$ for $i=0,\ldots,r-1$ and

\begin{align}\label{eq 10.5.3}
h(w,X_W^iw)=h(w',{X'}_W^i w') \mbox{ for } i=0,\ldots,m-1.
\end{align}

\noindent Moreover, by definition of $\Lambda'$, $\left(w,X_W w,\ldots,X_W^{m-1} w\right)$ and $\left(w',X'_W w',\ldots, {X'}_W^{m-1}w'\right)$ are basis of $W_{\overline{F}}$. Let $g$ be the unique $\overline{E}$-linear automorphism of $W_{\overline{F}}$ sending $X_W^i w$ to ${X'}_W^i w'$ for all $i=0,\ldots,m-1$. By \ref{eq 10.5.3}, we have $g\in U(W)$ and we easily check that $gXg^{-1}=X'$. It is also easy to see that $g$ is the only element of $U(W)$ with this property. Hence, we have proved that the map \ref{eq 10.5.1} induces a bijection from $U(W)\times \Lambda'$ to $\mathcal{Y}$ and we have constructed the inverse, which is obviously a morphism of algebraic varieties. This proves the claim \ref{eq 10.5.2}. $\blacksquare$ 

\vspace{3mm}

\begin{cor}\label{corollary 10.5.1}
We have an inequality

$$\sigma_G(t) \ll \sigma_{H\backslash G}(t) \sigma_{\Sigma'}(X)$$

\noindent for all $X\in \Sigma'$ and all $t\in G_X$.
\end{cor}

\noindent \ul{Proof}: First, we prove that

\begin{align}\label{eq 10.5.4}
\sigma_G(h)\ll \sigma_{\Sigma'}(X) +\sigma_{\mathfrak{g}}(hXh^{-1})
\end{align}

\noindent for all $h\in H$, $X\in \Sigma'$. From the previous proposition, we know that

$$H\times \Sigma'\to \Sigma'\times \Sigma'$$
$$(h,X)\mapsto (X,hXh^{-1})$$

\noindent is a closed immersion. Hence, we have

$$\sigma_G(h)\ll\sigma_{\Sigma'}(X)+\sigma_{\Sigma'}(hXh^{-1})$$

\noindent for all $X\in \Sigma'$ and $h\in H$. Moreover, since $\Sigma'$ is the principal Zariski open subset of $\Sigma$ defined by the non-vanishing of the polynomial $Q$ and since $Q$ is $G$-invariant, we have

$$\sigma_{\Sigma'}(hXh^{-1})\sim \sigma_{\mathfrak{g}}(hXh^{-1})+\log(2+\lvert Q(X)\rvert^{-1})$$

\noindent for all $h\in H$, $X\in \Sigma'$. Combining this with the inequality

$$\log\left(2+\lvert Q(X)\rvert^{-1}\right)\ll \sigma_{\Sigma'}(X)$$

\noindent for all $X\in \Sigma'$, we get \ref{eq 10.5.4}. We now deduce the following inequality:

\begin{align}\label{eq 10.5.5}
\sigma_G(h)\ll \sigma_{\Sigma'}(X)+\sigma_G(ht)
\end{align}

\noindent  for all $h\in H$, $X\in \Sigma'$, $t\in G_X$. Indeed, since $(g,X)\in G\times \mathfrak{g}\mapsto gXg^{-1}\in \mathfrak{g}$ is a regular map, we have

$$\sigma_{\mathfrak{g}}(hXh^{-1})=\sigma_{\mathfrak{g}}(htX(ht)^{-1})\ll  \sigma_G(ht)+\sigma_{\mathfrak{g}}(X)$$

\noindent for all $h\in H$, $X\in \Sigma$ and $t\in G_X$. Combining this with \ref{eq 10.5.4} and the inequality $\sigma_{\mathfrak{g}}(X)\ll \sigma_{\Sigma'}(X)$ for all $X\in \Sigma'$, we get \ref{eq 10.5.5}.

\vspace{2mm}

\noindent We are now in position to prove the lemma. By \ref{eq 10.5.5}, we have the following chain of inequalities

$$\sigma_G(t)=\sigma_G(h^{-1}ht)\ll \sigma_G(h)+\sigma_G(ht)\ll \sigma_{\Sigma'}(X)+\sigma_{G}(ht)\ll \sigma_{\Sigma'}(X)\sigma_{G}(ht)$$

\noindent for all $h\in H$, $X\in \Sigma'$, $t\in G_X$. By Lemma \ref{lemma 6.2.1}(i), taking the infimum over $h\in H$ gives the desired result. $\blacksquare$

\subsection{Borel subalgebras and $\Sigma'$}\label{section 10.6}

\begin{prop}\label{proposition 10.6.1}
Let $X\in \Sigma'$ and $\mathfrak{b}$ be a Borel subalgebra of $\mathfrak{g}$ (defined over $\overline{F}$) containing $X$, then

$$\mathfrak{b}\oplus \mathfrak{h}=\mathfrak{g}$$
\end{prop}

\noindent \ul{Proof}: Let $X\in \Sigma'$ and $\mathfrak{b}\subseteq \mathfrak{g}$ be a Borel subalgebra containing $X$. By Proposition \ref{proposition 10.2.1}, up to $N(F)$-conjugation we may assume that $X\in \Lambda'$ and we will assume this is so henceforth. Write $X=(X_W,X_V)$ with $X_W\in \mathfrak{u}(W)$ and $X_V\in \mathfrak{u}(V)$. By definition of $\Lambda$, we have a decomposition

\begin{align}\label{eq 10.6.1}
\displaystyle X_V=\Xi-X_W+c(z_r,w)+\sum_{i=0}^r \lambda_i c(z_i,\eta z_i)+\sum_{i=0}^{r-1} \mu_i c(z_i,z_{i+1})
\end{align}

\noindent where $w\in W_{\overline{F}}$, $\lambda_i\in \overline{F}$, $0\leqslant i\leqslant r$ and $\mu_i\in \overline{F}$, $0\leqslant i\leqslant r-1$.

\vspace{2mm}

\noindent It is easy to check that $\dim(\mathfrak{b})+\dim(\mathfrak{h})=\dim(\mathfrak{g})$ so that it suffices to prove

\begin{align}\label{eq 10.6.2}
\mathfrak{b}\cap \mathfrak{h}=0
\end{align}

\noindent There exist Borel subalgebras $\mathfrak{b}_W$ and $\mathfrak{b}_V$ of $\mathfrak{u}(W)$ and $\mathfrak{u}(V)$ respectively such that $X_W\in \mathfrak{b}_W$, $X_V\in\mathfrak{b}_V$ and $\mathfrak{b}=\mathfrak{b}_W\times \mathfrak{b}_V$. Then obviously \ref{eq 10.6.2} is equivalent to

\begin{align}\label{eq 10.6.3}
\left(\mathfrak{b}_W+\mathfrak{n}\right)\cap \mathfrak{b}_V=0
\end{align}

\noindent Fix an $F$-embedding $E\hookrightarrow \overline{F}$ and set $\overline{V}=V\otimes_E \overline{F}$, $\overline{W}=W\otimes_E \overline{F}$. Denote by $\overline{V}^*$ and $\overline{W}^*$ the $\overline{F}$-dual of $\overline{V}$ and $\overline{W}$ respectively. Then, we have isomorphisms of $\overline{F}$-vector space

$$V_{\overline{F}}\simeq \overline{V}\oplus \overline{V}^*$$
$$W_{\overline{F}}\simeq \overline{W}\oplus \overline{W}^*$$

\noindent sending $v\otimes_F \lambda$ and $w\otimes_F \lambda$ to $(v\otimes_E \lambda, h(v,.)\otimes_E\lambda)$ and $(w\otimes_E \lambda, h(w,.)\otimes_E\lambda)$ respectively. For $v\in V$, we will denote by $\overline{v}$ and $\overline{v}^*$ the image of $v$ in $\overline{V}$ and $\overline{V}^*$ respectively. Also, if $U$ is a subspace of $V$ we will set $\overline{U}=U\otimes_E \overline{F}$ and see it as a subspace of $\overline{V}$. We will adopt similar notation with respect to $\overline{W}$. We have an isomorphism

$$\mathfrak{u}(V)_{\overline{F}}\simeq \mathfrak{gl}(\overline{V})$$

\noindent which sends $X\in \mathfrak{u}(V)_{\overline{F}}$ to its restriction to $\overline{V}$ (the inverse is given by mapping $X\in \mathfrak{gl}(\overline{V})$ to the endomorphism of $V_{\overline{F}}$ acting as $X$ on $\overline{V}$ and as $-{}^t X$ on $\overline{V}^*$). Similarly, we have an isomorphism

$$\mathfrak{u}(W)_{\overline{F}}\simeq \mathfrak{gl}(\overline{W})$$

\noindent and we will use these isomorphisms as identifications. Then, $\mathfrak{b}_W$ is the stabilizer in $\mathfrak{gl}(\overline{W})$ of a complete flag

$$0=\overline{W}_0\subsetneq \overline{W}_1\subsetneq\ldots\subsetneq \overline{W}_m=\overline{W}$$

\noindent and similarly $\mathfrak{b}_V$ is the stabilizer in $\mathfrak{gl}(\overline{V})$ of a complete flag

$$\displaystyle \mathcal{F}:\;\;\; 0=\overline{V}_0\subsetneq \overline{V}_1\subsetneq\ldots\subsetneq \overline{V}_d=\overline{V}$$

\noindent Let us define another complete flag 

$$\displaystyle \mathcal{F}':\;\;\; 0=\overline{V}'_0\subsetneq \overline{V}'_1\subsetneq\ldots\subsetneq \overline{V}'_d=\overline{V}$$

\noindent of $\overline{V}$ by setting:

\vspace{2mm}

\begin{itemize}
\renewcommand{\labelitemi}{$\bullet$}
\item $\overline{V}'_i=<\overline{z}_r,\ldots,\overline{z}_{r-i+1}>$ for $i=1,\ldots,r+1$;

\item $\overline{V}'_{r++1+i}=\overline{Z}_+\oplus \overline{D}\oplus \overline{W}_i$ for $i=1,\ldots,m$;

\item $\overline{V}'_{r+m+1+i}=\overline{Z}_+\oplus \overline{V}_0\oplus <\overline{z}_{-1},\ldots,\overline{z}_{-i}>$ for $i=1,\ldots,r$.
\end{itemize}

\vspace{2mm}

\noindent For all $v\in \overline{V}$, let us denote by $V(X_V,v)$ the subspace of $\overline{V}$ generated by $v,X_Vv,X_V^2v,\ldots$. We will need the following lemma

\vspace{2mm}

\begin{lem}\label{lemma 10.6.1}
Let $1\leqslant i\leqslant d$, then we have

\begin{enumerate}[(i)]
\item For all $v\in \overline{V}'_i$ which is nonzero, $\overline{V}'_{i-1}+V(X_V,v)=\overline{V}$;

\item $\overline{V}'_i\cap \overline{V}_{d-i}=0$
\end{enumerate}
\end{lem}

\vspace{2mm}

\noindent\ul{Proof}: First we prove that (i) implies (ii). Indeed, if $v\in \overline{V}'_i\cap \overline{V}_{d-i}$ is nonzero, then by (i), we would have $\dim(V(X_V,v))\geqslant d+1-i$. But $V(X_V,v)\subset \overline{V}_{d-i}$ (since $v\in \overline{V}_{d-i}$ and $X_V\in \mathfrak{b}_V$ preserves $\overline{V}_{d-i}$), and so $\dim(V(X_V,v))\leqslant dim(\overline{V}_{d-i})=d-i$. This is a contradiction.

\vspace{2mm}

\noindent We now turn to the proof of (i). Let $v\in \overline{V}_i'$ be non-zero. Without loss of generality, we may assume that $v\in \overline{V}'_i\backslash \overline{V}'_{i-1}$ since otherwise the result with $i-1$ instead of $i$ is stronger. We assume this is so henceforth and it follows that

\begin{align}\label{eq 10.6.4}
\overline{V}'_{i-1}+V(X_V,v)=\overline{V}_i'+V(X_V,v)
\end{align}

\noindent Obviously $\overline{z}_r\in \overline{V}_i'+V(X_V,v)$ and so by Proposition \ref{proposition 10.4.1}, it suffices to show that $\overline{V}'_{i-1}+V(X_V,v)$ is $X_V$-stable. The subspaces $V(X_V,v)$ is $X_V$-stable almost by definition. Hence, we are left with proving that

\begin{align}\label{eq 10.6.5}
X_V\overline{V}'_{i-1}\subseteq \overline{V}'_{i}+V(X_V,v)
\end{align}

\noindent This is clear if $1\leqslant i\leqslant r+1$ or $r+m+2\leqslant i\leqslant d=2r+m+1$ since in this cases using the decomposition \ref{eq 10.6.1} we easily check that $X_V\overline{V}_{i-1}'\subseteq \overline{V}'_i$. It remains to show that \ref{eq 10.6.5} holds for $r+2\leqslant i\leqslant r+m+1$. In this cases, again using the decomposition \ref{eq 10.6.1}, we easily check that 

\begin{align}\label{eq 10.6.6}
X_Vv'\in \overline{V}'_i+\langle \overline{z}_0^*,v'\rangle X_V\overline{z}_0
\end{align}

\noindent for all $v'\in \overline{V}'_i$ (where $X_V\overline{z}_0=-\overline{z}_{-1}+2\nu\eta\lambda_0\overline{z}_0+\mu_0\nu\overline{z}_1$ if $r\geqslant 1$ and $X_V\overline{z}_0=\nu\overline{w}+2\nu\eta\lambda_0\overline{z}_0$ if $r=0$). Here, we have used the fact that $X_W\in \mathfrak{b}_W$ so that $\overline{W}_{i-r-2}$ and $\overline{W}_{i-r-1}$ are $X_W$-stable. As $v\in \overline{V}'_i$, it suffices to show that the existence of $k\geqslant 0$ such that $\langle \overline{z}^*_0, X_V^kv\rangle\neq 0$. By Proposition \ref{proposition 10.4.1}, the family

$$\overline{z}_r^*,{}^tX_V\overline{z}_r^*,{}^tX_V^2\overline{z}_r^*,\ldots$$

\noindent generates $\overline{V}^*$. Hence, since $v\neq 0$, there exists $k_0\geqslant 0$ such that $\langle {}^t X_V^{k_0}\overline{z}_r^*,v\rangle=\langle \overline{z}_r^*,X_V^{k_0}v\rangle\neq 0$. This already settles the case where $r=0$. In the case $r\geqslant 1$, since $\overline{V}'_i$ is included in the kernel of $\overline{z}_r^*$ this shows that the sequence $v,X_Vv,X_V^2v,\ldots$ eventually escapes from $\overline{V}'_i$ and by \ref{eq 10.6.6} this implies also the existence of $k\geqslant 0$ such that $\langle\overline{z}_0^*,X_V^kv\rangle\neq 0$. This ends the proof of \ref{eq 10.6.5} and of the lemma. $\blacksquare$ 

\vspace{2mm}

\noindent Let us now set $D_i=\overline{V}'_i\cap \overline{V}_{d+1-i}$ for $i=1,\ldots,d$. By the previous lemma, these are one dimensional subspaces of $\overline{V}$ and we have

$$\displaystyle \overline{V}=\bigoplus_{i=1}^d D_i$$

\noindent Let $Y\in  \big(\mathfrak{b}_W+\mathfrak{n}\big)\cap \mathfrak{b}_V$. We want to prove that $Y=0$ (to get \ref{eq 10.6.2}). Obviously $Y$ must stabilize the flags $\mathcal{F}$ and $\mathcal{F}'$ so that $Y$ stabilizes the lines $D_1,\ldots,D_d$. We claim that

\begin{align}\label{eq 10.6.7}
\displaystyle Y(D_i)=0 \mbox{ for all } i=1,\ldots,r+1 \mbox{ and all } i=r+m+2,\ldots,d.
\end{align}

\noindent Indeed, since $Y\in \mathfrak{b}_W+\mathfrak{n}$, we have $Y\overline{V}_i'\subseteq \overline{V}_{i-1}'$ for all $i=1,\ldots,r+1$ and all $i=r+m+2,\ldots,d$ and so $YD_i\subseteq \overline{V}'_{i-1}\cap \overline{V}_{d+1-i}=0$ (by the previous lemma).

\vspace{2mm}

\noindent To deduce that $Y=0$, it only remains to show that

\begin{align}\label{eq 10.6.8}
\displaystyle Y(D_i)=0 \mbox{ for all } i=r+2,\ldots,r+m+1.
\end{align}

\noindent Assume, by way of contradiction, that there exists $1\leqslant j\leqslant m$ such that $YD_{r+1+j}\neq 0$. Since $Y\in \mathfrak{b}_W+\mathfrak{n}$, we have $Y\overline{V}'_{r+1+j}\subseteq \overline{W}_j\oplus \overline{Z}_+$ so that $D_{r+1+j}=YD_{r+1+j}\subseteq \overline{W}_j\oplus \overline{Z}_+$. Let $v\in D_{r+1+j}$ be non-zero. We claim that

\begin{align}\label{eq 10.6.9}
\left(\overline{Z}_+\oplus \overline{W}_j\right)+ V(X_V,v)=\overline{V}
\end{align}

\noindent Indeed by the previous lemma, it suffices to prove that $\overline{z}_0\in \left(\overline{Z}_+\oplus \overline{W}_j\right)+ V(X_V,v)$. By the decomposition \ref{eq 10.6.1}, we easily check that

$$X_V\left(\overline{Z}_+\oplus\overline{W}_j\right)\subseteq \overline{Z}_+\oplus \overline{W}_j\oplus \overline{F}\overline{z}_0$$

\noindent so that we only need to check that the sequence $v,X_Vv,\ldots$ eventually escapes $\overline{Z}_+\oplus \overline{W}_j$. From Proposition \ref{proposition 10.4.1}, we know that there exists $k\geqslant 0$ such that $\langle \overline{z}_r^*, X_V^kv\rangle\neq 0$. Since $\overline{Z}_+\oplus \overline{W}_j$ is included in the kernel of $\overline{z}_r^*$, this proves \ref{eq 10.6.9}.

\vspace{2mm}

\noindent From \ref{eq 10.6.9}, we deduce that $\dim V(X_V,v)\geqslant 1+d-j-r$. On the other hand, we have $v\in \overline{V}_{d-r-j}$ (since $v\in \overline{D}_{r+1+j}$) and $X_V$ leaves $\overline{V}_{d-r-j}$ stable (since $X_V\in \mathfrak{b}_V$). As $\dim \overline{V}_{d-r-j}=d-r-j$ it is a contradiction. This ends the proof of \ref{eq 10.6.8} and of the proposition. $\blacksquare$

\vspace{2mm}

\noindent Recall that we have fixed a (classical) norm $\lvert.\rvert_{\mathfrak{g}}$ on $\mathfrak{g}$ and that for all $R>0$, $B(0,R)$ denotes the closed ball of radius $R$ centered at the origin in $\mathfrak{g}(F)$.

\begin{cor}\label{corollary 10.6.1}
There exists a $c>0$ such that, for all $\epsilon>0$ sufficiently small, all $X\in \Sigma'(F)$ and all parabolic subalgebras $\mathfrak{p}$ of $\mathfrak{g}$ defined over $F$ and containing $X$, we have

$$\displaystyle \exp\left[ B\left(0,\epsilon e^{-c\sigma_{\Sigma'}(X)}\right)\right]\subseteq H(F) \exp\left( B(0,\epsilon)\cap \mathfrak{p}(F)\right)$$
\end{cor}

\noindent\ul{Proof}: 
\begin{lem}
There exists a constant $\epsilon_0>0$ such that for every subspace $\mathcal{V}$ of $\mathfrak{g}(F)$ and all $c_0\geqslant 1$ satisfying

$$B(0,1)\subseteq B(0,c_0)\cap \mathfrak{h}(F)+B(0,c_0)\cap \mathcal{V}$$

\noindent we have

$$\displaystyle \exp\left( B(0,\frac{\epsilon}{2c_0})\right)\subseteq H(F).\exp\left( B(0,\epsilon)\cap\mathcal{V}\right)$$

\noindent for all $0<\epsilon<\epsilon_0$.
\end{lem}

\noindent \ul{Proof}: Use the Campbell-Hausdorff formula and successive approximations. $\blacksquare$

\vspace{2mm}

\noindent Because of this lemma, it suffices to prove the following:

\vspace{3mm}

\begin{num}
\item\label{eq 10.6.10} There exists $c_1>0$ such that for all $X\in \Sigma'(F)$ and all parabolic subalgebras $\mathfrak{p}$ of $\mathfrak{g}$ defined over $F$ and containing $X$, we have
$$\displaystyle B(0,1)\subseteq B\left(0,e^{c_1\sigma_{\Sigma'}(X)}\right)\cap \mathfrak{h}(F)+ B\left(0,e^{c_1\sigma_{\Sigma'}(X)}\right)\cap \mathfrak{p}(F)$$
\end{num}

\vspace{3mm}

\noindent Let us denote by $\mathcal{B}$ the variety of all Borel subalgebras of $\mathfrak{g}$. Let $\mathcal{X}$ be the closed subvariety of $\mathfrak{g}_{\reg}\times \mathcal{B}$ defined by

$$\displaystyle \mathcal{X}:=\{(X,\mathfrak{b})\in \mathfrak{g}_{\reg}\times \mathcal{B}; \;\;\; X\in \mathfrak{b}\}$$

\noindent We will denote by $p$ the natural projection $\mathcal{X}\to \mathfrak{g}_{\reg}$ and by $\mathcal{X}_\Sigma$ the inverse image by $p$ of $\Sigma'$. By the previous proposition, for all $(X,\mathfrak{b})\in \mathcal{X}_\Sigma$ we have $\mathfrak{b}\oplus \mathfrak{h}=\mathfrak{g}$ and we will denote by $p_{\mathfrak{h}}^{\mathfrak{b}}$ (resp.\ $p_{\mathfrak{b}}^{\mathfrak{h}}$) the projection with range $\mathfrak{b}$ (resp.\ $\mathfrak{h}$) and kernel $\mathfrak{h}$ (resp.\ $\mathfrak{b}$). Denote by $\lvert .\rvert$ the subordinate norm on $\End_{\overline{F}}(\mathfrak{g})$ coming from the norm we fixed on $\mathfrak{g}$. We first prove the following fact

\vspace{3mm}

\begin{num}
\item\label{eq 10.6.11} There exists $c_2>0$ such that
$$\lvert p_{\mathfrak{h}}^{\mathfrak{b}}\rvert+\lvert p_{\mathfrak{b}}^{\mathfrak{h}}\rvert \leqslant e^{c_2\sigma_{\Sigma'}(X)}$$

for all $(X,\mathfrak{b})\in \mathcal{X}_{\Sigma}$.
\end{num}

\vspace{3mm}

\noindent Since the map $(X,\mathfrak{b})\in \mathcal{X}_\Sigma\mapsto (p_{\mathfrak{h}}^{\mathfrak{b}},p_{\mathfrak{b}}^{\mathfrak{h}})\in \End_{\overline{F}}(\mathfrak{g})^2$ is regular, we have an inequality

$$\log\left(\lvert p_{\mathfrak{h}}^{\mathfrak{b}}\rvert+ \lvert p_{\mathfrak{b}}^{\mathfrak{h}}\rvert\right)\ll \sigma_{\mathcal{X}_\Sigma}(X,\mathfrak{b})$$

\noindent for all $(X,\mathfrak{b})\in \mathcal{X}_\Sigma$. Moreover, the morphism $p:\mathcal{X}\to \mathfrak{g}_{\reg}$ is finite \'etale and therefore, so is its restriction $p_{\Sigma}:\mathcal{X}_\Sigma\to \Sigma'$. Therefore, by Lemma \ref{lemma 1.2.1}, we also have an inequality

$$\sigma_{\mathcal{X}_\Sigma}(X,\mathfrak{b})\ll \sigma_{\Sigma'}(X)$$

\noindent for all $(X,\mathfrak{b})\in \mathcal{X}_\Sigma$. Combining this with the previous inequality, we get \ref{eq 10.6.11}.

\vspace{2mm}

\noindent We will also need the following fact, whose easy proof is left to the reader.

\vspace{3mm}

\begin{num}
\item\label{eq 10.6.12} There exists a finite Galois extension $K$ of $F$, contained in $\overline{F}$, such that for all $(X,\mathfrak{b})\in \mathcal{X}$ with $X\in \mathfrak{g}_{\reg}(F)$, the Borel subalgebra $\mathfrak{b}$ is defined over $K$.
\end{num}

\vspace{3mm}

\noindent We are now able to prove the corollary. Fix a Galois extension $K$ as in \ref{eq 10.6.12} and a constant $c_2$ as in \ref{eq 10.6.11}. For $R>0$, we will denote by $B_K(0,R)$ the closed ball of radius $R$ centered at the origin in $\mathfrak{g}(K)$ (recall that $\lvert .\rvert_{\mathfrak{g}}$ is defined on all of $\mathfrak{g}=\mathfrak{g}(\overline{F})$). Let $X$ and $\mathfrak{p}$ be as in the corollary and let $\mathfrak{b}\in \mathcal{B}$ be such that $X\in \mathfrak{b}$ and $\mathfrak{b}\subset \mathfrak{p}$. By \ref{eq 10.6.11} and \ref{eq 10.6.12}, we have

\begin{align}\label{eq 10.6.13}
\displaystyle B(0,1)\subseteq B_K \left(0,e^{c_2\sigma_{\Sigma'}(X)}\right)\cap \mathfrak{b}(K)\oplus B_K\left(0,e^{c_2\sigma_{\Sigma'}(X)}\right)\cap \mathfrak{h}(K)
\end{align}

\noindent Denote by $P$ the projection $\mathfrak{g}(K)\to \mathfrak{g}(F)$ given by

$$\displaystyle P(Y)=\frac{1}{d_K}\sum_{\sigma\in Gal(K/F)} \sigma(Y)$$

\noindent where $d_K=[K:F]$. Let $\lvert P\rvert$ be the subordinate norm of $P$ relative to the norms on $\mathfrak{g}(K)$ and $\mathfrak{g}(F)$ (obtained by restrictions of that on $\mathfrak{g}$). Then applying $P$ to the inclusion \ref{eq 10.6.13}, we get

$$\displaystyle B(0,1)\subseteq B_K \left(0,\lvert P\rvert e^{c_2\sigma_{\Sigma'}(X)}\right)\cap \mathfrak{b}(K)\oplus B_K\left(0,\lvert P\rvert e^{c_2\sigma_{\Sigma'}(X)}\right)\cap \mathfrak{h}(K)$$

\noindent Since $P(\mathfrak{b}(K))\subseteq \mathfrak{p}(F)$, we get \ref{eq 10.6.10}. $\blacksquare$

\subsection{The quotient $\Sigma'(F)/H(F)$}\label{section 10.7}

\noindent By Proposition \ref{proposition 10.2.1}, we know that $\Sigma'$ has a geometric quotient by $N$ and that $\Sigma'/N\simeq \Lambda'$. Because $H=N\rtimes U(W)$ and $U(W)$ is reductive, the geometric quotient of $\Sigma'$ by $H$ exists and we have $\Sigma'/H\simeq \Lambda'/U(W)$. Denote by $\mathfrak{g}'$ the non-vanishing locus of $Q$ in $\mathfrak{g}$ and by $\mathfrak{g}'/G$ the geometric quotient of $\mathfrak{g}'$ by $G$ for the adjoint action. The natural map $\Sigma'\to \mathfrak{g}'/G$ factors through the quotient $\Sigma'/H$ and we will denote by

$$\pi: \Sigma'/H\to \mathfrak{g}'/G$$

\noindent the induced morphism. We will also consider the $F$-analytic counterpart of this map:

$$\pi_F: \Sigma'(F)/H(F)\to \mathfrak{g}'(F)/G(F)$$

\noindent Recall that we put on $H(F)$ the Haar measure $\mu_H$ which lift the Haar measure $\mu_\mathfrak{h}$ on $\mathfrak{h}(F)$. Because $H(F)$ acts freely on $\Sigma'(F)$, we can define a measure $\mu_{\Sigma'/H}$ on $\Sigma'(F)/H(F)$ to be the quotient of (the restriction to $\Sigma'(F)$ of) $\mu_\Sigma$ by $\mu_H$. It is the unique measure on $\Sigma'(F)/H(F)$ such that

$$\displaystyle\int_{\Sigma(F)} \varphi(X)d\mu_\Sigma(X)=\int_{\Sigma'(F)/H(F)} \int_{H(F)} \varphi(h^{-1}Xh)dh d\mu_{\Sigma'/H}(X)$$

\noindent for all $\varphi\in C_c(\Sigma(F))$. Recall also that we have defined in Section \ref{section 1.8} a measure $dX$ on $\mathfrak{g}_{\reg}(F)/G(F)=\Gamma_{\reg}(\mathfrak{g})$. Moreover, $\mathfrak{g}'(F)/G(F)$ is an open subset of $\mathfrak{g}_{\reg}(F)/G(F)$ and we will still denote by $dX$ the restriction of this measure to $\mathfrak{g}'(F)/G(F)$.

\begin{prop}\label{proposition 10.7.1}
\begin{enumerate}[(i)]
\item $\pi$ is an isomorphism of algebraic varieties and $\pi_F$ is an open embedding of $F$-analytic spaces ;

\item $\pi_F$ sends the measure $d\mu_{\Sigma'/H}(X)$ to $D^G(X)^{1/2}dX$;

\item The natural projection $p: \Sigma'\to \Sigma'/H$ has the norm descent property.
\end{enumerate}
\end{prop}

\noindent\ul{Proof}:
\begin{enumerate}[(i)]
\item Both $\Sigma'/H$ and $\mathfrak{g}'/G$ are smooth. By Proposition \ref{proposition 10.5.1} $\pi$ and $\pi_F$ are injective. Moreover, using Proposition \ref{proposition 10.3.1} we see that $\pi$ is surjective and so $\pi$ is bijective. Since the tangent spaces at $X\in \Sigma'(F)/H(F)$ of $\Sigma'(F)/H(F)$ and $\Sigma'/H$ are the same, we only need to prove that $\pi$ is a local isomorphism (i.e. \'etale). Let $X\in \Sigma'$. The tangent space of $\Sigma'/H$ at $X$ is

$$T_X \Sigma'/H=\mathfrak{h}^\perp/\ad(X)(\mathfrak{h})$$

\noindent and the tangent space of $\mathfrak{g}'/G$ at $X$ is

$$T_X \left(\mathfrak{g}'/G\right)=\mathfrak{g}/\ad(X)(\mathfrak{g})$$

\noindent Modulo these identifications, the differential of $\pi$ at $X$ is the natural inclusion of $\mathfrak{h}^\perp/\ad(X)(\mathfrak{h})$ in $\mathfrak{g}/\ad(X)(\mathfrak{g})$. We want to prove that it is an isomorphism. Choose a Borel subalgebra $\mathfrak{b}$ of $\mathfrak{g}$ that contains $X$ and denote by $\mathfrak{u}$ its nilpotent radical. By Proposition \ref{proposition 10.6.1}, we have $\mathfrak{g}=\mathfrak{h}\oplus \mathfrak{b}$ so that $\mathfrak{g}=\mathfrak{h}^\perp\oplus \mathfrak{b}^\perp$. But, we have $\mathfrak{b}^\perp=\mathfrak{u}\subset \ad(X)(\mathfrak{g})$ and this establishes the surjectivity. On the other hand, from the equality $\mathfrak{g}=\mathfrak{h}\oplus \mathfrak{b}$, we deduce that

$$\ad(X)(\mathfrak{g})=\ad(X)(\mathfrak{h})\oplus \ad(X)(\mathfrak{b})=\ad(X)(\mathfrak{h})\oplus \mathfrak{u}$$

\noindent We just saw that $\mathfrak{g}=\mathfrak{h}^\perp\oplus \mathfrak{u}$ so that $\mathfrak{h}^\perp\cap \mathfrak{u}=0$. Hence, we have $\mathfrak{h}^\perp\cap \ad(X)(\mathfrak{g})=\ad(X)(\mathfrak{h})$ and this proves the injectivity.

\item In what follows, we will use heavily the notions defined at the end of Section \ref{section 1.6}. Let $X\in \Sigma'(F)$. Set $\mathfrak{g}_X=\Ker(\ad(X))$ and $\mathfrak{g}^X=\Ima(\ad(X))$. Then we have natural identifications (as above)

$$T_X \Sigma'(F)/H(F)=\mathfrak{h}^\perp(F)/\ad(X)(\mathfrak{h}(F))$$

$$T_X \mathfrak{g}'(F)/G(F)=\mathfrak{g}(F)/\mathfrak{g}^X(F)$$

\noindent and the tangent map of $\pi_F$ at $X$, that we will denote by $\iota_X$, is the natural inclusion (and this is an isomorphism). Let $F(X)\in \mathbb{R}_+$ be such that

$$\iota_{X,*} \left(\mu_{\mathfrak{h}}^\perp/\ad(X)_*\mu_{\mathfrak{h}}\right)=F(X)\;\; \mu_{\mathfrak{g}}/\mu_{\mathfrak{g}^X}$$

\noindent where $\mu_{\mathfrak{g}}$ and $\mu_{\mathfrak{g}^X}$ are the autodual measures with respect to $B(.,.)$ on $\mathfrak{g}(F)$ and $\mathfrak{g}^X(F)$ respectively. Then $\pi_F$ sends the measure $d\mu_{\Sigma'/H}(X)$ to $F(X)dX$ and so we have to prove that $F(X)=D^G(X)^{1/2}$. Fix a Borel subalgebra $\mathfrak{b}$ of $\mathfrak{g}$ containing $X$ and let $\mathfrak{u}$ be its nilpotent radical. Then we have $\mathfrak{g}=\mathfrak{h}\oplus \mathfrak{b}$, $\mathfrak{g}=\mathfrak{h}^\perp\oplus \mathfrak{u}$ and $\mathfrak{b}=\mathfrak{g}_X\oplus \mathfrak{u}$. Let $\mu_{\mathfrak{u}}$ be the unique measure on $\mathfrak{u}$ such that

\begin{align}\label{eq 10.7.1}
\mu_{\mathfrak{g}}=\mu_{\mathfrak{h}}^\perp\otimes \mu_{\mathfrak{u}}
\end{align}

\noindent Then by \ref{eq 1.6.2} and \ref{eq 1.6.3}, we have

$$\mu_{\mathfrak{g}}=\mu_{\mathfrak{h}}\otimes \mu_{\mathfrak{u}}^\perp$$

\noindent But, $\mathfrak{u}^\perp=\mathfrak{b}$ and it is easy to see that $\mu_{\mathfrak{u}}^\perp=\mu_{\mathfrak{g}_X}\otimes \mu_{\mathfrak{u}}$, where $\mu_{\mathfrak{g}_X}$ is the autodual measure on $\mathfrak{g}_X(F)$. Hence, we have

\begin{align}\label{eq 10.7.2}
\mu_{\mathfrak{g}}=\mu_{\mathfrak{h}}\otimes \mu_{\mathfrak{g}_X}\otimes \mu_{\mathfrak{u}}
\end{align}

\noindent Let $T$ be the endomorphism of $\mathfrak{g}$ that is equal to $\ad(X)$ on $\mathfrak{h}\oplus \mathfrak{u}$ and to $Id$ on $\mathfrak{g}_X$. Then we have $\det(T)=\det(\ad(X)_{|\mathfrak{g}/\mathfrak{g}_X})$ so that $\lvert \det(T)\rvert=D^G(X)$. Thus, using \ref{eq 10.7.2}, we have

\[\begin{aligned}
D^G(X)\mu_{\mathfrak{g}}=T_*\mu_{\mathfrak{g}} & =\ad(X)_*\mu_{\mathfrak{h}}\otimes \mu_{\mathfrak{g}_X}\otimes \ad(X)_*\mu_{\mathfrak{u}} \\
 & =D^G(X)^{1/2} \left(\ad(X)_*\mu_{\mathfrak{h}}\otimes \mu_{\mathfrak{g}_X}\otimes \mu_{\mathfrak{u}}\right)
\end{aligned}\]

\noindent As $\mu_{\mathfrak{g}}=\mu_{\mathfrak{g}_X}\otimes \mu_{\mathfrak{g}^X}$ (relative to the decomposition $\mathfrak{g}=\mathfrak{g}_X\oplus \mathfrak{g}^X$), this implies $\mu_{\mathfrak{g}^X}=D^G(X)^{1/2} \ad(X)_*\mu_{\mathfrak{h}}\otimes \mu_{\mathfrak{u}}$ (relative to the decomposition $\mathfrak{g}^X=\ad(X)(\mathfrak{h})\oplus \mathfrak{u}$). From this and \ref{eq 10.7.1}, we easily deduce that $F(X)=D^G(X)^{1/2}$.

\item By Proposition \ref{proposition 10.2.1}, it is sufficient to show that

$$\Lambda'\to \Lambda'/U(W)$$

\noindent has the norm descent property. Denote by $\Lambda_{Q_0}$ the non-vanishing locus of $Q_0$ in $\Lambda$ (where $Q_0\in F[\mathfrak{g}]^G$ is defined in Section \ref{section 10.3}). Then we have the following Cartesian diagram where horizontal maps are open immersions

$$\xymatrix{ \Lambda' \ar@{^{(}->}[r] \ar[d] & \Lambda_{Q_0} \ar[d] \\ \Lambda'/U(W) \ar@{^{(}->}[r] & \Lambda_{Q_0}/U(W) }$$

\noindent Thus, if we prove that $\Lambda_{Q_0}\to \Lambda_{Q_0}/U(W)$ has the norm descent property, we will be done. By definition of $\Lambda$ (cf.\ Section \ref{section 10.2}), we have an $U(W)$-equivariant isomorphism

$$\Lambda\simeq \mathfrak{u}(W)\times W\times \mathbb{A}^{2r+1}$$

\noindent where the action of $U(W)$ on the right hand side is the product of the adjoint action on $\mathfrak{u}(W)$, the natural action on $W$ and the trivial action on $\mathbb{A}^{2r+1}$. Denote by $\big(\mathfrak{u}(W)\times W\big)_0$ the open-Zariski subset of $\mathfrak{u}(W)\times W$ consisting of all pairs $(X,w)\in \mathfrak{u}(W)\times W$ such that $\left(w,Xw,X^2w,\ldots\right)$ generates $W$. Then $\Lambda_{Q_0}$ corresponds via the previous isomorphism to $\big(\mathfrak{u}(W)\times W\big)_0\times \mathbb{A}^{2r+1}$. Since $U(W)$ acts trivially on $\mathbb{A}^{2r+1}$, we are reduced to show that

$$\big(\mathfrak{u}(W)\times W\big)_0\to \big(\mathfrak{u}(W)\times W\big)_0/U(W)$$

\noindent has the norm descent property. Let $\mathcal{B}$ be the variety of basis of $W$ and let $Pol_m$ be the variety of monic polynomial $P\in \overline{E}[T]$ of degree $m$. Consider the action of $U(W)$ on $Pol_m\times \mathcal{B}$ which is trivial on $Pol_m$ and given by $g.(e_1,\ldots,e_m)=(ge_1,\ldots,ge_m)$ on $\mathcal{B}$. The map

$$\left(\mathfrak{u}(W)\times W\right)_0\to Pol_m\times \mathcal{B}$$
$$(X,w)\to \left(P_X,w,Xw,\ldots,X^{m-1}w\right)$$

\noindent is a $U(W)$-equivariant closed immersion. Passing to the quotient, we get a commutative diagram

$$\xymatrix{(\mathfrak{u}(W)\times W)_0 \ar@{^{(}->}[r] \ar[d] & Pol_m\times \mathcal{B} \ar[d] \\ (\mathfrak{u}(W)\times W)_0/U(W) \ar@{^{(}->}[r] & Pol_m\times \mathcal{B}/U(W)}$$

\noindent where horizontal maps are closed immersion (since $U(W)$ is reductive). Moreover the diagram is Cartesian (because all $U(W)$-orbits in $\mathcal{B}$ are closed and so the quotient separates all orbits). Thus, we are finally reduced to showing that $\mathcal{B}\to \mathcal{B}/U(W)$ has the norm descent property. Choosing a particular basis of $W$, this amounts to proving that

$$GL(W)\to GL(W)/U(W)$$

\noindent has the norm descent property. Since this map is $GL(W)$-equivariant for the action by left translation, by Lemma \ref{lemma 1.2.2}(i), it suffices to show the existence of a nonempty Zariski open subset of $GL(W)/U(W)$ over which the previous map has the norm descent property. Choose an orthogonal basis $\left(e_1,\ldots,e_m\right)$ of $W$ and denote by $B$ the standard Borel subgroup of $GL(W)$ relative to this basis. Then $B\cap U(W)=Z$ is the subtorus acting by multiplication by an element of $\Ker\; N_{E/F}$ on each $e_i$ (so that $Z\simeq \left(\Ker\; N_{E/F}\right)^m$). Denote $B\times^Z U(W)$ the quotient of $B\times U(W)$ by the $Z$-action given by $z.(b,g)=(bz^{-1},zg)$. Then the multiplication map

$$B\times^Z U(W)\to GL(W)$$
$$(b,g)\mapsto bg$$

\noindent is an open immersion. Thus, it suffices to prove that 

$$B\times^Z U(W)\to \left(B\times^Z U(W)\right)/U(W)=B/Z$$

\noindent has the norm descent property. Let $U$ be the unipotent radical of $B$ and $T$ the subtorus of $B$ stabilizing the lines $\langle e_1\rangle,\ldots,\langle e_m\rangle$. Then the previous map is isomorphic to

$$U\times \left(T\times^Z U(W)\right)\to U\times T/Z$$

\noindent and so we are reduced to showing that $T\times^Z U(W)\to T/Z$ has the norm descent property. Let $S\simeq \left(\mathbb{Z}/2\mathbb{Z}\right)^m$ be the subgroup of elements of $Z$ with eigenvalues $\pm 1$ and let $T_i=\left(\mathbb{G}_m\right)^m$ be the split part of the torus $T=\left(R_{E/F} \mathbb{G}_m\right)^m$. Then we have $T=T_i\times^S Z$ and thus $T\times^Z U(W)=T_i\times^S U(W)$ and $T/Z=T_i/S$. So, we need to prove that the map $T_i\times^S U(W)\to T_i/S$ has the norm descent property. For all $c\in \Ker\left( H^1(F,S)\to H^1(F,U(W))\right)$ let $g_c\in U(W)$ such that $c_\sigma=\sigma(g_c)g_c^{-1}$, for all $\sigma\in \Gamma_F$. Denote by $\mathcal{E}$ the set of all $g_c$ (this is a finite set). Let $t\in T_i$ and denote by $\overline{t}$ its image in $T_i/S$. Assume that $\overline{t}\in \Ima\left( \left(T_i\times^S U(W)\right)(F)\to T_i/S(F)\right)$. Then the $1$-cocycle $\sigma\in \Gamma_F\mapsto c_\sigma=t\sigma(t)^{-1}\in S$ splits in $U(W)$. So there exists (a unique) $g_t\in \mathcal{E}$ such that $c_\sigma=\sigma(g_t)g_t^{-1}$, for all $\sigma\in \Gamma_F$. Moreover, the element $(t,g_t)\in \left(T_i\times^S U(W)\right)(F)$ maps to $\overline{t}$ and we clearly have an inequality

$$\displaystyle \sigma_{T_i\times^S U(W)}(t,g_t)\ll \sigma_{T_i/S}(\overline{t})$$

\noindent for every such $t$. This proves that $T_i\times^S U(W)\to T_i/S$ has the norm descent property and this ends the proof of (iii). $\blacksquare$
\end{enumerate}

\subsection{Statement of the spectral expansion of $J^{\Lie}$}\label{section 10.8}

\noindent Let us define $\gls{GammaSigma}$ to be the subset of $\Gamma(\mathfrak{g})$ consisting of the conjugacy classes of the semi-simple parts of elements in $\Sigma(F)$. We equip this subset with the restriction of the measure defined on $\Gamma(\mathfrak{g})$. Thus, if $\mathcal{T}(G)$ is a set of representatives for the $G(F)$-conjugacy classes of maximal tori in $G$ and if for all $T\in\mathcal{T}(G)$ we denote by $\mathfrak{t}(F)_\Sigma$ the subset of elements $X\in\mathfrak{t}(F)$ whose conjugacy class belongs to $\Gamma(\Sigma)$, then we have

$$\displaystyle \int_{\Gamma(\Sigma)}\varphi(X)dX=\sum_{T\in\mathcal{T}(G)} \lvert W(G,T)\rvert^{-1}\int_{\mathfrak{t}(F)_\Sigma}\varphi(X) dX$$

\noindent for all $\varphi\in C_c(\Gamma(\Sigma))$. Recall that in Section \ref{section 7.2}, we have defined a continuous linear form $J^{\Lie}$ on $\mathcal{S}_{\scusp}(\mathfrak{g}(F))$. The purpose of the next 3 sections is to prove the following theorem.

\vspace{3mm}

\begin{theo}\label{theorem 10.8.1}
We have

$$\displaystyle J^{\Lie}(f)=\int_{\Gamma(\Sigma)}D^G(X)^{1/2} \widehat{\theta}_f(X) dX$$

\noindent for all $f\in \mathcal{S}_{\scusp}(\mathfrak{g}(F))$.
\end{theo}

\vspace{5mm}

\noindent Let $f\in\mathcal{S}_{\scusp}(\mathfrak{g}(F))$. By Theorem \ref{proposition 7.2.1}, Lemma \ref{lemma 5.2.2}(i) and \ref{eq 1.8.2}, both sides of the equality of the theorem are continuous in $f$. Hence, by \ref{eq 5.1.1}, we may assume that $\overline{\Supp(\widehat{f})^G}$ is compact modulo conjugation (this condition is automatic if $F$ is $p$-adic). We assume this is so henceforth.

\subsection{Introduction of a truncation}\label{section 10.9}

We fix a sequence $(\kappa_N)_{N\geqslant 1}$ of functions $\gls{kappaN}:H(F)\backslash G(F)\to [0,1]$ satisfying the following conditions:

\vspace{3mm}

\begin{num}
\item\label{eq 10.9.1} There exist $C_1,C_2>0$ such that for all $x\in H(F)\backslash G(F)$ and all $N\geqslant 1$, we have:
$$\sigma_{H\backslash G}(x)\leqslant C_1N\Rightarrow \kappa_N(x)=1$$
$$\kappa_N(x)\neq 0\Rightarrow \sigma_{H\backslash G}(x)\leqslant C_2N$$

\item\label{eq 10.9.2} If $F$ is $p$-adic, there exists an open-compact subgroup $K'\subset G(F)$ such that the function $\kappa_N$ is right-invariant by $K'$ for all $N\geqslant 1$.

\item\label{eq 10.9.3} If $F=\mathbb{R}$, the functions $\kappa_N$ are smooth and there exists a positive real number $C$ such that
$$\displaystyle \left\lvert \frac{d}{dt}\kappa_N(xe^{tX})_{\mid t=0}\right\rvert\leqslant C.\lvert X\rvert_{\mathfrak{g}}$$
for all $x\in H(F)\backslash G(F)$, all $X\in \mathfrak{g}(F)$ and all $N\geqslant 1$.
\end{num}

\vspace{3mm}

\noindent Such a sequence of truncation functions is easy to construct. Indeed, pick any sequence $(\kappa^0_N)_{N\geqslant 1}$ of measurable functions $\kappa^0_N:H(F)\backslash G(F)\to [0,1]$ that satisfy the condition \ref{eq 10.9.1} above, let $\varphi\in C_c^\infty(G(F))$ be any positive function such that $\displaystyle \int_{G(F)}\varphi(g)dg=1$, then the sequence of functions $\kappa_N=\kappa^0_N\ast \varphi$, $N\geqslant 1$, satisfies the conditions \ref{eq 10.9.1}, \ref{eq 10.9.2} and \ref{eq 10.9.3} above.

\vspace{2mm}

\noindent Set
$$\displaystyle \gls{JLieN}(f)=\int_{H(F)\backslash G(F)} \kappa_N(g)\int_{\mathfrak{h}(F)} f(g^{-1}Xg)\xi(X)dXdg$$
Then, by definition of $J^{\Lie}$, we have

\begin{align}\label{eq 10.9.4}
\displaystyle J^{\Lie}(f)=\lim\limits_{N\to \infty} J^{\Lie}_N(f)
\end{align}

\noindent By \ref{eq 10.1.2}, we have

\begin{align}\label{eq 10.9.5}
\displaystyle J^{\Lie}_N(f)=\int_{H(F)\backslash G(F)}\kappa_N(g) \int_{\Sigma(F)}\widehat{f}(g^{-1}Xg)d\mu_\Sigma(X)dg
\end{align}

\noindent Fix a set $\mathcal{T}(G)$ of representatives for the conjugacy classes of maximal tori in $G$. Recall that in Section \ref{section 10.3}, we have defined a $G$-invariant polynomial function $Q$ on $\mathfrak{g}$. For all $T\in \mathcal{T}(G)$, let us denote by $\mathfrak{t}'$ the principal Zariski-open subset

$$\mathfrak{t}'=\{X\in\mathfrak{t}; \;Q(X)\neq 0\}$$

\noindent and set $\gls{tF'}=\mathfrak{t}'(F)\cap \mathfrak{t}(F)_\Sigma$. Then, $\mathfrak{t}(F)'$ is exactly the subset of elements $X\in\mathfrak{t}(F)$ that are conjugate to some element in $\Sigma'(F)$. Let us fix, for all $T\in\mathcal{T}(G)$, two maps

$$X\in\mathfrak{t}(F)'\mapsto \gls{gammaX}\in G(F)$$

$$X\in\mathfrak{t}(F)'\mapsto \gls{XSigma}\in \Sigma'(F)$$

\noindent such that $\gamma_X^{-1}X\gamma_X=X_\Sigma$ for all $X\in \mathfrak{t}(F)'$. Then, by Proposition \ref{proposition 10.7.1}(i) and (ii), we have

$$\displaystyle \int_{\Sigma(F)}\widehat{f}(g^{-1}Xg)d\mu_\Sigma(X)=\sum_{T\in\mathcal{T}(G)} \lvert W(G,T)\rvert^{-1} \int_{\mathfrak{t}(F)'} D^G(X)^{1/2}\int_{H(F)} \widehat{f}(g^{-1}h^{-1}X_\Sigma hg)dhdX$$

\noindent for all $g\in G(F)$. By inserting this expression in \ref{eq 10.9.5} and switching two integrals, we get

$$\displaystyle J_N^{\Lie}(f)=\sum_{T\in\mathcal{T}(G)} \lvert W(G,T)\rvert^{-1} \int_{\mathfrak{t}(F)'} D^G(X)^{1/2}\int_{T(F)\backslash G(F)} \widehat{f}(g^{-1}Xg)\kappa_{N,X}(g)dgdX$$

\noindent for all $N\geqslant 1$, where we have set

$$\displaystyle \gls{kappaNX}(g)=\int_{T(F)} \kappa_N(\gamma_X^{-1}tg) dt$$

\noindent Define

\begin{align}\label{eq 10.9.6}
\displaystyle \gls{JLieNT}(f)=\int_{\mathfrak{t}(F)'} D^G(X)^{1/2}\int_{T(F)\backslash G(F)} \widehat{f}(g^{-1}Xg)\kappa_{N,X}(g)dgdX,\;\;\; \mbox{for } N\geqslant 1
\end{align}

\noindent for all $T\in \mathcal{T}(G)$, so that

\begin{align}\label{eq 10.9.7}
\displaystyle J_N^{\Lie}(f)=\sum_{T\in\mathcal{T}(G)} \left\lvert W(G,T)\right\rvert^{-1} J^{\Lie}_{N,T}(f),\;\;\; \mbox{for } N\geqslant 1
\end{align}

\noindent We fix from now on a torus $T\in \mathcal{T}(G)$. The previous formal manipulation (interchange of two integrals) will be justified by the next lemma proving the absolute convergence of \ref{eq 10.9.6}. But first we need to prove the following:

\vspace{3mm}

\begin{num}
\item\label{eq 10.9.8} We can choose the maps $X\in\mathfrak{t}(F)'\mapsto \gamma_X$ and $X\in \mathfrak{t}(F)'\mapsto X_\Sigma$ so that they satisfy inequalities
$$\displaystyle \sigma_{\Sigma'}(X_\Sigma)\ll \sigma_{\mathfrak{g}}(X)+\log\left(2+\lvert Q(X)\rvert^{-1}\right)$$
$$\displaystyle \sigma_G(\gamma_X)\ll \sigma_{\mathfrak{g}}(X)+\log\left(2+\lvert Q(X)\rvert^{-1}\right)$$
for all $X\in \mathfrak{t}(F)'$.
\end{num}

\vspace{3mm}

\noindent By Proposition \ref{proposition 10.7.1}(i) and (iii), we can choose the map $X\in\mathfrak{t}(F)'\to X_\Sigma\in \Sigma'(F)$ such that

\begin{align}\label{eq 10.9.9}
\sigma_{\Sigma'}(X_\Sigma)\ll \sigma_{\mathfrak{g}'/G}(X)
\end{align}

\noindent for all $X\in \mathfrak{t}(F)'$. Moreover, since $\mathfrak{g}'/G$ is the principal open subset of $\mathfrak{g}/G$ defined by $Q$, we have

\begin{align}\label{eq 10.9.10}
\displaystyle \sigma_{\mathfrak{g}'/G}(X)\sim \sigma_{\mathfrak{g}/G}(X)+\log\left(2+\lvert Q(X)\rvert^{-1}\right)
\end{align}

\noindent for all $X\in\mathfrak{t}(F)'$. But, since $T$ is a torus, we have $\sigma_{\mathfrak{g}/G}(X)\sim \sigma_{\mathfrak{g}}(X)$ for all $X\in \mathfrak{t}(F)$. The first inequality of \ref{eq 10.9.8} is now a consequence of \ref{eq 10.9.9} and \ref{eq 10.9.10}. On the other hand, by \ref{eq 1.2.1} and \ref{eq 1.2.2}, we can choose the map $X\mapsto \gamma_X$ such that

$$\sigma_G(\gamma_X)\ll \sigma_{\mathfrak{g}_{\reg}}(X_\Sigma)$$

\noindent for all $X\in\mathfrak{t}(F)'$. Since we have $\sigma_{\mathfrak{g}_{\reg}}(Y)\ll \sigma_{\Sigma'}(Y)$ for all $Y\in \Sigma'(F)$, the second inequality of \ref{eq 10.9.8} follows from the first. This ends the proof of \ref{eq 10.9.8}.

\vspace{2mm}

\noindent We will assume from now on that the maps $X\mapsto \gamma_X$ and $X\mapsto X_\Sigma$ satisfy the conditions of \ref{eq 10.9.8}. For all $\epsilon>0$, let us denote by $\gls{t'>}$ the set of elements $X\in \mathfrak{t}(F)'$ such that $\lvert Q(X)\rvert> \epsilon$. For all $C>0$, we will denote by $\mathbf{1}_{<C}$ the characteristic function of the set of $X\in\mathfrak{g}(F)$ such that $\sigma_{\mathfrak{g}}(X)<C$. For all $N\geqslant 1$ and all $\epsilon>0$, we define the two following expressions which are similar to \ref{eq 10.9.7}:

$$\displaystyle \lvert J\rvert_{N,T}^{\Lie}(f)=\int_{\mathfrak{t}(F)'} D^G(X)^{1/2}\int_{T(F)\backslash G(F)} \lvert\widehat{f}(g^{-1}Xg)\rvert\kappa_{N,X}(g)dg$$

$$\displaystyle \gls{JNTepsilonf}=\int_{\mathfrak{t}(F)'[>\epsilon]} D^G(X)^{1/2}\int_{T(F)\backslash G(F)}\mathbf{1}_{<\log(N)}(g^{-1}Xg) \widehat{f}(g^{-1}Xg)\kappa_{N,X}(g)dg$$

\noindent Let $\omega_T\subseteq \mathfrak{t}(F)$ be a relatively compact subset such that $\widehat{f}$ is zero on $\mathfrak{t}(F)-\omega_T$ (recall the assumption that $\overline{\Supp(\widehat{f})^G}$ is compact modulo conjugation). The following lemma proves in particular the absolute convergence of \ref{eq 10.9.7}.

\begin{lem}\label{lemma 10.8.1}
\begin{enumerate}[(i)]
\item There exist $k>0$ such that
$$\kappa_{N,X}(g)\ll N^k \log\left(2+\lvert Q(X)\rvert^{-1}\right)^k \sigma_{\mathfrak{g}}(g^{-1}Xg)^k$$
for all $X\in \mathfrak{t}(F)'\cap \omega_T$, all $N\geqslant 1$ and all $g\in G(F)$.

\item There exist $k>0$ such that
$$\lvert J\rvert^{\Lie}_{N,T}(f)\ll N^k$$
for all $N\geqslant 1$.

\item For $b>0$ large enough, we have
$$\left\lvert J^{\Lie}_{N,T}(f)-J^{\Lie}_{N,T,N^{-b}}(f)\right\rvert\ll N^{-1}$$
for all $N\geqslant 1$.
\end{enumerate}
\end{lem}

\noindent\ul{Proof}: (i) implies easily (ii) and (iii) by \ref{eq 1.2.3} and \ref{eq 1.2.4}. By the property \ref{eq 10.9.1} of $\kappa_N$, we have an inequality
$$\sigma_{H\backslash G}(g)\ll N$$
for all $N\geqslant 1$ and all $g\in G(F)$ such that $\kappa_N(g)> 0$. It follows that

\begin{align}\label{eq 10.9.11}
\sigma_{H\backslash G} (\gamma_X^{-1}t\gamma_X)\ll N\sigma_G(g)\sigma_G(\gamma_X)
\end{align}

\noindent for all $N\geqslant 1$, all $X\in \mathfrak{t}(F)'$, all $t\in T(F)$ and all $g\in G(F)$ such that $\kappa_N(\gamma_X^{-1}tg)>0$. On the other hand, since $\gamma_X^{-1}t\gamma_X\in G_{X_\Sigma}(F)$, for all $X\in\mathfrak{t}(F)'$ and all $t\in T(F)$, by Corollary \ref{corollary 10.5.1}, we have
$$\sigma_G(\gamma_X^{-1}t\gamma_X)\ll \sigma_{H\backslash G}(\gamma_X^{-1}t\gamma_X)\sigma_{\Sigma'}(X_\Sigma)$$
for all $X\in \mathfrak{t}(F)'$ and all $t\in T(F)$. Combining this with \ref{eq 10.9.8} and \ref{eq 10.9.11}, and since $\omega_T$ is bounded, we get
$$\sigma_{G}(t)\ll N\sigma_G(g) \log\left(2+\lvert Q(X)\rvert^{-1}\right)$$
for all $N\geqslant 1$, all $X\in \mathfrak{t}(F)'\cap \omega_T$, all $t\in T(F)$ and all $g\in G(F)$ such that $\kappa_N(\gamma_X^{-1}tg)>0$. The function $\kappa_N$ is nonnegative and bounded above by $1$. Hence, it follows from the previous inequality and the definition of $\kappa_{N,X}$ that there exists $c_0>0$ such that
$$\displaystyle \kappa_{N,X}(g)\leqslant \vol\{t\in T(F);\; \sigma_G(t)\leqslant  c_0 N \sigma_G(g)\log\left(2+\lvert Q(X)\rvert^{-1}\right)\}$$
for all $N\geqslant 1$, all $X\in \mathfrak{t}(F)'\cap \omega_T$ and all $g\in G(F)$. It is easy to see that there exists $k>0$ such that
$$\displaystyle \vol\{ t\in T(F); \sigma_G(t)\leqslant M\}\ll M^k$$
for all $M\geqslant 1$. Hence, we get
$$\displaystyle \kappa_{N,X}(g)\ll N^k\log\left(2+\lvert Q(X)\rvert^{-1}\right)^k \sigma_G(g)^k$$
for all $N\geqslant 1$, all $X\in \mathfrak{t}(F)'\cap \omega_T$ and all $g\in G(F)$. Since the function $\kappa_{N,X}$ is invariant by left translation by $T(F)$, we may replace $g$ by $tg$ for any $t\in T(F)$ in the right hand side of the inequality above. By \ref{eq 1.2.1}, taking the infimum over $T(F)$ gives the inequality

\begin{align}\label{eq 10.9.12}
\kappa_{N,X}(g)\ll N^k\log\left(2+\lvert Q(X)\rvert^{-1}\right)^k \sigma_{T\backslash G}(g)^k
\end{align}

\noindent for all $N\geqslant 1$, all $X\in \mathfrak{t}(F)'\cap \omega_T$ and all $g\in G(F)$. By \ref{eq 1.2.2}, we have an inequality

\begin{align}\label{eq 10.9.13}
\displaystyle \sigma_{T\backslash G}(g)\ll \sigma_{\mathfrak{g}}(g^{-1}Xg)\log\left(2+D^G(X)^{-1}\right)
\end{align}

\noindent for all $g\in G(F)$ and all $X\in \mathfrak{t}_{\reg}(F)$. But since the polynomial $d^G$ divides $Q$, we also have

\begin{align}\label{eq 10.9.14}
\displaystyle \log\left(2+D^G(X)^{-1}\right)\ll \log\left(2+\lvert Q(X)\rvert^{-1}\right)
\end{align}

\noindent for all $X\in\mathfrak{t}(F)'\cap\omega_T$. The point (i) now follows from the combination of \ref{eq 10.9.12}, \ref{eq 10.9.13} and \ref{eq 10.9.14}. $\blacksquare$

\vspace{3mm}

\noindent We fix for the moment a positive integer $b>0$ satisfying (iii) of the previous lemma. In Section \ref{section 10.11}, we shall assume (as we may) that $b$ has been chosen sufficiently large (how large will be made precise in \S \ref{section 10.11}).

\subsection{Change of truncation}\label{section 10.10}

\noindent Set $\gls{MT}=Z_G(A_T)$. It is a Levi subgroup of $G$. Fix a minimal Levi subgroup $M_{\mini}$ of $G$ included in $M_T$, a minimal parabolic subgroup $P_{\mini}$ having $M_{\mini}$ as a Levi component and $K$ a maximal compact subgroup of $G(F)$ which is special in the $p$-adic case. We use this compact subgroup to define the functions $H_Q$ for all $Q\in \mathcal{F}(M_{\mini})$. Let $\Delta_{\mini}$ be the set of roots of $A_{\mini}$ in $P_{\mini}$. Let $Y\in \mathcal{A}_{P_{\mini}}^+$ and define $Y_P$, for $P\in \mathcal{P}(M_{\mini})$, by $Y_P=w\cdot Y$ where $w$ is the only element in the Weyl group $W(G,M_{\mini})$ such that $wP_{\mini}w^{-1}=P$. Then $(Y_P)_{P\in\mathcal{P}(M_{\mini})}$ is a positive $(G,M_{\mini})$-orthogonal set. By the general constructions of Section \ref{section 1.9}, this determines a positive $(G,M_T)$-orthogonal set $(Y_{P_T})_{P_T\in \mathcal{P}(M_T)}$. For all $g\in G(F)$, we define another $(G,M_T)$-orthogonal set $\mathcal{Y}(g)=\big(\mathcal{Y}(g)_{P_T}\big)_{P_T\in \mathcal{P}(M_T)}$ by setting

$$\mathcal{Y}(g)_{P_T}=Y_{P_T}-H_{\overline{P}_T}(g)$$

\noindent for all $P_T\in \mathcal{P}(M_T)$ and where $\overline{P}_T$ denote the parabolic subgroup opposite to $P_T$ with respect to $M_T$. We will need the following

\vspace{3mm}

\begin{num}
\item\label{eq 10.10.15} There exists $c>0$ such that for all $g\in G(F)$ and $Y\in \mathcal{A}_{P_{\mini}}^+$ satisfying
$$\sigma(g)\leqslant c \inf_{\alpha\in \Delta_{\mini}} \alpha(Y)$$
the $(G,M_T)$-orthogonal set $\mathcal{Y}(g)$ is positive.
\end{num}

\vspace{3mm}

\noindent For $\mathcal{Y}=(\mathcal{Y}_{P_T})_{P_T\in \mathcal{P}(M_T)}$ a positive $(G,M_T)$-orthogonal set and $Q=LU_Q\in \mathcal{F}(M_T)$,we will denote by $\sigma_{M_T}^Q(., \mathcal{Y})$ and (resp.\ $\tau_Q$) the characteristic function in $\mathcal{A}_{M_T}$ of the sum of $\mathcal{A}_L$ and of the convex hull of the family $(\mathcal{Y}_{P_T})_{P_T\subset Q}$ (resp.\ the characteristic function of $\mathcal{A}_{M_T}^L+\mathcal{A}_Q^+$). Then we have (see \cite{LW} Lemme 1.8.4(3))

\begin{align}\label{eq 10.10.16}
\displaystyle \sum_{Q\in \mathcal{F}(M_T)} \sigma_{M_T}^Q(\zeta, \mathcal{Y}) \tau_Q(\zeta-\mathcal{Y}_Q)=1
\end{align}

\noindent for all $\zeta\in \mathcal{A}_{M_T}$.

\vspace{2mm}

\noindent For all $Y\in \mathcal{A}_{P_{\mini}}^+$, we define a function $\gls{tildevY}$ on $G(F)$ by

$$\displaystyle \tilde{v}(Y,g)=\int_{T(F)} \sigma_{M_T}^G(H_{M_T}(t), \mathcal{Y}(g)) dt$$

\noindent Set

$$\displaystyle \gls{JYTN-bf}=\int_{\mathfrak{t}(F)'[> N^{-b}]} D^G(X)^{1/2} \int_{T(F)\backslash G(F)} \mathbf{1}_{<\log(N)}(g^{-1}Xg) \widehat{f}(g^{-1}Xg) \tilde{v}(Y,g) dg dX$$

\noindent for all $N\geqslant 1$ and $Y\in \mathcal{A}_{P_{\mini}}^+$.

\begin{prop}\label{proposition 10.8.1}
There exist $c_1,c_2>0$ such that
$$\left\lvert J^{\Lie}_{N,T,N^{-b}}(f)-J_{Y,T,N^{-b}}(f)\right\rvert \ll N^{-1}$$
for all $N\geqslant 1$ and all $Y\in \mathcal{A}_{P_{\mini}}^+$ that satisfy the following two conditions

\begin{align}\label{eq 10.10.17}
\displaystyle c_1 \log(N)\leqslant \inf_{\alpha\in \Delta_{\mini}} \alpha(Y)
\end{align}

\begin{align}\label{eq 10.10.18}
\sup_{\alpha\in \Delta_{\mini}} \alpha(Y)\leqslant c_2 N
\end{align}
\end{prop}

\noindent\ul{Proof}: For all $N\geqslant 1$, let us denote by $\mathcal{A}_N$ the subset of $\left(\mathfrak{t}(F)'\cap\omega_T\right)\times T(F)\backslash G(F)$ consisting of pairs $(X,g)$ such that $\lvert Q(X)\rvert>N^{-b}$ and $\sigma_{\mathfrak{g}}(g^{-1}Xg)<\log(N)$. Then, we have

$$\displaystyle J^{\Lie}_{N,T,N^{-b}}(f)=\int_{\mathcal{A}_N} D^G(X)^{1/2}\widehat{f}(g^{-1}Xg)\kappa_{N,X}(g) dXdg$$

$$\displaystyle J_{Y,T,N^{-b}}(f)=\int_{\mathcal{A}_N} D^G(X)^{1/2}\widehat{f}(g^{-1}Xg)\tilde{v}(Y,g) dXdg$$

\noindent for all $N\geqslant 1$ and all $Y\in \mathcal{A}_{P_{\mini}}^+$. Let $c_1$, $c_2$ be positive real numbers. We will prove that the inequality of the proposition is valid for all $N\geqslant 1$ and all $Y\in \mathcal{A}_{P_{\mini}}^+$ that satisfies the inequalities \ref{eq 10.10.17} and \ref{eq 10.10.18} as long as $c_1$ is large enough and $c_2$ is small enough. We note the following

\vspace{3mm}

\begin{num}
\item\label{eq 10.10.19} We have an inequality $\sigma_{T\backslash G}(g)\ll \log(N)$ for all $N\geqslant 1$ and all $(X,g)\in\mathcal{A}_N$.
\end{num}

\vspace{3mm}

\noindent Indeed, this follows from \ref{eq 1.2.2} and the fact that $d^G$ divides $Q$. In particular, by \ref{eq 10.10.15}, if $c_1$ is sufficiently large, the $(G,M_T)$-family $\mathcal{Y}(g)$ is positive orthogonal for all $N\geqslant 1$, all $(X,g)\in \mathcal{A}_N$ and all $Y\in\mathcal{A}_{P_{\mini}}^+$ that satisfies \ref{eq 10.10.17}. We will henceforth assume that $c_1$ is at least that sufficiently large. Hence, for all $Q\in \mathcal{F}(M_T)$, we can set

$$\displaystyle \kappa^Y_{N,X,Q}(g):=\int_{T(F)} \kappa_N(\gamma_X^{-1}tg) \sigma_{M_T}^Q(H_{M_T}(t),\mathcal{Y}(g)) \tau_Q(H_{M_T}(t)-\mathcal{Y}(g)_Q) dt$$

\noindent for all $N\geqslant 1$, all $(X,g)\in\mathcal{A}_N$ and all $Y\in\mathcal{A}_{P_{\mini}}^+$ that satisfies \ref{eq 10.10.17}. By \ref{eq 10.10.16}, we have the decomposition

$$\displaystyle \kappa_{N,X}(g)=\sum_{Q\in \mathcal{F}(M_T)} \kappa^Y_{N,X,Q}(g)$$

\noindent for all $N,X,g$ and $Y$ as before. Obviously the functions $\kappa^Y_{N,X,Q}$ are left invariant by $T(F)$ and so we have accordingly a decomposition

$$\displaystyle J^{\Lie}_{N,T,N^{-b}}(f)=\sum_{Q\in \mathcal{F}(M_T)} J^{Q,Y}_{N,T,N^{-b}}(f)$$

\noindent for all $N\geqslant 1$ and all $Y\in \mathcal{A}_{P_{\mini}}^+$ satisfying \ref{eq 10.10.17}, where we have set

$$\displaystyle  J^{Q,Y}_{N,T,N^{-b}}(f)=\int_{\mathcal{A}_N} D^G(X)^{1/2} \widehat{f}(g^{-1}Xg) \kappa^Y_{N,X,Q}(g) dg dX$$

\noindent The proposition will now follows from the two following facts

\vspace{3mm}

\begin{num}
\item\label{eq 10.10.20} If $c_2$ is sufficiently small, there exists $N_0\geqslant 1$ such that
$$J^{G,Y}_{N,T,N^{-b}}(f)=J_{Y,T,N^{-b}}(f)$$
for all $N\geqslant N_0$ and all $Y\in \mathcal{A}_{P_{\mini}}^+$ satisfying \ref{eq 10.10.18}.
\end{num}

\vspace{3mm}

\begin{num}
\item\label{eq 10.10.21}Let $Q\in \mathcal{F}(M_T)$, $Q\neq G$. If $c_1$ is large enough, then we have an inequality
$$\left\lvert J^{Q,Y}_{N,T,N^{-b}}(f)\right\rvert \ll N^{-1}$$
for all $N\geqslant 1$ and all $Y\in \mathcal{A}_{P_{\mini}}^+$ satisfying \ref{eq 10.10.17}.
\end{num}

\vspace{3mm}

\noindent First we prove \ref{eq 10.10.20}. Obviously, we only need to prove that for $c_2$ small enough, $N$ large enough and $Y\in \mathcal{A}_{P_{\mini}}^+$ satisfying \ref{eq 10.10.18}, we have
$$\kappa^Y_{N,X,G}(g)=\tilde{v}(Y,g)$$
for all $(X,g)\in \mathcal{A}_N$. Expanding the definitions, it is certainly enough to prove that

\begin{align}\label{eq 10.10.22}
\displaystyle \sigma^G_{M_T}(H_{M_T}(t), \mathcal{Y}(g))=1\Rightarrow \kappa_N(\gamma_X^{-1}tg)=1
\end{align}

\noindent for all $(X,g)\in \mathcal{A}_N$ and all $t\in T(F)$. We have an inequality
$$\displaystyle \sigma(tg)\ll \sup_{\alpha\in \Delta_{\mini}} \alpha(Y)+\sigma_{T\backslash G}(g)$$
for all $Y\in \mathcal{A}_{P_{\mini}}^+$, all $g\in G(F)$ and all $t\in T(F)$ satisfying $\sigma_{M_T}^G(H_{M_T}(t), \mathcal{Y}(g))=1$. By \ref{eq 10.9.8}, we have $\sigma(\gamma_X)\ll \log(N)$ for all $N\geqslant 2$ and all $X\in \mathfrak{t}(F)'[>N^{-b}]\cap \omega_T$. Combining these two facts with \ref{eq 10.10.19}, we get

\begin{align}\label{eq 10.10.23}
\displaystyle \sigma_{M_T}^G(H_{M_T}(t), \mathcal{Y}(g))=1\Rightarrow \sigma(\gamma_X^{-1}tg)\ll \sup_{\alpha\in \Delta_{\mini}} \alpha(Y)+ \log(N)
\end{align}

\noindent for all $N\geqslant 2$, all $(X,g)\in\mathcal{A}_N$, all $t\in T(F)$ and all $Y\in \mathcal{A}_{P_{\mini}}^+$. Moreover, by the property \ref{eq 10.9.1} of $\kappa_N$, there exists $C_0>0$ such that for all $N\geqslant 1$ and all $\gamma\in G(F)$ we have $\sigma(\gamma)\leqslant C_0N\Rightarrow \kappa_N(\gamma)=1$. Hence, \ref{eq 10.10.22} follows from \ref{eq 10.10.23} when $c_2$ is small enough and $N$ large enough and this ends the proof of \ref{eq 10.10.20}.

\vspace{2mm}

\noindent We now move on to the proof of \ref{eq 10.10.21}. Fix a proper parabolic subgroup $Q=LU_Q\in \mathcal{F}(M_T)$ and denote by $\overline{Q}=LU_{\overline{Q}}$ the opposite parabolic subgroup with respect to $L$ (the only Levi component of $Q$ containing $M_T$). We have the Iwasawa decomposition $G(F)=L(F)U_{\overline{Q}}(F)K$ and for suitable choices of Haar measures we have $dg=dldudk$. For all $N\geqslant 1$, let $\mathcal{B}_N$ be the set of quadruple $(X,l,u,k)\in \mathfrak{t}(F)'\times \left(T(F)\backslash L(F)\right)\times U_{\overline{Q}}(F)\times K$ such that $(X,luk)\in \mathcal{A}_N$. Then we have
$$\displaystyle J^{Q,Y}_{N,T,N^{-b}}(f)=\int_{\mathcal{B}_N} D^G(X)^{1/2} \widehat{f}(k^{-1}u^{-1}l^{-1}Xluk)\kappa^Y_{N,X,Q}(luk) dldudkdX$$
for all $N\geqslant 1$ and all $Y\in \mathcal{A}_{P_{\mini}}^+$. We claim the following:

\vspace{3mm}

\begin{num}
\item\label{eq 10.10.24} If $c_1$ is large enough, we have
$$\lvert\kappa^Y_{N,X,Q}(luk)-\kappa^Y_{N,X,Q}(lk)\rvert\ll N^{-1}$$
for all $N\geqslant 1$, all $(X,l,u,k)\in \mathcal{B}_N$ and all $Y\in \mathcal{A}_{P_{\mini}}^+$ satisfying \ref{eq 10.10.17}.
\end{num}

\vspace{3mm}

\noindent We will postpone the proof of \ref{eq 10.10.24} and show how to deduce \ref{eq 10.10.21} from it. Assume that $c_1$ is large enough so that \ref{eq 10.10.24} holds. Then we have

\[\begin{aligned}
 & \displaystyle \left\lvert J^{Q,Y}_{N,T,N^{-b}}(f) -\int_{\mathcal{B}_N} D^G(X)^{1/2} \widehat{f}(k^{-1}u^{-1}l^{-1}Xluk) \kappa^Y_{N,X,Q}(lk)dldudkdX\right\rvert  \\
  & \ll N^{-1} \int_{\mathcal{B}_N} D^G(X)^{1/2} \left\lvert\widehat{f}(k^{-1}u^{-1}l^{-1}Xluk)\right\rvert dXdldudk \\
  & =  N^{-1} \int_{\mathcal{A}_N} D^G(X)^{1/2} \left\lvert \widehat{f}(g^{-1}Xg)\right\rvert dXdg \\
 & \ll  N^{-1} \int_{\mathfrak{t}(F)} J_G(X, \lvert \widehat{f}\rvert)dX
\end{aligned}\]

\noindent for all $N\geqslant 1$ and $Y\in \mathcal{A}_{P_{\mini}}^+$ satisfying \ref{eq 10.10.17}. The last integral above is convergent. Consequently, in order to prove \ref{eq 10.10.21}, we only need to prove an inequality

\begin{align}\label{eq 10.10.25}
\displaystyle \left\lvert\int_{\mathcal{B}_N} D^G(X)^{1/2} \widehat{f}(k^{-1}u^{-1}l^{-1}Xluk) \kappa^Y_{N,X,Q}(lk)dldudkdX\right\rvert\ll N^{-1}
\end{align}

\noindent for all $N\geqslant 1$ and all $Y\in \mathcal{A}_{P_{\mini}}^+$. Obviously, we have $\sigma_{T\backslash G}(lk)\ll \sigma_{T\backslash G}(luk)$ for all $l\in L(F)$, all $u\in U_{\overline{Q}}(F)$ and all $k\in K$. By \ref{eq 10.10.19}, it follows that there exists $c>0$ such that $\sigma_{T\backslash G}(lk)< c\log(N)$ for all $N\geqslant 2$ and all $(X,l,u,k)\in \mathcal{B}_N$. Let us denote, for all $N\geqslant 2$, by $\mathcal{C}_N$ the set of triples $(X,l,k)\in \left(\mathfrak{t}(F)'\cap\omega_T\right)\times \big(T(F)\backslash L(F)\big)\times K$ such that $\lvert Q(X)\rvert> N^{-b}$ and $\sigma_{T\backslash G}(lk)<c \log(N)$. By what we just said, we have $\mathcal{B}_N\subset \mathcal{C}_N\times U_{\overline{Q}}(F)$. Since $\widehat{f}$ is strongly cuspidal, we have
$$\displaystyle \int_{\mathcal{C}_N\times U_{\overline{Q}}(F)} D^G(X)^{1/2} \widehat{f}(k^{-1}u^{-1}l^{-1}Xluk) \kappa^Y_{N,X,Q}(lk) dudldkdX=0$$
Let us set $\mathcal{D}_N=\left(\mathcal{C}_N\times U_{\overline{Q}}(F)\right)\backslash \mathcal{B}_N$ for all $N\geqslant 2$. Note that by definition of $\mathcal{B}_N$, we have

\begin{align}\label{eq 10.10.26}
\displaystyle \sigma_{\mathfrak{g}}(k^{-1}u^{-1}l^{-1}Xluk)\geqslant \log(N)
\end{align}

\noindent for all $N\geqslant 2$ and all $(X,l,u,k)\in\mathcal{D}_N$. From the vanishing of the above integral, we deduce

\begin{align}\label{eq 10.10.27}
\displaystyle & \int_{\mathcal{B}_N} D^G(X)^{1/2} \widehat{f}(k^{-1}u^{-1}l^{-1}Xluk) \kappa^Y_{N,X,Q}(lk)dldudkdX \\
\nonumber & = -\int_{\mathcal{D}_N} D^G(X)^{1/2} \widehat{f}(k^{-1}u^{-1}l^{-1}Xluk) \kappa^Y_{N,X,Q}(lk)dldudkdX
\end{align}

\noindent for all $N\geqslant 2$ and all $Y\in \mathcal{A}_{P_{\mini}}^+$. Obviously, we have $\kappa^Y_{N,X,Q}\leqslant \kappa_{N,X}$. Using Lemma \ref{lemma 10.8.1}(i), it follows that there exists $k>0$ such that
$$\kappa^Y_{N,X,Q}(lk)\ll N^k$$
for all $N\geqslant 2$, all $(X,l,k)\in \mathcal{C}_N$ and all $Y\in \mathcal{A}_{P_{\mini}}^+$. Hence, we see that the integral \ref{eq 10.10.27}, is essentially bounded (in absolute value) by
$$\displaystyle N^{k}\int_{\mathcal{D}_N} D^G(X)^{1/2} \lvert\widehat{f}(k^{-1}u^{-1}l^{-1}Xluk)\rvert dldudkdX$$
for all $N\geqslant 2$ and all $Y\in \mathcal{A}_{P_{\mini}}^+$. But, since $\widehat{f}$ is a Schwartz function, we easily deduce from \ref{eq 10.10.26} that the last integral above is essentially bounded by $N^{-1-k}$ for all $N\geqslant 2$. This proves \ref{eq 10.10.25} and ends the proof of \ref{eq 10.10.21}.

\vspace{2mm}

\noindent It remains to prove the crucial point \ref{eq 10.10.24}. By \ref{eq 10.10.19}, there exists $c>0$ such that for all $N\geqslant 2$ and all $(X,l,u,k)\in\mathcal{B}_N$ up to translating $l$ by an element of $T(F)$, we have $\sigma(lul^{-1})<c\log(N)$ and $\sigma(lk)<c\log(N)$. Since $\kappa^Y_{N,X,Q}$ is left invariant by $T(F)$, it suffices to prove the following

\vspace{3mm}

\begin{num}
\item\label{eq 10.10.28} If $c_1$ is sufficiently large, we have 
$$\displaystyle \left\lvert \kappa^Y_{N,X,Q}(ug)-\kappa^Y_{N,X,Q}(g)\right\rvert\ll N^{-1}$$
for all $N\geqslant 2$, all $X\in \mathfrak{t}(F)'[>N^{-b}]\cap \omega_T$, all $u\in U_{\overline{Q}}(F)$ and all $g\in G(F)$ satisfying $\sigma(u)<c\log(N)$, $\sigma(g)<c\log(N)$ and all $Y\in \mathcal{A}_{P_{\mini}}^+$ satisfying \ref{eq 10.10.17}.
\end{num}

\vspace{3mm}

\noindent Let $N$, $X$, $u$, $g$ and $Y$ as above. We will prove that the inequality \ref{eq 10.10.28} holds provided $c_1$ is large enough. First, if $c_1$ is large enough, we see by \ref{eq 10.10.15} that the $(G,M_T)$-families $\mathcal{Y}(ug)$ and $\mathcal{Y}(g)$ are positive orthogonal. For $g'\in G(F)$, the function $\sigma_{M_T}^Q(.,\mathcal{Y}(g'))\tau_Q(.-\mathcal{Y}(g')_Q)$ only depends on $\mathcal{Y}(g')_{P_T}$ for $P_T\in \mathcal{P}(M_T)$ with $P_T\subset Q$ and these terms are invariant by left translation of $g'$ by $U_{\overline{Q}}(F)$. Thus, we have
$$\displaystyle \sigma_{M_T}^Q(.,\mathcal{Y}(ug))\tau_Q(.-\mathcal{Y}(ug)_Q)=\sigma_{M_T}^Q(.,\mathcal{Y}(g))\tau_Q(.-\mathcal{Y}(g)_Q)$$
Remembering the definition of $\kappa^Y_{N,X,Q}$, we deduce that

\[\begin{aligned}
\displaystyle \left\lvert \kappa^Y_{N,X,Q}(ug)-\kappa^Y_{N,X,Q}(g)\right\rvert\leqslant \int_{T(F)} \left\lvert \kappa_N(\gamma_X^{-1}tug)-\kappa(\gamma_X^{-1}tg)\right\rvert & \sigma_{M_T}^Q(H_{M_T}(t),\mathcal{Y}(g)) \\
 & \tau_Q(H_{M_T}(t)-\mathcal{Y}(g)_Q)dt
\end{aligned}\]

\noindent By the property \ref{eq 10.9.1} of $\kappa_N$, there exists $C_1>0$ such that for all $t\in T(F)$ we have
$$\displaystyle \left\lvert \kappa_N(\gamma_X^{-1}tug)-\kappa(\gamma_X^{-1}tg)\right\rvert\neq 0\Rightarrow \sigma(t)\leqslant C_1\left(N+\sigma(\gamma_X)+\sigma(u)+\sigma(g)\right)$$
By the hypothesis made on $g$, $u$ and $X$ and by \ref{eq 10.9.8}, this last condition implies $\sigma(t)\leqslant C_2 N$ for some bigger constant $C_2>C_1$. Moreover, there exists $C_3>0$ and $k>0$ such that the volume of the subset $\{t\in T(F);\;\sigma(t)\leqslant C_2N\}$ is bounded by $C_3N^{k}$. Hence, we will be done if we can establish the following

\vspace{3mm}

\begin{num}
\item\label{eq 10.10.29} Provided $c_1$ is sufficiently large, for all $t\in T(F)$ satisfying 
$$\displaystyle \sigma_{M_T}^Q(H_{M_T}(t),\mathcal{Y}(g))\tau_Q(H_{M_T}(t)-\mathcal{Y}(g)_Q)=1$$
we have
$$\displaystyle \left\lvert \kappa_N(\gamma_X^{-1}tug)-\kappa_N(\gamma_X^{-1}tg)\right\rvert\leqslant N^{-k-1}$$
\end{num}

\vspace{3mm}

\noindent Fix $t\in T(F)$ such that
$$\displaystyle \sigma_{M_T}^Q(H_{M_T}(t),\mathcal{Y}(g))\tau_Q(H_{M_T}(t)-\mathcal{Y}(g)_Q)=1$$
Let $\Sigma_Q^+$ be the set of roots of $A_T$ in $U_Q$ and consider it as a subset of $\mathcal{A}_{M_T}^*=X^*(A_T)\otimes \mathbb{R}$. There exist positive constants $C_4$ and $C_5$ such that
$$\displaystyle \langle \beta,H_{M_T}(t)\rangle\geqslant C_4 \inf_{\alpha\in \Delta_{\mini}} \alpha(Y')-C_5 \sigma_G(g'),\;\;\;\mbox{for all } \beta\in \Sigma_Q^+$$
for all $Y'\in \mathcal{A}_{P_{\mini}}^+$, all $g'\in G(F)$ and all $t\in T(F)$ such that
$$\displaystyle \sigma_{M_T}^Q(H_{M_T}(t),\mathcal{Y}'(g'))\tau_Q(H_{M_T}(t)-\mathcal{Y}'(g')_Q)=1$$
Hence, by the assumptions on $g$ and $Y$, we have

\begin{align}\label{eq 10.10.30}
\displaystyle \langle \beta,H_{M_T}(t)\rangle\geqslant (C_4c_1-C_5c)\log(N),\;\;\;\mbox{for all } \beta\in \Sigma_Q^+
\end{align}

\noindent Let $e$ be a positive real number that we will assume sufficiently large in what follows. By \ref{eq 10.10.30} and the assumption on $u$, if $c_1$ is sufficiently large, we have
$$\displaystyle tut^{-1}\in \exp\left(B\left(0,N^{-4e}\right)\right)$$
Hence, by \ref{eq 10.9.8} and the assumption on $X$, if $e$ is sufficiently large, we have
$$\displaystyle \gamma_X^{-1}tut^{-1}\gamma_X\in \exp\left(B\left(0,N^{-3e}\right)\right)$$
Let $P_t\in \mathcal{P}(M_T)$ be a parabolic subgroup such that $H_{M_T}(t)\in \overline{\mathcal{A}_{P_t}^+}$ (where the bar denotes the closure). Recall that $\gamma_X^{-1}X\gamma_X=X_\Sigma\in \Sigma'(F)$. By Corollary \ref{corollary 10.6.1} and \ref{eq 10.9.8} again, if we choose $e$ sufficiently large, we will have
$$\displaystyle \gamma_X^{-1}tut^{-1}\gamma_X \in  H(F) \gamma_X^{-1}\exp\left( B\left(0,N^{-2e}\right)\cap \mathfrak{p}_t(F)\right)\gamma_X$$
where $\mathfrak{p}_t=\Lie(P_t)$. Hence, we may write $\gamma_X^{-1}tut^{-1}\gamma_X=h\gamma_X^{-1}e^{X_P}\gamma_X$ with $h\in H(F)$ and $X_P\in B\left(0,N^{-2e}\right)\cap \mathfrak{p}_t(F)$. By left invariance of $\kappa_N$ by $H(F)$, we will have

\begin{align}\label{eq 10.10.31}
\displaystyle \kappa_N(\gamma_X^{-1}tug)=\kappa_N(\gamma_X^{-1}e^{X_P}tg)
\end{align}

\noindent Since we choose $P_t\in \mathcal{P}(M_T)$ so that $H_{M_T}(t)\in \overline{\mathcal{A}_{P_t}^+}$, there exists $c_5>0$ independent of $t$ such that $t^{-1}\left(B(0,1)\cap \mathfrak{p}_t(F)\right)t\subseteq B(0,c_5)\cap \mathfrak{p}_t(F)$. Hence, we get
$$\displaystyle t^{-1}X_Pt\in B\left(0,c_5N^{-2e}\right)\cap \mathfrak{p}_t(F)$$
It follows by the assumption on $g$ that if $e$ is large enough, we have
$$\displaystyle g^{-1}t^{-1}X_Ptg\in B\left(0,N^{-e}\right)$$
By the conditions \ref{eq 10.9.2} and \ref{eq 10.9.3} that we imposed on $\kappa_N$, for $e$ sufficiently large we have

$$\displaystyle \mbox{(In the non-Archimedean case)}\;\;\; \kappa_N(\gamma_X^{-1}e^{X_P}tg)=\kappa_N(\gamma_X^{-1}tg)$$

$$\displaystyle \mbox{(In the Archimedean case)}\;\;\; \left\lvert \kappa_N(\gamma_X^{-1}tg)- \kappa_N(\gamma_X^{-1}e^{X_P}tg)\right\rvert\leqslant N^{-1-k}$$

\noindent Combining this with \ref{eq 10.10.31}, we get, if $c_1$ is large enough,

$$\displaystyle \mbox{(In the non-Archimedean case)}\;\;\; \kappa_N(\gamma_X^{-1}tug)=\kappa_N(\gamma_X^{-1}tg)$$

$$\displaystyle \mbox{(In the Archimedean case)}\;\;\; \left\lvert \kappa_N(\gamma_X^{-1}tg)- \kappa_N(\gamma_X^{-1}tug)\right\rvert\leqslant N^{-1-k}$$

\noindent This proves \ref{eq 10.10.29} from which the claim \ref{eq 10.10.24} follows and this ends the proof of the proposition. $\blacksquare$

\subsection{End of the proof of Theorem \ref{theorem 10.8.1}}\label{section 10.11}

\noindent We are now in position to finish the proof of Theorem \ref{theorem 10.8.1}. Let us set
$$\displaystyle J_{Y,T}(f)=\int_{\mathfrak{t}(F)'} D^G(X)^{1/2}\int_{T(F)\backslash G(F)} \widehat{f}(g^{-1}Xg) \widetilde{v}(Y,g)dgdX$$
for all $Y\in \mathcal{A}_{P_{\mini}}^+$. Obviously, we can find $k>0$ such that
$$\displaystyle \widetilde{v}(Y,g)\ll (1+\lvert Y\rvert)^k \sigma_{T\backslash G}(g)^k$$
for all $Y\in \mathcal{A}_{\mini}^+$ and all $g\in G(F)$. Using \ref{eq 1.2.2}, \ref{eq 1.2.3} and \ref{eq 1.2.4}, it follows that there exists $\epsilon>0$ such that we have
$$\displaystyle \left\lvert J_{Y,T}(f)-J_{Y,T,N^{-b}}(f)\right\rvert \ll N^{k-b\epsilon}$$
for all $N\geqslant 2$ and all $Y\in \mathcal{A}_{\mini}^+$ satisfying the inequality \ref{eq 10.10.18} of Proposition \ref{proposition 10.8.1}. Hence, combining this inequality with Proposition \ref{proposition 10.8.1}, we see that if we choose $b$ large enough we have

\begin{align}\label{eq 10.11.32}
\displaystyle \left\lvert J_{Y,T}(f)-J^{\Lie}_{N,T}(f)\right\rvert\ll N^{-1}
\end{align}

\noindent for all $N\geqslant 2$ and all $Y\in\mathcal{A}_{\mini}^+$ that satisfies the inequalities \ref{eq 10.10.17} and \ref{eq 10.10.18} of Proposition \ref{proposition 10.8.1}. Arthur has computed the functions $Y\in \mathcal{A}_{\mini}^+\mapsto \widetilde{v}(Y,g)$ (cf.\ \cite{A1}, p.46). More precisely, the result is the following: for every $g\in G(F)$ the function $Y\in \mathcal{A}_{\mini}^+\mapsto \widetilde{v}(Y,g)$ for $Y$ in a certain lattice $\mathcal{R}$ is a sum of functions of the form $q_\zeta(Y,g)e^{\zeta(Y)}$ where $q_\zeta(.,g)$ is a polynomial in $Y$ and $\zeta\in \Hom(\mathcal{R},2i\pi\mathbb{R}/2i\pi\mathbb{Z})$. Such functions are linearly independent, and it follows from \ref{eq 10.11.32} that we have

\begin{align}\label{eq 10.11.33}
\displaystyle \lim\limits_{N\to \infty} J^{\Lie}_{N,T}(f)=\int_{\mathfrak{t}(F)'}D^G(X)^{1/2}\int_{T(F)\backslash G(F)} \widehat{f}(g^{-1}Xg)q_0(0,g) dgdX
\end{align}

\noindent Moreover, by \cite{A1} (6.6), we have
$$\displaystyle q_0(0,g)=(-1)^{a_{M_T}}\nu(T)^{-1}\sum_{Q\in \mathcal{F}(M_T)} c'_Q v_{M_T}^Q(g)$$
where the $c'_Q$ are certain constant with $c'_G=1$ and the functions $v_{M_T}^Q(g)$ are the one introduced in Section \ref{section 1.10} (recall that we fixed a maximal compact subgroup $K$ that is special in the $p$-adic case). We remind the reader that the factor $\nu(T)$ is the quotient between the Haar measure we fixed on $T(F)$ (i.e. the autodual one) and another natural Haar measure on $T(F)$ (cf.\ Section \ref{section 1.6}). It is present in the formula above because we are using different normalizations of measures than the one used by Arthur in \cite{A1}. By definition of the weighted orbital integrals $J^Q_{M_T}(X,.)$, \ref{eq 10.11.33} becomes
$$\displaystyle \lim\limits_{N\to \infty} J^{\Lie}_{N,T}(f)=(-1)^{a_{M_T}}\nu(T)^{-1}\sum_{Q\in \mathcal{F}(M_T)} c'_Q\int_{\mathfrak{t}(F)'}D^G(X)^{1/2}J_{M_T}^Q(X,\widehat{f})dX$$
Since $\widehat{f}$ is strongly cuspidal, by Lemma \ref{lemma 5.2.1}(i) we have $J_{M_T}^Q(X,\widehat{f})=0$ for all $X\in\mathfrak{t}_{\reg}(F)$ and all $Q\in\mathcal{F}(M_T)$ such that $Q\neq G$. So finally we obtain

\[\begin{aligned}
\displaystyle \lim\limits_{N\to \infty} J^{\Lie}_{N,T}(f) & =(-1)^{a_{M_T}}\nu(T)^{-1}\int_{\mathfrak{t}(F)'}D^G(X)^{1/2}J_{M_T}^G(X,\widehat{f})dX \\
 & =\int_{\mathfrak{t}(F)'}D^G(X)^{1/2} \theta_{\widehat{f}}(X) dX
\end{aligned}\]

\noindent Summing this last equality over $T\in\mathcal{T}(G)$, we obtain, by \ref{eq 10.9.7} and \ref{eq 10.9.4},

\[\begin{aligned}
\displaystyle J^{\Lie}(f) & =\sum_{T\in\mathcal{T}(G)}\left\lvert W(G,T)\right\rvert^{-1} \int_{\mathfrak{t}(F)'}D^G(X)^{1/2} \theta_{\widehat{f}}(X) dX \\
 & =\int_{\Gamma(\Sigma)}D^G(X)^{1/2} \theta_{\widehat{f}}(X) dX
\end{aligned}\]

\noindent But, by Proposition \ref{proposition 5.6.1}(i) we have $\theta_{\widehat{f}}=\widehat{\theta}_f$. This ends the proof of Theorem \ref{theorem 10.8.1}. $\blacksquare$

\section{Geometric expansions and a formula for the multiplicity}\label{section 11}

This is the last chapter on the proof of the local simple trace formulas for GGP triples. More precisely, we will establish geometric expansions for the distributions $J$ (Theorem \ref{theorem 11.2.1}) and $J^{\Lie}$ (Theorem \ref{theorem 11.2.3}) as well as a certain integral formula for the multiplicity $m(\pi)$ (Theorem \ref{theorem 11.2.2}) when the representation $\pi$ is tempered. These three results will be proved together in a common inductive proof which is scattered over Sections \ref{section 11.3} to \ref{section 11.7}. On the other hand, Sections \ref{section 11.0}, \ref{section 11.1} and \ref{section 11.1bis} contains preliminary material for the statement of the main theorems (and their proof). In more details, in Section \ref{section 11.0} we introduce some spaces of semi-simple conjugacy classes which are then used in Section \ref{section 11.1} to define certain linear forms $m_{\geom}(.)$ and $m_{\geom}^{\Lie}(.)$ on spaces of quasi-characters for the group $G(F)$ and its Lie algebra respectively. These linear forms are the main ingredients in the formulation of the three theorems of Section \ref{section 11.2}. Finally, in Section \ref{section 11.1bis} we record a result pertaining to the compatibility of the linear form $m_{\geom}(.)$ with parabolic induction.

\subsection{Some spaces of conjugacy classes}\label{section 11.0}

In this section, we introduce certain spaces of semi-simple conjugacy classes in $G(F)$, $G_x(F)$ (for $x\in H_{\ssi}(F)$) and $\mathfrak{g}(F)$ to be denoted $\Gamma(G,H)$, $\Gamma(G_x,H_x)$  and $\Gamma^{\Lie}(G,H)$ respectively. These will be needed in the next section to define certain continuous linear forms on the spaces of quasi-characters $QC(G(F))$ and $SQC(\mathfrak{g}(F))$ which are in turn the main ingredients entering in the statement of the geometric expansions of $J(f)$ and $J^{\Lie}(f)$ as well as the formula for the multiplicity $m(\pi)$ (see Section \ref{section 11.2}). 

\vspace{2mm}

\noindent Let $x\in H_{\ssi}(F)$. We first give an explicit description of the triple $\gls{GxHxxix}$ where $\xi_x=\xi_{\mid H_x(F)}$. Up to conjugation, we may assume that $x\in U(W)_{\ssi}(F)$. Denote by $\gls{W'_x}$ and $\gls{V'_x}$ the kernel of $1-x$ in $W$ and $V$ respectively and by $\gls{W''_x}$ the image of $1-x$. We then have the orthogonal decompositions $W=W'_x\oplus^\perp W''_x$ and $V=V'_x\oplus^\perp W''_x$. Set $\gls{H'_x}=U(W'_x)\ltimes N_x$ (where $N_x$ is the centralizer of $x$ in $N$), $\gls{G'_x}=U(W'_x)\times U(V'_x)$, $\gls{H''_x}=U(W''_x)_x$ and $\gls{G''_x}=U(W''_x)_x\times U(W''_x)_x$. We have natural decompositions 
$$\displaystyle U(V)_x=U(V'_x)\times U(W''_x)_x,\;\; U(W)_x=U(W'_x)\times U(W''_x)_x\mbox{ and } H_x=U(W)_x\ltimes N_x$$
Moreover, we easily check that $U(W''_x)_x$ commutes with $N_x$. Hence, we also have natural decompositions

\begin{align}\label{eq 11.1.1}
\displaystyle G_x=G'_x\times G''_x\mbox{ and } H_x=H'_x\times H''_x
\end{align}

\noindent the inclusions $H_x\subset G_x$ being the product of the two inclusions $H'_x\subset G'_x$ and $H''_x\subset G''_x$. It is clear that $\xi_x$ is trivial on $H''_x$, so that we get a decomposition
$$\displaystyle (G_x,H_x,\xi_x)=(G'_x,H'_x,\xi'_x)\times (G''_x,H''_x,1)$$
where $\xi'_x=\xi_{\mid H'_x}$ and the product of triples is obviously defined. Note that the triple $(G'_x,H'_x,\xi'_x)$ coincides with the GGP triple associated to the admissible pair $(V'_x,W'_x)$. The second triple $(G''_x,H''_x,1)$ is also of a particular shape: the group $G''_x$ is the product of two copies of $H''_x$ and the inclusion $H''_x\subset G''_x$ is the diagonal one. We shall call such a triple an {\em Arthur triple}. Finally, note that although we have assumed $x\in U(W)_{\ssi}(F)$, there is a decomposition similar to \ref{eq 11.1.1} for any $x\in H_{\ssi}(F)$ (just conjugated $x$ inside $H(F)$ to an element in $U(W)_{\ssi}(F)$) and that if $x,y\in H_{\ssi}(F)$ are $H(F)$-conjugate there are natural isomorphisms of triples
$$(G'_x,H'_x,\xi'_x)\simeq (G'_y,H'_y,\xi'_y) \mbox{ and } (G''_x,H''_x,1)\simeq (G''_y,H''_y,1)$$
well-defined up to inner automorphisms (by $H'_x(F)$ and $H''_x(F)$ respectively).

\vspace{2mm}

\noindent Let $x\in H_{\ssi}(F)$. As in Section \ref{section 1.8}, we denote by $\Gamma(H)$, $\Gamma(H_x)$, $\Gamma(G)$ and $\Gamma(G_x)$ the sets of semi-simple conjugacy classes in $H(F)$, $H_x(F)$, $G(F)$ and $G_x(F)$ respectively and we equip them with topologies. Then, we have the following

\vspace{3mm}

\begin{num}
\item\label{eq 11.1.2} The natural maps $\Gamma(H_x)\to \Gamma(G_x)$ and $\Gamma(H)\to \Gamma(G)$ are closed embeddings. Moreover, if $\Omega_x\subseteq G_x(F)$ is a sufficiently small $G$-good open neighborhood of $x$ (see Section \ref{section 3.2} for this notion), the following diagram (where we identify $\Omega_x$ with its image in $\Gamma(G_x)$) is Cartesian
$$\xymatrix{
\Gamma(H_x)\cap \Omega_x \ar[d] \ar[r] & \Omega_x\ar[d]\\
\Gamma(H) \ar[r] & \Gamma(G)}.$$
\end{num}

\vspace{3mm}

\noindent Indeed, we see easily using the above descriptions of both $H_x$ and $G_x$ that the two maps $\Gamma(H_x)\to \Gamma(G_x)$ and $\Gamma(H)\to \Gamma(G)$ are injective. Since these maps are continuous and proper (see \S \ref{section 1.8}) and $\Gamma(H_x)$, $\Gamma(G_x)$, $\Gamma(H)$, $\Gamma(G)$ are all Hausdorff and locally compact, it follows that $\Gamma(H_x)\to \Gamma(G_x)$ and $\Gamma(H)\to \Gamma(G)$ are closed embeddings. Let $\Omega_x\subseteq G_x(F)$ be a $G$-good open neighborhood of $x$. We show now that the above diagram is Cartesian provided $\Omega_x$ is sufficiently small. This amounts to proving that if $y\in \Omega_{x,ss}$ is $G(F)$-conjugate to an element of $H(F)$ then $y$ is $G_x(F)$-conjugate to an element of $H_x(F)$. Let $y$ be such an element and let us fix $\Omega_x'\subseteq G_x(F)$ another $G$-good open neighborhood of $x$. In what follows, we will assume (as we may) that $\Omega_x\subseteq \Omega'_x$. Since $G_x=Z_G(x)$ (because $G_{\der}$ is simply-connected), by definition of a $G$-good open subset it suffices to show that if $\Omega_x$ is sufficiently small then $y$ is $G(F)$-conjugate to an element in $H_x(F)\cap \Omega'_x$. We easily check that $H_x(F)\cap\Omega'_x$ is a $H$-good open neighborhood of $x$. Hence, the map $\Gamma(H_x)\cap \Omega'_x\to \Gamma(H)$ is injective and has open image. Similarly, the map $\Gamma(G_x)\cap\Omega_x\to \Gamma(G)$ is injective and has open image. Moreover, as $\Omega_x$ runs through the $G$-good open neighborhoods of $x$ the subsets $\Gamma(G_x)\cap \Omega_x$ form a basis of open neighborhoods of $x$ in $\Gamma(G)$. Hence, since $\Gamma(H)\to \Gamma(G)$ is a closed embedding, the subsets $\Gamma(H)\cap\left(\Gamma(G_x)\cap \Omega_x\right)$, as $\Omega_x$ runs through the $G$-good open neighborhoods of $x$, form a basis of open neighborhoods of $x$ in $\Gamma(H)$. It follows that for $\Omega_x$ sufficiently small we have $\Gamma(H)\cap\left(\Gamma(G_x)\cap \Omega_x\right)\subseteq \Gamma(H_x)\cap \Omega'_x$ and this implies the claim.

\vspace{2mm}

\noindent We now define a subset $\gls{GammaGH}$ of $\Gamma(H)$ as follows: $x\in \Gamma(G,H)$ if and only if $H''_x$ is an anisotropic torus (and hence $G''_x$ also). By \ref{eq 11.1.2}, we may also see $\Gamma(G,H)$ as a subset of $\Gamma(G)$. Notice that $\Gamma(G,H)$ is a subset of $\Gamma_{\elli}(G)$ that contains $1$. We now equip $\Gamma(G,H)$ with a topology, which is finer than the one induced from $\Gamma(G)$, and a measure. For this, we need to give a more concrete description of $\Gamma(G,H)$. Consider the following set $\gls{undercalT}$ of subtori of $U(W)$: $T\in\underline{\mathcal{T}}$ if and only if there exists a non-degenerate subspace $W''\subset W$ (possibly $W''=0$) such that $T$ is a maximal elliptic subtorus of $U(W'')$. For such a torus $T$, let us denote by $\gls{Tnat}$ the open Zariski subset of elements $t\in T$ which are regular in $U(W'')$ acting without the eigenvalue $1$ on $W''$. Then, $\Gamma(G,H)$ is the set of conjugacy classes that meet

$$\displaystyle \bigcup_{T\in\underline{\mathcal{T}}} T_{\natural}(F)$$

\noindent Indeed, for all $x\in \Gamma(G,H)$ (identified with one of its representatives in $U(W)_{\ssi}(F)$), we have $H''_x\in\underline{\mathcal{T}}$ and $x\in (H''_x)_{\natural}(F)$ whereas on the other hand if $x\in T_\natural(F)$ for some $T\in\underline{\mathcal{T}}$, then $H''_x=T$. For $T\in \underline{\mathcal{T}}$, the non-degenerate subspace $W''\subset W$ such that $T$ is a maximal torus of $U(W'')$ is unique (since we have $W''=W''_x$ for all $x\in T_\natural(F)$) and we shall denote it by $\gls{W''_T}$. We will also set $\gls{W(T)}$ for the Weyl group $W(U(W''_T),T)$. Let us now fix a set of representatives $\gls{calT}$ for the $U(W)(F)$-conjugacy classes in $\underline{\mathcal{T}}$. Then, we have a natural bijection

\begin{align}\label{eq 11.1.3}
\displaystyle \Gamma(G,H)\simeq \bigsqcup_{T\in\mathcal{T}} T_\natural(F)/W(T)
\end{align}

\noindent Indeed, the map that associates to an element of the right hand side its conjugacy class is a surjection onto $\Gamma(G,H)$. That this map is injective is an easy application of Witt's theorem. Now, the right hand side of \ref{eq 11.1.3} has a natural topology and we transfer it to $\Gamma(G,H)$. Moreover, we equip $\Gamma(G,H)$ with the unique regular Borel measure such that
$$\displaystyle \int_{\Gamma(G,H)}\varphi(x)dx=\sum_{T\in\mathcal{T}} \lvert W(T)\rvert^{-1}\nu(T)\int_{T(F)} \varphi(t)dt$$
for all $\varphi\in C_c(\Gamma(G,H))$. Recall that $\nu(T)$ is the only positive factor such that the total mass of $T(F)$ for the measure $\nu(T)dt$ is one. Note that $1$ is an atom for this measure whose mass is equal to $1$ (this corresponds to the contribution of the trivial torus in the formula above).

\vspace{2mm}

\noindent More generally, for all $x\in H_{\ssi}(F)$ we may construct a subset $\gls{GammaGxHx}$ of $\Gamma(G_x)$ which is equipped with its own topology and measure as follows. By \ref{eq 11.1.1}, we have a decomposition $\Gamma(G_x)=\Gamma(G'_x)\times \Gamma(G''_x)$. Since the triple $(G'_x,H'_x,\xi_x)$ is a GGP triple, the previous construction provides us with a space $\Gamma(G'_x,H'_x)$ of semi-simple conjugacy classes in $G'_x(F)$. On the other hand, we define $\Gamma(G''_x,H''_x)$ to be the image of $\Gamma_{ani}(H''_x)$ (the set of anisotropic conjugacy classes in $H''_x(F)$, cf.\ Section \ref{section 1.8}) by the natural inclusion $\Gamma(H''_x)\subset \Gamma(G''_x)$. In Section \ref{section 1.8}, we already equipped $\Gamma(G''_x,H''_x)=\Gamma_{ani}(H''_x)$ with a topology and a measure. We now set
$$\Gamma(G_x,H_x)=\Gamma(G'_x,H'_x)\times \Gamma(G''_x,H''_x)$$
and we equip this set with the product of the topologies and the measures defined on $\Gamma(G'_x,H'_x)$ and $\Gamma(G''_x,H''_x)$. Note that $\Gamma(G_x,H_x)=\emptyset$ unless $x\in G(F)_{\elli}$ (because otherwise $\Gamma_{ani}(H''_x)=\emptyset$). The following lemma establish a link between $\Gamma(G_x,H_x)$ and $\Gamma(G,H)$:

\vspace{3mm}

\begin{lem}\label{lemma 11.0.1}
Let $\Omega_x\subseteq G_x(F)$ be a $G$-good open neighborhood of $x$ and set $\Omega=\Omega_x^G$. Then, if $\Omega_x$ is sufficiently small, the restriction of the natural map $\Gamma(G_x)\to \Gamma(G)$ to $\Omega_x\cap \Gamma(G_x,H_x)$ induces an isomorphism of topological spaces
$$\Omega_x\cap \Gamma(G_x,H_x)\simeq \Omega\cap \Gamma(G,H)$$
preserving measures.
\end{lem}

\vspace{3mm}

\noindent\ul{Proof}: First, note that the restriction of the natural map $\Gamma(G_x)\to \Gamma(G)$ to $\Omega_x\cap \Gamma(G_x,H_x)$ is injective by definition of a $G$-good open subset. As a first step towards the proof of the lemma, we show the following

\vspace{3mm}

\begin{num}
\item\label{eq 11.1.17} If $\Omega_x$ is sufficiently small, the image of $\Omega_x\cap \Gamma(G_x,H_x)$ in $\Gamma(G)$ is $\Omega\cap\Gamma(G,H)$.
\end{num}

\vspace{3mm}

\noindent By \ref{eq 11.1.2}, if $\Omega_x$ is sufficiently small, every conjugacy class in $\Omega\cap \Gamma(G,H)$ has a representative in $\Omega_x\cap H_{x,ss}(F)$. Moreover, every conjugacy class in $\Omega_x\cap \Gamma(G_x,H_x)$ also have a representative in $\Omega_x\cap H_{x,ss}(F)$. Let $y\in\Omega_x\cap H_{x,ss}(F)$. We need only prove that the $G(F)$-conjugacy class of $y$ belongs to $\Gamma(G,H)$ if and only if the $G_x(F)$-conjugacy class of $y$ belongs to $\Gamma(G_x,H_x)$. Let $y'\in H'_x(F)$ and $y''\in H''_x(F)$ for the components of $y$ relative to the decomposition $H_x=H'_x\times H''_x$. Since $\Omega_x$ is a $G$-good open neighborhood of $x$, we have
$$G_y=(G_x)_y=(G'_x)_{y'}\times (G''_x)_{y''}$$
Moreover, we have the decomposition
$$(G'_x)_y=(G'_x)'_{y'}\times (G'_x)''_{y'}$$
and it follows that
$$G_y=(G'_x)_{y'}\times \left[ (G'_x)''_{y'}\times (G''_x)_{y''}\right]$$
We easily check that this corresponds to the decomposition $G_y=G'_y\times G''_y$, that is
$$G'_y=(G'_x)'_{y'} \mbox{ and } G''_y=(G'_x)''_{y'}\times (G''_x)_{y''}$$
By definition, $y$ belongs to $\Gamma(G,H)$ if and only if $G''_y$ is an anisotropic torus whereas $y$ belongs to $\Gamma(G_x,H_x)$ if and only if both $(G'_x)''_{y'}$ and $(G''_x)_{y''}$ are anisotropic tori. By the last equality above, these two conditions are equivalent. This ends the proof of \ref{eq 11.1.17}.

\vspace{2mm}

\noindent To finish the proof of the lemma, it only remains to show that the induced map 
$$\Omega_x\cap \Gamma(G_x,H_x)\to \Omega\cap \Gamma(G,H)$$
is locally a topological isomorphism that preserves measures. Let $y\in \Omega_x\cap \Gamma(G_x,H_x)$ and identify it with one of its representative in $H_{x,ss}(F)$. Set $T=G''_y$. We introduce as before the components $y'$ and $y''$ of $y$ relative to the decomposition $H_x=H'_x\times H''_x$. Then, $T$ is an anisotropic torus and we have a decomposition
$$T=T'\times T''$$
where $T'=(G'_x)''_{y'}$ and $T''=(G''_x)_{y''}$. We have $\Gamma(G_x,H_x)=\Gamma(G'_x,H'_x)\times \Gamma(G''_x,H''_x)$. By definition of the topological structures and measures on both $\Gamma(G'_x,H'_x)$ and $\Gamma(G''_x,H''_x)$, the two maps
$$t\in T'(F)\mapsto ty'\in \Gamma(G'_x,H'_x)$$
$$t\in T''(F)\mapsto ty''\in \Gamma(G''_x,H''_x)$$
are locally near $1$ topological isomorphisms that preserve measures (where we equip $T'(F)$ and $T''(F)$ with the unique Haar measures of total mass $1$). Hence, the map
$$t\in T(F)\mapsto ty\in \Gamma(G_x,H_x)$$
is also locally near $1\in T(F)$ a topological isomorphism that preserves measures (again equipping $T(F)$ with the Haar measure of total mass $1$). On the other hand, by definition of the topology and the measure on $\Gamma(G,H)$, the map
$$t\in T(F)\mapsto ty\in \Gamma(G,H)$$
has exactly the same property. Hence the map
$$\Omega_x\cap \Gamma(G_x,H_x)\to \Omega\cap \Gamma(G,H)$$
is locally near $y$ a topological isomorphism that preserves measures. This ends the proof of the lemma.$\blacksquare$

\vspace{3mm}

\noindent We also define a subset $\gls{GammaLieGH}$ of $\Gamma(\mathfrak{g})$, again equipped with a topology and a measure, as follows. For all $X\in\mathfrak{u}(W)_{\ssi}(F)$, we have decompositions

\begin{align}\label{eq 11.1.4}
\displaystyle G_X=G'_X\times G''_X,\;\;\; H_X=H'_X\times H''_X
\end{align}

\noindent where this time
$$\displaystyle \gls{G'_X}=U(W'_X)\times U(V'_X),\;\; \gls{G''_X}=U(W''_X)_X\times U(W''_X)_X,\;\; \gls{H'_X}=U(W'_X)\ltimes N_X\;\; \gls{H''_X}=U(W''_X)_X$$
for $\gls{W'_X}$, $\gls{V'_X}$ the kernels of $X$ acting on $W$ and $V$ respectively, $\gls{W''_X}$ the image of $X$ in $W$ and $N_X$ the centralizer of $X$ in $N$. Again, the decompositions \ref{eq 11.1.4} still hold for every $X\in\mathfrak{h}_{\ssi}(F)$ and they depend on the choice of representative in the conjugacy class of $X$ only up to an inner automorphism. We now define $\Gamma^{\Lie}(G,H)$ to be the set of semi-simple conjugacy classes $X\in \Gamma(\mathfrak{h})$ such that $H''_X$ is an anisotropic torus. The obvious analog of \ref{eq 11.1.2} for the Lie algebra allows us to identify $\Gamma^{\Lie}(G,H)$ with a subset of $\Gamma(\mathfrak{g})$. Notice that $\Gamma^{\Lie}(G,H)$ is a subset of $\Gamma_{\elli}(\mathfrak{g})$ that contains $0$. Moreover, fixing a set of tori $\mathcal{T}$ as before, we have a natural identification

\begin{align}\label{eq 11.1.5}
\displaystyle \Gamma^{\Lie}(G,H)=\bigsqcup_{T\in\mathcal{T}} \mathfrak{t}_\natural(F)/W(T)
\end{align}

\noindent where for $T\in \mathcal{T}$, $\gls{tnat}$ denotes the Zariski open subset consisting of elements $X\in\mathfrak{t}$ that are regular in $\mathfrak{u}(W''_T)$ and acting without the eigenvalue $0$ on $W''_T$. By the identification \ref{eq 11.1.5}, $\Gamma^{\Lie}(G,H)$ inherits a natural topology. Moreover, we equip $\Gamma^{\Lie}(G,H)$ with the unique regular Borel measure such that
$$\displaystyle\int_{\Gamma^{\Lie}(G,H)}\varphi(X) dX=\sum_{T\in\mathcal{T}} \lvert W(T)\rvert^{-1}\nu(T)\int_{\mathfrak{t}(F)}\varphi(X)dX$$
for all $\varphi\in C_c(\Gamma^{\Lie}(G,H))$. Note that $0$ is an atom for this measure whose associated mass is $1$. The following lemma establish a link between $\Gamma^{\Lie}(G,H)$ and $\Gamma(G,H)$:

\vspace{3mm}

\begin{lem}\label{lemma 11.0.2}
Let $\omega\subseteq \mathfrak{g}(F)$ be a $G$-excellent open neighborhood of $0$ (see Section \ref{section 3.3} for this notion) and set $\Omega=\exp(\omega)$. Then, the exponential map induces a topological isomorphism
$$\omega\cap \Gamma^{\Lie}(G,H)\simeq \Omega\cap\Gamma(G,H)$$
preserving measures.
\end{lem}

\vspace{3mm}

\noindent\ul{Proof}: Since $\omega\subseteq \mathfrak{g}(F)$ is a $G$-excellent open subset, the exponential map induces a bijection
$$\omega\cap \Gamma(\mathfrak{g})\simeq \Omega\cap \Gamma(G)$$
Moreover, this bijection restricts to a bijection between $\omega\cap \Gamma^{\Lie}(G,H)$ and $\Omega\cap \Gamma(G,H)$ since for $X\in \omega\cap \Gamma(\mathfrak{h})$ we have $G''_X=G''_{e^X}$. To finish the proof of the lemma, it only remains to show that this bijection is locally a topological isomorphism that preserves measure. Let $X\in\omega\cap \Gamma^{\Lie}(G,H)$ and set $T=G''_X$. We have the following commutative diagram
$$\xymatrix{
Y\in\mathfrak{t}(F) \ar[d]^{\exp} \ar@{|->}[r] & Y+X\in \omega\cap\Gamma^{\Lie}(G,H)\ar[d]^{\exp}\\
t\in T(F) \ar@{|->}[r] & te^X\in \Omega\cap \Gamma(G,H)}$$
where the maps at the top and on the left are locally near $0$ topological isomorphisms that preserve measures and the map at the bottom is locally near $1$ a topological isomorphism that preserves measures. This shows that the map on the right is locally near $X$ a topological isomorphism that preserves measure. This ends the proof of the lemma. $\blacksquare$

\subsection{The linear forms $m_{\geom}$ and $m_{\geom}^{\Lie}$}\label{section 11.1}

\noindent In this section, we define continuous linear forms $m_{\geom}$ and $m_{\geom}^{\Lie}$ on $QC(G(F))$ and $SQC(\mathfrak{g}(F))$ respectively. These linear forms are the main ingredients in the formulation of the three theorems of Section \ref{section 11.2}. They will be precisely defined in the proposition below but first we need to introduce some determinant functions. We keep the notation introduced in the previous section. First, we set
$$\displaystyle \gls{Deltax}=D^G(x)D^H(x)^{-2}$$
for all $x\in H_{\ssi}(F)$ where we recall that $D^G(x)=\lvert \det(1-\Ad(x))_{\mid \mathfrak{g}/\mathfrak{g}_x}\rvert$ and $D^H(x)=\lvert \det(1-\Ad(x))_{\mid \mathfrak{h}/\mathfrak{h}_x}\rvert$ (see Section \ref{section 1.1}). Then, we easily check that

\begin{align}\label{eq 11.1.6}
\displaystyle \Delta(x)=\left\lvert N_{E/F}(\det(1-x)_{\mid W''_x})\right\rvert
\end{align}

\noindent for all $x\in H_{\ssi}(F)$. Similarly, we define
$$\displaystyle \gls{DeltaX}=D^G(X)D^H(X)^{-2}$$
for all $X\in\mathfrak{h}_{\ssi}(F)$ and we have the equality

\begin{align}\label{eq 11.1.7}
\displaystyle \Delta(X)=\left\lvert N_{E/F}(\det(X_{\mid W''_X}))\right\rvert
\end{align}

\noindent for all $X\in\mathfrak{h}_{\ssi}(F)$. Let $\omega\subseteq \mathfrak{g}(F)$ be a $G$-excellent open neighborhood of $0$. Then, we easily check that $\omega\cap\mathfrak{h}(F)$ is a $H$-excellent open neighborhood of $0$ (cf. the remark at the end of Section \ref{section 3.3}). We may thus set
$$\displaystyle \gls{jGHX}=j^H(X)^2j^G(X)^{-1}$$
for all $X\in \omega\cap\mathfrak{h}(F)$. By \ref{eq 3.3.1}, we have

\begin{align}\label{eq 11.1.8}
\displaystyle j^H_G(X)=\Delta(X)\Delta(e^X)^{-1}
\end{align}

\noindent for all $X\in\omega\cap\mathfrak{h}_{\ssi}(F)$. Note that $j^H_G$ is a smooth, positive and $H(F)$-invariant function on $\omega\cap \mathfrak{h}(F)$. It actually extends (not uniquely although) to a smooth, positive and $G(F)$-invariant function on $\omega$. This can be seen as follows. We can embed the groups $H_1=U(W)$ and $H_2=U(V)$ into GGP triples $(G_1,H_1,\xi_1)$ and $(G_2,H_2,\xi_2)$. Then the function $X=(X_W,X_V)\in \omega\mapsto j_{G_1}^{H_1}(X_W)^{1/2}j_{G_2}^{H_2}(X_V)^{1/2}$ is easily seen to be such an extension (using for example the equality \ref{eq 11.1.8} above). We will always assume that such an extension has been chosen and we will still denote it by $j_G^H$.

\vspace{2mm}

\noindent Let $x\in H_{\ssi}(F)$. Then, we define
$$\displaystyle \gls{Deltaxy}=D^{G_x}(y)D^{H_x}(y)^{-2}$$
for all $y\in H_{x,ss}(F)$. On the other hand, since the triple $(G'_x,H'_x,\xi'_x)$ is a GGP triple, the previous construction yields a function $\Delta^{G'_x}$ on $H'_{x,ss}(F)$. We easily check that

\begin{align}\label{eq 11.1.9}
\displaystyle \Delta_x(y)=\Delta^{G'_x}(y')
\end{align}

\noindent for all $y=(y',y'')\in H_{x,ss}(F)=H'_{x,ss}(F)\times H''_{x,ss}(F)$.

\vspace{2mm}

\noindent Let $\Omega_x\subseteq G_x(F)$ be a $G$-good open neighborhood of $x$. Then, it is easy to see that $\Omega_x\cap H(F)\subseteq H_x(F)$ is a $H$-good open neighborhood of $x$ (cf. the remark at the end of Section \ref{section 3.2}). This allow us to set
$$\displaystyle \gls{etaGxHy}=\eta^H_x(y)^2 \eta^G_x(y)^{-1}$$
for all $y\in \Omega_x\cap H(F)$. By \ref{eq 3.2.4}, we have

\begin{align}\label{eq 11.1.10}
\displaystyle \eta_{G,x}^H(y)=\Delta_x(y)\Delta(y)^{-1}
\end{align}

\noindent for all $y\in \Omega_x\cap H_{\ssi}(F)$. Note that $\eta_{G,x}^H$ is a smooth, positive and $H_x(F)$-invariant function on $\Omega_x\cap H(F)$. It actually extends (not uniquely although) to a smooth, positive and $G_x(F)$-invariant function on $\Omega_x$. We will always still denote by $\eta_{G,x}^H$ such an extension.

\vspace{2mm}

\noindent The definitions of the distributions $m_{\geom}$ and $m_{\geom}^{\Lie}$ are contained in the following proposition:

\begin{prop}\label{proposition 11.1.1}
\begin{enumerate}[(i)]
\item Let $\theta\in QC(G(F))$. Then, for all $s\in\mathbb{C}$ such that $Re(s)>0$  the integral
$$\displaystyle\int_{\Gamma(G,H)}D^G(x)^{1/2}c_\theta(x) \Delta(x)^{s-1/2}dx$$
is absolutely convergent and the limit
$$\displaystyle \gls{mgeomtheta}:=\lim\limits_{s\to 0^+} \int_{\Gamma(G,H)}D^G(x)^{1/2}c_\theta(x) \Delta(x)^{s-1/2}dx$$
exists. Similarly, for all $x\in H_{\ssi}(F)$ and for all $\theta_x\in QC(G_x(F))$, the integral
$$\displaystyle\int_{\Gamma(G_x,H_x)}D^{G_x}(y)^{1/2}c_{\theta_x}(y) \Delta_x(y)^{s-1/2}dy$$
is absolutely convergent for all $s\in \mathbb{C}$ such that $Re(s)>0$ and the limit
$$\displaystyle \gls{mxgeomthetax}:=\lim\limits_{s\to 0^+} \int_{\Gamma(G_x,H_x)}D^{G_x}(y)^{1/2}c_{\theta_x}(y) \Delta_x(y)^{s-1/2}dy$$
exists. Moreover, $m_{\geom}$ is a continuous linear form on $QC(G(F))$ and for all $x\in H_{\ssi}(F)$, $m_{x,\geom}$ is a continuous linear form on $QC(G_x(F))$.

\item Let $x\in H_{\ssi}(F)$ and let $\Omega_x\subseteq G_x(F)$ be a $G$-good open neighborhood of $x$ and set $\Omega=\Omega_x^G$. Then, if $\Omega_x$ is sufficiently small, we have
$$\displaystyle m_{\geom}(\theta)=m_{x,\geom}((\eta_{G,x}^H)^{1/2}\theta_{x,\Omega_x})$$
for all $\theta\in QC_c(\Omega)$.

\item Let $\theta\in QC_c(\mathfrak{g}(F))$. Then, for all $s\in\mathbb{C}$ such that $Re(s)>0$  the integral
$$\displaystyle\int_{\Gamma^{\Lie}(G,H)}D^G(X)^{1/2}c_\theta(X) \Delta(X)^{s-1/2}dX$$
is absolutely convergent and the limit
$$\displaystyle \gls{mgeomLietheta}:=\lim\limits_{s\to 0^+} \int_{\Gamma^{\Lie}(G,H)}D^G(X)^{1/2}c_\theta(X) \Delta(X)^{s-1/2}dX$$
exists. Moreover, $m^{\Lie}_{\geom}$ is a continuous linear form on $QC_c(\mathfrak{g}(F))$ that extends continuously to $SQC(\mathfrak{g}(F))$ and we have
$$\displaystyle m^{\Lie}_{\geom}(\theta_\lambda)=\lvert \lambda\rvert^{\delta(G)/2} m_{\geom}^{\Lie}(\theta)$$
for all $\theta\in SQC(\mathfrak{g}(F))$ and all $\lambda\in F^\times$ (recall that $\theta_\lambda(X)=\theta(\lambda^{-1}X)$ for all $X\in\mathfrak{g}_{\reg}(F)$).

\item Let $\omega\subseteq \mathfrak{g}(F)$ be a $G$-excellent open neighborhood of $0$ and set $\Omega=\exp(\omega)$. Then, we have
$$\displaystyle m_{\geom}(\theta)=m^{\Lie}_{\geom}((j_G^H)^{1/2}\theta_\omega)$$
for all $\theta\in QC_c(\Omega)$.
\end{enumerate}
\end{prop}

\vspace{2mm}

\noindent\ul{Remark}: By Proposition \ref{proposition 4.5.1} 1.(i), in the integral defining $m_{\geom}(\theta)$ above only the conjugacy classes $x\in \Gamma(G,H)$ such that $G_x$ is quasi-split contribute. This means that we could have replaced $\Gamma(G,H)$ by the, usually smaller, set $\Gamma_{\quasid}(G,H,\xi)$ consisting of conjugacy classes $x\in \Gamma(G,H)$ such that $G_x$ is quasi-split. Of course, a similar remark applies to $m_{x,\geom}(\theta_x)$ and $m_{\geom}^{\Lie}(\theta)$.

\vspace{2mm}

\noindent\ul{Proof}:

\begin{enumerate}[(i)]
\item We first show the following

\vspace{3mm}

\begin{num}
\item\label{eq 11.1.11} For all $\theta\in QC(G(F))$ and all $s\in\mathbb{C}$ such that $Re(s)>0$, the integral
$$\displaystyle m_{\geom,s}(\theta):=\int_{\Gamma(G,H)}D^G(x)^{1/2}c_\theta(x)\Delta(x)^{s-1/2}dx$$
is absolutely convergent. Moreover, $m_{\geom,s}$ defines a continuous linear form on $QC(G(F))$ (for all $Re(s)>0$).
\end{num}

\vspace{3mm}

\noindent Indeed, since $\Gamma(G,H)$ is compact modulo conjugation and the function $(D^G)^{1/2}c_\theta$ is locally bounded by a continuous semi-norm on $QC(G(F))$ for all $\theta\in QC(G(F))$ (see Proposition \ref{proposition 4.5.1}), it is sufficient to show that the integral

$$\displaystyle \int_{\Gamma(G,H)} \Delta(x)^{s-1/2}dx$$

\noindent is absolutely convergent for $Re(s)>0$. By definition of the measure on $\Gamma(G,H)$ and \ref{eq 11.1.6}, this is a straightforward application of Lemma \ref{lemma B.1.2}(i). Similarly, we prove that

\vspace{3mm}

\begin{num}
\item\label{eq 11.1.12} For all $x\in H_{\ssi}(F)$, all $\theta_x\in QC(G_x(F))$ and all $s\in\mathbb{C}$ such that $Re(s)>0$, the integral
$$\displaystyle m_{x,\geom,s}(\theta_x):=\int_{\Gamma(G_x,H_x)}D^{G_x}(y)^{1/2} c_{\theta_x}(y)\Delta_x(y)^{s-1/2}dy$$
is absolutely convergent. Moreover, $m_{x,\geom,s}$ defines a continuous linear form on $QC(G_x(F))$ (for all $Re(s)>0$).
\end{num}

\vspace{3mm}

\noindent Assume for a moment that the limits

\begin{align}\label{eq 11.1.13}
m_{\geom}(\theta)=\lim\limits_{s\to 0^+} m_{\geom,s}(\theta)
\end{align}

\noindent and

\begin{align}\label{eq 11.1.14}
m_{x,\geom}(\theta_x)=\lim\limits_{s\to 0^+} m_{x,\geom,s}(\theta_x)
\end{align}

\noindent exist for all $\theta\in QC(G(F))$, all $x\in H_{\ssi}(F)$ and all $\theta_x\in QC(G_x(F))$. Then, by the uniform boundedness principle (cf.\ Appendix \ref{section A.1}), $m_{\geom}$ and $m_{\geom,x}$ ($x\in H_{\ssi}(F)$) will automatically be continuous linear forms on $QC(G(F))$ and $QC(G_x(F))$ respectively. Hence, it only remains to show that the limits \ref{eq 11.1.13} and \ref{eq 11.1.14} always exist. We prove this by induction on $\dim(G)$ i.e. we assume that the result holds for every GGP triple $(G',H',\xi')$ with $\dim(G')<\dim(G)$ (the result is trivial when $\dim(G)=1$) an we will show that it holds for $(G,H,\xi)$. We first show the following

\vspace{3mm}

\begin{num}
\item\label{eq 11.1.15} Let $x\in H_{\ssi}(F)$ and assume that $x\neq 1$. Then the limit \ref{eq 11.1.14} exists for all $\theta_x\in QC(G_x(F))$.
\end{num}

\vspace{3mm}

\noindent Since $G_x=G'_x\times G''_x$, by Proposition \ref{proposition 4.4.1}(v) we have $QC(G_x(F))=QC(G'_x(F))\widehat{\otimes}_p QC(G''_x(F))$. By \ref{eq A.5.3}, it thus suffices to show that the limit \ref{eq 11.1.14} exists for every quasi-character $\theta_x$ of the form $\theta_x=\theta'_x\otimes\theta''_x$ where $\theta'_x\in QC(G'_x(F))$ and $\theta''_x\in QC(G''_x(F))$. Fix such a quasi-character. Since $\Gamma(G_x,H_x)=\Gamma(G'_x,H'_x)\times \Gamma(G''_x,H''_x)$, using \ref{eq 11.1.9} we have
$$\displaystyle m_{x,\geom,s}(\theta_x)=\int_{\Gamma(G'_x,H'_x)} D^{G'_x}(y)^{1/2}c_{\theta'_x}(y)\Delta^{G'_x}(y)^{s-1/2} dy\times \int_{\Gamma(G''_x,H''_x)}D^{G''_x}(y)^{1/2}c_{\theta''_x}(y)dy$$
for all $s\in\mathbb{C}$ such that $Re(s)>0$. Recall that the triple $(G'_x,H'_x,\xi'_x)$ is a GGP triple. For $Re(s)>0$, let us denote by $m_{\geom,s}^{G'_x}$ the distribution on $QC(G'_x(F))$ defined the same way as $m_{\geom,s}$ but for this GGP triple instead of $(G,H,\xi)$. Then, the first integral above is equal to $m_{\geom,s}^{G'_x}(\theta'_x)$. Since $\dim(G'_x)<\dim(G)$, the induction hypothesis tells us that $m_{\geom,s}^{G'_x}(\theta'_x)$ has a limit as $s\to 0^+$ and this ends the proof of \ref{eq 11.1.15}.

\vspace{2mm}

\noindent We are now left with proving that the limit \ref{eq 11.1.13} exists for all $\theta\in QC(G(F))$. First, from Lemma \ref{lemma 11.0.1} we deduce that

\vspace{3mm}

\begin{num}
\item\label{eq 11.1.18} Let $x\in H_{\ssi}(F)$, $\Omega_x\subseteq G_x(F)$ be a $G$-good open neighborhood of $x$ and set $\Omega=\Omega_x^G$. Then, if $\Omega_x$ is sufficiently small, we have the equality
$$m_{\geom,s}(\theta)=m_{x,\geom,s}((\eta^H_{x,G})^{1/2-s}\theta_{x,\Omega_x})$$
for all $\theta\in QC_c(\Omega)$ and all $s\in\mathbb{C}$ such that $Re(s)>0$.
\end{num} 

\vspace{3mm}

\noindent Indeed, let $\theta\in QC_c(\Omega)$, by Lemma \ref{lemma 11.0.1} we have

\[\begin{aligned}
\displaystyle m_{\geom,s}(\theta) & =\int_{\Omega\cap \Gamma(G,H)}D^G(y)^{1/2}c_\theta(y)\Delta(y)^{s-1/2}dy \\
 & =\int_{\Omega_x\cap \Gamma(G_x,H_x)}D^G(y)^{1/2}c_\theta(y) \Delta(y)^{s-1/2}dy
\end{aligned}\]

\noindent for all $s\in\mathbb{C}$ such that $Re(s)>0$. For all $y\in \Omega_x\cap \Gamma(G_x,H_x)$, we have $D^G(y)^{1/2}c_\theta(y)=D^{G_x}(y)^{1/2}c_{\theta_{x,\Omega_x}}(y)$ (Proposition \ref{proposition 4.5.1}.1(iv)) and by \ref{eq 11.1.10} we also have $\Delta(y)=\eta^H_{x,G}(y)^{-1}\Delta_x(y)$. The equality \ref{eq 11.1.18} follows.

\vspace{2mm}

\noindent We now prove

\vspace{3mm}

\begin{num}
\item\label{eq 11.1.19} Let $\theta\in QC(G(F))$ and assume that $1\notin \Supp(\theta)$. Then, the limit \ref{eq 11.1.13} exists.
\end{num}

\vspace{3mm}

\noindent Indeed, since $\Gamma(G,H)$ is a subset of $\Gamma(H)$ which is compact modulo conjugation, by an invariant partition of unity process (Proposition \ref{proposition 3.1.1}(ii)), we are immediately reduced to proving \ref{eq 11.1.19} for $\theta\in QC_c(\Omega)$ where $\Omega$ is of the form $\Omega=\Omega_x^G$ for some $x\in H_{\ssi}(F)$ different from $1$ and some $\Omega_x\subseteq G_x(F)$ a $G$-good open neighborhood of $x$ that we can take as small as we want. In particular by \ref{eq 11.1.18}, if we take $\Omega_x$ sufficiently small, we have
$$\displaystyle m_{\geom,s}(\theta)=m_{x,\geom,s}((\eta_{x,G}^H)^{1/2-s}\theta_{x,\Omega_x})$$
for all $Re(s)>0$. By \ref{eq 11.1.15}, we know that the continuous linear forms $m_{x,\geom,s}$ converge point-wise to a continuous linear form $m_{x,\geom}$ on $QC(G_x(F))$ as $s\to 0$. By the uniform boundedness principle, it implies that $m_{x,\geom,s}$ converges uniformly on compact subsets of $QC(G_x(F))$. In particular, if we can show that the function $s\mapsto (\eta_{x,G}^H)^{1/2-s}\theta_{x,\Omega_x}$ has a limit in $QC(G_x(F))$ as $s\to 0$, then we will be done by the above equality. Clearly, $(\eta_{x,G}^H)^{1/2-s}$ converges to $(\eta_{x,G}^H)^{1/2}$ in $C^\infty(\Omega_x)^{G_x}$. Hence, by Proposition \ref{proposition 4.4.1}(iv), we have
$$\displaystyle \lim\limits_{s\to 0^+}(\eta_{x,G}^H)^{1/2-s}\theta_{x,\Omega_x}=(\eta_{x,G}^H)^{1/2}\theta_{x,\Omega_x}$$
in $QC_c(\Omega_x)$ and so also in $QC(G_x(F))$. This ends the proof of \ref{eq 11.1.19}.

\vspace{2mm}

\noindent Let $\omega\subseteq \mathfrak{g}(F)$ be a $G$-excellent open neighborhood of $0$ and set $\Omega=\exp(\omega)$. By \ref{eq 11.1.19}, and since there exists $\varphi\in C^\infty(G(F))^G$ such that $\Supp(\varphi)\subseteq \Omega$ and $\varphi=1$ in a neighborhood of $1$, we are left with proving

\vspace{3mm}

\begin{num}
\item\label{eq 11.1.20} The limit \ref{eq 11.1.13} exists for all $\theta\in QC_c(\Omega)$
\end{num}

\vspace{3mm}

\noindent As a first step towards the proof of \ref{eq 11.1.20}, we claim the following

\vspace{3mm}

\begin{num}
\item\label{eq 11.1.21} For all $\theta\in QC_c(\mathfrak{g}(F))$ and all $s\in\mathbb{C}$ such that $Re(s)>0$, the integral
$$\displaystyle m_{\geom,s}^{\Lie}(\theta)=\int_{\Gamma^{\Lie}(G,H)}D^G(X)^{1/2}c_{\theta}(X) \Delta(X)^{s-1/2}dX$$
converges absolutely. Moreover, $m_{\geom,s}^{\Lie}$ defines, for $Re(s)>0$, a continuous linear form on $QC_c(\mathfrak{g}(F))$, and we have
$$m_{\geom,s}(\theta)=m_{\geom,s}^{\Lie}((j_G^H)^{1/2-s}\theta_\omega)$$
for all $\theta\in QC_c(\Omega)$ and all $s\in\mathbb{C}$ such that $Re(s)>0$.
\end{num} 

\vspace{3mm}

\noindent Indeed, the first part of the claim can be proved in a way similar to \ref{eq 11.1.11} using Lemma \ref{lemma B.1.2}(ii) whereas the second part of the claim is proved as \ref{eq 11.1.18} using Lemma \ref{lemma 11.0.2} instead of Lemma \ref{lemma 11.0.1}.

\vspace{2mm}

\noindent Using again the uniform boundedness principle (cf.\ Appendix \ref{section A.1}), if the limit

\begin{align}\label{eq 11.1.23}
m_{\geom}^{\Lie}(\theta):=\lim\limits_{s\to 0^+}m^{\Lie}_{\geom,s}(\theta)
\end{align}

\noindent exists for all $\theta\in QC_c(\mathfrak{g}(F))$, then $m_{\geom}^{\Lie}$ will automatically be a continuous linear form on $QC_c(\mathfrak{g}(F))$ and the linear forms $m^{\Lie}_{\geom,s}$ will converge uniformly on compacta to $m_{\geom}^{\Lie}$ as $s\to 0^+$. Hence, by \ref{eq 11.1.21}, to prove \ref{eq 11.1.20} we only need to show the existence of the limit \ref{eq 11.1.23} for all $\theta\in QC_c(\mathfrak{g}(F))$. As a first step, we show

\vspace{3mm}

\begin{num}
\item\label{eq 11.1.24} For all $\theta\in QC_c(\omega)$ such that $0\notin \Supp(\theta)$, the limit \ref{eq 11.1.23} exists.
\end{num}

\vspace{3mm}

\noindent Let $L\subseteq \omega$ be a closed invariant neighborhood of $0$. Then, it suffices to show that the limit \ref{eq 11.1.23} exists for all $\theta\in QC_c(\omega-L)$. By \ref{eq 11.1.21}, we have
$$\displaystyle m_{\geom,s}^{\Lie}(\theta)=m_{\geom,s}((j^H_G\circ\exp)^{s-1/2}\theta_{\Omega})$$
for all $\theta\in QC_c(\omega)$ and $Re(s)>0$. Hence, it suffices to prove that the limit

\begin{align}\label{eq 11.1.25}
\displaystyle \lim\limits_{s\to 0^+}m_{\geom,s}((j^H_G\circ\exp)^{s-1/2}\theta)
\end{align}

\noindent exists for all $\theta\in QC_c(\Omega-e^L)$. Since $e^L$ is a closed invariant neighborhood of $1$ in $G(F)$, we already know, by \ref{eq 11.1.19}, that the limit \ref{eq 11.1.13} exists for all $\theta\in QC_c(\Omega-e^L)$. But, again by the uniform boundedness principle and Proposition \ref{proposition 4.4.1}(iv), it follows that the linear forms $m_{\geom,s}$ converge uniformly on compact subsets of $QC_c(\Omega-L)$ as $s\to 0^+$, hence the limit \ref{eq 11.1.25} exists for all $\theta\in QC_c(\Omega-e^L)$.

\vspace{2mm}

\noindent We now extend slightly \ref{eq 11.1.24} and prove

\vspace{3mm}

\begin{num}
\item\label{eq 11.1.26} For all $\theta\in QC_c(\mathfrak{g}(F))$ such that $0\notin \Supp(\theta)$, the limit \ref{eq 11.1.23} exists.
\end{num} 

\vspace{3mm}

\noindent Recall that for $\theta\in QC(\mathfrak{g}(F))$ and $\lambda\in F^\times$, the quasi-character $\theta_\lambda$ is defined by $\theta_\lambda(X)=\theta(\lambda^{-1}X)$, $X\in\mathfrak{g}_{\reg}(F)$. Of course, if $\theta\in QC_c(\mathfrak{g}(F))$ is such that $0\notin \Supp(\theta)$ then the quasi-character $\theta_\lambda$ has the same property for all $\lambda\in F^\times$. Moreover, if $\lvert \lambda\rvert$ is sufficiently small then $\theta_\lambda\in QC_c(\omega)$. Hence, \ref{eq 11.1.26} will follow from \ref{eq 11.1.24} once we have established the following

\vspace{3mm}

\begin{num}
\item\label{eq 11.1.27} Denote by $M_\lambda$, $\lambda\in F^\times$ the endomorphism of $QC_c(\mathfrak{g}(F))$ given by $M_\lambda(\theta)=\lvert \lambda\rvert^{-\delta(G)/2}\theta_\lambda$. There exists a positive integer $d>0$ such that for all $\theta\in QC_c(\mathfrak{g}(F))$ and all $\lambda\in F^\times$, we have
$$\displaystyle \lim\limits_{s\to 0^+}m_{\geom,s}^{\Lie}\left[(M_\lambda-1)^d\theta\right]=0$$
\end{num}

\vspace{3mm}

\noindent Let $d>0$ be an integer. For all $X\in \Gamma^{\Lie}(G,H)$, let us denote by $d_X$ the dimension of the torus $G''_X$. By definition of the measure on $\Gamma^{\Lie}(G,H)$ and Proposition \ref{proposition 4.5.1} 2.(iv), it is easy to see that
$$\displaystyle m^{\Lie}_{\geom,s}\left[(M_\lambda-1)^d\theta\right]=\int_{\Gamma^{\Lie}(G,H)} (\lvert \lambda\rvert^{2sd_X}-1)^d D^G(X)^{1/2}c_\theta(X) \Delta(X)^{s-1/2}dX$$
for all $\theta\in QC_c(\mathfrak{g}(F))$, all $\lambda\in F^\times$ and all $Re(s)>0$. By this and the definition of the measure on $\Gamma^{\Lie}(G,H)$, to establish \ref{eq 11.1.27} it is now sufficient to show the following

\vspace{3mm}

\begin{num}
\item\label{eq 11.1.28} For all $T\in \mathcal{T}$ there exists a positive integer $d>0$ such that for all $\theta\in QC_c(\mathfrak{g}(F))$, we have
$$\displaystyle \lim\limits_{s\to 0^+}s^d\int_{\mathfrak{t}(F)}D^G(X)^{1/2}c_\theta(X)\Delta(X)^{s-1/2} dX=0$$
\end{num}

\vspace{3mm}

\noindent The point \ref{eq 11.1.28} is a straightforward consequence of Lemma \ref{lemma B.1.2}(ii) and of the fact that for all $\theta\in QC_c(\mathfrak{g}(F))$ the function $(D^G)^{1/2}c_\theta$ is locally bounded and compactly supported modulo conjugation. This ends the proof of \ref{eq 11.1.27} and hence of \ref{eq 11.1.26}.

\vspace{2mm}

\noindent We now again improve \ref{eq 11.1.26} slightly. More precisely, we prove

\vspace{3mm}

\begin{num}
\item\label{eq 11.1.29} Let $\theta\in QC_c(\mathfrak{g}(F))$ and assume that $c_{\theta,\mathcal{O}}(0)=0$ for all $\mathcal{O}\in \Nil_{\reg}(\mathfrak{g})$. Then, the limit \ref{eq 11.1.23} exists.
\end{num}

\vspace{3mm}

\noindent Choose $d$ so that \ref{eq 11.1.27} holds. Fix $\lambda\in F^\times$ such that $\lvert \lambda\rvert\neq 1$. Let $\theta\in QC_c(\mathfrak{g}(F))$ and assume that $c_{\theta,\mathcal{O}}(0)=0$ for all $\mathcal{O}\in \Nil_{\reg}(\mathfrak{g})$. Then, by Proposition \ref{proposition 4.6.1}(i) we may find $\theta_1\in QC_c(\mathfrak{g}(F))$ and $\theta_2\in QC_c(\mathfrak{g}(F))$ such that $\theta=(M_\lambda-1)^d\theta_1+\theta_2$ and $\theta_2$ is supported away from $0$. We can now deduce \ref{eq 11.1.29} from \ref{eq 11.1.26} and \ref{eq 11.1.27}.

\vspace{2mm}

\noindent If $G$ is not quasi-split, then \ref{eq 11.1.29} already shows that the limit \ref{eq 11.1.23} always exists (as in this case $\Nil_{\reg}(\mathfrak{g})=\emptyset$), ending the proof of (i). We assume henceforth that $G$ is quasi-split. Then, $G$ has two regular nilpotent orbits $\mathcal{O}_+, \mathcal{O}_-\in \Nil_{\reg}(\mathfrak{g})$ (this follows from the description of regular nilpotent orbits of unitary groups given in Section \ref{section 6.1}). Let $\varphi\in C^\infty(\mathfrak{g}(F))^G$ be an invariant smooth function which is compactly supported modulo conjugation and equals $1$ in some neighborhood of $0$. Set $\theta_+=\varphi\widehat{j}(\mathcal{O}_+,.)$ and $\theta_-=\varphi\widehat{j}(\mathcal{O}_-,.)$. These are compactly supported quasi-characters on $\mathfrak{g}(F)$ and by \ref{eq 11.1.29}, it only remains to show that the two limits

\begin{align}\label{eq 11.1.30}
\displaystyle \lim\limits_{s\to 0^+}m_{\geom,s}^{\Lie}(\theta_+) \mbox{ and } \lim\limits_{s\to 0^+} m_{\geom,s}^{\Lie}(\theta_-)
\end{align}

\noindent exist. The quasi-character $\widehat{j}(\mathcal{O}_+,.)+\widehat{j}(\mathcal{O}_-,.)$ is parabolically induced from a maximal torus of a Borel subgroup of $G$ (cf.\ Chapter \ref{section 3.4} for the notion of induction of quasi-characters and more particularly \ref{eq 3.4.6} for the case at hand). Hence the support of $\theta_++\theta_-$ intersects $\Gamma_{\elli}(G)$ only at the center of $\mathfrak{g}(F)$. It follows from \ref{eq 1.7.5} that the support of $\theta_++\theta_-$ intersects $\Gamma^{\Lie}(G,H)$ only at $0$. Hence, we have

\begin{align}\label{eq 11.1.31}
\displaystyle m_{\geom,s}^{\Lie}(\theta_+)+m_{\geom,s}^{\Lie}(\theta_-)=c_{\theta_+}(0)+c_{\theta_-}(0)=1
\end{align}

\noindent for all $Re(s)>0$. On the other hand, let $\lambda\in N(E^\times)-F^\times$ and $d$ as in \ref{eq 11.1.27}. Then, multiplication by $\lambda$ exchanges the two orbits $\mathcal{O}_+$ and $\mathcal{O}_-$ (this follows from the description of regular nilpotent orbits of unitary groups given in Section \ref{section 6.1}). It follows that $(M_\lambda-1)^d\theta_+$ coincides near $0$ with a nonzero multiple of $\theta_+-\theta_-$. Hence, by \ref{eq 11.1.26} and \ref{eq 11.1.27}, the limit
$$\displaystyle \lim\limits_{s\to 0^+}m_{\geom,s}^{\Lie}(\theta_+-\theta_-)$$
exists. Combining this with \ref{eq 11.1.31}, this shows that the two limits \ref{eq 11.1.30} exist and this ends the proof of (i).

\item This follows from \ref{eq 11.1.18}.

\item We have already proved the first part during the proof of (i). For all $T\in\mathcal{T}$, let us introduce the distribution
$$\displaystyle m_{\geom,s}^{\Lie,T}(\theta)=\lvert W(T)\rvert^{-1}\int_{\mathfrak{t}(F)}D^G(X)^{1/2}c_\theta(X) \Delta(X)^{s-1/2}dX,\;\;\; \theta\in QC_c(\mathfrak{g}(F))$$
for $Re(s)>0$. Again using Lemma \ref{lemma B.1.2}(ii), this integral is absolutely convergent and defines a continuous linear form on $QC_c(\mathfrak{g}(F))$. Moreover, if we follow the proof of (i) closely, then we see that it actually also works ``torus by torus" so that the limit
$$\lim\limits_{s\to 0^+}m_{\geom,s}^{\Lie,T}(\theta)$$
exists for all $\theta\in QC_c(\mathfrak{g}(F))$ and all $T\in\mathcal{T}$. Since, we have
$$\displaystyle m^{\Lie}_{\geom,s}(\theta_\lambda)=\sum_{T\in\mathcal{T}} \lvert \lambda\rvert^{\delta(G)/2+2s\dim(T)}m_{\geom,s}^{\Lie,T}(\theta)$$
for all $\theta\in QC_c(\mathfrak{g}(F))$, all $\lambda\in F^\times$ and all $Re(s)>0$, it now easily follows that
$$m_{\geom}^{\Lie}(\theta_\lambda)=\lvert \lambda\rvert^{\delta(G)/2}m_{\geom}^{\Lie}(\theta)$$
for all $\theta\in QC_c(\mathfrak{g}(F))$ and all $\lambda\in F^\times$. The fact that $m_{\geom}^{\Lie}$ extends continuously to $SQC(\mathfrak{g}(F))$ now follows from Proposition \ref{proposition 4.6.1}(ii).

\item This follows from \ref{eq 11.1.21}. $\blacksquare$
\end{enumerate}

\subsection{Geometric multiplicity and parabolic induction}\label{section 11.1bis}

Let $L$ be a Levi subgroup  of $G$. Then, as in Section \ref{section 8.4}, we may decompose $L$ as a product

$$L=L^{GL}\times \widetilde{G}$$

\noindent where $L^{GL}$ is a product of general linear groups over $E$ and $\widetilde{G}$ belongs to a GGP triple $(\widetilde{G},\widetilde{H},\widetilde{\xi})$ which is well-defined up to $\widetilde{G}(F)$-conjugation. In particular, Proposition \ref{proposition 11.1.1} applied to this GGP triple provides us with a continuous linear form $m_{\geom}^{\widetilde{G}}$ on $QC(\widetilde{G}(F))$. We define a continuous linear form $m_{\geom}^L$ on $QC(L(F))=QC(L^{GL}(F))\widehat{\otimes}_p QC(\widetilde{G}(F))$ by setting
$$\displaystyle \gls{mgeomL}(\theta^{GL}\otimes \widetilde{\theta})=m_{\geom}^{\widetilde{G}}(\widetilde{\theta})c_{\theta^{GL}}(1)$$
for all $\theta^{GL}\in QC(L^{GL}(F))$ and all $\widetilde{\theta}\in QC(\widetilde{G}(F))$. The next lemma is precisely \cite{Beu1} Lemme 17.2.1 in the $p$-adic case and moreover the proof of {\it loc. cit.} adapts verbatim to the real case once we replace references of {\it loc. cit.} to Lemme 2.3 of \cite{Wa4} by references to Proposition \ref{proposition 4.7.1}(ii).

\vspace{2mm}

\begin{lem}\label{lemma 11.1.1}
Let $\theta^L\in QC(L(F))$ and set $\theta=i_L^G(\theta^L)$. Then, we have

$$m_{\geom}(\theta)=m^L_{\geom}(\theta^L)$$
\end{lem}

\subsection{Statement of three theorems}\label{section 11.2}

\noindent Set

$$\displaystyle \gls{Jgeomf}=m_{\geom}(\theta_f),\;\;\;\mbox{for all } f\in \mathcal{C}_{\scusp}(G(F))$$

$$\displaystyle \gls{mgeompi}=m_{\geom}(\theta_\pi),\;\;\;\mbox{for all } \pi\in R_{\tempe}(G)$$

$$\displaystyle \gls{JgeomLief}=m^{\Lie}_{\geom}(\theta_f),\;\;\;\mbox{for all } f\in \mathcal{S}_{\scusp}(\mathfrak{g}(F))$$

\noindent Recall that in Chapter \ref{section 7} we have defined two continuous linear forms $J(.)$ and $J^{\Lie}(.)$ on $\mathcal{C}_{\scusp}(G(F))$ and $\mathcal{S}_{\scusp}(G(F))$ respectively. Recall also that in Section \ref{section 6.3}, we have defined a multiplicity $m(\pi)$ for all $\pi\in \Temp(G)$. The goal of the next 5 sections is to prove the following three theorems.

\vspace{3mm}

\begin{theo}\label{theorem 11.2.1}
We have

$$J(f)=J_{\geom}(f)$$

\noindent for all $f\in \mathcal{C}_{\scusp}(G(F))$.
\end{theo}

\vspace{2mm}

\begin{theo}\label{theorem 11.2.2}
We have

$$m(\pi)=m_{\geom}(\pi)$$

\noindent for all $\pi\in \Temp(G)$.
\end{theo}

\vspace{2mm}

\begin{theo}\label{theorem 11.2.3}
We have

$$J^{\Lie}(f)=J_{\geom}^{\Lie}(f)$$

\noindent for all $f\in \mathcal{S}_{\scusp}(\mathfrak{g}(F))$.
\end{theo}

\vspace{2mm}

\noindent The proof is by induction on $\dim(G)$ (the case $\dim(G)=1$ being obvious). Hence, until the end of Section \ref{section 11.7}, we make the following induction hypothesis

\vspace{4mm}

\hspace{5mm} {\bf (HYP)} Theorem \ref{theorem 11.2.1}, Theorem \ref{theorem 11.2.2} and Theorem \ref{theorem 11.2.3} hold for all GGP triples

\vspace{0.5pt}

\hspace{19mm} $(G',H',\xi')$ such that $\dim(G')<\dim(G)$.

\subsection{Equivalence of Theorem \ref{theorem 11.2.1} and Theorem \ref{theorem 11.2.2}}\label{section 11.3}

\noindent Recall that $R_{\tempe}(G)$ stands for the space of complex virtual tempered representations of $G(F)$. In Section \ref{section 2.7}, we have defined the subspace $R_{ind}(G)$ of ``properly induced" virtual representations and also a set $\mathcal{X}_{\elli}(G)$ of virtual representations well-defined up to multiplication by a scalar of module $1$ (the set of elliptic representations).

\vspace{2mm}

\begin{prop}\label{proposition 11.3.1}
Assume the induction hypothesis {\bf (HYP)}. Then,
\begin{enumerate}[(i)]
\item Let $\pi\in R_{ind}(G)$. Then, we have

$$m(\pi)=m_{\geom}(\pi)$$

\item For all $f\in \mathcal{C}_{\scusp}(G(F))$, we have the equality

$$\displaystyle J(f)=J_{\geom}(f)+\sum_{\pi\in \mathcal{X}_{\elli}(G)} D(\pi) \widehat{\theta}_f(\pi) \left(m(\overline{\pi})-m_{\geom}(\overline{\pi})\right)$$

\noindent the sum in the right hand side being absolutely convergent.

\item Theorem \ref{theorem 11.2.1} and Theorem \ref{theorem 11.2.2} are equivalent.

\item There exists a unique continuous linear form $\gls{Jqc}$ on $QC(G(F))$ such that

\begin{itemize}
\renewcommand{\labelitemi}{$\bullet$}

\item $J(f)=J_{\qc}(\theta_f)$ for all $f\in \mathcal{C}_{\scusp}(G(F))$;

\item $\Supp(J_{\qc})\subseteq G(F)_{\elli}$.
\end{itemize}
\end{enumerate}
\end{prop}

\vspace{2mm}

\noindent\ul{Proof}:

\begin{enumerate}[(i)]

\item Let $\pi\in R_{ind}(G)$. By linearity, we may assume that there exists a proper parabolic subgroup $Q=LU_Q$ of $G$ and a representation $\sigma\in \Temp(L)$ such that $\pi=i_Q^G(\sigma)$. Then, as in Section \ref{section 8.4}, we may decompose $L$ as a product

$$L=L^{GL}\times \widetilde{G}$$

\noindent where $L^{GL}$ is a product of general linear groups over $E$ and $\widetilde{G}$ belongs to a GGP triple $(\widetilde{G},\widetilde{H},\widetilde{\xi})$. We have accordingly a decomposition

$$\sigma=\sigma^{GL}\boxtimes \widetilde{\sigma}$$

\noindent of $\sigma$, where $\sigma^{GL}\in \Temp(L^{GL})$ and $\widetilde{\sigma}\in \Temp(\widetilde{G})$. By Corollary \ref{corollary 8.6.1}(i), we have

$$m(\pi)=m(\widetilde{\sigma})$$

\noindent On the other hand, by \ref{eq 3.4.3} and Lemma \ref{lemma 11.1.1}, we have

$$m_{\geom}(\pi)=m_{\geom}^{\widetilde{G}}(\widetilde{\sigma})c_{\sigma^{GL}}(1)$$

\noindent Finally, by the induction hypothesis {\bf (HYP)} applied to the GGP triple $(\widetilde{G},\widetilde{H},\widetilde{\xi})$, we have

$$m(\widetilde{\sigma})=m_{\geom}^{\widetilde{G}}(\widetilde{\sigma})$$

\noindent It only remains to see that $c_{\sigma^{GL}}(1)=1$. But this follows from Proposition \ref{proposition 4.8.1}(i), since every tempered representation of a general linear group admits a Whittaker model (cf.\ \cite{Ze} Theorem 9.7).

\item Let $f\in \mathcal{C}_{\scusp}(G(F))$. The sum on the right hand side of the identity is absolutely convergent by Lemma \ref{lemma 5.4.2}, Proposition \ref{proposition 4.8.1}(ii) and \ref{eq 2.7.2}. By Theorem \ref{theorem 9.1.1}, we have the equality

\begin{align}\label{eq 11.3.1}
\displaystyle J(f) & =\int_{\mathcal{X}(G)}D(\pi)\widehat{\theta}_f(\pi)m(\overline{\pi})d\pi \\
\nonumber & =\int_{\mathcal{X}_{ind}(G)}D(\pi)\widehat{\theta}_f(\pi)m(\overline{\pi})d\pi+\sum_{\pi\in \mathcal{X}_{\elli}(G)}D(\pi)\widehat{\theta}_f(\pi)m(\overline{\pi})
\end{align}

\noindent Since the (virtual) representations in $\mathcal{X}_{ind}(G)$ are properly induced, by (i), we have
$$m(\overline{\pi})=m_{\geom}(\overline{\pi})$$
for all $\pi\in\mathcal{X}_{ind}(G)$. Hence, \ref{eq 11.3.1} becomes
$$\displaystyle J(f)=\int_{\mathcal{X}(G)}D(\pi)\widehat{\theta}_f(\pi)m_{\geom}(\overline{\pi})d\pi+\sum_{\pi\in\mathcal{X}_{\elli}(G)} D(\pi) \widehat{\theta}_f(\pi) \left(m(\overline{\pi})-m_{\geom}(\overline{\pi})\right)$$
Finally, by Proposition \ref{proposition 5.6.1}(ii), we have the equality
$$\displaystyle J_{\geom}(f)=\int_{\mathcal{X}(G)}D(\pi)\widehat{\theta}_f(\pi)m_{\geom}(\overline{\pi})d\pi$$
and (ii) follows.

\item It is obvious from (ii) that Theorem \ref{theorem 11.2.2} implies Theorem \ref{theorem 11.2.1}. For the converse, it suffices to use (i), Corollary \ref{corollary 5.7.1}(iv) and \ref{eq 2.7.1}.

\item The unicity follows from Corollary \ref{corollary 5.7.1}(i). For the existence, it suffices to set
$$\displaystyle J_{\qc}(\theta)=m_{\geom}(\theta)+\sum_{\pi\in\mathcal{X}_{\elli}(G)}D(\pi)\left(m(\overline{\pi})-m_{\geom}(\overline{\pi})\right) \int_{\Gamma_{\elli}(G)} D^G(x)\theta(x)\theta_\pi(x)dx$$
for all $\theta\in QC(G(F))$. That it defines a continuous linear form follows from Corollary \ref{corollary 5.7.1}(ii) and (iii), Proposition \ref{proposition 4.8.1}(ii) and \ref{eq 2.7.2}. $\blacksquare$
\end{enumerate}

\subsection{Semi-simple descent and the support of $J_{\qc}-m_{\geom}$}\label{section 11.4}

\begin{prop}\label{proposition 11.4.1}
Assume the induction hypothesis {\bf (HYP)}. Let $\theta\in QC(G(F))$ and assume that $1\notin \Supp(\theta)$. Then, we have

$$J_{\qc}(\theta)=m_{\geom}(\theta)$$
\end{prop}

\vspace{2mm}

\noindent\ul{Proof}: Since $J_{\qc}$ and $m_{\geom}$ are both supported in $\Gamma_{\elli}(G)$, by a partition of unity process, we only need to prove the equality of the proposition for $\theta\in QC_c(\Omega)$ where $\Omega$ is a completely $G(F)$-invariant open subset of $G(F)$ of the form $\Omega_x^G$ for some $x\in G(F)_{\elli}$, $x\neq 1$, and some $G$-good open neighborhood $\Omega_x\subseteq G_x(F)$ of $x$. Moreover, we may take $\Omega_x$ as small as we want. In particular, we will assume that $\Omega_x$ is relatively compact modulo conjugation.

\vspace{2mm}

\noindent Assume first that $x$ is not conjugate to any element of $H_{\ssi}(F)$. Then, since $\Gamma(H)$ is closed in $\Gamma(G)$ (by \ref{eq 11.1.2}), if $\Omega_x$ is chosen sufficiently small, we would have $\Omega\cap \Gamma(H)=\emptyset$. In this case, both sides of the equality are easily seen to be zero (for all $\theta\in QC_c(\Omega)$).

\vspace{2mm}

\noindent Assume now that $x$ is conjugate to some element of $H_{\ssi}(F)$. We may as well assume that $x\in H_{\ssi}(F)$. Then, we have the decompositions

$$G_x=G'_x\times G''_x,\;\;\; G''_x=H''_x\times H''_x$$

\noindent Shrinking $\Omega_x$ if necessary, we may assume that $\Omega_x$ decomposes as a product

$$\Omega=\Omega'_x\times\left(\Omega''_x\times \Omega''_x\right)$$

\noindent where $\Omega'_x\subseteq G'_x(F)$ (resp.\ $\Omega_x''\subseteq H''_x(F)$) is open and completely $G'_x(F)$-invariant (resp.\ completely $H''_x(F)$-invariant). Note that $x$ is elliptic in both $G'_x$ and $H''_x$. Hence, by Corollary \ref{corollary 5.7.1}(i), shrinking $\Omega_x$ further if necessary, we may assume, and we do, that the linear maps

$$f'_x\in\mathcal{S}_{\scusp}(\Omega'_x)\mapsto \theta_{f'_x}\in QC_c(\Omega'_x)$$

$$f''_x\in \mathcal{S}_{\scusp}(\Omega''_x)\mapsto \theta_{f''_x}\in QC_c(\Omega''_x)$$

\noindent have dense image. Since $QC_c(\Omega)=QC_c(\Omega'_x)\widehat{\otimes}_p QC_c(\Omega''_x)\widehat{\otimes}_p QC_c(\Omega''_x)$ (Proposition \ref{proposition 4.4.1}(v)) and $J_{\qc}$ and $m_{\geom}$ are continuous linear forms on $QC_c(\Omega)$, we only need to prove the equality of the proposition for quasi-characters $\theta\in QC_c(\Omega)$ such that $\theta_{x,\Omega_x}=\theta_{f_x}$ for some $f_x\in \mathcal{S}_{\scusp}(\Omega_x)$ which further decomposes as a tensor product $f_x=f'_x\otimes \left(f''_{x,1}\otimes f''_{x,2}\right)$ where $f'_x\in \mathcal{S}_{\scusp}(\Omega'_x)$ and $f''_{x,1},f''_{x,2}\in\mathcal{S}_{\scusp}(\Omega''_x)$. So fix a quasi-character $\theta\in QC_c(\Omega)$ with this property and fix functions $f_x$, $f'_x$, $f''_{x,1}$ and $f''_{x,2}$ as before.

\vspace{2mm}

\noindent Consider a map as in Proposition \ref{proposition 5.7.1}, and set $f=\widetilde{f_x}\in \mathcal{S}_{\scusp}(\Omega)$. Then, we have

\begin{align}\label{eq 11.4.1}
J(f)=J_{\qc}(\theta_f)=J_{\qc}(\theta)
\end{align}

\noindent (Notice that here $Z_G(x)=G_x$ since $G_{\der}$ is simply-connected). Let us denote by $J^{G'_x}$ the continuous linear form on $\mathcal{C}_{\scusp}(G'_x(F))$ associated to the GGP triple $(G'_x,H'_x,\xi'_x)$. Also, let us denote by $J^{A,H''_x}$ the continuous bilinear form on $\mathcal{C}_{\scusp}(H''_x(F))$ introduced in Section \ref{section 5.5} (where we replace $G$ by $H''_x$). Recall that we have defined a $H_x(F)$-invariant smooth and positive function $\eta_{x,G}^H$ on $\Omega_x\cap H_x(F)=(\Omega'_x\cap H'_x(F))\times \Omega''_x$. It is easy to see, using the formulas \ref{eq 11.1.9}, \ref{eq 11.1.10} and \ref{eq 11.1.6}, that this function factorizes through the projection $\Omega_x\cap H_x(F)\to \Omega''_x$. We shall identify $\eta_{x,G}^H$ with the function it defines on $\Omega''_x$. Let us show the following

\vspace{3mm}

\begin{num}
\item\label{eq 11.4.2} If $\Omega_x$ is sufficiently small, we have 
$$\displaystyle J(f)=J^{G'_x}(f'_x)J^{A,H''_x}((\eta^H_{x,G})^{1/2}f''_{x,1},f''_{x,2}).$$
\end{num}

\vspace{3mm}

\noindent The intersection $\Omega_x\cap H_x(F)\subseteq H_x(F)$ is a $H$-good open neighborhood of $x$ (cf. the remark at the end of Section \ref{section 3.2}). Moreover, by \ref{eq 11.1.2} if $\Omega_x$ is sufficiently small, we have $\Omega\cap H(F)=\left(\Omega_x\cap H_x(F)\right)^H$. We henceforth assume $\Omega_x$ that sufficiently small. Then, by \ref{eq 3.2.5}, we have

\begin{align}\label{eq 11.4.3}
\displaystyle J(f) & =\int_{H(F)\backslash G(F)} \int_{H(F)} {}^gf(h)\xi(h) dhdg \\
\nonumber & =\int_{H(F)\backslash G(F)} \int_{H_x(F)\backslash H(F)}\int_{H_x(F)} \eta_x^H(h_x)\; {}^{hg}f(h_x)\xi_x(h_x)dh_xdhdg \\
\nonumber & =\int_{H(F)\backslash G(F)} \int_{H_x(F)\backslash H(F)}\int_{H_x(F)} \eta_{x,G}^H(h_x)^{1/2}({}^{hg}f)_{x,\Omega_x}(h_x)\xi_x(h_x)dh_xdhdg
\end{align}

\noindent Assume one moment that the exterior double integral above is absolutely convergent. Then, we would have

\[\begin{aligned}
\displaystyle J(f) & =\int_{H_x(F)\backslash G(F)}\int_{H_x(F)} \eta_{x,G}^H(h_x)^{1/2}({}^gf)_{x,\Omega_x}(h_x)\xi_x(h_x)dh_xdg \\
 & =\int_{G_x(F)\backslash G(F)}\int_{H_x(F)\backslash G_x(F)}\int_{H_x(F)} \eta_{x,G}^H(h_x)^{1/2}({}^gf)_{x,\Omega_x}(g_x^{-1}h_xg_x)\xi_x(h_x)dh_xdg_xdg
\end{aligned}\]

\noindent Introduce a function $\alpha$ on $G_x(F)\backslash G(F)$ as in Proposition \ref{proposition 5.7.1}. Let $g\in G(F)$. Up to translating $g$ by an element of $G_x(F)$, we may assume that $({}^gf)_{x,\Omega_x}=\alpha(g)f_x$. Then, the interior integral above decomposes as

\[\begin{aligned}
\displaystyle \alpha(g) & \int_{H_x(F)\backslash G_x(F)}\int_{H_x(F)} \eta_{x,G}^H(h_x)^{1/2}f_x(g_x^{-1}h_xg_x)\xi_x(h_x)dh_xdg_xdg \\
 & =\alpha(g) \int_{H'_x(F)\backslash G'_x(F)}\int_{H'_x(F)} f'_x({g'_x}^{-1}h'_xg'_x)\xi'_x(h'_x)dh'_xdg'_x \\
 & \times \int_{H''_x(F)} \int_{H''_x(F)} \eta_{x,G}^H(h''_x)^{1/2}f''_{x,1}({g''_x}^{-1}h''_xg''_x)f''_{x,2}(h''_x)dh''_xdg''_x
\end{aligned}\]

\noindent We recognize the two integrals above: the first one is $J^{G'_x}(f'_x)$ and the second one is $J^{A,H''_x}((\delta_{x,G}^H)^{1/2}f''_{x,1},f''_{x,2})$ (Note that the center of $H''_x(F)$ is compact since $x$ is elliptic). By Theorem \ref{proposition 7.1.1}(ii) and Theorem \ref{theorem 5.5.1}(ii) and since the function $\alpha$ is compactly supported, this shows that the exterior double integral of the last line of \ref{eq 11.4.3} is absolutely convergent. Moreover, since we have
$$\displaystyle \int_{G_x(F)\backslash G(F)} \alpha(g)dg=1$$
this also proves \ref{eq 11.4.2}.

\vspace{2mm}

\noindent We assume from now on that $\Omega_x$ is sufficiently small so that \ref{eq 11.4.2} holds. By the induction hypothesis {\bf (HYP)}, we have
$$J^{G'_x}(f'_x)=m_{\geom}^{G'_x}(\theta_{f'_x})$$
On the other hand, by Theorem \ref{theorem 5.5.1}(iv), we have
$$\displaystyle J_A^{H''_x}((\eta_{x,G}^H)^{1/2}f''_{x,1},f''_{x,2})=\int_{\Gamma_{ani}(H''_x)} \eta_{G,x}^H(y)^{1/2}D^{G''_x}(y)^{1/2}\theta_{f''_{x,1}}(y) \theta_{f''_{x,2}}(y) dy$$
(Notice that here both $f''_{x,1}$ and $f''_{x,2}$ are strongly cuspidal, hence the terms corresponding to $\Gamma(H''_x)-\Gamma_{ani}(H''_x)$ vanish). Moreover, it is easy to check that
$$\displaystyle m_{\geom,x}((\eta_{G,x}^H)^{1/2}\theta_{f_x})=m_{\geom}^{G'_x}(\theta_{f'_x})\times \int_{\Gamma_{ani}(H''_x)} \eta_{x,G}^H(y)^{1/2}D^{G''_x}(y)^{1/2}\theta_{f''_{x,1}}(y) \theta_{f''_{x,2}}(y) dy$$
Hence, by \ref{eq 11.4.1}, \ref{eq 11.4.2} and Proposition \ref{proposition 11.1.1}(ii), we have
$$J_{\qc}(\theta)=J(f)=m_{\geom,x}((\eta_{G,x}^H)^{1/2}\theta_{f_x})=m_{\geom}(\theta)$$
This ends the proof of the proposition. $\blacksquare$

\subsection{Descent to the Lie algebra and equivalence of Theorem \ref{theorem 11.2.1} and Theorem \ref{theorem 11.2.3}}\label{section 11.5}

\noindent Let $\omega\subseteq \mathfrak{g}(F)$ be a $G(F)$-excellent open neighborhood of $0$ and set $\Omega=\exp(\omega)$. Recall that for any quasi-character $\theta\in QC(\mathfrak{g}(F))$ and all $\lambda\in F^\times$, $\theta_\lambda$ denotes the quasi-character given by $\theta_\lambda(X)=\theta(\lambda^{-1}X)$ for all $X\in\mathfrak{g}_{\reg}(F)$.

\begin{prop}\label{proposition 11.5.1}
Assume the induction hypothesis {\bf (HYP)}. Then,
\begin{enumerate}[(i)]
\item For all $f\in \mathcal{S}_{\scusp}(\Omega)$, we have

$$J(f)=J^{\Lie}((j_G^H)^{1/2}f_\omega)$$

\item There exists a unique continuous linear form $\gls{JqcLie}$ on $SQC(\mathfrak{g}(F))$ such that

$$J^{\Lie}(f)=J_{\qc}^{\Lie}(\theta_f)$$

\noindent for all $f\in \mathcal{S}_{\scusp}(\mathfrak{g}(F))$. Moreover, we have

$$J_{\qc}^{\Lie}(\theta_{\lambda})=\lvert \lambda\rvert^{\delta(G)/2}J_{\qc}^{\Lie}(\theta)$$

\noindent for all $\theta\in SQC(\mathfrak{g}(F))$ and all $\lambda\in F^\times$.

\item Theorem \ref{theorem 11.2.1} and Theorem \ref{theorem 11.2.3} are equivalent.

\item Let $\theta\in SQC(\mathfrak{g}(F))$ and assume that $0\notin \Supp(\theta)$. Then, we have

$$J^{\Lie}_{\qc}(\theta)=m^{\Lie}_{\geom}(\theta)$$
\end{enumerate}
\end{prop}

\vspace{2mm}

\noindent\ul{Proof}:

\begin{enumerate}[(i)]
\item The intersection $\omega_{\mathfrak{h}}=\omega\cap\mathfrak{h}(F)\subseteq \mathfrak{h}(F)$ is a $H$-excellent open neighborhood of $0$ (cf. the remark at the end of Section \ref{section 3.3}). Recall that
$$j_G^H(X)=j^H(X)^{2}j^G(X)^{-1}$$
for all $X\in \omega_{\mathfrak{h}}$. Hence, by \ref{eq 3.3.2}, we have

\[\begin{aligned}
\displaystyle \int_{H(F)} f(g^{-1}hg)\xi(h)dh & =\int_{\omega_{\mathfrak{h}}} j^H(X)f(g^{-1}e^X g) \xi(X)dX \\
 & =\int_{\mathfrak{h}(F)} j_G^H(X)^{1/2} f_{\omega}(g^{-1}Xg)\xi(X)dX
\end{aligned}\]

\noindent for all $f\in \mathcal{S}(\Omega)$ and all $g\in G(F)$. The identity of the proposition follows immediately.

\item The unicity follows from Proposition \ref{proposition 5.6.1}(i). Let us prove the existence. Set

\begin{align}\label{eq 11.5.1}
\displaystyle J_{\qc}^{\Lie}(\theta)=\int_{\Gamma(\Sigma)} D^G(X)^{1/2}\widehat{\theta}(X)dX
\end{align}

\noindent for all $\theta\in SQC(\mathfrak{g}(F))$, where $\Gamma(\Sigma)$ is defined as in Section \ref{section 10.8}. By \ref{eq 1.8.2}, the integral above is absolutely convergent and $J_{\qc}^{\Lie}$ is a continuous linear form on $SQC(\mathfrak{g}(F))$. Moreover, by Theorem \ref{theorem 10.8.1}, we have
$$J^{\Lie}(f)=J_{\qc}^{\Lie}(\theta_f)$$
for all $f\in\mathcal{S}_{\scusp}(\mathfrak{g}(F))$. This shows the existence. The last claim is easy to check using the formula \ref{eq 11.5.1}.

\item By (i) and Proposition \ref{proposition 11.1.1}(iv), it is clear that Theorem \ref{theorem 11.2.3} implies Theorem \ref{theorem 11.2.1}. Let us prove the converse. Assume that Theorem \ref{theorem 11.2.1} holds. Then by (i) and Proposition \ref{proposition 11.1.1}(iv), we have the equality
$$J^{\Lie}(f)=J_{\geom}^{\Lie}(f)$$
for all $f\in \mathcal{S}_{\scusp}(\omega)$. Hence, by Proposition \ref{proposition 5.6.1}(i), we also have
$$J^{\Lie}_{\qc}(\theta)=m^{\Lie}_{\geom}(\theta)$$
for all $\theta\in QC_c(\omega)$. By the homogeneity properties of $J^{\Lie}_{\qc}$ and $m^{\Lie}_{\geom}$ (cf.\ Proposition \ref{proposition 11.1.1}(iii)), this last equality extends to $QC_c(\mathfrak{g}(F))$. By the density of $QC_c(\mathfrak{g}(F))$ in $SQC(\mathfrak{g}(F))$ (Lemma \ref{lemma 4.2.2}(v)), this identity is even true for all $\theta\in SQC(\mathfrak{g}(F))$. This implies Theorem \ref{theorem 11.2.3}.

\item First, we show that we may assume $\theta\in QC_c(\mathfrak{g}(F))$. By Lemma \ref{lemma 4.2.2}(v), we may find a sequence $(\theta_n)_{n\geqslant 1}$ in $QC_c(\mathfrak{g}(F))$ such that
$$\lim\limits_{n\to\infty}\theta_n=\theta$$
in $SQC(\mathfrak{g}(F))$. Let $\varphi\in C^\infty(\mathfrak{g}(F))^G$ be an invariant compactly supported modulo conjugation function which is equal to $1$ in a neighborhood of $0$ and such that $(1-\varphi)\theta=\theta$. By Lemma \ref{lemma 4.2.2}(iv) and the closed graph theorem, multiplication by $(1-\varphi)$ induces a continuous endomorphism of $SQC(\mathfrak{g}(F))$. Hence, we have
$$\lim\limits_{n\to\infty}(1-\varphi)\theta_n=\theta$$
and each of the quasi-characters $(1-\varphi)\theta_n$, $n\geqslant 1$, is supported away from $0$. By continuity of the linear forms $J_{\qc}^{\Lie}$ and $m_{\geom}^{\Lie}$, this shows that we may assume that $\theta\in QC_c(\mathfrak{g}(F))$. By the homogeneity properties of $J^{\Lie}_{\qc}$ and $m_{\geom}^{\Lie}$ (cf.\ Proposition \ref{proposition 11.1.1}(iii)), we may even assume that $\Supp(\theta)\subseteq \omega$. Finally, using Proposition \ref{proposition 5.6.1}(i), we only need to prove the equality for $\theta=\theta_f$ where $f\in \mathcal{S}_{\scusp}(\omega)$. It is then a consequence of (i) and Proposition \ref{proposition 11.1.1}(iv). $\blacksquare$

\end{enumerate}

\subsection{A first approximation of $J_{\qc}^{\Lie}-m_{\geom}^{\Lie}$}\label{section 11.6}

\begin{prop}\label{proposition 11.6.1}
Assume the induction hypothesis {\bf (HYP)}. Then, there exists a constant $c\in \mathbb{C}$ such that

$$J_{\qc}^{\Lie}(\theta)-m_{\geom}^{\Lie}(\theta)=c.c_\theta(0)$$

\noindent for all $\theta\in SQC(\mathfrak{g}(F))$.
\end{prop}

\vspace{2mm}

\noindent\ul{Proof}: We first prove the following weaker result

\vspace{3mm}

\begin{num}
\item\label{eq 11.6.1} There exists constants $c_{\mathcal{O}}$, $\mathcal{O}\in \Nil_{\reg}(\mathfrak{g})$, such that
$$\displaystyle J_{\qc}^{\Lie}(\theta)-m_{\geom}^{\Lie}(\theta)=\sum_{\mathcal{O}\in \Nil_{\reg}(\mathfrak{g})}c_{\mathcal{O}} c_{\theta,\mathcal{O}}(0)$$
for all $\theta\in SQC(\mathfrak{g}(F))$.
\end{num}

\vspace{3mm}

\noindent Let $\theta\in SQC(\mathfrak{g}(F))$ be such that $c_{\theta,\mathcal{O}}(0)=0$ for all $\mathcal{O}\in \Nil_{\reg}(\mathfrak{g})$. We want to show that $J_{\qc}^{\Lie}(\theta)=m_{\geom}^{\Lie}(\theta)$. Let $\lambda\in F^\times$ be such that $\lvert \lambda\rvert\neq 1$. Denote by $M_\lambda$ the operator on $SQC(\mathfrak{g}(F))$ given by $M_\lambda \theta=\lvert\lambda\rvert^{-\delta(G)/2}\theta_\lambda$. Then by Proposition \ref{proposition 4.6.1}(i), we may find $\theta_1,\theta_2\in SQC(\mathfrak{g}(F))$ such that $\theta=(M_\lambda-1)\theta_1+\theta_2$ and $\theta_2$ is supported away from $0$. Then, by Proposition \ref{proposition 11.5.1}(iv), we have $J_{\qc}^{\Lie}(\theta_2)=m_{\geom}^{\Lie}(\theta_2)$. On the other hand, by the homogeneity property of $J_{\qc}^{\Lie}$ and $m_{\geom}^{\Lie}$ (cf.\ Proposition \ref{proposition 11.5.1}(ii) and Proposition \ref{proposition 11.1.1}(iii)), we also have $J_{\qc}^{\Lie}((M_\lambda-1)\theta_1)=m_{\geom}^{\Lie}((M_\lambda-1)\theta_1)=0$. This proves \ref{eq 11.6.1}.

\vspace{2mm}

\noindent To ends the proof of the proposition, it remains to show that the coefficients $c_\mathcal{O}$, for $\mathcal{O}\in \Nil_{\reg}(\mathfrak{g})$, are all equal. If $G$ is not quasi-split, there is nothing to prove (as then $\Nil_{\reg}(\mathfrak{g})=\emptyset$). Assume now that $G$ is quasi-split. Then $G$ has two nilpotent orbits and multiplication by any element $\lambda\in F^\times-N(E^\times)$ exchanges the two orbits (this follows from the description of regular nilpotent orbits of unitary groups given in Section \ref{section 6.1}). The results now follows from \ref{eq 1.7.5} and the homogeneity property of $J_{\qc}^{\Lie}$ and $m_{\geom}^{\Lie}$. $\blacksquare$

\subsection{End of the proof}\label{section 11.7}

\noindent By Proposition \ref{proposition 11.3.1}(iii), Proposition \ref{proposition 11.5.1}(ii) and (iii), in order to finish the proof of Theorem \ref{theorem 11.2.1}, Theorem \ref{theorem 11.2.2} and Theorem \ref{theorem 11.2.3}, it only remains to show that the coefficients $c$ of Proposition \ref{proposition 11.6.1} is zero. If $G$ is not quasi-split there is nothing to prove. Assume now that $G$ is quasi-split. Fix a Borel subgroup $B\subset G$ and a maximal torus $T_{\quasid}\subset B$ (both defined over $F$). Denote by $\Gamma_{\quasid}(\mathfrak{g})$ the subset of $\Gamma(\mathfrak{g})$ consisting of the conjugacy classes that meet $\mathfrak{t}_{\quasid}(F)$. Recall that in Section \ref{section 10.8}, we have defined a subset $\Gamma(\Sigma)\subset \Gamma(\mathfrak{g})$. It consists in the conjugacy classes of the semi-simple parts of elements in the affine subspace $\Sigma(F)\subset \mathfrak{g}(F)$ defined in Section \ref{section 10.1}. We claim that

\begin{align}\label{eq 11.7.1}
\Gamma_{\quasid}(\mathfrak{g})\subseteq \Gamma(\Sigma)
\end{align}

\noindent Up to $G(F)$-conjugation, we may assume that $B$ is a good Borel subgroup (cf.\ Section \ref{section 6.4}). Then, we have $\mathfrak{h}\oplus \mathfrak{b}=\mathfrak{g}$ and it follows that 

\begin{align}\label{eq 11.7.2}
\mathfrak{h}^\perp\oplus \mathfrak{u}=\mathfrak{g}
\end{align}

\noindent where $\mathfrak{u}$ denotes the nilpotent radical of $\mathfrak{b}$. Recall that $\Sigma=\Xi+\mathfrak{h}^\perp$. From \ref{eq 11.7.2}, we easily deduce that the restriction of the natural projection $\mathfrak{b}\to \mathfrak{t}_{\quasid}$ to $\Sigma\cap\mathfrak{b}$ induces an affine isomorphism $\Sigma\cap\mathfrak{b}\simeq\mathfrak{t}_{\quasid}$ and this clearly implies \ref{eq 11.7.1}.

\vspace{2mm}

\noindent Let $\theta_0\in C_c^\infty(\mathfrak{t}_{\quasid,\reg}(F))$ be $W(G,T_{\quasid})$-invariant and such that

\begin{align}\label{eq 11.7.3}
\displaystyle \int_{\mathfrak{t}_{\quasid}(F)} D^G(X)^{1/2}\theta_0(X)dX\neq 0
\end{align}

\noindent We may extend $\theta_0$ to a smooth invariant function on $\mathfrak{g}_{\reg}(F)$, still denoted by $\theta_0$, which is zero outside $\mathfrak{t}_{\quasid,\reg}(F)^G$. Obviously, $\theta_0$ is a compactly supported quasi-character. Consider its Fourier transform $\theta=\widehat{\theta}_0$. By Proposition \ref{proposition 4.1.1}(iii), Lemma \ref{lemma 4.2.2}(iii) and \ref{eq 3.4.5}, $\theta$ is supported in $\Gamma_{\quasid}(\mathfrak{g})$. Since $\Gamma_{\quasid}(\mathfrak{g})\cap \Gamma(G,H)=\{1\}$, by definition of $m_{\geom}^{\Lie}$, we have
$$m_{\geom}^{\Lie}(\theta)=c_\theta(0)$$
On the other hand, by Proposition \ref{proposition 4.1.1}(iii), Lemma \ref{lemma 4.2.2}(iii) and Proposition \ref{proposition 4.5.1}.2(v), we have

\begin{align}\label{eq 11.7.4}
\displaystyle c_\theta(0)=\int_{\Gamma(\mathfrak{g})} D^G(X)^{1/2}\theta_0(X)c_{\widehat{j}(X,.)}(0)dX=\int_{\Gamma_{\quasid}(\mathfrak{g})} D^G(X)^{1/2}\theta_0(X) dX
\end{align}

\noindent By definition of $J_{\qc}^{\Lie}$ and \ref{eq 11.7.1}, this last term is also equal to $J_{\qc}^{\Lie}(\theta)$. Hence, we have $J_{\qc}^{\Lie}(\theta)=m_{\geom}^{\Lie}(\theta)$. Combining \ref{eq 11.7.3} with \ref{eq 11.7.4}, we see that $c_\theta(0)\neq 0$ and it follows that the constant $c$ of Proposition \ref{proposition 11.6.1} is zero. $\blacksquare$

\section{An application to the Gan-Gross-Prasad conjecture}\label{section 12}

In this chapter, we fix, as we did in Chapters \ref{section 6} to \ref{section 11}, an admissible pair of hermitian spaces $(V,W)$ and we denote by $(G,H,\xi)$ the corresponding GGP triple but we now make the following additional assumption

\begin{center}
\ul{$G$ and $H$ are quasi-split}
\end{center}

\noindent The goal of this chapter is to give an application of the multiplicity formula of Theorem \ref{theorem 11.2.2} to the so-called local Gan-Gross-Prasad conjecture. Roughly speaking, this application states that there is exactly one distinguished representation $\pi$ (i.e. one such that $m(\pi)=1$) in every (extended) tempered $L$-packet of $G(F)$. In the $p$-adic case this result was already proved by the author \cite{Beu1} and the idea of the proof goes back to Waldspurger \cite{Wa1}. As explained in the introduction, it is based on showing the existence of many cancellations when we sum the multiplicities over an extended tempered $L$-packet. For this, we use character identities between the stable characters associated to a tempered $L$-packet on $G(F)$ and its (pure) inner forms. In the $p$-adic case there are only two relevant pure inner forms to consider: $G$ itself together with a nonquasi-split one $G'$ and the character identity involves the sign $-1$. In the real case, there are far more pure inner forms to keep track of ($m+1$ precisely where $m$ is the rank of $H$) and the signs involved in the character identities alternate. To get a uniform treatment, we will use Kottwitz \cite{Kott2} general definition for these signs. We recall this in Section \ref{section 12.1} where we also introduce a certain notion of stable conjugacy for semi-simple conjugacy classes (that we call {\em strongly stable conjugacy}) which differs from the generally accepted one but will be the one relevant for us. Then in Section \ref{section 12.2} we describe the so-called {\em pure inner forms} of a GGP triple. These are needed to state the main result where we consider a GGP triple together with all its pure inner forms at the same time. In Section \ref{section 12.3}, we state our main requirement on tempered $L$-packets under the form of three hypothesis. The first two assumptions ((STAB) and (TRANS)) pertain to the aforementioned character identities. The third one, (WHITT), concerns the existence and unicity of a generic representation (with respect to a fixed Whittaker datum) in each tempered $L$-packet of $G$. That these hypothesis are satisfied in the Archimedean case follows from work of Shelstad, Kostant and Vogan. In the $p$-adic case, we give precise references to the recent work of Mok \cite{Mok} and Kaletha-Minguez-Shin-White \cite{KMSW} on the local Langlands correspondence for unitary groups where these assumptions are proved. The main result of this chapter (Theorem \ref{theorem 12.4.1}) is stated in Section \ref{section 12.4}. As a preparation for its proof, we study in Section \ref{section 12.5} strongly stable conjugacy classes inside the space of conjugacy class $\Gamma(G,H)$ introduced in Section \ref{section 11.0} (strictly speaking we consider the disjoint union of $\Gamma(G_\alpha,H_\alpha)$ over all the pure inner forms $(G_\alpha,H_\alpha,\xi_\alpha)$ of $(G,H,\xi)$). Finally, the proof of Theorem \ref{theorem 12.4.1} is given in Section \ref{section 12.6}.

\subsection{Strongly stable conjugacy classes, transfer between pure inner forms and the Kottwitz sign}\label{section 12.1}

\noindent In this section, we recall some definitions and facts that will be needed later. These considerations are general and we can forget about the GGP triple that we fixed. So, let $G$ be any connected reductive group defined over $F$. Recall that two regular elements $x,y\in G_{\reg}(F)$ are said to be {\em stably conjugate} if there exists $g\in G(\overline{F})$ such that $y=gxg^{-1}$ and $g^{-1}\sigma(g)\in G_x$ for all $\sigma\in \Gamma_F=Gal(\overline{F}/F)$. We will need to extend this definition to more general semi-simple elements. The usual notion of stable conjugacy for semi-simple elements (cf. for example \cite{Kott1}) is too weak for our purpose. The definition that we will adopt is as follows. We will say that two semi-simple elements $x,y\in G_{\ssi}(F)$ are {\em strongly stably conjugate} and we will write

$$x\gls{stabconj} y$$

\noindent if there exists $g\in G(\overline{F})$ such that $y=gxg^{-1}$ and the isomorphism $\Ad(g):G_x\simeq G_y$ is defined over $F$. This last condition has the following concrete interpretation: it means that the $1$-cocycle $\sigma\in \Gamma_F\mapsto g^{-1}\sigma(g)$ takes its values in $Z(G_x)$ the center of $G_x$ (this is because $Z(G_x)$ coincides with the centralizer of $G_x$ in $Z_G(x)$ since $G_x$ contains a maximal torus of $G$ which is its own centralizer). Moreover, for $x\in G_{\ssi}(F)$ the set of $G(F)$-conjugacy classes inside the strong stable conjugacy class of $x$ is easily seen to be in natural bijection with

$$\displaystyle \Ima\left(H^1(F,Z(G_x))\to H^1(F,Z_G(x))\right)\cap \Ker\left(H^1(F,Z_G(x))\to H^1(F,G)\right)$$

\noindent We now recall the notion of pure inner forms. A {\em pure inner form} for $G$ is defined formally as a triple $(G',\psi,c)$ where

\vspace{2mm}

\begin{itemize}
\renewcommand{\labelitemi}{$\bullet$}
\item $G'$ is a connected reductive group defined over $F$;
\item $\psi: G_{\overline{F}}\simeq G'_{\overline{F}}$ is an isomorphism defined over $\overline{F}$;
\item $c:\sigma\in \Gamma_F\to c_\sigma\in G(\overline{F})$ is a $1$-cocycle such that $\psi^{-1}{}^\sigma\psi=\Ad(c_\sigma)$ for all $\sigma\in \Gamma_F$.
\end{itemize}

\vspace{2mm}

\noindent There is a natural notion of isomorphism between pure inner forms the equivalence classes of which are naturally in bijection with $H^1(F,G)$ (the isomorphism class of $(G',\psi,c)$ being parametrized by the image of $c$ in $H^1(F,G)$). Moreover, inside an equivalence class of pure inner forms $(G',\psi,c)$, the group $G'$ is well-defined up to $G'(F)$-conjugacy. We will always assume fixed for all $\alpha\in H^1(F,G)$ a pure inner form in the class of $\alpha$ that we will denote by $(G_\alpha,\psi_\alpha,c_\alpha)$ or simply by $\gls{Galpha}$ if no confusion arises.

\vspace{2mm}

\noindent Let $(G',\psi,c)$ be a pure inner form of $G$. Then, we will say that two semi-simple elements $x\in G_{\ssi}(F)$ and $y\in G'_{\ssi}(F)$ are {\em strongly stably conjugate} and we will write

$$x\gls{stabconj} y$$

\noindent (this extends the previous notation) if there exists $g\in G(F)$ such that $y=\psi(gxg^{-1})$ and the isomorphism $\psi\circ \Ad(g):G_x\simeq G_y$ is defined over $F$. Again, the last condition has an interpretation in terms of cohomological classes: it means that the $1$-cocycle $\sigma\in \Gamma_F\mapsto g^{-1}c_\sigma \sigma(g)$ takes its values in $Z(G_x)$. For $x\in G_{\ssi}(F)$ the set of semi-simple conjugacy classes in $G'(F)$ that are strongly stably conjugate to $x$ is naturally in bijection with

$$\displaystyle \Ima\left(H^1(F,Z(G_x))\to H^1(F,Z_G(x))\right)\cap p_x^{-1}(\alpha)$$

\noindent where $\alpha\in H^1(F,G)$ parametrizes the equivalence class of the pure inner form $(G',\psi,c)$ and $p_x$ denotes the natural map $H^1(F,Z_G(x))\to H^1(F,G)$. We will need the following fact

\vspace{3mm}

\begin{num}
\item\label{eq 12.1.1} Let $y\in G'_{\ssi}(F)$ and assume that $G$ and $G'_y$ are both quasi-split. Then, the set
$$\{x\in G_{\ssi}(F);\; x\sim_{\stab} y\}$$
is non-empty.
\end{num}

\vspace{3mm}

\noindent Indeed, since $G'_y$ is quasi-split, we can fix a Borel subgroup $B_y$ of $G'_y$ and a maximal torus $T_y\subset B_y$ both defined over $F$. Consider the embedding $\iota=\psi^{-1}_{\mid T_y}:T_{y,\overline{F}}\hookrightarrow G_{\overline{F}}$. Since $\psi$ is an inner form, it is easy to see that for all $\sigma\in \Gamma_F$ the embedding ${}^\sigma \iota$ is conjugate to $\iota$. As $G$ is quasi-split, by a result of Kottwitz (\cite{Kott1} Corollary 2.2), there exists $g\in G(\overline{F})$ such that $\Ad(g)\circ \iota$ is defined over $F$. Set $x=g\iota(y)g^{-1}=g\psi^{-1}(y)g^{-1}$ and $T_x=g\iota(T_y)g^{-1}$. We claim that $x$ and $y$ are strongly stably conjugate. To see this, we first note that the $1$-cocycle $\sigma\in \Gamma_F\mapsto g^{-1}c_\sigma \sigma(g)$ takes its values in $T_x$ (as $\Ad(g)\circ \iota$ is defined over $F$). Hence, we only need to show that the map $H^1(F,Z(G_x))\to H^1(F,T_x)$ is surjective or, what amounts to the same, that the map $H^1(F,Z(G'_y))\to H^1(F,T_y)$ is surjective. Denote by $\Delta$ the set of simple roots of $T_y$ in $B_y$ (a priori these roots are not defined over $F$ but $\Gamma_F$ acts on them). We have the equality
$$\displaystyle X^*_{\overline{F}}\left(T_y/Z(G_y)\right)=\bigoplus_{\alpha\in \Delta} \mathbb{Z}\alpha$$
and the Galois action on $X^*_{\overline{F}}\left(T_y/Z(G_y)\right)$ permutes the basis $\Delta$. It follows that the torus $T_y/Z(G_y)$ is a finite product of tori of the form $R_{F'/F}\mathbb{G}_m$ where $F'$ is a finite extension of $F$ and $R_{F'/F}$ denotes the functor of restriction of scalars. In particular, by Hilbert 90 the group $H^1(F,T_y/Z(G_y))$ is trivial and it immediately follows that the morphism $H^1(F,Z(G'_y))\to H^1(F,T_y)$ is surjective.

\vspace{2mm}

\noindent We continue to consider a pure inner form $(G',\psi,c)$ of $G$. We say of a quasi-character $\theta$ on $G(F)$ that it is {\em stable} if for all regular elements $x,y\in G_{\reg}(F)$ that are stably conjugate we have $\theta(x)=\theta(y)$. Let $\theta$ and $\theta'$ be stable quasi-characters on $G(F)$ and $G'(F)$ respectively and assume moreover that $G$ is quasi-split. Then, we say that $\theta'$ is a {\em transfer} of $\theta$ if for all regular points $x\in G_{\reg}(F)$ and $y\in G'_{\reg}(F)$ that are stably conjugate we have $\theta'(y)=\theta(x)$. Note that if $\theta'$ is a transfer of $\theta$ the quasi-character $\theta'$ is entirely determined by $\theta$ (this is because every regular element in $G'(F)$ is stably conjugate to some element of $G(F)$ for example by the point \ref{eq 12.1.1} above). We will need the following:

\vspace{3mm}

\begin{num}
\item\label{eq 12.1.2} Let $\theta$ and $\theta'$ be stable quasi-characters on $G(F)$ and $G'(F)$ respectively and assume that $\theta'$ is a transfer of $\theta$. Then, for all $x\in G_{\ssi}(F)$ and $y\in G'_{\ssi}(F)$ that are strongly stably conjugate we have
$$c_{\theta'}(y)=c_\theta(x)$$
\end{num}

\vspace{3mm}

\noindent Let $x\in G_{\ssi}(F)$ and $y\in G'_{\ssi}(F)$ be two strongly stably conjugate semi-simple elements. Choose $g\in G(\overline{F})$ such that $y=\psi(gxg^{-1})$ and the isomorphism $\psi\circ \Ad(g):G_x\simeq G_y$ is defined over $F$. We will denote by $\iota$ this isomorphism. If $G_x$ and $G'_y$ are not quasi-split there is nothing to prove since by Proposition \ref{proposition 4.5.1}(i) both sides of the equality we want to establish are equal to zero. Assume now that the groups $G_x$, $G'_y$ are quasi-split. Let $B_x$ be a Borel subgroup of $G_x$ and $T_x\subset B_x$ be a maximal torus, both defined over $F$. Set $B_y=\iota(B_x)$ and $T_y=\iota(T_x)$. Then, $B_y$ is a Borel subgroup of $G'_y$ and $T_y\subset B_y$ is a maximal torus. By Proposition \ref{proposition 4.5.1}(ii), we have
$$\displaystyle D^G(x)^{1/2}c_\theta(x)=\lvert W(G_x,T_x)\rvert^{-1} \lim\limits_{x'\in T_x(F)\to x} D^G(x')^{1/2}\theta(x')$$
and

\[\begin{aligned}
\displaystyle D^{G'}(y)^{1/2}c_{\theta'}(y) & =\lvert W(G'_y,T_y)\rvert^{-1} \lim\limits_{y'\in T_y(F)\to y} D^{G'}(y')^{1/2}\theta'(y') \\
 & =\lvert W(G'_y,T_y)\rvert^{-1}\lim\limits_{x'\in T_x(F)\to x} D^{G'}(\iota(x'))^{1/2}\theta'(\iota(x'))
\end{aligned}\]

\noindent For all $x'\in T_x(F)\cap G_{\reg}(F)$, the elements $x'$ and $\iota(x')$ are stably conjugate. Hence, since $\theta'$ is a transfer of $\theta$, we have $\theta'(\iota(x'))=\theta(x')$ for all $x'\in T_x(F)\cap G_{\reg}(F)$. On the other hand, we also have $D^G(x)=D^{G'}(y)$, $\lvert W(G_x,T_x)\rvert=\lvert W(G'_y,T_y)\rvert$ and $D^{G'}(\iota(x'))=D^G(x')$ for all $x'\in T_x(F)\cap G_{\reg}(F)$. Consequently, the two formulas above imply the equality $c_{\theta'}(y)=c_\theta(x)$.

\vspace{2mm}

\noindent Let us assume henceforth that $G$ is quasi-split. Following Kottwitz \cite{Kott2}, we may associate to any class of pure inner forms $\alpha\in H^1(F,G)$ a sign $\gls{eGalpha}$ (as the notation suggests, this sign actually only depends on the isomorphism class of the group $G_\alpha$). Let $\gls{Br2F}=H^2(F,\{\pm 1\})=\{\pm 1\}$ be the $2$-torsion subgroup of the Bauer group of $F$. The sign $e(G_\alpha)$ will more naturally be an element of $Br_2(F)$. To define it, we need to introduce a canonical algebraic central extension

\begin{align}\label{eq 12.1.3}
\displaystyle 1\to\{\pm 1\}\to \gls{Gtilde}\to G\to 1
\end{align}

\noindent of $G$ by $\{\pm 1\}$. Recall that a quasi-split connected group over $F$ is classified up to conjugation by its (canonical) based root datum $\Psi_0(G)=\left(X_G,\Delta_G,X_G^\vee,\Delta_G^\vee\right)$ together with the natural action of $\Gamma_F$ on $\Psi_0(G)$. For any Borel pair $(B,T)$ of $G$ that is defined over $F$, we have a canonical $\Gamma_F$-equivariant isomorphism $\Psi_0(G)\simeq \left(X^*(T),\Delta(T,B),X_*(T),\Delta(T,B)^\vee\right)$ where $\Delta(T,B)\subseteq X^*(T)$ denotes the set of simple roots of $T$ in $B$ and $\Delta(T,B)^\vee\subseteq X_*(T)$ denotes the corresponding sets of simple coroots. Fix such a Borel pair and set
$$\displaystyle \rho=\frac{1}{2}\sum_{\beta\in R(G,T)}\beta\in X^*(T)\otimes \mathbb{Q}$$
for the half sum of the roots of $T$ in $B$. The image of $\rho$ in $X_G\otimes \mathbb{Q}$ doesn't depend on the particular Borel pair $(B,T)$ chosen and we shall still denote by $\rho$ this image. Consider now the following based root datum

\begin{align}\label{eq 12.1.4}
(\widetilde{X}_G,\Delta_G,\widetilde{X}_G^\vee,\Delta_G^\vee)
\end{align}

\noindent where $\widetilde{X}_G=X_G+\mathbb{Z}\rho\subseteq X_G\otimes \mathbb{Q}$ and $\widetilde{X}_G=\{\lambda^\vee\in X_G^\vee;\; \langle \lambda^\vee,\rho\rangle\in \mathbb{Z}\}$. Note that we have $\Delta_G^\vee\subseteq \widetilde{X}_G^\vee$ since $\langle \alpha^\vee,\rho\rangle=1$ for all $\alpha^\vee\in \Delta_G^\vee$. The based root datum \ref{eq 12.1.4}, with its natural $\Gamma_F$-action, is the based root datum of a unique quasi-split group $\widetilde{G}_0$ over $F$ well-defined up to conjugacy. Moreover, we have a natural central isogeny $\widetilde{G}_0\to G$, well-defined up to $G(F)$-conjugacy, whose kernel is either trivial or $\{\pm 1\}$ (depending on whether $\rho$ belongs to $X_G$ or not). If the kernel is $\{\pm 1\}$, we set $\widetilde{G}=\widetilde{G}_0$ otherwise we simply set $\widetilde{G}=G\times \{\pm 1\}$. In any case, we obtain a short exact sequence like \ref{eq 12.1.3} well-defined up to $G(F)$-conjugacy. The last term of the long exact sequence associated to \ref{eq 12.1.3} yields a canonical map

\begin{align}\label{eq 12.1.5}
H^1(F,G)\to H^2(F,\{\pm 1\})=Br_2(F)\simeq \{\pm 1\}
\end{align}

\noindent We now define the sign $e(G_\alpha)$, for $\alpha\in H^1(F,G)$, simply to be the image of $\alpha$ by this map. We will need the following fact

\vspace{3mm}

\begin{num}
\item\label{eq 12.1.6} Let $T$ be a (not necessarily maximal) subtorus of $G$. Then, the composite of \ref{eq 12.1.5} with the natural map $H^1(F,T)\to H^1(F,G)$ is a group morphism $H^1(F,T)\to Br_2(F)$. Moreover, if $T$ is anisotropic this morphism is onto if and only if the inverse image $\widetilde{T}$ of $T$ in $\widetilde{G}$ is a torus (i.e., is connected).
\end{num}

\vspace{3mm}

\noindent The first part is obvious since the map $H^1(F,T)\to Br_2(F)=H^2(F,\{\pm 1\})$ is a connecting map of the long exact sequence associated to the short exact sequence

\begin{align}\label{eq 12.1.7}
1\to \{\pm 1\}\to \widetilde{T}\to T\to 1
\end{align}

\noindent of algebraic abelian groups. To see why the second part of the claim is true, we first note that if $\widetilde{T}$ is not a torus then the short exact sequence \ref{eq 12.1.7} splits so that the connecting map $H^1(F,T)\to H^2(F,\{\pm 1\})$ is trivial. On the other hand, if $\widetilde{T}$ is a torus and $T$ is anisotropic then $\widetilde{T}$ is also anisotropic. From the long exact sequence associated to \ref{eq 12.1.7}, we can extract the following short exact sequence
$$H^1(F,T)\to H^2(F,\{\pm 1\})\to H^2(F,\widetilde{T})$$
But, by Tate-Nakayama we have $H^2(F,\widetilde{T})=0$ and it immediately follows that the morphism $H^1(F,T)\to H^2(F,\{\pm 1\})$ is onto.

\subsection{Pure inner forms of a GGP triple}\label{section 12.2}

\noindent Let $V$ be a hermitian space. We have the following explicit description of the pure inner forms of $U(V)$. The cohomology set $H^1(F,U(V))$ naturally classifies the isomorphism classes of hermitian spaces of the same dimension as $V$. Let $\alpha\in H^1(F,U(V))$ and choose a hermitian space $V_\alpha$ in the isomorphism class corresponding to $\alpha$. Set $V_{\overline{F}}=V\otimes_F \overline{F}$ and $V_{\alpha,\overline{F}}=V_\alpha\otimes_F \overline{F}$. Fix an isomorphism $\phi_\alpha: V_{\overline{F}}\simeq V_{\alpha,\overline{F}}$ of $\overline{E}$-hermitian spaces. Then, the triple $(U(V_\alpha),\psi_\alpha,c_\alpha)$, where $\psi_\alpha$ is the isomorphism $U(V)_{\overline{F}}\simeq U(V_\alpha)_{\overline{F}}$ given by $\psi_\alpha(g)=\phi_\alpha\circ g\circ \phi_\alpha^{-1}$ and $c_\alpha$ is the $1$-cocycle given by $\sigma\in \Gamma_F\mapsto \phi_\alpha^{-1}{}^\sigma \phi_\alpha$, is a pure inner form of $U(V)$ in the class of $\alpha$. Moreover, the $2$-cover $\gls{U(V)tilde}$ of $U(V)$ that has been defined for general reductive groups at the end of the previous section admits the following explicit description:

\vspace{2mm}

\begin{itemize}
\item If $\dim(V)$ is odd, then $\widetilde{U(V)}=U(V)\times \{\pm 1\}$;

\item If $\dim(V)$ is even then $\widetilde{U(V)}=\{(g,z)\in U(V)\times \Ker N_{E/F}; \; \det(g)=z^2\}$.
\end{itemize}

\vspace{2mm}

\noindent We now return to the GGP triple $(G,H,\xi)$ that we have fixed. Recall that this GGP triple comes from an admissible pair $(V,W)$ of hermitian spaces and that we are assuming in this chapter that the groups $G$ and $H$ are quasi-split. Let $\alpha\in H^1(F,H)$. We are going to associate to $\alpha$ a new GGP triple $\gls{GalphaHalphaxialpha}$ well-defined up to conjugacy. Since $H^1(F,H)=H^1(F,U(W))$, to the cohomology class $\alpha$ corresponds an isomorphism class of hermitian spaces of the same dimension as $W$. Let $\gls{Walpha}$ be a hermitian space in this isomorphism class and set $\gls{Valpha}=W_\alpha\oplus^\perp Z$ (recall that $Z$ is the orthogonal complement of $W$ in $V$). Then, the pair $(V_\alpha,W_\alpha)$ is easily seen to be admissible and hence there is a GGP triple $(G_\alpha,H_\alpha,\xi_\alpha)$ associated to it. Of course, this GGP triple is well-defined up to conjugacy. We call such a GGP triple a {\em pure inner form} of $(G,H,\xi)$. By definition, these pure inner forms are parametrized by $H^1(F,H)$. Note that for all $\alpha\in H^1(F,H)$, $G_\alpha$ is a pure inner form of $G$ in the class corresponding to the image of $\alpha$ in $H^1(F,G)$ and that the natural map $H^1(F,H)\to H^1(F,G)$ is injective.

\subsection{The local Langlands correspondence}\label{section 12.3}

\noindent In this section, we recall the local Langlands correspondence in a form that will be suitable for us. Let $G$ be a quasi-split connected reductive group over $F$ and denotes by $\gls{LG}=\widehat{G}(\mathbb{C})\rtimes W_F$ its Langlands dual, where $\gls{W_F}$ denotes the Weil group of $F$. Recall that a {\em Langlands parameter} for $G$ is a homomorphism from the group

$$\displaystyle \gls{L_F}=\left\{
    \begin{array}{ll}
        W_F\times SL_2(\mathbb{C}) & \mbox{ if } F \mbox{ is } p-\mbox{adic} \\
        W_F & \mbox{ if } F=\mathbb{R}
    \end{array}
\right.
$$

\noindent to ${}^L G$ satisfying the usual conditions of continuity, semi-simplicity, algebraicity and compatibility with the projection ${}^L G\to W_F$. A Langlands parameter $\varphi$ is said to be {\em tempered} if $\varphi(W_F)$ is bounded. By the hypothetical local Langlands correspondence, a tempered Langlands parameter $\varphi$ for $G$ should give rise to a finite set $\gls{PiGphi}$, called a {\em $L$-packet}, of (isomorphism classes of) tempered representations of $G(F)$. Actually, such a parameter $\varphi$ should also give rise to tempered $L$-packets $\gls{PiGalphaphi}\subseteq \Temp(G_\alpha)$ for all $\alpha\in H^1(F,G)$ (we warn the reader that in this formulation of the local Langlands correspondence, some of the $L$-packets $\Pi^{G_\alpha}(\varphi)$ may be empty). These families of $L$-packets should of course satisfy some conditions. Among them, we expect the following properties to hold for every tempered Langlands parameter $\varphi$ of $G$:

\vspace{5mm}

\hspace{3mm}(STAB) For all $\alpha\in H^1(F,G)$, the character

$$\displaystyle \gls{thetaalphaphi}=\sum_{\pi\in \Pi^{G_\alpha}(\varphi)}\theta_\pi$$

\hspace{18mm} is stable

\vspace{3mm}

\noindent see \S \ref{section 12.1} for the definition of stable; also notice that our different notion of ``strongly stable conjugate" does not affect this property since it only involves the values of $\theta_{\alpha,\varphi}$ at regular semi-simple elements (for which the two notions of stable conjugacy coincide). For $\alpha=1\in H^1(F,G)$, in which case $G_\alpha=G$, we shall simply set $\gls{thetaphi}=\theta_{1,\varphi}$.

\vspace{2mm}

\hspace{3mm}(TRANS) For all $\alpha\in H^1(F,G)$, the stable character $\theta_{\alpha,\varphi}$ is the transfer of $e(G_\alpha)\theta_{\varphi}$

\vspace{0.4pt}

\hspace{20mm} where $e(G_\alpha)$ is the Kottwitz sign whose definition has been recalled in

\vspace{0.4pt}

\hspace{20mm} Section \ref{section 12.1}.

\vspace{5mm}

\hspace{3mm}(WHITT) For every $\mathcal{O}\in \Nil_{\reg}(\mathfrak{g})$, there exists exactly one representation in the 

\vspace{0.4pt}

\hspace{20mm} $L$-packet $\Pi^G(\varphi)$ admitting a Whittaker model of type $\mathcal{O}$.

\vspace{3mm}

\noindent (cf.\ Section \ref{section 4.8} for the bijection between $\Nil_{\reg}(\mathfrak{g})$ and the set of types of Whittaker models). 

\vspace{2mm}

\noindent Notice that these conditions are far from characterizing the compositions of the $L$-packets uniquely. However, by the linear independence of characters, conditions (STAB) and (TRANS) uniquely characterize the $L$-packets $\Pi^{G_\alpha}(\varphi)$, $\alpha\in H^1(F,G)$, in terms of $\Pi^G(\varphi)$.

\vspace{2mm}

\noindent When $F=\mathbb{R}$, the local Langlands correspondence has been constructed by Langlands himself \cite{Lan} building on previous results of Harish-Chandra. This correspondence indeed satisfies the three conditions stated above. That (STAB) and (TRANS) hold is a consequence of early work of Shelstad (\cite{She} Lemma 5.2 and Theorem 6.3). The property (WHITT) for its part, follows from results of Kostant (\cite{Kost} Theorem 6.7.2) and Vogan (\cite{Vog} Theorem 6.2).

\noindent When $F$ is $p$-adic, the local Langlands correspondence is known in a variety of cases. In particular, for unitary groups, which are our main concern, the existence of the Langlands correspondence is now fully established thanks to Mok \cite{Mok} and Kaletha-Minguez-Shin-White \cite{KMSW} both building up on previous work of Arthur who dealt with orthogonal and symplectic groups \cite{A7}. That the tempered $L$-packets constructed in these references verify the conditions (STAB) and (TRANS) follows from \cite{Mok} Theorem 3.2.1 (a) and \cite{KMSW} Proposition 1.5.2. Moreover, the $L$-packets on the quasi-split form $G$ satisfy condition (WHITT) by \cite{Mok} Corollary 9.2.4.

\subsection{The theorem}\label{section 12.4}

\noindent Recall that we fixed a GGP triple $(G,H,\xi)$ with the requirement that $G$ and $H$ be quasi-split. Also, we have defined in Section \ref{section 12.2} the pure inner forms $(G_\alpha,H_\alpha,\xi_\alpha)$ of $(G,H,\xi)$. These are also GGP triples, they are parametrized by $H^1(F,H)$ and $G_\alpha$ is a pure inner form of $G$ corresponding to the image of $\alpha$ in $H^1(F,G)$ via the natural map $H^1(F,H)\to H^1(F,G)$.

\vspace{2mm}

\noindent As we said, the local Langlands correspondence, in the form we stated it, is known for unitary groups. It is a fortiori known for the product of two such groups and hence for $G$.

\vspace{2mm}

\noindent The purpose of this chapter is to show the following theorem. It has already been shown in \cite{Beu1} (th\'eor\`eme 18.4.1) in the $p$-adic case. The proof we present here is essentially the same as the one given in \cite{Beu1} (which itself follows closely the proof of th\'eor\`eme 13.3 of \cite{Wa1})  but here we are treating both the $p$-adic and the real case at the same time and it requires more care since in the real case there are usually far more pure inner forms to keep track of (in the $p$-adic case there are only two unless $\dim(W)=0$).

\begin{theo}\label{theorem 12.4.1}
Let $\varphi$ be a tempered Langlands parameter for $G$. Then, there exists a unique representation $\pi$ in the disjoint union of $L$-packets

$$\displaystyle \bigsqcup_{\alpha\in H^1(F,H)} \Pi^{G_\alpha}(\varphi)$$

\noindent such that $m(\pi)=1$.
\end{theo}

\subsection{Stable conjugacy classes inside $\Gamma(G,H)$}\label{section 12.5}

\noindent Recall that in Section \ref{section 11.0}, we have defined a set $\Gamma(G,H)$ of semi-simple conjugacy classes in $G(F)$. It consists in the $G(F)$-conjugacy classes of elements $x\in U(W)_{\ssi}(F)$ such that

$$\displaystyle \gls{T_x}:=U(W''_x)_x$$

\noindent is an anisotropic torus (where we recall that $W''_x$ denotes the image of $x-1$ in $W$). Two elements $x,x'\in U(W)_{\ssi}(F)$ are $G(F)$-conjugate if and only if they are $U(W)(F)$-conjugate and moreover if it is so, any element $g\in U(W)(F)$ conjugating $x$ to $x'$ induces an isomorphism

$$\displaystyle U(W''_x)_x\simeq U(W''_{x'})_{x'}$$

\noindent Moreover, this isomorphism depends on the choice of $g$ only up to an inner automorphism. From this it follows that any conjugacy class $x\in \Gamma(G,H)$ determines the anisotropic torus $T_x$ up to a unique isomorphism so that we can speak of ``the torus" $T_x$ associated to $x$.

\vspace{2mm}

\noindent These considerations of course apply verbatim to the pure inner forms $(G_\alpha,H_\alpha,\xi_\alpha)$, $\alpha\in H^1(F,H)$, of the GGP triple $(G,H,\xi)$ that were introduced in Section \ref{section 12.2}. In particular, for all $\alpha\in H^1(F,H)$, we have a set $\Gamma(G_\alpha,H_\alpha)$ of semi-simple conjugacy classes in $G_\alpha(F)$ and to any $y\in \Gamma(G_\alpha,H_\alpha)$ is associated an anisotropic torus $T_y$.

\vspace{2mm}

\begin{prop}\label{proposition 12.5.1}
\begin{enumerate}[(i)]
\item Let $\alpha\in H^1(F,H)$ and $y\in \Gamma(G_\alpha,H_\alpha)$ be such that $G_{\alpha,y}$ is quasi-split. Then, the set

$$\{x\in \Gamma(G,H);\; x\sim_{\stab} y\}$$

\noindent is non-empty.

\item Let $\alpha\in H^1(F,H)$, $x\in \Gamma(G,H)$ and $y\in \Gamma(G_\alpha,H_\alpha)$ be such that $x\sim_{\stab} y$. Choose $g\in G_\alpha(\overline{F})$ such that $g\psi_\alpha(x)g^{-1}=y$ and $\Ad(g)\circ \psi_\alpha:G_x\simeq G_y$ is defined over $F$. Then, $\Ad(g)\circ \psi_\alpha$ restricts to an isomorphism

$$T_x\simeq T_y$$

\noindent that is independent of the choice of $g$.

\item Let $x\in \Gamma(G,H)$. Then, for all $\alpha\in H^1(F,H)$ there exists a natural bijection between the set

$$\{y\in \Gamma(G_\alpha,H_\alpha);\; x\sim_{\stab} y\}$$

\noindent and the set

$$q_x^{-1}(\alpha)$$

\noindent where $q_x$ denotes the natural map $H^1(F,T_x)\to H^1(F,G)$.

\item Let $x\in \Gamma(G,H)$, $x\neq 1$. Then, the composition of the map $\alpha\in H^1(F,G)\mapsto e(G_\alpha)\in Br_2(F)$ with the natural map $H^1(F,T_x)\to H^1(F,G)$ gives a surjective morphism of groups $H^1(F,T_x)\to Br_2(F)$.
\end{enumerate}
\end{prop}

\vspace{3mm}

\noindent\ul{Proof}: For all $\alpha\in H^1(F,H)=H^1(F,U(W))$, let us fix a hermitian space $W_\alpha$ of the same dimension as $W$ and in the isomorphism class defined by $\alpha$. Set $V_\alpha=W_\alpha\oplus^\perp Z$. Then, we may assume that $G_\alpha=U(W_\alpha)\times U(V_\alpha)$. Moreover, we may also fix the other parts of the data $(G_\alpha,\psi_\alpha,c_\alpha)$ of a pure inner form of $G$ in the class of $\alpha$ as follows. Choose an isomorphism $\phi_\alpha^W:W_{\overline{F}}\simeq W_{\alpha,\overline{F}}$ of $\overline{E}$-hermitian spaces, where we have set as usual $W_{\overline{F}}=W\otimes_F \overline{F}$ and $W_{\alpha,\overline{F}}=W_\alpha\otimes_F \overline{F}$, and extend it to an isomorphism $\phi_\alpha^V:V_{\overline{F}}\simeq V_{\alpha,\overline{F}}$ that is the identity on $Z_{\overline{F}}$. Then, we may take $\psi_\alpha$ to be the isomorphism

$$G_{\overline{F}}=U(W_{\overline{F}})\times U(V_{\overline{F}})\simeq G_{\alpha,\overline{F}}= U(W_{\alpha,\overline{F}})\times U(V_{\alpha,\overline{F}})$$

\noindent given by

$$(g_W,g_V)\mapsto \left(\phi_\alpha^W\circ g_W\circ (\phi_\alpha^W)^{-1}, \phi_\alpha^V\circ g_V\circ (\phi_\alpha^V)^{-1}\right)$$

\noindent and we may take the $1$-cocycle $c_\alpha$ to be given by

$$\sigma\in \Gamma_F\mapsto \left((\phi_\alpha^W)^{-1}{}^\sigma \phi_\alpha^W, (\phi^V_\alpha)^{-1}{}^\sigma \phi_\alpha^V\right)\in U(W_{\overline{F}})\times U(V_{\overline{F}})=G_{\overline{F}}$$

\noindent Notice that if we do so, then the isomorphism $\psi_\alpha$ sends $H_{\overline{F}}$ onto $H_{\alpha,\overline{F}}$.

\begin{enumerate}[(i)]
\item Let $y\in \Gamma(G_\alpha,H_\alpha)$ and assume that $G_{\alpha,y}$ is quasi-split. Identify $y$ with one of its representatives in $U(W_\alpha)$. Writing the decomposition \ref{eq 11.1.1} for $G_{\alpha,y}$, we have

$$G_{\alpha,y}=G'_{\alpha,y}\times G''_{\alpha,y}$$

\noindent where

$$G'_{\alpha,y}=U(W'_\alpha)\times U(V'_\alpha),\;\; G''_{\alpha,y}=U(W''_\alpha)_y\times U(W''_\alpha)_y$$

\noindent where $W'_\alpha$ is the kernel of $y-1$ in $W_\alpha$, $V'_\alpha=W'_\alpha\oplus Z$ and $W''_\alpha$ is the image of $y-1$ in $W_\alpha$. By definition of $\Gamma(G_\alpha,H_\alpha)$, $T_y=U(W''_\alpha)_y$ is an anisotropic torus and in particular $y$ is a regular element of $U(W''_\alpha)$ without the eigenvalue $1$. Since $G_{\alpha,y}$ is quasi-split, we see that the unitary groups of both $W'_\alpha$ and $V'_\alpha=W'_\alpha\oplus Z$ are quasi-split. As $Z$ is odd-dimensional, by Witt's theorem these two conditions are easily seen to imply that the hermitian space $W'_\alpha$ embeds in $W$. Fix such an embedding $W'_\alpha\hookrightarrow W$ and let us denote by $W''$ the orthogonal complement of $W'_\alpha$ in $W$. Then, $U(W''_\alpha)$ is a pure inner form of $U(W'')$ and moreover $U(W'')$ is quasi-split. Hence, by \ref{eq 12.1.1}, there exists a regular element $x\in U(W'')(F)$ that is stably conjugate to $y$. In particular, $x$ doesn't have the eigenvalue $1$ when acting on $W''$ and moreover we have an isomorphism $U(W'')_x\simeq U(W''_\alpha)_y$ which is defined over $F$. Thus, $U(W'')_x$ is an anisotropic torus and it follows that $x\in \Gamma(G,H)$. This ends the proof of (i).

\item Let us identify $x$ and $y$ with representatives in $U(W)(F)$ and $U(W_\alpha)(F)$ respectively. Then, $\psi_\alpha(x)$ and $y$ are $G_\alpha(\overline{F})$-conjugate. However, as $G_\alpha=U(W_\alpha)\times U(V_\alpha)$, two elements in $U(W_\alpha)(\overline{F})$ are $G_\alpha(\overline{F})$-conjugate if and only if they are $U(W_\alpha)(\overline{F})$-conjugate. Hence, there exists $h\in U(W_\alpha)(\overline{F})$ such that $y=h\psi_\alpha(x)h^{-1}$. Obviously, $h\circ \phi_\alpha$ sends $W''_x$ to $W''_{\alpha,y}$ and so $\Ad(h)\circ \psi_\alpha$ induces an isomorphism $T_x=U(W''_x)_x\simeq U(W''_{\alpha,y})_y=T_y$. Now every element $g\in G_\alpha(\overline{F})$ such that $g\psi_\alpha(x)g^{-1}=y$ may be written $g=g_yh$ for some element $g_y\in G_{\alpha,y}(\overline{F})$ (since $G_\alpha$ has a derived subgroup which is simply-connected, here we have $G_{\alpha,y}=Z_{G_\alpha}(y)$) and as $T_y$ is contained in the center of $G_{\alpha,y}$, the restriction of $\Ad(g)\circ \psi_\alpha$ to $T_x$ will coincide with the isomorphism $\Ad(h)\circ \psi_\alpha:T_x\simeq T_y$.

\item Recall that the set

$$\displaystyle \{y\in \Gamma(G_\alpha);\; x\sim_{\stab} y\}$$

\noindent is naturally in bijection with $\Ima\left(H^1(F,Z(G_x))\to H^1(F,G_x)\right)\cap p_x^{-1}(\alpha)$ where $p_x$ denotes the natural map $H^1(F,G_x)\to H^1(F,G)$ (since the derived subgroup of $G_\alpha$ is simply connected we have $Z_G(x)=G_x$). This bijection is given as follows: for $y\in \Gamma(G_\alpha)$ such that $x\sim_{\stab}y$ choose $g\in G(\overline{F})$ so that $\psi_\alpha(gxg^{-1})=y$. Then, the cohomological class associated to $y$ is the image of the $1$-cocycle $\sigma\in \Gamma_F\mapsto g^{-1}c_{\alpha,\sigma}\sigma(g)$ in $H^1(F,G_x)$ which by definition belongs to $\Ima\left(H^1(F,Z(G_x))\to H^1(F,G_x)\right)\cap p_x^{-1}(\alpha)$. To obtain the desired bijection, it suffices to show that the image of the subset

$$\displaystyle \{y\in \Gamma(G_\alpha,H_\alpha);\; x\sim_{\stab} y\}$$

\noindent by this bijection is exactly $H^1(F,T_x)\cap p_x^{-1}(\alpha)$ (notice that the restriction of the map $H^1(F,Z(G_x))\to H^1(F,G_x)$ to $H^1(F,T_x)$, where the embedding $T_x\hookrightarrow Z(G_x)$ being induced from $H\hookrightarrow G$ is the `diagonal' one, is injective as the latter is a direct summand of $H^1(F,G_x)$). Let $y\in \Gamma(G_\alpha,H_\alpha)$ be such that $x\sim_{\stab} y$ and identify $y$ with one of its representatives in $H_\alpha(F)$. Since two semi-simple elements of $H_\alpha(\overline{F})$ are $G_\alpha(\overline{F})$-conjugate if and only if they are $H_\alpha(\overline{F})$-conjugate and as $\psi_\alpha$ sends $H(\overline{F})$ to $H_\alpha(\overline{F})$, we may find $h\in H(\overline{F})$ such that $y=\psi_\alpha(hxh^{-1})$ and it follows that the cohomology class associated to $y$ lies in

\begin{align}\label{eq 12.5.1}
\displaystyle \Ima\left(H^1(F,H_x)\to H^1(F,G_x)\right)\cap \Ima\left(H^1(F,Z(G_x))\to H^1(F,G_x)\right)\cap p_x^{-1}(\alpha)
\end{align}

\noindent Conversely, assume that $y\in \Gamma(G_\alpha)$ is strongly stably conjugate to $x$ and that its associated cohomology class belongs to \ref{eq 12.5.1}. Then, since the map $H^1(F,H)\to H^1(F,G)$ is injective, we may find $h\in H(\overline{F})$ such that $\psi_\alpha(hxh^{-1})\in H_\alpha(F)$ is in the conjugacy class of $y$. From this, it is easy to infer that $y\in \Gamma(G_\alpha,H_\alpha)$. Thus, we have proved that the set
$$\displaystyle \{y\in \Gamma(G_\alpha,H_\alpha);\; x\sim_{\stab} y\}$$
is in bijection with \ref{eq 12.5.1}. To conclude, we only need to prove the equality

\begin{align}\label{eq 12.5.2}
\displaystyle \Ima\left(H^1(F,H_x)\to H^1(F,G_x)\right)\cap \Ima\left(H^1(F,Z(G_x))\to H^1(F,G_x)\right)=H^1(F,T_x)
\end{align}

\noindent Using the decompositions
$$G_x=G'_x\times G''_x,\;\;\; H_x=H'_x\times T_x$$
$$G''_x=T_x\times T_x$$
it is easy to deduce that

\[\begin{aligned}
\displaystyle & \Ima\left(H^1(F,H_x)\to H^1(F,G_x)\right)\cap \Ima\left(H^1(F,Z(G_x))\to H^1(F,G_x)\right)= \\
 & \left[\Ima\left(H^1(F,Z(G'_x))\to H^1(F,G'_x)\right)\cap \Ima\left(H^1(F,H'_x)\to H^1(F,G'_x)\right)\right]\times H^1(F,T_x)
\end{aligned}\]

\noindent Hence, to get \ref{eq 12.5.2} it suffices to see that

\begin{align}\label{eq 12.5.3}
\displaystyle \Ima\left(H^1(F,Z(G'_x))\to H^1(F,G'_x)\right)\cap \Ima\left(H^1(F,H'_x)\to H^1(F,G'_x)\right)=\{1\}
\end{align}

\noindent Now recall that (see Section \ref{section 11.0})
$$G'_x=U(W'_x)\times U(V'_x)$$
$$H'_x=U(W'_x)\ltimes N_x$$
where $W'_x$ is the kernel of $x-1$ in $W$, $V'_x=W'_x\oplus Z$ and $N_x$ is the centralizer of $x$ in $N$. Thus, we have
$$H^1(F,G'_x)=H^1(F,U(W'_x))\times H^1(F,U(V'_x))$$
$$H^1(F,H'_x)=H^1(F,U(W'_x))$$
and the map $H^1(F,H'_x)\to H^1(F,G'_x)$ is the product of the two maps
$$H^1(F,U(W'_x))\to H^1(F,U(W'_x)),\;\;\; H^1(F,U(W'_x))\to H^1(F,U(V'_x))$$
which are both injective. So, to get \ref{eq 12.5.3} we only need to show that

\begin{align}\label{eq 12.5.4}
\displaystyle \Ima\left(H^1\left(F,Z(U(W'_x))\right)\to H^1(F,U(V'_x))\right)\cap \Ima\left(H^1\left(F,Z(U(V'_x))\right)\to H^1(F,U(V'_x))\right)=\{1\}
\end{align}

\noindent Recall that $H^1(F,U(V'_x))$ classifies the (isomorphism classes of) hermitian spaces of the same dimension as $V'_x$. Let $\delta\in F^\times\smallsetminus N_{E/F}(E^\times)$. Then, the group $H^1\left(F,Z(U(W'_x))\right)$ contains only one nontrivial element whose image in $H^1(F,U(V'_x))$ corresponds to the hermitian space $\delta W'_x\oplus Z$ (where $\delta W'_x$ denotes the hermitian space obtained from $W'_x$ by multiplying its hermitian form by $\delta$). Similarly, $H^1\left(F,Z(U(V'_x))\right)$ has only one nontrivial element whose image in $H^1(F,U(V'_x))$ corresponds to the hermitian space $\delta V'_x$. As $Z$ is odd dimensional, we have $\delta Z\not\simeq Z$ (the two hermitian spaces have distinct discriminants) so that by Witt's theorem we also have $\delta V'_x=\delta W'_x\oplus \delta Z\not\simeq \delta W'_x\oplus Z$. This proves \ref{eq 12.5.4} and ends the proof of (iii).

\item Let us denote by $\widetilde{G}$ the $2$-cover of $G$ defined at the end of Section \ref{section 12.1} and let $\widetilde{T}_x$ be the inverse image of $T_x$ in this $2$-cover. Then, by \ref{eq 12.1.6}, it suffices to check that $\widetilde{T}_x$ is connected. By the precise description of $\widetilde{U(V)}$ and $\widetilde{U(W)}$ given at the beginning of Section \ref{section 12.2} and since exactly one of the hermitian spaces $V$ and $W$ is even-dimensional, we have
$$\widetilde{T}_x=\{(t,z)\in T_x\times \Ker N_{E/F}; \det(t)=z^2\}$$
Thus, we need to show that the determinant has no square-root in the character group of $T_x$ but this is obvious since over the algebraic closure there exists an isomorphism
$$(T_x)_{\overline{F}}\simeq \mathbb{G}_m^\ell$$
$$t\mapsto (t_i)_{1\leqslant i\leqslant \ell}$$
for some integer $\ell\geqslant 1$ such that
$$\displaystyle \det(t)=\prod_{i=1}^\ell t_i$$
for all $t\in T_x(\overline{F})$. $\blacksquare$
\end{enumerate}

\subsection{Proof of Theorem \ref{theorem 12.4.1}}\label{section 12.6}

\noindent Let $\varphi$ be a tempered Langlands parameter for $G$. We want to show that the sum

\begin{align}\label{eq 12.6.1}
\displaystyle \sum_{\alpha\in H^1(F,H)} \sum_{\pi\in \Pi^{G_\alpha}(\varphi)}m(\pi)
\end{align}

\noindent is equal to $1$. Let $\alpha\in H^1(F,H)$. By Theorem \ref{theorem 11.2.2}, we have
$$\displaystyle m(\pi)=\lim\limits_{s\to 0^+} \int_{\Gamma(G_\alpha,H_\alpha)}c_\pi(y) D^{G_\alpha}(y)^{1/2} \Delta(y)^{s-1/2}dy$$
for all $\pi\in \Pi^{G_\alpha}(\varphi)$. Summing this equality over the $L$-packet $\Pi^{G_\alpha}(\varphi)$, we deduce that

\begin{align}\label{eq 12.6.2}
\displaystyle \sum_{\pi\in \Pi^{G_\alpha}(\varphi)}m(\pi)=\lim\limits_{s\to 0^+} \int_{\Gamma(G_\alpha,H_\alpha)} c_{\varphi,\alpha}(y) D^{G_\alpha}(y)^{1/2}\Delta(y)^{s-1/2}dy
\end{align}

\noindent where we have set
$$\displaystyle c_{\varphi,\alpha}=\sum_{\pi\in \Pi^{G_\alpha}(\varphi)}c_\pi$$
Denote by $\Gamma_{\stab}(G_\alpha,H_\alpha)$ the set of strongly stable conjugacy classes in $\Gamma(G_\alpha,H_\alpha)$. We endow this set with the quotient topology and with the unique measure such that the projection map $\Gamma(G_\alpha,H_\alpha)\to \Gamma_{\stab}(G_\alpha,H_\alpha)$ is locally measure-preserving. By the condition (STAB) of Section \ref{section 12.3}, the character $\displaystyle \theta_{\varphi,\alpha}=\sum_{\pi\in \Pi^{G_\alpha}(\varphi)}\theta_\pi$ is stable and it follows from \ref{eq 12.1.2} that the function $y\mapsto c_{\varphi,\alpha}(y)$ is constant on the fibers of the map $\Gamma(G_\alpha,H_\alpha)\to \Gamma_{\stab}(G_\alpha,H_\alpha)$. Since it is also trivially true for the functions $D^{G_\alpha}$ and $\Delta$, we may rewrite the right hand side of \ref{eq 12.6.2} as an integral over $\Gamma_{\stab}(G_\alpha,H_\alpha)$ to obtain

\begin{align}\label{eq 12.6.3}
\displaystyle \sum_{\pi\in \Pi^{G_\alpha}(\varphi)} m(\pi)=\lim\limits_{s\to 0^+} \int_{\Gamma_{\stab}(G_\alpha,H_\alpha)} \lvert p_{\alpha,\stab}^{-1}(y)\rvert c_{\varphi,\alpha}(y) D^{G_\alpha}(y)^{1/2} \Delta(y)^{s-1/2}dy
\end{align}

\noindent where $p_{\alpha,\stab}:\Gamma(G_\alpha,H_\alpha)\to \Gamma_{\stab}(G_\alpha,H_\alpha)$ is the natural projection.

\vspace{2mm}

\noindent Let us introduce the subset $\Gamma_{\stab}^{\quasid}(G_\alpha,H_\alpha)\subseteq \Gamma_{\stab}(G_\alpha,H_\alpha)$ of elements $x\in \Gamma_{\stab}(G_\alpha,H_\alpha)$ such that $G_{\alpha,x}$ is quasi-split. This subset is open and closed in $\Gamma_{\stab}(G_\alpha,H_\alpha)$ and by Proposition \ref{proposition 4.5.1} 1.(i), the function $c_{\varphi,\alpha}$ vanishes on the complement $\Gamma_{\stab}(G_\alpha,H_\alpha)\smallsetminus \Gamma_{\stab}^{\quasid}(G_\alpha,H_\alpha)$ so that \ref{eq 12.6.3} becomes

\begin{align}\label{eq 12.6.4}
\displaystyle \sum_{\pi\in \Pi^{G_\alpha}(\varphi)}m(\pi)=\lim\limits_{s\to 0^+} \int_{\Gamma_{\stab}^{\quasid}(G_\alpha,H_\alpha)} \lvert p_{\alpha,\stab}^{-1}(y)\rvert c_{\varphi,\alpha}(y) D^{G_\alpha}(y)^{1/2}\Delta(y)^{s-1/2}dy
\end{align}

\noindent By Proposition \ref{proposition 12.5.1}(i), we have an injection $\Gamma_{\stab}^{\quasid}(G_\alpha,H_\alpha)\hookrightarrow \Gamma_{\stab}^{\quasid}(G,H)$ such that if $y\mapsto x$ then $y$ and $x$ are strongly stably conjugate. For all $y\in \Gamma_{\stab}^{\quasid}(G_\alpha,H_\alpha)$, denoting by $x$ its image in $\Gamma_{\stab}^{\quasid}(G,H)$, we have the following commutative diagram

\begin{center}
\begin{tikzpicture}
\node (A) at (0,0) {$T_y(F)$};
\node (B) at (3,0) {$T_x(F)$};
\node (C) at (0,-2) {$\Gamma_{\stab}^{\quasid}(G_\alpha,H_\alpha)$};
\node (D) at (3,-2) {$\Gamma_{\stab}^{\quasid}(G,H)$};
\draw[transparent] (A) edge node[opacity=1] {\resizebox{1.5cm}{0.3cm}{$\simeq$}} (B);
\draw[->] (A) edge (C);
\draw[->] (B) edge (D);
\draw[right hook->] (C) edge (D);
\end{tikzpicture}
\end{center}

\noindent where the two vertical arrows are only defined in some neighborhood of $1$, given by $t\mapsto ty$, $t\mapsto tx$ and are both locally preserving measures (when $T_y(F)$ and $T_x(F)$ are both equipped with their unique Haar measure of total mass one) and the top vertical arrow is the restriction to the $F$-points of the isomorphism provided by Proposition \ref{proposition 12.5.1}(ii). From this diagram, we easily infer that the embedding $\Gamma_{\stab}^{\quasid}(G_\alpha,H_\alpha)\hookrightarrow \Gamma_{\stab}^{\quasid}(G,H)$ preserves measures. Moreover, by \ref{eq 12.1.2} and the condition (TRANS) of Section \ref{section 12.3}, if $y\in \Gamma_{\stab}^{\quasid}(G_\alpha,H_\alpha)$ maps to $x\in \Gamma_{\stab}^{\quasid}(G,H)$ we have the equality $c_{\varphi,\alpha}(y)=e(G_\alpha)c_\varphi(x)$ (where we have set $c_{\varphi}=c_{\varphi,1}$). Using these two facts, we may now express the right hand side of \ref{eq 12.6.4} as an integral over $\Gamma_{\stab}^{\quasid}(G,H)$. More precisely, the result is the following formula
$$\displaystyle \sum_{\pi\in \Pi^{G_\alpha}(\varphi)} m(\pi)=\lim\limits_{s\to 0^+} \int_{\Gamma_{\stab}^{\quasid}(G,H)} e(G_\alpha) n_\alpha(x) c_\varphi(x) D^G(x)^{1/2} \Delta(x)^{s-1/2} dx$$
where we have set
$$\displaystyle n_\alpha(x)=\left\lvert \{y\in \Gamma(G_\alpha,H_\alpha); y\sim_{\stab} x\}\right\rvert$$
for all $x\in \Gamma_{\stab}^{\quasid}(G,H)$. Summing the above equality over $\alpha\in H^1(F,H)$, we get

\begin{align}\label{eq 12.6.5}
\displaystyle \sum_{\alpha\in H^1(F,H)}\sum_{\pi\in \Pi^{G_\alpha}(\varphi)} m(\pi)=\lim\limits_{s\to 0^+} \int_{\Gamma_{\stab}^{\quasid}(G,H)} \sum_{\alpha\in H^1(F,H)} e(G_\alpha)n_\alpha(x) c_\varphi(x) D^G(x)^{1/2} \Delta(x)^{s-1/2} dx
\end{align}

\noindent Let $x\in \Gamma_{\stab}^{\quasid}(G,H)$ and consider the inner sum

\begin{align}\label{eq 12.6.6}
\displaystyle \sum_{\alpha\in H^1(F,H)} e(G_\alpha)n_\alpha(x)
\end{align}

\noindent By Proposition \ref{proposition 12.5.1}(iii), this sum equals
$$\displaystyle \sum_{\beta\in H^1(F,T_x)} e(G_\beta)$$
Moreover, by Proposition \ref{proposition 12.5.1}(iv), the map $\beta\in H^1(F,T_x)\mapsto e(G_\beta)\in Br_2(F)\simeq \{\pm 1\}$ is a group homomorphism which is non-trivial for $x\neq 1$. We deduce from this that the sum \ref{eq 12.6.6} is zero unless $x=1$. Returning to the formula \ref{eq 12.6.5} and taking into account these cancellations, we see that the right hand side of \ref{eq 12.6.5} reduces to the contribution of $1\in \Gamma(G,H)$ and thus we get

\begin{align}\label{eq 12.6.7}
\displaystyle \sum_{\alpha\in H^1(F,H)}\sum_{\pi\in \Pi^{G_\alpha}(\varphi)} m(\pi)=c_\varphi(1)
\end{align}

\noindent By Proposition \ref{proposition 4.8.1}(i), the term $c_\varphi(1)$ has the following representation-theoretic interpretation: it equals the number of representations in $\Pi^G(\varphi)$ admitting a Whittaker model, a representation being counted as many times as the number of different types of Whittaker model it has, divided by the number of different type of Whittaker models for $G$. By the condition (WHITT) of Section \ref{section 12.3}, this number is $1$. It follows that the left hand side of \ref{eq 12.6.7} is also equal to $1$ and we are done. $\blacksquare$

\bigskip

\appendix
\section{Topological vector spaces}\label{section A}

\noindent In this appendix, we collect some facts about topological vector spaces that will be use constantly throughout the paper. We will only consider Hausdorff locally convex topological vector spaces over $\mathbb{C}$. To abbreviate we will call them topological vector spaces. Let $E$ be a topological vector space. Recall that a subset $B\subseteq E$ is said to be {\em bounded} if for every neighborhood $U$ of the origin there exists $\lambda>0$ such that $B\subseteq \lambda U$. An equivalent condition is that every continuous semi-norm on $E$ is bounded on $B$. By the Banach-Steinhaus theorem, a subset $B\subseteq E$ is bounded if and only if it is weakly bounded, meaning that for any continuous linear form $e'$ on $E$ the set $\{\langle b,e'\rangle,\; b\in B\}$ is bounded. If $B\subseteq E$ is a bounded convex and radial (that is $\lambda B\subseteq B$ for all $\lvert \lambda\rvert\leqslant 1$) subset of $E$ then we will denote by $E_B$ the subspace spanned by $B$ equipped with the norm

$$\displaystyle q_B(e)=\inf\{t\geqslant 0; \; e\in tB\},\;\;\; e\in E_B$$

\noindent We say that $E$ is {\em quasi-complete} if every bounded closed subset is complete. Of course complete implies quasi-complete. If $E$ is quasi-complete, then for every closed bounded convex and radial subset $B\subseteq E$, the space $E_B$ is a Banach space. We will denote by $E'$ the topological dual of $E$ that we will always endow it with the {\em strong topology} (that is the topology of uniform convergence on bounded subsets). If $E$ is quasi-complete, a subset $B\subseteq E'$ is bounded if and only if for all $e\in E$ the set $\{\langle e,b\rangle;\; b\in B\}$ is bounded (again by the Banach-Steinhaus theorem). More generally, if $F$ is another topological vector space then we will equip the space $\Hom(E,F)$ of continuous linear maps from $E$ to $F$ with the strong topology. Thus, a generating family of semi-norms for the topology on $\Hom(E,F)$ is given by

$$p_B(T)=\sup_{e\in B} p(Te),\;\;\; T\in \Hom(E,F)$$

\noindent where $B$ runs through the bounded subsets of $E$ and $p$ runs through a generating family of semi-norms for $F$.

\subsection{LF spaces}\label{section A.1}

\noindent Let $(E_i,f_{ij})$, $i\in I$, be a direct system of topological vector spaces (the connecting linear maps $f_{ij}$ being continuous). Consider the direct limit in the category of vector spaces

$$E=\varinjlim_{I} E_i$$

\noindent We will in general endow $E$ with the direct limit topology that is the finest locally convex topology on $E$ such that all the natural maps $E_i\to E$ are continuous. If $F$ is another topological space, a linear map $E\to F$ is continuous if and only if all the induced linear maps $E_i\to F$, $i\in I$ are continuous. If $I$ is at most countable, the $E_i$, $i\in I$, are Fr\'echet spaces and $E$ is Hausdorff (this is not automatic) then we will call $E$ an {\em LF space}. If $E$ is an LF space then we can write it as the direct limit of a sequence $(E_n)_{n\geqslant 0}$ of Fr\'echet spaces where the connecting maps are just inclusions $E_n\subseteq E_{n+1}$ so that we have

$$\displaystyle E=\bigcup_{n\geqslant 0}E_n$$

\noindent A bounded subset $B$ of $E$ is not always included in one of the $E_n$ but this is the case if $B$ is convex radial and complete in which case there exists $n\geqslant 0$ such that $B$ is included and bounded in $E_n$ (\cite{Gr1} Corollaire IV.2 p.17). Hence, in particular, if $E$ is quasi-complete then every bounded subset is included and bounded in some $E_n$.

\vspace{2mm}

\noindent LF spaces share many of the good properties of Fr\'echet spaces. First of all, an LF space is barrelled (\cite{Tr} Corollary 3 of Proposition 33.2), hence the Banach-Steinhaus theorem (aka uniform boundedness principle) applies to them (\cite{Tr} Theorem 33.1). In particular, if $E$ is an LF space, $F$ any topological vector space and $(u_n)_{n\geqslant 1}$ is a sequence of continuous linear maps from $E$ to $F$ that converges pointwise to a linear map $u:E\to F$, then $u$ is continuous and the sequence $(u_n)_{n\geqslant 1}$ converges to $u$ uniformly on compact subsets of $E$. Less known is the fact the open mapping theorem and the closed graph theorem continue to hold for LF spaces. This result is due to to Grothendieck (cf.\ \cite{Gr1} Theorem 4.B) and will be constantly used throughout this paper. We state it in the next proposition.

\begin{prop}\label{proposition A.1.1}
Let $E$ and $F$ be LF spaces. Then we have the following
\begin{enumerate}[(i)]
\item (Open mapping theorem) If $f:E\to F$ is a continuous surjective linear map, then it is open.

\item (Closed Graph theorem) A linear map $f:E\to F$ is continuous if and only if its graph is closed. 
\end{enumerate}
\end{prop}

\subsection{Vector-valued integrals}\label{section A.2}

\noindent Let $E$ be a topological vector space and let $X$ be a Hausdorff locally compact topological space equipped with a regular Borel measure $\mu$. Let $\varphi:X\to E$ be a continuous function.  We say that $\varphi$ is {\em weakly integrable} (with respect to $\mu$) if for all $e'\in E'$ the function $x\in X\mapsto \langle \varphi(x),e'\rangle$ is absolutely integrable. If this is so, there exists a unique vector $I(\varphi)$ in $(E')^*$, the algebraic dual of $E'$, such that

$$\displaystyle \langle I(\varphi),e'\rangle=\int_X \langle \varphi(x),e'\rangle d\mu(x),\;\;\; \mbox{for all } e'\in E'$$

\noindent We call $I(\varphi)$ the {\em weak integral} of $\varphi$ and we usually denote it simply by

$$\displaystyle \int_X \varphi(x) d\mu(x)$$

\noindent Note that this integral doesn't necessarily belong to $E$ or the completion of $E$.

\vspace{2mm}

\noindent We will say that $\varphi$ is {\em absolutely integrable} if for any continuous semi-norm $p$ on $E$ the function $x\in X\mapsto p(\varphi(x))$ is integrable. Of course, absolutely integrable implies weakly integrable. If $\varphi$ is absolutely integrable and $E$ is quasi-complete then the integral $\displaystyle \int_X \varphi(x) d\mu(x)$ belongs to $E$. In general, $\displaystyle \int_X \varphi(x) d\mu(x)$ belongs to $\widehat{E}$, the completion of $E$. Of course, if $\varphi$ is compactly supported then it is absolutely integrable (recall that we assume $\varphi$ to be continuous).

\subsection{Smooth maps with values in topological vector spaces}\label{section A.3}

\noindent Let $E$ be a topological vector space and $M$ be a real smooth manifold. Let $k\in \mathbb{N}\cup \{\infty\}$. Let $\varphi$ be a map from $M$ to $E$. We say that $\varphi$ is {\em strongly $C^k$} if it admits derivatives of all orders up to $k$ (in the classical sense) and that all these derivatives are continuous. We say that $\varphi$ is {\em weakly $C^k$} if for all $e'\in E'$ the map $m\in M\mapsto \langle \varphi(m),e'\rangle$ is of class $C^k$. Obviously, if $\varphi$ is strongly $C^k$ then it is also weakly $C^k$. A function that is weakly $C^\infty$ will also be called {\em smooth}. The next proposition summarizes the main properties of strongly and weakly $C^k$ functions. The first point of the proposition shows that for quasi-complete spaces there is not much difference between weakly $C^k$ and strongly $C^k$ maps. The second point implies that the notion of smooth functions with values in $E$ only depends on the bornology of $E$ (that is the family of its bounded subset). This is a very important property and the subsequent points of the proposition follow easily from it. In lack of a reference, we include a proof.

\vspace{2mm}

\begin{prop}\label{proposition A.3.1}

Let $k\in \mathbb{N}\cup\{\infty\}$. Then

\begin{enumerate}[(i)]

\item If $E$ is quasi-complete and $\varphi:M\to E$ is weakly $C^{k+1}$ then $\varphi$ is strongly $C^k$ .

\item A map $\varphi:M\to E$ is smooth if and only if for all $k\geqslant 0$ and every relatively compact open subset $\Omega\subset M$, there exists a bounded subset $B_k\subset E$ such that $\varphi_{\mid \Omega}$ factorizes through $E_{B_k}$ and such that the induced map $\Omega\to E_{B_k}$ is weakly $C^k$;

\item Let $A\subset E'$ be a subset such that every subset $B\subset E$ that is weakly-$A$-bounded (meaning that for all $e'\in A$ the set $\{\langle b,e'\rangle,\; b\in B\}$ is bounded) is bounded. Then, a map $\varphi:M\to E$ is smooth if and only if for all $a\in A$ the function $m\in M\mapsto \langle \varphi(m),a\rangle$ is smooth.

\item Assume again that $E$ is quasi-complete. Let $F$ be another topological vector space. Then a map $\varphi:M\to \Hom(E,F)$ is smooth if and only if for all $e\in E$ and $f'\in F'$, the map $m\in M\mapsto \langle \varphi(m)(e),f'\rangle$ is smooth.

\item Let $F$ and $G$ be two other topological vector spaces and assume that $E$ is quasi-complete. Let $A:E\times F\to G$ be a separately continuous bilinear map. Then, if $\varphi_1:M\to E$ and $\varphi_2:M\to F$ are smooth functions, the function

$$m\in M\mapsto A(\varphi_1(m),\varphi_2(m))\in G$$

\noindent is also smooth.
\end{enumerate}

\end{prop}

\vspace{2mm}

\noindent\ul{Proof}:

\begin{enumerate}[(i)]
\item The question being local, we may assume that $M=\mathbb{R}^n$. Assume first that $\varphi$ is weakly $C^1$. Then for all $e'\in E'$ and all $x,y\in \mathbb{R}^n$, we have

\[\begin{aligned}
\displaystyle \left\lvert \langle \varphi(y),e'\rangle-\langle \varphi(x),e'\rangle\right\rvert & =\left\lvert \int_0^1 \partial(y-x)\left(\langle \varphi(.),e'\rangle\right)\left((1-t)x+ty\right) dt\right\rvert \\
 & \leqslant \lVert y-x\rVert \sup_{z\in [x,y]} \left\lvert\partial\left(\frac{y-x}{\lVert y-x\rVert}\right)\left(\langle \varphi(.),e'\rangle\right)(z)\right\rvert
\end{aligned}\]

\noindent This shows that for every compact subset $\mathcal{K}\subseteq \mathbb{R}^n$ the family

$$\displaystyle \{ \frac{\varphi(y)-\varphi(x)}{\lVert y-x\rVert},\; x,y\in \mathcal{K}\}$$

\noindent is weakly bounded hence bounded so that $\varphi$ is Lipschitz hence continuous.

\vspace{2mm}

\noindent Assume now that $\varphi$ is weakly $C^2$. Let $u,x\in \mathbb{R}^n$ and $e'\in E'$. Then by Rolle's theorem, for all $t\in[-1,1]\smallsetminus \{0\}$ there exists $s_t\in [-1,1]$ with $\lvert s_t\rvert<\lvert t\rvert$ such that

$$\displaystyle \left\langle \frac{\varphi(x+tu)-\varphi(x)}{t},e'\right\rangle=\partial(u)\left(\left\langle \varphi(.),e'\right\rangle\right)(x+s_tu)$$

\noindent Hence, for all $t,t'\in [-1,1]\smallsetminus\{0\}$, we have

\[\begin{aligned}
\displaystyle \left\lvert\left\langle \frac{\varphi(x+tu)-\varphi(x)}{t}-\frac{\varphi(x+t'u)-\varphi(x)}{t'},e'\right\rangle\right\rvert & =\left\lvert\int_{s_{t'}}^{s_t}\partial(u^2)\left(\left\langle \varphi(.),e'\right\rangle\right)(x+su)ds\right\rvert \\
 & \leqslant \lvert s_t-s_{t'}\rvert \sup_{s\in [-1,1]}\left\lvert \partial(u^2)\left(\left\langle \varphi(.),e'\right\rangle\right)(x+su)\right\rvert \\
  & \leqslant \left(\lvert t\rvert+\lvert t'\rvert\right) \sup_{s\in [-1,1]}\left\lvert \partial(u^2)\left(\left\langle \varphi(.),e'\right\rangle\right)(x+su)\right\rvert
\end{aligned}\]

\noindent It follows that the family

$$\left\{ \left(\lvert t\rvert+\lvert t'\rvert\right)^{-1}\left(\frac{\varphi(x+tu)-\varphi(x)}{t}-\frac{\varphi(x+t'u)-\varphi(x)}{t'}\right),\;\; t,t'\in [-1,1]\smallsetminus\{0\} \right\}$$

\noindent is weakly bounded, hence bounded. Since $E$ is quasi-complete, this immediately implies that the limit

$$\displaystyle \lim\limits_{t\to 0} \frac{\varphi(x+tu)-\varphi(x)}{t}$$

\noindent exists in $E$. We just prove that the function $\varphi$ is strongly derivable everywhere in every direction. But for $u\in \mathbb{R}^n$, $\partial(u)\varphi$ is of course weakly $C^1$, hence continuous by the first part of the proof. This shows that $\varphi$ is indeed strongly $C^1$.

\vspace{2mm}

\noindent The general case follows by induction. Indeed, assume that the result is true for $k\geqslant 1$. Let $\varphi:\mathbb{R}^n\to E$ be a weakly $C^{k+2}$ map. Then, of course $\varphi$ is weakly $C^2$ hence strongly $C^1$ by what we just saw. Moreover, for every $u\in \mathbb{R}^n$, the function $\partial(u)\varphi$ is obviously weakly $C^{k+1}$ so that by the induction hypothesis it is also strongly $C^k$. This proves that $\varphi$ is in fact strongly $C^{k+1}$ and this ends the proof of (i).

\item Up to replacing $E$ by its completion, we may of course assume $E$ complete. Then, the proof of (i) shows that it suffices to take

$$B_k=\{(D\varphi)(m); \; m\in \overline{\Omega},\; D\in \Diff^\infty_{\leqslant k+1}(M)\}$$

\noindent for all $k\geqslant 0$.

\item As we already explained, the true meaning of (ii) is that the notion of smooth function with values in a topological vector space only depends on the bornology of that space. Hence, here it suffices to notice that the weak-$A$-topology on $E$ defines the same bornology as the original topology on $E$ (this is exactly the assumption made on $A$).

\item Again, this is because the topology defined by the semi-norms $T\mapsto \lvert \langle Te,f'\rangle\rvert$, $e\in E$ $f\in F'$, defines the same bornology on $\Hom(E,F)$ as the strong topology.

\item We are immediately reduced to the case when $G=\mathbb{C}$. Let $k\geqslant 0$ and $\Omega\subset M$ be a relatively compact open subset. Then by (ii), there exists bounded subsets $B_k\subseteq E$ and $B'_k\subseteq F$ such that the restrictions of $\varphi_1$ and $\varphi_2$ to $\Omega$ factorize through $E_{B_k}$ and $F_{B'_k}$ respectively and induce weakly $C^{k+1}$ maps from $\Omega$ into these spaces. Since $E$ is quasi-complete, the bilinear map $A$ restricted to $E_{B_k}\times F_{B'_k}$ induces a continuous bilinear form $E_{B_k}\times F_{B'_k}\to \mathbb{C}$. Hence, this bilinear form extends continuously to $\widehat{E}_{B_k}\times \widehat{F}_{B'_k}$, where $\widehat{E}_{B_k}$ and $\widehat{F}_{B'_k}$ denote the completion of $E_{B_k}$ and $F_{B'_k}$ respectively (these are Banach spaces). Since the maps $\varphi_1:\Omega\to \widehat{E}_{B_k}$ and $\varphi_2:\Omega\to \widehat{F}_{B'_k}$ are weakly $C^{k+1}$, by (i) they are also strongly $C^k$. It immediately follows that the map $m\in \Omega\mapsto A(\varphi_1(m),\varphi_2(m))$ is of class $C^k$. This of course implies (v). $\blacksquare$
\end{enumerate}

\vspace{2mm}

\noindent We will denote by $C^\infty(M,E)$ the space of all smooth maps from $M$ to $E$. We equip $C^\infty(M,E)$ with a topology as follows. If $E$ is complete then we endow $C^\infty(M,E)$ with the topology defined by the semi-norms

$$p_{D,\mathcal{K},q}(\varphi)=\sup_{m\in \mathcal{K}} q\left[(D\varphi)(m)\right]$$

\noindent where $\mathcal{K}$ runs through the compact subsets of $M$, $q$ runs through the continuous semi-norms on $E$ and $D$ runs through $\Diff^\infty(M)$ the space of smooth differential operators on $M$ (note that by the point (i) of the last proposition, since $E$ is complete, $\Diff^\infty(M)$ acts on $C^\infty(M,E)$). In the general case, we equip $C^\infty(M,E)$ with the subspace topology coming from the inclusion $C^\infty(M,E)\subset C^\infty(M,\widehat{E})$ where $\widehat{E}$ is the completion of $E$.

\subsection{Holomorphic maps with values in topological vector spaces}\label{section A.4}

\noindent Let $E$ be a topological vector space and $M$ a complex analytic manifold. Let $\varphi:M\to E$ be a map. We say that $\varphi$ is {\em holomorphic} if for all $e'\in E'$ the function $m\in M\mapsto \langle \varphi(m),e'\rangle$ is holomorphic and we say that $\varphi$ is {\em strongly holomorphic} if it admits complex derivatives of all orders. Obviously, if $\varphi$ is strongly holomorphic then $\varphi$ is holomorphic. The analog of Proposition \ref{proposition A.3.1} for holomorphic functions is true (\cite{Gr2} Theorem 1) and is summarized in the next proposition.

\vspace{2mm}

\begin{prop}\label{proposition A.4.1}
\begin{enumerate}[(i)]

\item If $E$ is quasi-complete and $\varphi:M\to E$ is holomorphic if and only if it is strongly holomorphic.

\item A map $\varphi:M\to E$ is holomorphic if and only if for every relatively compact open subset $\Omega\subset M$ there exists a bounded convex and radial subset $B\subset E$ such that $\varphi_{\mid \Omega}$ factorizes through $E_{B}$ and the induced map $\Omega\to E_{B}$ is holomorphic;

\item Let $A\subset E'$ be a subset such that every subset $B\subset E$ that is weakly-$A$-bounded (meaning that for all $e'\in A$ the set $\{\langle b,e'\rangle,\; b\in B\}$ is bounded) is bounded. Then, a map $\varphi:M\to E$ is holomorphic if and only if for all $a\in A$ the function $m\in M\mapsto \langle \varphi(m),a\rangle$ is holomorphic.

\item Assume again that $E$ is quasi-complete. Let $F$ be another topological vector space. Then a map $\varphi:M\to \Hom(E,F)$ is holomorphic if and only if for all $e\in E$ and $f'\in F'$, the map $m\in M\mapsto \langle \varphi(m)(e),f'\rangle$ is holomorphic.
\end{enumerate}
\end{prop}

\vspace{3mm}

\noindent Let $X$ be a compact real smooth manifold. Denote by $C^\infty(X)$ the space of smooth complex-valued functions on $X$. It is naturally a Fr\'echet space and we will denote by $C^{-\infty}(X)$ its topological dual. For each integer $k\geqslant 0$, we also have the space $C^k(X)$ of $C^k$ complex-valued functions on $X$. It is naturally a Banach space and its dual $C^{-k}(X)$ is also a Banach space. We have a natural continuous inclusion $C^{-k}(X)\subseteq C^{-\infty}(X)$. Let $M$ be a complex analytic manifold. In Section \ref{section B.3}, we will need the following fact

\vspace{3mm}

\begin{num}
\item\label{eq A.4.1} Let $\varphi:M\to C^{-\infty}(X)$ be holomorphic. Then, for every relatively compact open subset $\Omega\subseteq M$, there exists an integer $k\geqslant 0$ such that the map $\varphi_{\mid\Omega}$ factorizes through $C^{-k}(X)$ and induces an holomorphic map $\Omega\to C^{-k}(X)$.
\end{num}

\vspace{3mm}

\noindent This will follow from the point (ii) of the last proposition if we can prove that every bounded subset $B\subseteq C^{-\infty}(X)$ is contained and bounded in some $C^{-k}(X)$, $k\geqslant 0$. Note that we have a natural continuous and bijective linear map

\begin{align}\label{eq A.4.2}
\displaystyle \varinjlim_{k} C^{-k}(X)\to C^{-\infty}(X)
\end{align}

\noindent Since $C^{-\infty}(X)$ is the strong dual of a Fr\'echet space, it is complete (\cite{Bour} Proposition IV.3.2). Hence, if we can prove that \ref{eq A.4.2} is a topological isomorphism then we will be done (cf.\ Section \ref{section A.1}). By the open mapping theorem, it suffices to show that $C^{-\infty}(X)$ is an LF space. This follows from the fact that $C^\infty(X)$ is a nuclear space (cf. next section), hence it is reflexive and the strong dual of a reflexive Fr\'echet space is an LF space (\cite{Bour} IV.23 Proposition 4).

\subsection{Completed projective tensor product, nuclear spaces}\label{section A.5}

\noindent Let $E$ and $F$ be two topological vector spaces. For all continuous semi-norms $p$ and $q$ on $E$ and $F$ respectively, define a semi-norm $p\otimes q$ on $E\otimes F$ by

$$\displaystyle (p\otimes q)(v)=\inf\{\sum_{i=1}^n p(e_i)q(f_i);\; v=\sum_{i=1}^n e_i\otimes f_i\},\;\; v\in E\otimes F$$

\noindent where the infimum is taken over all decompositions $v=\sum_{i=1}^n e_i\otimes f_i$ of $v$ where $e_i\in E$ and $f_i\in F$. If we choose two families of semi-norms $(p_i)_{i\in I}$ and $(q_j)_{j\in J}$ that generate the topology on $E$ and $F$ respectively, then the family of semi-norms $(p_i\otimes q_j)_{(i,j)\in I\times J}$ defines a topology on $E\otimes F$ which is independent of the two chosen families $(p_i)_{i\in I}$ and $(q_j)_{j\in J}$. We shall call it the {\em projective topology} and we will denote by $E\gls{otimesp} F$ the completion of $E\otimes F$ for this topology. We call $E\widehat{\otimes}_p F$ the {\em completed projective tensor product} of $E$ and $F$. It satisfies the following universal property (\cite{Tr} Proposition 43.8): for any complete topological vector space $G$, the map that associates to $T\in \Hom(E\widehat{\otimes}_p F,G)$ the bilinear map $(e,f)\in E\times F\mapsto T(e\otimes f)$ induces a bijection

$$\displaystyle \Hom(E\widehat{\otimes}_p F,G)\simeq B(E,F;G)$$

\noindent where $B(E,F;G)$ denotes the space of all continuous bilinear maps $E\times F\to G$. We shall need the following

\vspace{3mm}

\begin{num}
\item\label{eq A.5.1} Let $E$ and $F$ be Fr\'echet spaces, $G$ be any topological vector space, $M$ be a real smooth manifold and $\varphi$ be a map $M\to \Hom(E\widehat{\otimes}_p F,G)$. Then, $\varphi$ is smooth if and only if for all $e\in E$, all $f\in F$ and all $g'\in G'$, the map
$$\displaystyle m\in M\mapsto \langle \varphi(m)(e\otimes f),g'\rangle$$
is smooth.
\end{num}

\vspace{3mm}

\noindent By Proposition \ref{proposition A.3.1}(iii), it suffices to see that a subset $B\subseteq \Hom(E\widehat{\otimes}_pF,G)$ is bounded if and only if for all $e\in E$, $f\in F$ and $g'\in G'$ the set $\{\langle b(e\otimes f),g'\rangle;\; b\in B\}$ is bounded. We are immediately reduced to the case where $G=\mathbb{C}$. Let $B$ be a set of continuous linear forms on $E\widehat{\otimes}_p F$, that we will identify with a set of continuous bilinear forms on $E\times F$, and assume that the set $\{b(e,f),\; b\in B\}$ is bounded for all $e\in E$ and $f\in F$. We want to show that $B$ is bounded in $\left(E\widehat{\otimes}_p F\right)'$. This amounts to proving that the set $\{b(v);\; b\in B\}$ is bounded for all $v\in E\widehat{\otimes}_p F$. Since $E$ and $F$ are Fr\'echet spaces, by \cite{Tr} Theorem 45.1, every $v\in E\widehat{\otimes}_p F$ is the sum of an absolutely convergent series

\begin{align}\label{eq A.5.2}
\displaystyle v=\sum_{n=0}^\infty \lambda_n x_n\otimes y_n
\end{align}

\noindent where $(\lambda_n)_{n\geqslant 0}$ is a sequence of complex numbers such that $\displaystyle \sum_{n=0}^\infty \lvert \lambda_n\rvert<\infty$ and $(x_n)_{n\geqslant 0}$ (resp.\ $(y_n)_{n\geqslant 0}$) is a sequence converging to $0$ in $E$ (resp.\ in $F$). Thus, we only need to show that the set $\{b(x_n,y_n);\; b\in B\; n\geqslant 0\}$ is bounded. But this follows from the usual uniform boundedness principle (for Banach spaces) using the fact that $E$ and $F$ are quasi-complete. This ends the proof of \ref{eq A.5.1}.

\vspace{2mm}

\noindent Another property of the projective tensor product that we shall need is the following.

\vspace{3mm}

\begin{num}
\item\label{eq A.5.3} Let $E$ and $F$ be Fr\'echet spaces. Then a sequence $(\lambda_n)_{n\geqslant 0}$ of continuous linear forms on $E\widehat{\otimes}_p F$ converges pointwise if and only if for all $e\in E$ and all $f\in F$ the sequence $(\lambda_n(e\otimes f))_{n\geqslant 0}$ converges.
\end{num}

\vspace{3mm}

\noindent Indeed, by the universal property of the projective tensor product, the $\lambda_n$ correspond to continuous bilinear forms $B_n\colon E\times F\to \mathbb{C}$. Now, if $B_n$ converges pointwise, since $E$ and $F$ are Fr\'echet spaces, by the Banach-Steinhaus theorem (and \cite{Tr} Theorem 34.1) the sequence $(B_n)_{n\geqslant 0}$ is equicontinuous and thus so does the sequence $(\lambda_n)_{n\geqslant 0}$. As it is converging on a dense subspace of $E\widehat{\otimes}_p F$, this sequence is therefore converging everywhere.

\vspace{2mm}

\noindent {\em Nuclear spaces} are a class of topological vector spaces that have many of the good properties of finite dimensional vector spaces. For the precise definition of nuclear spaces we refer the reader to \cite{Tr}. Examples of nuclear spaces are $C^\infty(M)$ or $C_c^\infty(M)$, where $M$ is a real smooth manifold. Any subspace and any quotient by a closed subspace of a nuclear space is nuclear (\cite{Tr} Proposition 50.1). A nuclear Fr\'echet space is a Montel space (\cite{Tr} Proposition 50.2), hence is reflexive. Important for us will be the following description of the completed projective tensor product of two nuclear spaces of functions (\cite{Gr1} Theorem 13). For every set $S$, let us denote by $\mathcal{F}(S)$ the space of all complex-valued functions on $S$. We shall endow this space with the topology of pointwise convergence.

\vspace{3mm}

\begin{prop}{(Grothendieck's weak-strong principle)}\label{proposition A.5.1}
Let $S$ and $T$ be two sets and let $E$ and $F$ be subpaces of $\mathcal{F}(S)$ and $\mathcal{F}(T)$ respectively. Assume that both $E$ and $F$ are equipped with locally convex topologies that are finer than the topology of pointwise convergence and that turn both $E$ and $F$ into nuclear LF spaces. Then, the natural bilinear map
$$E\times F\to \mathcal{F}(S\times T)$$
is continuous and extends to an injective continuous linear map
$$E\widehat{\otimes}_p F\to \mathcal{F}(S\times T)$$
whose image consists in the functions $\varphi:S\times T\to \mathbb{C}$ which satisfies

\begin{itemize}
\renewcommand{\labelitemi}{$\bullet$}
\item For all $s\in S$, the function $t\in T\mapsto \varphi(s,t)$ belongs to $F$;

\item For all $f'\in F'$, the function $s\in S\mapsto \langle \varphi(s,.),f'\rangle$ belongs to $E$.
\end{itemize}
\end{prop}

\vspace{2mm}

\noindent Finally, the next lemma will be used in Section \ref{section 4.2} in order to prove that spaces of quasi-characters are nuclear. Before stating it, we introduce a notation. Let $V\subseteq U\subseteq \mathbb{R}^n$ be open subsets. We will denote by $\gls{CbinfVU}$ the space of all smooth functions $\varphi:V\to \mathbb{C}$ such that for all $u\in S(\mathbb{R}^n)$ the function $\partial(u)\varphi$ (extended by $0$ outside $V$) is locally bounded on $U$. We endow this space with the topology defined by the semi-norms

$$\displaystyle q_{\mathcal{K},u}(\varphi)=\sup_{x\in \mathcal{K}} \left\lvert (\partial(u)\varphi)(x)\right\rvert,\;\;\; \varphi\in C^\infty_b(V,U)$$

\noindent where $\mathcal{K}$ runs through the compact subsets of $U$ and $u$ runs through $S(\mathbb{R}^n)$. With this topology, $C^\infty_b(V,U)$ is easily seen to be a Fr\'echet space.

\begin{lem}\label{lemma A.5.1}
Assume that the following condition is satisfied: For all $x\in U$, there exists an open neighborhood $U_x\subseteq U$ of $x$ such that $U_x\cap V$ can be written as a finite union of open convex subsets. Then, $C^\infty_b(V,U)$ is nuclear.
\end{lem}

\vspace{2mm}

\noindent\ul{Proof}: Fix for every $x\in U$ an open neighborhood $U_x\subseteq U$ satisfying the condition of the lemma. Set $V_x=U_x\cap V$ for all $x\in U$. Then, the natural restriction maps $C^\infty_b(V,U)\to C^\infty_b(V_x,U_x)$ induce a closed embedding

$$\displaystyle C_b^\infty(V,U)\hookrightarrow \prod_{x\in U} C_b^\infty(V_x,U_x)$$

\noindent Hence, we only need to prove that the spaces $C_b^\infty(V_x,U_x)$ are nuclear. We may thus assume that $V$ itself can be written as a finite union of open convex subsets. Let $V=V_1\cup V_2\cup\ldots\cup V_d$, $d\geqslant 1$, be such a presentation. The restriction maps again induce a closed embedding

$$\displaystyle C_b^\infty(V,U)\hookrightarrow \bigoplus_{i=1}^d C_b^\infty(V_i,U)$$

\noindent from which it follows that we may assume, without loss of generality, that $V$ is convex. If $V$ is convex, then for all $u\in S(\mathbb{R}^n)$ and all $\varphi\in C_b^\infty(V,U)$, by the mean value theorem, the function $\partial(u)\varphi$ extends continuously to $cl_U(V)$ (the closure of $V$ in $U$). By Whitney extension theorem \cite{Wh}, it follows that every function $\varphi\in C_b^\infty(V,U)$ extends to a smooth function on $U$. Hence, the restriction map

$$C^\infty(U)\to C_b^\infty(V,U)$$

\noindent is surjective. This linear map is obviously continuous and both $C^\infty(U)$ and $C_b^\infty(V,U)$ are Fr\'echet spaces. Consequently, by the open mapping theorem, $C_b^\infty(V,U)$ is a quotient of $C^\infty(U)$ and the result follows since $C^\infty(U)$ is nuclear. $\blacksquare$

\section{Some estimates}\label{section B}

\subsection{Three lemmas}\label{section B.1}

\begin{lem}\label{lemma B.1.1}
Let $\mathcal{U}\subseteq F^\times$ be a compact neighborhood of $1$. Then, for all $\delta>0$ there exists $\delta'>0$ such that

$$\displaystyle \int_{\mathcal{U}} \left(1+\left\lvert tx-t^{-1}y\right\rvert\right)^{-\delta}dt\ll \left(1+\lvert x\rvert\right)^{-\delta'}\left(1+\lvert y\rvert\right)^{-\delta'}$$

$$\displaystyle \int_{\mathcal{U}} \left(1+\left\lvert tx-y\right\rvert\right)^{-\delta}dt \ll \left(1+\lvert x\rvert\right)^{-\delta'}\left(1+\lvert y\rvert\right)^{-\delta'}$$

\noindent for all $x,y\in F$.
\end{lem}

\vspace{3mm}

\begin{lem}\label{lemma B.1.2}
\begin{enumerate}[(i)]
\item Let $T$ be a torus over $F$ and let $\chi_1,\ldots,\chi_n$ be linearly independent elements of $X^*(T)_{\overline{F}}$. Then, for all $Re(s)>0$, the function

$$\displaystyle t\mapsto \left\lvert \chi_1(t)-1\right\rvert^{s-1}\ldots\left\lvert \chi_n(t)-1\right\rvert^{s-1}$$

\noindent is locally integrable over $T(F)$.

\item Let $V$ be a finite dimensional $F$-vector space and let $\lambda_1,\ldots,\lambda_n$ be linearly independent elements in $V_{\overline{F}}^*$. Then, for all $Re(s)>0$, the function

$$\displaystyle v\mapsto \left\lvert \lambda_1(v)\right\rvert^{s-1}\ldots\left\lvert \lambda_n(v)\right\rvert^{s-1}$$

\noindent is locally integrable over $V$. Moreover, for all $d>n$ and every compact subset $\mathcal{K}_V\subseteq V$, we have

$$\displaystyle \lim\limits_{s\to 0^+} s^d\int_{\mathcal{K}_V} \left\lvert \lambda_1(v)\right\rvert^{s-1}\ldots\left\lvert \lambda_n(v)\right\rvert^{s-1} dv=0$$

\end{enumerate}
\end{lem}

\vspace{3mm}

\begin{lem}\label{lemma B.1.3}
Let $H$ be an algebraic group over $F$, $\sigma$ a log-norm on $H$ (cf.\ Section \ref{section 1.2}) and $d_Lh$ a left Haar measure on $H(F)$. For all $b>0$, let us denote by $\mathbf{1}_{\sigma<b}$ the characteristic function of $\{h\in H(F);\; \sigma(h)<b\}$. Then, there exists $R>0$ such that

$$\displaystyle \int_{H(F)} \mathbf{1}_{\sigma<b}(h)d_Lh\ll e^{Rb}$$

\noindent for all $b>0$.
\end{lem}

\subsection{Asymptotics of tempered Whittaker functions for general linear groups}\label{section B.2}

\noindent Let $V$ be an $F$-vector space of finite dimension $d$ and set $G=GL(V)$. Let $(e_1,\ldots,e_d)$ be a basis of $V$ and $(B,T)$ be the standard Borel pair of $G$ with respect to this basis. We have an isomorphism

$$T\simeq \left(\mathbb{G}_m\right)^d$$
$$t\mapsto (t_i)_{1\leqslant i\leqslant d}$$

\noindent where $t_i$, $1\leqslant i\leqslant d$, denotes the eigenvalue of $t$ acting on $e_i$. Set

$$\displaystyle \widetilde{T}:=\{t\in T; \; t_d=1\}$$

\noindent Let $N$ be the unipotent radical of $B$ and let $\xi$ be a generic character on $N(F)$, for example

$$\displaystyle \xi(n)=\psi\left(\sum_{i=1}^{d-1} \langle ne_{i+1},e_i^*\rangle\right),\;\;\; n\in N(F)$$

\noindent where $(e_1^*,\ldots,e_d^*)$ denotes the dual basis of $(e_1,\ldots,e_d)$.

\vspace{2mm}

\noindent Let $\pi\in \Temp(G)$ and fix a Whittaker model for $\pi$

$$\displaystyle \pi^\infty\hookrightarrow C^\infty\left(N(F)\backslash G(F),\xi\right)$$
$$v\mapsto W_v$$

\noindent Hence, there exists a nonzero continuous linear form

$$\ell:\pi^\infty\to \mathbb{C}$$

\noindent such that $\ell\circ \pi(n)=\xi(n)\ell$ for all $n\in N(F)$ and so that

$$W_v(g)=\ell(\pi(g)v)$$

\noindent for all $v\in \pi^\infty$ and all $g\in G(F)$.

\vspace{3mm}

\begin{lem}\label{lemma B.2.1}
For all $R>0$, there exists a continuous semi-norm $\nu_R$ on $\pi^\infty$ such that

$$\displaystyle \left\lvert W_v(\widetilde{t})\right\rvert \leqslant \nu_R(v)\Xi^G(\widetilde{t}) \prod_{i=1}^{d-1} \max\left(1,\lvert \widetilde{t}_i\rvert\right)^{-R}$$

\noindent for all $v\in \pi^\infty$ and all $\widetilde{t}\in\widetilde{T}(F)$.
\end{lem}

\vspace{3mm}

\noindent\ul{Proof}: Assume first that $F$ is $p$-adic. Then, we have the following stronger inequality whose proof is classical (cf. for example \cite{Wa4} Lemme 3.7 (i))

\vspace{3mm}

\begin{num}
\item\label{eq B.2.1} For all $v\in \pi^\infty$, there exists $c>0$ such that
$$\displaystyle \left\lvert W_v(\widetilde{t})\right\rvert \ll \Xi^{G}(\widetilde{t}) \prod_{i=1}^{d-1} \mathbf{1}_{]0,c]}(\left\lvert \widetilde{t}_i\right\rvert)$$
for all $\widetilde{t}\in \widetilde{T}(F)$ and where $\mathbf{1}_{]0,c]}$ denotes the characteristic function of the interval $]0,c]$.
\end{num}

\vspace{3mm}

\noindent We henceforth assume that $F=\mathbb{R}$. In this case, the inequality of the lemma is a consequence of the two following facts:

\vspace{3mm}

\begin{num}
\item\label{eq B.2.2} There exists $R_0>0$ and a continuous semi-norm $\nu$ on $\pi^\infty$ such that
$$\displaystyle \left\lvert W_v(\widetilde{t})\right\rvert \leqslant \nu(v) \Xi^G(\widetilde{t}) \prod_{i=1}^{d-1} \max\left(1,\lvert \widetilde{t}_i\rvert\right)^{R_0}$$
for all $v\in \pi^\infty$ and all $\widetilde{t}\in \widetilde{T}(\mathbb{R})$.
\end{num} 

\vspace{3mm}

\begin{num}
\item\label{eq B.2.3} For all $1\leqslant i\leqslant d-1$, there exists $u_i\in \mathcal{U}(\mathfrak{n})$ such that
$$\displaystyle W_v(\widetilde{t})=\widetilde{t}_i^{-1} W_{\pi(u_i)v}(\widetilde{t})$$
for all $v\in \pi^\infty$ and all $\widetilde{t}\in \widetilde{T}(\mathbb{R})$.
\end{num}

\vspace{3mm}

\noindent Let us identify $\mathfrak{g}$ with $\mathfrak{gl}_d$ using the basis $e_1,\ldots,e_d$ and let us denote for $1\leqslant i\leqslant d-1$ by $X_i\in \mathfrak{g}$ the matrix with a $1$ at the crossing of the $i$th row and $(i+1)$th column and zeros everywhere else. Set $u_i'=X_iX_{i+1}\ldots X_{d-1}$ and $u_i=d\xi(u'_i)^{-1} u'_i$ for $1\leqslant i\leqslant d-1$. Here $d\xi:\mathcal{U}(\mathfrak{n})$ is the natural extension of the character $d\xi:\mathfrak{n}(\mathbb{R})\to \mathbb{C}$ obtained by differentiating $\xi$ at the origin. Note that the elements $u_i$ are well-defined since by the hypothesis that $\xi$ is generic we have $d\xi(u'_i)\neq 0$ for all $1\leqslant i\leqslant d-1$. Moreover, it is easy to check that for all $1\leqslant i\leqslant d-1$, $u_i$ satisfies the claim \ref{eq B.2.3} We are thus only left with proving \ref{eq B.2.2}.

\vspace{2mm}

\noindent In what follows, we fix a positive integer $k$ that we assume sufficiently large throughout. Denote by $\overline{B}=T\overline{N}$ the Borel subgroup opposite to $B$ with respect to $T$. Let $Y_1,\ldots,Y_b$ be a basis of $\overline{\mathfrak{b}}(\mathbb{R})$ and set $\Delta_{\overline{B}}=Y_1^2+\ldots+Y_b^2\in \mathcal{U}(\overline{\mathfrak{b}})$. Then, by elliptic regularity (cf.\ \ref{eq 2.1.2}), there exists functions $\varphi_{\overline{B}}^1\in C_c^{2k-\dim(\overline{B})-1}(\overline{B}(\mathbb{R}))$ and $\varphi_{\overline{B}}^2\in C_c^\infty(\overline{B}(\mathbb{R}))$ such that
$$\displaystyle \varphi_{\overline{B}}^1\ast \Delta_{\overline{B}}^k+\varphi_{\overline{B}}^2=\delta_1^{\overline{B}}$$
Applying $\pi$ to this equality, we get
$$\displaystyle W_v(\widetilde{t})=\ell\left(\pi(\widetilde{t})v\right)=\ell\left(\pi(\varphi_{\overline{B}}^1)\pi(\Delta_{\overline{B}}^k)\pi(\widetilde{t})v\right)+\ell\left(\pi(\varphi_{\overline{B}}^2)\pi(\widetilde{t})v\right)$$
for all $v\in \pi^\infty$ and all $\widetilde{t}\in \widetilde{T}(\mathbb{R})$. Let $\varphi_N\in C_c^\infty(N(\mathbb{R}))$ be any function such that $\int_{N(\mathbb{R})} \varphi_N(n)\xi(n)dn=1$. Then we have $\ell=\ell\circ \pi(\varphi_N)$. Plugging this into the last equality, we get

\begin{align}\label{eq B.2.4}
\displaystyle W_v(\widetilde{t})=\ell\left(\pi(\varphi^1)\pi(\Delta_{\overline{B}}^k)\pi(\widetilde{t})v\right)+\ell\left(\pi(\varphi^2)\pi(\widetilde{t})v\right)
\end{align}

\noindent for all $v\in \pi^\infty$ and all $\widetilde{t}\in \widetilde{T}(\mathbb{R})$ where we have set $\varphi^1=\varphi_N\ast \varphi_{\overline{B}}^1$ and $\varphi^2=\varphi_N\ast \varphi_{\overline{B}}^2$. Note $\varphi^1$ and $\varphi^2$ both belong to $C_c^{2k-\dim(\overline{B})-1}(G(\mathbb{R}))$. It follows that for $k$ sufficiently large, the two vectors $\ell\circ \pi(\varphi^1)$, $\ell\circ \pi(\varphi^2)\in \pi^{-\infty}$ actually belong to $\overline{\pi}^\infty$. Assuming $k$ to be that sufficiently large, by \ref{eq 2.2.6} we get the existence of a continuous semi-norm $\nu_0$ on $\pi^\infty$ such that
$$\displaystyle \left\lvert\ell\left(\pi(\varphi^i)\pi(g)v\right)\right\rvert \leqslant \nu_0(v) \Xi^G(g)$$
for all $v\in \pi^\infty$, all $g\in G(\mathbb{R})$ and all $i\in \{1,2\}$. Combining this inequality with \ref{eq B.2.4}, we get
$$\displaystyle \left\lvert W_v(\widetilde{t})\right\rvert \leqslant \left(\nu_0\left(\pi(\widetilde{t}^{-1}\Delta_{\overline{B}}^k\widetilde{t})v\right)+\nu_0(v)\right)\Xi^G(\widetilde{t})$$
for all $v\in \pi^\infty$ and all $\widetilde{t}\in \widetilde{T}(\mathbb{R})$. To ends the proof of \ref{eq B.2.2}, it suffices to notice that $\Delta_{\overline{B}}^k$ is a sum of eigenvectors for the adjoint action of $\widetilde{T}(\mathbb{R})$ on $\mathcal{U}(\overline{\mathfrak{b}})$ and that any such eigenvector has an associated eigen-character of the form $\widetilde{t}\in \widetilde{T}(\mathbb{R})\mapsto \widetilde{t}_1^{-n_1}\ldots \widetilde{t}_{d-1}^{-n_{d-1}}$ for some nonnegative integers $n_1,\ldots,n_{d-1}$. $\blacksquare$

\subsection{Unipotent estimates}\label{section B.3}

\noindent Let us fix the following

\vspace{3mm}

\begin{itemize}
\renewcommand{\labelitemi}{$\bullet$}

\item $G$ a connected reductive group over $F$;

\item $P_{\mini}=M_{\mini}N_{\mini}$ and $\overline{P}_{\mini}=M_{\mini}\overline{N}_{\mini}$ two opposite minimal parabolic subgroups of $G$;

\item $A_{\mini}=A_{M_{\mini}}$ denotes the maximal split subtorus of $M_{\mini}$;

\item $\lambda_{\mini}:N_{\mini}\to \mathbb{G}_a$ is a non-degenerate additive character;

\item $\xi_{\mini}=\psi\circ \lambda_{min,F}:N_{\mini}(F)\to \mathbb{C}^\times$ where $\psi:F\to \mathbb{C}^\times$ is a continuous unitary character;

\item $N'_{\mini}=\Ker(\lambda_{\mini})$.

\item $K$ is a maximal compact subgroup of $G(F)$ that is special in the $p$-adic case. We denote by $m_{\overline{P}_{\mini}}:G(F)\to M_{\mini}(F)$ any map such that $m_{\overline{P}_{\mini}}(g)^{-1}g\in \overline{N}_{\mini}(F)K$ for all $g\in G(F)$.
\end{itemize}

\vspace{3mm}

\noindent The purpose of this section is to show the following estimate

\begin{prop}\label{proposition B.3.1}
There exists $\epsilon>0$ such that the integral

$$\displaystyle \int_{N'_{\mini}(F)} \delta_{\overline{P}_{\mini}}(m_{\overline{P}_{\mini}}(n'n))^{1/2-\epsilon} dn'$$

\noindent is absolutely convergent for all $n\in N_{\mini}(F)$ and is bounded uniformly in $n$.
\end{prop}

\noindent To prove this estimate, we will use the holomorphic continuation of the Jacquet integral. Let us recall what it means. For all $s\in \mathbb{C}$, we introduce the smooth normalized induced representation

$$\pi_s^\infty=i_{\overline{P}_{\mini}}^G(\delta_{\overline{P}_{\mini}}^s)^\infty$$

\noindent By restriction to $K$, all the spaces underlying the representations $\pi_s^\infty$ (for $s\in \mathbb{C}$) become isomorphic to $C^\infty(K_{\mini}\backslash K)$, where $K_{\mini}=K\cap \overline{P}_{\mini}(F)$. We will use these isomorphisms as identifications and for $e\in C^\infty(K_{\mini}\backslash K)$, $s\in \mathbb{C}$, we will denote by $e_s$ the corresponding vector in the space of $\pi^\infty_s$. Now for $Re(s)>0$, we may define the following functional

$$\Lambda_s:C^\infty(K_{\mini}\backslash K)\to \mathbb{C}$$
$$\displaystyle e\mapsto \int_{N_{\mini}(F)} e_s(n)\overline{\xi(n)} dn$$

\noindent (the integral is absolutely convergent). This functional is called the Jacquet integral. The space $C^\infty(K_{\mini}\backslash K)$ is naturally a topological vector space: if $F$ is Archimedean then it has a structure of Fr\'echet space whereas if $F$ is $p$-adic we equip it with the finest locally convex topology. Then $\Lambda_s$, for $Re(s)>0$, is a continuous linear form hence it belongs to the topological dual of $C^\infty(K_{\mini}\backslash K)$ that we will denote by $C^{-\infty}(K_{\mini}\backslash K)$. The holomorphic continuation of the Jacquet integral now means the following

\vspace{3mm}

\hspace{6mm} The map $s\in \{Re>0\}\mapsto \Lambda_s\in C^{-\infty}(K_{\mini}\backslash K)$ is holomorphic and admits an

\hspace{6mm} holomorphic continuation to $\mathbb{C}$.

\vspace{3mm}

\noindent\ul{Proof of Proposition }\ref{proposition B.3.1}: We may assume without loss of generality that $G$ is adjoint. This implies the existence of a one-parameter subgroup

$$a:\mathbb{G}_m\to A_{\mini}$$
$$x\mapsto a(x)$$

\noindent such that $\lambda_{\mini}\left(a(x)na(x)^{-1}\right)=x\lambda(n)$ for all $x\in \mathbb{G}_m$ and all $n\in N$. Note that $a\in X_*(A_{\mini})$ is in the positive chamber corresponding to $P_{\mini}$. We start by proving the following

\vspace{3mm}

\begin{num}
\item\label{eq B.3.1} For all $e\in C^\infty(K_{\mini}\backslash K)$, all $n\in N_{\mini}(F)$ and all $s\in \{Re>0\}$, we have
$$\displaystyle \int_{N'_{\mini}(F)} e_s(n'n)dn'=\int_{F} \psi(x\lambda_{\mini}(n)) \Lambda_s\left(\pi_s(a(x))e\right) \delta_{P_{\mini}}(a(x))^{s-1/2} dx$$
where $dx$ denotes some additive Haar measure on $F$ and both integrals are absolutely convergent.
\end{num}

\vspace{3mm} 

\noindent Fix $e\in C^\infty(K_{\mini}\backslash K)$ and $s\in \mathbb{C}$ such that $Re(s)>0$. Using $\lambda_{min,F}$, we may identify $N'_{\mini}(F)\backslash N_{\mini}(F)$ with $F$. Then both sides of \ref{eq B.3.1} may be seen as a functions on $F$:

$$\displaystyle \varphi_1:y\in F\mapsto \int_{N'_{\mini}(F)} e_s(n'y)dn'$$

$$\displaystyle \varphi_2:y\in F\mapsto \int_F \psi(xy) \Lambda_s\left(\pi_s(a(x))e\right) \delta_{P_{\mini}}(a(x))^{s-1/2} dx$$

\noindent The integral defining the $\varphi_1$ is absolutely convergent and the resulting function is integrable over $F$. Moreover, in the $p$-adic case we obtain a uniformly smooth function on $F$ whereas in the Archimedean case we obtain a function on $F$ that is smooth with all its derivatives integrable. All of these easily follow from the following two facts

\vspace{3mm}

\begin{num}
\item\label{eq B.3.2} For every compact $C\subseteq N_{\mini}(F)$, we have
$$\lvert e_s(n'n)\rvert\ll \delta_{\overline{P}_{\mini}}\left(m_{\overline{P}_{\mini}}(n')\right)^{1/2+Re(s)}$$
for all $n'\in N_{\mini}(F)$ and all $n\in C$.
\end{num} 

\vspace{3mm}

\begin{num}
\item\label{eq B.3.3} The integral
$$\int_{N_{\mini}(F)}\delta_{\overline{P}_{\mini}}\left(m_{\overline{P}_{\mini}}(n)\right)^{1/2+Re(s)} dn$$
is absolutely convergent.
\end{num}

\vspace{3mm}

\noindent In both cases, this implies that $\varphi_1$ admits a Fourier transform

$$\displaystyle \widehat{\varphi}_1(x)=\int_F \varphi_1(y) \overline{\psi(xy)} dy$$

\noindent (the Haar measure $dy$ being the quotient of the Haar measure on $N_{\mini}(F)$ and on $N'_{\mini}(F)$) which is integrable over $F$ and that for $dx$ the Haar measure dual to $dy$, we have

$$\displaystyle \varphi_1(y)=\int_F \widehat{\varphi}_1(x)\psi(xy) dx$$

\noindent for all $y\in F$. To prove \ref{eq B.3.1} it is thus sufficient to establish that $\Lambda_s\left(\pi_s(a(x))e\right) \delta_{P_{\mini}}(a(x))^{s-1/2}=\widehat{\varphi}_1(x)$ for all $x\in F^\times$. We have

\[\begin{aligned}
\displaystyle \Lambda_s\left(\pi_s(a(x))e\right) & =\int_{N_{\mini}(F)} e_s(na(x))\overline{\xi_{\mini}(n)}dn \\
 & =\delta_{\overline{P}_{\mini}}(a(x))^{1/2+s}\int_{N_{\mini}(F)} e_s(a(x)^{-1}na(x)) \overline{\xi_{\mini}(n)} dn \\
 & =\delta_{P_{\mini}}(a(x))^{1/2-s}\int_{N_{\mini}(F)} e_s(n)\overline{\xi_{\mini}(a(x)na(x)^{-1})} dn \\
 & =\delta_{P_{\mini}}(a(x))^{1/2-s}\int_F \varphi_1(y)\overline{\psi(xy)}dy
\end{aligned}\]

\noindent for all $x\in F^\times$, where at the third line we made the variable change $n\mapsto a(x)na(x)^{-1}$. This proves the equality $\Lambda_s\left(\pi_s(a(x))e\right) \delta_{P_{\mini}}(a(x))^{s-1/2}=\widehat{\varphi}_1(x)$ for all $x\in F^\times$ and ends the proof of \ref{eq B.3.1}.

\vspace{2mm}

\noindent We will now prove the following

\vspace{3mm}

\begin{num}
\item\label{eq B.3.4} There exists $\delta>0$ such that for all $e\in C^\infty(K_{\mini}\backslash K)$ the integral
$$\displaystyle \int_F \Lambda_s\left(\pi_s(a(x))e\right) \delta_{P_{\mini}}(a(x))^{s-1/2} dx$$
is absolutely and locally uniformly convergent for all $s$ in $\{Re>-\delta\}$.
\end{num}

\vspace{3mm}

\noindent First we show how to deduce the proposition from this last claim. It implies in particular that the right hand side of \ref{eq B.3.1} admits an holomorphic continuation to some half plan $\{Re>-\delta\}$, $\delta>0$. Consequently, the left hand side also admits an holomorphic continuation to such an half plane. Let us consider the case where $e=e_0\in C^\infty(K_{\mini}\backslash K)$ is the constant function equal to $1$. Then the left hand side of \ref{eq B.3.1} is, for $Re(s)>0$,

$$\displaystyle \int_{N'_{\mini}(F)} \delta_{\overline{P}_{\mini}}(m_{\overline{P}_{\mini}}(n'n))^{1/2+s}dn'$$

\noindent Since the integrand is positive for $s$ real, this implies that the integral is still absolutely convergent in the half plane $\{Re>-\delta\}$. Hence, for $\epsilon<\delta$, the integral of the proposition is convergent. Moreover, it is equal to the integral of the right hand side of \ref{eq B.3.1} evaluated at $s=-\epsilon$. By \ref{eq B.3.4}, this integral is also absolutely convergent. Since the absolute value of the integrand is independent of $n\in N(F)$, this shows the uniform boundedness of the proposition. Hence, we are left with establishing \ref{eq B.3.4}.

\vspace{2mm}

\noindent Fix $e\in C^\infty(K_{\mini}\backslash K)$. We now split \ref{eq B.3.4} into the two following estimates:

\vspace{3mm}

\begin{num}
\item\label{eq B.3.5} The integral
$$\displaystyle \int_{\{\lvert x\rvert\geqslant 1\}} \Lambda_s\left(\pi_s(a(x))e\right) \delta_{P_{\mini}}(a(x))^{s-1/2} dx$$
is absolutely convergent locally uniformly in $s$ for all $s\in \mathbb{C}$.
\end{num}

\vspace{3mm}

\begin{num}
\item\label{eq B.3.6} There exists $\delta>0$ such that the integral
$$\displaystyle \int_{\{\lvert x\rvert <1\}} \Lambda_s\left(\pi_s(a(x))e\right) \delta_{P_{\mini}}(a(x))^{s-1/2} dx$$
is absolutely convergent locally uniformly in $s$ for all $s\in\{Re>-\delta\}$.
\end{num}

\vspace{3mm}

\noindent When $F=\mathbb{R}$, for each integer $k\geqslant 0$ the (continuous) dual $C^{-k}(K_{\mini}\backslash K)$ of the space $C^k(K_{\mini}\backslash K)$ of functions continuously derivable up to order $k$ is a Banach space which naturally embeds in $C^{-\infty}(K_{\mini}\backslash K)$ and by \ref{eq A.4.1} we have the following:

\vspace{3mm}

\begin{num}
\item\label{eq B.3.7} Assume $F=\mathbb{R}$. Then, for every relatively compact open subset $\Omega\subseteq \mathbb{C}$, there exists an integer $k\geqslant 0$ such that the map $s\in \Omega\mapsto \Lambda_s$ factors through $C^{-k}(K_{\mini}\backslash K)$ and defines an holomorphic function into that Banach space.
\end{num}

\vspace{3mm}

\noindent We now prove \ref{eq B.3.5}. By assumption, there exists $X\in \mathfrak{n}_{\mini}$ such that $\xi_{\mini}(e^X)\neq 1$ and $a(x)Xa(x)^{-1}=xX$ for all $x\in F^\times$. We now separate the proof according to whether $F$ is $p$-adic or real

\vspace{3mm}

\begin{itemize}
\renewcommand{\labelitemi}{$\bullet$}

\item First assume that $F$ is a $p$-adic field. It is easy to see that $x\in F^\times\mapsto \Lambda_s\left(\pi_s(a(x))e\right)$ is bounded on compact subsets of $F^\times$ locally uniformly in $s$. Hence, it is sufficient to establish the following

\vspace{3mm}

\begin{num}
\item\label{eq B.3.8} There exists $c\geqslant 1$ such that
$$\Lambda_s\left(\pi_s(a(x))e\right)=0$$
for all $s\in \mathbb{C}$ and all $x\in F^\times$ satisfying $\lvert x\rvert\geqslant c$.
\end{num}

\vspace{3mm}

\noindent Indeed, for all $s\in \mathbb{C}$ and all $x\in F^\times$, we have

\[\begin{aligned}
\xi_{\mini}(e^X)\Lambda_s\left(\pi_s(a(x))e\right) & =\Lambda_s\left(\pi_s(e^X)\pi_s(a(x))e\right) \\
 & =\Lambda_s\left(\pi_s(a(x))\pi_s(e^{x^{-1}X})e\right)
\end{aligned}\]

\noindent and for $x$ sufficiently large, $e^{x^{-1}X}$ is in $K$ and stabilizes $e$. Since $\xi_{\mini}(e^X)\neq 1$, this shows the vanishing \ref{eq B.3.8}.

\item Let us now treat the case where $F=\mathbb{R}$. Fix $\Omega$ a compact subset of $\mathbb{C}$. By \ref{eq B.3.7}, there exists an integer $k\geqslant 0$ such that
$$\displaystyle \left\lvert \Lambda_s(e')\right\rvert\ll \lVert e'\rVert_{C^k}$$
for all $e'\in C^\infty(K_{\mini}\backslash K)$ and all $s\in \Omega$. Moreover, it is easy to see that there exists a positive integer $N_0>0$ such that
$$\displaystyle \lVert \pi_s(a(x))e'\rVert_{C^k}\ll \lvert x\rvert^{N_0} \lVert e'\rVert_{C^k}$$
for all $e'\in  C^\infty(K_{\mini}\backslash K)$, all $s\in \Omega$ and all $x\in F^\times$ such that $\lvert x\rvert\geqslant 1$. Hence, we have

\begin{align}\label{eq B.3.9}
\displaystyle \left\lvert \Lambda_s\left(\pi_s(a(x))e'\right)\right\rvert\ll \lvert x\rvert^{N_0} \lVert e'\rVert_{C^k}
\end{align}

\noindent for all $e'\in  C^\infty(K_{\mini}\backslash K)$, all $s\in \Omega$ and all $x\in F^\times$ such that $\lvert x\rvert\geqslant 1$. Also, there exists a positive integer $N_1$ such that

\begin{align}\label{eq B.3.10}
\displaystyle \left\lvert \delta_{P_{\mini}}(a(x))^{s-1/2}\right\rvert\ll \lvert x\rvert^{N_1}
\end{align}

\noindent for all $s\in \Omega$ and all $x\in F^\times$ such that $\lvert x\rvert\geqslant 1$. \\

Consider the element $X\in \mathfrak{n}_{\mini}(F)$ previously introduced. Up to a scaling, we may assume that $d\xi_{\mini}(X)=1$. Then, for every positive integer $N_2$, we have

\[\begin{aligned}
\displaystyle \Lambda_s\left(\pi_s(a(x))e\right) & =d\xi_{\mini}(X)^{N_2} \Lambda_s\left(\pi_s(a(x))e\right)\\
 & =\Lambda_s\left(\pi_s(X^{N_2})\pi_s(a(x))e\right) \\
 & =\lvert x\rvert^{-N_2}\Lambda_s\left(\pi_s(a(x))\pi_s(X^{N_2})e\right)
\end{aligned}\]

\noindent Since the family $\left(\pi_s(X^{N_2})e\right)_{s\in \Omega}$ is bounded in $C^\infty(K_{\mini}\backslash K)$, combining the previous equality with \ref{eq B.3.9}, we get that for every integer $N_2>0$, we have an inequality
$$\displaystyle \left\lvert \Lambda_s\left(\pi_s(a(x))e\right)\right\rvert\ll \lvert x\rvert^{N_0-N_2}$$
for all $s\in \Omega$ and all $x\in F^\times$ such that $\lvert x\rvert\geqslant 1$. Combining this further with \ref{eq B.3.10}, we get an inequality
$$\displaystyle \left\lvert \delta_{P_{\mini}}(a(x))^{s-1/2}\Lambda_s\left(\pi_s(a(x))e'\right)\right\rvert\ll \lvert x\rvert^{-2}$$
for all $s\in \Omega$ and all $x\in F^\times$ such that $\lvert x\rvert\geqslant 1$. This ends the proof of \ref{eq B.3.5}.
\end{itemize}

\vspace{4mm}

\noindent We now prove \ref{eq B.3.6}. We will deduce it from the following claim:

\vspace{3mm}

\begin{num}
\item\label{eq B.3.11} There exists $d>0$ such that we have an inequality which is uniform locally in $s$
$$\left\lvert \Lambda_s\left(\pi_s(a(x))e\right)\right\rvert\ll \delta_{P_{\mini}}(a(x))^{1/2-\lvert Re(s)\rvert}\sigma(a(x))^d$$
for all $x\in F^\times$ such that $\lvert x\rvert<1$.
\end{num}

\vspace{3mm}

\noindent That \ref{eq B.3.11} implies \ref{eq B.3.6} is clear since for $\delta>0$ sufficiently small the function $x\mapsto \delta_{P_{\mini}}(a(x))^{-\delta}\sigma(a(x))^d$ is locally integrable on $F$.

\vspace{2mm}

\noindent It only remains to prove \ref{eq B.3.11}. Define a bilinear form $\langle.,.\rangle_K$ on $C(K_{\mini}\backslash K)$ (space of all complex-valued continuous functions on $K_{\mini}\backslash K$) by
$$\displaystyle \langle e,e'\rangle_K=\int_K e(k)e'(k)dk,\;\;\; e,e'\in C(K_{\mini}\backslash K)$$
This pairing induces a continuous embedding $C(K_{\mini}\backslash K)\subset C^{-\infty}(K_{\mini}\backslash K)$ and we have $\langle \pi_s(g)e,e'\rangle_K=\langle e,\pi_{-s}(g^{-1})e'\rangle_K$ for all $e,e'\in C(K_{\mini}\backslash K)$, all $s\in \mathbb{C}$ and all $g\in G(F)$. Moreover, the space $C(K_{\mini}\backslash K)$ is a Banach space when equipped with the norm
$$\displaystyle \lVert e\rVert_{\infty}=\sup_{k\in K}\; \lvert e(k)\rvert,\;\;\; e\in C(K_{\infty}\backslash K)$$
We have

\vspace{3mm}

\begin{num}
\item\label{eq B.3.12} There exists $d>0$, such that
$$\left\lvert \langle e,\pi_s(a(x))e'\rangle_K\right\rvert\ll \delta_{P_{\mini}}(a(x))^{1/2-\lvert Re(s)\rvert} \sigma(a(x))^d \lVert e\rVert_{\infty} \lVert e'\rVert_{\infty}$$
for all $s\in \mathbb{C}$, all $e,e'\in C(K_{\mini}\backslash K)$ and all $x\in F^{\times}$ such that $\lvert x \rvert<1$.
\end{num}

\vspace{3mm}

\noindent Indeed, let $e,e'\in C(K_{\mini}\backslash K)$, unraveling the definitions we have
$$\displaystyle \langle e,\pi_s(a(x))e'\rangle_K=\int_K \delta_{\overline{P}_{\mini}}\left(m_{\overline{P}_{\mini}}(ka(x))\right)^{1/2+s} e'\left(k_{\overline{P}_{\mini}}(ka(x))\right) e(k) dk$$
Hence,
$$\displaystyle \left\lvert \langle e,\pi_s(a(x))e'\rangle_K\right\rvert\leqslant \lVert e\rVert_{\infty}\lVert e'\rVert_{\infty} \sup_{k\in K} \left[\delta_{\overline{P}_{\mini}}\left(m_{\overline{P}_{\mini}}(ka(x))\right)^{Re(s)}\right] \int_K \delta_{\overline{P}_{\mini}}\left(m_{\overline{P}_{\mini}}(ka(x))\right)^{1/2}dk$$
for all $s\in \mathbb{C}$ and all $x\in F^\times$. The integral above is (by definition) $\Xi^G(a(x))$ and since $a(.)$ is in the positive chamber relative to $P_{\mini}$, by Proposition \ref{proposition 1.5.1}(i) there exists a $d_0>0$ such that $\Xi^G(a(x))\ll \delta_{P_{\mini}}(a(x))^{1/2}\sigma(a(x))^{d_0}$ for all $x\in F^\times$ with $\lvert x\rvert<1$. On the other hand, it easily follows from Proposition \ref{proposition 1.5.1}(i) and (ii) that there exists $d_1>0$ such that
$$\displaystyle \sup_{k\in K} \left[\delta_{\overline{P}_{\mini}}\left(m_{\overline{P}_{\mini}}(ka(x))\right)^{Re(s)}\right]\ll \delta_{P_{\mini}}(a(x))^{-\lvert Re(s)\rvert}\sigma(a(x))^{d_1}$$ 
for all $s\in \mathbb{C}$ and all $x\in F^\times$ such that $\lvert x \rvert <1$. The estimate \ref{eq B.3.12} follows.

\vspace{2mm}

\noindent We now prove \ref{eq B.3.11} distinguishing again the case where $F$ is $p$-adic from the case where $F=\mathbb{R}$.

\vspace{5mm}

\begin{itemize}
\renewcommand{\labelitemi}{$\bullet$}

\item Assume first that $F$ is a $p$-adic field. Since $a\in X_*(A_{\mini})$ is in the positive chamber for $P_{\mini}$, there exists a compact-open subgroup $K_{\overline{P}}\subseteq K_{\mini}\cap \overline{P}_{\mini}(F)$ fixing $e$ and such that $a(x)^{-1}K_{\overline{P}}a(x)\subseteq K_{\overline{P}}$ for all $x\in F^\times$ with $\lvert x \rvert<1$. Fix such a subgroup and let $K'\subseteq K\cap \Ker(\xi_{\mini})K_{\overline{P}}$ be a compact-open subgroup. Then, we have
$$\displaystyle \Lambda_s\left(\pi_s(a(x))e\right)=\Lambda_s\left(\pi(e_{K'})\pi_s(a(x))e\right)$$
for all $s\in \mathbb{C}$ and all $x\in F^\times$ such that $\lvert x\rvert<1$, where $e_{K'}=\vol(K')^{-1}\mathbf{1}_{K'}$. Now the function $s\mapsto \Lambda_s\circ \pi(e_{K'})$ actually takes value in $C^\infty(K_{\mini}\backslash K)$ and we have
$$\displaystyle \Lambda_s\left(\pi(e_{K'})e'\right)=\langle \Lambda_s\circ \pi(e_{K'}),e'\rangle_K$$
for all $e'\in C^\infty(K_{\mini}\backslash K)$ and all $s\in \mathbb{C}$. By \ref{eq B.3.12}, it is thus sufficient to show that the map $s\mapsto \Lambda_s\circ \pi(e_{K'})\in C(K_{\mini}\backslash K)$ is locally bounded (as a map into $C(K_{\mini}\backslash K)$). Obviously, this map is continuous as a map into $C^{-\infty}(K_{\mini}\backslash K)$ (because $s\mapsto \Lambda_s$ is continuous). Moreover, the map $s\mapsto \Lambda_s\circ \pi(e_{K'})$ in fact takes value in the finite dimensional subspace of functions invariant on the right by $K'$. Since the topologies induced by either $C^{-\infty}(K_{\mini}\backslash K)$ or $C(K_{\mini}\backslash K)$ on that subspace are the same, we conclude that the map $s\mapsto \Lambda_s\circ \pi(e_{K'})\in C(K_{\mini}\backslash K)$ is continuous and we are done.

\item Assume now that $F=\mathbb{R}$. Fix a basis $X_1,\ldots X_N$ of $\overline{\mathfrak{p}}_{\mini}(F)$ and set
$$\displaystyle \Delta_{\overline{P}}=X_1^2+\ldots+X_N^2\in \mathcal{U}\left(\overline{\mathfrak{p}}(F)\right)$$
By elliptic regularity (cf.\ \ref{eq 2.1.2}), for every positive integer $m$ such that $2m>\dim(\overline{P}_{\mini})$, there exist functions $\varphi^1_{\overline{P}}\in C_c^{2m-dim\overline{P}_{\mini}-1}\left(\overline{P}_{\mini}(F)\right)$ and $\varphi^2_{\overline{P}}\in C_c^\infty\left(\overline{P}_{\mini}(F)\right)$ such that
$$\displaystyle \varphi^1_{\overline{P}}\ast \Delta_{\overline{P}}^m+\varphi_{\overline{P}}^2=\delta_e^{\overline{P}_{\mini}}$$
Hence, for all $s\in \mathbb{C}$ we have
$$\displaystyle \pi_s(\varphi^1_{\overline{P}})\pi_s(\Delta_{\overline{P}}^m)+\pi_s(\varphi_{\overline{P}}^2)=Id$$
Choose a function $\varphi_N\in C_c^\infty\left(N_{\mini}(F)\right)$ such that $\displaystyle \int_{N_{\mini}(F)} \varphi_N(n)\xi_{\mini}(n) dn=1$. Then, we have

\begin{align}\label{eq B.3.13}
\displaystyle \Lambda_s(e')=\Lambda_s\left( \pi_s(\varphi^1)\pi_s(\Delta_{\overline{P}}^m)e'\right)+\Lambda_s\left(\pi_s(\varphi^2)e'\right)
\end{align}

\noindent for all $e'\in C^\infty(K_{\mini}\backslash K)$ and all $s\in \mathbb{C}$, where $\varphi^1=\varphi_N\ast \varphi_{\overline{P}}^1\in C_c^{2m-dim\overline{P}_{\mini}-1}\left(G(F)\right)$ and $\varphi^2=\varphi_N\ast \varphi^2_{\overline{P}}\in C_c^\infty\left(G(F)\right)$.

\vspace{2mm}

\noindent Fix $\Omega\subseteq \mathbb{C}$ a compact subset. By \ref{eq B.3.7}, there exists an integer $k\geqslant 0$ such that
$$\displaystyle s\mapsto \Lambda_s\in C^{-k}(K_{\mini}\backslash K)$$
is holomorphic. If $m$ is sufficiently large, we will have $\varphi^i\in C_c^k(G(F))$ and hence $\Lambda_s\circ \pi_s(\varphi^i)\in C(K_{\mini}\backslash K)$ for $i=1,2$ and for all $s\in \Omega$. Henceforth, we will assume that $m$ is that sufficiently large. By \ref{eq B.3.13}, we have
$$\displaystyle \Lambda_s\left(\pi_s(a(x))e\right)=\left\langle \Lambda_s\circ \pi_s(\varphi^1),\pi_s(a(x))\pi_s\left(a(x)^{-1}\Delta_{\overline{P}}^ma(x)\right)e\right\rangle_K+\left\langle \Lambda_s\circ \pi_s(\varphi^2), \pi_s(a(x))e\right\rangle_K$$
for all $s\in \Omega$ and all $x\in F^\times$. Since $a(.)$ is in the positive chamber relative to $P_{\mini}$, the function $x\in F^\times\mapsto a(x)^{-1}\Delta_{\overline{P}}^ma(x)$ is a finite sum of terms of the form $x\mapsto x^\ell D$ where $\ell\geqslant 0$ is an integer and $D\in \mathcal{U}(\overline{\mathfrak{p}}_{\mini}(F))$. Since for all $D\in \mathcal{U}(\mathfrak{g}(F))$ the map $s\mapsto \pi_s(D)e$ is locally bounded, by \ref{eq B.3.12} we are reduced to showing that the maps
$$\displaystyle s\in \Omega\mapsto \Lambda_s\circ \pi_s(\varphi^i)\in C(K_{\mini}\backslash K)\;\;\; (i=1,2)$$
are bounded. Fix $i$ to be $1$ or $2$ and set $\varphi=\varphi^i$. In any case, we have $\varphi\in C_c^k(G(F))$. Then, for all $s\in \mathbb{C}$ the operator $\pi_s(\varphi)$ is given by a kernel function $\pi_s(\varphi)(.,.)$ on $K_{\mini}\backslash K\times K_{\mini}\backslash K$ i.e., we have
$$\displaystyle \left(\pi_s(\varphi)e'\right)(k)=\int_K \pi_s(\varphi)(k,k')e(k')dk'$$
for all $e'\in C^\infty(K_{\mini}\backslash K)$ and all $k\in K$, where
$$\displaystyle \pi_s(\varphi)(k,k')=\int_{\overline{P}_{\mini}(F)} \varphi(k^{-1}\overline{p}_{\mini}k') \delta_{\overline{P}_{\mini}}(\overline{p}_{\mini})^{1/2+s}d_L\overline{p}_{\mini}$$
For all $s\in \mathbb{C}$, we have $\pi_s(\varphi)(.,.)\in C^k(K_{\mini}\backslash K\times K_{\mini}\backslash K)$ and it is not hard to see that the map $s\mapsto \pi_s(\varphi)(.,.)\in C^k(K_{\mini}\backslash K\times K_{\mini}\backslash K)$ is holomorphic hence continuous. Now for all $s\in \Omega$ and all $k\in K$, we have
$$\displaystyle \left(\Lambda_s\circ \pi_s(\varphi)\right)(k)=\Lambda_s\left(\pi_s(\varphi)(k,.)\right)$$
Since $s\mapsto \Lambda_s\in C^{-k}(K_{\mini}\backslash K)$ is continuous, we deduce that

\[\begin{aligned}
\displaystyle \left\lvert \left(\Lambda_s\circ\pi_s(\varphi)\right)(k)\right\rvert & \ll \left\lVert \pi_s(\varphi)(k,.)\right\rVert_{C^k} \\
 & \ll \left\lVert \pi_s(\varphi)(.,.)\right\rVert_{C^k}
\end{aligned}\]
for all $s\in \Omega$ and all $k\in K$ (the first norm above is the norm on $C^k(K_{\mini}\backslash K)$ whereas the second norm is the norm on $C^k(K_{\mini}\backslash K\times K_{\mini}\backslash K)$). Since the last norm above is bounded on $\Omega$ (because $s\mapsto \pi_s(\varphi)(.,.)$ is continuous), this proves that the map $s\in \Omega\mapsto \Lambda_s\circ\pi_s(\varphi)\in C^k(K_{\mini}\backslash K)$ is bounded and ends the proof of \ref{eq B.3.11}. $\blacksquare$
\end{itemize}

\begin{cor}\label{corollary B.3.1}
Let $P=MN\supseteq P_{\mini}=M_{\mini}N_{\mini}$ be a parabolic subgroup of $G$ containing $P_{\mini}$. Then, for all $\delta>0$ there exists $\epsilon>0$ such that the integral
$$\displaystyle \int_{N(F)} \Xi^G(n) e^{\epsilon \sigma(n)} \left(1+\lvert \lambda_{\mini}(n)\rvert\right)^{-\delta} dn$$
is absolutely convergent.
\end{cor}

\vspace{3mm}

\noindent\ul{Proof}: First we do the case where $P=P_{\mini}$. Fix $\epsilon_0>0$ such that the conclusion of Proposition \ref{proposition B.3.1} holds for $\epsilon=\epsilon_0$. Then, the integral
$$\displaystyle \int_{N_{\mini}(F)} \delta_{\overline{P}_{\mini}}(m_{\overline{P}_{\mini}}(n))^{1/2-\epsilon_0}\left(1+\lvert \lambda_{\mini}(n)\rvert\right)^{-2} dn$$
is absolutely convergent. By Proposition \ref{proposition 1.5.1}(ii), there exists $d>0$ such that $\Xi^G(g)\ll\delta_{\overline{P}_{\mini}}(m_{\overline{P}_{\mini}}(g))^{1/2}\sigma(g)^d$ for all $g\in G(F)$. Moreover, by the \cite{Wa2} Lemme II.3.4 (in the $p$-adic case) and \cite{Wall} Lemma 4.A.2.3 (in the real case), there exists $c_1>0$ such that $e^{c_1\sigma(n)}\ll \delta_{\overline{P}_{\mini}}(m_{\overline{P}_{\mini}}(n))^{-1}$ for all $n\in N_{\mini}(F)$. It follows that for $\epsilon< \epsilon_0c_1$, the integral
$$\displaystyle \int_{N_{\mini}(F)} \Xi^G(n) e^{\epsilon  \sigma(n)} \left(1+\lvert \lambda_{\mini}(n)\rvert\right)^{-2} dn$$
is absolutely convergent. This establishes the corollary for $\delta\geqslant 2$. Let $0<\delta<2$ and set $p=\frac{2}{\delta}$, $q=\frac{2}{2-\delta}$. By H\"older inequality, for all $\epsilon>0$ we have

\[\begin{aligned}
\displaystyle \int_{N_{\mini}(F)} \Xi^G(n) e^{\epsilon \sigma(n)} \left(1+\lvert \lambda_{\mini}(n)\rvert\right)^{-\delta} dn & =\int_{N_{\mini}(F)} \Xi^G(n)^{1/p+1/q} e^{2\epsilon \sigma(n)-\epsilon \sigma(n)} \left(1+\lvert \lambda_{\mini}(n)\rvert\right)^{-\delta} dn \\
 & \leqslant \left(\int_{N_{\mini}(F)} \Xi^G(n) e^{2\epsilon p\sigma(n)} \left(1+\lvert \lambda_{\mini}(n)\rvert\right)^{-2} dn\right)^{1/p} \\
 & \times \left(\int_{N_{\mini}(F)} \Xi^G(n) e^{-\epsilon q \sigma(n)}dn\right)^{1/q}
\end{aligned}\]

\noindent By what we just saw, the first integral above is absolutely convergent if $\epsilon$ is sufficiently small. On the other hand, the second integral above is always absolutely convergent by Proposition \ref{proposition 1.5.1}(iv). Hence, the integral
$$\displaystyle \int_{N_{\mini}(F)} \Xi^G(n) e^{\epsilon \sigma(n)} \left(1+\lvert \lambda_{\mini}(n)\rvert\right)^{-\delta} dn$$
is absolutely convergent for $\epsilon>0$ sufficiently small. This settles the case $P=P_{\mini}$. We may deduce the general case from this particular case as follows. Let $\delta>0$ and choose $\epsilon>0$ such that the conclusion of the corollary holds for $P=P_{\mini}$. Then, the integral
$$\displaystyle \int_{N(F)\backslash N_{\mini}(F)} \int_{N(F)} \Xi^G(nn')e^{\epsilon \sigma(nn')} \left(1+\lvert \lambda_{\mini}(nn')\rvert\right)^{-\delta} dndn'$$
is absolutely convergent. By Fubini, it follows that there exists $n'\in N_{\mini}(F)$ such that the inner integral
$$\displaystyle \int_{N(F)} \Xi^G(nn')e^{\epsilon \sigma(nn')} \left(1+\lvert \lambda_{\mini}(nn')\rvert\right)^{-\delta} dn$$
is also absolutely convergent. Up to translation by $N(F)$, we may assume that $\lambda_{\mini}(n')=0$. Moreover, we have inequalities 
$$\Xi^G(n)\ll \Xi^G(nn'),\;\;\; e^{\epsilon \sigma(n)}\ll e^{\epsilon \sigma(nn')}$$
for all $n\in N(F)$. It follows that the same integral without the $n'$ is also absolutely convergent and this ends the proof of the corollary $\blacksquare$

\bigskip

\phantomsection
\clearpage
\addcontentsline{toc}{section}{\protect\numberline{}Bibliography}

\phantomsection
\cleardoublepage
\addcontentsline{toc}{section}{\indexname}

\printnoidxglossary[type=index,style=mcolindex,sort=use]
\end{document}